\documentclass[11pt]{amsart}

\usepackage[T1]{fontenc}
\usepackage{lmodern,tikz-cd}
\usepackage[expansion=false]{microtype}
\usepackage{amsmath,amssymb,mathtools,mathrsfs,fullpage}
\usepackage{stmaryrd}
\SetSymbolFont{stmry}{bold}{U}{stmry}{m}{n}
\usepackage[shortlabels]{enumitem}
\usepackage[pagebackref,colorlinks=true,linkcolor=blue,citecolor=blue,urlcolor=blue]{hyperref}
\renewcommand*{\backref}[1]{}
\renewcommand*{\backrefalt}[4]{%
  \ifnum#1>0\quad#2\fi}

\hypersetup{pdftitle={Loops and mirrors},pdfauthor={Nero Budur}}

\newtheorem{theorem}{Theorem}[section]
\newtheorem{proposition}[theorem]{Proposition}
\newtheorem{lemma}[theorem]{Lemma}
\newtheorem{corollary}[theorem]{Corollary}
\newtheorem{conjecture}[theorem]{Conjecture}

\theoremstyle{definition}
\newtheorem{definition}[theorem]{Definition}
\newtheorem{example}[theorem]{Example}
\newtheorem{remark}[theorem]{Remark}
\newtheorem{convention}[theorem]{Convention}

\AddToHook{cmd/@thm/before}{\setlist[enumerate]{font=\bfseries}}

\newcommand{\A}{\mathbb A}
\def\bA{\mathbb A}
\def\bQ{\mathbb Q}
\newcommand{\C}{\mathbb C}
\newcommand{\F}{\mathbb Q}
\newcommand{\N}{\mathbb N}
\newcommand{\Pj}{\mathbb P}
\newcommand{\R}{\mathbb R}
\newcommand{\Z}{\mathbb Z}
\newcommand{\cL}{\mathcal L}
\newcommand{\vct}{\mathbf v}
\newcommand{\uct}{\mathbf u}
\newcommand{\wct}{\mathbf w}
\newcommand{\zero}{\mathbf 0}
\newcommand{\ord}{\operatorname{ord}}
\newcommand{\Spec}{\operatorname{Spec}}
\DeclareMathOperator{\Cot}{Cot}
\newcommand{\gr}{\operatorname{gr}}
\newcommand{\hocolim}{\operatorname*{hocolim}}
\newcommand{\Gm}{\mathbb G_m}

\newcommand{\Hom}{\operatorname{Hom}}

\newcommand{\supp}{\operatorname{supp}}
\newcommand{\Tot}{\operatorname{Tot}}

\newcommand{\fin}{\mathrm{fin}}
\def\sL{\mathscr L}

\title{Loops and mirrors}
\author{Nero Budur}
\address{Department of Mathematics, KU Leuven, Celestijnenlaan 200B,
3001 Leuven, Belgium; YMSC, Tsinghua University, 100084 Beijing, China;
BCAM, Mazarredo 14, 48009 Bilbao, Spain}
\email{nero.budur@kuleuven.be}
\date{\today\ at \begingroup
  \count255=\time
  \divide\count255 by 60
  \ifnum\count255<10 0\fi\number\count255:%
  \multiply\count255 by 60
  \advance\count255 by -\time
  \multiply\count255 by -1
  \ifnum\count255<10 0\fi\number\count255\endgroup}
\makeatletter
\renewcommand*\l@section{\@tocline{1}{0pt}{1pc}{}{\small}}
\makeatother

\begin{document}

\begin{abstract}
For a smooth affine complex variety $X$ we construct a $\bQ$-algebra $LH^*(X)$, loop cohomology, in terms of the formal-loop space of $X$. We conjecture that it is isomorphic to the symplectic cohomology of $X$.
The degree-zero loop cohomology subalgebra is finitely generated for polarized log Calabi-Yau $X$, constructing an intrinsic mirror affine coordinate ring. Its polyvector algebra is identified with the full loop cohomology algebra for a positive log Calabi--Yau pair with maximal boundary and smooth loop mirror, proving a form of mirror symmetry. In this case we further identify the degree-zero loop cohomology algebra with the Gross--Siebert intrinsic mirror algebra. Assuming in addition that $X$ is a surface, we prove that the full loop and symplectic cohomology algebras agree.
\end{abstract}

\maketitle
\begingroup\small
\tableofcontents
\endgroup
\clearpage

\section{Introduction}
\addtocontents{toc}{\protect\setcounter{tocdepth}{1}}

String theory has been a major source of geometric structures
in mathematics. For a smooth manifold, the free loop space
 provides a natural setting for
some of these structures: Chas--Sullivan string topology
studies operations on its homology, while Hamiltonian Floer
theory develops an infinite-dimensional Morse theory for
action functionals on loop spaces. In algebraic geometry,
related structures are often indirectly constructed through moduli spaces
of stable maps and their virtual fundamental classes, as replacements of moduli of Floer trajectories. 

We construct a direct counterpart of string topology and Floer theory via formal loops on a smooth affine complex variety $X$. A formal loop, and henceforth simply a {\it loop}, on the complex affine space $\bA^n$ is an assignment of a Laurent power series in $\C(\!( t )\!)$ for each of the $n$ coordinates. A loop on a subvariety of $\bA^n$ is a loop on $\bA^n$ satisfying the equations of the subvariety.
 From the loop space $\sL X$ of a smooth affine complex variety $X$ of dimension $> 0$, we
construct a  $\bQ$-algebra $LH^*(X)$ which we call formal-loop cohomology, henceforth simply {\it loop cohomology}. It is endowed with a  multiplication, which we call fusion product,
 defined by collisions of formal loops,  and with a BV operator induced
by loop rotation.
The  grading  is $\Z/2$. If a nowhere-vanishing algebraic volume form $\Omega$ on $X$ is chosen,
it updates to a
 $\Z$-graded version
 $LH^*_\Omega(X).$
 We conjecture that there is an isomorphism between loop cohomology  and symplectic
cohomology $SH^*(X):=SH^*(X,\bQ)$ with the fusion product corresponding to the pair-of-pants product. 

A
polarized simple-normal-crossings compactification $X=M\setminus D$
gives a strongly convergent multiplicative spectral sequence with first page the
logarithmic cohomology of $(M,D)$ and abutment
$LH^*_\Omega(X)$. This an analog of the symplectic spectral sequence of
 Ganatra and Pomerleano \cite{GP}. 
We prove loop analogs of degeneracy criteria of \cite{GP, PomerleanoIntrinsic}.  

For a polarized log Calabi--Yau complement $X=M\setminus D$, the algebra $SH^0(X)\otimes\C$ is supposed be the global regular functions on the mirror of $X$. The degree-zero loop algebra $LH^0_\Omega(X)$ is finitely generated, geometrically reduced, and comes with
a distinguished additive basis.  If $D$ has zero-dimensional strata and $\check X:=\Spec (LH^0_\Omega(X)\otimes \C)$ is smooth, we prove that the polyvector algebra of $\check X$ recovers the full loop algebra $LH^*_\Omega(X)\otimes\C$. This is one form of intrinsic mirror symmetry. We identify in this case  $LH^0_\Omega(X)$ with the intrinsic mirror algebra of Gross--Siebert \cite{GSIntrinsic}, although our construction in terms of loops is very different. Assuming in addition that $X$ is a surface, we prove in this case that $LH^*_\Omega(X;\C)\cong SH^*(X;\C)$ as unital graded algebras.

We give now a more detailed version of these results. 

\subsection{Loop cohomology and its algebra} 
A \emph{polarized SNC compactification} of $X$ is a triple
$(M,D,A)$ consisting of a smooth projective complex variety $M$,
an identification $X=M\setminus D$ with the complement of a
reduced simple-normal-crossings divisor
$D=\bigcup_{i\in I}D_i$, and an ample divisor
\[
 A=\sum_{i\in I}k_iD_i,\qquad k_i\in\Z_{>0}.
\]
Here the $D_i$ are the irreducible components of $D$; in
particular, the support of $A$ is exactly $D$. 

Fix such a polarized SNC compactification $(M,D,A)$. A loop in $X$ gives an arc in $M$ with contact orders $\vct=(v_i)_{i\in I}$ along the boundary $D$, and total weighted contact
$w(\vct)=\sum_i k_iv_i$ with $A$. The contact filtration, or the pole-stage filtration, $F_P\mathscr L X$ consists of loops whose
extension to $M$ has weighted boundary contact at most $P$.

Ideally, one  constructs a loop complex $B$ giving $LH^*(X)$ from
Verdier dualizing complexes on finite-type varieties obtained via coefficient-truncation of these infinite-dimensional pole
stages, using proper pushforward along their closed inclusions and normalized
smooth pullback under coefficient-truncation maps. However, the  coefficient-truncation maps need not always be smooth or have affine-space fibers.
This issue is bypassed using the
finite-factor theorem of Drinfeld \cite{Drinfeld}, which provides suitable local
finite-factor presentations of $F_P\sL X$ by products $Y\times \bA^J$ of a finite-type scheme times an affine space with
countably many coordinates. 
Comparison morphisms between such charts are realized in zig-zags on finite-factors by smooth morphisms with affine-space fibers.
The dimension of the finite-factors $Y$ is normalized in a compatible way by Drinfeld's volume indices, which we show to admit an explicit formula in terms of contact orders at infinity plus shifts in terms of $\Omega$. Using normalized
smooth pullbacks under finite-factor refinements, a complex $B_P$ is defined for each pole stage $F_P\sL X$. This is not really new, it is the Betti version with $\bQ$-coefficients of the complex  with \(\mathbb Z_\ell\)-coefficients introduced by Drinfeld before taking the linear dual in \cite[(6.5)]{Drinfeld}. It should be thought of as a renormalized Borel-Moore complex, and its cohomology as a renormalized Borel-Moore {\it homology}. While \cite{Drinfeld} works with ordinary derived category of complexes over $\bQ$, we work with the derived $\infty$-category, so that in absence of $\Omega$ we can still construct $B$ with a $\Z/2$-grading. Then $B$ is defined as the homotopy colimit of $B_P$ along the closed inclusions $$F_0\sL X\subset\ldots\subset F_P\sL X\subset F_{P+1}\sL X\subset\ldots\subset \sL X=\varinjlim_{P}F_P\sL X$$ and  $LH:=H(B)=\varinjlim_{P}H(B_P)$ is the colimit of the renormalized Borel-Moore homology of these poles stages. 

We then use the families of loops from the factorization monoid structure on $\sL X$ constructed by Kapranov-Vasserot \cite{KV} to model collisions of two or more loops. In a simplified 1-parameter case, such a model fits into a Cartesian diagram
\[
\begin{tikzcd}[column sep=large, row sep=large]
\Gm \times F_P\sL X \times F_Q\sL X \times \bA^J
  \arrow[r, hook, "j"]
  \arrow[d, "\operatorname{pr}_1"']
&
\mathcal Y
  \arrow[d, "\pi"]
&
F_{P+Q}\sL X \times \bA^I
  \arrow[l, hook', "i"']
  \arrow[d]
\\
\Gm
  \arrow[r, hook, "j_0"']
&
\A^1
&
\{0\}
  \arrow[l, hook', "i_0"]
\end{tikzcd}
\]
where $I$ and $J$ are countable. Specialization of Borel-Moore cycles on the finite-factor models of these collision families, renormalized for our infinite-dimensional setup using volume indices, define the fusion product first at the level of pole stages $B_P\otimes B_Q\to B_{P+Q}$, then via homotopy colimit on $B$, and finally by taking cohomology  on $LH$. 
Rotating the source of a loop via
$t\mapsto \lambda t$ for $\lambda\in\C^*$ preserves the pole stages $F_P\sL X$. We extract from this action a degree-$-1$ operator $\Delta$ at each pole stage, which gives us one on $B$, and it induces one on $LH$.
Properties of the fusion product and of the operator $\Delta$, like associativity and BV relations, are proved by analizing the renormalized Borel-Moore specialization associated with collisions with three or more input loops.

\begin{theorem}
\label{thm:summary-fusion}
For every smooth affine complex variety $X$, loop cohomology
$LH^*(X)$ is a well-defined parity-graded
$\bQ$-vector space with a canonical fusion product
\[
 \mu:LH^i(X)\otimes LH^j(X)
                  \longrightarrow LH^{i+j}(X).
\]
This product is associative, unital, and graded-commutative.
It is endowed with an operator
$\Delta:LH^j(X)\to LH^{j-1}(X)$
which makes $(LH^*,\mu,\Delta)$ a unital BV algebra.

Each polarized SNC compactification $(M,D,A)$ gives an exhaustive
nonnegative contact filtration $F_\bullet LH^*$, with
\[
 F_{-1}LH^*=0,\qquad
 \mu(F_PLH^*,F_QLH^*)\subseteq F_{P+Q}LH^*,\qquad
 \Delta(F_PLH^*)\subseteq F_PLH^*.
\]
There is a natural map
$\iota_X:H^*(X;\bQ)\to LH^*(X)$, with ordinary
cohomological degrees read modulo two, such that
\[
 F_0LH^*=\operatorname{im}(\iota_X),\qquad
 1_{LH}=\iota_X(1).
\]
After forgetting the filtration, the BV algebra is
independent of the polarized compactification. 

A chosen nowhere-vanishing algebraic volume form $\Omega$ on $X$ gives
an integer-graded algebra $LH^*_\Omega(X)$ with the
same properties, whose parity reduction recovers $LH^*(X)$.
\end{theorem}

The construction of the loop complex $B$ and of the loop cohomology is done in Section~\ref{sec:loop-complex}. Section~\ref{sec:fusion-construction} constructs the fusion product, and Section~\ref{sec:incidence} proves the associative unital graded-commutative version of the above theorem. Section \ref{sec:bv} introduces the BV operator and proves the rest of the theorem.

The terminology loop cohomology has appeared in different contexts, e.g. for topological loops and derived analogs of them. We stress that our terminology is short for formal-loop cohomology and that this is a genuinely new structure.

\subsection{The multiplicative contact spectral sequence} The successive differences $F_{P+1}\sL X\setminus F_P\sL X$ are topological disjoint unions of loci of arcs on $M$ with fixed contact with $D$. The latter have simple topological description since the boundary $D$ is SNC. This leads to:

\begin{theorem}
\label{thm:summary}
For each polarized SNC compactification $(M,D,A)$, the pole-stage
filtration has a spectral sequence of unital BV algebras
converging strongly as BV algebras:
\begin{equation}\label{eq:loop-ss}
 E_1^{p,q}=H^{p,q}_{\log}(M,D)
       \quad\Longrightarrow\quad LH^{p+q}(X).
\end{equation}
Here $q$ and total degree are read modulo two, and the abutment
has the filtration of Definition~\ref{def:loop-cohomology}.
A chosen nowhere-vanishing algebraic volume form $\Omega$
gives the same assertions in integer grading with abutment
$LH^*_\Omega(X)$.
\end{theorem}

Here the first-page is the logarithmic cohomology algebra of \cite{GP}, defined in Section \ref{sec:setup}. The logarithmic BV operator is defined in  Subsection \ref{subsec:logarithmic-bv}. The theorem is proved in
Section~\ref{sec:spectral-sequence}. This is a loop analog of the spectral sequence of \cite[Theorem~1.1]{GP} for symplectic cohomology and naturally leads to:

\begin{conjecture}[Loop--Floer Conjecture]
\label{conj:intro-multiplicative}
There is an isomorphism of filtered graded BV algebras
\[
 \Phi:(LH^*(X),\mu,\Delta)
       \xrightarrow{\sim}(SH^*(X;\F),\star,\Delta_{SH}),
\]
together with an isomorphism of their multiplicative spectral
sequences commuting with the induced BV operators on every page.
Its first-page map is the canonical logarithmic identification,
and its limit map is $\gr\Phi$. 
The statement holds in parity and in the integer grading specified
by any chosen volume form $\Omega$.
\end{conjecture}

The symplectic cohomology of $SH^*(X;\F)$ of $X$
with its associated Liouville structure, originates
in the constructions of Floer--Hofer
\cite{FloerHoferSH} and Cieliebak--Floer--Hofer
\cite{CieliebakFloerHoferSH}; the Liouville-domain version used
here is due to Viterbo \cite{ViterboSH}.
Its product $\star$ counts Floer pair-of-pants solutions,
and $\Delta_{SH}$ is the degree-$-1$ operator obtained by rotating
an asymptotic marker; see \cite[Section~8(a)]{SeidelBiased}.

A more general conjecture would
ask for a filtered derived equivalence between the loop and
symplectic cochain complexes intertwining their products
and BV operators. Conjecture \ref{conj:intro-multiplicative} is the analog of the Arc-Floer Conjecture of \cite{BBLN} between arc contact loci and fixed-point Floer theory of monodromy iterates of a germ of holomorphic function; see also \cite{BodegaLorenzo}. Theorem \ref{thm:summary} is a richer loop analog of the main result on arcs of \cite{BBLN}.

There are computational simplifications in presence of 
global finite-factor presentations of pole stages, see  Subsection~\ref{subsec:global-finite-presentations} for $LH$ additively, Subsection~\ref{subsec:global-binary-fusion} for the fusion product on $LH$, and Subsection~\ref{subsec:global-finite-bv} for
the BV operator. These express the computations
in terms of finite-dimensional Borel--Moore homology.
A K\"unneth theorem and deformation invariance for $LH$ are proved in 
Section~\ref{sec:conclusion} and provide further reductions. 
We work out examples of loop cohomology algebras in
Part~\ref{part:examples}. 

\subsection{Degeneration in all degrees}
Our main degeneration criterion for the above spectral sequence is the loop counterpart of
  \cite[Corollary~1.2]{GP}.

\begin{theorem}[Degeneration for proportional boundary classes]
\label{thm:intro-proportional}
Let $(M,D,A)$ be a polarized SNC compactification of $X$, with
$D=\bigcup_{i=1}^rD_i$. Suppose that a $\Z$-Cartier divisor
$H$ and positive integers $n_i$ satisfy $D_i\sim n_iH$ for every $i$.
If some $n_i>1$, or if all $n_i=1$ and
$D_1\cap\cdots\cap D_r=\varnothing$, the contact BV spectral
sequence \eqref{eq:loop-ss} degenerates at $E_1$ in the sense
of Definition~\ref{def:bv-spectral-degeneration}.
This holds in parity grading and in the integer grading specified
by any nowhere-vanishing algebraic volume form $\Omega$ on $X$.
\end{theorem}

The theorem determines the associated graded ring, but products
on $LH$ may still have terms of smaller contact weight. 
Its proof, given also in
Section~\ref{sec:spectral-sequence} starts with noting that the proportionality $D_i\sim n_iH$ makes the
class of the contact vector $\vct$ in $\Z^r/\Z(n_1,\ldots,n_r)$
locally constant on the loop space. The first-page algebra $E_1$
is generated by contact vectors with entries $0$ or $1$.
Under the hypotheses of Theorem~\ref{thm:intro-proportional},
subtracting any positive multiple of $(n_1,\ldots,n_r)$
from such a realized vector produces a negative entry.
Thus these generators have no possible lower-contact targets
for a differential and the Leibniz rule for the differentials of the multiplicative spectral sequence then forces degeneration.

\subsection{Foundations}
Part~\ref{part:general} of this paper is dedicated to the construction of loop cohomology and proving the above results. Let us explain a bit the length of this foundational part. As mentioned, formal loop spaces are studied here through local finite-factor presentations of pole stages $F_P\sL X$ by
finite-dimensional schemes with auxiliary affine coordinates. 
Coherent descent glues their normalized complexes and the maps
induced by geometric constructions, retaining the homotopies
and higher compatibilities required on multiple overlaps.
This allows Borel-Moore specialization, fusion, and source rotation to
define global operations independent of the chosen finite-factor presentations. This requires simultaneous finite-dimensional models for the
relevant local diagrams. On finite-type varieties one could use the enhanced six functor formalism \cite{ScholzeSixFunctors} as a black-box. However, we also need to carry one operation, $p^{-1}$, where $p$ is the smooth projection away from countably many auxilliary affine coordinates, which falls outside this formalism. 
To handle this, under smooth changes of finite model,
the normalization must cancel the additional dimension shifts,
and the resulting comparisons must intertwine the maps of
complexes compatibly with base change and further refinement. Much of Part~\ref{part:general} establishes these geometric compatibilities,
allowing the functorial sheaf constructions to supply the
higher coherences needed for gluing. 

The foundations for this hybrid finite-infinite dimensional setup are not fixed or standardized. This would be desirable in the future, but for now we laboriously carry similar arguments from one step to another when it comes to gluing in enhanced derived cateories and deal exclusively in our context with the technical difficulties regarding the normalization process.
Regarding the latter, in order to normalize Borel-Moore specialization to our infinite-dimensional setup and thus define the fusion product on loop cohomology, we needed to extend Drinfeld's finite-factor theorem to a relative setup, and give a precise formula in terms of contact orders along the boundary at infinity for his volume index used in the normalization process.

\subsection{Degree zero for log Calabi--Yau pairs}
\label{subsec:intro-positive-lcy}
A \emph{positive log Calabi--Yau pair} is a smooth connected
projective variety $M$ with reduced SNC boundary $D$ such that
$K_M+D$ is trivial and $D$ supports an ample divisor with strictly
positive coefficients; it is {\it polarized} if such an ample divisor $A$  has been chosen.
Its complement $X=M\setminus D$ is affine.
The contact spectral sequence has a further consequence for such
a pair. Fix a volume form with
$\operatorname{div}_M\Omega=-D$.

\begin{theorem}[Degree-zero degeneration for log Calabi--Yau pairs]
\label{thm:intro-lcy-zero}
For every polarized log Calabi--Yau pair $(M,D,A)$,  the spectral sequence of
Theorem~\ref{thm:summary} degenerates at $E_1$ in total degree
zero, and
\[
 \gr_FLH^0_\Omega(X)\cong H^0_{\log}(M,D)
\]
as rings graded by the contact weight
$w(\vct)=\sum_i k_iv_i$ associated with the polarization
$A=\sum_i k_iD_i$
(Definition~\ref{def:log-cohomology}), and
$LH^j_\Omega(X)=0$ for $j<0$.
\end{theorem}

The proof is given in Section~\ref{sec:log-cy-zero} and it almost automatic at this point.
Recall from the construction of $LH$ that volume indices are used to normalize dimensions of the finite factors in a nice atlas for poles stages. The volume-index function admits an explicit expression in terms of local dimension and a term in terms of contact orders at infinity contributed to by $\Omega$, see Proposition \ref{prop:general-factor-contact-index}.
The log Calabi--Yau condition makes this last term automatically zero, so it
gives
equidimensionality: on a finite-factor chart with constant
volume index $d$, every nonempty closed pole-stage factor
is pure of dimension $d$, with irreducible components given
by reduced contact-piece closures. Their fundamental classes
have normalized degree zero and glue to lifts of all
degree-zero first-page classes. These lifts force outgoing
differentials to vanish, while the absence of negative degrees
excludes incoming differentials. 

Theorem~\ref{thm:intro-lcy-zero} is the loop counterpart of the
Floer degree-zero degeneration theorem of Pomerleano \cite[Theorem~1.2]{PomerleanoIntrinsic}, which completed partial degree-zero degeneration results of \cite{GP,Pascaleff}.  The loop and Floer  degree-zero degeneration
results give, in this full positive log Calabi--Yau range,
\[
 \gr_FLH^0_\Omega(X)
 \cong H^0_{\log}(M,D)
 \cong\gr_F SH^0_\Omega(X;\F)
\]
as contact-weight-graded unital $\F$-algebras.
A full degree-zero ring comparison does not
follow from it since lower-weight product terms remain to be
compared. 

\subsection{The loop theta basis and the intrinsic loop affine mirror}\label{subsLa} Write $D=\bigcup_{i=1}^rD_i$,
and, for a boundary contact vector $\vct\in\Z_{\ge0}^r$, put
\[
 I(\vct)=\{i:v_i>0\},\qquad
 D_I^\circ=\left(\bigcap_{i\in I}D_i\right)
              \setminus\bigcup_{j\notin I}D_j
              \quad(I\ne\varnothing).
\]
Write $\pi_0(T)$ for the set of connected components of $T$, and put
\[
 \Lambda_D=\{0\}\sqcup
 \{(\vct,C):\vct\in\Z_{\ge0}^r\setminus\{0\},\quad
                    C\in\pi_0(D_{I(\vct)}^\circ)\},
\]
where $0$ labels arcs in $X$. For an ample divisor $A=\sum_i k_iD_i$ with
$k_i\in\Z_{>0}$ put 
$$ w_A(\vct,C)=\sum_i k_iv_i,
 \qquad w_A(0)=0.$$

\begin{theorem}[The distinguished theta basis]
\label{thm:intro-theta-basis}
Let $(M,D)$ be a positive log Calabi--Yau pair, put
$X=M\setminus D$, and choose $\Omega$ with
$\operatorname{div}_M\Omega=-D$.  There is a canonical basis
$\{\vartheta_p\}_{p\in\Lambda_D}$ such that
\[
 LH^0_\Omega(X)=\bigoplus_{p\in\Lambda_D}\F\vartheta_p,
 \qquad \vartheta_0=1.
\]
For the fixed pair and volume form, the classes are independent
of the boundary polarization. 
For every ample divisor $A=\sum_i k_iD_i$ with
$k_i\in\Z_{>0}$ and every $P\in\Z_{\ge0}$, they satisfy
\[
 F^A_PLH^0_\Omega(X)
   =\bigoplus_{w_A(p)\le P}\F\vartheta_p.
\]
Under Theorem~\ref{thm:intro-lcy-zero}, their symbols are the
constant contact classes $e_p\in H^0_{\log}(M,D)$ of
Definition~\ref{def:lcy-constant-contact-classes}:
\[
 \operatorname{in}_{F^A}(\vartheta_p)=e_p
 \quad\text{in weight }w_A(p).
\]
Here $F^A$ denotes the contact filtration determined by $A$, and
$\operatorname{in}_{F^A}(\vartheta_p)$ is the image of
$\vartheta_p$ in
$F^A_{w_A(p)}LH^0_\Omega(X)/F^A_{w_A(p)-1}LH^0_\Omega(X)$.
\end{theorem}

We call $\{\vartheta_p\}$ the \emph{(loop) theta basis}. The theorem is proved in Section~\ref{sec:log-cy-zero} and uses the same construction  behind Theorem \ref{thm:intro-lcy-zero} mentioned above. So its proof is also an automatic consequence of the volume-index formula.
For each contact label $p$, consider the finite-factor charts of
$F^A_{w_A(p)}\sL X$. In each chart, sum the normalized fundamental
classes of the reduced closures of all contact pieces labeled
by $p$, each with coefficient one. These are degree-zero classes by the equidimensionality mentioned above, 
agree under smooth comparisons, and descend to a class at that
pole stage. The image is in $LH^0_\Omega(X)$ is $\vartheta_p$.

Put
\[
 R_X=LH^0_\Omega(X),\qquad
 \check X=\Spec(R_X\otimes_{\F}\C).
\]
The algebra $H^0_{\log}(M,D)$ is identified already in \cite{GP, PomerleanoIntrinsic} with the boundary face algebra, with a separate face for every connected
component of a boundary intersection; the so-called Stanley-Reisner algebra.
The above results imply that the contact filtration gives a flat Rees degeneration of $R_X$ to the
Stanley-Reisner algebra. Properties of the latter together with standard filtered commutative algebra imply that $R_X$ is finitely generated and geometrically reduced,
and that
\[
 \dim \check X
 =\max\{\lvert I\rvert:D_I=\bigcap_{i\in I}D_i\ne\varnothing\}.
\]
 Thus this $\check X$
has dimension $n=\dim X$ exactly under the \emph{maximal boundary}
condition: $D$ has a zero-dimensional stratum. So finite-generation is a triviality for us, whereas for other approaches to intrinsic mirror constructions it is a sensitive point, equivalent to showing polynomiality of the multiplication in the theta bases. 

\begin{theorem}[Calabi--Yau property of the loop mirror]
\label{mir:thm-intro-cy}
\label{mir:thm-loop-cy}
For a positive log Calabi--Yau pair $(M,D)$ with maximal boundary,
$\check X$ is a connected semi-log-canonical Gorenstein
scheme of pure dimension $n$, with
$
 \omega_{\check X}\simeq
 \mathcal O_{\check X}.
$
In particular its smooth locus has a nowhere-vanishing algebraic
$n$-form.  If $\nu:\check X^\nu\to \check X$ is the normalization
and $C$ its reduced conductor divisor, then $(\check X^\nu,C)$ is log
canonical and $K_{\check X^\nu}+C\sim0$. If $\check X$ is normal,
it is therefore an integral log-canonical Calabi--Yau variety.
\end{theorem}


This is proved in Section \ref{mir:sec-affine-scheme}. 
It also almost automatic. It follows from properties of the dual complex of a log Calabi-Yau pair via standard filtered commutative algebra once the isomorphism from Theorem \ref{thm:intro-lcy-zero} is known, there is no additional loop input needed for it.
After blowing up boundary strata, the contact filtration gives
a flat degeneration of $\check X$ to affine spaces glued along
coordinate subspaces according to the boundary dual complex.
The homology-sphere properties of this complex and its links
make the central fiber Gorenstein and semi-log-canonical.
Compatibly oriented logarithmic volume forms have opposite
residues along shared coordinate hyperplanes and glue to
trivialize its dualizing sheaf. The graded Rees construction
lifts this trivialization and these singularity properties
to $\check X$. Under the hypotheses of the theorem, we expect, but do not show, that $\check X$ is also normal.
In the expected mirror picture,
$\check X$ is the mirror itself when the full mirror is
affine, and its affinization otherwise. 

\subsection{The polyvector algebra of a smooth loop mirror}
Another intrinsic mirror certificate is the following.

\begin{theorem}[The global polyvector comparison]
\label{thm:intro-global-polyvectors}
Let $(M,D)$ be a positive log Calabi--Yau pair of dimension $n\ge1$.
Put $X=M\setminus D$ and choose $\Omega$ with
$\operatorname{div}_M\Omega=-D$. Suppose $D$ has a zero-dimensional
stratum and $\check X$ is smooth. There is a canonical isomorphism of graded BV algebras
\[
 \sigma:LH^*_\Omega(X)\otimes_{\F}\C
 \xrightarrow{\sim}
 \Gamma\bigl(\check X,\wedge^\bullet T_{\check X}\bigr).
\]
 In particular, the complexified loop
algebra is generated over $R_X\otimes_{\F}\C$ by its degree-one part;
its degree-$j$ part is a projective module of rank $\binom nj$ over
that ring, and vanishes outside $0\le j\le n$.
\end{theorem}

On the right-hand side the BV operator is the divergence $\Delta_\eta$ associated to the mirror volume form $\eta$ of Definition \ref{globpv:def-mirror-volume}.  The proof again uses no new loop computations, so it is in some sense also automatic. It is essentially an algebraic consequence of the filtered BV spectral sequence and degree-zero degeneration from above. The BV bracket lets degree-one loop classes
differentiate functions on $\check X$, while iterated brackets
define higher polyvectors. The de Rham class $[\Omega]$,
mapped into loop cohomology through the arc locus, produces
a top polyvector. Boundary residues and the Rees degeneration
show that it is nowhere vanishing, and its contractions
generate all vector fields. Since source rotation annihilates
this class, the loop BV operator becomes divergence with
respect to the dual mirror volume form. The theorem is addressed in Section \ref{globpv:comparison}.

\subsection{Comparison with the Gross--Siebert mirror}
The theta basis of Theorem~\ref{thm:intro-theta-basis}
defines  a linear  function, called the \emph{loop trace},
\[
 \epsilon:R_X\longrightarrow\F,\qquad
 \epsilon(\vartheta_p)=\delta_{p,0}.
\]
Here $p=0\in \Lambda_D$ labels the arc component as in \ref{subsLa}  and $\vartheta_0=1$. So $\epsilon$ amounts to taking the ``constant coefficient'' in the theta basis.

 Gross and Siebert \cite{GSIntrinsic}
construct an algebra $R^{\mathrm{GS}}_{M,D}$ with theta basis
$\vartheta_p^{\mathrm{GS}}$, using punctured logarithmic curves.
For a positive pair its structure constants are polynomials in
 variables $q^\beta$ with $\beta$ in a curve-class monoid $Q$, see Section~\ref{mir:sec-gs-comparison}. Define the base chage
\[
 R^{\mathrm{GS},1}_{M,D}
 =R^{\mathrm{GS}}_{M,D}\otimes_{\F[Q],\,q^\beta\mapsto1}\F.
\]

\begin{theorem}[The Gross--Siebert comparison and logarithmic descendants]
\label{mir:thm-intro-gs}\label{mir:conj-intro-enumerative}
Let $(M,D)$ be a positive log Calabi--Yau pair with maximal boundary,
let $X=M\setminus D$ be connected, and suppose its loop mirror
$\check X$ is smooth. Choose $\Omega$ with
$\operatorname{div}_M\Omega=-D$.
The pairing $B(a,b)=\epsilon(ab)$ on $R_X$ is nondegenerate, and
the theta-basis identification is a canonical isomorphism of
unital algebras with trace:
\[
 \Phi_{\mathrm{GS}}:R_X\xrightarrow{\ \sim\ }
 R^{\mathrm{GS},1}_{M,D},\qquad
 \vartheta_p\longmapsto\vartheta_p^{\mathrm{GS}}.
\]
For every $d\ge2$ and nonzero contact labels $p_1,\ldots,p_d$,
\[
 \epsilon(\vartheta_{p_1}\cdots\vartheta_{p_d})
 =\sum_\beta
 \int_{[\overline{\mathcal M}^{\log}_{0,d+1}
 ((M,D),\beta;p_1,\ldots,p_d,0)_x]^{\mathrm{vir}}}
 \psi_{\mathrm{out}}^{d-2}.
\]
The sum is finite. Here $x\in X$ is general, the inputs are labeled,
the output has contact zero and maps to $x$, and
$\psi_{\mathrm{out}}$ is its stable-map cotangent class; the definition of the right-hand side is recalled in Section \ref{mir:sec-gs-comparison}.
\end{theorem}

The proof of nondegeneracy of the two-input loop trace is given in Section \ref{mir:sec-trace-potentials}. It uses first the log Calabi-Yau condition as an obstruction to $\bA^1$-uniruledness to show that arcs with sufficiently generic centers on $X$ have arc neighborhoods in each pole stage $F_P\sL X$. Together with the equidimensionality of the finite-factors mentioned above, it implies that $\epsilon\circ \Delta=0$ on $LH^1_\Omega(X)$. General algebra facts then kick in: Theorem~\ref{thm:intro-global-polyvectors}  implies that the radical of the two-input pairing is stable under all derivations, and differential-simplicity of smooth affine varieties implies triviality of the radical.

The proof of the identification with the Gross-Siebert algebra is given in Section~\ref{mir:sec-gs-comparison}. As in \cite{JohnstonComparison}, both algebras are reconstructed from their respective two-input and three-input traces if the two-input trace is nondegenerate. To compare the two-input and three-input loop traces with theirs, we extract global curves
from loop collisions, which we explain now in a simplified context. Collide marked points
$z=sc_i$ on the $z$-line, where the distinct constants $c_i$ specify a general
marked configuration. Maps from the formal neighborhood of the union of these marked points are the loops we set to collide as $s$ becomes $0$. In our collision family, such an ordered tuple is a point, and letting $s$ go to $0$ we follow a collision branch whose limit at $s=0$ is another loop. Assume this output loop is an arc; this is the relevant type of output loop to look at if we want to compute loop traces. The collision
branch then gives a rational map from a formal surface to $M$. This is the loop version of a {\it wedge}, a fundamental tool for Nash-type problems describing adjacencies between arc contact loci, see for example \cite{EmbeddedNash,BBLN,BobadillaPePereiraNash,deFernexDocampoNash}. Blowing up the ideal $(s,z)$ on this surface exposes an
exceptional projective line with coordinate $w=z/s$:
the input markings lie at $w=c_i$, and the disk carrying the
output arc attaches at $w=\infty$. Resolving the map and
making the source semistable, after a finite base change,
produces a proper tree of rational curves. Replacing its
attachment to the output disk by a marking and stabilizing
gives a global logarithmic map with the prescribed input
contacts and output contact zero at the arc's center.
The loop
traces are shown to be counts of these curves. Using Johnston's enumerativity
theorem \cite{JohnstonComparison} and the 
punctured log trace identities of \cite{GSIntrinsic} we identify the two-input
and three-input traces: the loop versions agree with the Gross-Siebert versions. This comparison does not require
smoothness of $\check X$.  The smoothness hypothesis is only used as explained above to show nondegeneracy of the two-input loop trace. The equivalence with the Gross-Sibert two-input trace shows that the latter is also nondegenerate.

\subsection{Loop and symplectic cohomology for log Calabi--Yau surfaces}
\label{subsec:intro-surface-sh}
The polyvector and the Gross-Sibert comparisons above, together with existing results for surfaces, prove the following case of the
Loop-Floer Conjecture~\ref{conj:intro-multiplicative}.

\begin{theorem}[Surface loop--Floer comparison]
\label{thm:surface-loop-floer}
With the assumptions as in Theorem \ref{mir:thm-intro-gs}, assume in addition that $\dim M=2$.
There is an isomorphism of unital $\Z$-graded $\C$-algebras
\[
 LH^*_\Omega(X)\otimes_{\F}\C
       \xrightarrow{\ \sim\ }SH^*_\Omega(X;\C)
\]
taking the fusion product to the pair-of-pants product. Here the
standard finite-type Weinstein completion of $X$ is graded by
$\Omega|_X$. Both algebras identify with
$\Gamma(\check X,\bigwedge^\bullet T_{\check X})$, with exterior
multiplication and $j$-vectors in degree $j$.
\end{theorem}

The proof is given in Section~\ref{sec:surface-loop-floer}. 
A collection of results  on mirror symmetry  \cite{Varolgunes, HackingKeating,LaiZhou, JohnstonComparison} 
imply that the specialization $q^\beta=1$ of the complexified Gross-Siebert ring give $SH^0(X;\C)$ and that its polyvector algebra gives $SH^*(X;\C)$. Thus Theorem \ref{thm:surface-loop-floer} follows by direct application of 
Theorems~\ref{thm:intro-global-polyvectors} and~\ref{mir:thm-intro-gs}.

\subsection{How this project arose.}  For a germ of a holomorphic function $f$ with an isolated singularity, McLean \cite{McLean} constructed a spectral sequence converging to the fix-point Floer cohomology $HF^*(\phi^m_f,+)$ of iterates of the monodromy, with $E_1$ in terms of a log resolution of $f$. In \cite{BBLN}, we constructed a spectral sequence with the same $E_1$ page, converging to $H^*_c(\mathscr{X}_m(f),\Z)$, where $\mathscr{X}_m(f)$ are the (jet truncations of) the locus of formal-arcs with contact order $m$ along $f$ at the singularity. We posed in \cite{BBLN} the Arc--Floer conjecture, that the two types of cohomologies agree. One of the most exciting prospects, if this conjecture would be true, is that there should be a pair-of-pants product on contact loci. There are however not many calculations of pair-of-pants products on $HF^*(\phi^m_f,+)$ available, making an algebraic counterpart difficult to guess.

Symplectic cohomology is a much more studied Floer algebra on the other hand and we hoped it sheds more light on the previous product. There was a striking similarity with the arc contact loci that fired the spark for this project: Ganatra and Pomerleano \cite{GP} constructed a multiplicative spectral sequence converging to $SH^*(X)$ with $E_1$  similar to the one from \cite{McLean, BBLN}. This suggested us a connection with formal loops, since arc contact at infinity is pole order for loops on $X$. Work of Bouaziz \cite{BouazizCurves} on $\sL X$ for punctured curves strengthened this suggestion. There was also indirect evidence describing quantum cohomology (symplectic cohomology is the open affine counterpart of quantum cohomology for projective manifolds) in terms of formal-loop spaces: toric Fano manifolds \cite{arkhipov2006toric},  partial flag varieties \cite{lam2010quantum}, and the rather vague \cite[Open Question E]{kato2025loop}. Also, vertex algebras were constructed in \cite{KV}  from formal-loop spaces, the indirect connection between vertex algebras and Floer theory being that both are shadows of string theory, or analysis on the usual free loop space. 
With these clues in hand, we started exploring formal-loop spaces a few years ago. AI gave it an unexpected boost.

\subsection{Next steps in this project.} The following are within reach:
\begin{enumerate}[(i), wide=0pt]
\item \emph{$A_\infty$--structures on the loop complex and loop cohomology.} The multi-collision families in this paper provide a natural departure point for this.
\item \emph{De Rham loop cohomology.} In fact, Drinfeld even suggests the existence of a motivic version for the normalized dualizing complexes at least on pole stages. Endowing $LH^*$ with mixed Hodge structures is part of this. This leads to computational simplifications for the algebra structure on loop cohomology. 
\item \emph{Intrinsic loop mirrors.} The de Rham version of the fusion product can be used to address the normality Conjecture \ref{mir:conj:loop-normality} for affine loop mirrors of maximal boundary log Calabi-Yau's, which in turn would extend the coverage of the polyvector comparison of Theorem \ref{thm:intro-global-polyvectors}. On the other hand,
the comparison with the Gross-Siebert intrinsic mirror construction in this paper is only under a specialization of the latter and under the smooth affine maximal boundary  absolute case. Extending the coefficient rings and addressing the relative mirror constructions is also a  natural next step.
\item \emph{Formal-loop space counterpart of the small $QH^*(M)$.} Anticipating a positive answer to the Loop-Floer Conjecture \ref{conj:intro-multiplicative}, one can guess at least in some cases from the symplectic counterparts a deformation of the loop complex of $X$ whose cohomology gives the small quantum cohomology ring of the compactification $M$.  
\end{enumerate} 

These tasks  seem AI-easy at this point, but we have not verified the outputs. The Loop-Floer Conjecture remains difficult in general.

\subsection{Organization}
\label{subsec:intro-results}
Part~\ref{part:general}, \emph{Loop cohomology}, begins with
logarithmic cohomology in Section~\ref{sec:setup}.
Section~\ref{sec:formal-loops} constructs the formal-loop contact loci
and normalized dualizing complexes on pole stages.
Section~\ref{sec:loop-complex} defines the filtered loop complexes
$B$ and $B_\Omega$ and their cohomology $LH$ and $LH_\Omega$, proves independence of
finite-factor presentations, and gives computations from global finite factors.
Section~\ref{sec:topological-preparations} consists of general topological preparations for the fusion product, which at the core relies on Borel-Moore specialization for finite factors of poles stages. 
Section~\ref{sec:fusion-preparations} defines the collision
families we eventually use for the fusion product, and specifies what would make them fusion-ready.
Section~\ref{sec:direct-finite-factors} starts the process of making the collision families fusion-ready. It proves a relative version of Drinfeld's
finite-factor theorem in the context of our collision families. Relative Borel-Moore complexes are needed, as prescribed in the previous sections, to normalize to our infinite-dimensional context the Borel-Moore specialization for finite factors.
Section~\ref{sec:fusion-ready-families}
constructs the compatible relative volume indices, completing the normalization and the process of making the collision families
fusion-ready.
Section~\ref{sec:relative-complexes} begins the process of extracting a fusion product from the fusion-ready collision families. This section is only concerned with
factorization properties of the relative Borel-Moore complexes. Factorization is what turns specialization into a multiplication.
Section~\ref{sec:fusion-construction} finishes this process, puts everything together and constructs the fusion product. Section~\ref{sec:incidence} computes its leading
terms and proves its unit, ring, independence, and grading properties.
Section~\ref{sec:bv} proves the BV identities, gives the operator
on global finite models, and identifies its associated-graded
action on logarithmic cohomology.
Section~\ref{sec:spectral-sequence} constructs the contact spectral
sequence and proves proportional-boundary degeneration. 
Section~\ref{sec:conclusion} concludes Part~\ref{part:general}
with the K\"unneth theorem and deformation invariance for loop
cohomology.

Part~\ref{part:mirror}, \emph{Mirror symmetry via loop spaces},
begins by proving Theorems~\ref{thm:intro-lcy-zero} and
\ref{thm:intro-theta-basis} together, with an explicit construction
of the theta classes. Section~\ref{mir:sec-affine-scheme}
constructs the affine mirror, proves finite generation and the
dimension formula, and proves Theorem~\ref{mir:thm-intro-cy}
on its Gorenstein semi-log-canonical Calabi--Yau property. Section \ref{globpv:comparison} 
then establishes the global polyvector certificate for smooth
mirrors of positive log Calabi--Yau pairs with maximal boundary.
Section~\ref{mir:sec-trace-potentials} defines the loop trace,
proves that it vanishes on the image of the BV operator, and proves
nondegeneracy of its product pairing for a smooth loop mirror.
Section~\ref{mir:sec-gs-comparison} proves
Theorem~\ref{mir:thm-intro-gs}. It first explains how a collision
branch produces a global curve, recalls Gross--Siebert's polynomial
algebra theorem, and defines the curve-class specialization. It then
compares two-input and three-input loop traces with fixed-source
logarithmic curve counts. The local proof separates the construction
of collision normal forms from the computation of branch
multiplicities. Nondegeneracy reconstructs the algebra product,
and the Gross--Siebert--Johnston trace identity gives all logarithmic
point descendants.
Section~\ref{sec:surface-loop-floer} concludes Part~\ref{part:mirror}
with the proof of Theorem~\ref{thm:surface-loop-floer}.

Part~\ref{part:examples} consists of some examples of  calculations for loop cohomology and compares them when possible to symplectic counterparts. It treats affine spaces, tori, curves,
projective hypersurface complements, and the hyperbola complement.

Appendix~\ref{app:derived-colimits} records the derived-complex,
coherent-descent, homotopy-limit, and homotopy-colimit conventions, which were used in Part \ref{part:general}.


\subsection*{AI use}  AI was an essential collaborator in the development of this project, which was conceived, guided, and curated by the author.  The author verified the AI-generated content.

\subsection*{Acknowledgements} We thank J.F. de Bobadilla, R. Docampo, M. Hablicsek, S. Johnston, C.Y. Mak, L. Mihalcea, M. Pe Pereira, H. Ruddat for useful conversations. 
This work was supported by KU Leuven Grant Methusalem
METH/21/03, FWO Grant G0B3123N. LLMs available with a ChatGPT Pro subscription were used.

\addtocontents{toc}{\protect\setcounter{tocdepth}{1}}
\clearpage
\part{Loop  cohomology}\label{part:general}

\section{Logarithmic cohomology}\label{sec:setup}

This section recalls the definition from \cite{GP} of the logarithmic cohomology ring that appears
as the first page of the spectral sequence in
Theorem~\ref{thm:summary}. 

\subsection{Polarized compactifications}
\label{subsec:polarized-compactifications}
Let $M$ be a smooth projective complex variety of dimension $n>0$, let
$D=\bigcup_{i\in I}D_i$ be a reduced simple-normal-crossings divisor,
and put $X=M\setminus D$. Fix an ample divisor
\begin{equation}\label{eq:polarized-setup}
 A=\sum_{i\in I}k_iD_i,\qquad k_i\in\Z_{>0}.
\end{equation}
Write $\mathbf{k}=(k_i)_{i\in I}$.
A sufficiently large multiple of $A$ is very ample, and its defining
section realizes $X$ as an affine open. Every smooth affine complex
variety admits such a compactification.

For $\vct\in\N^I$, where $\N=\Z_{\geq0}$, write
\[
 I(\vct)=\{i:v_i>0\},\qquad |\vct|=\sum_i v_i,
 \qquad w(\vct)=\sum_i k_iv_i.
\]
For a meromorphic volume form
$\Omega$ on $M$, nowhere zero on $X$, 
writing
\[
 \operatorname{div}(\Omega)=-\sum_{i\in I}a_iD_i,
 \qquad a_i\in\Z,
\]
we define
\begin{equation}\label{eq:volume-grading}
 \mu_\Omega(\vct)=2\sum_i(1-a_i)v_i.
\end{equation}
Such $\Omega$ is available precisely when the
canonical divisor has a representative supported on $D$.

Without \(\Omega\), we read cohomological degrees modulo two. We use the Koszul sign rule in both gradings.

\subsection{Normal faces and restriction maps}
For $J\subset I$, set
\[
 D_J=\bigcap_{j\in J}D_j,\qquad
 D_J^\circ=D_J\setminus\bigcup_{i\notin J}D_i,
 \qquad D_\varnothing=M.
\]
Let $\rho:M^{\log}\to M$ be the real oriented blowup along $D$.
Locally it replaces each boundary coordinate $z_i$ by
$z_i=r_i e^{\sqrt{-1}\theta_i}$ with $r_i\geq0$; hence
$M^{\log}$ is a compact manifold with corners and interior $X$.
Define
\[
 S_J^{\log}=\rho^{-1}(D_J),\qquad
 S_J^\circ=\rho^{-1}(D_J^\circ).
\]
The map $S_J^\circ\to D_J^\circ$ is the product of the unit normal
circle bundles to $D_j$, $j\in J$, with their natural twisting.
In particular $S_\varnothing^\circ=X$ and
$S_\varnothing^{\log}=M^{\log}$.

\begin{lemma}\label{lem:face-restriction}
The inclusion $j_J:S_J^\circ\hookrightarrow S_J^{\log}$ is a
homotopy equivalence. For $J\subset L$, the face inclusion
$\iota_{JL}:S_L^{\log}\hookrightarrow S_J^{\log}$ therefore defines
a canonical homomorphism of cohomology rings
\[
 r_{JL}^*=j_L^*\iota_{JL}^*(j_J^*)^{-1}:
 H^*(S_J^\circ)\longrightarrow H^*(S_L^\circ).
\]
These maps satisfy $r_{LK}^*r_{JL}^*=r_{JK}^*$ for
$J\subset L\subset K$.
\end{lemma}
\begin{proof}
The interior $S_J^\circ$ of the manifold with corners
$S_J^{\log}$ is homotopy equivalent to it: compatible boundary
collars push all boundary-defining coordinates to positive values,
and interpolation gives the required homotopies. The displayed
formula uses only the inverse of the induced cohomology isomorphism,
so is independent of the collars. Composition follows from composition
of the face inclusions.
\end{proof}

\subsection{The logarithmic ring}
\begin{definition}[Logarithmic cohomology]\label{def:log-cohomology}
Following \cite{GP}, define the logarithmic ring with the following
gradings.

\emph{Parity grading.}
For $\epsilon\in\Z/2$, put
\begin{equation}\label{eq:log-cohomology-definition}
 H^{\epsilon}_{\log}(M,D)
 :=\bigoplus_{\vct\in\N^I}
       \bigoplus_{r\equiv\epsilon\ (2)}
       H^r(S_{I(\vct)}^\circ;\F)t^{\vct}.
\end{equation}
The symbols $t^{\vct}$ record contact multiplicities. Each element
has finite contact support, and empty faces contribute zero.
For $\alpha\in H^r(S_J^\circ;\F)$ and
$\beta\in H^s(S_K^\circ;\F)$, where $J=I(\vct)$ and
$K=I(\uct)$, set $L=J\cup K$ and define
\begin{equation}\label{eq:end-product}
 (\alpha t^{\vct})(\beta t^{\uct})
 =\bigl(r_{JL}^*\alpha\smile r_{KL}^*\beta\bigr)t^{\vct+\uct}.
\end{equation}
This product is zero when $D_L^\circ=\varnothing$.
The parity of $\alpha t^{\vct}$ is $\bar r$, its weight is
$w(\vct)$, and its spectral-sequence bidegree is
$(p,q)=(-w(\vct),\overline{r+w(\vct)})\in\Z\times\Z/2$.
No volume form is required.

\emph{Integer grading.}
Given a nowhere-zero algebraic volume form $\Omega$ on $X$,
put, for $d\in\Z$,
\[
 H^d_{\log}(M,D)
 :=\bigoplus_{\vct\in\N^I}
       H^{d-\mu_\Omega(\vct)}(S_{I(\vct)}^\circ;\F)t^{\vct},
\]
where $\mu_\Omega$ is defined in \eqref{eq:volume-grading}.
Use the same product \eqref{eq:end-product} and weight.
Thus $\alpha\in H^r(S_{I(\vct)}^\circ;\F)$ gives
$|\alpha t^{\vct}|=r+\mu_\Omega(\vct)$ and bidegree
\begin{equation}\label{eq:bidegree}
 (p,q)=\bigl(-w(\vct),r+\mu_\Omega(\vct)+w(\vct)\bigr)
       \in\Z^2.
\end{equation}
The function $\mu_\Omega$ is additive and even, so the product
preserves integer degrees and reduction modulo two recovers
the parity-graded ring.

In either case write $H^*_{\log}(M,D)$ for the resulting graded
ring and $H^{p,q}_{\log}(M,D)$ for its spectral-sequence
bigrading. We suppress the polarization and, in integer grading,
the chosen volume form from the notation.
\end{definition}

\begin{proposition}[\cite{GP}]\label{prop:end-ring} 
The product \eqref{eq:end-product} makes $H^*_{\log}(M,D)$ an
associative, unital, graded-commutative ring. It adds bidegrees, and
its unit is $1\in H^0(X)$ in the zero-contact summand.
\end{proposition}

\addtocontents{toc}{\protect\setcounter{tocdepth}{2}}
\section{Formal loops and normalized dualizing complexes}\label{sec:formal-loops}

In this section we define the analytic topology on formal
loop spaces, their pole stages, and their contact strata.
Drinfeld's finite-factor charts then allow us to construct
normalized dualizing complexes on the pole stages. We compute  the volume indices behind this normalization process in terms of contact orders at infinity.

\subsection{Formal loops and their analytic topology}\label{subsec:loop-topology}
Write $\mathcal D=\Spec\C\llbracket t \rrbracket$ and
$\mathcal D^*=\Spec\C(\!( t )\!)$. For an affine variety $Y$ we use
\[
 \mathscr L Y(R)=Y(R(\!( t )\!)),\qquad
 \mathscr L^+Y(R)=Y(R\llbracket t \rrbracket)
\]
for formal loops and  arcs. The  notation $\sL X$ and $\sL ^+X$ also denotes their
spaces of complex points when a topology is discussed.

\begin{definition}\label{par:coefficient-schemes}\label{def:laurent-coefficients}
Fix a closed immersion
\[
 Y=\Spec\C[x_1,\ldots,x_m]/(f_1,\ldots,f_l)\subset\A^m.
\]
For $b\ge0$, introduce indeterminates $x_{i,j}$, $1\le i\le m$,
$j\ge-b$, and put
\begin{equation}\label{eq:coefficient-ring}
 \begin{split}
 C_b(Y)&=\C[x_{i,j}:j\ge-b]/I_b,\\
 I_b&=\bigl([t^a]f_\nu(\textstyle\sum_{j\ge-b}x_{1,j}t^j,
                 \ldots,\sum_{j\ge-b}x_{m,j}t^j):\nu,a\bigr),
 \end{split}
\end{equation}
where $[t^a]g(t)$ means the coefficient of $t^a$ in $g(t)$.
Each displayed coefficient is a polynomial in finitely many variables.
The affine scheme $\Spec C_b(Y)$ represents exactly the functor
\[
 R\longmapsto\{\gamma\in Y(R(\!( t )\!)):x_i(\gamma)\in t^{-b}R\llbracket t \rrbracket
                         \text{ for all }i\}.
\]
We call this the \emph{Laurent coefficient scheme of bound $b$ for
this immersion}, or its \emph{Laurent stage of bound $b$}.
This bound refers to the specified coordinate functions $x_i$.
The scheme is usually not of finite type. For $b'\ge b$,
its closed immersion into $\Spec C_{b'}(Y)$ sets
$x_{i,j}=0$ for $-b'\le j<-b$. The union over $b$ represents
$\mathscr L Y$. A polynomial change of affine generators gives the
same ind-scheme by cofinality of these bounds.

 Let $C_{b,N}(Y)$ be the
subalgebra of $C_b(Y)$ generated by the classes of $x_{i,j}$ with
$-b\le j\le N$. Then
\begin{equation}\label{eq:literal-coefficient-limit}
 \Spec C_b(Y)=\varprojlim_N\Spec C_{b,N}(Y).
\end{equation}
The defining ideal of $C_{b,N}(Y)$ is the contraction of $I_b$ to
this finite polynomial ring, so it includes all relations among the
retained coefficients. The transitions in
\eqref{eq:literal-coefficient-limit} need not be smooth or affine
bundles, so this presentation alone supplies no smooth-pullback
isomorphisms. The finite factors used for normalized dualizing
complexes are introduced in
Subsection~\ref{subsec:drinfeld-finite-factors}.
\end{definition}

Give $\Spec C_b(Y)(\C)$ the inverse-limit topology from
\eqref{eq:literal-coefficient-limit}, using the usual analytic
topology on the complex points of each finite-type scheme. Give
$\mathscr L Y(\C)$ the final topology for the inclusions of these
spaces as $b$ increases. This is the \emph{ind--pro analytic topology}.
\begin{convention}[Topology on formal loops]
\label{conv:loop-analytic-topology}
Throughout this article, the \emph{analytic topology on
$\mathscr LY(\C)$} means the ind--pro analytic topology just
defined. Unless another topology is explicitly specified, all
topological assertions about formal loops, including singular
cohomology, use this topology. 
\end{convention}

\begin{example}
For $X=\A^1$,
\[
 \mathscr LX(\C)=\C(\!( t )\!)=\bigcup_{b\ge0}E_b,
 \qquad
 E_b:=t^{-b}\C\llbracket t \rrbracket
     \cong\prod_{j\ge-b}\C.
\]
Each $E_b$ has the product topology, so a basic neighborhood
restricts only finitely many coefficients. The union has the
final topology:
\[
 U\subset\C(\!( t )\!)\text{ is open}
 \iff U\cap E_b\text{ is open for every }b.
\]
A sequence converges precisely when its pole orders are uniformly
bounded and every coefficient converges in $\C$. Indeed, boundedness
places the sequence and its coefficientwise limit in a single product
stage. Conversely, coefficient functions are continuous by the final
topology. If $f_m\to f$ had unbounded pole orders, one could choose a
subsequence $f_{m_q}$ with strictly increasing pole orders $b_q$, all
larger than the pole order of $f$. Put
$a_q=[t^{-b_q}]f_{m_q}\ne0$. Then
\[
 V=\{g\in\C(\!( t )\!):|[t^{-b_q}]g|<|a_q|/2\text{ for every }q\}
\]
is a neighborhood of $f$: its intersection with each $E_b$ imposes
only finitely many inequalities. It excludes every $f_{m_q}$, a
contradiction. For example,
\[
 \frac1n t^{-1}\longrightarrow0,\qquad
 n!t^n\longrightarrow0,\qquad
 \sum_{j=0}^{n}j!t^j\longrightarrow\sum_{j\ge0}j!t^j.
\]
In contrast, $n^{-1}t^{-n}$ does not converge to zero, although
each fixed coefficient does. Indeed,
\[
 U=\left\{f\in\C(\!( t )\!):
       |[t^{-j}]f|<\frac1{2j}\text{ for every }j\ge1\right\}
\]
is an open neighborhood of zero excluding every $n^{-1}t^{-n}$:
its intersection with $E_b$ imposes only $b$ inequalities.
This topology also differs from the $t$-adic topology, since
the constant series $1/n$ converges to zero here.
\end{example}

The ind-representability assertion below is due to \cite{KV};
we prove the assertions about the complex-point topology directly.

\begin{proposition}
\label{lp:prop:topology}
Let $Y$ be an affine complex variety. The analytic topology on $\sL Y$ is independent of the  closed immersion of $Y$ into affine
space, and every morphism $Y\to Z$ of affine varieties induces a
continuous map $\mathscr L Y\to\mathscr L Z$.
The analytic topology on $\sL Y$ is compactly generated and Hausdorff.
For each chosen immersion and every $b\ge0$, the inclusion
\[
 (\Spec C_b(Y))(\C)\longrightarrow\mathscr L Y
\]
is a closed topological embedding, with the inverse-limit topology
on its source. For every $f\in\mathcal O(Y)$ and $j\in\Z$, the
function $\gamma\mapsto[t^j]f(\gamma(t))$ is continuous.
The functor $R\mapsto Y(R(\!( t )\!))$ is represented by the ind-scheme
$\varinjlim_b\Spec C_b(Y)$, whose terms are affine schemes,
usually not of finite type.
\end{proposition}
\begin{proof}
The ind-representability assertion and its coefficient-equation
construction are \cite[Proposition~1.5.1(a)]{KV}. We verify the
additional assertions about complex-point topology. Polynomial
substitution takes each Laurent bound continuously into a sufficiently
large bound, since every output coefficient is polynomial in finitely
many input coefficients. This proves independence of generators,
functoriality, and continuity of the displayed coefficient functions.

Each $(\Spec C_b(Y))(\C)$ is a closed subspace of a countable
product of copies of $\C$, hence metrizable and compactly generated.
Consequently a subset of the final space that is closed under all
tests by maps from compact Hausdorff spaces has closed inverse image
in every Laurent stage, and is closed. This proves compact generation.
The coefficient functions separate points, proving Hausdorffness.
The image of $(\Spec C_b(Y))(\C)$ is cut out by vanishing of
all $[t^j]x_i$ with $j<-b$; its retained coefficient functions
recover its inverse-limit topology. Thus the inclusion is a closed
topological embedding.
\end{proof}

So a convergent sequence of loops in our topology has a unique limit, due to Hausdorffness. Below, paths, homotopy and cohomology always use the analytic topology of $\sL Y$. 

When discussing algebraic irreducibility in $\sL Y$, we instead give each
$\Spec C_b(Y)$ its Zariski topology and use the final
\emph{ind-Zariski topology}. 

Passing to a reduced Laurent stage, that is, replacing $C_b(Y)$ by
its quotient by all nilpotent elements,
does not change the underlying Zariski space or the complex points
with their analytic topology.

\begin{lemma}\label{lp:lem:topological-products}\label{lp:lem:arc-contraction}
For affine varieties $Y,Z$,
\[
 \mathscr L(Y\times Z)\cong\mathscr L Y\times\mathscr L Z,
 \qquad \mathscr L^+ Y\sim Y(\C).
\]
The first map is a natural homeomorphism for the analytic topology
and the compactly generated product; it also comes from an
isomorphism of the corresponding functors. The second equivalence is
a homotopy equivalence for the inverse-limit analytic topology on
$\mathscr L^+Y(\C)$. More precisely, center evaluation and the
inclusion of constant arcs are inverse up to the strong deformation
retraction $\gamma(t)\mapsto\gamma(st)$, $s\in[0,1]$, which is
center evaluation at $s=0$ and the identity at $s=1$.
\end{lemma}
\begin{proof}
Use the union of the two generating sets. Products preserve
colimits in compactly generated spaces, and the diagonal in
$\N^2$ is cofinal. For arcs, substitution multiplies the
$j$th coefficient by $s^j$, continuously; at $s=0$ it gives
the constant arc at the center.
\end{proof}

The first assertion is not necessarily true in the usual product topology: for $Y=Z=\bA^1$, $U=\{(f,g): |[t^0](fg)|<1\}\subset \C(\!( t )\!)^{\times 2}$ is open in $\sL \bA^2$, but not open in the usual product topology on $\sL \bA^1\times\sL \bA^1$.

\subsection{Contact strata and logarithmic cohomology}\label{sec:ends}

Let
\[
 \cL_\ell(M)(R)=\operatorname{Hom}_{\C}
       \bigl(\Spec R[t]/(t^{\ell+1}),M\bigr),\qquad
 \cL^+(M)=\varprojlim_\ell\cL_\ell(M).
\]
Since $M$ is quasi-compact and quasi-separated, Bhatt's
algebraization theorem \cite[Corollary~1.2]{Bhatt} gives natural
bijections, for every $\C$-algebra $R$,
\[
 M(R\llbracket t \rrbracket)\xrightarrow{\ \sim\ }
 \varprojlim_\ell M(R[t]/(t^{\ell+1}))
 =\cL^+(M)(R).
\]
Thus $\cL^+(M)$ represents the formal-arc functor also for this
non-affine $M$. On complex points we give
\[
 \cL^+(M)(\C)=\varprojlim_\ell\cL_\ell(M)(\C)
\]
the inverse-limit topology of the usual analytic topologies on
the finite-type jet schemes. In particular, a map to this arc
space is continuous if and only if all its jet truncations are
continuous.

Write $j_\ell:\cL^+(M)\to\cL_\ell(M)$ for jet truncation.
For $\vct\in\N^I$, put $J=I(\vct)$ and
$h(\vct)=\max(\{v_i:i\in I\}\cup\{0\})$.
The \emph{contact locus} $\mathscr P_{\vct}$ is the locally
closed subscheme of $\cL^+(M)$ whose arcs have center in $D_J^\circ$
and contact order $v_i$ with every $D_i$. Precisely, in an SNC
coordinate neighborhood with $D_i=(z_i=0)$ for $i\in J$, write
$z_i(\gamma(t))=\sum_{a\ge0}z_{i,a}t^a$. The defining conditions are
\begin{equation}\label{eq:contact-sector-definition}
 z_{i,a}=0\quad(0\le a<v_i),\qquad
 z_{i,v_i}\ne0\quad(i\in J),\qquad
 \gamma(0)\in D_J^\circ.
\end{equation}
Here $z_{i,v_i}\ne0$ means the principal open subscheme on which this
coefficient is invertible. Thus, for an arc over a $\C$-algebra $R$,
the condition is $z_i(\gamma(t))=t^{v_i}u_i(t)$ with
$u_i(t)\in R\llbracket t \rrbracket^\times$. These local definitions glue because
changing a local equation of $D_i$ multiplies it by a unit. In
particular, $\mathscr P_{\zero}=\cL^+(X)$. If $D_J^\circ$ is empty,
then $\mathscr P_{\vct}$ is empty.

For $\ell\ge h(\vct)$, we use the notation $j_\ell(\mathscr P_{\vct})$
for the locally closed subscheme of $\cL_\ell(M)$ defined by the same
equations and open conditions \eqref{eq:contact-sector-definition}.
It is the jet truncation of the contact locus: smoothness of $M$ gives
surjective truncations, and
\begin{equation}\label{eq:contact-sector-truncation}
 \mathscr P_{\vct}
 =j_\ell^{-1}\bigl(j_\ell(\mathscr P_{\vct})\bigr)
 =\varprojlim_{m\ge h(\vct)}j_m(\mathscr P_{\vct}).
\end{equation}
The first equality is an equality of subschemes of $\cL^+(M)$.

Properness of $M$ identifies $\mathscr LX(\C)$ with the arcs in
$\cL^+(M)(\C)$ whose generic point lies in $X$: a loop
$\gamma:\Spec\C(\!( t )\!)\to X$ extends uniquely to an arc
$\bar\gamma:\Spec\C\llbracket t \rrbracket\to M$. Under this bijection we regard
$\mathscr P_{\vct}(\C)$ as a subspace of $\mathscr LX(\C)$, with
its induced analytic topology, and abbreviate it to
$\mathscr P_{\vct}$ in topological formulas. This identification
concerns complex points only.

\begin{definition}[Pole stages]
\label{def:pole-stages}
For $P\in\Z_{\ge0}$, define the subspace
\[
 F_P\mathscr LX
 :=\bigcup_{\substack{\vct\in\N^I\\w(\vct)\le P}}
       \mathscr P_{\vct}(\C)
 \subset\mathscr LX(\C),\qquad
 F_{-1}\mathscr LX=\varnothing,
\]
with the induced analytic topology. Equivalently,
\[
 F_P\mathscr LX
 =\cL^+(X)(\C)
 \;\cup\!
 \bigcup_{\substack{\vct\ne\zero,\ w(\vct)\le P\\
                    D_{I(\vct)}^\circ\ne\varnothing}}
       \mathscr P_{\vct}(\C).
\]
Indeed, $I(\vct)=\varnothing$ precisely when $\vct=\zero$,
and $D_\varnothing^\circ=X$, so
$\mathscr P_{\zero}=\cL^+(X)$. A contact vector with
$D_{I(\vct)}^\circ=\varnothing$ contributes the empty set.
In particular, $F_0\mathscr LX=\cL^+(X)(\C)$.
\end{definition}

Next proposition is similar to
 \cite[Proposition~3.2]{EmbeddedNash}. 

\begin{proposition}
\label{prop:stable-strata}
Let $\vct\in\N^I$, $J=I(\vct)$, and $\ell\ge h(\vct)$.
Assume $D_J^\circ\ne\varnothing$. The $\ell$-jet contact locus $j_\ell(\mathscr P_{\vct})$ is an iterated
Zariski-locally-trivial affine fibration over
\[
 \underset{D_J^\circ}{\prod_{j\in J}}
 \bigl(N_{D_j/M}|_{D_J^\circ}\setminus0\bigr),
\]
 the fiber product over $D_J^\circ$ of the punctured normal
line bundles indexed by $j\in J$. For $J=\varnothing$, this
base is the empty fiber product over $D_\varnothing^\circ=X$,
namely $X$; the jet contact locus is then $\cL_\ell(X)$.
If $D_J^\circ=\varnothing$, the contact locus and its jet truncations
are empty and contribute zero cohomology. In the nonempty case,
\[
 H^*\bigl(j_\ell(\mathscr P_{\vct});\F\bigr)
 \cong H^*(S_J^\circ;\F),
\]
canonically and compatibly with further truncation. The map
$j_{\ell+1}(\mathscr P_{\vct})\to j_\ell(\mathscr P_{\vct})$
is an affine $n$-bundle, where $n=\dim_\C X$.
\end{proposition}

\begin{proof}
Write $T_J$ for the displayed product of punctured normal bundles.
The center and the leading normal coefficients define a morphism
$j_\ell(\mathscr P_{\vct})\to T_J$: replacing a local equation
$z_j$ of $D_j$ by $u_jz_j$ multiplies its leading coefficient
by $u_j$ evaluated at the center, which is precisely the normal-bundle
transition rule. On an affine \'etale SNC coordinate neighborhood,
write, modulo $t^{\ell+1}$,
\[
 z_j(t)=t^{v_j}\Bigl(c_j+\sum_{a=1}^{\ell-v_j}c_{j,a}t^a\Bigr)
 \quad(j\in J),\qquad
 z_k(t)=z_k(0)+\sum_{a=1}^{\ell}b_{k,a}t^a
 \quad(k\notin J),
\]
with $c_j$ invertible. Formal \'etaleness identifies the jet over
its center with these coordinate expansions.
For $a=1,\ldots,\ell$, retain all tangential coefficients $b_{k,a}$
and those normal coefficients $c_{j,a}$ for which $a\le\ell-v_j$.
Here $a$ is the relative order in the normal expansion: $c_{j,a}$
is the coefficient of $t^{v_j+a}$ in $z_j(t)$. Under an SNC
coordinate change, the new tangential coefficients depend linearly
on the $b_{k,a}$, with linear part the invertible tangential Jacobian
along $D_J^\circ$. Any normal coefficient entering their remaining
terms has relative order at most $a-v_j<a$, so is already retained.
Expanding $z'_j=u_j(z)z_j$, the new $c_{j,a}$ has coefficient
$u_j(\gamma(0))$ on the old $c_{j,a}$; its only other dependence on
newly retained coefficients is linear in the $b_{k,a}$. All other
terms use previously retained coefficients. Thus the transition at
each step is affine with block triangular, invertible linear part.
The coordinate-change identities give the cocycle identities for
these affine transitions, so each step is a torsor under a
finite-rank vector bundle and hence a Zariski-locally-trivial affine
fibration. Separately, increasing the absolute jet order from $\ell$
to $\ell+1$, with $\ell\ge h(\vct)$, adds exactly $n$ unrestricted
coefficients. Their transitions have linear part the full Jacobian
of the coordinate change at the center, proving the stated affine
$n$-bundle. Analytically these affine fibrations have contractible
fibers. The quotient of each punctured normal line by positive real
scaling is its oriented normal circle, so the leading-normal map
induces the asserted canonical cohomology isomorphism. It commutes
with truncation.
\end{proof}

The following proposition is an analogue for loops, in the SNC
setting, of \cite[Lemmas~2.3 and~2.5]{BBLN}:

\begin{proposition}\label{prop:loop-sectors}
For every integer $P\ge0$, $F_P\mathscr L X$ is a closed subspace of
$\mathscr L X(\C)$ in the analytic topology. More precisely,
choose $r>0$ such that $rA$ is very ample and a basis
$s_0,s_1,\ldots,s_m$ of $H^0(M,\mathcal O_M(rA))$ with
$\operatorname{div}(s_0)=rA$. For the resulting closed immersion
$X\hookrightarrow\A^m$ with coordinates $s_j/s_0$, one has
\[
 F_P\mathscr L X=(\Spec C_{rP}(X))(\C)
\]
as subspaces of $\mathscr L X(\C)$. For the following algebraic
assertion, give $F_P\mathscr LX$ the Laurent coefficient-scheme structure
$\Spec C_{rP}(X)$ and put $F_{-1}\mathscr LX=\varnothing$.
Restriction of arcs to the punctured disk induces an isomorphism
of reduced schemes
\[
 \bigl(F_P\mathscr L X\setminus F_{P-1}\mathscr L X\bigr)_{\mathrm{red}}
 \cong\coprod_{\substack{w(\vct)=P\\D_{I(\vct)}^\circ\ne\varnothing}}
 \mathscr P_{\vct}.
\]
The contact schemes on the right are reduced. On complex points
this gives a finite topological coproduct in the analytic topology.
On every nonempty contact
locus, the subspace topologies induced by $\mathscr LX(\C)$ and by
$\cL^+(M)(\C)$ coincide with the inverse-limit analytic topology
of its contact-jet truncations. Sending an arc to its leading
normal vectors and then to their classes modulo positive real
scaling defines a canonical homotopy equivalence
\[
 \mathscr P_{\vct}\longrightarrow S_{I(\vct)}^\circ.
\]
For $\vct=\zero$, this map is center evaluation
$\cL^+(X)(\C)\to X(\C)$.
The continuous bijection
$\coprod_{\vct}\mathscr P_{\vct}\to\mathscr L X$
need not be a homeomorphism.
\end{proposition}

\begin{proof}
\emph{The Laurent bound.}
Write $g_i=s_i/s_0$ for $1\le i\le m$, and let $\bar\gamma$
be the extension to $M$ of a loop $\gamma$. If its contact vector
is $\vct$, then in a local trivialization of $\mathcal O_M(rA)$
at its center,
\[
 \ord_t\bar\gamma^*s_0=r\sum_i k_iv_i=r\,w(\vct).
\]
The sections $s_0,\ldots,s_m$ have no common zero, so at least
one pulled-back section has order zero. Consequently
\[
 \max\bigl(\{0\}\cup
       \{-\ord_t g_i(\gamma):1\le i\le m\}\bigr)
       =r\,w(\vct).
\]
Thus $w(\vct)\le P$ is equivalent to the Laurent bound $rP$
for every $g_i$. Proposition~\ref{lp:prop:topology} identifies
this Laurent stage with a closed topologically embedded subspace
of $\mathscr LX$, proving the first assertion.

\emph{Algebraic identification and analytic topology.}
Fix $P>0$, put $b=rP$, and let $U_P$ be the open complement
of $\Spec C_{r(P-1)}(X)$ in $\Spec C_b(X)$.
The Laurent-bound calculation above applies over every residue
field: the maximal pole order is always a multiple of $r$.
Consequently $U_P$ is also the open complement of
$\Spec C_{b-1}(X)$, and the principal opens where some
$[t^{-b}]g_j$ is invertible cover it.
On such an open, write the universal series as
$g_i(t)=t^{-b}u_i(t)$, with $u_i\in R\llbracket t \rrbracket$ and
$u_j(0)\in R^\times$, where $R$ is its coordinate ring.
They define an arc with projective coordinates
\[
 [t^b:u_1(t):\cdots:u_m(t)]
\]
and affine coordinates in the chart $s_j\ne0$
\[
 y_0(t)=\frac{t^b}{u_j(t)},\qquad
 y_i(t)=\frac{u_i(t)}{u_j(t)}\quad(i\ne j),
 \qquad y_j(t)=1.
\]
These are precisely $s_0/s_j$ and $s_i/s_j$ evaluated on the
arc. The homogeneous equations of $M$ hold after inverting $t$,
hence in $R\llbracket t \rrbracket$ because $R\llbracket t \rrbracket\to R(\!( t )\!)$ is injective.
Each coefficient of the displayed series is a regular function
of finitely many Laurent coefficients and $u_j(0)^{-1}$, so
the formulas define morphisms of the coefficient and arc schemes.
The arc has total contact exactly $b$ with $rA$: the pulled-back
local equation satisfies
\[
 y_0(t)=t^b a(t),\qquad a(t)\in R\llbracket t \rrbracket^\times.
\]
This equality is imposed over $R$, not merely after passage to
its residue fields.
Conversely, on the arc scheme defined by this total-contact
condition, write $y_0=t^b a(t)$ with $a(0)$ invertible and set
\[
 g_i(\gamma(t))=\frac{y_i(t)}{y_0(t)}.
\]
The inverse-series coefficients are again regular functions of
finitely many arc coefficients and $a(0)^{-1}$.
These inverse constructions agree on overlaps and identify $U_P$
algebraically with the scheme of total contact exactly $b$ just described.

To decompose its reduction, fix a residue-field arc with contact
$\vct$ and choose an SNC coordinate neighborhood of its center
meeting no $D_i$ with $i\notin I(\vct)$.
Consider the open subset where the center stays in this neighborhood
and every coefficient $z_{i,v_i}$, $i\in I(\vct)$, is invertible.
Every residue-field arc in this open has contacts $v_i'\le v_i$,
and the total-contact condition gives
\[
 \sum_i k_i v_i'=P=\sum_i k_i v_i,
\]
so positivity of the $k_i$ forces $\vct'=\vct$.
Thus every lower coefficient $z_{i,a}$, $a<v_i$, vanishes at
every prime of this open and therefore vanishes on its reduction.
Its reduced scheme consequently factors through the contact scheme
$\mathscr P_{\vct}$ defined in
\eqref{eq:contact-sector-definition}. Conversely, the part of
$\mathscr P_{\vct}$ centered in this neighborhood maps to this
open total-contact locus. In \'etale SNC
coordinates its jet schemes are smooth and its omitted coefficients
are free affine coordinates, so $\mathscr P_{\vct}$ is reduced.
The two factorizations are inverse. The finitely many contact
pieces are open and closed and give the asserted reduced scheme
isomorphism. For $P=0$, one has already
$\Spec C_0(X)=\cL^+(X)$ as schemes, and this scheme is reduced
because $X$ is smooth.

All coefficient formulas above and their inverses are continuous
on complex points in the inverse-limit analytic topologies.
Reduction leaves these topologies unchanged. This proves the
topological coproduct assertion and the equality of the subspace
topologies on each contact locus; use also
\eqref{eq:contact-sector-truncation} for its contact-jet topology.
Reduction cannot in general be omitted from the algebraic assertion:
over $\C[\epsilon]/(\epsilon^2)$, the series
$z_1=t+\epsilon$ and $z_2=t-\epsilon$ have product $t^2$
but neither has individual contact order exactly one in the
above sense: neither series belongs to $tR\llbracket t \rrbracket$ for
$R=\C[\epsilon]/(\epsilon^2)$, although both reduce to $t$
over the residue field.

\emph{The homotopy equivalence.}
Put $J=I(\vct)$ and, for this proof, write
\[
 T_J:=\underset{D_J^\circ}{\prod_{j\in J}}
       \bigl(N_{D_j/M}|_{D_J^\circ}\setminus0\bigr).
\]
For $J=\varnothing$ set $T_J=X$. In SNC coordinates an arc in
$\mathscr P_{\vct}$ has expansions
\[
 z_j(t)=t^{v_j}(c_j+c_{j,1}t+\cdots),\qquad c_j\ne0
 \quad(j\in J).
\]
Its center and the leading coefficients $c_j$ define a map
$p_{\vct}:\mathscr P_{\vct}\to T_J$. If a local equation
$z_j$ is replaced by $a_jz_j$ for a unit $a_j$, its leading
coefficient is multiplied by $a_j(\gamma(0))$; this is precisely
the change of coordinates for the normal vector. Thus
$p_{\vct}$ is globally defined.

Over an analytic SNC coordinate neighborhood, fixing the center
and these leading normal vectors leaves all higher normal
coefficients and all positive-degree tangential coefficients free.
It follows that $p_{\vct}$ is locally a product projection with
fiber a countable product of copies of $\C$ (possibly a point).
Indeed, analytic coordinate changes and their inverses give
compatible maps on all finite jet spaces; each output coefficient
uses only finitely many input coefficients and derivatives evaluated
at the center, which vary continuously with the center. Both maps
are therefore continuous for the inverse-limit topology, giving
the stated local trivializations without requiring convergence of
any formal series. The fiber is contractible by coefficientwise
scalar multiplication. Since $T_J$ is a finite-dimensional
paracompact manifold, the trivializations admit a subordinate
partition of unity. Each local product projection is a fibrewise
homotopy equivalence to the identity bundle over its base, using
this contraction. The local-to-global theorem for fibrewise
homotopy equivalences over a numerable cover
\cite[Theorem~13.3.3]{Dieck} therefore makes $p_{\vct}$ a fibrewise
homotopy equivalence and, in particular, a homotopy equivalence.

A point of the real oriented normal circle is a nonzero normal
vector modulo positive real scaling. Hence the definition of the
real oriented blowup gives a canonical identification
\[
 T_J/(\R_{>0})^J=S_J^\circ.
\]
The quotient map is a locally trivial bundle with contractible
fiber $(\R_{>0})^J$, so the same theorem makes it a homotopy
equivalence. Its composite with $p_{\vct}$ is the canonical map
in the statement. For $J=\varnothing$, it is center evaluation,
with the explicit deformation retraction
$\gamma(t)\mapsto\gamma(st)$ of
Lemma~\ref{lp:lem:arc-contraction}.

Finally, take $X=\A^1$, $M=\Pj^1$, $D=\{\infty\}$, and $A=D$.
The family $\gamma_s(t)=s/t$ is continuous in the Laurent stage
of bound one, including at $s=0$. Its contact is one for $s\ne0$
and zero for $s=0$. Its lift to the disjoint union of contact
loci is therefore not continuous. This proves the final assertion.
\end{proof}

We will also call the closed subspaces $F_P\mathscr L X$ the
\emph{pole stages for $A$}, or simply the \emph{pole stages} when
$A$ is fixed. Their intrinsic bound is $\ord_t\bar\gamma^*A\le P$;
the corresponding Laurent bound for the coordinates above is $rP$.
\begin{corollary}
\label{cor:contact-log-bigrading}
Fix a volume form $\Omega$. For every bidegree $(p,q)$ there is a
canonical isomorphism
\begin{equation}\label{eq:contact-log-bigrading}
 H^{p,q}_{\log}(M,D)
 \cong
 \bigoplus_{\substack{\vct\in\N^I,\ w(\vct)=-p\\
                       D_{I(\vct)}^\circ\ne\varnothing}}
 H^{p+q-\mu_\Omega(\vct)}(\mathscr P_{\vct};\F).
\end{equation}
In the parity theory, omit $\mu_\Omega$ and read total degrees
modulo two. Under these identifications the logarithmic product is the common-face product \eqref{eq:end-product}. 
\end{corollary}
\begin{proof}
Proposition~\ref{prop:loop-sectors} gives canonical cohomology
isomorphisms $H^*(\mathscr P_{\vct};\F)\cong
H^*(S_{I(\vct)}^\circ;\F)$. By
Definition~\ref{def:log-cohomology}, a class in
$H^k(\mathscr P_{\vct};\F)$ therefore corresponds to a class
$\alpha t^{\vct}$ of weight $w(\vct)$ and total degree
$k+\mu_\Omega(\vct)$. This gives the displayed formula and its
parity version. Transporting the restriction maps between the
normal faces gives the stated common-face product.
\end{proof}


\subsection{Drinfeld's finite-factor theorem}
\label{subsec:drinfeld-finite-factors}

Write $\A^J=\Spec\C[u_j:j\in J]$, with $J$ finite or countable.
A \emph{nice chart} of a $\C$-scheme $S$ is a finitely presented
\'etale morphism of $\C$-schemes $U\to S$, together with an
isomorphism of $\C$-schemes
\[
 U\cong Y\times_{\C}\A^J,
\]
where $Y$ is an affine scheme of finite type over $\C$, possibly
singular.  Here \'etale means smooth of relative dimension zero.
Finite presentation is relative: for an affine chart map
$\Spec B\to\Spec A$ it means
$B=A[x_1,\ldots,x_m]/(f_1,\ldots,f_r)$ with $m,r$ finite;
it does not require $A$ or $B$ to be of finite type over $\C$.

An \emph{atlas} is a jointly surjective family of these
chart morphisms. We call $S$ \emph{locally nice} if it has such an
atlas, and \emph{nice} if it itself has the displayed product
decomposition. 
The projection $\pi:U\to Y$ is called a \emph{finite factor}.
Retaining a finite subset $J_0\subset J$ replaces this
factor by $Y\times\A^{J_0}$, and the projection back to $Y$ is an
oriented affine-bundle refinement of relative dimension $|J_0|$.

A closed quasi-compact subscheme $S$ of a loop ind-scheme is
\emph{reasonable} if its ideal in every larger closed quasi-compact
subscheme is finitely generated. 

We use the following theorem of \cite[Section~6]{Drinfeld} to
construct $B$. We will prove in Theorem~\ref{thm:relative-finite-factors}
a stronger version, as explained in
Remark~\ref{rem:direct-recovers-drinfeld}.

\begin{theorem}[Drinfeld]
\label{thm:drinfeld-finite-factors}
Let $X$ be a smooth affine variety over $\C$, and let
$\mathscr LX$ be the loop ind-scheme over $\C$. Every reasonable closed quasi-compact subscheme $S\subset\mathscr LX$
admits a jointly surjective family of finitely presented \'etale
maps $U\to S$ with $U\cong Y\times\A^J$, where $Y$ is affine of
finite type over $\C$. Thus $S$ is locally nice.
If $S$ is also quasi-separated, it admits a finite filtration
\[
 S=S_0\supset S_1\supset\cdots\supset S_m=\varnothing
\]
by finitely presented closed subschemes such that each successive
difference $S_i\setminus S_{i+1}$ is isomorphic to
$Y_i\times\A^{J_i}$, with $Y_i$ affine of finite type over $\C$
and $J_i$ finite or countable.
\end{theorem}

Here local niceness follows from \cite[Theorem~6.3 and
Section~6.3.7]{Drinfeld}, and the finite filtration is given in
\cite[Section~6.11]{Drinfeld}. His Nisnevich charts are in
particular \'etale charts; taking affine opens of the finite-type
factors, and affine locally closed pieces of the finite-type
strata, gives the affine formulation above. Countability follows
from the coefficient presentation in our setting. Indeed, $S$ is
contained in a bounded Laurent stage: on a finite affine cover,
take a common bound for the finitely many coordinate functions of
$X$. Its coordinate algebra is consequently countably generated
over $\C$, as are those of its finitely presented affine charts
and affine locally closed strata. Such an algebra has at most
countable dimension as a $\C$-vector space. If
$U=\Spec A[u_j:j\in J]$ is nonempty, the elements $u_j$ are
$\C$-linearly independent, so $J$ is at most countable.

Chart maps cannot in general be replaced by Zariski open embeddings \cite[Section~6.13]{Drinfeld}.
The pole stages $F_P\sL X$ are affine, hence quasi-compact and
quasi-separated. Their reasonableness is checked at the beginning
of the proof of Proposition~\ref{prop:normalized-pole-sheaves},
so both conclusions of the theorem will apply to them.

\begin{remark}
\label{rem:finite-nice-atlas}
Every quasi-compact locally nice scheme $S$ admits a nice
\'etale atlas with finitely many charts. Indeed, the images
of the chart maps are open and cover $S$, so quasi-compactness
gives a finite subcover. This applies to every $S$ in
Theorem~\ref{thm:drinfeld-finite-factors}, in particular to
the pole stages $F_P\mathscr LX$.
\end{remark}

\begin{lemma}[Nice overlaps]
\label{lem:nice-cech-overlaps}
Let $S$ be a separated complex scheme and let
$\{u_i:U_i\to S\}_{i\in I}$ be finitely presented \'etale
maps with each $U_i$ nice. For every finite tuple
$(i_0,\ldots,i_n)$, the scheme
\[
 U_{i_0\cdots i_n}
 :=U_{i_0}\times_S\cdots\times_S U_{i_n}
\]
is affine and nice. 
\end{lemma}
\begin{proof}
Since $S$ is separated, the natural morphism
\[
 U_{i_0\cdots i_n}
 \longrightarrow
 U_{i_0}\times_{\C}\cdots\times_{\C}U_{i_n}
\]
is a closed immersion. Its target is affine, so its source
$W:=U_{i_0\cdots i_n}$ is affine.
The projection $W\to U_{i_0}$ is finitely presented and
\'etale, being a finite composite of base changes of chart maps.
Write $U_{i_0}=\Spec A[x_\alpha:\alpha\in J]$, where $A$
is a finitely generated $\C$-algebra. Finite presentation gives
\[
 \mathcal O(W)
 \cong A[x_\alpha:\alpha\in J][t_1,\ldots,t_m]/
       (f_1,\ldots,f_r).
\]
Choose a finite subset $J_0\subset J$ containing every index
$\alpha$ for which $x_\alpha$ occurs in one of the $f_\ell$,
and put
\[
 A_0:=A[x_\alpha:\alpha\in J_0][t_1,\ldots,t_m]/
       (f_1,\ldots,f_r).
\]
Then $A_0$ is finitely generated over $\C$ and
\[
 \mathcal O(W)\cong A_0[x_\alpha:\alpha\in J\setminus J_0],
 \qquad W\cong\Spec A_0\times\A^{J\setminus J_0}.
\]
Thus $W$ is nice. 
\end{proof}

All face and degeneracy maps between these overlaps are also
finitely presented \'etale. Faces are composites of base changes
of chart maps. A degeneracy is a base change of the diagonal of
such a map. This diagonal is an open immersion because the map
is \'etale, and is quasi-compact because the map is finitely
presented; hence the diagonal is a finitely presented \'etale map.

\subsection{Volume indices}
\label{subsec:volume-index-preparations}

We specify the finite-factor comparisons and compatible finite-factor
presentations needed to define volume indices on the pole stages.

\emph{Comparison of finite factors.} Recall the definition of finite factors from \ref{subsec:drinfeld-finite-factors}.
For two finite factors $\pi:U\to Y$ and $\pi':U\to Y'$, we say
that $(Y',\pi')$ \emph{dominates} $(Y,\pi)$ if there is a morphism
$f:Y'\to Y$ such that $\pi=f\circ\pi'$. We call $f$ a
\emph{comparison morphism}, or a \emph{domination morphism}, between
these finite factors. 
 A \emph{common finite factor dominating both} is a
finite factor $q:U\to W$ with morphisms
\[
 p:W\to Y,\qquad p':W\to Y',\qquad
 \pi=p\circ q,\quad \pi'=p'\circ q.
\]
The following comparison lemma is from
\cite[Section~6.8]{Drinfeld}.

\begin{lemma}\label{lem:finite-factor-comparison}
Any two finite factors of the same nice chart $U$ admit a
common finite factor dominating both. If a comparison morphism
between two such finite factors exists, it is unique, surjective,
and smooth. 
\end{lemma}
\begin{proof}
Write $Y=\Spec A$ and $U=\Spec A[u_j:j\in J]$ using the first
finite-factor presentation. Since $Y'$ is affine of finite type, the
pullbacks under $\pi'$ of finitely many generators of its coordinate
ring involve only a finite subset $J_0\subset J$. Consequently
$\pi'$ factors through
\[
 W=\Spec A[u_j:j\in J_0]=Y\times\A^{J_0}.
\]
The projection $q:U\to W$ is itself a finite factor, since
$U\cong W\times\A^{J\setminus J_0}$. Its projection to $Y$ and
the descended map to $Y'$ give $p$ and $p'$ as required.

Now suppose that $\pi=f\pi'$. The product projection $\pi'$ has a section $s:Y'\to U$, so necessarily $f=\pi s$; this proves uniqueness. Surjectivity follows from that of $\pi=f\pi'$.
To prove smoothness, consider a nilpotent extension of affine
schemes $T_0\hookrightarrow T$, a map $T\to Y$, and a compatible
map $T_0\to Y'$. Compose the latter with $s$. The resulting map
$T_0\to U$ lifts over $Y$ to $T\to U$: in the finite-factor presentation
$U=Y\times\A^J$, each affine coordinate can be lifted across the
nilpotent extension. Composing this lift with $\pi'$ gives the
required lift to $Y'$. Thus $f$ is formally smooth. Since $Y$ and
$Y'$ are affine schemes of finite type over $\C$, $f$ is finitely
presented, hence smooth. This also proves smoothness of $p$ and $p'$.
\end{proof}

Fix the polarized compactification of Section~\ref{sec:setup}.
For its closed pole-stage diagram, choose the coordinates
$g_j=s_j/s_0$ associated with $rA$ and use
\begin{equation}\label{eq:pole-stage-scheme}
 F_P\mathscr LX=\Spec C_{rP}(X)
\end{equation}
as the scheme model of $F_P\mathscr LX$.
Its analytic complex-point space is the pole stage identified in
Proposition~\ref{prop:loop-sectors}. We use the same notation for
the scheme and its analytic complex-point space.

\emph{Retaining data in a finite factor.}
Let $\pi:U\cong Y\times\A^J\to Y$ be a finite-factor
projection. A morphism $h:U\to H$, with $H$ affine of finite
type over $\C$, is \emph{retained in $Y$} if
$h=\bar h\pi$ for a morphism $\bar h:Y\to H$.
This includes regular functions, by taking $H=\A^1$.
A finitely presented closed subscheme $C\subset U$ is
\emph{retained in $Y$} if
\[
 C=U\times_Y C_Y
\]
scheme-theoretically for a closed subscheme $C_Y\subset Y$.
The descended morphism $\bar h$, or closed subscheme $C_Y$,
is unique because $\pi$ is faithfully flat.
Given any initial finite factor and a specified finite collection
of these data, adjoining the finitely many affine coordinates
occurring in pullbacks of algebra generators and in generators
of the defining ideals retains the whole collection.
Additional finite data are accommodated by further enlargement.
For open or locally closed subschemes specified by finitely many
equations and nonvanishing conditions, retention has the same
inverse-image meaning; these conditions descend after a finite
enlargement.
If such a finite-type factor is not affine, use a finite affine open
cover and the induced product charts.

\emph{Finite-factor presentations of closed diagrams.}
Let $S_a$ be the schemes in a finite diagram whose arrows
$S_a\hookrightarrow S_b$ are finitely presented closed immersions.
The principal example is a finite chain of pole schemes
$F_0\mathscr LX\hookrightarrow\cdots\hookrightarrow F_P\mathscr LX$.
Choose nice atlases and finite factors on their charts and iterated
overlaps. If the $S_a$ are separated, the overlaps are nice by
Lemma~\ref{lem:nice-cech-overlaps} and need no further chart covers.
For a chart $U_b\to S_b$, put $U_a=S_a\times_{S_b}U_b$.
Enlarge any initial finite factor of $U_b$ to retain the finitely
many closed subschemes $U_a$. A \emph{finite-factor presentation}
of these pulled-back inclusions is then
\begin{equation}\label{eq:presented-closed-diagram}
 (Y_a\hookrightarrow Y_b)\times\A^J,
 \qquad Y_a\hookrightarrow Y_b\text{ closed and of finite type},
\end{equation}
where $U_b\cong Y_b\times\A^J$ and
$U_a=U_b\times_{Y_b}Y_a$. Thus the inclusions use the same
omitted coordinates, and their squares with the finite-factor
projections are Cartesian. Their open complements are
$(Y_b\setminus Y_a)\times\A^J$, with affine covers as above.

\emph{Finite-factor presentations of chart maps.}
For a finitely presented \'etale map $v:V\to U$ between nice
charts, a \emph{finite-factor presentation of $v$} consists of
finite-factor projections $q:V\to Z$, $\pi:U\to Y$, and a
morphism $f:Z\to Y$ such that
\[
                         fq=\pi v.
\]
The morphism $f$ is determined by $q$ and $\pi$.
Given $\pi$ and any initial finite factor of $V$, such a $q$
exists by retaining $\pi v$. The resulting $f$ is smooth:
as in the proof of Lemma~\ref{lem:finite-factor-comparison},
a section of $q$ transfers the infinitesimal lifting property
of the formally smooth morphism $\pi v$ to $f$, and $f$ is
finitely presented because $Y$ and $Z$ are of finite type over
$\C$. The square with $v$ need not be Cartesian.

 The following is the normalization
associated with a volume form in \cite[Proposition~3.12 and
Section~6]{Drinfeld}. 

\begin{proposition}[Drinfeld]
\label{prop:finite-factor-normalization}
A nowhere-vanishing algebraic volume form $\Omega$ on $X$ supplies
a preferred locally constant function $\delta_Y:Y\to\Z$ for each
finite factor $\pi:U\to Y$ of a nice chart of a pole stage $F_P\sL X$.
The notation suppresses its dependence on $\pi$ and $\Omega$.
It has the following properties:
\begin{enumerate}[(i)]
\item If $\pi=f\pi'$ for a comparison morphism $f:Y'\to Y$,
then $\delta_{Y'}=f^*\delta_Y+e$, where $e$ is the relative complex
dimension of $f$.
\item In a finite-factor presentation of a closed inclusion
$(Y_a\hookrightarrow Y_b)\times\A^J$, one has
$\delta_{Y_a}=\delta_{Y_b}|_{Y_a}$.
\item For a finite-factor presentation
$(q:V\to Z,\ \pi:U\to Y,\ f:Z\to Y)$ of a finitely
presented \'etale chart map $v:V\to U$, one has
\[
 \delta_Z=f^*\delta_Y+e,
\]
where $e:Z\to\Z$ is the locally constant relative complex
dimension of $f$.
\end{enumerate}
The identities in \textup{(i)--(iii)} are compatible with composition
and with restriction to the closed subschemes and open complements
in the finite-factor presentations.
\end{proposition}

Subsection~\ref{subsec:contact-volume-indices} states the contact
formulas for $\delta_Y$; Subsection~\ref{subsec:volume-index-proofs}
proves them together with this proposition.
Proposition~\ref{prop:direct-relative-volume-indices} constructs their
relative extension on the collision families.

\subsection{Coherent descent}
\label{subsec:coherent-descent}
We use finite factors and volume indices to normalize Verdier
dualizing complexes. We construct their comparison maps and trace
maps as coherent diagrams in derived $\infty$-categories, and then
apply analytic descent. The diagrams include the homotopies for
composition and their higher compatibilities on chart overlaps.

\emph{Compatible finite-factor presentations.}
For a finite  string of finitely presented \'etale maps
between nice charts
\[
 U_0\xrightarrow{v_0}U_1\xrightarrow{v_1}\cdots
                         \xrightarrow{v_{m-1}}U_m,
\]
fix, at each occurrence $U_i$, a finite collection $\mathcal C_i$
of closed subschemes $C\hookrightarrow U_i$ whose inclusions
are finitely presented. Include in these collections their
scheme-theoretic inverse images under the maps in the string,
so that
\[
 C'\in\mathcal C_{i+1}
 \quad\Longrightarrow\quad
 U_i\times_{U_{i+1}}C'\in\mathcal C_i.
\]
Only finitely many subschemes are added by this requirement.
If no closed subschemes are specified, take all $\mathcal C_i$
empty.

A closed diagram as in \emph{Finite-factor presentations of closed diagrams}
in Subsection~\ref{subsec:volume-index-preparations} supplies
such collections as follows. Fix one ambient scheme $S$ and
finitely many finitely presented closed immersions
$S_a\hookrightarrow S$, $1\le a\le r$.
Suppose the charts in the string have structure maps
$u_i:U_i\to S$ and the chart maps are over $S$, that is,
\[
 u_{i+1}v_i=u_i.
\]
Set
\[
 C_{a,i}:=U_i\times_S S_a\hookrightarrow U_i,
 \qquad
 \mathcal C_i:=\{C_{a,i}:1\le a\le r\}.
\]
These collections already satisfy the required pullback condition,
since
\[
 U_i\times_{U_{i+1}}C_{a,i+1}
 \cong U_i\times_S S_a=C_{a,i}
\]
as closed subschemes of $U_i$.
Thus the string compares charts over the fixed scheme $S$,
and the collections $\mathcal C_i$ record the chosen closed
subschemes of $S$ on those charts.
For a finite chain of pole stages, take $S$ to be its largest
stage and the $S_a$ to be its specified closed substages.
For a general finite closed diagram, apply this construction
for each fixed target $S_b$, using the specified closed
immersions $S_a\hookrightarrow S_b$.

Then
a \emph{compatible finite-factor presentation} $\mathcal P$ of this
string with these closed subschemes consists of finite-factor
projections $\pi_i:U_i\to Y_i$, morphisms
$f_i:Y_i\to Y_{i+1}$ satisfying
\[
 f_i\pi_i=\pi_{i+1}v_i,
\]
and, for every $C\in\mathcal C_i$, a closed subscheme
$C_{Y_i}\hookrightarrow Y_i$ satisfying
\[
 C=U_i\times_{Y_i}C_{Y_i}
\]
scheme-theoretically. Thus each specified closed inclusion has
a finite-factor presentation as in
\eqref{eq:presented-closed-diagram}.
A factor is chosen for each \emph{occurrence} of a chart;
repeated occurrences need not have the same factor.
A \emph{refinement of $\mathcal P$} is another compatible finite-factor
presentation $\widetilde{\mathcal P}$ of the same string and the
same collections $\mathcal C_i$, with projections
$\widetilde\pi_i:U_i\to\widetilde Y_i$ and descended maps
$\widetilde f_i:\widetilde Y_i\to\widetilde Y_{i+1}$,
together with comparison morphisms $r_i:\widetilde Y_i\to Y_i$
satisfying
\[
 \pi_i=r_i\widetilde\pi_i,\qquad
 r_{i+1}\widetilde f_i=f_i r_i,
\]
and the following condition on the closed factors: for every
$C\in\mathcal C_i$, the corresponding closed factors satisfy
\[
 C_{\widetilde Y_i}=\widetilde Y_i\times_{Y_i}C_{Y_i}.
\]
A compatible finite-factor presentation induces finite-factor presentations of
shorter strings as follows. Retain any nonempty subsequence
of chart occurrences, together with their projections $\pi_i$
and collections $\mathcal C_i$, and compose the intervening
chart maps and descended maps $f_i$.
Repeating an occurrence $U_i$ inserts $\mathrm{id}_{U_i}$;
the induced finite-factor presentation repeats $\pi_i$ and $\mathcal C_i$
and inserts $\mathrm{id}_{Y_i}$ between the repeated finite
factors.

The following lemma extends Drinfeld's category of finite-factor
presentations \cite[Section~6.8]{Drinfeld} to finite strings of
charts with prescribed closed subschemes. It supplies common
refinements of finite-factor presentations of the same string and the indexing
categories and cofinal comparisons used in the coherent descent
construction of Proposition~\ref{prop:normalized-pole-sheaves}.
For small and filtered categories, nerves, and cofinal functors,
see Appendix~\ref{app:derived-colimits}; coherent diagrams are
defined in Subsection~\ref{app:coherent-diagrams}.

\begin{lemma}[Categories of finite-factor presentations]
\label{lem:finite-presentation-categories}
For a nice chart $U$, let $\mathsf F(U)$ be the poset of finitely
generated $\C$-subalgebras $A\subset\Gamma(U,\mathcal O_U)$ for
which the induced map $\pi_A:U\to\Spec A$ is a finite-factor
projection, ordered by inclusion. Let
\[
 \sigma=(U_0\xrightarrow{v_0}U_1\longrightarrow\cdots
                         \xrightarrow{v_{m-1}}U_m)
\]
be a string as above with prescribed closed collections
$\mathcal C_i$. Define $\mathsf P_\sigma$ to be the poset of tuples
$A=(A_0,\ldots,A_m)$ such that $A_i\in\mathsf F(U_i)$,
\[
 v_i^*A_{i+1}\subset A_i\qquad(0\le i<m),
\]
and every $C\in\mathcal C_i$ is retained in $Y_i=\Spec A_i$.
The order is componentwise inclusion, directed toward refinements.
The tuple determines the projections $\pi_i=\pi_{A_i}$ and the
smooth descended maps $f_i:Y_i\to Y_{i+1}$ induced by $v_i^*$.
\begin{enumerate}[(i),wide=0pt]
\item The poset $\mathsf P_\sigma$ is small, nonempty, and filtered.
More precisely, for any finite collection of prescribed finite
factors at each occurrence $U_i$, there is a tuple in
$\mathsf P_\sigma$ dominating every prescribed factor at that
occurrence. In particular, $\mathsf F(U)$ is small and filtered.
If subcollections $\mathcal C_i'\subset\mathcal C_i$ still
satisfy the pullback condition, the inclusion from the
poset of finite-factor presentations for $(\mathcal C_i)_i$ into that for
$(\mathcal C_i')_i$ is cofinal.
\item Write $[m]=\{0\to\cdots\to m\}$. For an order-preserving
map $\alpha:[r]\to[m]$, let $\sigma\alpha$ be the string with
occurrences $U_{\alpha(j)}$, the corresponding composite or
identity maps, and closed collections $\mathcal C_{\alpha(j)}$.
The functor
\[
 R_\alpha:\mathsf P_\sigma\longrightarrow\mathsf P_{\sigma\alpha},
 \qquad R_\alpha(A)=(A_{\alpha(j)})_{0\le j\le r},
\]
is cofinal on nerves. For every order-preserving
$\beta:[s]\to[r]$, one has
\[
 R_\beta R_\alpha=R_{\alpha\circ\beta}.
\]
\item If finite sets of regular functions
$G_i\subset\Gamma(U_i,\mathcal O_{U_i})$ are additionally
specified, put $G=(G_0,\ldots,G_m)$ and let
$\mathsf P_\sigma(G)$ be the full subposet
defined by $G_i\subset A_i$ for all $i$.
Its inclusion into $\mathsf P_\sigma$ is cofinal, and
\textup{(i)--(ii)} remain valid for these subposets when the
functions prescribed at occurrence $j$ of $\sigma\alpha$ are
$G_{\alpha(j)}$.
\end{enumerate}
\end{lemma}

\begin{proof}
A finite-factor projection is faithfully flat, so its pullback
identifies the coordinate ring of its target with an object of
$\mathsf F(U)$. Conversely, each such subalgebra specifies the
projection. The condition $v_i^*A_{i+1}\subset A_i$ supplies
the descended map $f_i$; it is smooth by
\emph{Finite-factor presentations of chart maps} in
Subsection~\ref{subsec:volume-index-preparations}.
Componentwise inclusions give comparison morphisms commuting
with the descended chart maps. Their closed factors satisfy
the required inverse-image identities because those identities
hold after the faithfully flat finite-factor projections.
Thus the order on $\mathsf P_\sigma$ is precisely refinement
of the corresponding finite-factor presentations.

For \textup{(i)}, start at $U_m$ with a factor dominating its
prescribed finite collection, using
Lemma~\ref{lem:finite-factor-comparison}, and enlarge it to retain
every member of $\mathcal C_m$. Proceed backwards. At $U_i$,
choose a factor dominating its prescribed factors and enlarge
it to retain $\mathcal C_i$ and the pullbacks of generators of
the algebra just chosen at $U_{i+1}$. Each enlargement uses
only finitely many affine coordinates: the closed immersions
in $\mathcal C_i$ have finitely generated defining ideals,
and the next factor has a finitely generated coordinate ring.
This is the finite-retention procedure of
Subsection~\ref{subsec:volume-index-preparations}.
It gives a tuple in $\mathsf P_\sigma$ dominating all the
prescribed factors. Taking empty prescribed collections proves
nonemptiness; prescribing the factors of finitely many tuples
proves filteredness. Smallness follows because the algebras
are subalgebras of fixed rings. The case of a one-occurrence
string with no closed subschemes gives the assertion for
$\mathsf F(U)$. For the inclusion obtained by forgetting closed
subschemes, fix a tuple satisfying the remaining conditions.
The tuples retaining all the original closed subschemes and
dominating this tuple form a nonempty filtered poset by the
same construction. Its nerve is contractible, which proves
cofinality of the inclusion.

For \textup{(ii)}, the pullback condition on the closed collections
is preserved under composition and identity maps, so
$R_\alpha$ is well defined. For $B=(B_0,\ldots,B_r)$ in
$\mathsf P_{\sigma\alpha}$, its comma category with $R_\alpha$
is the poset
\[
 (B\downarrow R_\alpha)
 =\{A\in\mathsf P_\sigma:
           B_j\subset A_{\alpha(j)}\text{ for every }j\}.
\]
Prescribe each $B_j$ at occurrence $\alpha(j)$ of $\sigma$.
Part~\textup{(i)} makes this poset nonempty. The same construction,
also prescribing the factors of any finite collection of its
objects, gives a common upper bound in it. It is therefore
filtered and has contractible nerve, which proves cofinality.
The composition identity follows by selecting the same tuple
entries on both sides.

For \textup{(iii)}, add the finitely many functions in $G_i$
to the data retained at $U_i$ in the backwards construction.
It still dominates every finite collection of prescribed
factors. In particular, for each $A\in\mathsf P_\sigma$, the
poset of objects of $\mathsf P_\sigma(G)$ dominating $A$ is
nonempty and filtered. This proves cofinality of the inclusion.
With the induced lists $G_{\alpha(j)}$, the comma-category
argument in \textup{(ii)} applies unchanged.
\end{proof}

\emph{Comparing different atlases.}
Lemma~\ref{lem:finite-presentation-categories} compares finite-factor presentations
of a fixed string of charts. To compare constructions using two
nice \'etale atlases $W_0\to S$ and $W'_0\to S$, a
\emph{common atlas refinement} means a nice \'etale atlas
$V_0\to S$ with maps to both atlases over $S$, given on each
chart $V$ by finitely presented \'etale maps
$v_\alpha:V\to U_\alpha$, where $U_1$ and $U_2$ are charts
of $W_0$ and $W'_0$, respectively.
Such a refinement is obtained by taking a nice \'etale atlas
of $W_0\times_S W'_0$.
If $S$ is separated, Lemma~\ref{lem:nice-cech-overlaps}
allows us to take the charts $U_1\times_S U_2$ themselves.

For prescribed finite-factor projections
$\pi_\alpha:U_\alpha\to Y_\alpha$, a
\emph{common refinement of the presented atlases}
also includes, on each such $V$, a finite-factor projection
$q:V\to Z$ and morphisms $p_\alpha:Z\to Y_\alpha$ satisfying
\[
 p_\alpha q=\pi_\alpha v_\alpha,\qquad \alpha=1,2.
\]
These data exist by Lemma~\ref{lem:finite-presentation-categories}(iii),
applied to the one-occurrence string $V$ with the pullbacks
of generators of both coordinate rings prescribed.
The maps $p_\alpha$ are smooth by
\emph{Finite-factor presentations of chart maps} in
Subsection~\ref{subsec:volume-index-preparations}.

If a specified finitely presented closed substage
$S_a\hookrightarrow S$ is presented on $U_\alpha$ by
$Y_{\alpha,a}\hookrightarrow Y_\alpha$, its pullback to $V$
is presented by
\[
 Z_a=Z\times_{Y_1}Y_{1,a}
     =Z\times_{Y_2}Y_{2,a}
\]
as closed subschemes of $Z$.
The equality follows after the faithfully flat pullback by $q$.
Thus the comparison squares for closed factors are Cartesian,
and the maps $Z_a\to Y_{\alpha,a}$ are smooth;
their open complements are compared by restriction.
The squares presenting $v_\alpha$ need not be Cartesian.

For an infinite pole-stage diagram, these constructions are
applied to finite chains of stages and finite strings of chart
maps, with further refinements comparing the choices.
The coherent comparisons of normalized complexes, including
independence of the atlas, are constructed in the proof of
Proposition~\ref{prop:normalized-pole-sheaves}.

\emph{Analytic realization.}
For an affine complex scheme $S=\Spec A$, give $S(\C)$ the initial
topology for the evaluation maps $h:S(\C)\to\C$, $h\in A$, where
the target $\C$ is endowed with the analytic topology.
For the pole schemes this is exactly their analytic topology from
Section~\ref{subsec:loop-topology}; for a nice chart it is the
product topology on $Y(\C)\times\C^J$.

Every finitely presented \'etale chart induces a local
homeomorphism $U(\C)\to S(\C)$. Indeed, write $A$ as the filtered
union of its finitely generated complex subalgebras. Finite
presentation and \'etaleness descend to one of them, giving
\[
 U=U_0\times_{S_0}S,\qquad
 U(\C)\cong U_0(\C)\times_{S_0(\C)}S(\C),
\]
where $U_0\to S_0$ is \'etale between finite-type complex schemes;
see \cite[Tags~01ZM and~07RP]{Stacks}. The second identification
is a homeomorphism: evaluation on a tensor-product algebra gives
precisely the fiber-product topology. The analytic implicit-function
theorem applied to an \'etale presentation makes
$U_0(\C)\to S_0(\C)$ a local analytic isomorphism, and local
homeomorphisms remain such under topological base change.
The same argument works locally when $S$ is not affine.

A jointly surjective \'etale atlas covers all complex points,
since a nonempty \'etale fiber over $\C$ has a $\C$-point.
Restricting each chart to neighborhoods on which its map is a
homeomorphism therefore gives an ordinary analytic open cover.

\emph{Sheaf categories and finite-dimensional operations.} The derived $\infty$-category   obtained
from cochain complexes by making quasi-isomorphisms equivalences is denoted $\mathcal D(\bQ)$, cf. Appendix \ref{app:complexes-duality}.
For a topological space $T$, let $\mathcal D(T;\bQ)$ be the
derived $\infty$-category obtained from cochain complexes of
sheaves of $\bQ$-vector spaces by inverting quasi-isomorphisms.
It is equivalently the category of hypersheaves with values in
$\mathcal D(\bQ)$; see \cite[Proposition~7.1 and its proof]{ScholzeSixFunctors}.
Write $D(T;\bQ)=h\mathcal D(T;\bQ)$ for its ordinary derived
category, and $\mathcal D_{\mathrm{per}}(T;\bQ)$ for the
corresponding category of two-periodic sheaf complexes. No
constructibility condition is imposed. For a complex scheme $Y$
we abbreviate $\mathcal D(Y(\C);\bQ)$ to $\mathcal D(Y;\bQ)$,
and similarly in parity grading.

For the periodic category we use hypersheaves valued in
$\mathcal D(\Lambda)$, the derived $\infty$-category of
differential graded modules over $\Lambda=\bQ[u,u^{-1}]$,
with $|u|=2$ and $du=0$; see
\emph{Sheaf complexes and perfect duality} in
Subsection~\ref{app:complexes-duality}.

On finite-type separated complex schemes and the locally closed
subanalytic subsets used below, fix the enhanced topological
six-functor formalism of \cite[Theorem~7.4, Remark~7.5, and
Proposition~7.7]{ScholzeSixFunctors}. These analytic spaces are
locally compact Hausdorff, paracompact, and of finite covering
dimension; hence the sheaf and hypersheaf versions agree by
\cite[Proposition~7.1]{ScholzeSixFunctors}. We use the following
operations and natural transformations in these enhanced categories.
For a map $f:V\to Y$, the adjunctions are
\[
 \begin{gathered}
 f^{-1}:\mathcal D(Y;\bQ)\longrightarrow\mathcal D(V;\bQ),
       \qquad f^{-1}\dashv Rf_*,\\
 Rf_!: \mathcal D(V;\bQ)\longrightarrow\mathcal D(Y;\bQ),
       \qquad Rf_!\dashv f^!.
 \end{gathered}
\]
Here $f^{-1}$ is ordinary sheaf inverse image, $Rf_!$ is derived
direct image with proper support, and $f^!$ is its right adjoint,
called exceptional inverse image. In particular, their units
and counits include
\[
 E\longrightarrow Rf_*f^{-1}E,
 \qquad Rf_!f^!E\longrightarrow E.
\]
If $f$ is proper, $Rf_!=Rf_*$. For $g:Y'\to Y$, put
$V'=V\times_Y Y'$, and let $f':V'\to Y'$ and $g':V'\to V$
be the projections. Proper base change is the equivalence
\[
 g^{-1}Rf_*\simeq Rf'_*g'^{-1}\qquad(f\text{ proper}).
\]
For a closed embedding $i:Z\hookrightarrow Y$ and its open
complement $j:O\hookrightarrow Y$, localization is the natural
cofiber sequence
\[
 i_*i^!E\longrightarrow E\longrightarrow Rj_*j^{-1}E.
\]
Cofibers and cofiber sequences are defined in
\emph{Cofiber sequences of complexes} in
Subsection~\ref{app:homotopy-limits-colimits}; their images in $D(Y;\bQ)$ are
distinguished triangles.

For a smooth complex morphism $f:V\to Y$ of relative dimension
$e$, let $\operatorname{or}_f$ be the orientation local system
of its real relative tangent bundle. Its complex structure gives
the canonical orientation, hence
$\operatorname{or}_f\simeq\bQ_V$. The resulting
\emph{purity transformation} for $f$ is the following
equivalence in $\mathcal D(V;\bQ)$, natural in
$E\in\mathcal D(Y;\bQ)$:
\begin{equation}\label{eq:oriented-smooth-purity}
 f^!E\simeq f^{-1}E\otimes\operatorname{or}_f[2e]
       \simeq f^{-1}E[2e].
\end{equation}
For smooth $g:W\to V$ of relative dimension $d$, the exact
sequence of relative tangent bundles identifies
$\operatorname{or}_{f\circ g}$ with
$\operatorname{or}_g\otimes g^{-1}\operatorname{or}_f$.
Thus the purity equivalences identify transitivity
$(f\circ g)^!\simeq g^!f^!$ with the shift $[2(d+e)]$;
they also respect base change. The enhanced formalism supplies
these transformations with their composition homotopies and
all higher compatibilities.

For an arbitrary continuous map $a:T'\to T$, ordinary sheaf
inverse image is exact and induces
\[
 a^{-1}:\mathcal D(T;\bQ)\longrightarrow\mathcal D(T';\bQ).
\]
In the hypersheaf description this is hypersheafified inverse
image. In particular, a nice chart will use below only
\[
 p^{-1}:\mathcal D(Y(\C);\bQ)\longrightarrow
      \mathcal D(Y(\C)\times\C^J;\bQ),
 \qquad p:Y(\C)\times\C^J\to Y(\C),
\]
to transport finite-factor complexes and their maps; no
exceptional inverse image $p^!$ is used. For infinite $J$, this
product need not be locally compact, so the six-functor
theorem is not applied to it. Lemma~\ref{lem:infinite-product-descent}
supplies the additional product comparisons and descent needed
on these infinite-dimensional spaces.

On arbitrary topological spaces, we also use the right adjoint
$Ra_*$ of $a^{-1}$ and open--closed localization.
More precisely, let $i:Z\hookrightarrow T$ be a closed embedding
and $j:T\setminus Z\hookrightarrow T$ its open complement.
The exact, fully faithful functor $i_*$ has a right adjoint
$i^!$, given, for $E\in\mathcal D(T;\bQ)$, by
\[
 i^!E\simeq
 i^{-1}\operatorname{fib}
 \bigl(E\longrightarrow Rj_*j^{-1}E\bigr),
\]
where $\operatorname{fib}$ denotes the fiber in the derived
$\infty$-category. Indeed, this fiber restricts to zero on
$T\setminus Z$, so it lies in the essential image of $i_*$;
the adjunctions $i^{-1}\dashv i_*$ and $j^{-1}\dashv Rj_*$
then identify the displayed functor as the right adjoint of $i_*$.
Consequently, the localization sequence
\[
 i_*i^!E\longrightarrow E\longrightarrow Rj_*j^{-1}E
\]
holds without local compactness or finite-dimensionality.
The same constructions apply with periodic coefficients.
Where the six-functor formalism above applies, these functors
agree with its functors by uniqueness of right adjoints.

\emph{Representatives of diagrams.}
We construct diagrams and their morphisms in the derived
$\infty$-categories, including their composition homotopies and
higher compatibilities. Propositions~\ref{prop:diagram-model-comparison}
and~\ref{prop:simultaneous-diagram-representatives} in
Appendix~\ref{app:diagram-representatives} allow joint
representatives by complexes of sheaves and cochain maps,
in either grading, after replacing complexes when necessary.
In these strict representatives, the composition and identity
relations hold as equalities of cochain maps; the comparison
with the original diagram is an equivalence of entire coherent
diagrams. We use the same assertion for diagrams of complexes
of vector spaces after derived global sections.


\begin{lemma}[Product pullback and descent]
\label{lem:infinite-product-descent}
All equivalences and functors below are in the derived
$\infty$-categories just defined, and limits are limits of
$\infty$-categories. Both assertions also hold with
two-periodic coefficients.
\begin{enumerate}[(i),wide=0pt]
\item Let $Y$ be a topological space and let
$p:Y\times\C^J\to Y$, with $J$ finite or countable.
For an open embedding $j:O\hookrightarrow Y$ and a closed
embedding $i:Z\hookrightarrow Y$, write $j',i'$ for their
products with $\C^J$ and $p_O,p_Z$ for the projections.
For $K\in\mathcal D(O;\bQ)$ and $L\in\mathcal D(Z;\bQ)$,
the canonical comparisons are equivalences:
\begin{equation}\label{eq:infinite-product-base-change}
 p^{-1}Rj_*K\simeq Rj'_*p_O^{-1}K,
 \qquad p^{-1}i_*L\simeq i'_*p_Z^{-1}L.
\end{equation}
For $E\in\mathcal D(Y;\bQ)$, the adjunction unit
$E\to Rp_*p^{-1}E$ is an equivalence. Moreover,
$p_Z^{-1}i^!E\simeq i'^!p^{-1}E$, where $i^!$ and $i'^!$
are right adjoints of the closed direct images. These
comparisons are natural in the complexes and compatible
with successive product projections.

\item Let $S$ be a topological space and let $u:U\to S$ be a
surjective local homeomorphism. Write $U_\bullet\to S$ for its
augmented \v Cech nerve, indexed by $\Delta^{\mathrm{op}}$,
as in Definition~\ref{def:cech-nerve}.
For every small ordinary category $A$, restriction induces
an equivalence
\begin{equation}\label{eq:diagram-cech-descent}
 \operatorname{Fun}(N(A),\mathcal D(S;\bQ))
 \simeq
 \lim_{[n]\in\Delta}
 \operatorname{Fun}(N(A),\mathcal D(U_n;\bQ)).
\end{equation}
In particular, objects, arrows ($A=[1]$), and their
coherent diagrams descend. Cofiber sequences, in the sense of
\emph{Cofiber sequences of complexes} in
Subsection~\ref{app:homotopy-limits-colimits}, descend as cofiber sequences.
\end{enumerate}
\end{lemma}
\begin{proof}
(i) A basic neighborhood in $\C^J$ can be chosen to be a convex
cylinder $W$: finitely many coordinates lie in open discs and
the others are unrestricted. Such a $W$ is weakly contractible
and locally weakly contractible. The homotopy-invariance theorem
\cite[Theorem~0.6]{HainePortaTeyssier} applies to an arbitrary
topological base and gives, for every open $V\subset Y$,
\[
 R\Gamma((V\cap O)\times W,p_O^{-1}K)
 \simeq R\Gamma(V\cap O,K).
\]
Taking stalks over these product neighborhoods proves the first
comparison in \eqref{eq:infinite-product-base-change}. The
second follows from the stalkwise description of closed direct
image. Applying the same theorem to $\C^J$ gives the unit
equivalence. For the assertion about $i^!$, use the open
complement $j:Y\setminus Z\hookrightarrow Y$ and pull back
the localization cofiber sequence
$i_*i^!E\to E\to Rj_*j^{-1}E$.
The two product comparisons identify it with the localization
sequence on $Y\times\C^J$; full faithfulness of $i'_*$ gives
the result. These are the canonical base-change transformations,
so they have the composition compatibilities of those transformations.
The homotopy-invariance theorem allows arbitrary presentable (in the sense of \cite[p.~12]{HA})
coefficient categories, including thus periodic complexes.

(ii) Hypersheaf descent for a surjective local homeomorphism
gives $\mathcal D(S;\bQ)\simeq\lim_\Delta\mathcal D(U_n;\bQ)$;
see \cite[Observations~3.19--3.20]{HainePortaTeyssier}.
Applying $\operatorname{Fun}(N(A),-)$, which preserves limits,
gives \eqref{eq:diagram-cech-descent}. Restriction to the charts
is exact and jointly conservative; hence it also detects
cofiber sequences.
\end{proof}

\subsection{Normalized dualizing complexes on pole stages}
\label{subsec:finite-presentations}

The following is a Betti adaptation with $\bQ$ coefficients of
\cite[Section~6.8]{Drinfeld}. We use the enhanced sheaf categories
and operations of Subsection~\ref{subsec:coherent-descent}, whereas Drinfeld uses BBD-type gluing; see also Remark \ref{rmkBBD}.

For a finite factor $Y$, let
$\omega_Y=a_Y^!\bQ\in\mathcal D(Y;\bQ)$ be its Verdier
dualizing complex, where $a_Y:Y(\C)\to\mathrm{pt}$.
We write $\operatorname{Per}_2$ for direct-sum two-periodization,
as defined in \emph{Parity grading and periodization} in
Subsection~\ref{app:complexes-duality}; the same
construction applies to complexes of sheaves.

We also use the enhanced six-functor formalism with coefficients
in $\mathcal D(\Lambda)$; write $a_{Y,\Lambda}^!$ for its
exceptional inverse image along $a_Y$. Its underlying rational
functor is $a_Y^!$ applied to the underlying complex: the
projection formula makes extension of scalars commute with
$Ra_{Y,!}$, and the assertion follows by adjunction.
A \emph{Euclidean neighborhood retract} is a space homeomorphic
to a retract of an open subset of some $\mathbb R^N$.
By local triangulation, every point of $Y(\C)$ has an open
neighborhood $V$ of this kind, including when $Y$ is singular.
For $a_V:V\to\mathrm{pt}$,
\cite[Propositions~7.8 and~6.5(2)]{ScholzeSixFunctors} give
the natural equivalence
\[
 a_V^!\bQ\otimes a_V^{-1}E\simeq a_V^!E
 \qquad(E\in\mathcal D(\bQ)).
\]
For $E=\Lambda$, this is the underlying rational map of
the canonical $\Lambda$-linear coefficient-change map on $V$.
Globally, this map is adjoint to the trace map
\[
 Ra_{Y,!}(\omega_Y\otimes_{\bQ}\Lambda_Y)
 \simeq Ra_{Y,!}\omega_Y\otimes_{\bQ}\Lambda
 \longrightarrow\Lambda,
\]
where the last arrow is the counit
$Ra_{Y,!}a_Y^!\bQ\to\bQ$ tensored with $\Lambda$.
Since these neighborhoods cover $Y(\C)$ and forgetting the
action detects equivalences, this shows
\[
 \operatorname{Per}_2\omega_Y
 \simeq\omega_Y\otimes_{\bQ}\Lambda_Y
 \xrightarrow{\sim}a_{Y,\Lambda}^!\Lambda.
\]
Thus the periodic finite-factor complexes are dualizing objects
in this coefficient formalism. Their traces, localization
sequences, and smooth-purity comparisons below are taken in
that formalism.

\begin{proposition}
\label{prop:normalized-pole-sheaves}
For every polarized SNC pair $(M,D,A)$ and $P\ge0$, put
$S_P=F_P\mathscr LX(\C)$. There is a canonical normalized
dualizing complex
$\omega_P^{\mathrm{per}}\in\mathcal D_{\mathrm{per}}(S_P;\bQ)$.
A choice of $\Omega$ also gives
$\omega_{P,\Omega}\in\mathcal D(S_P;\bQ)$.
For a nice chart $h:U\to F_P\mathscr LX$ and a finite factor
$\pi:U\to Y$, their local descriptions are
\begin{equation}\label{eq:nice-chart-dualizing}
 h^{-1}\omega_{P,\Omega}
 \simeq\pi^{-1}\!\left(\omega_Y[-2\delta_Y]\right),
\end{equation}
and
\[
 h^{-1}\omega_P^{\mathrm{per}}
 \simeq\pi^{-1}\operatorname{Per}_2\omega_Y.
\]
These equivalences lie in $\mathcal D(U;\bQ)$ and
$\mathcal D_{\mathrm{per}}(U;\bQ)$, respectively; inverse images
are taken on analytic complex-point spaces.
The variable shift is taken componentwise: on a connected component
$Y_a$ where $\delta_Y=d$, it is $\omega_{Y_a}[-2d]$.
These descriptions are canonically independent of the finite factors
and the nice atlases, compatibly with refinements and chart overlaps.
When $\Omega$ exists,
\[
          \operatorname{Per}_2\omega_{P,\Omega}
                    \simeq\omega_P^{\mathrm{per}}.
\]
The parity object requires no volume form and no choice of an index
modulo $2$.

For every closed stage inclusion $i_{P,Q}:S_P\hookrightarrow S_Q$,
$P\le Q$, proper trace gives canonical maps in
$\mathcal D(S_Q;\bQ)$,
\[
 (i_{P,Q})_*\omega_{P,\Omega}\longrightarrow\omega_{Q,\Omega}
\]
and their counterparts in $\mathcal D_{\mathrm{per}}(S_Q;\bQ)$,
compatible with composition and finite-factor presentations as coherent
diagrams. For $P\ge1$, write
$i:F_{P-1}\mathscr LX\hookrightarrow F_P\mathscr LX$ and
$j:U_P=F_P\mathscr LX\setminus F_{P-1}\mathscr LX
\hookrightarrow F_P\mathscr LX$. These maps fit into the
localization cofiber sequence in $\mathcal D(S_P;\bQ)$,
\[
 i_*\omega_{P-1,\Omega}\longrightarrow\omega_{P,\Omega}
 \longrightarrow Rj_*j^{-1}\omega_{P,\Omega}\xrightarrow{+1},
\]
and the analogous sequence in $\mathcal D_{\mathrm{per}}(S_P;\bQ)$.
The displayed triangle denotes the image of this cofiber
sequence in the ordinary derived category.
\end{proposition}

\begin{proof}
\emph{Finite-factor presentations of the pole schemes.}
We use the individual finite-factor presentations of
Subsection~\ref{subsec:volume-index-preparations} and the compatible
finite-factor presentations and enhanced sheaf operations of
Subsection~\ref{subsec:coherent-descent}.
In the affine coordinates $g_j=s_j/s_0$ let $f_\nu$ be defining
polynomials of $X$. The ring of $F_P\mathscr L X$ is precisely
\eqref{eq:coefficient-ring} with lower coefficient index $-rP$.
Inside the ring for $F_Q\mathscr L X$, $Q\ge P$, its ideal is
\[
 (g_{j,a}:1\le j\le m,\ -rQ\le a<-rP).
\]
It is finitely generated. Every quasi-compact closed subscheme
of $\mathscr LX$ is contained in a bounded stage: choose a finite
affine cover and a common Laurent bound for the finitely many
affine coordinate functions on its members. In a closed subscheme
of that stage containing $F_P\mathscr LX$, the ideal of
$F_P\mathscr LX$ is the image of the displayed finitely generated
ideal. This proves reasonableness. It also verifies the local
finite-generation condition in Drinfeld's definition for a
containing closed subscheme that is not quasi-compact: each of
its affine opens has a bounded Laurent presentation. Increase
that bound to at least $rP$; the same coefficient equations then
generate the ideal of its intersection with $F_P\mathscr LX$.
Theorem~\ref{thm:drinfeld-finite-factors}, applied to the smooth
affine $\C$-variety $X$, therefore gives the required nice atlases.

For each finite chain of pole inclusions, apply the
finite-data enlargement of \emph{Retaining data in a finite factor}
in Subsection~\ref{subsec:volume-index-preparations} to the pullbacks
of its negative-coefficient equations on a nice chart
$U=\Spec A[u_j:j\in J]$ of the largest stage.
For a finite set $J_0\subset J$ containing the coordinates in
those equations, the factor $Y=\Spec A[u_j:j\in J_0]$ presents
the chain as $(Y_a\subset Y)\times\A^{J\setminus J_0}$.
The overlaps are nice by Lemma~\ref{lem:nice-cech-overlaps}.
Present their finite strings by
Lemma~\ref{lem:finite-presentation-categories}, prescribing the
pulled-back closed chain at each occurrence. The descended
chart maps are smooth, and the closed comparison squares are
Cartesian. Equality of descended composites follows
after pullback along the faithfully flat source projection.
We now use the categories of finite-factor presentations of
Lemma~\ref{lem:finite-presentation-categories} to construct
the coherent maps of their normalized complexes.

\emph{Enhanced finite-factor comparisons.}
Write $\omega_Y^{\mathrm n}=\omega_Y[-2\delta_Y]\in\mathcal D(Y;\bQ)$, with
shifts taken on the open-and-closed level sets of $\delta_Y$.
The superscript $\mathrm n$ records the specified index normalization.
Its dependence on the finite-factor projection and on $\Omega$
is suppressed, as in $\delta_Y$.
For a smooth comparison $f:Y'\to Y$ of relative complex
dimension $e$, the purity transformation
\eqref{eq:oriented-smooth-purity} and
Proposition~\ref{prop:finite-factor-normalization}(i) give
\begin{equation}\label{eq:enhanced-normalized-purity}
 \omega_{Y'}^{\mathrm n}\simeq f^!\omega_Y[-2\delta_{Y'}]
 \simeq f^{-1}\omega_Y[2e-2\delta_{Y'}]
 \simeq f^{-1}\omega_Y^{\mathrm n}.
\end{equation}
We call \eqref{eq:enhanced-normalized-purity} the
\emph{normalized purity equivalence} for $f$. Ordinary inverse
image under $\pi':U\to Y'$ gives the comparison on the chart
\[
 \pi'^{-1}\omega_{Y'}^{\mathrm n}\simeq\pi'^{-1}f^{-1}\omega_Y^{\mathrm n}
       \simeq\pi^{-1}\omega_Y^{\mathrm n},\qquad \pi=f\pi'.
\]
For the small filtered poset
$\mathsf F(U)$ of Lemma~\ref{lem:finite-presentation-categories}(i),
these pulled-back normalized purity equivalences, inverted when
necessary, define a functor
\[
 \mathsf F(U)\longrightarrow\mathcal D(U;\bQ),
 \qquad A\longmapsto\pi_A^{-1}\omega_{\Spec A}^{\mathrm n},
\]
which sends every arrow to an equivalence. Composition and
higher compatibilities are inherited from the purity
transformations of Subsection~\ref{subsec:coherent-descent}
and the compatible volume-index identities.

Its colimit can be computed by an explicit sequence. Fix
$U=Y\times\A^J$ and finite subsets
$J_0\subseteq J_1\subseteq\cdots$ with union $J$; put
$Y_n=Y\times\A^{J_n}$ and let $\pi_n:U\to Y_n$ be the
projection. Every finite factor $\pi':U\to Y'$ factors
through some $\pi_n$, since generators of $\mathcal O(Y')$
involve finitely many coordinates. Consequently the chain
$\{\mathcal O(Y_n)\}_n$ is cofinal in $\mathsf F(U)$, and
the chart object is
\[
 \omega_U^{\mathrm n}:=\operatorname*{colim}_{n\ge0}\pi_n^{-1}\omega_{Y_n}^{\mathrm n}
       \quad\text{in }\mathcal D(U;\bQ).
\]
The colimit cone identifies $\omega_U^{\mathrm n}$ with each finite-factor
term. Common dominating factors compare different initial
finite-factor presentations and exhaustions. This system is countable;
it need not have a finite cofinal subsystem, even when the
atlas has one chart.

We now construct the chart comparisons together with their
higher compatibilities. Fix a finite string
\[
 \sigma=(U_0\xrightarrow{v_0}U_1\xrightarrow{v_1}\cdots
                    \xrightarrow{v_{m-1}}U_m)
\]
of finitely presented \'etale maps between nice charts, and use
its poset of finite-factor presentations $\mathsf P_\sigma$ from
Lemma~\ref{lem:finite-presentation-categories}.
When presenting a finite chain of pole stages, use chart maps
over its largest stage and choose the collections $\mathcal C_i$
by the pullback construction in \emph{Compatible finite-factor presentations}
in Subsection~\ref{subsec:coherent-descent}.
If finite sets $G_i$ of regular functions on $U_i$ are additionally
specified, the constructions below also apply on the cofinal
full subposets $\mathsf P_\sigma(G)$ of part~(iii) of that lemma,
with the induced lists of functions under omission and repetition.

For each tuple, smooth purity and
Proposition~\ref{prop:finite-factor-normalization}(iii) give
the normalized purity equivalences
\[
                         \omega_{Y_i}^{\mathrm n}\simeq f_i^{-1}\omega_{Y_{i+1}}^{\mathrm n}.
\]
Pull these equivalences to $U_0$, using the composites
$U_0\to U_i$. They form a diagram of $m$ composable equivalences
in $\mathcal D(U_0;\bQ)$. An enlargement of tuples induces
smooth comparison maps between their factors. Although the
squares for these maps need not be Cartesian, their two
composites agree. More explicitly, write $Y_i'$ for the enlarged
factors, $h_i:Y_i'\to Y_i$ for their comparison maps, and
$f_i':Y_i'\to Y_{i+1}'$ for the enlarged adjacent maps. Then
$f_i h_i=h_{i+1}f_i'$, and transitivity gives
\[
 h_i^!f_i^!\simeq(f_i h_i)^!
 \simeq(h_{i+1}f_i')^!\simeq(f_i')^!h_{i+1}^!.
\]
Transitivity of complex orientations identifies the corresponding
normalized purity equivalences, with the composition homotopies
and higher compatibilities inherited from the enhanced formalism. This comparison
uses transitivity, without requiring a base-change equivalence
for a non-Cartesian square. The entire
diagrams form a functor from $\mathsf P_\sigma$ to
$\operatorname{Fun}(N([m]),\mathcal D(U_0;\bQ))$.

The indexing of all these diagrams is explicit. Its objects
are pairs $(\sigma,A)$, with $A\in\mathsf P_\sigma$. A morphism
$(\sigma,A)\to(\tau,B)$, where $\sigma$ is indexed by $[m]$
and $\tau$ by $[n]$, consists of an order-preserving map
$\theta:[m]\to[n]$ such that $\sigma=\tau\theta$, and the
componentwise inclusions $A_i\subset B_{\theta(i)}$ in the
coordinate rings of the corresponding chart occurrences.
Composition composes the order maps and the inclusions. Thus
this category records both refinements and the omissions and
repetitions of occurrences. For such a morphism, inverse image
along the chart composite from the first chart of $\tau$ to
its $\theta(0)$-th chart, followed by the factor comparisons,
compares the $[m]$-diagram for $(\sigma,A)$ with the
$\theta$-restriction of the $[n]$-diagram for $(\tau,B)$.
These transformations compose coherently by the transitivity
of the normalized purity equivalences just described. In particular, this specifies the
comparison when the first vertex is omitted.

For each fixed $\sigma$, take the colimit over its refinement
poset $\mathsf P_\sigma$. Evaluation and inverse image preserve
colimits. Lemma~\ref{lem:finite-presentation-categories}(ii)
gives cofinality of every $R_\theta$. For a one-occurrence
subchain, its poset of finite-factor presentations is also cofinal in
$\mathsf F(U_i)$ by part~(i), forgetting the closed conditions.
Thus restriction identifies the terms with the corresponding
pullbacks of $\omega_{U_i}^{\mathrm n}$.
Consequently the induced colimit comparisons identify faces
and degeneracies with the diagrams for the corresponding
strings. These comparisons are induced by the morphisms of
the displayed indexing category. Functoriality of colimits
therefore preserves their compositions and all the higher
compatibilities inherited from the normalized purity equivalences. This constructs
the coherent overlap diagram used in descent.

\emph{Descent on the pole stage.}
Put $S=F_P\mathscr LX(\C)$ and choose a finite nice atlas
$U=\coprod_i U_i\to S$ as in
Remark~\ref{rem:finite-nice-atlas}. Its analytification is a
surjective local homeomorphism by \emph{Analytic realization} in
Subsection~\ref{subsec:coherent-descent}. On its \v Cech nerve the
preceding construction gives the complexes
$\pi_i^{-1}\omega_{Y_i}^{\mathrm n}$ and their coherent overlap equivalences.
By Lemma~\ref{lem:nice-cech-overlaps}, each \v Cech degree
has finitely many nice overlap schemes $U_{i_0\cdots i_n}$;
all degrees are included. On triple
overlaps, the images of the comparison equivalences in the
ordinary derived category $D(U_{ijk};\bQ)$ satisfy
\begin{equation}\label{eqComp}
 \alpha_{jk}\alpha_{ij}=\alpha_{ik}.
\end{equation}
This equality is the image of the composition homotopy in the
enhanced diagram just constructed. Descent uses that diagram,
including its higher compatibilities, rather than only
\eqref{eqComp}. Lemma~\ref{lem:infinite-product-descent}(ii)
descends it to $\omega_{P,\Omega}\in\mathcal D(S;\bQ)$.
For a second atlas, choose a common refinement of the presented
atlases as in \emph{Comparing different atlases} in
Subsection~\ref{subsec:coherent-descent}. Apply the preceding
coherent comparison and descent constructions to its chart
maps to both atlases. The resulting comparison of the descended
objects proves the asserted independence.

For the parity theory use the enhanced formalism with
periodic coefficients and the canonical even-shift
identifications
\[
 \operatorname{Per}_2\omega_{Y'}
 \simeq f^{-1}\operatorname{Per}_2\omega_Y.
\]
They compose because smooth relative dimensions add and
complex orientations compose. Thus the same descent constructs
$\omega_P^{\mathrm{per}}$ without an integer normalization.
When $\Omega$ exists, the canonical local equivalences
$\operatorname{Per}_2\omega_Y^{\mathrm n}\simeq\operatorname{Per}_2\omega_Y$
commute with these comparisons. For each chart map $h:U\to S$,
exact inverse image commutes with direct sums, so
\[
 h^{-1}\operatorname{Per}_2\omega_{P,\Omega}
 \simeq\operatorname{Per}_2h^{-1}\omega_{P,\Omega}.
\]
Hence these local equivalences descend from
$\operatorname{Per}_2\omega_{P,\Omega}$ itself and give the
stated periodization equivalence. This argument does not exchange
periodization with an infinite \v Cech totalization.

\emph{Descent of proper trace and localization.}
On an ambient-stage chart, present the pulled-back closed pole
inclusion by $i:Z\hookrightarrow Y$
with $\delta_Z=i^*\delta_Y$. Transitivity of exceptional inverse
image and its adjunction counit give the arrow
\[
 i_*\omega_Z^{\mathrm n}\simeq i_*i^!\omega_Y^{\mathrm n}\longrightarrow \omega_Y^{\mathrm n}
\]
in $\mathcal D(Y;\bQ)$. Let $f:Y'\to Y$ be a smooth
comparison arising from a refinement or a chart map. The
finite-factor presentation of the same closed substage gives a Cartesian square with
$Z'=Y'\times_Y Z$, $i':Z'\hookrightarrow Y'$, and
$g:Z'\to Z$. Its normalized purity equivalences and proper
base change identify
\[
 f^{-1}i_*\omega_Z^{\mathrm n}\simeq i'_*g^{-1}\omega_Z^{\mathrm n}\simeq i'_*\omega_{Z'}^{\mathrm n},
 \qquad f^{-1}\omega_Y^{\mathrm n}\simeq \omega_{Y'}^{\mathrm n}.
\]
Under these identifications, the pullback of the counit
$i_*i^!\omega_Y^{\mathrm n}\to \omega_Y^{\mathrm n}$ is the counit
$i'_*i'^!\omega_{Y'}^{\mathrm n}\to \omega_{Y'}^{\mathrm n}$: this is the compatibility of the
adjunction counit with the Cartesian base-change transformation
and smooth purity. Thus the \emph{whole arrow}, including its
source and target, is compatible with comparison. Transitivity
for nested closed embeddings and higher compatibilities are
those of these same enhanced transformations. Pullback to
$Y\times\C^J$ preserves the arrows and their coherent diagrams,
and \eqref{eq:infinite-product-base-change} identifies the
source with the closed direct image from $Z\times\C^J$.

Apply Lemma~\ref{lem:finite-presentation-categories} to the
chart strings with the pullbacks of the finite closed chain
as their prescribed closed collections. Taking
colimits in the corresponding category of arrows constructs
coherent overlap comparisons for the trace arrows. Their
sources identify with the closed direct images of the already
descended normalized objects on the substages, by the closed
base-change comparison in
Lemma~\ref{lem:infinite-product-descent}(i). Apply part~(ii) of
that lemma with $A=[1]$ to descend each entire arrow. For a
finite chain of closed stages, the same construction in the
category of diagrams preserves the compositional counit maps
and their higher compatibilities. Restriction to smaller chains
is cofinal. If the largest ambient stage remains in the chain,
forgetting closed substages gives an inclusion of posets of finite-factor
presentations, cofinal by Lemma~\ref{lem:finite-presentation-categories}(i).
If the ambient stage is removed, one must
also compare finite factors on a closed chart with restrictions of
ambient finite factors. Write the ambient coordinate algebra as
$R=A[x_j:j\in J]$, with $A$ of finite type, and retain the ideal
of the closed chart so that it has the form $I=I_0R$ for an ideal
$I_0\subset A$. For any finite-factor algebra $B\subset R/I$,
lift finitely many algebra generators of $B$ to $R$. These lifts
belong to $A'=A[x_j:j\in F]$ for a finite subset $F\subset J$.
Polynomial extension gives
\[
 I\cap A'=I_0A',\qquad
 R/I\cong(A'/I_0A')[x_j:j\in J\setminus F].
\]
Thus $A'/I_0A'$ is a finite-factor algebra of the closed chart
containing $B$. Lifting generators for any finite collection of
prescribed factors, and then propagating enlargements backwards
along the chart string, makes the corresponding comma categories
nonempty and filtered. This proves cofinality also when the ambient
stage changes. The construction gives the asserted
coherent diagram of pole maps for all $P\le Q$.

On a finite factor, with $j:Y\setminus Z\hookrightarrow Y$,
localization is the functorial cofiber sequence
\[
 i_*\omega_Z^{\mathrm n}\longrightarrow \omega_Y^{\mathrm n}
       \longrightarrow Rj_*j^{-1}\omega_Y^{\mathrm n}\xrightarrow{+1}.
\]
Exact inverse image and
\eqref{eq:infinite-product-base-change} identify its pullback
with the asserted localization sequence on the nice chart.
Its compatibility with refinement follows from that of the
trace arrow and functoriality of cofibers. The last assertion
of Lemma~\ref{lem:infinite-product-descent}(ii) therefore descends
this sequence to $S$. All these constructions apply in the
periodic coefficient category as well. At no point is a
six-functor theorem applied to an infinite-dimensional chart.
\end{proof}

\begin{definition}[Pole-stage loop complexes]
\label{def:pole-stage-complexes}
For $P\ge0$, define
\[
 \begin{aligned}
 B_{\Omega,P}
   &:=R\Gamma(F_P\mathscr LX,\omega_{P,\Omega}),\\
 B_P
   &:=R\Gamma(F_P\mathscr LX,\omega_P^{\mathrm{per}}),
 \end{aligned}
\]
These are objects of $\mathcal D(\bQ)$ and
$\mathcal D_{\mathrm{per}}(\bQ)$, respectively. They are the
\emph{pole-stage loop complexes}, also called the
\emph{normalized Borel--Moore complexes}, in the integer and
parity gradings, respectively. Put $B_{\Omega,-1}=B_{-1}=0$. 

For $0\le P\le Q$, the transition map is the composite
\[
\begin{aligned}
 B_{\Omega,P}
 &\simeq
 R\Gamma\!\left(F_Q\mathscr LX,
               (i_{P,Q})_*\omega_{P,\Omega}\right)\to R\Gamma(F_Q\mathscr LX,\omega_{Q,\Omega})
 =B_{\Omega,Q}.
\end{aligned}
\]
Here $i_{P,Q}:F_P\mathscr LX\hookrightarrow F_Q\mathscr LX$
is the closed stage inclusion. The equivalence follows from
composition of derived direct images along
\[
 F_P\mathscr LX\hookrightarrow F_Q\mathscr LX
 \longrightarrow\mathrm{pt},
\]
with sheaf operations taken on the analytic complex-point
spaces. The arrow is obtained by applying
$R\Gamma(F_Q\mathscr LX,-)$ to the proper-trace map
\[
 (i_{P,Q})_*\omega_{P,\Omega}
 \longrightarrow\omega_{Q,\Omega}
\]
of Proposition~\ref{prop:normalized-pole-sheaves}.
Replacing $\omega_{P,\Omega}$ and $\omega_{Q,\Omega}$ by
$\omega_P^{\mathrm{per}}$ and $\omega_Q^{\mathrm{per}}$
gives the periodic transition map
\[
 B_P\longrightarrow B_Q.
\]
The coherence of the proper-trace maps and of composition
of derived direct images makes these complexes and transition
maps into coherent diagrams in the derived $\infty$-categories
of complexes in the respective gradings. Their descent models
and coherent transition maps are described explicitly in
Remark~\ref{rem:coherent-pole-descent}.
\end{definition}

\begin{remark}\label{rmkBBD}
The sheaf complex $\omega_{P,\Omega}$ is the Betti version with
$\bQ$ coefficients of Drinfeld's renormalized dualizing complex
\cite[Section~6.8]{Drinfeld}. Consequently,
$B_{\Omega,P}$ is the Betti version of the expression
on the right-hand side of \cite[(6.5)]{Drinfeld} before taking
its linear dual. Drinfeld glues in the ordinary derived category
using negative self-Ext vanishing and BBD gluing
\cite[Section~6.8]{Drinfeld}. We can do the same, but we chose to use coherent descent because
for a nonzero two-periodic object the relevant negative self-Ext might not vanish.
Gluing first in integer degrees and then periodizing requires a compatible integer-valued dimension theory, supplied here by \(\Omega\). Our parity theory is constructed without choosing such a dimension theory.
\end{remark}

\begin{remark}[Descent models and coherent pole-stage diagrams]
\label{rem:coherent-pole-descent}
In the volume grading, fix $P$ and choose a finite nice \'etale atlas
$\{h_i:U_i\to F_P\mathscr LX\}$ with finite-factor projections
$\pi_i:U_i\to Y_i$. The chart complexes have the identifications
\[
 \omega_{U_i}^{\mathrm n}\simeq\pi_i^{-1}\omega_{Y_i}^{\mathrm n}.
\]
For $\mathbf i=(i_0,\ldots,i_q)$, let
\[
 U_{\mathbf i}:=
 U_{i_0}\times_{F_P\mathscr LX}\cdots\times_{F_P\mathscr LX}U_{i_q},
 \qquad p_0:U_{\mathbf i}\to U_{i_0}
\]
be the overlap and its first projection. The coherent overlap
comparisons of Proposition~\ref{prop:normalized-pole-sheaves}
give a cosimplicial diagram $C_P^\bullet$ in $\mathcal D(\bQ)$,
with terms
\[
 C_P^q:=\prod_{i_0,\ldots,i_q}
 R\Gamma\bigl(U_{\mathbf i},p_0^{-1}\omega_{U_{i_0}}^{\mathrm n}\bigr).
\]
Its coface and codegeneracy maps come from the projections and
diagonals between overlaps, using those comparisons when the
chosen local complex changes. Each product is finite, but all
\v Cech degrees $q\ge0$ are included. Analytic descent gives
\[
 B_{\Omega,P}
 \xrightarrow{\sim}
 \operatorname*{holim}_{[q]\in\Delta}C_P^q,
\]
where $\Delta$ is the category of nonempty finite ordered sets
and order-preserving maps; this is the homotopy totalization
of Appendix~\ref{app:totalization-realization}.

Put $K_P:=\operatorname*{holim}_{[q]\in\Delta}C_P^q$.
To describe the transition $K_P\to K_Q$ for $P\le Q$,
write $S_R=(F_R\mathscr LX)^{\mathrm{an}}$ for $R=P,Q$.
Choose a finite nice \'etale atlas of $F_Q\mathscr LX$, with
finite factors presenting the inverse image of $F_P\mathscr LX$.
Let $U_Q^q$ be the disjoint union of the analytic realizations
of its $(q+1)$-fold overlaps, put
\[
 U_P^q:=S_P\times_{S_Q}U_Q^q,
\]
and write $h_{R,q}:U_R^q\to S_R$ for the structure maps.
These compatible atlases give cosimplicial section diagrams
\[
 \widetilde C_R^q
   :=R\Gamma(U_R^q,h_{R,q}^{-1}\omega_{R,\Omega}),
   \qquad R=P,Q.
\]
The closed base-change and proper-trace comparisons of
Proposition~\ref{prop:normalized-pole-sheaves} give
\[
 \begin{aligned}
 \widetilde C_P^q
   &\simeq R\Gamma\!\left(
       U_Q^q,h_{Q,q}^{-1}(i_{P,Q})_*\omega_{P,\Omega}\right)\\
   &\longrightarrow
       R\Gamma(U_Q^q,h_{Q,q}^{-1}\omega_{Q,\Omega})
       =\widetilde C_Q^q.
 \end{aligned}
\]
Compatibility with the overlap projections and diagonals makes
these maps a morphism
$\widetilde C_P^\bullet\to\widetilde C_Q^\bullet$
of cosimplicial diagrams. Totalization and the coherent
common-refinement comparisons with the originally chosen
models give
\[
 t_{P,Q}:K_P
   \simeq\operatorname*{holim}_{[q]\in\Delta}\widetilde C_P^q
   \longrightarrow
   \operatorname*{holim}_{[q]\in\Delta}\widetilde C_Q^q
   \simeq K_Q.
\]
Under the descent equivalences $B_{\Omega,R}\simeq K_R$,
this is the transition map of
Definition~\ref{def:pole-stage-complexes}: the first equivalence
in the degreewise formula uses closed base change and composition
of derived direct images, and the arrow is induced by proper trace.

For a finite chain of pole bounds, use an atlas of its largest
stage and present all the pulled-back closed stages on its
finite factors. Transitivity of proper trace and composition
of derived direct images supply specified homotopies
\[
 t_{Q,R}\circ t_{P,Q}\simeq t_{P,R}
 \qquad(P\le Q\le R),\qquad t_{P,P}=\mathrm{id}_{K_P}.
\]
For $P\le Q\le R\le S$, the two resulting homotopies from
$t_{R,S}\circ t_{Q,R}\circ t_{P,Q}$ to $t_{P,S}$ are related
by a specified higher homotopy. The construction supplies the
corresponding higher homotopies for every finite chain,
compatibly with omission and repetition of indices. The cofinal
comparisons in the proof of
Proposition~\ref{prop:normalized-pole-sheaves} give this
compatibility also when the largest stage is omitted.
Thus the equivalences $B_{\Omega,P}\simeq K_P$ identify the
entire pole-stage diagrams.

For a common refinement of the descent models, write
$\widetilde K_P$ and $\widetilde t_{P,Q}$ for the resulting
totalizations and transition maps. The refinement comparisons
give equivalences $r_P:\widetilde K_P\xrightarrow{\sim}K_P$
and homotopies
\[
 r_Q\circ\widetilde t_{P,Q}\simeq t_{P,Q}\circ r_P,
\]
compatible with the composition homotopies, their higher
coherences, and successive refinements. They therefore give
equivalences of the entire pole-stage diagrams.
Lemma~\ref{lem:totalization-realization-representatives} supplies
simultaneous cochain representatives of these totalizations and
their comparison maps.

For the parity grading, replace $\omega_{U_i}^{\mathrm n}$ by
$\pi_i^{-1}\operatorname{Per}_2\omega_{Y_i}$ and perform the
same construction in the derived $\infty$-category of
two-periodic complexes. It gives a cosimplicial diagram
$C_{P,\mathrm{per}}^\bullet$ with
\[
 B_P
 \simeq\operatorname*{holim}_{[q]\in\Delta}C_{P,\mathrm{per}}^q,
\]
with the same composition and refinement coherences, now in
$\mathcal D_{\mathrm{per}}(\bQ)$.
\end{remark}

\subsection{Computing volume indices from contact loci}
\label{subsec:contact-volume-indices}

Fix $\Omega$ and write $\operatorname{div}(\Omega)=-\sum_i a_iD_i$
as in Section~\ref{subsec:polarized-compactifications}. We will prove the next two propositions together with
Proposition~\ref{prop:finite-factor-normalization} below.

\begin{proposition}
\label{prop:contact-volume-indices}
Let $\mathscr P_{\vct,W}$ be the nonempty locus of contact
arcs centered in an affine \'etale SNC coordinate neighborhood
$W\subset M$. For $\ell\ge h(\vct)$, jet truncation is a
finite-factor projection
\[
 \pi_{\vct,\ell}:\mathscr P_{\vct,W}\longrightarrow
 Y_{\vct,\ell}:=j_\ell(\mathscr P_{\vct,W}).
\]
Its target is smooth of dimension
$d_{\vct,\ell}=n(\ell+1)-|\vct|$, and its volume index is
\begin{equation}\label{eq:volume-dimension-index}
 \delta_{\vct,\ell}=n(\ell+1)-\sum_i a_iv_i.
\end{equation}
More generally, every finite factor $Z$ of a nice \'etale chart
of $\mathscr P_\vct$ is smooth and satisfies
\begin{equation}\label{eq:contact-factor-volume-index}
 \delta_Z(z)=\dim_z Z+\tfrac12\mu_\Omega(\vct).
\end{equation}
Here $\dim_z Z$ is the dimension of the smooth component through $z$.
\end{proposition}

\begin{proposition}[Volume index of a general finite factor]
\label{prop:general-factor-contact-index}
Let $h:U\to F_P\mathscr LX$ be a nice chart with finite factor
$\pi:U\to Y$, and let $u\in U(\C)$ correspond to a loop of
contact $\vct$. Write $U\simeq Y\times\A^J$ with
$\pi=\operatorname{pr}_Y$. For a suitable finite subset
$J_0\subset J$, put $e=|J_0|$ and $J'=J\setminus J_0$.
There is an affine open neighborhood $Y'\subset Y\times\A^{J_0}$
of the image of $u$ such that its inverse image $U'\subset U$
under projection to $Y\times\A^{J_0}$ has the finite-factor
presentation
\[
 U'\simeq Y'\times\A^{J'},\qquad
 \pi':U'\to Y',\quad \pi'=\operatorname{pr}_{Y'},\qquad
 \pi|_{U'}=f\pi',
\]
where $f:Y'\to Y$ is the restriction of the projection
$Y\times\A^{J_0}\to Y$. The inverse image of $\mathscr P_\vct$ in
$U'$ is $Z\times\A^{J'}$ for a smooth locally closed reduced
subscheme $Z\subset Y'$. For $z=\pi'(u)$ one has
\begin{equation}\label{eq:contact-volume-general-factor}
 \delta_Y(\pi(u))=\dim_z Z+\tfrac12\mu_\Omega(\vct)-e.
\end{equation}
This value is independent of the finite-factor presentation
compatible with $\pi$ used to obtain $Z$ and $e$. 
\end{proposition}

The following lemmas give explicit coefficient proofs of the
cotangent-quotient description and specialization invariance of indices
underlying Drinfeld's construction
\cite[Lemma~6.5 and its corollary; Proposition~3.12]{Drinfeld}.
Subsection~\ref{subsec:volume-index-proofs} then computes the
resulting volume-normalized indices on contact loci and relates
them to the finite factors used here.

\begin{lemma}
\label{lem:elementary-coefficient-differentials}
Let $X\subset\A^m$ be smooth affine of dimension $n>0$, let
$S\subset\mathscr LX$ be a reasonable closed quasi-compact subscheme, and let
$U=Y\times\A^J\to S$ be a nice chart with projection $\pi$.
For an affine test scheme $T=\Spec A\to U$, write $\gamma_T$
for its family of loops. Define its \emph{Laurent cotangent module}
in the parameter $t$ by
\[
 \Cot_t(\gamma_T):=\gamma_T^*\Omega_X^1
            \otimes_{A(\!( t )\!)} A(\!( t )\!)\,dt,
\]
where $dt$ is a formal generator. The \emph{coefficient-residue map}
$q_{\pi,T}$ associated with $\pi$ and $T\to U$ is the natural
$A$-linear surjection
\begin{equation}\label{eq:elementary-residue-map}
 \begin{aligned}
 q_{\pi,T}:\Cot_t(\gamma_T)&\longrightarrow
                 (\Omega^1_{U/Y})|_T=\bigoplus_{j\in J}A\,du_j,\\
 a(t)\,dg\otimes dt&\longmapsto
      \operatorname{Res}_{t=0}a(t)\,(d_{U/Y}g(\gamma_U(t)))|_T\,dt.
 \end{aligned}
\end{equation}
Here $a(t)\in A(\!(t)\!)$ and $g\in\mathcal O(X)$.
The Laurent residue extracts the coefficient of $t^{-1}dt$. The differential $d_{U/Y}$ is applied coefficientwise to the
Laurent series $g(\gamma_U(t))$ on $U$, before restriction to $T$.
Both Laurent series have lower exponent bounds, so the residue
uses finitely many coefficients and belongs to the displayed
direct sum.
Whenever $\Cot_t(\gamma_T)$ is free over $A(\!( t )\!)$, choose a basis.
In these coordinates, the coefficient-residue kernel
$L_T:=\ker q_{\pi,T}$ satisfies
\begin{equation}\label{eq:elementary-kernel-bounds}
 t^NA\llbracket t \rrbracket^n\subset L_T\subset t^{-N}A\llbracket t \rrbracket^n
\end{equation}
for some $N\ge0$. These are inclusions of $A$-modules;
$L_T$ need not be stable under multiplication by $A\llbracket t \rrbracket$.
For any $A$-algebra $A'$, put $T'=\Spec A'$ and let
$\gamma_{T'}$ be the pullback of $\gamma_T$. The natural map
$\Cot_t(\gamma_T)\to\Cot_t(\gamma_{T'})$ sends $L_T$ into $L_{T'}$.
\end{lemma}

\begin{proof}
Choose $b$ so that $S$ is a closed subscheme of $\Spec C_b(X)$.
Such a bound exists: on each affine open of a finite affine cover
of $S$, its loop family has bounded poles in the finitely many
coordinates of $X$; take the maximum of these bounds. The local
factorizations glue because coefficient stages are monomorphisms.
Write
\[
 X=\Spec\C[x_1,\ldots,x_m]/(f_1,\ldots,f_l),
 \qquad g_i=x_i|_X.
\]
Set
\[
 P_b=\C[x_{i,k}:1\le i\le m,\ k\ge-b],
 \qquad G_i(t)=\sum_{k\ge-b}x_{i,k}t^k,
 \qquad F_{\nu,a}=[t^a]f_\nu(G_1(t),\ldots,G_m(t)).
\]
Thus $C_b(X)=P_b/(F_{\nu,a}:\nu,a)$, as in
\eqref{eq:coefficient-ring}. Let
$\varphi:C_b(X)\to A$ correspond to the composite
$T\to U\to S\to\Spec C_b(X)$.
Put $B=A(\!( t )\!)$, and write
\[
 \gamma_i(t)=\sum_{k\ge-b}\varphi(x_{i,k})t^k,
 \qquad
 v_\nu(t)=
 \left(
 \frac{\partial f_\nu}{\partial x_i}(\gamma(t))
 \right)_{i=1}^m\in B^m.
\]
The conormal presentation, tensored with $A$, gives
\[
 (\Omega^1_{\Spec C_b(X)/\C})|_T
 \cong
 \frac{\displaystyle\bigoplus_{i=1}^m
                  \bigoplus_{k\ge-b}A\,dx_{i,k}}
      {\displaystyle\sum_{\nu,a}A\,r_{\nu,a}},
 \qquad
 r_{\nu,a}
 =
 \sum_{i=1}^m\sum_{k\ge-b}
 [t^{a-k}]
 \frac{\partial f_\nu}{\partial x_i}(\gamma(t))\,dx_{i,k}.
\]
Each $r_{\nu,a}$ is a finite sum. If $e_1,\ldots,e_m$
is the standard basis of $B^m$, the $A$-linear isomorphism
\[
 \bigoplus_{i,k\ge-b}A\,dx_{i,k}
 \xrightarrow{\sim} B^m/t^bA\llbracket t \rrbracket^m,
 \qquad
 dx_{i,k}\longmapsto[t^{-k-1}e_i],
\]
sends $r_{\nu,a}$ to $[t^{-a-1}v_\nu(t)]$.

Put $J_\gamma=\sum_{\nu=1}^l B\,v_\nu(t)$.
Choose $c\ge0$ such that every $v_\nu(t)$ belongs to
$t^{-c}A\llbracket t \rrbracket^m$. For any $h\in B$, choose a Laurent
polynomial $h_0$ with $h-h_0\in t^{b+c}A\llbracket t \rrbracket$.
Then $(h-h_0)v_\nu\in t^bA\llbracket t \rrbracket^m$. Hence the image
of $J_\gamma$ in $B^m/t^bA\llbracket t \rrbracket^m$ is exactly the
$A$-span of the classes $[t^{-a-1}v_\nu]$.
Combining this calculation with the cotangent presentation
of $X$, in which $e_i$ maps to $dg_i\otimes dt$, gives
the following identifications of $A$-modules:
\begin{equation}\label{eq:elementary-coefficient-quotient}
 (\Omega^1_{\Spec C_b(X)/\C})|_T
 \cong B^m/(J_\gamma+t^bA\llbracket t \rrbracket^m),
 \qquad \Cot_t(\gamma_T)\cong B^m/J_\gamma.
\end{equation}
At a field-valued loop $\gamma$ in $\Spec C_b(X)$,
\eqref{eq:elementary-coefficient-quotient} identifies the cotangent
fiber of this stage with a quotient of $\Cot_t(\gamma)$.
With $t$ fixed, $\Cot_t(\gamma)$ identifies with the fiber at $\gamma$
of the Tate cotangent sheaf of $\sL X$ described by Drinfeld
\cite[p.~289 and p.~293]{Drinfeld}.

Reasonableness gives
\[
 S=\Spec\bigl(C_b(X)/(h_1,\ldots,h_q)\bigr).
\]
Put $R_Y=\Gamma(Y,\mathcal O_Y)$ and
$R_U=\Gamma(U,\mathcal O_U)=R_Y[u_j:j\in J]$.
Choose lifts $\widetilde h_\rho\in P_b$ of the $h_\rho$ and set
\[
 w_\rho(t)=\sum_{i,k}
 \left.\frac{\partial\widetilde h_\rho}{\partial x_{i,k}}\right|_U
 t^{-k-1}e_i\in R_U[t,t^{-1}]^m.
\]
These sums are finite. Since $U\to S$ is \'etale,
\[
 \Omega^1_{R_U/\C}
 \cong R_U\otimes_{C_b(X)/(h_1,\ldots,h_q)}
                  \Omega^1_{(C_b(X)/(h_1,\ldots,h_q))/\C}.
\]
Choose algebra generators $y_1,\ldots,y_d$ of $R_Y$ over $\C$.
The conormal presentation and this isomorphism allow us to write
\[
 dy_\lambda=\sum_{i,k}c_{\lambda,i,k}\,dx_{i,k}
       \quad\text{in }\Omega^1_{R_U/\C},
 \qquad c_{\lambda,i,k}\in R_U,
\]
with only finitely many nonzero terms for each $\lambda$.
Here the coefficient differentials are pulled back to $U$.
Set $z_\lambda(t)=\sum_{i,k}c_{\lambda,i,k}t^{-k-1}e_i$.
All these choices are made over $U$. Write $w_{\rho,T}$ and
$z_{\lambda,T}$ for their coefficientwise images in $B^m$.
The conormal sequence for $S\subset\Spec C_b(X)$ and the
relative cotangent sequence for $R_Y\to R_U$ now give
\[
 (\Omega^1_{U/Y})|_T
 \cong
 \frac{B^m}{J_\gamma+t^bA\llbracket t \rrbracket^m+
              \sum_{\rho=1}^q A w_{\rho,T}+
              \sum_{\lambda=1}^d A z_{\lambda,T}}.
\]
These are cokernel presentations, so tensoring them with $A$
requires no flatness assumption. Under the coefficient-residue
identification above, the quotient map from $\Cot_t(\gamma_T)$ is exactly \eqref{eq:elementary-residue-map}. In particular it
is surjective. Writing $p:B^m\twoheadrightarrow \Cot_t(\gamma_T)$, the kernel
of $q_{\pi,T}$ is
\[
 L_T=p(t^bA\llbracket t \rrbracket^m)
       +\sum_{\rho=1}^q A\,p(w_{\rho,T})
       +\sum_{\lambda=1}^d A\,p(z_{\lambda,T}).
\]

Suppose $\Cot_t(\gamma_T)$ is free and use the chosen basis to identify it
with $B^n$. Choose a $B$-linear section $\sigma:B^n\to B^m$
of $p$. Choose $c_0\ge0$ so that the matrices of $p$ and
$\sigma$, and the finitely many vectors $p(w_{\rho,T})$ and
$p(z_{\lambda,T})$, have all their entries in $t^{-c_0}A\llbracket t \rrbracket$.
Then
\[
 t^{b+c_0}A\llbracket t \rrbracket^n
 \subset p(t^bA\llbracket t \rrbracket^m)
 \subset L_T\subset t^{-c_0}A\llbracket t \rrbracket^n.
\]
The first inclusion follows by applying $\sigma$ and using
$p\sigma=1$; the last follows from the displayed formula for
$L_T$. Thus $N=b+c_0$ satisfies
\eqref{eq:elementary-kernel-bounds}.

Finally, a homomorphism $A\to A'$ induces the coefficientwise
map $A(\!( t )\!)\to A'(\!( t )\!)$ and the identification
\[
 \Cot_t(\gamma_{T'})\cong \Cot_t(\gamma_T)\otimes_{A(\!( t )\!)}A'(\!( t )\!).
\]
The fixed vectors $w_\rho,z_\lambda$ and the formula for
$q_{\pi,T}$ specialize under this map. Consequently
\[
 q_{\pi,T'}(\xi|_{T'})=q_{\pi,T}(\xi)\otimes1
 \quad(\xi\in \Cot_t(\gamma_T)),
\]
which proves the asserted map $L_T\to L_{T'}$.
\end{proof}

The same construction of $q_{\pi,T}$ and the same proof apply
to a chosen nice chart with finite-factor projection $\pi$ of a
finitely presented locally closed subscheme of a reasonable
closed subscheme of $\mathscr LX$.
Indeed, locally factor its immersion as a closed immersion followed
by an open immersion. The additional finitely many equations add
conormal relations to the calculation above, while open restriction
localizes its coefficients. The formulas for $q_{\pi,T}$ agree
on overlaps and remain compatible with pullback.

\begin{definition}[Pointwise coefficient index]
\label{def:pointwise-coefficient-index}
Keep the notation and hypotheses of
Lemma~\ref{lem:elementary-coefficient-differentials}, and fix a
nowhere-vanishing volume form $\Omega$ on $X$. For a
field-valued point $x:\Spec k\to U$, let $\gamma_x$ be its loop,
write $q_{\pi,x}$ for the corresponding coefficient-residue map, and put
$L_x:=\ker q_{\pi,x}$.
The form $\Omega$ defines
$\varphi_\Omega:\det\Cot_t(\gamma_x)\xrightarrow{\sim}k(\!( t )\!)$ by
$\varphi_\Omega(\gamma_x^*\Omega\otimes dt^{\otimes n})=1$.
Choose a basis $e_1,\ldots,e_n$ with
$\varphi_\Omega(e_1\wedge\cdots\wedge e_n)=1$ and an integer
$N\ge0$ satisfying \eqref{eq:elementary-kernel-bounds}.
The \emph{pointwise coefficient index} of $\pi$ at $x$,
normalized by $\Omega$, is
\begin{equation}\label{eq:elementary-field-index}
 c_\pi(x):=
 \dim_k\bigl(L_x/t^Nk\llbracket t \rrbracket^n\bigr)-nN.
\end{equation}
\end{definition}
This is the value on $L_x$ of the dimension theory normalized
by $\varphi_\Omega$ in \cite[Sections~3.5 and~3.9]{Drinfeld}.
We verify its well-definedness directly.
Increasing $N$ adds the same dimension to both terms. This
number is also independent of the volume-normalized basis:
the transition matrix has determinant $1$, and the elementary
divisor theorem over $k\llbracket t \rrbracket$ gives relative dimension zero
between its two reference power-series modules. Explicitly,
for an arbitrary transition matrix $G\in\mathrm{GL}_n(k(\!( t )\!))$
that relative dimension is $-\ord_t\det G$, as is seen by
reducing $G$ to a diagonal matrix with entries powers of $t$.
Comparison inside a common smaller power-series module then
proves the asserted independence for any kernel $L_x$.

The pointwise coefficient index $c_\pi(x)$ is also unchanged
by extension of the field of definition. Let $k\subset K$ be
a field extension, let $x_k=x:\Spec k\to U$, and let
$x_K:\Spec K\to U$ be its pullback. For $R=k,K$, let $\gamma_R$
be the loop of $x_R$ and put $L_R=L_{x_R}$. Under
\[
 \Cot_t(\gamma_K)\cong
 \Cot_t(\gamma_k)\otimes_{k(\!( t )\!)}K(\!( t )\!),
\]
the chosen volume-normalized basis extends to a volume-normalized
basis. Choose a common $N\ge0$
such that
\[
 H_R:=t^NR\llbracket t \rrbracket^n\subset L_R\subset t^{-N}R\llbracket t \rrbracket^n
 \qquad(R=k,K).
\]
Such an $N$ exists by
\eqref{eq:elementary-kernel-bounds}, since increasing $N$
preserves both inclusions. The monomial bases give
\[
 (\Cot_t(\gamma_k)/H_k)\otimes_kK\cong \Cot_t(\gamma_K)/H_K,
 \qquad
 \left(\bigoplus_{j\in J}k\,du_j\right)\otimes_kK
 \cong\bigoplus_{j\in J}K\,du_j.
\]
Under these identifications, the quotient map induced by
$q_{\pi,x_k}$ becomes the quotient map induced by
$q_{\pi,x_K}$ after tensoring with $K$, by the coefficientwise
compatibility in
Lemma~\ref{lem:elementary-coefficient-differentials}.
Tensoring the kernel sequence of the induced map on $\Cot_t(\gamma_k)/H_k$
with the flat $k$-module $K$ therefore gives
\[
 (L_k/H_k)\otimes_kK\cong L_K/H_K.
\]
The bounds make these kernels finite-dimensional, whence
$c_\pi(x_K)=c_\pi(x_k)$. Only the quotients by $H_R$ are
identified by ordinary scalar extension: their elements have
finite Laurent-polynomial representatives with exponents
less than $N$. No identification of
$k\llbracket t \rrbracket\otimes_kK$ with $K\llbracket t \rrbracket$ is used. Thus $c_\pi$
defines a function on the underlying scheme points of $U$,
whose value can be computed over the residue field or any
extension of that field.

\begin{lemma}
\label{lem:coefficient-index-specialization}
In the setup of Definition~\ref{def:pointwise-coefficient-index},
let $\Spec k\llbracket s \rrbracket\to U$ be a family, where
$k$ is an algebraically closed extension of $\C$.
The pointwise coefficient index $c_\pi$ has equal values at
the generic and special points of this family.
\end{lemma}

\begin{proof}
Put $A=k\llbracket s \rrbracket$. The ring
\[
 B=A(\!( t )\!)=k\llbracket s,t \rrbracket[1/t]
\]
is a one-dimensional Noetherian UFD, hence a PID: the regular
local ring $k\llbracket s,t \rrbracket$ is factorial of dimension two, and its
unique height-two prime contains $t$. Let $\gamma_A$ be the loop
over $A$ given by the family. Smoothness of $X$ makes
$\Cot_t(\gamma_A)$ finite projective over $B$, hence free. The volume form
trivializes its determinant, so choose a volume-normalized basis
$e_1,\ldots,e_n$.

For $R=A$, $k$, or $k(\!( s )\!)$, write $\gamma_R$ for the
loop induced by $\gamma_A$ and $L_R=\ker q_{\pi,\Spec R}$.
For $F=k$ or $k(\!( s )\!)$,
\[
 \Cot_t(\gamma_F)\cong \Cot_t(\gamma_A)\otimes_{A(\!( t )\!)}F(\!( t )\!),
\]
so the chosen basis induces a volume-normalized basis of $\Cot_t(\gamma_F)$.
Apply Lemma~\ref{lem:elementary-coefficient-differentials}
separately over these three rings and choose one $N$ large
enough that
\[
 H_R:=t^NR\llbracket t \rrbracket^n\subset L_R\subset t^{-N}R\llbracket t \rrbracket^n
              \qquad(R=A,k,k(\!( s )\!)).
\]
The coefficient-residue map $q_{\pi,\Spec A}$ induces a
split exact sequence
\[
 0\longrightarrow L_A/H_A\longrightarrow \Cot_t(\gamma_A)/H_A
   \xrightarrow{\bar q_{\pi,A}}
       \bigoplus_{j\in J}A\,du_j\longrightarrow0,
\]
because its last module is free. Moreover,
\[
 \Cot_t(\gamma_A)/H_A=\bigoplus_{i=1}^n\bigoplus_{r<N}A[t^re_i].
\]
For either field $F$, the same monomial bases identify
\[
 (\Cot_t(\gamma_A)/H_A)\otimes_A F\cong \Cot_t(\gamma_F)/H_F,
 \qquad
 \left(\bigoplus_{j\in J}A\,du_j\right)\otimes_A F
       \cong\bigoplus_{j\in J}F\,du_j.
\]
Under these identifications,
$\bar q_{\pi,A}\otimes_A F=\bar q_{\pi,F}$: on each basis
vector $[t^re_i]$ this is the coefficientwise compatibility
of \eqref{eq:elementary-residue-map}, whose evaluation is a
finite sum.
Tensoring the split exact sequence therefore gives
\[
                 (L_A/H_A)\otimes_A F\cong L_F/H_F.
\]
As in the field-extension argument above, ordinary scalar
extension is used only after quotienting by $H_R$.

The module $L_A/H_A$ is a direct summand of the free module
$\Cot_t(\gamma_A)/H_A$, hence projective. It is also a submodule of
$t^{-N}A\llbracket t \rrbracket^n/H_A\cong A^{2nN}$, so it is finitely generated
because $A$ is Noetherian. Since $A$ is local, it is finite
free, of some rank $d$. For the generic and special points
$x_\eta$ and $x_0$, the preceding fiber identifications and
\eqref{eq:elementary-field-index} give
$
                    c_\pi(x_\eta)=d-nN=c_\pi(x_0).
$
\end{proof}

\subsection{Proof of the volume-index propositions}
\label{subsec:volume-index-proofs}
We prove Propositions~\ref{prop:finite-factor-normalization},
\ref{prop:contact-volume-indices}, and
\ref{prop:general-factor-contact-index}.
In Steps~1--3 we work with the pointwise coefficient indices
$c_\pi$ of Definition~\ref{def:pointwise-coefficient-index}.
Step~4 constructs the locally constant functions $\delta_Y$
with $c_\pi=\pi^*\delta_Y$ and completes the three propositions.

\emph{Step 1: comparisons of the pointwise coefficient indices $c_\pi$.}
We establish the pointwise counterparts of
Proposition~\ref{prop:finite-factor-normalization}(i)--(iii)
and verify that the construction of $c_\pi$ applies to the
contact loci.
For a smooth comparison $\pi=f\pi'$ of relative dimension $e$,
where $f:Y'\to Y$, the finite-factor presentation for $\pi'$ gives
an exact sequence
\[
 0\longrightarrow (\pi')^*\Omega^1_{Y'/Y}
 \longrightarrow\Omega^1_{U/Y}
 \longrightarrow\Omega^1_{U/Y'}\longrightarrow0.
\]
At a field-valued point $x\in U$, the coefficient-residue maps
satisfy
\[
 q_{\pi',x}=
 \bigl((\Omega^1_{U/Y})|_x\longrightarrow
             (\Omega^1_{U/Y'})|_x\bigr)\circ q_{\pi,x}.
\]
Put $L_\pi:=\ker q_{\pi,x}$ and
$L_{\pi'}:=\ker q_{\pi',x}$.
Surjectivity of $q_{\pi,x}$ and $q_{\pi',x}$ gives
\[
 L_\pi\subset L_{\pi'},\qquad
 L_{\pi'}/L_\pi\cong
 ((\pi')^*\Omega^1_{Y'/Y})|_x.
\]
The quotient has dimension $e$, so
\eqref{eq:elementary-field-index} gives
$c_{\pi'}(x)=c_\pi(x)+e$. Relative dimensions here and below
are the ranks of the corresponding relative cotangent sheaves
at the specified points.
In a finite-factor presentation
$(Y_a\hookrightarrow Y_b)\times\A^J$, let $\pi_a,\pi_b$
be the respective product projections. For a field-valued
point $x$ of $Y_a\times\A^J$, write $x_b$ for its image in
$Y_b\times\A^J$. The relative cotangent quotients have the
same omitted-coordinate basis $du_j$, and the coefficient-residue
maps satisfy $q_{\pi_a,x}=q_{\pi_b,x_b}$ under this
identification. Hence their kernels agree and
$c_{\pi_a}(x)=c_{\pi_b}(x_b)$. Open restrictions preserve
this calculation.

For a finitely presented \'etale chart change $v:V\to U$,
its finite presentation uses only finitely many coordinates of
$U$. Retaining these coordinates descends it to a morphism
$V_0\to Y\times\A^e$; this morphism is \'etale because its
base change along the faithfully flat polynomial extension
adjoining the remaining coordinates is $v$. Thus we have
$V\cong V_0\times\A^{J'}$ over
$U\cong (Y\times\A^e)\times\A^{J'}$.
The coefficient-residue maps for the projections to $V_0$
and $Y\times\A^e$ have the same omitted-coordinate basis
and commute with \'etale pullback. Writing $g:V_0\to Y$ for
the composite projection, we obtain $c_{V_0}=v^*c_\pi+e$,
where $c_{V_0}$ denotes the pointwise coefficient index of
the projection $V\to V_0$.
Now let $(\rho:V\to Z,\ \pi:U\to Y,\ f:Z\to Y)$ be
another finite-factor presentation of $v$.
Choose a common finite factor $V\to W$ with smooth comparisons
$a:W\to Z$ and $b:W\to V_0$. The product projection
$V\to W$ has a section, so $fa=gb$. At corresponding points put
$d_a=\operatorname{reldim}(a)$ and
$d_b=\operatorname{reldim}(b)$. The comparison just proved and
additivity of smooth relative dimension give
\[
 c_W=c_\rho+d_a=c_{V_0}+d_b,
 \qquad d_a+\operatorname{reldim}(f)=d_b+e.
\]
Hence $c_\rho=v^*c_\pi+\operatorname{reldim}(f)$.
These statements use only one chart map at a time and do not
presuppose Proposition~\ref{prop:finite-factor-normalization}.

To apply Lemma~\ref{lem:elementary-coefficient-differentials}
to the contact schemes, fix
a nice chart $h:U=Y\times\A^J\to F_P\mathscr LX$, with
$\pi=\operatorname{pr}_Y$. Retain the finitely many affine
coordinates occurring in the equations of all pulled-back pole
substages. Thus, for some $e\ge0$,
\[
 U=T\times\A^{J'},\qquad T=Y\times\A^e,\qquad
 h^{-1}(F_p\mathscr LX)=T_p\times\A^{J'}
       \quad(0\le p\le P),
\]
where $T_p\subset T$ are closed subschemes, $T_P=T$, and
$T_{-1}=\varnothing$. Put
$L_p=(T_p\setminus T_{p-1})_{\mathrm{red}}$.
For every ring $R$,
\[
 \sqrt{(0)}_{\,R[u_j:j\in J']}
      =\sqrt{(0)}_{\,R}[u_j:j\in J'].
\]
Since each $T_p$ is of finite type, its nilradical is finitely
generated. Hence the reduced pulled-back layer is
$L_p\times\A^{J'}$, and its immersion into $U$ is finitely
presented. By Proposition~\ref{prop:loop-sectors}, it is the
disjoint union of the clopen subschemes
\[
 U_{\vct}:=U\times_{F_P\mathscr LX}\mathscr P_{\vct},
       \qquad w(\vct)=p.
\]
Indeed, these base changes are reduced because they are
\'etale over the reduced contact schemes. The layer is
quasi-compact, so each clopen piece is quasi-compact and its
immersion into $U$ is finitely presented. This property can
be checked on the \'etale charts, so the contact immersions
$\mathscr P_{\vct}\hookrightarrow F_P\mathscr LX$ are finitely
presented as well. The extension of
Lemma~\ref{lem:elementary-coefficient-differentials} to finitely
presented locally closed subschemes therefore applies.

\emph{Step 2: the contact-jet calculation.}
We prove the finite-factor and smoothness assertions of
Proposition~\ref{prop:contact-volume-indices} and establish
its two index formulas for the pointwise coefficient indices
$c_\pi$.
Let $x:\Spec k\to\mathscr P_{\vct,W}$ be a field-valued point.
Near the center of its contact arc, choose \'etale SNC coordinates
$z_1,\ldots,z_n$, avoiding boundary components outside
$I(\vct)$. The contact equations fix $|\vct|$ coefficients
of an $\ell$-jet and invert its specified leading coefficients.
Thus the jet factor is smooth of dimension
$n(\ell+1)-|\vct|$, and \'etale lifting makes all higher
coefficients free. Write
\[
 \Omega=b\prod_{i\in I(\vct)}z_i^{-a_i}
                   dz_1\wedge\cdots\wedge dz_n,
 \qquad b\in\mathcal O(W)^\times.
\]
The infinitesimal deformations fixing this jet are
$\bigoplus_i t^{\ell+1}k\llbracket t \rrbracket\partial_{z_i}$.
For the contact-jet projection $\pi_{\vct,\ell}$, the
coefficient-residue kernel is therefore
\[
 L_x=\ker q_{\pi_{\vct,\ell},x}
     =\bigoplus_i t^{-\ell-1}k\llbracket t \rrbracket\,dz_i\otimes dt.
\]
Indeed, a Laurent series has zero Laurent residue against every
$t^{\ell+1}c(t)$ precisely when its exponents are at least
$-\ell-1$. The determinant of this displayed basis, measured
by $\varphi_\Omega$, has order
$-n(\ell+1)+\sum_i a_iv_i$. The elementary-divisor calculation
following \eqref{eq:elementary-field-index} gives
\[
 c_{\pi_{\vct,\ell}}(x)=n(\ell+1)-\sum_i a_iv_i,
\]
the pointwise formula asserted in
\eqref{eq:volume-dimension-index}.

Every finitely presented \'etale chart over a contact scheme
descends locally to a sufficiently high contact-jet factor:
its finitely many equations and inverses use finitely many
arc coefficients. An arbitrary finite factor is compared with
such a factor through common smooth refinements by
Lemma~\ref{lem:finite-factor-comparison}. Its smoothness follows
from the infinitesimal lifting property of the contact
scheme: the zero section of the finite-factor presentation transfers
this property to its finite factor. Under a smooth refinement,
both dimension and the pointwise coefficient index $c_\pi$
increase by its relative dimension. This gives
\eqref{eq:contact-factor-volume-index} pointwise and completes
the structural assertions of
Proposition~\ref{prop:contact-volume-indices}.

\emph{Step 3: general factors and constructibility.}
We construct the contact-factor presentation of
Proposition~\ref{prop:general-factor-contact-index}, establish
its formula for $c_\pi$, and express $c_\pi$ on a finite
refinement as the pullback of a constructible function.
Keep the enlarged finite-factor presentation $U=T\times\A^{J'}$ above.
On an affine open $\Spec R\subset L_p$, a clopen contact
piece is defined by an idempotent of $R[u_j:j\in J']$.
Every such idempotent belongs to $R$: its nonconstant
coefficients vanish modulo every prime of $R$, hence vanish
because $R$ is reduced. Thus the contact piece descends uniquely
to $\Spec R$.
Uniqueness makes these descended pieces agree on overlaps.
We obtain reduced locally closed subschemes $Z_{\vct}\subset T$
such that
\[
 L_p=\coprod_{w(\vct)=p}Z_{\vct},\qquad
 U_{\vct}=Z_{\vct}\times\A^{J'}.
\]
Each $Z_{\vct}$ is smooth by Step~2, applied on its affine
open subsets. For $u\in U_{\vct}$ with image $z\in T$, choose
an affine open neighborhood $Y'\subset T$ of $z$ and put
$U'=Y'\times\A^{J'}$ and $Z=Z_{\vct}\cap Y'$.
This gives the finite-factor presentation in
Proposition~\ref{prop:general-factor-contact-index}.

By Step~1, restriction to this contact piece and removal of
the $e$ retained coordinates give
\[
 c_\pi(u)=\dim_z Z+\tfrac12\mu_\Omega(\vct)-e.
\]
This proves \eqref{eq:contact-volume-general-factor}
pointwise, independently of the finite-factor presentation.
The finitely many smooth schemes $Z_{\vct}$ have finitely
many connected components. The rule
\[
 c_T(z)=\dim_z Z_{\vct}+\tfrac12\mu_\Omega(\vct)-e
       \qquad(z\in Z_{\vct})
\]
therefore defines a finite-valued constructible function
$c_T:T\to\Z$. Writing $\pi_T:U\to T$ for the refined
projection, we have
\[
 c_\pi=\pi_T^*c_T,\qquad
 c_{\pi_T}=\pi_T^*c_T+e.
\]

\emph{Step 4: local constancy and descent.}
We prove local constancy and descend $c_\pi$ to $\delta_Y$ on
the original finite factor, completing
Proposition~\ref{prop:finite-factor-normalization} and the
volume-index formulas of
Propositions~\ref{prop:contact-volume-indices} and
\ref{prop:general-factor-contact-index}.
Let $s_0:T\to U=T\times\A^{J'}$ be the zero section
of the remaining affine tail. Every morphism
$\operatorname{Spec}\C\llbracket s \rrbracket\to T$ can be composed with $s_0$;
Lemma~\ref{lem:coefficient-index-specialization} therefore gives
equal generic and special values of $c_T$ along this morphism.
To prove local constancy, suppose a constructible level set
$C_m=c_T^{-1}(m)$ is not closed. Since $T$ is of finite type
over $\C$, there is a closed point
$x_0\in\overline{C_m}\setminus C_m$ and an integral curve
through $x_0$ whose generic point lies in $C_m$. One obtains
such a curve by taking successive hyperplane sections through
$x_0$ in an irreducible component of $\overline{C_m}$ on
which $C_m$ contains a dense open subset, and choosing a curve
component through $x_0$ not contained in its complement.
Normalize the curve and complete at a point above $x_0$. Its completed local ring is $\C\llbracket s \rrbracket$. The resulting morphism
$\operatorname{Spec}\C\llbracket s \rrbracket\to T$ has generic value $m$:
completion extends the function field of the curve, and the
pointwise coefficient index $c_\pi$ is invariant under this
field extension.
Its special value is different from $m$, a contradiction.
Thus every $C_m$ is closed. There are finitely many nonempty
level sets, so each is also open.

Write $p:T=Y\times\A^e\to Y$ for the projection and
$i_0:Y\to T$ for its zero section. Each fiber of $p$ is an
affine space over a residue field and hence connected, so $c_T$
is constant on that fiber. Consequently
$\delta_Y:=i_0^*c_T$ is locally constant and
$c_T=p^*\delta_Y$. Surjectivity of $p$ proves uniqueness.
We have $c_\pi=\pi^*\delta_Y$. Step~1 now proves
Proposition~\ref{prop:finite-factor-normalization}(i)--(iii)
and their compatibilities; Steps~2 and~3 give
Propositions~\ref{prop:contact-volume-indices} and
\ref{prop:general-factor-contact-index}, respectively.


\emph{Reduction of finite models.}
For a nice chart $U=Y\times\A^J\to F_P\mathscr LX$, one has
\[
 U_{\mathrm{red}}=Y_{\mathrm{red}}\times\A^J.
\]
Give $Y_{\mathrm{red}}$ the inherited index
$\delta_{Y_{\mathrm{red}}}:=\delta_Y|_{Y_{\mathrm{red}}}$.
Write $\pi_{\mathrm{red}}:U_{\mathrm{red}}\to Y_{\mathrm{red}}$
for the induced projection. At every field-valued point $x$,
the coefficient-residue maps $q_{\pi_{\mathrm{red}},x}$ and
$q_{\pi,x}$ agree under the identification of their
omitted-coordinate cotangent modules. Thus Step~1 gives
$c_{\pi_{\mathrm{red}}}=\pi_{\mathrm{red}}^*\delta_{Y_{\mathrm{red}}}$.
Identifying $Y_{\mathrm{red}}(\C)$ with $Y(\C)$ gives
\[
 \omega_{Y_{\mathrm{red}}}[-2\delta_{Y_{\mathrm{red}}}]
       \simeq\omega_Y[-2\delta_Y],
\]
also for singular $Y$. Every morphism in the chosen finite-factor
presentation diagram induces the same analytic map after reduction,
so these identifications preserve
proper traces, localization sequences, and their compositions.
For a smooth comparison $f:Y'\to Y$, one has
$Y'_{\mathrm{red}}=Y'\times_Y Y_{\mathrm{red}}$; hence its
reduction is smooth of the same relative dimension, and the
normalized smooth-pullback comparisons agree as well.
Base-change comparisons are those of the original scheme-theoretic
finite-factor presentation diagrams, transported through these identifications of
analytic spaces. $\hfill\Box$

\begin{example}[Affine line]\label{ex:affine-line-volume-index}
Take $X=\A^1$ with coordinate $x$, $\Omega=dx$, and the standard
compactification $(\Pj^1,\{\infty\})$ polarized by $\{\infty\}$.
The pole stage $F_P\sL X$ consists of
$x(t)=\sum_{j\ge-P}x_jt^j$. Retaining $x_{-P},\ldots,x_N$,
with $N\ge0$, gives the finite factor $Y_N=\A^{P+N+1}$.
At a loop of contact $q\le P$, its contact locus is already
presented by a locally closed subscheme $Z\subset Y_N$:
the coefficients below order $-q$ vanish, and $x_{-q}\ne0$
if $q>0$. For $q=0$ all negative coefficients vanish and
there is no nonvanishing condition on $x_0$.
Thus $\dim Z=q+N+1$ and no additional coordinates are needed,
so $e=0$. Since $dx$ has pole order two at infinity,
$\mu_\Omega(q)/2=-q$. Proposition~\ref{prop:general-factor-contact-index}
therefore gives, at every point,
\[
                  \delta_{Y_N}=(q+N+1)-q=N+1.
\]
Since $\dim Y_N=P+N+1$, this yields
\[
 \omega_{Y_N}[-2\delta_{Y_N}]\simeq\bQ_{Y_N}[2P],
 \qquad B_{\Omega,P}\simeq\bQ[2P].
\]
\end{example}

\section{The loop complex and loop cohomology}\label{sec:loop-complex}

This section constructs the filtered loop complex and loop
cohomology underlying Theorem~\ref{thm:summary-fusion}. We assemble
the pole-stage complexes by a homotopy colimit and prove independence
of various choices for a fixed polarized compactification.
We compute
the associated graded objects, providing the first-page calculation for
Theorem~\ref{thm:summary}. We also explain changes of volume grading
and, under explicit compatibility hypotheses, give a computational tool
in presence of global finite-factor presentations. 

\subsection{The loop complex}
\label{subsec:loop-definition}
\label{subsec:well-definedness}

\begin{definition}[Loop complex]\label{def:loop-complex}
For every polarized SNC pair $(M,D,A)$, the pole-stage diagrams of
Definition~\ref{def:pole-stage-complexes} define the loop complexes
of $\bQ$-vector spaces
\begin{equation}\label{eq:loop-definition}
 \begin{aligned}
 B_{\Omega}
   &:=\hocolim_P B_{\Omega,P},\\
 B
   &:=\hocolim_P B_P,
 \end{aligned}
\end{equation}
in the integer and parity gradings, respectively. Their pole
filtrations are the filtered objects
\[
 P\longmapsto B_{\Omega,P},\qquad
 P\longmapsto B_P
\]
in the sense of Definition~\ref{def:filtered-derived-object}: these are functors from
$\Z_{\ge0}$ to the respective derived $\infty$-categories,
with the transition maps of
Definition~\ref{def:pole-stage-complexes} and the composition
and refinement coherences described explicitly in
Remark~\ref{rem:coherent-pole-descent}. Equation~\eqref{eq:loop-definition}
makes these filtrations exhaustive. The objects $B_P$ and
$B_{\Omega,P}$ come with their canonical maps to $B$ and
$B_\Omega$. Theorem~\ref{thm:canonical-complex} realizes these
filtered objects by cochain complexes $\mathbf B$ and
$\mathbf B_\Omega$ with filtrations by subcomplexes
$F_P\mathbf B$ and $F_P\mathbf B_\Omega$. Homotopy
colimits are recalled in
Appendix~\ref{app:homotopy-limits-colimits}. \end{definition}

\begin{theorem}
\label{thm:canonical-complex}
The loop objects $B_\Omega$ and $B$
admit representatives by, respectively, a $\Z$-graded cochain
complex $\mathbf B_\Omega$ and a two-periodic cochain complex
$\mathbf B$ of $\bQ$-vector spaces, with exhaustive nonnegative
filtrations by subcomplexes. Choose together with them equivalences of the
entire stage diagrams
\[
 (F_P\mathbf B_\Omega)_{P\ge0}\simeq(B_{\Omega,P})_{P\ge0},
 \qquad
 (F_P\mathbf B)_{P\ge0}\simeq(B_P)_{P\ge0},
\]
in $\operatorname{Fun}(N(\Z_{\ge0}),\mathcal D(\bQ))$ and
$\operatorname{Fun}(N(\Z_{\ge0}),\mathcal D_{\mathrm{per}}(\bQ))$,
respectively, in the sense of
Definition~\ref{def:filtered-derived-object}. They induce
$\mathbf B_\Omega\simeq B_\Omega$
and $\mathbf B\simeq B$ in the respective derived categories.
In particular,
\[
 F_0\mathbf B_\Omega\simeq C^*(X;\bQ),\qquad
 F_0\mathbf B\simeq\operatorname{Per}_2C^*(X;\bQ).
\]
For the fixed polarized compactification, these representatives
are independent of the choices of nice charts and compatible
finite-factor presentations in
Proposition~\ref{prop:normalized-pole-sheaves}, up to a zigzag
of filtration-preserving cochain maps that are quasi-isomorphisms
on every stage.

Put $F_{-1}\mathbf B=F_{-1}\mathbf B_\Omega=0$. The literal
quotient complexes
\[
 \gr_P\mathbf B:=F_P\mathbf B/F_{P-1}\mathbf B,
 \qquad
 \gr_P\mathbf B_\Omega:=F_P\mathbf B_\Omega/F_{P-1}\mathbf B_\Omega
\]
are the weight-$P$ terms of the associated-graded complexes
\[
 \gr_F\mathbf B:=\bigoplus_{P\ge0}\gr_P\mathbf B,
 \qquad
 \gr_F\mathbf B_\Omega:=\bigoplus_{P\ge0}\gr_P\mathbf B_\Omega.
\]
The specified comparisons of filtered diagrams induce
\[
 \gr_P\mathbf B\simeq\operatorname{cofib}(B_{P-1}\to B_P),
 \qquad
 \gr_P\mathbf B_\Omega
 \simeq\operatorname{cofib}(B_{\Omega,P-1}\to B_{\Omega,P}).
\]
There is an
equivalence in $\mathcal D(\bQ)$
\begin{equation}\label{eq:loop-layer}
 \gr_P\mathbf B_{\Omega}
 \simeq\bigoplus_{w(\vct)=P}
 C^{*-\mu_\Omega(\vct)}(S_{I(\vct)}^{\log};\F).
\end{equation}
The parity counterpart lies in $\mathcal D_{\mathrm{per}}(\bQ)$
and has summands $\operatorname{Per}_2C^*(S_{I(\vct)}^{\log};\F)$;
it does not require $\Omega$. When $\Omega$ exists,
$\operatorname{Per}_2B_\Omega\simeq B$
in $\mathcal D_{\mathrm{per}}(\bQ)$.
\end{theorem}

\begin{proof}
\emph{Filtered models and independence of finite-factor presentations.}
Let $\mathcal C=\mathcal D(\bQ)$ in the volume grading and
$\mathcal C=\mathcal D_{\mathrm{per}}(\bQ)$ in the parity grading.
Write $C_\bullet=(C_P)_{P\ge0}$ for the corresponding coherent
pole-stage diagram, so $C_P=B_{\Omega,P}$ in the volume grading
and $C_P=B_P$ in the parity grading. Put $C=B_\Omega$ or $C=B$,
respectively. Propositions~\ref{prop:diagram-model-comparison}
and~\ref{prop:simultaneous-diagram-representatives} give a
strict sequence of cochain complexes
\[
 \mathbf K_0\xrightarrow{u_0}\mathbf K_1\xrightarrow{u_1}\mathbf K_2
 \xrightarrow{u_2}\cdots
\]
and an equivalence $\mathbf K_\bullet\simeq C_\bullet$ in
$\operatorname{Fun}(N(\Z_{\ge0}),\mathcal C)$.
Here $u_P:\mathbf K_P\to\mathbf K_{P+1}$ are cochain maps;
they need not be inclusions. Comparisons are represented together
with their source and target diagrams, using
$\Z_{\ge0}\times[1]$. Model homotopy limits and colimits compute
those in $\mathcal C$ by
\cite[Propositions~1.3.4.23--1.3.4.25]{HA}.
Put
\[
 T=\operatorname{Tel}(\mathbf K_\bullet)
  :=\operatorname{Cone}\!\left(
       1-s:\bigoplus_{p\ge0}\mathbf K_p\longrightarrow\bigoplus_{p\ge0}\mathbf K_p
     \right),
\]
where $s$ sends the $p$th summand to the $(p+1)$st by $u_p$.
Its finite initial telescopes are the subcomplexes
\[
 T_P:=\operatorname{Cone}\!\left(
       1-s:\bigoplus_{p<P}\mathbf K_p\longrightarrow\bigoplus_{p\le P}\mathbf K_p
     \right),\qquad T=\bigcup_{P\ge0}T_P.
\]
Proposition~\ref{prop:sequential-telescope} gives natural
quasi-isomorphisms $\varepsilon_P:T_P\to\mathbf K_P$ satisfying
\[
 \varepsilon_{P+1}\circ(T_P\hookrightarrow T_{P+1})
       =u_P\circ\varepsilon_P.
\]
Hence $T_\bullet\simeq\mathbf K_\bullet\simeq C_\bullet$ as filtered
objects in $\operatorname{Fun}(N(\Z_{\ge0}),\mathcal C)$.
The compatibility of model homotopy colimits with colimits in
$\mathcal C$, recalled in
Proposition~\ref{prop:hocolim-properties}, then gives
\[
 T\simeq\hocolim_P\mathbf K_P\simeq\hocolim_P C_P=C.
\]
Set $\mathbf B_\Omega:=T$ and $F_P\mathbf B_\Omega:=T_P$
in the volume grading, and $\mathbf B:=T$ and
$F_P\mathbf B:=T_P$ in the parity grading. These are equalities
of cochain complexes. The preceding equivalences give the specified
comparisons of their entire filtrations with the original
pole-stage diagrams. The construction uses homotopy colimits after derived
global sections; no interchange of $R\Gamma$ with these
colimits is required.

Represent the coherent refinement equivalences of
Remark~\ref{rem:coherent-pole-descent} together with their
pole-stage diagrams.
The natural maps $\varepsilon_P$ identify the resulting
telescope filtrations stagewise. Objectwise quasi-isomorphisms
of strict sequences induce quasi-isomorphisms of every finite
initial telescope, so these comparisons give the asserted
zigzags of filtered quasi-isomorphisms.

For $P\ge1$, let $U_P=F_P\mathscr LX\setminus F_{P-1}\mathscr LX$
and write $j:U_P\hookrightarrow F_P\mathscr LX$. As a right
adjoint between stable categories, derived global sections
preserve finite limits and hence finite colimits. Applying them
to the localization cofiber sequence of
Proposition~\ref{prop:normalized-pole-sheaves} gives a cofiber
sequence in $\mathcal D(\bQ)$,
\[
 B_{\Omega,P-1}\longrightarrow B_{\Omega,P}
 \longrightarrow R\Gamma(U_P,j^{-1}\omega_{P,\Omega}).
\]
The specified comparison of filtered diagrams identifies this
cofiber with $\gr_P\mathbf B_\Omega$. This is an equivalence
in the derived category; the latter notation denotes the literal
quotient of the chosen subcomplexes.
The same argument gives the parity filtration.

\emph{The contact grading.}
Proposition~\ref{prop:loop-sectors} identifies $(U_P)_{\mathrm{red}}$
with the finite disjoint union of the contact schemes of weight
$P$. Reduction leaves their analytic spaces unchanged, and its
compatibility with the normalizing indices was verified in
Subsection~\ref{subsec:volume-index-proofs}. On a contact-jet
finite factor $Y=Y_{\vct,\ell}$ of
Proposition~\ref{prop:contact-volume-indices}, complex orientation
and the index formula give
\[
 \omega_Y\simeq\bQ_Y[2\dim Y],\qquad
 2(\delta_Y-\dim Y)=\mu_\Omega(\vct),\qquad
 \omega_Y[-2\delta_Y]\simeq\bQ_Y[-\mu_\Omega(\vct)].
\]
These identifications respect the actual comparison maps.
Indeed, for a smooth comparison $f:Y'\to Y$ of relative
dimension $e$, both dimension and the volume index increase by
$e$. Oriented smooth purity therefore identifies the normalized
comparison with the identity on
$f^{-1}\bQ_Y[-\mu_\Omega(\vct)]$. These are the normalized purity
equivalences of Proposition~\ref{prop:normalized-pole-sheaves};
their composition coherences identify them with the canonical
pullbacks of the constant complex.
Thus Lemma~\ref{lem:infinite-product-descent}(ii) descends these
identifications of the restricted normalized sheaf with the constant complex
$\bQ_{\mathscr P_{\vct}}[-\mu_\Omega(\vct)]$ on the entire
contact locus. Independently, the same argument with periodic
coefficients identifies its periodic normalized complex with
$\operatorname{Per}_2\bQ_{\mathscr P_{\vct}}$; this uses no
commutation of periodization with derived global sections.

Write $J=I(\vct)$ and put
\[
 T_J:=\underset{D_J^\circ}{\prod_{j\in J}}
       \bigl(N_{D_j/M}|_{D_J^\circ}\setminus0\bigr),
 \qquad T_\varnothing:=X.
\]
Here $N_{D_j/M}$ is the normal line bundle of $D_j$ in $M$.
The leading-normal map
$p_{\vct}:\mathscr P_{\vct}(\C)\to T_J(\C)$ of
Proposition~\ref{prop:loop-sectors} is locally a product
projection with fiber a finite or countable product of copies
of $\C$. Lemma~\ref{lem:infinite-product-descent}(i), applied
on these local trivializations, makes its adjunction unit an
equivalence for constant rational and periodic coefficients.
Hence the derived global sections on the contact locus are
those on $T_J(\C)$. The homotopy equivalence
$T_J(\C)\to S_J^\circ$ and
Lemma~\ref{lem:face-restriction} identify the rational result
with
\[
 C^{*-\mu_\Omega(\vct)}(S_J^\circ;\F)
 \simeq C^{*-\mu_\Omega(\vct)}(S_J^{\log};\F).
\]
The adjunction units and normalized purity equivalences are
natural under jet truncation, so these identifications commute
with its comparison maps.
The finite contact decomposition proves
\eqref{eq:loop-layer}. For $\vct=\zero$, the leading-normal
map is center evaluation to $X$, and the same argument gives
$B_{\Omega,0}\simeq C^*(X;\bQ)$. The specified comparison
$F_0\mathbf B_\Omega\simeq B_{\Omega,0}$ gives the asserted
arc-stage equivalence for the cochain model.

\emph{Comparison with periodization.}
Recall $\Lambda=\bQ[u,u^{-1}]$, $|u|=2$ and $du=0$, from
Subsection~\ref{subsec:coherent-descent}. On every contact
piece, the natural coefficient-change map on derived global
sections is an equivalence. To see this without exchanging
periodization with an infinite totalization, use the preceding
product comparison to compute on $T_J(\C)$. This space is
homotopy equivalent to the compact manifold with corners
$S_J^{\log}$, which admits a finite triangulation. Let $C_J$
be the bounded complex of finite-dimensional rational cellular
chains of such a triangulation, with the cellular $k$-chains
placed in cohomological degree $-k$. The coefficient-change map is
then represented by the isomorphism
\[
 \operatorname{Hom}_{\bQ}(C_J,\bQ)\otimes_{\bQ}\Lambda
 \xrightarrow{\sim}
 \operatorname{Hom}_{\bQ}(C_J,\Lambda).
\]
The map sends $\phi\otimes\lambda$ to
$c\mapsto\phi(c)\lambda$; all degrees of $\Lambda$ are even.
The finite number and finite dimensions of the terms of $C_J$
justify the tensor--Hom interchange. The cellular model is
used only to verify that the natural coefficient-change map
is an equivalence. This proves the stated
parity description of every contact summand and, for $J=\varnothing$,
gives $B_0\simeq\operatorname{Per}_2C^*(X;\bQ)$.
The comparison $F_0\mathbf B\simeq B_0$ gives the corresponding
statement for the chosen cochain model.

When $\Omega$ exists, put $S_P=F_P\mathscr LX(\C)$ for $P\ge0$ and let
$a_P:S_P\to\mathrm{pt}$ be its structure map. Define
\[
 \alpha_P:
 \operatorname{Per}_2R\Gamma(S_P,\omega_{P,\Omega})
 \longrightarrow R\Gamma(S_P,\omega_P^{\mathrm{per}})
\]
by adjunction from the periodization of the evaluation map
$a_P^{-1}R\Gamma(S_P,\omega_{P,\Omega})\to\omega_{P,\Omega}$,
followed by the coefficient-change equivalence
$\operatorname{Per}_2\omega_{P,\Omega}\simeq
\omega_P^{\mathrm{per}}$ of
Proposition~\ref{prop:normalized-pole-sheaves}.
These coefficient changes respect the proper traces:
on each finite closed diagram, the dualizing coefficient-change
map is defined by tensoring the adjunction counit with $\Lambda$,
as in Subsection~\ref{subsec:finite-presentations}, and the
normalization indices restrict to the closed factor.
Thus it carries the normalized rational trace to the periodic
trace. The comparison and descent construction of
Proposition~\ref{prop:normalized-pole-sheaves} preserves this
compatibility and its composition coherences. Consequently
the $\alpha_P$ form a morphism of the pole diagrams.
Functoriality of the cofiber of the closed trace gives the
localization comparison. On the open complement its map is
the coefficient-change map defined by the same evaluation
adjunction; equivalently, it is induced by the canonical
exchange map $\operatorname{Per}_2Rj_*\to Rj_*\operatorname{Per}_2$.
No commutation with an arbitrary open direct image is needed.

The zero-contact calculation proves that $\alpha_0$ is an
equivalence. For $P\ge1$, apply the exact functor
$\operatorname{Per}_2$ to the integer localization sequence
above and compare it with the periodic sequence. Its left
vertical map is $\alpha_{P-1}$, its middle vertical map is
$\alpha_P$, and its right vertical map is the finite direct
sum of the contact-piece coefficient changes just proved to
be equivalences. Induction on $P$ proves that every $\alpha_P$
is an equivalence. Finally, periodization preserves the pole
homotopy colimit by \emph{Parity grading and periodization} in
Subsection~\ref{app:complexes-duality}; hence
$\operatorname{Per}_2B_\Omega\simeq B$.

The stages start at $P=0$, and their homotopy colimit is the
complex in \eqref{eq:loop-definition}; thus the filtration is
nonnegative and exhaustive. By
Proposition~\ref{prop:hocolim-properties}, every cohomology
class lies in the image of the cohomology of some finite pole
stage; no completion in contact weight is taken.
\end{proof}

\subsection{Loop cohomology}
\label{subsec:loop-cohomology}

\begin{definition}[Loop cohomology]\label{def:loop-cohomology}
The cohomology of the loop complexes is denoted by
\begin{equation}\label{eq:lh-notation}
 \begin{aligned}
 LH^{\epsilon}(X;\F)
   &:=H^{\epsilon}(B),\\
 LH_{\Omega}^q(X;\F)
   &:=H^q(B_{\Omega}).
 \end{aligned}
\end{equation}
Here $\epsilon\in\Z/2$ and $q\in\Z$.  We usually suppress $\F$.
The specified comparisons of Theorem~\ref{thm:canonical-complex}
identify these groups with $H^*(\mathbf B)$ and
$H^*(\mathbf B_\Omega)$ in the respective gradings.
The \emph{contact filtration} on loop cohomology is
\begin{equation}\label{eq:lh-contact-filtration}
 \begin{aligned}
 F_PLH^{\epsilon}
   &:=\operatorname{im}\bigl(H^{\epsilon}(B_P)
                                  \longrightarrow H^{\epsilon}(B)\bigr),\\
 F_PLH_\Omega^q
   &:=\operatorname{im}\bigl(H^q(B_{\Omega,P})
                                  \longrightarrow H^q(B_\Omega)\bigr),
 \end{aligned}
\end{equation}
with $F_{-1}LH^*=F_{-1}LH_\Omega^*=0$. Equivalently, these
are the images of $H^*(F_P\mathbf B)\to H^*(\mathbf B)$ and
$H^*(F_P\mathbf B_\Omega)\to H^*(\mathbf B_\Omega)$ under
the specified comparisons.
\end{definition}

\begin{remark}[Cohomology of the pole colimit]
\label{rem:pole-cohomology-colimit}
The canonical maps from the pole stages induce isomorphisms
\begin{equation}\label{eq:pole-cohomology-colimit}
 \varinjlim_P H^{\epsilon}(B_P)\xrightarrow{\sim}LH^{\epsilon},
 \qquad
 \varinjlim_P H^q(B_{\Omega,P})\xrightarrow{\sim}LH_\Omega^q.
\end{equation}
Indeed, Definition~\ref{def:loop-complex} expresses each loop
complex as the filtered homotopy colimit of its pole stages,
and filtered colimits of vector spaces are exact; apply
Proposition~\ref{prop:hocolim-properties}(iii), also in parity.
Consequently every cohomology class is represented at a finite
pole stage, and the filtrations in
\eqref{eq:lh-contact-filtration} are exhaustive. A stage class
maps to zero in loop cohomology if and only if its image is
zero at some later stage.

\end{remark}

\begin{definition}[Finite-factor representatives]
\label{def:finite-factor-representatives}
Let $\alpha\in LH_\Omega^q(X)$ and choose a pole-stage lift
$\alpha_P\in H^q(B_{\Omega,P})$ mapping to $\alpha$.
Let $h:U\to F_P\mathscr LX$ be a nice chart with finite-factor
projection $\pi:U\to Y$, working componentwise where the
volume index $\delta_Y=\delta$ is constant.
We say that
\[
 \beta\in H^{BM}_{2\delta-q}(Y;\F)
       =H^qR\Gamma(Y,\omega_Y[-2\delta])
\]
\emph{represents $\alpha$ on this chart with respect to
$\alpha_P$} if
\[
 h^*\alpha_P=\pi^*\beta
 \quad\text{in}\quad
 H^qR\Gamma(U,h^{-1}\omega_{P,\Omega}).
\]
Here both pullbacks are induced by ordinary sheaf inverse
image, with the target of $\pi^*$ identified by
\eqref{eq:nice-chart-dualizing}.
Lemma~\ref{lem:infinite-product-descent}(i) makes $\pi^*$
an isomorphism; hence $\beta$ is uniquely determined by
$\alpha_P$ and the chart. The chosen pole-stage lift is
part of this terminology.
We also call $\beta$ the finite-factor representative of
the pole-stage class $\alpha_P$ on this chart.

In parity grading, use $\alpha\in LH^\epsilon(X)$,
$\alpha_P\in H^\epsilon(B_P)$, and
\[
 \beta\in H^\epsilon R\Gamma
       (Y,\operatorname{Per}_2\omega_Y)
 \cong\bigoplus_{j\equiv\epsilon\; (2)}H_j^{BM}(Y;\F),
\]
with the periodic chart comparison.
\end{definition}

\begin{remark}[Change of volume form]\label{rem:volume-change}
Suppose $\Omega'=u\Omega$ with $u\in\mathcal O(X)^\times$.
Put $\nu_u(\gamma)=\ord_t(u(\gamma))$. On a contact locus,
\[
 \mu_{\Omega'}(\vct)-\mu_\Omega(\vct)
 =2\sum_i\ord_{D_i}(u)v_i
 =2\nu_u(\gamma),\qquad\gamma\in\mathscr P_{\vct}.
\]
The function $\nu_u$ is locally constant on the full Laurent-loop space.
Indeed, if $\nu_u(\gamma)=m$, the conditions
$[t^m]u(\gamma')\ne0$ and $[t^{-m}]u(\gamma')^{-1}\ne0$
define an open neighborhood of $\gamma$, by
Proposition~\ref{lp:prop:topology}. On this neighborhood they imply,
respectively, $\nu_u(\gamma')\le m$ and
$-\nu_u(\gamma')\le-m$, hence $\nu_u(\gamma')=m$.
Its values on each pole stage form a finite set: polynomial
expressions for $u$ and $u^{-1}$ in the affine coordinates of $X$
give a common bound on their poles there, and therefore upper
and lower bounds on $\nu_u$. The stage
complexes therefore split into summands indexed by $\nu_u=m$.
Write $(B_{\Omega,P})_{u=m}$ for the summand on that
open-and-closed locus, taking it to be zero when the locus is
empty. The pole-stage maps preserve these loci. Define
\[
 (B_\Omega)_{u=m}
 :=\hocolim_P(B_{\Omega,P})_{u=m}.
\]
Since homotopy colimits commute with coproducts,
\[
 B_\Omega
 \simeq\hocolim_P\bigoplus_{m\in\Z}(B_{\Omega,P})_{u=m}
 \simeq\bigoplus_{m\in\Z}(B_\Omega)_{u=m}.
\]
Thus the resulting decomposition is a direct sum; each stage
has only finitely many nonzero summands. Use the same subscript
for their cohomology. On the locus $\nu_u=m$, the finite-factor
indices satisfy $\delta_{\Omega'}=\delta_\Omega+m$, so
\[
 \begin{aligned}
 (B_{\Omega'})_{u=m}
       &\simeq (B_\Omega)_{u=m}[-2m],\\
 LH_{\Omega'}^q(X)_{u=m}
       &\cong LH_\Omega^{q-2m}(X)_{u=m}.
 \end{aligned}
\]
The same shift applies to every pole-stage object
$(B_{\Omega,P})_{u=m}$ and preserves its transition maps.
These are decompositions and comparisons of the derived diagrams;
they do not require a literal direct-sum decomposition of the
chosen cochain model $\mathbf B_\Omega$. Thus a class of degree $q$ for $\Omega$ has degree
$q+2m$ for $\Omega'$. The shift may vary between these summands;
it need not be a single shift on the whole complex. It is always
even, so the parity-graded complex $B$ and its cohomology $LH$
are unchanged.
\end{remark}

\subsection{Global finite-factor presentations}
\label{subsec:global-finite-presentations}
For a finite-dimensional complex algebraic variety $Y$, recall that the
Borel--Moore homology is
\begin{equation}\label{eq:bm-duality}
 H_j^{BM}(Y;\F)
 :=H^{-j}R\Gamma(Y,\omega_Y)
 \cong H_c^j(Y;\F)^\vee.
\end{equation}
Here $\omega_Y$ is the Verdier dualizing complex on $Y(\C)$,
and $V^\vee=\operatorname{Hom}_{\F}(V,\F)$.

\begin{proposition}
\label{prop:global-finite-presentations}
Suppose that the pole schemes admit compatible global finite-factor
presentations
\[
 F_P\mathscr LX=\varprojlim_N Z_{P,N}
\]
in the sense of Subsections~\ref{subsec:volume-index-preparations}
and~\ref{subsec:coherent-descent}:
each projection $\pi_{P,N}:F_P\mathscr LX\to Z_{P,N}$
is a finite-factor projection, with a product decomposition
$F_P\mathscr LX\cong Z_{P,N}\times\A^{J_{P,N}}$
for a finite or countable set $J_{P,N}$. Require
closed inclusions $Z_{P,N}\hookrightarrow Z_{P+1,N}$,
affine-bundle projections $Z_{P,N+1}\to Z_{P,N}$, and Cartesian
squares. The finite-dimensional schemes $Z_{P,N}$ may be singular
or nonreduced. For the integer-graded statements, fix a
nowhere-vanishing volume form $\Omega$ and let $\delta_{P,N}$
be the locally constant volume index of the projection to
$Z_{P,N}$. Work at sufficiently large cutoffs for the normalization
and each specified closed diagram to descend. Then:
\begin{enumerate}[(i),wide=0pt]
\item \emph{The filtered loop complex.}
The finite factors compute the derived pole stages and the loop object:
\[
 \begin{aligned}
 B_{\Omega,P}
   &\simeq R\Gamma(Z_{P,N},
                    \omega_{Z_{P,N}}[-2\delta_{P,N}]),\\
 B_\Omega
   &\simeq\hocolim_P\hocolim_N
          R\Gamma(Z_{P,N},\omega_{Z_{P,N}}[-2\delta_{P,N}]).
 \end{aligned}
\]
The affine-bundle projections induce normalized oriented pullback
quasi-isomorphisms between the finite models of the same stage
$B_{\Omega,P}$. The proper traces of the closed
inclusions represent the filtration maps
\[
 B_{\Omega,P}\longrightarrow
 B_{\Omega,P+1}.
\]
These identifications are coherent under refinement and composition,
so they identify the entire original pole-stage diagram. Its
comparison with $(F_P\mathbf B_\Omega)_P$ is the specified
equivalence of pole-stage diagrams in Theorem~\ref{thm:canonical-complex}. The $N$-maps in the
homotopy colimit are the normalized pullbacks, and the $P$-maps are
the proper traces.

\item \emph{The associated graded complexes.}
Put $U_{P,N}=Z_{P,N}\setminus Z_{P-1,N}$ for $P\ge1$, and
$U_{0,N}=Z_{0,N}$. Then
\[
 \gr_P\mathbf B_\Omega\simeq
 R\Gamma(U_{P,N},
   \omega_{U_{P,N}}[-2\delta_{P,N}|_{U_{P,N}}]).
\]
Without a volume form, the assertions in \textup{(i)} and
\textup{(ii)} hold for $B$ after two-periodization,
with the even normalization shifts omitted.

\item \emph{Loop cohomology.}
There are isomorphisms
\begin{equation}\label{eq:bm-limits}
 \begin{aligned}
 LH^q_{\Omega}(X)
  &\cong\varinjlim_P\varinjlim_N
    \bigoplus_{C\in\pi_0(Z_{P,N})}
    H^{BM}_{2\delta_{P,N}(C)-q}(C;\F),\\
 LH^{\epsilon}(X)
  &\cong\varinjlim_P\varinjlim_N
    \bigoplus_{j\equiv\epsilon\ (2)}
    H^{BM}_j(Z_{P,N};\F).
 \end{aligned}
\end{equation}
The $P$-maps are proper pushforwards, and the $N$-maps are the
oriented affine-bundle pullbacks
\[
 p^!:H_j^{BM}(Z_{P,N};\F)
       \xrightarrow{\sim}H_{j+2e}^{BM}(Z_{P,N+1};\F),
\]
where $e$ is the complex relative dimension on the component in
question. They commute and preserve the normalized degree, or
parity. The contact filtration on $LH$ consists of the images
of the fixed-$P$ groups in \eqref{eq:bm-limits}; these groups
need not inject into the limit. Moreover,
\begin{equation}\label{eq:coefficient-layer}
 H^q(\gr_P\mathbf B_\Omega)\cong
 \bigoplus_{C\in\pi_0(U_{P,N})}
 H^{BM}_{2\delta_{P,N}(C)-q}(C;\F).
\end{equation}
All Borel--Moore groups here are those of the analytic
complex-point spaces. For a smooth connected component $C$ of
complex dimension $d$, its normalized term is
$H^{q+2d-2\delta_{P,N}(C)}(C;\F)$, and is ordinary cohomology
in degree $q$ when $\delta_{P,N}(C)=d$.
\end{enumerate}
\end{proposition}

\begin{proof}
For \textup{(i)}, the global finite-factor hypothesis and
Proposition~\ref{prop:normalized-pole-sheaves} identify the
pole-stage sheaf with the inverse image of
$\omega_{Z_{P,N}}^{\mathrm n}=\omega_{Z_{P,N}}[-2\delta_{P,N}]$ under the finite-factor
projection $\pi:F_P\mathscr LX(\C)\to Z_{P,N}(\C)$.
By the global finite-factor hypothesis, this projection is the
projection from $Z_{P,N}(\C)\times\C^J$ for a finite or countable
set $J$. Lemma~\ref{lem:infinite-product-descent}(i) makes the
adjunction unit
\[
 \omega_{Z_{P,N}}^{\mathrm n}\longrightarrow R\pi_*\pi^{-1}\omega_{Z_{P,N}}^{\mathrm n}
\]
an equivalence. Applying derived global sections proves the
first formula without an additional nice-chart atlas.

For an affine-bundle projection $p$ of complex relative dimension
$e$, the indices satisfy
$\delta_{P,N+1}=p^*\delta_{P,N}+e$. The normalized purity
equivalence \eqref{eq:enhanced-normalized-purity} is therefore
\[
 \omega_{Z_{P,N+1}}[-2\delta_{P,N+1}]
 \simeq p^{-1}\!\left(\omega_{Z_{P,N}}[-2\delta_{P,N}]\right).
\]
An affine bundle is analytically locally a product with $\C^e$.
On each trivializing open subset,
Lemma~\ref{lem:infinite-product-descent}(i) therefore identifies
the adjunction unit $\omega_{Z_{P,N}}^{\mathrm n}\to Rp_*p^{-1}\omega_{Z_{P,N}}^{\mathrm n}$ as an equivalence.
Equivalences of sheaf complexes are detected on an open cover,
so the unit is an equivalence globally. Applying derived global
sections to this unit and then using the inverse of the displayed
normalized purity equivalence gives the normalized pullback
quasi-isomorphism.
Along a closed inclusion the indices
restrict to each other, so its normalized proper trace represents
the pole-stage transition map. The coherent proper--smooth
base-change and composition compatibilities in
Proposition~\ref{prop:normalized-pole-sheaves} identify these maps
under increasing $N$. Thus the finite models recover the entire
coherent pole-stage diagram of
Remark~\ref{rem:coherent-pole-descent}. Taking its homotopy colimit, as in
Definition~\ref{def:loop-complex}, proves the second formula in
\textup{(i)}. The inner homotopy colimit identifies the different
coefficient cutoffs; the outer one assembles the filtered stages.

For \textup{(ii)}, apply the dualizing localization triangle to
$Z_{P-1,N}\hookrightarrow Z_{P,N}$ with its common normalization
shift and take derived global sections. The resulting cofiber
is the displayed complex on $U_{P,N}$. The comparison of the
entire pole-stage diagram with $(F_P\mathbf B_\Omega)_P$
identifies this cofiber with the quotient
$F_P\mathbf B_\Omega/F_{P-1}\mathbf B_\Omega$.

For the periodic assertions, let $Z$ be a finite factor or one
of these open complements, and put $C_Z:=R\Gamma_c(Z;\bQ)$.
This complex is perfect: it is equivalent to a bounded complex
of finite-dimensional rational vector spaces. Recall
$\Lambda=\bQ[u,u^{-1}]$, with $|u|=2$ and $du=0$.
In the unfolded differential graded $\Lambda$-module model,
Verdier duality and the dualizing coefficient-change equivalence
of Subsection~\ref{subsec:finite-presentations} identify the
natural comparison as
\[
 \begin{aligned}
 \operatorname{Per}_2R\Gamma(Z,\omega_Z)
 &\simeq R\operatorname{Hom}_{\bQ}(C_Z,\bQ)
                      \otimes_{\bQ}\Lambda\\
 &\xrightarrow{\sim}R\operatorname{Hom}_{\bQ}(C_Z,\Lambda)\\
 &\simeq R\Gamma(Z,\operatorname{Per}_2\omega_Z).
 \end{aligned}
\]
The first and last equivalences are the duality adjunctions
for the map $Z(\C)\to\mathrm{pt}$; the middle equivalence
uses perfectness of $C_Z$. These are the canonical evaluation
and adjunction comparisons, so their compatibility with proper
traces and oriented smooth pullbacks is the one constructed in
Proposition~\ref{prop:normalized-pole-sheaves}. No smoothness
or reducedness of $Z$, and no volume form, is needed.
Apply the product, purity, and trace arguments of \textup{(i)}
in the periodic coefficient category. The displayed equivalence
identifies the finite models with the indicated periodizations;
the even purity shifts have their canonical periodic
identifications. Periodization preserves the localization
cofiber sequences and the homotopy colimits, proving the
periodic assertions in \textup{(i)} and \textup{(ii)}.

For \textup{(iii)}, cohomology commutes with filtered homotopy
colimits by Proposition~\ref{prop:hocolim-properties}. Apply
\eqref{eq:bm-duality} to the finite factors in \textup{(i)} and
\textup{(ii)}. Their compactly supported cohomology is finite
dimensional and bounded, so its dual and periodization give
exactly the displayed Borel--Moore groups. Under these
identifications, proper trace becomes proper pushforward and
normalized pullback becomes the oriented affine-bundle pullback.
The index change by $e$ cancels its homological degree change
by $2e$, while proper--smooth base change gives commutativity.
The filtration on cohomology is, by definition, the image of the
cohomology of each filtered stage; no injectivity of the
proper-trace maps is asserted. Finally, smooth duality gives
$\omega_C\simeq\F_C[2d]$, proving the ordinary-cohomology
formula. The layer formula is the finite-model expression for
the contact cohomology of Theorem~\ref{thm:canonical-complex}.
\end{proof}

\subsection{Examples}
\label{subsec:loop-examples}

\begin{example}[Affine line -- continued]
For the affine line with $\Omega=dx$, the calculation in Example \ref{ex:affine-line-volume-index} gives
$B_{\Omega,P}\simeq\bQ[2P]$. Its transition maps
are represented by proper traces for the coordinate hyperplanes
$\A^{P+N+1}\hookrightarrow\A^{P+N+2}$. Every such pole-stage morphism is zero
in the derived category, since
\[
 \operatorname{Hom}_{D(\bQ)}
       (\bQ[2P],\bQ[2P+2])=0.
\]
This is a statement about the derived morphisms
$B_{\Omega,P}\to B_{\Omega,P+1}$, represented by the
inclusions $F_P\mathbf B_\Omega\hookrightarrow
F_{P+1}\mathbf B_\Omega$ under the specified comparisons.
Cohomology commutes with the filtered homotopy colimit, so
\begin{equation}\label{eqBA1}
 B_{\Omega}\simeq0,
 \qquad LH_{\Omega}^q(\A^1)=0
 \quad\text{for every }q.
\end{equation}

Let us illustrate two more points by embedding the same affine
line differently, as the parabola $X=\{y=x^2\}\subset\A^2$, still with
$\Omega=dx$. This embedding comes from the sections $1,x,x^2$
of $\mathcal O_{\Pj^1}(2\infty)$: we keep $A=\{\infty\}$
and take $r=2$ in \eqref{eq:pole-stage-scheme}.
Consider the Laurent coefficient scheme $S=F_1\mathscr LX=\Spec C_2(X)$,
where both $x$ and $y$ have pole order at most two. Over a complex
algebra $R$, write
\[
 x(t)=\varepsilon t^{-2}+at^{-1}+\sum_{j\ge0}x_jt^j.
\]
The coefficients of $t^{-4}$ and $t^{-3}$ in $x^2$ give
$\varepsilon^2=\varepsilon a=0$. All remaining $y$-coefficients
are determined by $y=x^2$, so elimination gives
\[
 S\cong\Spec
 \frac{\C[\varepsilon,a,x_0,x_1,\ldots]}
      {(\varepsilon^2,\varepsilon a)}.
\]
For $N\ge0$, put
\[
 E_N=\frac{\C[\varepsilon,a,x_0,\ldots,x_N]}
              {(\varepsilon^2,\varepsilon a)},
 \qquad Y_N=\Spec E_N.
\]
Then $S\cong Y_N\times\Spec\C[x_{N+1},x_{N+2},\ldots]$ is a
global nice chart, and $Y_{N+1}=Y_N\times\A^1\to Y_N$ is an
affine-line projection. We compute its index directly using
Proposition~\ref{prop:general-factor-contact-index}. The contact-one
and contact-zero loci have reduced finite factors
$\{\varepsilon=0,\ a\ne0\}\subset Y_N$ and
$\{\varepsilon=a=0\}\subset Y_N$, respectively, of dimensions
$N+2$ and $N+1$. Both already descend to $Y_N$, so the proposition
applies with $e=0$. Since $\mu_\Omega(q)/2=-q$ for contact order
$q\in\{0,1\}$, it gives
\[
 \delta_{Y_N}=
 \begin{cases}
  (N+2)-1=N+1,&a\ne0,\\
  (N+1)+0=N+1,&a=0.
 \end{cases}
\]
Thus $\delta_{Y_N}=N+1$ at every complex point. Since
$Y_N(\C)=\C^{N+2}$, we obtain
\[
 R\Gamma(Y_N,\omega_{Y_N}[-2\delta_{Y_N}])\simeq\bQ[2]
\]
which agrees with the standard $P=1$ calculation above.
 Note that $S_{\mathrm{red}}$ is the standard pole-one stage of the affine line. The first point we wanted to make with this is that the finite factors in our construction of the loop complex can be nonreduced.

The second point is that a natural coefficient presentation of $\sL X$ can
fail the affine-bundle hypothesis of
Proposition~\ref{prop:global-finite-presentations}. Consider the
literal cutoffs
\[
 Z_{P,N}:=\Spec C_{2P,N}(X),
 \qquad F_P\mathscr LX=\varprojlim_N Z_{P,N},
\]
where $C_{2P,N}(X)$ retains the coefficients of both $x$ and $y$
through order $N$, as in Definition~\ref{def:laurent-coefficients}.
We check the failure at $P=1$. On reduction, set $b=x_0$,
$c=x_1$, and $d=x_2$. The formulas $y_0=b^2+2ac$ and
$y_1=2bc+2ad$ give
\[
 \begin{aligned}
 \mathcal O((Z_{1,0})_{\mathrm{red}})
   &=\C[a,b,ac]\cong\C[a,b,s_0],\\
 \mathcal O((Z_{1,1})_{\mathrm{red}})
   &=\C[a,b,c,ad]\cong\C[a,b,c,s_1].
 \end{aligned}
\]
These subalgebras lie in $\mathcal O(S_{\mathrm{red}})$;
reduction commutes with the inclusions by contraction of the
nilradical. The reduced transition map
$(Z_{1,1})_{\mathrm{red}}\to(Z_{1,0})_{\mathrm{red}}$ is
\[
 \A^4\longrightarrow\A^3,\qquad
 (a,b,c,s_1)\longmapsto(a,b,s_0=ac).
\]
Its fiber is $\A^1$ when $a\ne0$, is $\A^2$ when $a=s_0=0$,
and is empty when $a=0$ and $s_0\ne0$. It is therefore not an
affine bundle. An affine bundle remains an affine bundle after
reduction, so the original map $Z_{1,1}\to Z_{1,0}$ cannot be
one either. The failure persists at every cutoff: from
\[
 y_N=2ax_{N+1}+\sum_{i=0}^Nx_ix_{N-i}
\]
one gets
$\mathcal O((Z_{1,N})_{\mathrm{red}})
 =\C[a,x_0,\ldots,x_N,ax_{N+1}]$,
and the next reduced transition has the same fiber jump.
Nevertheless, inside $\mathcal O(S)$ one has
\[
 E_N\subset\mathcal O(Z_{1,N})\subset E_{N+2},
\]
so both systems present the same scheme $S$. Thus the system $Z_{P,N}$ does not satisfy
Proposition~\ref{prop:global-finite-presentations}, although it is
a coefficient presentation of the pole stages. 
\end{example}

\begin{example}[The multiplicative group]
\label{ex:cstar-loop-complex}
Take $X=\C^*$ with coordinate $x$, compactification
$M=\Pj^1$, boundary $D=\{0,\infty\}$, polarization $A=D$,
and volume form $\Omega=dx/x$. The coordinates $x,x^{-1}$
come from $\mathcal O_M(A)$, so $r=1$ and the pole scheme
$F_P\mathscr LX$ bounds both Laurent series by $-P$.
Every complex loop has a unique expression
\[
 x(t)=t^v(a_0+a_1t+a_2t^2+\cdots),\qquad a_0\ne0.
\]
Its winding $v$ satisfies $-P\le v\le P$ in $F_P\mathscr LX$.
Positive and negative windings have contact $v$ at zero and
$-v$ at infinity, respectively, and winding zero is the arc
locus of $X$. Since $\Omega$ has simple poles at both boundary
points, every contact shift $\mu_\Omega$ is zero.

The reduction of this pole scheme is therefore a finite disjoint
union of $\Gm\times\A^{\Z_{>0}}$. We retain the nilpotents to exhibit actual
global finite factors. Write $y=x^{-1}$,
$x_j=[t^j]x(t)$, and $y_j=[t^j]y(t)$, and fix a winding
$v\in[-P,P]\cap\Z$. The principal open where $x_vy_{-v}$ is invertible
is precisely the winding-$v$ part: over every residue field,
both coefficients are nonzero exactly for winding $v$. As $v$ ranges over the integers $-P,\ldots,P$, the principal
opens $D(x_vy_{-v})$ form a disjoint clopen decomposition
of the pole scheme.
 Put $m=P+v$ and $n=P-v$, and let $T_{P,v}$
be the finite-type scheme of monic polynomial pairs
\[
 q(t)=t^m+q_{m-1}t^{m-1}+\cdots+q_0,\qquad
 p(t)=t^n+p_{n-1}t^{n-1}+\cdots+p_0,
 \qquad q(t)p(t)=t^{2P}.
\]
Its only closed point is $(q,p)=(t^m,t^n)$, so $T_{P,v}$
is a finite scheme supported at one point; it may be nonreduced.

\noindent\emph{Claim.} The entire winding-$v$ scheme, that is, the open-and-closed subscheme
$D(x_vy_{-v})\subset F_P\mathscr LX$, has the finite-factor presentation
\[
 D(x_vy_{-v})
 \simeq T_{P,v}\times\Gm\times\A^{\Z_{>0}}.
\]
\begin{proof}
For every complex algebra $R$, the map from right to
left sends a pair $(q,p)\in T_{P,v}(R)$ and coefficients
$a_0\in R^\times$, $a_j\in R$ for $j\ge1$, to the Laurent
series
\begin{equation}\label{eq:cstar-full-sector-factor}
 x=t^{-P}q(t)a(t),\qquad y=t^{-P}p(t)a(t)^{-1},
 \qquad a(t)=a_0+a_1t+a_2t^2+\cdots.
\end{equation}

Here is an elementary justification over an arbitrary coefficient
ring $R$. The coefficients of $t^Px$ below degree $m$, and of
$t^Py$ below degree $n$, generate a nilpotent ideal $I$: they
vanish in every residue field, and there are finitely many of
them. Put $f=t^Px$. Modulo $I$, this series is $t^m$ times a unit.
Starting with this factorization, suppose $f=qa$ modulo $I^k$.
Divide the error by the unit $a$ and write
$(f-qa)/a=r+t^ms$, where $r\in I^k[t]$ has degree less than
$m$ and $s\in I^k\llbracket t \rrbracket$; this splitting is unique.
Since $q\equiv t^m\pmod I$ and $k\ge1$,
\[
 (q+r)a(1+s)\equiv qa+a(r+t^ms)
                  =f\pmod{I^{k+1}}.
\]
Thus replacing $q$ by $q+r$ and $a$ by $a(1+s)$ corrects the
factorization modulo $I^{k+1}$, preserving monicity and the
invertibility of the constant coefficient of $a$.
Nilpotence makes this terminate. For uniqueness, suppose
$f=qa=q'a'$ are two such factorizations, where $q,q'$ are
monic of degree $m$, their lower coefficients are nilpotent,
and $a,a'\in R\llbracket t \rrbracket^\times$. Enlarge $I$ by the finitely many
lower coefficients of $q$ and $q'$; it remains nilpotent.
Modulo $I$, both monic factors equal $t^m$, so $a=a'$ modulo
$I$, since multiplication by $t^m$ is injective on
$(R/I)\llbracket t \rrbracket$. If the two factorizations agree coefficientwise
modulo $I^k$, then
\[
 0=(q-q')+t^m(a-a')/a\pmod{I^{k+1}}.
\]
The first summand has degree less than $m$, whereas the second
has no terms of degree less than $m$. Both therefore vanish
modulo $I^{k+1}$. Induction proves $q=q'$ and $a=a'$.
Applying this uniqueness statement to the product of the
factorizations of $t^Px$ and $t^Py$ gives $qp=t^{2P}$ and
makes their unit factors inverse. Uniqueness also makes these
factorizations compatible with every homomorphism of
coefficient rings. Conversely, for an $R$-point $(q,p)$ of
$T_{P,v}$ and $a\in R\llbracket t \rrbracket^\times$, the product formula gives
$xy=1$ and $x_v,y_{-v}\in R^\times$: modulo the nilpotent
ideal generated by the lower coefficients of $q,p$, those two
coefficients are $a_0$ and $a_0^{-1}$. Thus
\eqref{eq:cstar-full-sector-factor} defines the asserted
isomorphism of schemes.
\end{proof}

For $N\ge0$, retain $a_0,\ldots,a_N$ and set
\[
 Z_{P,N}=\coprod_{v=-P}^{P}
       \bigl(T_{P,v}\times\Gm\times\A^N\bigr).
\]
Then $F_P\mathscr LX=\varprojlim_N Z_{P,N}$, with an
independent affine coefficient tail. Increasing $N$ adds one
affine coordinate. Increasing $P$ sends $(q,p)$ to $(tq,tp)$
and retains $a(t)$; this gives closed inclusions
$Z_{P,N}\hookrightarrow Z_{P+1,N}$ and Cartesian squares.
On each component the reduced contact locus has finite factor
$\Gm\times\A^N$, of dimension $N+1$. It already descends to
$Z_{P,N}$, so Proposition~\ref{prop:general-factor-contact-index},
with $e=0$ and $\mu_\Omega=0$, gives
\[
 \delta_{P,N}=N+1.
\]
The complex-point space of this component is the smooth space
$\C^*\times\C^N$, whose dualizing complex is
$\F_{\C^*\times\C^N}[2N+2]$. 
Proposition~\ref{prop:global-finite-presentations}
therefore gives
\[
 B_{\Omega,P}
   \simeq\bigoplus_{v=-P}^{P}C^*(\C^*;\F).
\]
The $N$-maps are the normalized affine-line pullbacks. On every
existing winding, the $P$-map retains $a_0,\ldots,a_N$ and
induces the identity of $\C^*\times\C^N$, since both
$T_{P,v}(\C)$ and $T_{P+1,v}(\C)$ consist of one point.
The dualizing complexes and proper traces used here belong
to these analytic spaces. The restricted trace is therefore
the adjunction counit for the identity map, hence the identity;
the common normalization shift $-2(N+1)$ preserves this
identification. The two new windings add summands. Consequently
\[
 B_\Omega
   \simeq\bigoplus_{v\in\Z}C^*(\C^*;\F),\qquad
 LH^j_\Omega(\C^*)\cong
 \begin{cases}
  \displaystyle\bigoplus_{v\in\Z}\F,&j=0,1,\\
  0,&j\ne0,1.
 \end{cases}
\]
These are direct sums: the homotopy colimit adds finitely many
windings at each stage. The contact filtration consists precisely
of the summands with $|v|\le P$. 

For future use, define
$e_v\in LH_\Omega^0(\C^*)$ and
$\eta_v\in LH_\Omega^1(\C^*)$ as follows.
The preceding equivalence identifies the winding-$v$ summand
of $LH_\Omega^j(\C^*)$ with $H^j(\C^*;\F)$, where the
$\C^*$ coordinate is the leading coefficient
$a_0=[t^v]x$, also when $v<0$.
Let $\theta\in H^1(\C^*;\F)$ be the unique class satisfying
\[
 \langle\theta,[S^1]\rangle=1,
\]
where $[S^1]\in H_1(\C^*;\F)$ is the fundamental class of
the counterclockwise oriented unit circle.
We define $e_v$ and $\eta_v$ to have winding-$v$ components
$1\in H^0(\C^*;\F)$ and $\theta$, respectively, and zero
components in every other winding.

 Two-periodization gives
$B\simeq
\bigoplus_{v\in\Z}\operatorname{Per}_2C^*(\C^*;\F)$ and
$LH^{\bar0}(\C^*)\cong
LH^{\bar1}(\C^*)\cong\bigoplus_{v\in\Z}\F$.
Finally, every other algebraic volume form is
$\Omega'=c x^k\,dx/x$ with $c\in\C^*$ and $k\in\Z$.
Remark~\ref{rem:volume-change} gives the winding-$v$ summand
$C^*(\C^*;\F)[-2kv]$, with its two classes in degrees
$2kv$ and $2kv+1$.
\end{example}

\section{Topological preparations for the fusion product}
\label{sec:topological-preparations}

We recall first the classical specialization maps for finite-dimensional families. The goal then will be to normalize and glue them in our infinite-dimensional loop space context. We give in this section the general framework for this construction. We return to the loop space context in the next section to apply the construction from here.

\subsection{Sheaf operations and coherent diagrams}
\label{conv:section-five-sheaf-operations}
We use the sheaf categories and enhanced operations fixed in
Subsection~\ref{subsec:coherent-descent}. The finite-dimensional
spaces here are analytic spaces of separated complex schemes of
finite type and their locally closed subanalytic subsets, and
the paracompact locally compact Hausdorff parameter spaces of
finite covering dimension specified below. The same topological
formalism applies to the latter, in particular to manifolds with
corners and locally finite finite-dimensional polyhedra.
For such a space $T$, put
\[
 a_T:T\to\mathrm{pt},\qquad
 \omega_T:=a_T^!\F,\qquad
 \mathbb D_T E:=R\mathcal Hom_T(E,\omega_T).
\]
For a continuous map $h:M\to T$ between these spaces, the
projection formula and $h^!\omega_T\simeq\omega_M$ give
\begin{equation}\label{eq:duality-inverse-image}
 \mathbb D_Mh^{-1}E
 =R\mathcal Hom_M(h^{-1}E,h^!\omega_T)
 \simeq h^!R\mathcal Hom_T(E,\omega_T)
 =h^!\mathbb D_TE.
\end{equation}
Indeed, the middle equivalence is the right-adjoint form of
$Rh_!(V\otimes h^{-1}E)\simeq Rh_!V\otimes E$.
It requires neither constructibility nor biduality.

Adjunction, localization, purity, and their comparison maps are
taken as natural transformations in $\mathcal D(-;\F)$,
with their composition and base-change homotopies. The statements
in the ordinary categories $D(-;\F)$ below are their images
under passage to homotopy categories. A coherent diagram means
a functor $N(A)\to\mathcal D(-;\F)$, as in
Appendix~\ref{app:coherent-diagrams}; cochain representatives
are chosen for whole diagrams as specified there. Use
$\mathcal D_{\mathrm{per}}(-;\F)$ for the corresponding parity
constructions. No constructibility condition is imposed unless
stated.

\subsection{Specialization of cycles}
\label{subsec:classical-specialization}
We recall specialization first on dualizing complexes, then as the
dual of a restriction map, and finally on Borel--Moore homology.

\begin{lemma}[Specialization on dualizing complexes]
\label{lp:lem:bm-collision-specialization}
Let $\pi:\mathcal Y\to\A^1_s$ be a morphism of separated
complex schemes of finite type. Set
\[
 \mathcal Y_0=\mathcal Y\times_{\A^1_s}\{0\},\qquad
 \mathcal Y^\times=\mathcal Y\times_{\A^1_s}\Gm,
 \qquad \Gm=\Spec\C[s,s^{-1}].
\]
Write $i:\mathcal Y_0\hookrightarrow\mathcal Y$ and
$j:\mathcal Y^\times\hookrightarrow\mathcal Y$ for the inclusions.
Let $L=\R_{>0}\subset\C^*$, oriented from infinity toward zero,
and put $\mathcal Y_L=\pi^{-1}(L)$, with closed inclusion
$\ell:\mathcal Y_L\hookrightarrow\mathcal Y^\times$.

The following statements are in $D(\mathcal Y;\F)$ for the
analytic topology, with their enhancements specified as in
Subsection~\ref{conv:section-five-sheaf-operations}. The intermediate sheaves involving
$L$ are constructible for real semialgebraic stratifications. Put
\[
 \mathcal K_{\mathcal Y,\rho}
       :=Rj_*\ell_*\omega_{\mathcal Y_L}[-1].
\]
The subscript $\rho$ records the chosen inward-oriented ray $L$.
Here $Rj_*\ell_*=(Rj_*)\circ\ell_*\simeq R(j\circ\ell)_*$,
since the closed embedding $\ell$ has exact direct image.
The proper trace for $\ell$ and the boundary of the dualizing
localization triangle define a morphism in
$D(\mathcal Y;\F)$,
\begin{equation}\label{eqSp}
 \operatorname{Sp}^{\omega}_{\mathcal Y}:
 \mathcal K_{\mathcal Y,\rho}
 \xrightarrow{\ Rj_*\operatorname{Tr}_{\ell}[-1]\ }
 Rj_*\omega_{\mathcal Y^\times}[-1]
 \xrightarrow{\ \partial_\omega[-1]\ }
 i_*\omega_{\mathcal Y_0}.
\end{equation}
Equivalently, put $H=[0,\infty)$ and $\mathcal Y_H=\pi^{-1}(H)$,
with inclusions
\[
 b:\mathcal Y_0\hookrightarrow\mathcal Y_H,
 \qquad k:\mathcal Y_L\hookrightarrow\mathcal Y_H,
 \qquad a:\mathcal Y_H\hookrightarrow\mathcal Y.
\]
If $\partial_H:Rk_*\omega_{\mathcal Y_L}\to
b_*\omega_{\mathcal Y_0}[1]$ is the boundary for
$\mathcal Y_0\subset\mathcal Y_H$, then
\begin{equation}\label{eq:specialization-half-ray-boundary}
 \operatorname{Sp}^{\omega}_{\mathcal Y}
       =a_*\bigl(\partial_H[-1]\bigr),\qquad
 a_*Rk_*\simeq Rj_*\ell_*,\quad a_*b_*=i_*.
\end{equation}
Both descriptions define the specialization morphism without
assuming or choosing a product trivialization of the punctured
family $\mathcal Y^\times\to\Gm$. Here $\rho$ denotes the image of the oriented fundamental class
$[L]$ in $H^{\mathrm{BM}}_1(\C^*;\F)$. 

If $U$ is a separated complex scheme of finite type and
$\vartheta:\Gm\times U\xrightarrow{\sim}\mathcal Y^\times$
is an isomorphism over $\Gm$, the orientation identifies
\[
 R\Gamma(\mathcal Y,\mathcal K_{\mathcal Y,\rho})
 \simeq R\Gamma(\mathcal Y_L,\omega_{\mathcal Y_L})[-1]
 \simeq R\Gamma(U,\omega_U).
\]
Thus derived global sections give a degree-zero morphism
\begin{equation}\label{eqSpG}
 \operatorname{Sp}_{\mathcal Y}:
 R\Gamma(U,\omega_U)\longrightarrow
 R\Gamma(\mathcal Y_0,\omega_{\mathcal Y_0}).
\end{equation}
For $\mathcal Y=\A^1\times U$ with its product trivialization,
this morphism is the identity under $\mathcal Y_0\cong U$.

These constructions have the following functoriality. Let
$f:\mathcal Y'\to\mathcal Y$ be a map over $\A^1_s$
between such families, with restrictions $f_0,f^\times,f_L$
to their closed fibers, punctured families, and positive rays.
If $f$ is proper, the trace for $f_L$ gives
$t_f:Rf_*\mathcal K_{\mathcal Y',\rho}\to
\mathcal K_{\mathcal Y,\rho}$, and
\[
 \operatorname{Sp}^{\omega}_{\mathcal Y}\circ t_f
   =i_*\operatorname{Tr}_{f_0}\circ
                   Rf_*\operatorname{Sp}^{\omega}_{\mathcal Y'}.
\]
The target of the second factor on the right is identified with
$i_*Rf_{0,*}\omega_{\mathcal Y'_0}$.
If $f$ is smooth of pure complex relative dimension $e$, the
Cartesian open, closed, and ray squares and oriented purity give
\[
 \begin{gathered}
 \mathcal K_{\mathcal Y',\rho}
       \simeq f^*\mathcal K_{\mathcal Y,\rho}[2e],\qquad
 i'_*\omega_{\mathcal Y'_0}
       \simeq f^*(i_*\omega_{\mathcal Y_0})[2e],\\
 \operatorname{Sp}^{\omega}_{\mathcal Y'}
       =f^*\operatorname{Sp}^{\omega}_{\mathcal Y}[2e]
       \quad\text{under these identifications}.
 \end{gathered}
\]
When product trivializations $\vartheta,\vartheta'$ are chosen and
$f^\times\vartheta'=\vartheta(\mathrm{id}_{\Gm}\times u)$
for $u:U'\to U$, the ray restriction is $f_L=\mathrm{id}_L\times u$;
taking sections gives the corresponding proper and smooth
comparisons for \eqref{eqSpG}.
\end{lemma}
\begin{proof}
$L$ is closed in $\C^*$, so $\ell$ is a proper closed
embedding. The adjunction counit defines
\[
 \operatorname{Tr}_{\ell}:
 \ell_*\omega_{\mathcal Y_L}
   \simeq\ell_*\ell^!\omega_{\mathcal Y^\times}
   \longrightarrow\omega_{\mathcal Y^\times}.
\]
The dualizing localization triangle
\[
 i_*\omega_{\mathcal Y_0}\longrightarrow\omega_{\mathcal Y}
 \longrightarrow Rj_*\omega_{\mathcal Y^\times}
 \xrightarrow{\partial_\omega}i_*\omega_{\mathcal Y_0}[1]
\]
then defines $\operatorname{Sp}^{\omega}_{\mathcal Y}$.
The proper trace $a_*\omega_{\mathcal Y_H}\to\omega_{\mathcal Y}$
extends to a morphism from the direct image of the dualizing
localization triangle on $\mathcal Y_H$ to that on $\mathcal Y$.
Its closed component is $\mathrm{id}_{i_*\omega_{\mathcal Y_0}}$;
its open component is $Rj_*\operatorname{Tr}_\ell$.
Naturality of the boundaries therefore gives
\[
 a_*\partial_H
   =\partial_\omega\circ Rj_*\operatorname{Tr}_\ell,
\]
which proves \eqref{eq:specialization-half-ray-boundary}.
When $\vartheta$ is given, it identifies $\mathcal Y_L$ with
$L\times U$. The orientation gives $\omega_L\simeq\F_L[1]$, and
$R\Gamma(L,\F_L)\simeq\F$ because $L$ is contractible.
All tensor products here and below are over $\F$. Tensor products
and permutations use the Koszul rule; for homogeneous maps and
cochains,
$$
 (f\otimes g)(a\otimes b)=(-1)^{|g||a|}f(a)\otimes g(b),
 \qquad
 d_{C\otimes D}=d_C\otimes\mathrm{id}_D+
                    \mathrm{id}_C\otimes d_D.
$$
In parity grading the degrees are read modulo $2$.
Product duality and K\"unneth therefore give
\[
 R\Gamma(\mathcal Y_L,\omega_{\mathcal Y_L})[-1]
 \simeq
 \bigl(R\Gamma(L,\omega_L)\otimes
       R\Gamma(U,\omega_U)\bigr)[-1]
 \simeq R\Gamma(U,\omega_U).
\]
This proves the asserted source identification and defines
$\operatorname{Sp}_{\mathcal Y}$.

For proper $f$, composition of direct images identifies
\[
 Rf_*\mathcal K_{\mathcal Y',\rho}
 \simeq Rj_*\ell_*Rf_{L,*}\omega_{\mathcal Y'_L}[-1].
\]
Applying the trace of $f_L$ defines $t_f$.
Composition of proper traces and naturality of the localization
triangle give the stated equality. For smooth $f$, extraordinary base change gives
\[
 \begin{gathered}
 f^!Rj_*\simeq Rj'_* (f^\times)^!,\qquad
 (f^\times)^!\ell_*\simeq\ell'_*f_L^!,\qquad
 f^!i_*\simeq i'_*f_0^!.
 \end{gathered}
\]
The first identity is dual to base change for extension by zero;
the other two use proper base change. Since extraordinary
inverse image preserves dualizing complexes, these identities
carry the trace and localization boundary to those of
$\mathcal Y'$. Purity $f^!\simeq f^*[2e]$ gives the asserted
smooth comparison.

For the product family the map is the tensor product of the
scalar boundary
$R\Gamma(L,\omega_L)[-1]\to\F$ with the identity on
$R\Gamma(U,\omega_U)$. The inward-oriented ray has boundary
$[0]$, so that scalar map sends $1$ to $1$.
Thus $\operatorname{Sp}_{\A^1\times U}=\mathrm{id}$.
No properness of $\pi$ or smoothness of $\mathcal Y_0$ is used.
\end{proof}

\begin{lemma}[Specialization as the dual of restriction]
\label{lem:specialization-restriction}
Assume that the product trivialization $\vartheta$ in
Lemma~\ref{lp:lem:bm-collision-specialization} is given. For its
half-ray $H$, write
$h:H\hookrightarrow\C$, $i_H:\{0\}\hookrightarrow H$ and
$j_H:L\hookrightarrow H$. Put
\[
 \begin{gathered}
 \mathcal Y_H=\pi^{-1}(H),\qquad
 \pi_H:\mathcal Y_H\longrightarrow H,\qquad
 E=h^*R\pi_!\underline\F_{\mathcal Y}
      \simeq R\pi_{H,!}\underline\F_{\mathcal Y_H},\\
 V_0=R\Gamma_c(\mathcal Y_0;\F),\qquad
 V=R\Gamma_c(U;\F).
 \end{gathered}
\]
These are derived functors on the analytic spaces. Proper-support
base change and $\vartheta$ identify
$i_H^*E\simeq V_0$ and $j_H^*E\simeq\underline V_L$.
Evaluation at zero and on the trivialized positive ray gives
isomorphisms in $D(\F)$,
\[
 \varepsilon_0:R\Gamma(H,E)\xrightarrow{\sim}V_0,
 \qquad
 \varepsilon_L:R\Gamma(L,j_H^*E)\xrightarrow{\sim}V.
\]
Define
\begin{equation}\label{eq:specialization-restriction}
 r_{\mathcal Y}:=
 \varepsilon_L\circ\operatorname{res}_{H,L}\circ\varepsilon_0^{-1}
 :V_0\longrightarrow V.
\end{equation}
Equivalently, $r_{\mathcal Y}$ is
$i_H^*E\to i_H^*Rj_{H,*}j_H^*E$ under the above
identifications. Under Verdier duality for the finite-type fibers,
\begin{equation}\label{eq:specialization-restriction-dual}
 \operatorname{Sp}_{\mathcal Y}=r_{\mathcal Y}^{\vee}:
 R\Gamma(U,\omega_U)\longrightarrow
 R\Gamma(\mathcal Y_0,\omega_{\mathcal Y_0}).
\end{equation}
This equality is natural for the proper and smooth comparisons
of Lemma~\ref{lp:lem:bm-collision-specialization}.
\end{lemma}
\begin{proof}
The triangle
$j_{H,!}\underline V_L\to\underline V_H\to i_{H,*}V\xrightarrow{+1}$
gives $R\Gamma(H,j_{H,!}\underline V_L)=0$.
Applying derived sections to
$j_{H,!}j_H^*E\to E\to i_{H,*}i_H^*E\xrightarrow{+1}$
therefore proves the assertion about $\varepsilon_0$.
Contractibility of $L$ proves the assertion about $\varepsilon_L$.
Naturality of evaluation identifies \eqref{eq:specialization-restriction}
with the stated adjunction morphism.

Compact-support localization gives
\begin{equation}\label{eq:half-ray-compact-localization}
 R\Gamma_c(L,\underline V_L)\longrightarrow R\Gamma_c(H,E)
 \longrightarrow V_0
 \xrightarrow{\partial_c}R\Gamma_c(L,\underline V_L)[1].
\end{equation}
For the constant complex $\underline V_H$, the middle term is
zero, since $R\Gamma_c(H;\F)=0$. Its boundary is consequently
an isomorphism
\[
 \beta_V:V\xrightarrow{\sim}R\Gamma_c(L,\underline V_L)[1].
\]
Apply naturality of localization to the adjunction morphism
$E\to Rj_{H,*}j_H^*E\simeq\underline V_H$.
Its map at zero is $r_{\mathcal Y}$ and its map on $L$ is
the identity. Hence
\begin{equation}\label{eq:restriction-compact-boundary}
 \partial_c=\beta_V\circ r_{\mathcal Y}.
\end{equation}
The inverse of $\beta_{\F}$ is the compact-support orientation
identification for $L$ oriented toward zero: the dual fundamental
ray has boundary $[0]$. Tensoring this identification with $V$,
with the ray as the first factor, fixes the same sign in every degree.

Let $F=\mathbb D_HE$. Verdier duality and proper-support base
change identify
\[
 F\simeq R\pi_{H,*}\omega_{\mathcal Y_H},\qquad
 i_H^!F\simeq V_0^\vee,\qquad
 j_H^*F\simeq\omega_L\otimes V^\vee.
\]
Dualizing \eqref{eq:half-ray-compact-localization} gives the
localization triangle for the dualizing complex of
$\mathcal Y_H$, after derived global sections.
By \eqref{eq:restriction-compact-boundary}, its shifted boundary,
preceded by the oriented source identification, is
\[
 V^\vee
 \xrightarrow{\sim}
 R\Gamma(L\times U,\omega_{L\times U})[-1]
 \xrightarrow{R\Gamma(\partial_H[-1])}
 R\Gamma(\mathcal Y_0,\omega_{\mathcal Y_0})
 =V_0^\vee,
\]
and equals $r_{\mathcal Y}^\vee$.

By \eqref{eq:specialization-half-ray-boundary}, this boundary
is also the derived-global-section map of
$\operatorname{Sp}^{\omega}_{\mathcal Y}$, proving
\eqref{eq:specialization-restriction-dual}.
All comparisons used above are natural adjunction, base-change,
and localization maps, proving the last assertion.
\end{proof}

\begin{lemma}[Specialization on Borel--Moore homology]
\label{lp:lem:bm-specialization-formula}
Assume that the product trivialization $\vartheta$ of
Lemma~\ref{lp:lem:bm-collision-specialization} is given. Taking $H^{-d}$ of $\operatorname{Sp}_{\mathcal Y}$ gives
the degree-preserving homomorphism
\begin{equation}\label{lp:eq:bm-sp}
 \operatorname{sp}_{\mathcal Y}:
 H_d^{\mathrm{BM}}(U;\F)\longrightarrow
 H_d^{\mathrm{BM}}(\mathcal Y_0;\F),
 \qquad \xi\longmapsto\partial\,\vartheta_*(\rho\boxtimes\xi).
\end{equation}
Here $\partial$ is the connecting homomorphism in
\[
 H_{d+1}^{\mathrm{BM}}(\mathcal Y;\F)\longrightarrow
 H_{d+1}^{\mathrm{BM}}(\mathcal Y^\times;\F)
 \xrightarrow{\partial}H_d^{\mathrm{BM}}(\mathcal Y_0;\F)
 \longrightarrow H_d^{\mathrm{BM}}(\mathcal Y;\F).
\]
The ray is placed first in $\Gm\times U$. We suppress
$\vartheta_*$ when using the specified trivialization.
For the maps $f,u$ of Lemma~\ref{lp:lem:bm-collision-specialization},
the induced proper pushforwards and smooth Gysin maps satisfy
\[
 \begin{aligned}
 f_{0,*}\operatorname{sp}_{\mathcal Y'}
     &=\operatorname{sp}_{\mathcal Y}u_* &&(f\text{ proper}),\\
 \operatorname{sp}_{\mathcal Y'}u^!
     &=f_0^!\operatorname{sp}_{\mathcal Y}
       &&(f\text{ smooth of relative dimension }e).
 \end{aligned}
\]
Both Gysin maps raise Borel--Moore degree by $2e$.
For the product family $\mathcal Y=\A^1\times U$ with its
product trivialization, specialization is the identity under
$\mathcal Y_0\cong U$.
\end{lemma}
\begin{proof}
Taking derived global sections of the localization
triangle and using
$H^{-d}R\Gamma(V,\omega_V)=H_d^{\mathrm{BM}}(V;\F)$
gives the displayed exact sequence. Under the source
identification of Lemma~\ref{lp:lem:bm-collision-specialization}, the trace for $\ell$ sends $\xi$ to
$\vartheta_*(\rho\boxtimes\xi)$. The shifted localization
boundary then gives \eqref{lp:eq:bm-sp}. The ray raises the
Borel--Moore degree by one and the boundary lowers it by one.
The proper and smooth identities follow from that lemma; the
smooth pullback induced by purity is $R\Gamma(V,\omega_V)\to
R\Gamma(V',\omega_{V'})[-2e]$, giving the degree shift $2e$.
The product-family assertion follows from the identity of
complexes in that lemma. Placing the ray last would instead
introduce $(-1)^d$ in the displayed Borel--Moore formula.
\end{proof}
\begin{lemma}[Cohomology actions and specialization]
\label{lem:specialization-module}
For a topological space $N$, use the sheaf category
$\mathcal D(N;\F)$ of Subsection~\ref{subsec:coherent-descent} and write
\[
 H^d(N;\F):=H^dR\Gamma(N,\underline\F_N)
   =\operatorname{Hom}_{D(N;\F)}
       (\underline\F_N,\underline\F_N[d]).
\]
Choose a morphism $c:\underline\F_N\to\underline\F_N[d]$
representing such a class. For $E\in\mathcal D(N;\F)$, define
\[
 \begin{aligned}
 m_{c,E}:E\simeq\underline\F_N\otimes E
  &\xrightarrow{\ c\otimes\mathrm{id}_E\ }
  \underline\F_N[d]\otimes E\simeq E[d],\\
 c\frown\xi&:=H^qR\Gamma(N,m_{c,E})(\xi)
 \quad\bigl(\xi\in H^qR\Gamma(N,E)\bigr).
 \end{aligned}
\]
Thus the cap action has cohomological degree $d$.
\begin{enumerate}[(i),wide=0pt]
\item \emph{Naturality.} For every morphism
$\tau:K\to T$ in $\mathcal D(N;\F)$, tensor functoriality gives
a natural homotopy
\begin{equation}\label{eq:specialization-module-arrow}
 \tau[d]\circ m_{c,K}\simeq m_{c,T}\circ\tau
       \qquad\text{in }\mathcal D(N;\F).
\end{equation}
It gives an equality after passage to $D(N;\F)$ or to cohomology.
The action uses the actual shift of its object: under the usual
shift identifications,
\[
 m_{c,E[r]}=(-1)^{dr}(m_{c,E})[r].
\]

\item \emph{Finite-dimensional specialization.}
Let $\pi:Y\to\A^1_s$ be a morphism of separated complex schemes
of finite type, as in Lemma~\ref{lp:lem:bm-collision-specialization}.
On analytic spaces set
\[
 H=[0,\infty),\quad L=(0,\infty),\quad
 Y_H=\pi^{-1}(H),\quad Y_L=\pi^{-1}(L),\quad Y_0=\pi^{-1}(0).
\]
Let $v:N\hookrightarrow Y_H$ be any open neighborhood of $Y_0$,
put $N^+=N\cap Y_L$, and denote the inclusions by
$i_N:Y_0\hookrightarrow N$ and $j_N:N^+\hookrightarrow N$.
All dualizing complexes in this part are ordinary dualizing
complexes, with no volume-index normalization. Define
\[
 \mathcal K_{N,\rho}=Rj_{N,*}\omega_{N^+}[-1],\qquad
 T_N=i_{N,*}\omega_{Y_0},\qquad
 \sigma_N:=\partial_N[-1]:\mathcal K_{N,\rho}\longrightarrow T_N,
\]
where $\partial_N$ is the boundary of
\[
 i_{N,*}\omega_{Y_0}\longrightarrow\omega_N
 \longrightarrow Rj_{N,*}\omega_{N^+}
 \xrightarrow{\partial_N}i_{N,*}\omega_{Y_0}[1].
\]
For $c\in H^d(N;\F)$, part~(i) applies to $\tau=\sigma_N$.
More explicitly, the projection formula identifies its two actions as
\[
 \begin{aligned}
 m_{c,\mathcal K_{N,\rho}}
   &=(-1)^d\bigl(Rj_{N,*}m_{j_N^{-1}c,\omega_{N^+}}\bigr)[-1],\\
 m_{c,T_N}&=i_{N,*}m_{i_N^{-1}c,\omega_{Y_0}}.
 \end{aligned}
\]
With $\operatorname{sp}_N=H^*R\Gamma(N,\sigma_N)$, this gives
\begin{equation}\label{eq:specialization-module}
 \operatorname{sp}_N(c\frown\xi)
   =(i_N^{-1}c)\frown\operatorname{sp}_N(\xi),
 \qquad \xi\in H^qR\Gamma(N,\mathcal K_{N,\rho}).
\end{equation}
Both sides have degree $q+d$. The left action is on $\mathcal K_{N,\rho}$,
including its $[-1]$ shift; it is not the unshifted action on
$Rj_{N,*}\omega_{N^+}$.

For $a:Y_H\hookrightarrow Y^{\mathrm{an}}$, one has
$a_*\sigma_{Y_H}=\operatorname{Sp}^{\omega}_Y$ of
Lemma~\ref{lp:lem:bm-collision-specialization}. Open restriction gives
$v^{-1}\sigma_{Y_H}\simeq\sigma_N$ and hence
\[
 R\Gamma(N,\sigma_N)\circ\operatorname{res}_v
   \simeq R\Gamma(Y_H,\sigma_{Y_H}),\qquad
 \operatorname{res}_v:R\Gamma(Y_H,\mathcal K_{Y_H,\rho})\longrightarrow R\Gamma(N,\mathcal K_{N,\rho}),
\]
with common target $R\Gamma(Y_0,\omega_{Y_0})$.
No equivalence of these source complexes is required.

\item \emph{Product interpretation.}
Suppose $N=Y_H$ and an isomorphism
$\vartheta:\Gm\times U\xrightarrow{\sim}Y^\times$ over $\Gm$
is given, with $U$ separated of finite type. Let
$\alpha:R\Gamma(U,\omega_U)\xrightarrow{\sim}R\Gamma(Y_H,\mathcal K_{Y_H,\rho})$
be the inward-ray identification, with the ray factor first.
For $c\in H^d(Y_H;\F)$, let $c_U$ be its pullback to
$\{s_0\}\times U$, for any $s_0>0$, and put $c_0=c|_{Y_0}$.
Then
\[
 \alpha[d]\circ R\Gamma(U,m_{c_U,\omega_U})
    \simeq R\Gamma(Y_H,m_{c,\mathcal K_{Y_H,\rho}})\circ\alpha,
 \qquad
 (\operatorname{Sp}_Y)_*(c_U\frown\xi)
     =c_0\frown(\operatorname{Sp}_Y)_*(\xi).
\]
Here the first comparison uses the product homotopy on $L\times U$;
$\operatorname{Sp}_Y$ is \eqref{eqSpG}, and
$\xi\in H^qR\Gamma(U,\omega_U)=H^{\mathrm{BM}}_{-q}(U;\F)$.
The last equality lies in $H^{\mathrm{BM}}_{-q-d}(Y_0;\F)$.
\end{enumerate}
\end{lemma}
\begin{proof}
For (i), both composites in
\eqref{eq:specialization-module-arrow} are identified with
\[
 \underline\F_N\otimes K
    \xrightarrow{\ c\otimes\tau\ }
 \underline\F_N[d]\otimes T\simeq T[d].
\]
The bifunctor $\otimes$ supplies this homotopy and its composition
compatibilities. Moving the degree-$d$ action past the shift $[r]$
gives $(-1)^{dr}$, proving the shift formula. This argument applies
on arbitrary $N$; it uses no finite-dimensional duality.

For (ii), the open and closed embeddings satisfy
$j_N^{-1}\omega_N\simeq\omega_{N^+}$ and
$i_N^!\omega_N\simeq\omega_{Y_0}$.
Their localization triangle defines $\sigma_N$.
For an inclusion $g$ and a sheaf complex $E$ on its source,
the projection formula for the invertible object
$\underline\F_N[d]$ identifies
\[
 m_{c,Rg_*E}\simeq Rg_*m_{g^{-1}c,E}.
\]
Together with (i) for $r=-1$, this proves the two displayed
action formulas. Apply (i) to $\sigma_N$ and take cohomology
to obtain \eqref{eq:specialization-module}.
The identity $a_*\sigma_{Y_H}=\operatorname{Sp}^{\omega}_Y$ is
\eqref{eq:specialization-half-ray-boundary}.
Open base change for $Rj_{N,*}$ and localization identifies its
restriction to $N$ with $\sigma_N$. Applying $R\Gamma$ gives the
restriction identity.

For (iii), contraction of $L$ to $s_0$ gives
$\vartheta_L^{-1}(c|_{Y_L})=\operatorname{pr}_U^{-1}c_U$
in cohomology, and a homotopy between representatives.
Let $\rho=[L]\in H^{\mathrm{BM}}_1(L;\F)$ have the inward
orientation. In raw Borel--Moore grading the localization boundary
and exterior product obey
\[
 \begin{aligned}
 \partial\bigl((c|_{Y_L})\frown\eta\bigr)
     &=(-1)^d c_0\frown\partial\eta,\\
 (\operatorname{pr}_U^{-1}c_U)\frown(\rho\boxtimes\xi)
     &=(-1)^d\rho\boxtimes(c_U\frown\xi).
 \end{aligned}
\]
The first sign is the shift formula in (i) applied to
$\sigma_{Y_H}=\partial[-1]$; the second is the Koszul sign
from the ray class of cohomological degree $-1$.
They cancel in
$(\operatorname{Sp}_Y)_*\xi=\partial(\rho\boxtimes\xi)$.
Equivalently, shifting the ray factor by $[-1]$ makes
$\alpha$ intertwine the displayed actions. This proves (iii).
\end{proof}

\subsection{Restriction on collars and coordinate open stars}\label{subsec:collar-restrictions}
We construct restriction maps on parameter collars and
between coordinate open stars, generalizing
\eqref{eq:specialization-restriction}, and prove their
compatibility with coordinate-face inclusions and
exterior products.

For a topological space $T$ and a real number $a>0$,
a \emph{parameter collar of $T$} is the product
$[0,a)\times T$, with distinguished \emph{zero face}
$\{0\}\times T$ and \emph{positive part} $(0,a)\times T$.
We also allow the truncated model $[0,a]\times T$ and
the half-line model $[0,\infty)\times T$, with the same
zero face and with positive part defined by requiring
the first coordinate to be positive.
In the restriction constructions below, we abbreviate
\emph{parameter collar} to \emph{collar}.
The product structure concerns the parameter space;
the identifications of sheaf complexes on its zero face
and positive part are specified separately.

For an integer $r\ge1$, put $[r]=\{1,\ldots,r\}$ and
$Q=[0,\infty)^r$. For $A\subseteq[r]$, define the
\emph{coordinate stratum} $Q_A$ and its \emph{open star}
$W_A$ by
\[
 \begin{aligned}
 Q_A&=\{x\in Q:x_j>0\Longleftrightarrow j\in A\},\\
 W_A&=\{x\in Q:x_j>0\text{ for every }j\in A\}
      =\bigcup_{A\subseteq B\subseteq[r]}Q_B.
 \end{aligned}
\]
Thus $W_A$ is open in the subspace topology on
$Q\subseteq\mathbb R^r$ and is the union of precisely
those coordinate strata whose closures in $Q$
contain $Q_A$. For the bounded cube $[0,1]^r$, the corresponding
open stars are $W_A\cap[0,1]^r$.

\begin{lemma}[Collar and open-star restriction maps]
\phantomsection\label{lem:parameter-chain-properties}
\leavevmode
\begin{enumerate}[(i),wide=0pt]
\item Let $T$ be a paracompact finite-dimensional smooth real manifold,
possibly with corners, or a locally finite finite-dimensional
polyhedron. Put $M=[0,1]\times T$, with zero-end inclusion $i:T\to M$
and complementary open inclusion $j:(0,1]\times T\to M$.
Let $C,D$ be perfect complexes of $\F$-vector spaces.
Suppose a bounded constructible complex $E$ on $M$ is supplied with
identifications
\[
 i^*E\simeq\underline{D^\vee}_T,\qquad
 j^*E\simeq\operatorname{pr}_T^*\underline{C^\vee}_T.
\]
There is a \emph{collar specialization}
\begin{equation}\label{eq:collar-specialization}
 \begin{aligned}
 \alpha_E:\omega_T\otimes C
 &\simeq\mathbb D_T(i^*Rj_*j^*E)\\
 &\xrightarrow{\mathbb D_T(i^*\eta_E)}
       \mathbb D_T(i^*E)
       \simeq\omega_T\otimes D,
 \end{aligned}
\end{equation}
where $\eta_E:E\to Rj_*j^*E$ is the adjunction unit.
We call
\[
 i^*\eta_E:i^*E\longrightarrow i^*Rj_*j^*E
\]
the \emph{collar restriction morphism}. Its source is identified
with $\underline{D^\vee}_T$ by hypothesis, and the positive-part
trivialization identifies its target with
$\underline{C^\vee}_T$. Its Verdier dual is the collar
specialization $\alpha_E$.
These constructions are natural for morphisms of $E$ and the
specified identifications.

\item  Let $Q$, $Q_A$, and $W_A$ be as defined above.
For each $A\subseteq[r]$, choose the point $v_A\in Q_A$
given by
\[
 (v_A)_j=
 \begin{cases}
 1,&j\in A,\\
 0,&j\notin A.
 \end{cases}
\]
Let $E$ be a bounded complex whose restriction to each $Q_A$
is locally constant. Evaluation at $v_A$ induces a natural equivalence
\begin{equation}\label{eq:orthant-star-stalk}
 \epsilon_A:G_E(A):=R\Gamma(W_A,E)\xrightarrow{\sim}E_{v_A},
 \qquad E_{v_A}:=i_{v_A}^*E.
\end{equation}
We call $\epsilon_A$ the \emph{open-star evaluation equivalence};
it evaluates sections on $W_A$ at $v_A$.

\item Let $\mathcal P([r])$ be the poset of subsets of $[r]$,
ordered by inclusion. For $A\subseteq A'$, the inclusion
$W_{A'}\subseteq W_A$ induces the \emph{open-star restriction map}
\[
 r_{A,A'}:
 G_E(A)=R\Gamma(W_A,E)
 \longrightarrow
 G_E(A')=R\Gamma(W_{A'},E).
\]
The arrow is restriction of derived sections along that open
inclusion. These maps form a functor
$G_E:N\mathcal P([r])\to\mathcal D(\F)$.
In particular, a chain $A_0\subseteq\cdots\subseteq A_n$
gives the simplex
$G_E|_{\Delta^n}:\Delta^n\to\mathcal D(\F)$, including its
composition homotopies, faces, and degeneracies.
If the stalks $E_{v_A}$ are perfect, linear duality gives
\[
 G_E^\vee:N\mathcal P([r])^{\mathrm{op}}\longrightarrow\mathcal D(\F),
 \qquad
 r_{A,A'}^\vee:G_E(A')^\vee\longrightarrow G_E(A)^\vee.
\]
For $r=1$, $r_{\varnothing,\{1\}}^\vee$ agrees with (i) for
$T=\mathrm{pt}$, using $C=E_1^\vee$ and $D=E_0^\vee$.

For $J\subseteq[r]$, let
$\iota_J:[0,\infty)^J\hookrightarrow Q$ set the other coordinates
equal to zero. For $A\subseteq J$, put
\[
 W_A^J:=\iota_J^{-1}(W_A)
       =\{x\in[0,\infty)^J:x_j>0
                           \text{ for every }j\in A\}.
\]
Restriction to the coordinate face induces the comparison
\[
 R\Gamma(W_A,E)
 \longrightarrow
 R\Gamma(W_A^J,\iota_J^*E).
\]
These maps assemble into an equivalence of diagrams
\[
 \left.G_E\right|_{N\mathcal P(J)}
       \xrightarrow{\sim}G_{\iota_J^*E}.
\]
They compare sections of $E$ with sections of its inverse image
on the coordinate face; the open-star restriction maps are
the arrows within these diagrams.
If $E'$ is bounded and face-constructible on $Q'=[0,\infty)^s$, the
exterior-product maps give an equivalence of diagrams
\[
 G_E(A)\otimes G_{E'}(B)
       \xrightarrow{\sim}G_{E\boxtimes E'}(A\sqcup B),
 \qquad A\subseteq[r],\ B\subseteq[s],
\]
where the two coordinate sets are regarded as disjoint.
These equivalences are natural in $E,E'$ and compatible with
successive face inclusions and iterated exterior products.
For $E,E'$ with perfect stalks, their linear duals have the
corresponding compatibilities.

\end{enumerate}
In the two-periodic version, replace bounded sheaf complexes by
two-periodic sheaf complexes with the same local-constancy and
trivialization conditions. Require perfect periodic complexes
in (i) and perfect stalks for the dual diagrams in (iii).
\end{lemma}
\begin{proof}
For (i), the positive-collar trivialization induces
$\underline{C^\vee}_T\to i^*Rj_*j^*E$.
For $t\in T$, take contractible neighborhoods $U$ of $t$ and
$0<\varepsilon<1$. The induced map on stalks is the equivalence
\[
 C^\vee\xrightarrow{\sim}
 \underset{\varepsilon,U}{\operatorname{colim}}
 R\Gamma\bigl((0,\varepsilon)\times U,
                  \underline{C^\vee}\bigr).
\]
Thus the displayed sheaf map is an equivalence. Dualizing
$i^*\eta_E$ gives \eqref{eq:collar-specialization}, since
$\mathbb D_T\underline{C^\vee}_T\simeq\omega_T\otimes C$
by perfectness. Naturality follows from the adjunction unit
and duality. This uses $\omega_T$ itself; on a space with
boundary it need not be a shifted constant sheaf.

For (ii), the homeomorphism
\[
 Q\longrightarrow\{t\in\Delta^r:t_0>0\},
 \qquad x\longmapsto
       \frac{(1,x_1,\ldots,x_r)}{1+\sum_jx_j}
\]
identifies $Q_B$ with the open simplex indexed by $\{0\}\cup B$.
The exit-path $\infty$-category of a simplicial complex is
equivalent to the nerve of its face poset
\cite[Theorem~A.6.10]{HA}. Its restriction to $W_A$ is therefore
\[
 \operatorname{Exit}(W_A)\simeq N\mathcal P_A,
 \qquad \mathcal P_A=\{B\subseteq[r]:A\subseteq B\}.
\]
An exit simplex has an increasing chain of strata. If its first
stratum belongs to this upper set, every stratum does; hence
restriction to their union gives the full exit subcategory on
these objects.
Here exit paths run from a face to a containing face, and the
exit category includes their homotopies of all orders.
The spaces are paracompact, locally contractible, and conically
stratified by a finite poset. The constructible-sheaf
classification \cite[Theorem~A.9.3]{HA} consequently applies.
It represents $E$ by a diagram $F_E:N\mathcal P([r])\to\mathcal D(\F)$
whose value at $B$ is $E_{v_B}$, and gives
\[
 R\Gamma(W_A,E)\simeq
       \operatorname*{holim}_{B\in\mathcal P_A}F_E(B)
       \xrightarrow{\sim}F_E(A)=E_{v_A}.
\]
The last arrow is evaluation at the initial object $A$ of
$\mathcal P_A$; under the classification it is precisely
restriction of sections to the stalk. This proves (ii).
The statement for complexes follows from the classification
for sheaves of spaces by stabilization and passage to
$\F$-module objects.

For (iii), the inclusions $W_{A'}\subseteq W_A$ and derived
sections define $G_E$ as a functor; no choices of maps in the
ordinary derived category are involved. Under the preceding
limit description, the restriction $r_{A,A'}$ is the map
$F_E(A)\to F_E(A')$. Hence the equivalences $\epsilon_A$
identify the entire diagrams, not merely their objects.
For $J\subseteq[r]$ and $A\subseteq J$, let $v_A^J$ denote
the point of $[0,\infty)^J$ with coordinates $1$ on $A$ and
$0$ on $J\setminus A$. With $W_A^J$ as in the statement,
$\iota_J^{-1}(W_A)=W_A^J$ and $\iota_J(v_A^J)=v_A$.
The diagram
\[
 \begin{array}{ccc}
 R\Gamma(W_A,E)&\longrightarrow&R\Gamma(W_A^J,\iota_J^*E)\\
 \big\downarrow\rlap{$\scriptstyle\epsilon_A$}&&
 \big\downarrow\rlap{$\scriptstyle\epsilon_A$}\\
 E_{v_A}&=&(\iota_J^*E)_{v_A^J}
 \end{array}
\]
commutes by naturality of restriction to a stalk. Both vertical
arrows are equivalences by (ii), proving the face assertion.
Likewise, evaluation identifies the exterior-product map with
\[
 E_{v_A}\otimes E'_{v_B}
       \xrightarrow{\sim}(E\boxtimes E')_{(v_A,v_B)}.
\]
The naturality of inverse image, exterior product, and sections
supplies these comparisons as transformations of diagrams,
compatible with their composites. Linear duality is a functor
on the opposite category; for perfect objects it also carries
tensor products to tensor products as in
Subsection~\ref{app:coherent-diagrams}.

For $r=1$, evaluation identifies
$R\Gamma(Q,E)\to R\Gamma(W_{\{1\}},E)$ with
$i_0^*\eta_E:E_0\to i_0^*Rj_*j^*E\simeq E_1$.
Its dual is the map in (i); restricting to $[0,1]$ makes the
collar identification explicit. Finally, the same argument
uses module objects over $\F[u,u^{-1}]$, $|u|=2$, in the
periodic case; it requires no integer-graded lift.
\end{proof}
\subsection{Finite-factor presentations over a base and relative volume indices}
\label{subsec:relative-volume-preparations}
Let $S$ be a quasi-compact separated locally nice complex scheme,
let $T$ be a separated complex scheme of finite type, and let
$f:S\to T$ be a morphism. A \emph{compatible finite-factor presentation
over $T$} is a compatible finite-factor presentation as defined in
\emph{Compatible finite-factor presentations} in
Subsection~\ref{subsec:coherent-descent}, with a
factorization on each chart $u:U=Y\times\A^J\to S$,
\[
 f u=g p,\qquad p:U\to Y,\qquad g:Y\to T,
\]
and with all finite-factor comparison maps over $T$. This is
what it means for a finite factor to \emph{retain the map to the
base}; it is the displayed factorization even when $T$ is not
affine. Specified finitely presented closed subfamilies are
presented by inverse images of closed subschemes of $Y$, as in
\eqref{eq:presented-closed-diagram}, with their open complements.
According to the application we have in mind, we will also need to choose from start a collection of base changes $T'\to T$. For each specified base change $T'\to T$ between separated
finite-type complex schemes, include the induced finite-factor presentation
\[
 U\times_TT'\cong(Y\times_TT')\times\A^J.
\]
If $T'\hookrightarrow T$ is closed, include its inverse-image
closed inclusion and open complement in the finite-factor presentation diagram.
The finite-type factors $Y\times_TT'$ and open complements in
$Y$ need not be affine. In that case take finite affine open
covers and use their products with $\A^J$. These are finitely
presented open refinements with the same omitted coordinates;
further affine refinements compare the covers and overlaps.

A \emph{common refinement over $T$} is a common refinement of
presented atlases as in \emph{Comparing different atlases} in
Subsection~\ref{subsec:coherent-descent}, whose comparison maps
commute with the maps to $T$, with the
Cartesian finite-factor presentations of the specified closed subfamilies,
and with the induced finite-factor presentations of the specified base changes.
The finite-factor presentations of the open complements are their restrictions.
For an infinite diagram, these conditions concern each finite
subdiagram and each finite composable string of chart maps;
factors are chosen for chart occurrences and compared by
further common refinements over the corresponding bases.

\begin{definition}[Relative volume indices]
\label{def:relative-volume-index}
\label{conv:relative-volume-index}
Let $S\to T$ have a specified compatible finite-factor presentation
as above, including its specified overlaps, closed subfamilies,
and base changes. A \emph{relative volume-index normalization over $T$} is a collection of locally constant functions
\[
                 \delta_{Y/T}:Y\longrightarrow\Z,
\]
one for each finite-factor projection $\pi:U\to Y\to T$ in
that finite-factor presentation, satisfying the following conditions.
\begin{enumerate}[(i),wide=0pt]
\item For a specified smooth refinement $p:Y'\to Y$ over $T$,
of relative complex dimension $e$, require
\[
                  \delta_{Y'/T}=p^*\delta_{Y/T}+e.
\]
The assignments agree on overlaps through this rule. Explicitly,
if an overlap has a common finite factor $W$ with smooth
comparison maps $a:W\to Y$ and $b:W\to Y'$, then
\[
 a^*\delta_{Y/T}+\dim(a)=b^*\delta_{Y'/T}+\dim(b).
\]
All these comparisons are required to respect composition.
\item For a specified base change $T'\to T$, let
$g:Y\times_TT'\to Y$ be the projection and require
\[
              \delta_{(Y\times_TT')/T'}=g^*\delta_{Y/T}.
\]
\item For every closed inclusion $i:Z\hookrightarrow Y$ included
in the specified finite-factor presentation diagram over $T$, require
\[
 \delta_{Z/T}=i^*\delta_{Y/T}.
\]
Here, if the ambient finite-factor projection is
$Y\times\A^J\to Y$, the induced projection on the closed chart
is $Z\times\A^J\to Z$, with the same omitted coordinates.
\end{enumerate}
Each function $\delta_{Y/T}$ is called the \emph{relative volume
index} of that finite-factor presentation. The notation suppresses
$\pi$, the whole family diagram, and its chosen relative volume-index
normalization;
it is not an invariant of the morphism $Y\to T$ alone. 

A further finitely presented closed subfamily may be included
in the diagram, together with its open complement, after enlarging
the finite factors to retain its defining equations. The new closed
factors carry the restricted indices of (iii). For each such factor
$i:Z\hookrightarrow Y$, its open complement
$j:V=Y\setminus Z\hookrightarrow Y$ carries
$\delta_{V/T}:=j^*\delta_{Y/T}$.
We call these their \emph{induced normalizations}; they introduce
no additional choice. Indeed, for a smooth comparison $p:Y'\to Y$
of relative dimension $e$, put $Z'=Y'\times_Y Z$ and write
$i':Z'\hookrightarrow Y'$ and $p_Z:Z'\to Z$ for the induced maps.
The morphism $p_Z$ is smooth of relative dimension $e$, and
\[
 \delta_{Z'/T}=i'^*\delta_{Y'/T}
             =p_Z^*\delta_{Z/T}+e.
\]
The overlap identities in (i) restrict to the corresponding
closed and open factors. Compatibility with the specified base
changes and successive restrictions follows from functoriality
of pullback. Thus these assignments extend the chosen normalization.

For two presented families $S_i\to T_i$ with chosen relative
volume-index normalizations, define the \emph{induced relative
volume-index normalization} on the product of their finite-factor presentations by
\[
 \delta_{Y_1\times Y_2/T_1\times T_2}
   :=\operatorname{pr}_1^*\delta_{Y_1/T_1}
     +\operatorname{pr}_2^*\delta_{Y_2/T_2}.
\]
Define the induced relative volume-index normalization on the
affine stabilization $S\times\A^e\to T$ with $e\in\Z_{\ge0}$ by
\[
 \delta_{Y\times\A^e/T}
       :=\operatorname{pr}_Y^*\delta_{Y/T}+e.
\]
For a finite or countable affine tail $\A^J$, apply this definition
to each finite subset $J_0\subset J$ of coordinates included in
the factor, with $e=|J_0|$.
Even when $J$ is infinite, every relative volume index is attached
to a finite factor and is integer-valued.
If one merely retains further coordinates in a finite factor
of the same chart, this formula already follows from (i).
The induced assignments satisfy (i)--(iii): relative dimensions
add under products, and sums of volume indices commute with
pullback. These formulas define normalizations on the new product
and stabilized finite-factor presentations; they impose no additional condition
on the original family diagram.
\end{definition}
\subsection{Normalized dualizing complexes on locally nice families}
\label{subsec:relative-compact-support}
Let $f:S\to T$ and its compatible finite-factor presentations be as in
Subsection~\ref{subsec:relative-volume-preparations}. Choose a finite
nice \'etale atlas $W_0=\coprod_iU_i\to S$. Its \v Cech nerve is
$W_n=W_0^{\times_S(n+1)}$. Lemma~\ref{lem:nice-cech-overlaps}
expresses each $W_n$ as a finite disjoint union of nice overlap
schemes; use their finite factors directly.
Use compatible finite-factor presentations over $T$ for finite composable
strings of maps between these overlap charts, with a factor for
each chart occurrence, as in
Subsection~\ref{subsec:relative-volume-preparations}. Common
refinements compare different strings, truncations, and atlases.

\begin{definition}[Normalized dualizing complex and its relative images]
\label{def:relative-normalized-dual-complexes}
Give the finite-factor presentations of $S/T$ compatible relative volume
indices as in Definition~\ref{def:relative-volume-index}. No base changes  $T'\to T$ are specified in this definition.
On a chart
$p:U=Y\times\A^J\to Y$, with $f|_U=g p$, put
\[
 \omega_{Y/T}^{\mathrm n}:=\omega_Y[-2\delta_{Y/T}],\qquad
 \omega_{U/T}^{\mathrm n}:=p^{-1}\omega_{Y/T}^{\mathrm n}.
\]
For the smooth finite-factor comparison maps, use the construction
of \eqref{eq:enhanced-normalized-purity} with the relative indices
and the comparison rule of Definition~\ref{def:relative-volume-index}(i).
We call the resulting equivalences the \emph{relative normalized
purity equivalences}. Their ordinary inverse images on charts,
with their composition homotopies, glue these chart complexes. Recall that $p^{-1}$ here is ordinary
sheaf inverse image for
$p:U^{\mathrm{an}}=Y^{\mathrm{an}}\times\C^J\to Y^{\mathrm{an}}$;
it transports both the finite-factor complexes and their maps.
On these chart spaces we use $p^{-1}$, $Rp_*$, and, for a closed
embedding $i$, the adjunction $i_*\dashv i^!$;
no exceptional inverse image $p^!$ is required.
Lemma~\ref{lem:infinite-product-descent} supplies the product
comparisons and descent. The resulting object on $S^{\mathrm{an}}$ is the
\emph{normalized dualizing complex}
$\omega^{\mathrm n}_{S/T}$.
The subscripts $Y/T$, $U/T$, and $S/T$ record the chosen
relative volume-index normalization; they do not indicate relative
Verdier duality. Define the \emph{relative normalized Borel--Moore complex}
by
\[
 \mathcal B_{S/T}:=Rf_*\omega^{\mathrm n}_{S/T}.
\]
Here $Rf_*$ is the derived direct image for the continuous map
$S^{\mathrm{an}}\to T^{\mathrm{an}}$.

There is also a compact-support companion, constructed on the
finite-dimensional base. Use the simplicial indexing category
$\Delta^{\mathrm{op}}$ of Definition~\ref{def:cech-nerve}.
For finite factors
$W_{n,\beta}=Y_{n,\beta}\times\A^{J_{n,\beta}}$ over
$g_{n,\beta}:Y_{n,\beta}\to T$, set
\begin{equation}\label{eq:relative-compact-factor}
 E_n=\bigoplus_\beta
 Rg_{n,\beta!}\underline\F_{Y_{n,\beta}}
                 [2\delta_{Y_{n,\beta}/T}],
 \qquad
 \mathcal E_{S/T}:=
 \operatorname*{hocolim}_{[n]\in\Delta^{\mathrm{op}}}E_n.
\end{equation}
The homotopy colimit is taken in
$\mathcal D(T^{\mathrm{an}};\F)$. The maps in this simplicial
diagram are normalized traces.
Explicitly, a chart map represented after refinement by a smooth
$v:Y'\to Y$ of relative dimension $e$ induces
\[
 R(gv)_!\underline\F_{Y'}[2\delta_{Y'/T}]
 \simeq Rg_!Rv_!v^!
       \bigl(\underline\F_Y[2\delta_{Y/T}]\bigr)
 \xrightarrow{\;Rg_!\operatorname{tr}_v\;}
 Rg_!\underline\F_Y[2\delta_{Y/T}].
\]
For the degeneracies, use the traces of the diagonal
open-and-closed immersions of the separated \'etale charts.
The enhanced adjunction and base-change transformations of
Subsection~\ref{conv:section-five-sheaf-operations} on the
finite factors give the simplicial identities with their higher
compatibilities. Thus \eqref{eq:relative-compact-factor} is a
homotopy realization, as defined in
Appendix~\ref{app:totalization-realization}, and is not an
unnormalized compact-support image of the infinite-dimensional $S$.
Coherent comparisons of these simplicial diagrams induce maps
of their realizations with simultaneous cochain representatives
by Lemma~\ref{lem:totalization-realization-representatives}.

Independently of any relative volume-index normalization, define
the parity versions by using $\operatorname{Per}_2\omega_Y$ and
$\operatorname{Per}_2\underline\F_Y$ on the finite factors.
The even shifts in purity and the trace comparisons disappear
under the canonical periodization identifications. If relative
volume indices have been chosen, periodizing the normalized
finite-factor terms gives these same comparisons.

For $\alpha\in H^qR\Gamma(S,\omega^{\mathrm n}_{S/T})$,
we extend the terminology of
Definition~\ref{def:finite-factor-representatives} as follows.
On a chart $h:U\to S$ with finite-factor projection $p:U\to Y$,
its {\it representative} is the unique class
$\beta\in H^qR\Gamma(Y,\omega_Y[-2\delta_{Y/T}])$
such that $h^*\alpha=p^*\beta$ under the chart comparison.
The same convention applies in parity grading.
\end{definition}

\begin{proposition}[Descent and relative images]
\label{lem:relative-compact-support}
Fix $f:S\to T$ as above and, in integer grading, its chosen
relative volume-index normalization. The objects of
Definition~\ref{def:relative-normalized-dual-complexes} are well
defined in the derived $\infty$-categories
\[
 \omega^{\mathrm n}_{S/T}\in\mathcal D(S^{\mathrm{an}};\F),
 \qquad
 \mathcal E_{S/T},\mathcal B_{S/T}
       \in\mathcal D(T^{\mathrm{an}};\F).
\]
Use the two-periodic categories in parity grading.

The objects are independent of the nice \'etale atlas and its
compatible finite-factor presentations, with the relative volume
indices transported by Definition~\ref{def:relative-volume-index}.
Common refinements over $T$, as defined in
Subsection~\ref{subsec:relative-volume-preparations}, induce
comparison equivalences of these objects and of the maps in the
corresponding diagrams of complexes.
Independence does not include changing the relative
volume-index normalization.

There is a natural equivalence
\[
 \mathcal B_{S/T}\simeq\mathbb D_T\mathcal E_{S/T}.
\]
The following assertions hold.
\begin{enumerate}[(i),wide=0pt]
\item If $i:Z\hookrightarrow S$ is a finitely presented closed
inclusion and $j:V=S\setminus Z\hookrightarrow S$, include this
closed--open diagram in the finite-factor presentations, with the
induced normalization of Definition~\ref{def:relative-volume-index}.
There are coherent identifications
\[
 i^!\omega^{\mathrm n}_{S/T}\simeq\omega^{\mathrm n}_{Z/T},
 \qquad
 j^{-1}\omega^{\mathrm n}_{S/T}\simeq\omega^{\mathrm n}_{V/T},
\]
and a localization triangle on $S^{\mathrm{an}}$
\[
 i_*\omega^{\mathrm n}_{Z/T}\longrightarrow
 \omega^{\mathrm n}_{S/T}\longrightarrow
 Rj_*\omega^{\mathrm n}_{V/T}\xrightarrow{+1}.
\]
Its direct image is the localization triangle for $\mathcal B$.
The dual finite-factor localization triangles give
$\mathcal E_{V/T}\to\mathcal E_{S/T}\to
\mathcal E_{Z/T}\xrightarrow{+1}$.

\item For a morphism $b:T'\to T$ between separated complex
schemes of finite type, put $S'=S\times_TT'$ and include the
induced finite-factor presentations with the normalization of
Definition~\ref{def:relative-volume-index}(ii). Then
\begin{equation}\label{eq:relative-dual-base-change}
 b^*\mathcal E_{S/T}\simeq\mathcal E_{S'/T'},
 \qquad
 b^!\mathcal B_{S/T}\simeq\mathcal B_{S'/T'}.
\end{equation}
The extraordinary inverse image in this formula is on the
finite-dimensional bases. For $i_t:\{t\}\hookrightarrow T$,
the induced chart is
\begin{equation}\label{eq:induced-fiber-presentation}
 U_t=Y_t\times\A^J,
 \qquad Y_t=Y\times_T\{t\},
 \qquad\delta_{Y_t/\{t\}}=\delta_{Y/T}|_{Y_t}.
\end{equation}
Thus $i_t^!\mathcal B_{S/T}\simeq
R\Gamma(S_t,\omega^{\mathrm n}_{S_t/\{t\}})$.

\item In integer grading, give products and affine stabilizations
the induced relative volume-index normalizations of
Definition~\ref{def:relative-volume-index}. Normalized dualizing
complexes then respect exterior products. The normalized dualizing
complex on an affine stabilization is the ordinary inverse image
of the original one. These assertions also hold in parity grading
and for a finite or countable affine tail. The relative complexes
$\mathcal E$ and $\mathcal B$ are unchanged by these stabilizations.

\item Suppose, in addition, that $S$ has a finite filtration
$S=S_0\supset\cdots\supset S_m=\varnothing$ by finitely
presented closed subschemes with nice differences. Include
finite-factor presentations of this filtration over $T$. Then the images
of $\mathcal E_{S/T}$ and $\mathcal B_{S/T}$ in the ordinary
derived category belong to $D^b_c(T^{\mathrm{an}};\F)$, and
$\mathcal E_{S/T}\simeq\mathbb D_T\mathcal B_{S/T}$.
For two such families, their relative complexes respect
exterior products with the induced relative volume-index
normalization. The corresponding parity assertions hold
with finite-dimensional stalk cohomology in each parity.
\end{enumerate}
All comparisons are compatible with composition, further
refinement, base change, and localization, with the corresponding
higher homotopies.
\end{proposition}

\begin{proof}
\emph{Normalized sheaves.}
For a smooth finite-factor comparison $v:Y'\to Y$, the relative
normalized purity equivalence of
Definition~\ref{def:relative-normalized-dual-complexes} is
\[
 \omega_{Y'}[-2\delta_{Y'/T}]
 \simeq v^{-1}\bigl(\omega_Y[-2\delta_{Y/T}]\bigr).
\]
These equivalences inherit the composition homotopies of the
smooth-purity transformations from
Subsection~\ref{conv:section-five-sheaf-operations} and hence form
the finite-factor diagrams. Exact inverse
image transports them to the product charts, and
Lemma~\ref{lem:infinite-product-descent} descends them to $S$.
Apply the same assertion to the arrow and cofiber-sequence
diagrams to descend the maps and localization sequences
simultaneously with their objects.

For the finitely presented closed inclusion $Z\hookrightarrow S$,
use the finite-factor presentation of Subsection~\ref{subsec:relative-volume-preparations}:
its pullback to a chart $U=Y\times\A^J$ is
$Z_Y\times\A^J$, with $Z_Y\hookrightarrow Y$ closed. The closed and open product
comparisons of Lemma~\ref{lem:infinite-product-descent} pull the
finite-factor localization triangle to $U$. They also identify
its sections-with-support term with
$p_Z^{-1}(\omega_{Z_Y}[-2\delta_{Y/T}|_{Z_Y}])$.
Descent gives (i). With the induced relative volume-index
normalizations of Definition~\ref{def:relative-volume-index},
the product and stabilization assertions on sheaves follow from
\[
 \omega_{Y_1\times Y_2}\simeq\omega_{Y_1}\boxtimes\omega_{Y_2},
 \qquad
 \omega_{Y\times\A^e}[-2\delta_{Y/T}-2e]
       \simeq\operatorname{pr}_Y^{-1}
                    (\omega_Y[-2\delta_{Y/T}]).
\]

\emph{Compact-support trace descent.}
For a refinement between two finite factors of one chart,
the smooth comparison $v:Y'\to Y$ has contractible analytic
fibers. Indeed, write the common analytic chart as both
$U=Y\times\C^J$ and $U=Y'\times\C^{J'}$, with projections
$p:U\to Y$, $p':U\to Y'$ satisfying $p=vp'$, and let
$s':Y'\to U$ be the zero section of $p'$. For $y\in Y$, put
$F_y=v^{-1}(y)$. Restriction gives maps
\[
 F_y\xrightarrow{s'|_{F_y}}p^{-1}(y)\cong\C^J
     \xrightarrow{p'|_{p^{-1}(y)}}F_y
\]
whose composite is the identity of $F_y$. In the displayed
coordinates $p^{-1}(y)\cong\C^J$, the map
\[
 F_y\times[0,1]\longrightarrow F_y,\qquad
 (z,\tau)\longmapsto p'\bigl((1-\tau)s'(z)\bigr)
\]
contracts $F_y$: its endpoints are the identity and the constant
map with value $p'(0)$, where $0$ is the origin of these fiber
coordinates. The fiber $F_y$ is a smooth complex manifold of dimension $e$,
where $e$ is the relative dimension of $v$ on the component
in question. Proper-support base change identifies the stalk
of the normalized trace with
\[
 R\Gamma_c(F_y;\F)[2e]\longrightarrow\F.
\]
The complex orientation and Poincar\'e duality identify this
map with integration on the compactly supported top-degree
class, an equivalence because $F_y$ is contractible. Thus
$Rv_!\underline\F_{Y'}[2e]\to\underline\F_Y$ is an equivalence.
These are the normalized refinement equivalences in the
compact-support diagram. Their coherence is the transitivity
of the enhanced trace, before passing to any homotopy category.

For a finite-type \'etale cover $V\to Y$, the stalk at $y$
of its augmented trace \v Cech diagram has degree-$n$ term
$\F[V_y^{\,n+1}]$, the vector space with basis $V_y^{\,n+1}$.
Faces omit a coordinate, and the augmentation sends each basis
element of $\F[V_y]$ to $1$. Choose $v_*\in V_y$, which is
nonempty. The extra degeneracy
\[
 (v_0,\ldots,v_n)\longmapsto(v_*,v_0,\ldots,v_n)
\]
contracts the augmented chain complex. Its realization is
consequently $\underline\F_Y$.
To compare two finite atlases, fix a nice chart
$U=Y\times\A^J$ of the first. The pullback of the second is
a surjective finitely presented \'etale cover $V_U\to U$.
After enlarging $Y$, this cover descends to $V\to Y$.
It is surjective because its pullback along the surjective
projection $U\to Y$ is surjective. Its whole \v Cech nerve is
then the product with $\A^J$ of the nerve of $V\to Y$.
The chosen finite factors in each degree compare with these
factors by the normalized refinement equivalences above.
The finite-type augmented trace calculation, followed by the
proper-support map to $T$, therefore gives the comparison
on $U$. The double \v Cech diagram of the two atlases proves
that their trace realizations agree. Further iterated diagrams
and the enhanced traces give the higher compatibilities.
Lemma~\ref{lem:totalization-realization-representatives} supplies
simultaneous cochain representatives of the resulting
realizations and their coherent comparison maps.
Localization of the
finite-factor constant sheaves proves the triangle for
$\mathcal E$ in (i).

\emph{Identification with the direct image of the descended sheaf.}
On $U=Y\times\C^J$, Lemma~\ref{lem:infinite-product-descent}
gives $Rp_*p^{-1}\omega_{Y/T}^{\mathrm n}\simeq
\omega_{Y/T}^{\mathrm n}$. Hence
\[
 R(f|_U)_*\omega_{U/T}^{\mathrm n}
 \simeq Rg_*\omega_{Y/T}^{\mathrm n}
 \simeq\mathbb D_T
          Rg_!\underline\F_Y[2\delta_{Y/T}].
\]
These equivalences identify the dual trace maps with the
restriction maps of the sheaf \v Cech diagram. Since $Rf_*$
preserves limits, descent yields
\[
 Rf_*\omega^{\mathrm n}_{S/T}
 \simeq\operatorname*{holim}_{[n]\in\Delta}\mathbb D_T E_n
 \simeq\mathbb D_T\mathcal E_{S/T}.
\]
The last equivalence uses that derived duality takes colimits
in its argument to limits; it does not assume biduality for
an arbitrary unbounded object.

\emph{Base change and finiteness.}
The finite-factor proper-support base-change equivalences
\[
 b^*Rg_!\underline\F_Y[2\delta_{Y/T}]
 \simeq Rg'_!\underline\F_{Y\times_TT'}
                             [2\delta_{Y/T}|_{Y\times_TT'}]
\]
commute with the trace diagrams and their realization. This
proves the first formula in (ii). The second is the composite
\[
 \begin{aligned}
 b^!\mathcal B_{S/T}
 &\simeq b^!\mathbb D_T\mathcal E_{S/T}
  \simeq\mathbb D_{T'}b^{-1}\mathcal E_{S/T}\\
 &\simeq\mathbb D_{T'}\mathcal E_{S\times_TT'/T'}
  \simeq\mathcal B_{S\times_TT'/T'}.
 \end{aligned}
\]
Here the middle duality comparison is
\eqref{eq:duality-inverse-image}; the other comparisons are
the first formula and the direct-image identifications above.
The induced finite-factor presentations and finite affine covers of
Subsection~\ref{subsec:relative-volume-preparations} show that
these closed restrictions and base changes preserve local niceness.

Under the additional filtration hypothesis, each nice
difference has a finite factor $Y_a\to T$. Its companion is
$R(g_a)_!\underline\F_{Y_a}[2\delta_{Y_a/T}]$, hence bounded
constructible. Induction with the localization triangles
proves the same for $\mathcal E_{S/T}$, and finite-dimensional
Verdier duality proves it for $\mathcal B_{S/T}$ and gives
biduality. Finite-type compact-support K\"unneth, followed by
the trace realizations, gives the product formula for
$\mathcal E$; constructible Verdier duality gives it for
$\mathcal B$. Affine stabilization follows from
$R\Gamma_c(\A^e;\F)[2e]\simeq\F$ and, for the sheaf and
ordinary-direct-image description, from product homotopy
invariance. The comparisons are transitive as coordinates
are adjoined. Periodize the bounded finite-factor constructions
to obtain the parity assertions.

\end{proof}
\subsection{Topological pullbacks and descent of specialization}
\label{subsec:specialization-arrow-descent}
First we extend the preceding normalized dualizing complexes to
pullbacks along topological parameter maps. We then descend the
finite-factor specialization morphism together with its source
and target, as an object of the arrow category
$\operatorname{Fun}([1],\mathcal D(S^{\mathrm{an}};\F))$.

\begin{lemma}[Pullback to topological parameter spaces]
\label{lem:topological-pullback-presentations}
Let $f:S\to T$ have the setup and compatible finite-factor presentations
of Definition~\ref{def:relative-normalized-dual-complexes}: $S$ is
quasi-compact, separated, and locally nice, and $T$ is a separated
complex scheme of finite type. Fix a relative volume-index
normalization in integer grading. Write
$\omega^{\mathrm n}_{S/T}$, $\mathcal B_{S/T}$, and
$\mathcal E_{S/T}$ for the objects of that definition; by
Proposition~\ref{lem:relative-compact-support},
$\mathcal B_{S/T}\simeq\mathbb D_T\mathcal E_{S/T}$.
Let $h:M\to T^{\mathrm{an}}$ be continuous, with $M$ paracompact,
locally compact Hausdorff, and of finite covering dimension.
\begin{enumerate}[(i),wide=0pt]
\item There is a natural equivalence in $\mathcal D(M;\F)$,
\[
 \mathbb D_M(h^{-1}\mathcal E_{S/T})
       \simeq h^!\mathcal B_{S/T}.
\]

\item Put $S_M=S^{\mathrm{an}}\times_{T^{\mathrm{an}}}M$, with
projection $f_M:S_M\to M$. For each finite factor $g:Y\to T$
in the finite-factor presentation, put
\[
 Y_M=Y^{\mathrm{an}}\times_{T^{\mathrm{an}}}M,
 \qquad b_Y:Y_M\to Y^{\mathrm{an}},
 \qquad\delta_M=b_Y^*\delta_{Y/T}.
\]
Assume that every $Y_M$ is paracompact, locally compact
Hausdorff, and of finite covering dimension. Assume also that
every pullback $v_M:Y'_M\to Y_M$ of a smooth finite-factor
comparison $v:Y'\to Y$ of complex relative dimension $e$ is
locally on its source an oriented product projection of real
relative dimension $2e$.
These hypotheses hold for semialgebraic rays, closed half-rays,
and coordinate faces with semialgebraic maps to $T^{\mathrm{an}}$.

For the pulled-back charts
$U_M=Y_M\times\C^J\to S_M$, let $p_M:U_M\to Y_M$ be the
projection. The complexes
\[
 p_M^{-1}(\omega_{Y_M}[-2\delta_M])
       \in\mathcal D(U_M;\F)
\]
descend to
$\omega^{\mathrm n}_{S_M/M}\in\mathcal D(S_M;\F)$.
Pullbacks of common refinements over $T$ give comparison
equivalences independent of the atlas and finite-factor presentations,
with the relative volume indices transported as in
Definition~\ref{def:relative-volume-index}. Their direct image
satisfies
\[
 R(f_M)_*\omega^{\mathrm n}_{S_M/M}
   \simeq\mathbb D_M(h^{-1}\mathcal E_{S/T})
   \simeq h^!\mathcal B_{S/T}
       \qquad\text{in }\mathcal D(M;\F).
\]
This is a sheaf construction on $S_M$; no algebraic local
niceness of this topological pullback is asserted.
\end{enumerate}
The assertions hold in the periodic derived categories using
periodized finite-factor complexes, without integer volume indices.
\end{lemma}
\begin{proof}
For (i), apply \eqref{eq:duality-inverse-image} to
$E=\mathcal E_{S/T}$ and use
$\mathbb D_T\mathcal E_{S/T}\simeq\mathcal B_{S/T}$ from
Proposition~\ref{lem:relative-compact-support}.

For (ii), the oriented product hypothesis gives
$v_M^!\simeq v_M^{-1}[2e]$. Together with
$\delta'_M=v_M^*\delta_M+e$, this gives
\[
 \omega_{Y'_M}[-2\delta'_M]
 \simeq v_M^{-1}(\omega_{Y_M}[-2\delta_M]).
\]
The pulled-back atlas is a jointly surjective family of local
homeomorphisms, and its overlaps are the pullbacks of the
chosen overlap charts over $T$. The displayed normalized purity
equivalences, with the composition and base-change compatibilities
inherited from oriented purity, therefore give a coherent descent
diagram on this atlas. Lemma~\ref{lem:infinite-product-descent}(ii) descends
it to $\omega^{\mathrm n}_{S_M/M}$. Applying the same descent
to pullbacks of common refinements proves independence.

Let $E_\bullet$ be the compact-support trace diagram of
\eqref{eq:relative-compact-factor}. Proper-support base change
identifies its pulled-back finite-factor terms by
\[
 h^{-1}Rg_!\underline\F_Y[2\delta_{Y/T}]
 \simeq Rg_{M,!}\underline\F_{Y_M}[2\delta_M].
\]
These equivalences commute with the normalized traces, so the
trace diagram of the pulled-back finite-factor presentation is
$E^M_\bullet\simeq h^{-1}E_\bullet$ and
\[
 \operatorname*{hocolim}_{[n]\in\Delta^{\mathrm{op}}}E_n^M
 \simeq h^{-1}\mathcal E_{S/T}.
\]
The chart-to-trace comparison in the proof of
Proposition~\ref{lem:relative-compact-support} applies here:
product homotopy invariance gives
$Rp_{M,*}p_M^{-1}K\simeq K$, and finite-dimensional duality gives
$Rg_{M,*}\omega_{Y_M}\simeq
\mathbb D_M Rg_{M,!}\underline\F_{Y_M}$.
Thus the product of the chart direct images in degree $n$ is
$\mathbb D_M E_n^M$, with restriction maps dual to the trace
maps. Descent and preservation
of limits by $R(f_M)_*$ yield
\[
 R(f_M)_*\omega^{\mathrm n}_{S_M/M}
 \simeq\operatorname*{holim}_{[n]\in\Delta}\mathbb D_M E_n^M
 \simeq\mathbb D_M(h^{-1}\mathcal E_{S/T}).
\]
Part~(i) gives the remaining equivalence.

For the stated semialgebraic parameter spaces, each $Y_M$ is
locally closed semialgebraic in finite-dimensional charts, hence
has the required topological properties. Smooth algebraic maps
are analytically locally product projections; their pullbacks
preserve these projections and their complex orientations.
This verifies the additional hypotheses.

For the parity assertions, we verify coefficient change both on the
parameter space and on the pulled-back finite factors. Let
$X=M$ or $X=Y_M$. The assertion is immediate if $X$ is empty; otherwise put
$d=\dim_{\mathrm{cov}}X$, where $\dim_{\mathrm{cov}}$ denotes
covering dimension, and write $a_X:X\to\mathrm{pt}$.
Put $\Lambda=\F[u,u^{-1}]$, with $|u|=2$ and $du=0$, and let
$\Lambda_X$ be its constant sheaf complex. Write $a_{X,\Lambda}^!$
for the right adjoint of $Ra_{X,!}$ with $\Lambda$ coefficients,
and put $\omega_{X,\Lambda}:=a_{X,\Lambda}^!\Lambda$.
Extension of scalars and
the projection formula give a canonical morphism
\[
 c_X:\operatorname{Per}_2\omega_X
       =\omega_X\otimes_{\F}\Lambda_X
       \longrightarrow\omega_{X,\Lambda}.
\]
It is adjoint to the trace
$Ra_{X,!}\omega_X\to\F$ tensored with $\Lambda$.
The underlying rational complex of $\omega_{X,\Lambda}$ is
$a_X^!\Lambda$, by the same adjunction argument used in
Subsection~\ref{subsec:finite-presentations}.

We verify that $c_X$ is an equivalence on a basis of open sets.
Given $x\in O\subset X$ with $O$ open, local compactness and
the Hausdorff property give a continuous function
$\varphi:X\to[0,1]$ with $\varphi(x)=1$ and compact support
contained in $O$. Put
\[
 U=\{\varphi>0\},\qquad K=\operatorname{supp}(\varphi),\qquad
 K_n=\{\varphi\ge1/n\}\quad(n\ge1).
\]
Then $U=\bigcup_{n\ge1}K_n$ is a relatively compact,
$\sigma$-compact open neighborhood of $x$. It is paracompact,
as every $\sigma$-compact locally compact Hausdorff space is.
Each $K_n$ is closed in $X$, so $\dim_{\mathrm{cov}}K_n\le d$;
the countable closed-sum theorem for covering dimension on
normal spaces gives $\dim_{\mathrm{cov}}U\le d$.
Thus $U$ has sheaf-cohomological dimension at most $d$.
The compact space $K$ has the same bound. If $j_U:U\hookrightarrow K$
is the open inclusion, then
\[
 C_U:=R\Gamma_c(U;\F)\simeq R\Gamma(K,j_{U,!}\underline\F_U)
                         \in D^{[0,d]}(\F).
\]
Here $D^{[a,b]}$ denotes complexes with cohomology zero outside
the integers $a,\ldots,b$. Verdier duality gives
\[
 R\Gamma(U,\omega_X|_U)\simeq R\operatorname{Hom}_{\F}(C_U,\F)
                              \in D^{[-d,0]}(\F).
\]
Taking stalks over this basis shows that $\omega_X$ is bounded,
with cohomology sheaves only in degrees $-d,\ldots,0$.

For a bounded sheaf complex $K'$ on $U$, the canonical map
\[
 \operatorname{Per}_2R\Gamma(U,K')\longrightarrow
 R\Gamma(U,\operatorname{Per}_2K')
\]
is an equivalence. Indeed, in the hypercohomology spectral
sequence the sheaf-cohomological degree lies between $0$ and $d$;
boundedness of $K'$ therefore makes only finitely many of its
even shifts contribute to any fixed total degree. Finite direct
sums commute with sheaf cohomology, which proves the comparison.
Also, boundedness of $C_U$ gives the equivalence
\[
 \operatorname{Per}_2R\operatorname{Hom}_{\F}(C_U,\F)
       \xrightarrow{\sim}R\operatorname{Hom}_{\F}(C_U,\Lambda).
\]
To check this, represent $C_U$ by a complex of vector spaces
concentrated in degrees $0,\ldots,d$. In total degree $n$,
\[
 \operatorname{Hom}_{\F}^{n}(C_U,\Lambda)
 =\bigoplus_{k=0}^{d}
   \operatorname{Hom}_{\F}(C_U^k,\F)\otimes_{\F}\Lambda^{n+k}.
\]
Each $\Lambda^{n+k}$ is zero or one-dimensional, and this is
the degree-$n$ term of
$\operatorname{Hom}_{\F}(C_U,\F)\otimes_{\F}\Lambda$.
The comparison respects differentials; all degrees of $\Lambda$
are even. Thus no finite-dimensionality of the $C_U^k$ is
required. Applying these two comparisons gives
\[
 \begin{aligned}
 R\Gamma(U,\operatorname{Per}_2\omega_X)
 &\simeq\operatorname{Per}_2R\operatorname{Hom}_{\F}(C_U,\F)\\
 &\simeq R\operatorname{Hom}_{\F}(C_U,\Lambda)
  \simeq R\Gamma(U,\omega_{X,\Lambda}).
 \end{aligned}
\]
These are the maps induced by $c_X$, and are natural in $U$.
They prove that $c_X$ is an equivalence. For a map $r:X'\to X$
in the chosen diagram, the proper-support trace comparisons commute
with $c_{X'}$ and $c_X$: after adjunction, both composites
are the trace for $a_{X'}=a_Xr$ tensored with $\Lambda$, by
composition of the adjunction counits. This identity is an
identity of the enhanced natural transformations, so it also
preserves their composition coherences. Oriented product purity
uses the relative orientation class, whose image under
$\F\to\Lambda$ is the same orientation class with extended
coefficients; hence its comparison also commutes with $c_X$.
Exactness of extension of scalars gives the corresponding
compatibility for localization boundaries. The preceding
constructions and descent can therefore
be performed directly with $\Lambda$ coefficients, using
$\operatorname{Per}_2\omega_X$ on every pulled-back finite factor.
This proves the parity assertions, including their comparisons,
without interchanging periodization with an infinite \v Cech limit.
\end{proof}

\begin{lemma}[Descent of the specialization arrow]
\label{lem:specialization-arrow-descent}
Let $S$ be a quasi-compact separated locally nice complex
scheme and $f:S\to\A^1_s$ a morphism with compatible finite-factor
presentations over $\A^1_s$, including finite-factor presentations of the closed--open diagram
\[
 S_0:=S\times_{\A^1}\{0\}\xrightarrow{i}S
          \xleftarrow{j}S^\times:=S\times_{\A^1}\Gm
\]
as in
Subsection~\ref{subsec:relative-volume-preparations}.
Use the normalized dualizing complexes of
Definition~\ref{def:relative-normalized-dual-complexes}, with a
chosen relative volume-index normalization in integer grading
on this closed--open diagram.
On analytic spaces put
\[
 H=[0,\infty),\quad L=(0,\infty),\quad
 S_H=S^{\mathrm{an}}\times_{\C}H,\quad
 S_L=S^{\mathrm{an}}\times_{\C}L,
 \qquad \ell:S_L\hookrightarrow(S^\times)^{\mathrm{an}}.
\]
Use the same letters $i,j$ for the analytic inclusions.
For $M=H,L$, let $\omega^{\mathrm n}_{S_M/M}$ be the complexes
of Lemma~\ref{lem:topological-pullback-presentations}(ii);
its hypotheses hold for the semialgebraic inclusions
$H,L\hookrightarrow\C$.
\begin{enumerate}[(i),wide=0pt]
\item There is a canonical morphism in $\mathcal D(S^{\mathrm{an}};\F)$,
\begin{equation}\label{eq:normalized-specialization-arrow}
 \operatorname{Sp}^{\mathrm n}_S:
 \mathcal K^{\mathrm n}_{S,\rho}
 :=Rj_*\ell_*\omega^{\mathrm n}_{S_L/L}[-1]
 \longrightarrow i_*\omega^{\mathrm n}_{S_0/\{0\}}.
\end{equation}
On each finite factor it is \eqref{eqSp}, shifted by
$-2\delta_{Y/\A^1}$, without a product assumption on $Y^\times$.
The symbol $\rho$ uses the inward-oriented ray fixed in
Lemma~\ref{lp:lem:bm-collision-specialization}.

\item The morphism in (i), including its source and target,
is compatible in the derived enhancement with smooth
finite-factor comparisons, common refinements over $\A^1$,
and proper traces for the finitely presented closed inclusions
of families included in these diagrams. If $v:V\hookrightarrow S$ is an open
neighborhood of $S_0$, with induced finite-factor presentations, then
\[
 v^{-1}\mathcal K^{\mathrm n}_{S,\rho}
       \simeq\mathcal K^{\mathrm n}_{V,\rho},\qquad
 v^{-1}\operatorname{Sp}^{\mathrm n}_S
       \simeq\operatorname{Sp}^{\mathrm n}_V.
\]
Consequently
\[
 R\Gamma(\operatorname{Sp}^{\mathrm n}_V)\circ\operatorname{res}_v
       \simeq R\Gamma(\operatorname{Sp}^{\mathrm n}_S),\qquad
 \operatorname{res}_v:
 R\Gamma(S^{\mathrm{an}},\mathcal K^{\mathrm n}_{S,\rho})
 \longrightarrow R\Gamma(V^{\mathrm{an}},\mathcal K^{\mathrm n}_{V,\rho}),
\]
under the common target
$R\Gamma(S_0^{\mathrm{an}},\omega^{\mathrm n}_{S_0/\{0\}})$.
No equivalence is asserted for $\operatorname{res}_v$.
\end{enumerate}
Both assertions hold in parity grading with the periodic
complexes of Definition~\ref{def:relative-normalized-dual-complexes},
without choosing relative volume indices.
\end{lemma}
\begin{proof}
\emph{Finite factors and their comparisons.}
For a finite factor $Y\to\A^1$, let
$a:Y_H\hookrightarrow Y$ be the closed inclusion and let
$\partial_{Y_H}:Rk_*\omega_{Y_L}\to b_*\omega_{Y_0}[1]$
be the localization boundary for
$b:Y_0\hookrightarrow Y_H$ and $k:Y_L\hookrightarrow Y_H$.
By \eqref{eq:specialization-half-ray-boundary},
\[
 \operatorname{Sp}^{\omega}_Y
   =a_*\bigl(\partial_{Y_H}[-1]\bigr):
 Rj_{Y,*}\ell_{Y,*}\omega_{Y_L}[-1]
       \longrightarrow i_{Y,*}\omega_{Y_0}.
\]
Shift the arrow by $-2\delta_{Y/\A^1}$.
The proper and smooth naturality of
Lemma~\ref{lp:lem:bm-collision-specialization} gives the comparisons
of these arrows. For a smooth finite-factor comparison
$g:Y'\to Y$ of relative complex dimension $e$, the purity shift
$[2e]$ cancels against
\[
 \delta_{Y'/\A^1}=g^*\delta_{Y/\A^1}+e.
\]
Thus normalized smooth comparisons are degree-zero equivalences
of arrows, and proper traces for the presented closed inclusions
are morphisms of arrows. Their compatibilities are the enhanced
naturality of the same localization boundaries and traces.

\emph{Pullback to a nice chart.}
For $p:U=Y\times\C^J\to Y$, write
$p_M:U_M\to Y_M$ for $M=0,L$, and
$\delta_M=\delta_{Y/\A^1}|_{Y_M}$. Apply $p^{-1}$ to the
normalized finite-factor arrow. The product comparisons of
Lemma~\ref{lem:infinite-product-descent} identify its endpoints as
\[
 \begin{aligned}
 p^{-1}Rj_{Y,*}\ell_{Y,*}
          (\omega_{Y_L}[-2\delta_L])[-1]
 &\simeq Rj_{U,*}\ell_{U,*}
          p_L^{-1}(\omega_{Y_L}[-2\delta_L])[-1],\\
 p^{-1}i_{Y,*}(\omega_{Y_0}[-2\delta_0])
 &\simeq i_{U,*}p_0^{-1}(\omega_{Y_0}[-2\delta_0]).
 \end{aligned}
\]
Here the open comparison is product base change for $Rj_*$;
the closed direct-image comparisons are exact.

\emph{Descent of the arrow.}
Use the atlas and \v Cech nerve chosen in
Definition~\ref{def:relative-normalized-dual-complexes}.
On its charts and all overlaps, the preceding construction gives
local arrows and comparison equivalences between them, including
the source and target comparisons just displayed. The specified
composition and base-change compatibilities make this an object
of the descent diagram of arrow categories.
Lemma~\ref{lem:infinite-product-descent}(ii), with $A=[1]$,
descends it to an arrow on $S^{\mathrm{an}}$.
Evaluation at $0,1\in[1]$ commutes with descent, so its endpoints
are respectively
$Rj_*\ell_*\omega^{\mathrm n}_{S_L/L}[-1]$ and
$i_*\omega^{\mathrm n}_{S_0/\{0\}}$. This proves (i).
Common refinements give equivalences of the local arrow diagrams;
presented closed inclusions give morphisms between them.
Their descent proves the corresponding assertions in (ii).

Let $v:V\hookrightarrow S$ be an arbitrary open neighborhood
of $S_0$. The open scheme $V$ need not be quasi-compact.
For $M=H,L$, put $V_M=V^{\mathrm{an}}\times_{\C}M$ and write
$v_M:V_M\hookrightarrow S_M$ for its open inclusion. Define
\[
 \omega^{\mathrm n}_{V_M/M}:=
       v_M^{-1}\omega^{\mathrm n}_{S_M/M}.
\]
These restrictions have the normalized finite-factor descriptions
on every induced nice chart; open descent identifies their
comparisons on overlaps. Thus the notation applies to $V$
without requiring a finite atlas of $V$.
Let $j_V:(V^\times)^{\mathrm{an}}\hookrightarrow V^{\mathrm{an}}$,
$\ell_V:V_L\hookrightarrow(V^\times)^{\mathrm{an}}$, and
$i_V:S_0^{\mathrm{an}}\hookrightarrow V^{\mathrm{an}}$ be the
induced inclusions, where $V^\times=V\times_{\A^1}\Gm$.
Set
\[
 \mathcal K^{\mathrm n}_{V,\rho}:=
       Rj_{V,*}\ell_{V,*}\omega^{\mathrm n}_{V_L/L}[-1].
\]
Open base change gives the canonical endpoint identifications
\[
 v^{-1}\mathcal K^{\mathrm n}_{S,\rho}\simeq
       \mathcal K^{\mathrm n}_{V,\rho},\qquad
 v^{-1}i_*\omega^{\mathrm n}_{S_0/\{0\}}\simeq
       i_{V,*}\omega^{\mathrm n}_{S_0/\{0\}}.
\]
Under these identifications, define
$\operatorname{Sp}^{\mathrm n}_V:=v^{-1}\operatorname{Sp}^{\mathrm n}_S$.
On each induced chart this is the normalized half-ray boundary,
by naturality of localization under open restriction. It therefore
agrees with the local specialization construction. Applying
derived global sections and the natural restriction map proves
the displayed identity in (ii).
Periodizing the same diagrams proves the parity assertions.
\end{proof}

\begin{corollary}[Cohomology actions on normalized specialization]
\label{cor:normalized-specialization-module}
Let $f:S\to\A^1_s$ satisfy the hypotheses of
Lemma~\ref{lem:specialization-arrow-descent}: $S$ is quasi-compact,
separated, and locally nice, and its compatible finite-factor presentations
include finite-factor presentations of $S_0\hookrightarrow S\hookleftarrow S^\times$.
In integer grading choose a relative volume-index normalization
$\{\delta_{Y/\A^1}\}$ on this finite-factor presentation diagram as in
Definition~\ref{def:relative-volume-index}. In parity grading use
the periodic complexes without choosing volume indices.

Put $S_H=f^{-1}(H)$ and $S_L=f^{-1}(L)$ on analytic spaces.
Let $v:N\hookrightarrow S_H$ be an open neighborhood of $S_0$,
and write
\[
 N^+=N\cap S_L,\quad
 i_N:S_0\hookrightarrow N,\quad j_N:N^+\hookrightarrow N,
 \qquad
 \omega_+^{\mathrm n}=\omega^{\mathrm n}_{S_L/L}|_{N^+},\quad
 \omega_0^{\mathrm n}=\omega^{\mathrm n}_{S_0/\{0\}}.
\]
Use the complexes constructed in
Definition~\ref{def:relative-normalized-dual-complexes} and
Lemma~\ref{lem:topological-pullback-presentations}.
The descended half-ray specialization restricts to an arrow
\[
 \sigma_N^{\mathrm n}:\mathcal K^{\mathrm n}_{N,\rho}:=Rj_{N,*}\omega_+^{\mathrm n}[-1]
      \longrightarrow T_N^{\mathrm n}:=i_{N,*}\omega_0^{\mathrm n}
      \quad\text{in }\mathcal D(N;\F),
\]
or in $\mathcal D_{\mathrm{per}}(N;\F)$, respectively.
For every $c:\underline\F_N\to\underline\F_N[d]$, one has
\[
 \sigma_N^{\mathrm n}[d]\circ m_{c,\mathcal K^{\mathrm n}_{N,\rho}}
     \simeq m_{c,T_N^{\mathrm n}}\circ\sigma_N^{\mathrm n}
\] with $m_{c,E}$ as in Lemma \ref{lem:specialization-module}.
Thus, for $\operatorname{sp}^{\mathrm n}_N
 =H^*R\Gamma(N,\sigma_N^{\mathrm n})$,
\begin{equation}\label{eq:normalized-specialization-module}
 \operatorname{sp}^{\mathrm n}_N(c\frown\xi)
   =(i_N^{-1}c)\frown\operatorname{sp}^{\mathrm n}_N(\xi),
 \qquad \xi\in H^qR\Gamma(N,\mathcal K^{\mathrm n}_{N,\rho}).
\end{equation}
The action and degrees are those of
Lemma~\ref{lem:specialization-module}. In particular, the source
action includes $[-1]$; the shifts $[-2\delta_{Y/\A^1}]$ are even.
The class $c$ is a class on $N$ itself: no finite-factor
representative of $c$ is required.
In parity grading the displayed action uses
$\operatorname{Per}_2(c):\underline\F_N^{\mathrm{per}}
\to\underline\F_N^{\mathrm{per}}[d]$, where
$\underline\F_N^{\mathrm{per}}:=\operatorname{Per}_2(\underline\F_N)$
is the monoidal unit of $\mathcal D_{\mathrm{per}}(N;\F)$;
degrees are then read modulo $2$.

The identities commute with open restriction and the specified
proper and smooth comparisons when the classes are pulled back
along those comparisons. They also commute with an inverse
normalized affine-bundle comparison wherever that comparison is
an equivalence. If a geometric product presentation, with normalized
complex $\omega_U^{\mathrm n}$ on its factor $U$, gives
\[
 N^+\simeq L\times U,\qquad
 \omega_+^{\mathrm n}\simeq\omega_L\boxtimes \omega_U^{\mathrm n},\qquad
 c|_{N^+}=\operatorname{pr}_U^{-1}c_U,
\]
use its normalized sheaf comparison and the ray orientation to
identify $R\Gamma(N,\mathcal K^{\mathrm n}_{N,\rho})$ with $R\Gamma(U,\omega_U^{\mathrm n})$.
Then \eqref{eq:normalized-specialization-module} is the input-class
formula of Lemma~\ref{lem:specialization-module}(iii) for these
normalized complexes.
\end{corollary}
\begin{proof}
Let $a:S_H\hookrightarrow S^{\mathrm{an}}$, $j_H:S_L\hookrightarrow S_H$,
and $i_H:S_0\hookrightarrow S_H$ be the inclusions.
The half-ray construction in the proof of
Lemma~\ref{lem:specialization-arrow-descent} identifies
$\operatorname{Sp}^{\mathrm n}_S$ with $a_*\sigma^{\mathrm n}_{S_H}$, where
\[
 \sigma^{\mathrm n}_{S_H}:
 Rj_{H,*}\omega^{\mathrm n}_{S_L/L}[-1]
      \longrightarrow i_{H,*}\omega_0^{\mathrm n}.
\]
Composition of direct images identifies these two endpoints,
after applying $a_*$, with the source and target of
$\operatorname{Sp}^{\mathrm n}_S$; full faithfulness of the closed
pushforward determines the displayed arrow. Open base change
identifies $v^{-1}\sigma^{\mathrm n}_{S_H}$ with
$\sigma_N^{\mathrm n}:\mathcal K^{\mathrm n}_{N,\rho}\to T_N^{\mathrm n}$.
Lemma~\ref{lem:specialization-module}(i), applied on $N$, gives
the module homotopy. Taking $R\Gamma$ and cohomology proves
\eqref{eq:normalized-specialization-module}.

These actions commute with inverse image by monoidality, and
with direct image by the projection formula for the invertible
object $\underline\F[d]$, as in the proof of
Lemma~\ref{lem:specialization-module}(ii).
Apply this naturality to the comparison arrows of
Lemma~\ref{lem:specialization-arrow-descent}(ii). For refinements,
pull $c$ to the nice charts and their common refinements; no
class on a finite factor is needed. Inverting an intertwining
equivalence preserves the identity. The product assertion is
the ray calculation of Lemma~\ref{lem:specialization-module}(iii)
with $\omega_U^{\mathrm n}$; the even volume shifts add no sign. Periodizing $c$
and using the periodic unit gives the same tensor argument in
$\mathcal D_{\mathrm{per}}$, including the source shift $[-1]$.
\end{proof}

\begin{remark}[Dual restriction and descent]
\label{rem:restriction-arrow-descent}
In the setup of Lemma~\ref{lem:specialization-arrow-descent},
the comparison between compact-support and dualizing localization
diagrams is defined before any product trivialization. We spell
out its consequence for the descended compact-support companion.
Write
\[
 h:H=[0,\infty)\hookrightarrow\C,\qquad
 i_H:\{0\}\hookrightarrow H,\qquad j_H:L=(0,\infty)\hookrightarrow H,
\]
and put $E_H=h^{-1}\mathcal E_{S/\A^1}$ and $V_0=i_H^{-1}E_H$.
Choose a trivialization $j_H^{-1}E_H\simeq\underline V_L$,
where $V$ is a complex of $\F$-vector spaces and
$\underline V_L$ denotes its constant sheaf. The constant-sheaf
localization triangle on $H$ gives
$R\Gamma(H,j_{H,!}\underline V_L)=0$. Hence evaluation gives
\[
 \varepsilon_0:R\Gamma(H,E_H)\xrightarrow{\sim}V_0,
 \qquad
 \varepsilon_L:R\Gamma(L,\underline V_L)\xrightarrow{\sim}V.
\]
Define the restriction morphism
\[
 r:=\varepsilon_L\circ\operatorname{res}_{H,L}
                   \circ\varepsilon_0^{-1}:V_0\longrightarrow V.
\]
Compact-support localization gives
\[
 R\Gamma_c(L,\underline V_L)\longrightarrow R\Gamma_c(H,E_H)
 \longrightarrow V_0\xrightarrow{\partial_c}
 R\Gamma_c(L,\underline V_L)[1].
\]
For the constant complex $\underline V_H$, its middle term is
zero, and its boundary is the inward-ray equivalence
\[
 \beta_V:V\xrightarrow{\sim}R\Gamma_c(L,\underline V_L)[1].
\]
The adjunction map
$E_H\to Rj_{H,*}j_H^{-1}E_H\simeq\underline V_H$
is $r$ at zero and the identity on $L$. Naturality of localization
therefore gives $\partial_c=\beta_V\circ r$, as in
Lemma~\ref{lem:specialization-restriction}.

For the trace terms $E_n$ of \eqref{eq:relative-compact-factor},
put $E_n^H=h^{-1}E_n$. Proper-support base change identifies
these with the pulled-back finite-factor terms, and
$E_H\simeq\operatorname*{hocolim}_{\Delta^{\mathrm{op}}}E_n^H$.
On each finite factor, duality identifies compact-support
localization with normalized dualizing localization; the
identifications commute with the normalized proper and smooth
comparisons. Proposition~\ref{lem:relative-compact-support}
and Lemma~\ref{lem:topological-pullback-presentations}, using
\[
 \mathbb D_HE_H\simeq
 \operatorname*{holim}_{[n]\in\Delta}\mathbb D_HE_n^H,
\]
identify the dual boundary with the descended half-ray boundary
of Lemma~\ref{lem:specialization-arrow-descent}. Writing
$W^\vee:=R\operatorname{Hom}_{\F}(W,\F)$ for linear duality,
the endpoint identifications are
\[
 \begin{aligned}
 V_0^\vee&\simeq R\Gamma(S_0^{\mathrm{an}},
                         \omega^{\mathrm n}_{S_0/\{0\}}),\\
 V^\vee&\simeq R\Gamma(S_L,
                         \omega^{\mathrm n}_{S_L/L})[-1].
 \end{aligned}
\]
Under these identifications and composition of direct images,
$r^\vee=R\Gamma(\operatorname{Sp}^{\mathrm n}_S)$.
The argument uses duality from colimits to limits and requires
neither biduality nor perfectness of $V$. With periodic
coefficients, put $\Lambda=\F[u,u^{-1}]$, $|u|=2$, $du=0$,
take $V\in\mathcal D(\Lambda)$, and use
$W^\vee=R\operatorname{Hom}_{\Lambda}(W,\Lambda)$;
the same localization and adjunction argument applies.
A product identification of the descended open family supplies
the stipulated trivialization; individual finite total factors
need not have product decompositions over the punctured base.
\end{remark}
\subsection{Parameterized morphisms and chain operations}
\label{subsec:parameter-chains}

To pass from parameterized sheaf morphisms to cochain operations,
we use the following singular-chain description. The conventions
for perfect complexes and periodization are those of
Appendix~\ref{app:complexes-duality}.

\begin{definition}[Parameter-chain representative]
\label{def:parameter-chain-representative}
An \emph{admissible parameter space} is a paracompact
finite-dimensional smooth real manifold, possibly with corners,
or a locally finite finite-dimensional polyhedron. Let $T$ be
such a space, let $C,D$ be perfect complexes of $\F$-vector
spaces, and let
\[
 \alpha:\omega_T\otimes C\longrightarrow\omega_T\otimes D
\]
be a morphism in $\mathcal D(T;\F)$, where
$\omega_T$ is the dualizing complex fixed in
Subsection~\ref{conv:section-five-sheaf-operations}, and
$\omega_T\otimes C$ abbreviates
$\omega_T\otimes\underline C_T$.
In parity grading, $C,D$ are perfect two-periodic complexes,
and the morphism belongs to $\mathcal D_{\mathrm{per}}(T;\F)$.
All parameter restrictions of $\alpha$ mean extraordinary
pullback, with the canonical identification $u^!\omega_T\simeq
\omega_{T'}$ for a continuous map $u:T'\to T$ between
admissible parameter spaces.

Write $S_T=C_{-*}(T;\F)$ for normalized singular chains,
with an $n$-simplex in cohomological degree $-n$ and differential
the singular boundary. A \emph{parameter-chain representative}
of $\alpha$ is a degree-zero cochain map
\begin{equation}\label{eq:parameter-chain-representative}
 \Theta_\alpha:S_T\otimes C\longrightarrow D
\end{equation}
representing $\alpha$ under the natural equivalence
\begin{equation}\label{eq:parameter-chain-adjunction}
 \begin{aligned}
 R\operatorname{Hom}_T(\omega_T\otimes C,\omega_T\otimes D)
 &\simeq R\Gamma\bigl(T,
             \underline{R\operatorname{Hom}_{\F}(C,D)}\bigr)\\
 &\simeq R\operatorname{Hom}_{\F}(S_T\otimes C,D).
 \end{aligned}
\end{equation}
The existence and compatibility properties of this equivalence
are proved in Lemma~\ref{lem:parameter-chain-operations}.
For a coherent diagram of relative morphisms, transport the
whole diagram through this equivalence before choosing cochain
representatives. The resulting comparisons include homotopies
and higher homotopies; they need not be equalities between
independently chosen representatives.

In parity grading, replace $S_T$ by its direct-sum
two-periodization $\operatorname{Per}_2S_T$ and use the derived
Hom of two-periodic complexes. In particular, linear duals and
Hom complexes use products, not direct sums, over the singular
chain degrees.
\end{definition}

\begin{lemma}[Parameter-chain operations]
\phantomsection\label{lem:parameter-chain-operations}
The equivalence \eqref{eq:parameter-chain-adjunction} is natural
in $C,D$ and in parameter restrictions. It has the following
properties.
\begin{enumerate}[wide=0pt]
\item[(i)] For a homogeneous $n$-chain $c$, write
      $\Theta_\alpha(c)(x)=\Theta_\alpha(c\otimes x)$. Then
\[
 |\Theta_\alpha(c)|=-n,
 \qquad
 d_{\operatorname{Hom}}\Theta_\alpha(c)
                  =\Theta_\alpha(\partial c).
\]
For $u:T'\to T$, the map
$\Theta_\alpha\circ(u_*\otimes\mathrm{id})$ represents
$u^!\alpha$. If $\gamma:[0,1]\to T$ joins $t_0$ to $t_1$,
then $H=\Theta_\alpha(\gamma)$ has degree $-1$ and satisfies
\[
 d_DH+Hd_C=
       \Theta_\alpha([t_1])-\Theta_\alpha([t_0]).
\]
Here the endpoint maps represent $i_{t_j}^!\alpha:C\to D$.

\item[(ii)] Suppose
$\beta:\omega_T\otimes D\to\omega_T\otimes E$ is another
relative morphism. Let
\[
 \Delta_T:S_T\longrightarrow S_T\otimes S_T
\]
be the Alexander--Whitney diagonal: on a singular $n$-simplex
$\sigma$ it is
\[
 \Delta_T\sigma=
    \sum_{p=0}^n
       \sigma|[0,\ldots,p]\otimes\sigma|[p,\ldots,n].
\]
The cochain map
\begin{equation}\label{eq:parameter-chain-composition}
 \Theta_\beta\circ
 (\mathrm{id}_{S_T}\otimes\Theta_\alpha)\circ
 (\Delta_T\otimes\mathrm{id}_C)
 :S_T\otimes C\longrightarrow E
\end{equation}
represents $\beta\circ\alpha$. Thus it is homotopic to any
chosen representative $\Theta_{\beta\circ\alpha}$.

\item[(iii)] For a relative morphism
$\beta:\omega_{T'}\otimes C'\to\omega_{T'}\otimes D'$, let
$\alpha\boxtimes\beta$ be its exterior product with $\alpha$,
using $\omega_{T\times T'}\simeq\omega_T\boxtimes\omega_{T'}$.
For the Eilenberg--Zilber chain cross product
$\operatorname{EZ}:S_T\otimes S_{T'}\to S_{T\times T'}$,
there is a comparison homotopy
\[
 \begin{aligned}
 &\Theta_{\alpha\boxtimes\beta}\circ
       (\operatorname{EZ}\otimes\mathrm{id}_{C\otimes C'})\\
 &\qquad\simeq
 (\Theta_\alpha\otimes\Theta_\beta)\circ
 \bigl(\mathrm{id}_{S_T}\otimes
       \tau_{S_{T'},C}\otimes\mathrm{id}_{C'}\bigr),
 \end{aligned}
\]
where $\tau$ is the Koszul interchange. These comparisons and
those of \textup{(i)--(ii)} can be constructed as maps of coherent
diagrams, with all higher compatibility homotopies.
\end{enumerate}
\end{lemma}
\begin{proof}
Verdier biduality for the constant sheaf identifies the scalar
action with an equivalence
\[
 \underline{\F}_T\xrightarrow{\sim}
 R\mathcal{H}om_T(\omega_T,\omega_T).
\]
Since $C,D$ are perfect, tensor--Hom adjunction therefore gives
\[
 R\mathcal{H}om_T(\omega_T\otimes C,\omega_T\otimes D)
 \simeq
 \underline{R\operatorname{Hom}_{\F}(C,D)}_T.
\]
Taking derived sections proves the first equivalence in
\eqref{eq:parameter-chain-adjunction}. Its inverse is induced
by tensoring morphisms between constant complexes with
$\omega_T$, so it preserves composition and its coherences.

Admissible parameter spaces are Hausdorff, paracompact, and
locally contractible. The natural comparison between derived
sections of the constant sheaf and singular cochains is
multiplicative, including its higher coherences; see
\cite[Theorem~2.16]{ChataurCirici}. Tensoring that comparison
with a perfect coefficient complex gives
\[
 R\Gamma(T,\underline M_T)
       \simeq\operatorname{Hom}_{\F}(S_T,M),
 \qquad M=R\operatorname{Hom}_{\F}(C,D).
\]
The extension to $M$ and to the coefficient pairing
$\operatorname{Hom}(D,E)\otimes\operatorname{Hom}(C,D)
\to\operatorname{Hom}(C,E)$ uses finite sums, shifts, cones,
and the Koszul rule. Tensor--Hom adjunction proves
\eqref{eq:parameter-chain-adjunction}.

For $u:T'\to T$, Verdier duality identifies extraordinary
pullback of the displayed relative morphisms with ordinary
pullback of the corresponding constant-coefficient morphisms.
The latter is singular-cochain pullback, hence is represented
under adjunction by precomposition with $u_*$. The chain-map
identity is
\[
 d_D\Theta_\alpha(c,x)
 =\Theta_\alpha(\partial c,x)
                   +(-1)^n\Theta_\alpha(c,d_Cx),
\]
which proves \textup{(i)}, including the endpoint formula.

Composition of constant-coefficient derived morphisms is the
cup product followed by composition of their coefficient maps.
The singular-cochain cup product is induced by $\Delta_T$;
its adjoint is exactly \eqref{eq:parameter-chain-composition}.
Exterior products instead use the chain cross product and the
Koszul interchange, giving \textup{(iii)}. The natural
multiplicative comparison just used is a comparison of whole
coherent diagrams, so its composition homotopies and higher
homotopies are transported with them.

For parity grading, put $\Lambda=\F[u,u^{-1}]$, with $|u|=2$
and $du=0$, and use $\Lambda$-linear derived Hom. The scalar-action
equivalence also holds after periodizing $\omega_T$. To check this
locally, choose a locally finite triangulation of $T$ and an open
star meeting only finitely many simplices. On this star, the constant
sheaf and its dualizing complex have bounded finite-dimensional
cellular data. Derived Hom between these data is computed by the
finite resolution indexed by chains of proper face inclusions
$\sigma_0\subsetneq\cdots\subsetneq\sigma_r$ of simplices,
with $r\ge0$.
Extension of coefficients from $\F$ to $\Lambda$ commutes with its
finite sums, cones, and finite-dimensional linear Hom. Applying
this to the scalar-action equivalence above, and then to stars in subdivisions
to compute stalks, gives
\[
 \underline\Lambda_T\xrightarrow{\sim}
 R\mathcal Hom_{\Lambda,T}
       (\operatorname{Per}_2\omega_T,
        \operatorname{Per}_2\omega_T).
\]
Here $R\mathcal Hom_{\Lambda,T}$ denotes derived sheaf Hom for
sheaves of dg $\Lambda$-modules. The preceding tensor--Hom argument
therefore applies to perfect periodic $C,D$, without choosing an
integer-graded lift of their relative morphism.

For the global comparison, let $d$ be the dimension of the chosen
triangulation, let $\mathcal C_p$ be its set of $p$-simplices,
choose an orientation of each simplex, and let $R_T$ be the
cellular chain complex in cohomological grading:
$R_T^{-p}=\bigoplus_{\sigma\in\mathcal C_p}\F$, for $0\le p\le d$,
with differential the cellular boundary. For a perfect periodic
coefficient complex $M$, choose a model with finite-dimensional
terms in both parities. Cellular cochains compute
$R\Gamma(T,\underline M_T)$, and, for $\epsilon\in\{0,1\}$, their terms are
\[
 \operatorname{Hom}_{\Lambda}
      (\operatorname{Per}_2R_T,M)^{\bar\epsilon}
   =\prod_{p=0}^{d}\prod_{\sigma\in\mathcal C_p}
                  M^{\overline{\epsilon-p}}.
\]
Thus there are finitely many cellular degrees, although each set
$\mathcal C_p$ may be infinite; cochains use products over its
cells. The comparison of cellular and singular chains gives
$\operatorname{Per}_2R_T\simeq\operatorname{Per}_2S_T$.
It identifies this cellular calculation with the two-periodic
tensor--Hom adjunction using $\operatorname{Per}_2S_T$.
Singular chains preserve the functorial restriction and product
comparisons above. No interchange of direct-sum periodization
with an arbitrary infinite product is used.
\end{proof}

\section{Collision families}
\label{sec:fusion-preparations}
\label{lp:sec:collision-strategy}\label{lp:sec:pairpants-bm}
In this section we specify the families modelling collisions between loops that we will use to define the fusion product.

\subsection{Volume compatibility for loop families}
\label{conv:loop-volume-compatibility}
Throughout this section, fix the polarized SNC compactification
$(M,D,A)$ of Subsection~\ref{subsec:polarized-compactifications},
with $X=M\setminus D$ smooth affine over $\C$ and
$A=\sum_\alpha k_\alpha D_\alpha$ ample. Fix the affine
embedding defined by the chosen very ample multiple $rA$ and the
resulting pole schemes $F_P\mathscr LX$ of
\eqref{eq:pole-stage-scheme}. Write
\[
 \iota_{P,P'}:F_P\mathscr LX\hookrightarrow F_{P'}\mathscr LX
 \qquad(P\le P')
\]
for their closed inclusions. In integer grading, also fix the
nowhere-vanishing algebraic volume form
$\Omega\in\Gamma(X,\Omega_X^n)$, $n=\dim X$, and its
pole-stage volume indices from
Subsection~\ref{subsec:volume-index-preparations}. The parity-graded
construction does not require $\Omega$.

When $\Omega$ is given, a relative volume-index normalization on a presented family $S\to T$
is \emph{compatible with $\Omega$} if the following condition holds.
Whenever a specified base change $T_0\to T$ has a product
chart in the chosen finite-factor presentation diagram of the form
\[
 U\simeq T_0\times\prod_{i=1}^k U_i,
 \qquad U_i\longrightarrow F_{P_i}\mathscr LX,
 \qquad \pi_i:U_i\longrightarrow Y_i,
\]
where each \(U_i\to F_{P_i}\mathscr LX\) is a chosen nice étale chart, with finite-factor projection \(\pi_i:U_i\to Y_i\), the product finite factor $Y=T_0\times\prod_iY_i$ has
\begin{equation}\label{eq:loop-product-volume-index}
 \delta_{Y/T_0}=\sum_{i=1}^k\operatorname{pr}_i^*\delta_{Y_i}.
\end{equation}
Here $\delta_{Y_i}$ is the absolute pole-stage volume index of
Proposition~\ref{prop:finite-factor-normalization}.
If the finite factor is enlarged by retaining $e<\infty$ further
coordinates of the chart, its volume index increases by $e$, as in
Definition~\ref{def:relative-volume-index}. Base coordinates contribute
zero. In particular $\delta_{Y/\mathrm{pt}}=\delta_Y$ on a pole-stage
finite-factor presentation. Equation~\eqref{eq:loop-product-volume-index} is a
requirement on the normalization, not an existence assertion.
Proposition~\ref{prop:direct-relative-volume-indices} will construct
such a normalization for the collision families.

\subsection{Collision families}
\label{subsec:factorization}

We recall the factorization construction of \cite{KV} and use it
to define the collision families. For affine $X$, the
full loop space $\sL X$ is denoted
$\widetilde{\mathcal L}(X)$ in \cite{KV}.


\begin{definition}[The multipoint loop functor]
\label{def:kv-multipoint-functor}
Fix the \emph{source} curve $\A^1_z$, and let $I$ be a nonempty
finite set. For a $\C$-algebra $R$ and marked points
$a=(a_i)_{i\in I}\in R^I$, set
\begin{equation}\label{eq:kv-source-algebras}
 \begin{split}
 D_a(z)&=\prod_{i\in I}(z-a_i),\\
 \mathcal O_a&=\varprojlim_N R[z]/(D_a^{N+1}),\qquad
 \mathcal K_a=\mathcal O_a[D_a^{-1}].
 \end{split}
\end{equation}
The \emph{multipoint loop functor} is
\begin{equation}\label{eq:kv-multipoint-functor}
 R\longmapsto
 \{(a,\gamma):a\in R^I,\quad
       \gamma:\mathcal O(X)\longrightarrow\mathcal K_a
                    \text{ a $\C$-algebra homomorphism}\}.
\end{equation}
Its map to $\A^I$ sends $(a,\gamma)$ to $a$; a homomorphism
of test algebras acts by applying it to the marked points and
the Laurent coefficients. We write $\mathscr L_I X$ for this
functor and for its representing ind-scheme, whose existence is
included in Theorem~\ref{thm:kv-factorization} below.
\end{definition}

\begin{definition}[Source algebras and Laurent coefficient schemes]
\label{def:kv-coefficient-presentation}
We define the Laurent coefficient schemes of $\mathscr L_I X$
over $\A^I$. Put $k=|I|$. Division by the monic polynomial
$D_a$ gives isomorphisms of the {\it completed source} (respectively, {\it punctured-source}) {\it algebra} 
\begin{equation}\label{eq:kv-monic-source}
 \mathcal O_a\cong R\llbracket q_0 \rrbracket[z]/(D_a(z)-q_0),\qquad
 \mathcal K_a\cong R(\!( q_0 )\!)[z]/(D_a(z)-q_0),
\end{equation}
where $q_0$ corresponds to $D_a(z)$. Indeed, for every $N\ge1$
substitution $q_0=D_a(z)$ gives
\[
 R[q_0,z]/(q_0^N,D_a(z)-q_0)\cong R[z]/(D_a^N).
\]
Monic division makes this a free $R[q_0]/(q_0^N)$-module
with basis $1,z,\ldots,z^{k-1}$. The transition maps preserve
this basis, so taking the inverse limit gives the first
isomorphism in \eqref{eq:kv-monic-source}; inverting $q_0$
gives the second. This argument applies to every test algebra
$R$, without a reducedness or Noetherian hypothesis.

The coefficient construction is performed over $T=\A^I$,
with coefficient ring $\mathcal O(T)=\C[a_i:i\in I]$, in place
of $\C$ in Definition~\ref{def:laurent-coefficients}.
The $a_i$ are base parameters; the variables introduced below
are Laurent coefficients of the coordinate functions of the map.
Choose an affine embedding of $X$ with coordinate functions
$g_1,\ldots,g_L$
and defining equations $f_1,\ldots,f_c$. For a uniform bound
$m\ge0$, write the universal coordinate series as
\[
 G_\ell=\sum_{\nu=0}^{k-1}z^\nu
                 \sum_{n\ge-m}c_{\ell,\nu,n}q_0^n.
\]
Substitute these series into the $f_\alpha$, reduce modulo
$D_a(z)-q_0$ to degree less than $k$ in $z$, and let $J_{X,m}$
be the ideal generated by every Laurent coefficient of every
resulting basis component. Set
\[
 C_{I,m}(X)=
 \C[a_i,c_{\ell,\nu,n}:i\in I,\ 1\le\ell\le L,
                    \ 0\le\nu<k,\ n\ge-m]/J_{X,m}.
\]
The affine scheme $\Spec C_{I,m}(X)\to\A^I$ is the
\emph{Laurent coefficient scheme of uniform bound $m$} for the
chosen affine generators. Each relation is a polynomial in
finitely many coefficients: the input series have a common
lower Laurent bound, each defining polynomial has finite degree,
and monic reduction introduces only nonnegative powers of $q_0$.
Thus only finitely many input exponents contribute to a fixed
output coefficient. This scheme represents exactly
the pairs $(a,\gamma)$ with
$g_\ell(\gamma)\in q_0^{-m}\mathcal O_a$ for every $\ell$;
the morphism to $\A^I$ is induced by the coordinate functions
$a_i$. There are generally infinitely many coefficient variables
and relations, so this scheme need not be of finite type.
Increasing $m$ gives closed coefficient-stage inclusions, and
their union represents $\mathscr L_I X$. This is the coefficient
description of \cite[Proposition~2.7.1]{KV}, applied to the
chosen affine generators. 
\end{definition}

For a surjection $\pi:I\twoheadrightarrow J$, write
$I_j=\pi^{-1}(j)$. Let $U_\pi\subset\A^I$ be the open locus on
which $a_i-a_{i'}$ is invertible whenever $\pi(i)\ne\pi(i')$;
points in the same group may coincide. Let
$\Delta_\pi:\A^J\hookrightarrow\A^I$ repeat the coordinate
$b_j$ at every $i\in I_j$.

\begin{theorem}[
{\cite[Theorem~2.2.4]{KV}}]
\label{thm:kv-factorization}
Let $X$ be a smooth affine complex variety. The functors
\eqref{eq:kv-multipoint-functor} are represented by ind-schemes
$\mathscr L_I X\to\A^I$. For a singleton,
$\mathscr L_{\{1\}}X\cong\A^1\times\mathscr LX$, using
$t=z-a_1$ as the local source parameter.
For every surjection $\pi:I\twoheadrightarrow J$ there are
canonical isomorphisms
\begin{equation}\label{eq:kv-factorization}
 \begin{aligned}
 \left.\mathscr L_I X\right|_{U_\pi}
   &\cong
   \left.\left(\prod_{j\in J}\mathscr L_{I_j}X\right)
                                                       \right|_{U_\pi},\\
 \Delta_\pi^*\mathscr L_I X&\cong\mathscr L_J X.
 \end{aligned}
\end{equation}
In the first line each factor is pulled back along the projection
to its group's coordinates; the second line is scheme-theoretic
base change. These isomorphisms are natural under relabeling and
satisfy the following identities: diagonal identifications compose
under successive surjections; separating groups and then separating
their subgroups agrees with separating all the subgroups at once;
and diagonal restriction within separated groups commutes with
their separation. The families are formally smooth over $\A^I$.
\end{theorem}
\begin{proof}
The source-algebra description in
Definition~\ref{def:kv-coefficient-presentation} identifies our
functor with the full affine loop functor of
\cite[Definition~2.3.2(b) and Proposition~2.7.1(a)--(b)]{KV}.
Its representing coefficient stages and their closed inclusions
were constructed in that definition. For a singleton, the
substitution $t=z-a_1$ gives $\mathcal O_a=R\llbracket t \rrbracket$ and
$\mathcal K_a=R(\!( t )\!)$, proving the stated product description.

We verify the factorization maps over an arbitrary test algebra
$R$. On $U_\pi$, put $D_j=\prod_{i\in I_j}(z-a_i)$, so
$D_a=\prod_{j\in J}D_j$. For $i\in I_j$ and $i'\in I_{j'}$ with $j\ne j'$,
the difference $(z-a_i)-(z-a_{i'})=a_{i'}-a_i$ is a unit.
Thus the ideals generated by the cross-group linear factors
are comaximal, and so are the product ideals $(D_j)$ and $(D_{j'})$.
For every integer $N\ge1$, the Chinese remainder theorem gives
\[
 R[z]/(D_a^N)\xrightarrow{\sim}
                    \prod_{j\in J}R[z]/(D_j^N).
\]
Taking inverse limits identifies $\mathcal O_a$ with the
product of the group completions. In the $D_j$-adic completion,
every $D_{j'}$ with $j'\ne j$ is a unit: its inverse modulo
$D_j$ lifts by a convergent geometric series. Localizing
therefore gives
$\mathcal K_a\simeq\prod_{j\in J}\mathcal K_{a|I_j}$.
Maps from $\mathcal O(X)$ into this finite product give the
separation isomorphism.

On the diagonal write $a_i=b_{\pi(i)}$, put $e_j=|I_j|$,
$e=\max_j e_j$, and set
\[
 D_J=\prod_{j\in J}(z-b_j),\qquad
 D_I=\prod_{j\in J}(z-b_j)^{e_j}.
\]
The divisibilities $D_J\mid D_I\mid D_J^e$ give
\[
 (D_J^{eN})\subseteq(D_I^N)\subseteq(D_J^N)
                         \qquad(N\ge1).
\]
Thus their powers are cofinal, their completions are canonically
identified, and inverting either polynomial gives the same
punctured completion. This argument allows the $b_j$ themselves
to coincide and proves the diagonal isomorphism.

The diagonal comparisons are induced by the identity on
$R[z]$, and the source-algebra separation maps are its
projections to the group completions. Successive groupings
therefore induce
the same projections to each final group; these projections
also commute with the diagonal comparisons. Relabeling permutes
the factors. This proves all the stated compatibility identities.
The constructions are natural in every $R$, including nonreduced
$R$, so they identify the represented functors over the specified
bases and hence give the asserted scheme-theoretic base changes.
This argument evaluates completed sources over each test algebra;
it does not identify their change of coefficients with an ordinary
tensor product of completed source algebras.

Finally, let $R\to R_0=R/\mathfrak n$ be a surjection with
$\mathfrak n^N=0$, and fix lifted marks $a\in R^I$, with image
$a_0\in R_0^I$. In the basis $1,z,\ldots,z^{|I|-1}$ of
Definition~\ref{def:kv-coefficient-presentation}, the map
$\mathcal K_a\to\mathcal K_{a_0}$ is coefficientwise surjective.
Its kernel consists of expansions whose coefficients belong to
$\mathfrak n$; its $N$th power is zero, since polynomial
reduction preserves membership in $\mathfrak n^N$.
This remains valid when $\mathfrak n$ is not finitely generated.
Smoothness and affineness of $X$ imply that $\mathcal O(X)$
is formally smooth over $\C$. Hence every
$\mathcal O(X)\to\mathcal K_{a_0}$ lifts to
$\mathcal O(X)\to\mathcal K_a$. The finitely many chosen affine
generators of $\mathcal O(X)$ have a common finite Laurent
bound after lifting, possibly larger than their original bound.
This proves formal smoothness of the full families over $\A^I$;
it does not assert formal smoothness of their fixed pole stages.
\end{proof}


Use the affine generators $g_\ell$ of the embedding by $rA$
fixed in Subsection~\ref{conv:loop-volume-compatibility},
as in \eqref{eq:pole-stage-scheme}. Thus the
one-point pole stage $F_P\mathscr LX$ requires
$t^{rP}g_\ell(\gamma)\in R\llbracket t \rrbracket$ for every $\ell$.

\begin{definition}[Pole-bounded multipoint collision family]
\label{def:kv-bounded-family}
Let $I$ be a nonempty finite set and
$\mathbf P=(P_i)_{i\in I}\in\Z_{\ge0}^{I}$. Put
\[
 H=\prod_i(z-a_i)^{rP_i},\qquad m=\max_i(rP_i).
\]
Use the coefficient algebra $C_{I,m}(X)$ and the universal
coordinate series $G_\ell$ of
Definition~\ref{def:kv-coefficient-presentation}, with the fixed
affine generators $g_\ell$. Let $J_H\subset C_{I,m}(X)$ be
generated by the negative $q_0$-coefficients of every basis
component of $HG_\ell$, after reduction modulo $D_a(z)-q_0$.
Define the pole-bounded multipoint loop family
\[
 \mathcal Z_{I,\mathbf P}
       :=\Spec\bigl(C_{I,m}(X)/J_H\bigr)\longrightarrow\A^I,
\]
where the projection records the $a_i$.
Equivalently, its $R$-points are exactly the pairs
$(a,\gamma)\in\mathscr L_I X(R)$ satisfying
\begin{equation}\label{eq:kv-marked-poles}
 \left(\prod_{i\in I}(z-a_i)^{rP_i}\right)
                  g_\ell(\gamma)\in\mathcal O_a
       \quad\text{for every }\ell.
\end{equation}
Indeed, $H\mid D_a^m$ puts every such $g_\ell(\gamma)$ in
$D_a^{-m}\mathcal O_a$, and membership of $Hg_\ell(\gamma)$
in $\mathcal O_a$ is exactly the vanishing imposed by $J_H$.
There are finitely many additional negative-coefficient
equations, although the ambient Laurent coefficient scheme
$\Spec C_{I,m}(X)$ need not be of finite type.
\end{definition}

\begin{lemma}
\label{lem:kv-pole-compatibility}
For the families of Definition~\ref{def:kv-bounded-family},
increasing the $P_i$ gives finitely presented closed pole
inclusions. For a surjection $\pi:I\twoheadrightarrow J$, put
$I_j=\pi^{-1}(j)$ and $P'_j=\sum_{i\in I_j}P_i$.
The isomorphisms of Theorem~\ref{thm:kv-factorization} restrict to
\begin{equation}\label{eq:kv-bounded-factorization}
 \begin{aligned}
 \left.\mathcal Z_{I,\mathbf P}\right|_{U_\pi}
 &\cong
 \left.\left(\prod_{j\in J}
       \mathcal Z_{I_j,\mathbf P|_{I_j}}\right)
                                                        \right|_{U_\pi},\\
 \Delta_\pi^*\mathcal Z_{I,\mathbf P}
 &\cong\mathcal Z_{J,\mathbf P'}.
 \end{aligned}
\end{equation}
All maps are compatible with pole inclusions, relabeling, and
successive separations and collisions.
\end{lemma}
\begin{proof}
Fix $\mathbf P\le\mathbf P'$ componentwise, put
$m'=\max_i(rP'_i)$, and view the larger family in
$\Spec C_{I,m'}(X)$. Write $k=|I|$ and
$H=\prod_i(z-a_i)^{rP_i}$ for the lower marked bound.
Reduction of each $Hz^\nu$, $0\le\nu<k$, modulo
$D_a(z)-q_0$ has coefficients in $\C[a_i:i\in I][q_0]$.
Hence no basis component of $HG_\ell$ has a Laurent exponent
below $-m'$. The lower bound is therefore imposed by its
coefficients at $q_0^{-m'},\ldots,q_0^{-1}$, at most $Lkm'$
equations in total. Each is polynomial in finitely many
coefficient variables. By the functorial description in
Definition~\ref{def:kv-bounded-family}, their vanishing defines
exactly $\mathcal Z_{I,\mathbf P}$ inside
$\mathcal Z_{I,\mathbf P'}$. Its ideal is finitely generated,
so the inclusion is a finitely presented closed immersion.

On $U_\pi$, the completed source decomposes by groups. On the
factor belonging to $I_j$, every $z-a_i$ with $i\notin I_j$
is a unit, so the condition reduces to precisely the bound for
that group. Under $\Delta_\pi$, all marked factors in $I_j$
become $z-b_j$, and their exponents add to $rP'_j$. This proves
both identities over arbitrary test algebras, hence for the
scheme-theoretic base changes. The compatibility identities
restrict from Theorem~\ref{thm:kv-factorization}.
\end{proof}

\begin{definition}[Collision families and data]
\label{def:geometric-collision-data}
For every nonempty finite set $I$ and
$\mathbf P\in\Z_{\ge0}^{I}$, consider the pole-bounded family
$\mathcal Z_{I,\mathbf P}\to\A^I$ of
Definition~\ref{def:kv-bounded-family}.
For a surjection $\pi:I\twoheadrightarrow J$, put
\[
 \operatorname{Conf}_J
   :=\{(b_j)\in\A^J:b_j\ne b_{j'}\text{ for }j\ne j'\},
 \qquad C_\pi:=\Delta_\pi(\operatorname{Conf}_J)\subset\A^I.
\]
The locally closed subschemes $C_\pi$ are the
\emph{collision strata}: $a_i=a_{i'}$ holds precisely when
$\pi(i)=\pi(i')$. Put $P'_j=\sum_{i\in\pi^{-1}(j)}P_i$.
The identifications of Lemma~\ref{lem:kv-pole-compatibility}
give
\[
 \left.\mathcal Z_{I,\mathbf P}\right|_{C_\pi}
       \cong C_\pi\times\prod_{j\in J}F_{P'_j}\mathscr LX,
\]
using $z-b_j$ as the local parameter at the $j$th marked point.

On the locus where the marked points $a_i$ are pairwise
distinct, the separation identification describes a point
of the family by loops
\[
 \gamma_i\in F_{P_i}\mathscr LX,\qquad i\in I,
\]
written in the local source coordinates $t_i=z-a_i$.
We call these the \emph{input loops}.
When all marked points coincide at $b$, the diagonal
identification describes a point of the collided fiber by
a single loop
\[
 \gamma_{\mathrm{out}}\in
 F_{\sum_{i\in I}P_i}\mathscr LX,
 \qquad t=z-b,
\]
called the \emph{output loop}.
The separation and total-diagonal identifications are
correspondingly called the \emph{input} and \emph{output
identifications}.

For a partial collision specified by $\pi:I\twoheadrightarrow J$,
the factor $F_{P'_j}\mathscr LX$ is the output for the group
$I_j=\pi^{-1}(j)$. These group outputs can serve as inputs
to a subsequent collision. Thus the terminology is relative
to the specified collision step. It does not assert that
an arbitrary tuple of input loops determines a unique
specialization in the collided fiber.

The \emph{collision data} consist of these pole-bounded
families, their closed pole inclusions, the relabeling
isomorphisms, and the separation and scheme-theoretic diagonal
identifications of Lemma~\ref{lem:kv-pole-compatibility}, with
pole bounds added on each diagonal group, together with the
binary restriction described below, with the compatibility
identities for successive groupings, further separation, and
diagonal restriction in that lemma.
The \emph{collision diagram associated with these data}
consists of these families, maps, and identifications, together
with the open complements of the closed pole and diagonal
inclusions, the collision-stratum restrictions obtained from
the separation opens and diagonals, and finite composites and
scheme-theoretic fiber products formed among these maps.
Thus the diagram organizes the collision data and
these induced constructions; it involves no additional independent
choices of families or identifications.

For $I=\{1,2\}$ and $P,Q\ge0$, write
\[
 \mathcal Z_{P,Q}
   :=\A^1_s\times_{\A^{\{1,2\}}}
                 \mathcal Z_{\{1,2\},(P,Q)},
 \qquad s\longmapsto(s,-s).
\]
Include this 1-parameter binary-slice restriction in the collision diagram,
together with its closed pole inclusions and its input and output
identifications over $\Gm$ and $\{0\}$, respectively.
The diagonal and separation open pull back to precisely these
two base loci.

A \emph{collision family} means a family occurring in the
collision diagram. In particular, the families
$\mathcal Z_{P,Q}\to\A^1_s$ are the
\emph{1-parameter binary collision families}; their specified pole maps
and fiber identifications form the \emph{1-parameter binary collision data}.
Their open and closed fibers are displayed in
Example~\ref{ex:kv-binary-family}.

For any other morphism $b:T'\to T$ of complex schemes and
any collision family $S\to T$, call
$T'\times_T S\to T'$ the \emph{pullback of that collision
family along $b$}. Such a pullback does not become an object
of the collision diagram by this terminology. In particular,
arbitrary point inclusions are not additional specified maps;
their pullbacks are simply fibers of the corresponding family.
The unbounded families $\mathscr L_I X\to\A^I$ are the
filtered unions of the pole-bounded collision families,
and their base changes are described as pullbacks of these
filtered unions.
\end{definition}

\begin{definition}[Fusion-ready families and data]
\label{def:fusion-ready-data}
A system of \emph{fusion-ready data} is the collision diagram of Definition~\ref{def:geometric-collision-data}
equipped with the following finite-factor presentation and comparison data.
\begin{enumerate}[(i),wide=0pt]
\item For every family $f:S\to T$ in the diagram, a finite nice
\'etale atlas and finite-factor maps over $T$,
\[
 u_\alpha:U_\alpha\longrightarrow S,\qquad
 U_\alpha\cong Y_\alpha\times\A^{J_\alpha},\qquad
 p_\alpha:U_\alpha\longrightarrow Y_\alpha,\qquad
 f u_\alpha=g_\alpha p_\alpha,
\]
where $g_\alpha:Y_\alpha\to T$, $Y_\alpha$ is affine of finite
type over $\C$, and $J_\alpha$ is finite or countable.
Also choose a finite filtration over $T$
\[
 S=S_0\supset S_1\supset\cdots\supset S_m=\varnothing,
 \qquad S_a\setminus S_{a+1}\cong Y_a\times\A^{J_a},
\]
where each $Y_a$ is affine of finite type over $\C$, each
$J_a$ is finite or countable, and all closed immersions in the
filtration are finitely presented. Specify maps $g_a:Y_a\to T$
such that, under these finite-factor presentations,
\[
 f|_{S_a\setminus S_{a+1}}=g_a\operatorname{pr}_{Y_a}.
\]

\item For every finite subdiagram and every finite list of its
\'etale overlap maps, choose compatible finite-factor presentations in
the sense of Subsection~\ref{subsec:relative-volume-preparations}.
Explicitly, for each specified arrow $h:S\to S'$ over
$b:T\to T'$, choose, on compatible chart refinements, maps $\widetilde h:U\to U'$
and $h_Y:Y\to Y'$ satisfying
\[
 u'\widetilde h=h u,\qquad
 p'\widetilde h=h_Yp,\qquad
 g'h_Y=bg.
\]
A specified closed inclusion and its open complement have the
finite-factor presentation of Subsection~\ref{subsec:volume-index-preparations}:
\[
 (Z\hookrightarrow Y\hookleftarrow Y\setminus Z)\times\A^J;
\]
a specified base change $b:T'\to T$ has the induced finite-factor presentation
\[
 U\times_TT'\cong(Y\times_TT')\times\A^J.
\]
For nonaffine factors, use the finite affine open refinements of
Subsection~\ref{subsec:relative-volume-preparations}. In integer
grading the relative volume indices restrict along these open
immersions, whose relative dimension is zero.
Include finite-factor presentations of the two isomorphisms
\eqref{eq:kv-bounded-factorization} and the 1-parameter binary-slice map
$\A^1_s\to\A^2$, $s\mapsto(s,-s)$.
For every bijection $\sigma:I\xrightarrow{\sim}I'$, also include
finite-factor presentations of
\[
 \begin{aligned}
 b_\sigma:\A^I&\xrightarrow{\sim}\A^{I'},
       &(b_\sigma(a))_{i'}&=a_{\sigma^{-1}(i')},\\
 r_\sigma:\mathcal Z_{I,\mathbf P}
       &\xrightarrow{\sim}\mathcal Z_{I',\sigma_*\mathbf P},
       &(\sigma_*\mathbf P)_{i'}&=P_{\sigma^{-1}(i')},
 \end{aligned}
\]
where $r_\sigma$ is the relabeling isomorphism over $b_\sigma$.
Require their descended maps to satisfy the stated identities
among these arrows and their composites. Separation uses product charts;
on $C_\pi$ their finite factors are
\[
 U=C_\pi\times\prod_{j\in J}U_j,
 \quad U_j=Y_j\times\A^{J_j}\longrightarrow F_{P'_j}\mathscr LX,
 \qquad Y=C_\pi\times\prod_{j\in J}Y_j.
\]
Comparisons with other charts are through common refinements.
In particular, two finite factors $p:U\to Y$ and $p':U\to Y'$
on a common chart have a common dominating factor
\[
 q:U\to W,\qquad p=a q,\quad p'=a'q,
 \qquad a:W\to Y,\quad a':W\to Y'.
\]
Larger finite subdiagrams and different choices are compared by
further such refinements, with the descended maps satisfying
the same equalities among their composites. As in \emph{Compatible
finite-factor presentations} in Subsection~\ref{subsec:coherent-descent},
finite factors are chosen for occurrences of charts in finite
composable strings; repeated occurrences are compared by
refinements. No single finite factor for the entire infinite
diagram is required.

\item In integer grading, choose locally constant functions
$\delta_{Y/T}:Y\to\Z$ on this finite-factor presentation diagram satisfying
Definition~\ref{def:relative-volume-index} and
Subsection~\ref{conv:loop-volume-compatibility}. In particular,
on the product finite factors in (ii) require
\[
 \delta_{C_\pi\times\prod_jY_j/C_\pi}
      =\sum_j\operatorname{pr}_j^*\delta_{Y_j}.
\]
In parity grading no volume form or integer volume-index
normalization is part of the data; use the two-periodic smooth
purity comparisons on the same finite-factor presentations.
\end{enumerate}
A family with its induced data is a \emph{fusion-ready
family}. In particular, the 1-parameter binary slice carries compatible
finite-factor presentations of the localization diagram
\[
 (\mathcal Z_{P,Q})_0\hookrightarrow\mathcal Z_{P,Q}
                         \hookleftarrow\mathcal Z_{P,Q}|_{\Gm},
\]
its closed pole inclusions, and its fiber identifications in
Example~\ref{ex:kv-binary-family}; in integer grading these
finite-factor presentations carry the relative volume indices in (iii).
The finite-factor presentation and normalization requirements above are part
of the data. Their existence for the collision diagram is proved
in Proposition~\ref{lem:all-multipoint-models}, using the relative
finite-factor and volume-index constructions of the next two sections.
\end{definition}

\begin{example}[The 1-parameter binary collision family]
\label{ex:kv-binary-family}
The slice $s\mapsto(s,-s)$ makes the source and the two
fiber identifications explicit. For a complex test algebra
$R$ and $s\in R$, Definition~\ref{def:kv-coefficient-presentation}
with $D_a(z)=z^2-s^2$ gives the first line below; the  identities (\ref{eq:kv-bounded-factorization}) of Lemma~\ref{lem:kv-pole-compatibility}
give the second:
\begin{equation}\label{eq:kv-binary-family}
 \begin{gathered}
 \mathcal O_s=R\llbracket q_0 \rrbracket[z]/(z^2-s^2-q_0),\qquad
 \mathcal K_s=\mathcal O_s[q_0^{-1}],\\
 \left.\mathcal Z_{P,Q}\right|_{\Gm}
       \cong F_P\mathscr LX\times F_Q\mathscr LX\times\Gm,
 \qquad
 (\mathcal Z_{P,Q})_0
       \cong F_{P+Q}\mathscr LX.
 \end{gathered}
\end{equation}
For invertible $s$, the two disk parameters are $z-s,z+s$;
the factors from the other marked point are units, so the
two pole bounds are $P,Q$. At $s=0$, the completed source
is identified with $R\llbracket z \rrbracket$ by
\[
 A(q_0)+zB(q_0)\longmapsto A(z^2)+zB(z^2),
 \qquad A,B\in R\llbracket q_0 \rrbracket.
\]
Both marked factors become $z$, giving the output bound $P+Q$.
Thus the first fiber formula is an
isomorphism over the whole punctured base $\Gm$, and the
second is the scheme-theoretic fiber over zero, over arbitrary
test algebras. For $(P,Q)\le(P',Q')$, the inherited closed pole
map restricts under these identifications to
\[
 (\iota_{P,P'}\times\iota_{Q,Q'})\times\mathrm{id}_{\Gm}
 \quad\text{and}\quad
 \iota_{P+Q,P'+Q'},
\]
respectively, by Lemma~\ref{lem:kv-pole-compatibility}.

For an affine presentation
$X=\Spec\C[x_1,\ldots,x_m]/(f_1,\ldots,f_l)$, we now write
the Laurent coefficient equations for maps
$\gamma:\Spec\mathcal K_s\to X$. Expand each coordinate
$x_i(\gamma)$ in the $R(\!( q_0 )\!)$-basis $1,z$ of $\mathcal K_s$.
Write
\[
 x_i(\gamma)=A_i(q_0)+zB_i(q_0),\qquad
 A_i,B_i\in R(\!( q_0 )\!),\qquad c=s^2+q_0,
\]
and put $\mathbf A=(A_i)_i$, $\mathbf b=(B_i)_i$.
Polynomial division gives the equations
\begin{equation}\label{eq:example-fusion-equations}
 \begin{aligned}
 f_j(\mathbf A+z\mathbf b)
 &\equiv f_{j,0}(\mathbf A,\mathbf b,c)
       +z f_{j,1}(\mathbf A,\mathbf b,c)\pmod{z^2-c},\\
 f_{j,0}(\mathbf A,\mathbf b,c)
 &=f_{j,1}(\mathbf A,\mathbf b,c)=0\qquad(1\le j\le l).
 \end{aligned}
\end{equation}
Here $f_{j,0},f_{j,1}$ are the two coefficients of the unique
remainder; after evaluation they lie in $R(\!( q_0 )\!)$. Their
vanishing means that every Laurent coefficient is zero; this is equivalent
to the tuple defining a map $\Spec\mathcal K_s\to X$. To obtain $\mathcal Z_{P,Q}$,
impose additionally the marked pole conditions on the fixed
pole generators $g_\ell$:
\[
 (z-s)^{rP}(z+s)^{rQ}g_\ell(\gamma)\in\mathcal O_s
                       \qquad\text{for every }\ell.
\]
These are the equations of Definitions~\ref{def:kv-coefficient-presentation}
and~\ref{def:kv-bounded-family} on the binary slice. 

If the affine presentation of $X$ includes a condition $g\ne0$,
with $g\in\C[x_1,\ldots,x_m]$, use
\[
 X=\Spec\bigl((\C[x_1,\ldots,x_m]/(f_1,\ldots,f_l))[g^{-1}]\bigr)
   \cong
   \Spec\C[x_1,\ldots,x_m,v]/(f_1,\ldots,f_l,vg-1).
\]
Here $v$ is an additional affine coordinate representing $g^{-1}$.
Apply the preceding Laurent coefficient construction to this
presentation, including the additional defining polynomial $vg-1$.
\end{example}

\section{Relative finite factors for collision families}
\label{sec:relative-finite-factors}
\label{sec:direct-finite-factors}

We begin prepping the collision families to become fusion-ready. The main goal of this section is to prove Theorem~\ref{thm:relative-finite-factors}, a relative version of Drinfeld's finite-factor theorem, constructing finite factors for our families.

\subsection{The relative finite-factor theorem}
\label{subsec:relative-finite-factor-theorem} 
Fix a smooth affine complex variety $X\subset\A^N$ of pure
dimension $n>0$. Let $A$ be a finitely generated complex algebra,
let $D(z)\in A[z]$ be monic of degree $k\ge1$, and, for each
$A$-algebra $R$, put
\begin{equation}\label{eq:direct-source-algebras}
 O_R=R\llbracket q \rrbracket[z]/(D(z)-q),\qquad E_R=O_R[1/q].
\end{equation}
The basis $1,z,\ldots,z^{k-1}$ identifies $O_R$ with
$R\llbracket q \rrbracket^k$. For $b\ge0$, let $S_b\to\Spec A$ be the affine
Laurent coefficient scheme representing
\[
 S_b(R)=\{\gamma\in X(E_R):
                 x_\ell(\gamma)\in q^{-b}O_R\ (1\le\ell\le N)\}.
\]
It is defined by substituting the coefficient series for the
$x_\ell$ into the equations of $X$ and equating each coefficient
to zero, as in Definition~\ref{def:laurent-coefficients}.

\begin{theorem}[Relative finite factors]
\label{thm:relative-finite-factors}
Every $S_b\to\Spec A$ above admits a finite Nisnevich atlas
of finitely presented \'etale maps
\[
 U_\alpha\longrightarrow S_b,\qquad
 U_\alpha\cong Y_\alpha\times_{\Spec A}\A^{J_\alpha}_A,
\]
where $Y_\alpha$ is affine of finite type over $A$ and
$J_\alpha$ is finite or countable. These finite-factor presentations retain
the map to $\Spec A$. The assertion is preserved by finite-type
base change and by finitely presented locally closed restrictions
of the total scheme. Each resulting family admits a finite
filtration by finitely presented closed subschemes with nice
successive differences, over the same base.

In particular, taking $A=\C[a_i:i\in I]$,
$D(z)=\prod_{i\in I}(z-a_i)$, and $q=q_0$ in
Definition~\ref{def:kv-coefficient-presentation}, every collision family $\mathcal Z_{I,\mathbf P}\to\A^I$ admits such an atlas
\begin{equation}\label{eq:direct-relative-nice-atlas}
 U_\alpha\longrightarrow\mathcal Z_{I,\mathbf P},
 \qquad
 U_\alpha\cong Y_\alpha\times_{\A^I}\A^{J_\alpha}_{\A^I},
\end{equation}
and such a finite filtration. The same conclusions hold after
any finite-type base change $T\to\A^I$ and any finitely presented
locally closed restriction.
\end{theorem}

A family $\{V_\alpha\to S\}$ is a \emph{Nisnevich covering}
if its maps are \'etale and
\[
 \forall x\in S\quad\exists\alpha,\ y\in V_\alpha:
 \quad y\longmapsto x,\qquad\kappa(x)\xrightarrow{\sim}\kappa(y).
\]
Here $\kappa(x)$ and $\kappa(y)$ are the residue fields at the
indicated points. A \emph{finite Nisnevich atlas} means such a covering with
finitely many charts. Over a Noetherian scheme, the displayed
condition is equivalent to the existence of a finite closed
filtration
$S=S_0\supset\cdots\supset S_m=\varnothing$
on whose successive differences the map
$\coprod_\alpha V_\alpha\to S$ has sections. The phrase
\emph{Nisnevich-locally} means after a covering of this kind;
see \cite[Section~3.2.5]{Drinfeld}.

\begin{remark}
\label{rem:direct-recovers-drinfeld}
Theorem~\ref{thm:relative-finite-factors} extends the finite-factor
conclusion of Drinfeld's theorem \cite[Section~6]{Drinfeld}, in our setting, from a field of coefficients to a finitely
generated complex coefficient ring. 
The theorem recovers
Theorem~\ref{thm:drinfeld-finite-factors}: set $A=\C$,
$D(z)=z$, so $E_R=R(\!( q )\!)$.
If $S\subset\mathscr LX$ is closed and quasi-compact, the finitely
many affine coordinates of $X$ have a common Laurent bound on
$S$, hence $S\subset S_b$ for some $b$. If $S$ is reasonable,
the ideal of this closed immersion is finitely generated, so
$S\hookrightarrow S_b$ is finitely presented. Applying the
restriction assertion gives the nice atlas and finite nice
filtration stated in Theorem~\ref{thm:drinfeld-finite-factors}.
No polarization or volume form is used in this deduction.
\end{remark}

The proof follows the strategy outlined in
\cite[Section~6.7, p.~293]{Drinfeld}: at its core lies the decomposition (\ref{eqTm}), whose generalization via Tate modules can be applied to tangent modules to obtain a linearized version of the finite-factor theorem; then we use an implicit-function
argument to construct finite-factor charts. To implement this, for maps given by
Laurent series we first record finitely many coefficients and
require the equations of $X$ to hold up to a sufficiently high
power of $q$. Smoothness of $X$ lets us correct such an approximate
map to an exact map while changing only higher coefficients.
Near this exact map, sufficiently high-order changes in tangent
directions can be chosen freely; successive corrections determine
the remaining changes needed to stay in $X$. After a suitable
Nisnevich cover of the finite approximation scheme, these tangent
parameters can be chosen as independent affine coordinates.
Each coefficient of the corrected map depends on only finitely
many of them. Consequently, requiring the corrected map to have
the recorded coefficients imposes equations on only finitely
many parameters. These equations define the finite factor, while
the unused parameters remain free and give the affine tail.

\subsection{Coefficient estimates and polynomial correction}

Write the ideal of the fixed embedding $X\subset\A^N$ as
$I_X=(f_1,\ldots,f_c)$.
Smoothness and affineness of $X$ give a lift
\[
 \mathcal O(X)\longrightarrow\C[x_1,\ldots,x_N]/I_X^2
\]
of the identity modulo $I_X$. Choose polynomial representatives
of the images of the classes of $\bar x_i\in \mathcal O(X)$ in $\C[x_1,\ldots,x_N]/I_X^2$. They define a polynomial map
$p:\A^N\to\A^N$ satisfying
\begin{equation}\label{eq:direct-polynomial-retraction}
 p(x)-x=G(x)f(x),\qquad
 f_\alpha(p(x))=
       \sum_{\beta,\gamma}B_{\alpha\beta\gamma}(x)
                          f_\beta(x)f_\gamma(x),
\end{equation}
for polynomial matrices $G$ and $B$.
Thus $p|_X=\mathrm{id}$ and $p$ induces a retraction
\[
 \Spec(\C[x_1,\ldots,x_N]/I_X^2)\longrightarrow X.
\]
Its iterates on $E_R^N$ are denoted
\[
 p^{\circ0}=\mathrm{id}_{E_R^N},\qquad
 p^{\circ(j+1)}=p\circ p^{\circ j}\quad(j\ge0).
\]
These iterates will ``correct'' the approximate solutions we define next.

For any finite $d\ge1$, give $E_R^d$ the linear topology with
neighborhood basis $\{q^mO_R^d:m\in\Z\}$ at zero. Thus
\[
 \begin{aligned}
 u\equiv v\pmod {q^m}
   &\ \Longleftrightarrow\ u-v\in q^mO_R^d,\\
 u_j\longrightarrow u
   &\ \Longleftrightarrow\
     \forall m\ \exists j_m\ \forall j\ge j_m:
                  u_j-u\in q^mO_R^d.
 \end{aligned}
\]
This is the separated, complete topology defined by the lattice
$O_R^d$; it is not the ideal-adic topology for $(q)\subset E_R$,
since $q$ is invertible in $E_R$. For $u=(u_1,\ldots,u_d)\in E_R^d$,
expand each component uniquely as
\[
 u_i=\sum_{j=0}^{k-1}u_{ij}(q)z^j,
 \qquad u_{ij}(q)\in R(\!( q )\!),\qquad 1\le i\le d.
\]
Define $\nu(u)$ to be the least exponent of $q$ having a nonzero
coefficient in any of the Laurent series $u_{ij}(q)$, with
$\nu(0)=+\infty$. For $0\le j<k$ and $m\in\Z$, set
\[
 [z^jq^m]u_i:=[q^m]u_{ij}(q),
\]
where $[q^m]$ extracts the coefficient of $q^m$ from a Laurent
series. Thus
$u\in q^mO_R^d$ exactly when $\nu(u)\ge m$.
The multiplication table of $O_R$ has coefficients in $A[q]$,
so $\nu(ab)\ge\nu(a)+\nu(b)$. All congruences and limits below
use this topology.

For the remainder of this subsection, $R$ ranges over $A$-algebras.
The fixed embedding $\iota:X\hookrightarrow\A^N$ induces
\[
 \iota_{E_R}:X(E_R)\hookrightarrow E_R^N,
 \qquad\gamma\longmapsto
       (\gamma(\bar x_1),\ldots,\gamma(\bar x_N)).
\]
Using this injection, identify $S_b(R)$ with
$\{u\in q^{-b}O_R^N:f(u)=0\}$.
For integers $b,r\ge0$, define the affine coefficient scheme
$W_{b,r}$ over $A$ by
\begin{equation}\label{eq:approximate-solution-locus}
 W_{b,r}(R):=
 \{u\in q^{-b}O_R^N:f(u)\in q^rO_R^c\}.
\end{equation}
Explicitly, in
$V_b:=\Spec A[u_{\ell,\nu,m}:1\le\ell\le N,
                   \ 0\le\nu<k,\ m\ge-b]$, write
\[
 u_\ell=\sum_{\nu=0}^{k-1}z^\nu
                   \sum_{m\ge-b}u_{\ell,\nu,m}q^m.
\]
Then $W_{b,r}\subset V_b$ is cut out by the coefficients
$[z^\nu q^m]f_\alpha(u)$ with $0\le\nu<k$ and $m<r$.
These are finitely many polynomial equations:
$f_\alpha(u)$ has a uniform lower Laurent bound, and each of
its coefficients uses finitely many $u_{\ell,\nu,m}$. 

So \(W_{b,r}\) parametrizes approximate solutions of the equations defining \(S_b\), with pole bound \(b\) and error divisible by \(q^r\). Increasing \(r\) requires a more accurate approximation. The first part of the next result says that the polynomial map $p$ can be used to correct these approximate solutions.

\begin{lemma}
\label{lem:direct-newton-retraction}
Fix a bound $b\ge0$. There is an integer $C\ge1$, independent
of $R$, with the following properties.
\begin{enumerate}[(i),wide=0pt]
\item For every integer $r>C$ and every $u\in W_{b,r}(R)$,
the limit
\[
 \rho_R(u):=\lim_{j\to\infty}p^{\circ j}(u)
\]
exists and satisfies
\[
 \rho_R(u)\in q^{-b}O_R^N,\qquad f(\rho_R(u))=0,
\]
as well as
\begin{equation}\label{eq:direct-retraction-error}
 \nu(\rho_R(u)-u)\ge\nu(f(u))-C.
\end{equation}
For each $\ell,\nu,m$, the coefficient
$[z^\nu q^m]\rho_R(u)_\ell$ is a polynomial over $A$ in
finitely many $u_{\ell',\nu',m'}$, independent of the test
algebra $R$.
For every $A$-algebra homomorphism $\varphi:R\to R'$, one has
\begin{equation}\label{eq:retraction-test-algebra-naturality}
 \varphi_*\rho_R(u)=\rho_{R'}(\varphi_*u),
\end{equation}
where $\varphi_*$ applies $\varphi$ to each coefficient.
Thus the maps $\rho_R$ define an $A$-scheme morphism
\[
 \rho_{b,r}:W_{b,r}\longrightarrow S_b,
 \qquad \rho_{b,r}|_{S_b}=\mathrm{id}_{S_b}.
\]
For $r'\ge r>C$, its restriction to $W_{b,r'}$ is
$\rho_{b,r'}$. We write $\rho(u)$ when the indices are understood.
\item If $b_0\in S_b(R)$, then
$P_{b_0}:=dp_{b_0}
=
\left(\frac{\partial p_i}{\partial x_j}(b_0)\right)_{i,j}
\in\operatorname{Mat}_N(E_R)$
satisfies $P_{b_0}^2=P_{b_0}$ and has image
\[
 M_{b_0}=b_0^*T_X\subset E_R^N.
\]
For a tangent displacement $v\in M_{b_0}$ with $\nu(v)>C$, put
\[
 \Phi_{b_0}(v)=\rho(b_0+v),\qquad
 H_{b_0}(v)=P_{b_0}(\Phi_{b_0}(v)-b_0)-v.
\]
Uniformly for $b_0\in S_b(R)$ and $v,w\in M_{b_0}$ with
$\nu(v),\nu(w)>C$, one has
\begin{align}
 \nu(H_{b_0}(v))&\ge2\nu(v)-C,\label{eq:direct-newton-quadratic}\\
 \nu(H_{b_0}(v)-H_{b_0}(w))
 &\ge\nu(v-w)+\min(\nu(v),\nu(w))-C.
 \label{eq:direct-newton-lipschitz}
\end{align}
For every fixed integer $h>C$, each Laurent coefficient of every component of \(\Phi_{b_0}(v)\) and \(H_{b_0}(v)\)
on $b_0\in S_b(R)$ and $v\in M_{b_0}\cap q^hO_R^N$ is a
polynomial over $A$ in finitely many coefficients of $b_0$ and
$v$. Both maps commute with coefficientwise application of
$A$-algebra homomorphisms.
\end{enumerate}
\end{lemma}

\begin{proof}
For every polynomial $F\in\C[x_1,\ldots,x_N]$, changing
one factor of a monomial at a time gives polynomials
$F_i\in\C[x_1,\ldots,x_N,y_1,\ldots,y_N]$ such that
\[
 F(u)-F(v)=\sum_{i=1}^N(u_i-v_i)F_i(u,v).
\]
On the ball $q^{-b}O_R^N$, polynomial evaluation has a uniform
lower $q$-bound. Fix $c_0\ge0$, independent of $R$, so that
the entries of $G,B,dp$, the polynomial difference coefficients
for $f,p,G,B,dp$, and the coefficients of the terms of degree
at least two in the Taylor expansions of $f$ and $p$ all have
order at least $-c_0$ on this ball. In particular,
\[
 \nu(p(u)-u)\ge\nu(f(u))-c_0,\qquad
 \nu(f(p(u)))\ge2\nu(f(u))-c_0.
\]
Choose $C\ge\max(1,3c_0)$. For $u\in W_{b,r}(R)$ with $r>C$,
put $u_j=p^{\circ j}(u)$. Induction gives
\[
 \nu(f(u_j))\ge 2^j(r-C)+C,\qquad
 \nu(u_{j+1}-u_j)\ge 2^j(r-C).
\]
All iterates remain in $q^{-b}O_R^N$, since every correction
has positive order. Completeness gives their limit $\rho_R(u)$,
and continuity of polynomial evaluation gives
$f(\rho_R(u))=0$. Using $\nu(f(u))$ in place of $r$ proves
\eqref{eq:direct-retraction-error}; if $f(u)=0$, then
$p(u)=u$ and the iteration is stationary.
Moreover,
\begin{equation}\label{eq:retraction-coefficient-stabilization}
 \nu(\rho_R(u)-p^{\circ j}(u))\ge 2^j(r-C).
\end{equation}
For a desired coefficient $[z^\nu q^m]$, choose $j$ with
$2^j(r-C)>m$. That coefficient of $\rho_R(u)$ is then the
corresponding coefficient of $p^{\circ j}(u)$, a fixed finite
polynomial iteration. It is therefore a polynomial over $A$
in finitely many input coefficients, uniformly on $W_{b,r}$.
Applying an $A$-algebra homomorphism to that polynomial proves
\eqref{eq:retraction-test-algebra-naturality}. The coefficient
polynomials define $\rho_{b,r}$ on affine coordinate rings;
the identities $f(\rho_R(u))=0$ make it factor through $S_b$.
The same iteration for all $r$ proves the restriction identity.

Differentiating \eqref{eq:direct-polynomial-retraction} at $b_0$
shows that $df_{b_0}P_{b_0}=0$ and that $P_{b_0}$ is the
identity on $\ker df_{b_0}=b_0^*T_X$. Hence it is the stated
idempotent. Taylor expansion gives
\[
 p(b_0+v)=b_0+P_{b_0}v+Q_{b_0}(v),
\]
where, uniformly on the bounded ball,
\[
 \nu(Q_{b_0}(v))\ge2\nu(v)-c_0,\qquad
 \nu(Q_{b_0}(v)-Q_{b_0}(w))
 \ge\nu(v-w)+\min(\nu(v),\nu(w))-c_0.
\]
For tangent vectors, $P_{b_0}v=v$. Put
$h=\min(\nu(v),\nu(w))$, $d=\nu(v-w)$, and
$x_j(v)=p^{\circ j}(b_0+v)$. Taylor expansion of $f$ gives
\[
 \nu(f(x_0(v)))\ge2h-c_0,\qquad
 \nu(f(x_0(v))-f(x_0(w)))\ge d+h-c_0.
\]
Set $r_0=2h-c_0$ and $r_{j+1}=2r_j-c_0$.
Since $h>C\ge3c_0$, simultaneous induction using
\eqref{eq:direct-polynomial-retraction} gives
\[
 \begin{split}
 \nu(f(x_j(v)))&\ge r_j,\\
 \nu(x_j(v)-x_j(w))&\ge d,\\
 \nu(f(x_j(v))-f(x_j(w)))&\ge d+r_j-h.
 \end{split}
\]
Indeed, the difference of the quadratic errors has order
at least
$\min(d+2r_j-h-c_0,d+2r_j-c_0)=d+r_{j+1}-h$.
The difference of the two increments
$G(x_j)f(x_j)$ has order at least $d+r_j-h-c_0$,
which exceeds $d$; this also preserves the second inequality.
Summing these increments and passing to the limit yields
\[
 \begin{split}
 \nu(\Phi_{b_0}(v)-b_0-v)&\ge2\nu(v)-2c_0,\\
 \nu\bigl((\Phi_{b_0}(v)-b_0-v)
       -(\Phi_{b_0}(w)-b_0-w)\bigr)&\ge d+h-2c_0.
 \end{split}
\]
The entries of $P_{b_0}$ have order at least $-c_0$ by
the initial choice of $c_0$. Multiplication by $P_{b_0}$ proves
\eqref{eq:direct-newton-quadratic}--
\eqref{eq:direct-newton-lipschitz}.
For fixed $h>C$, Taylor expansion puts $b_0+v$ in
$W_{b,2h-c_0}(R)$. Part~(i), applied to the coefficientwise
addition map $(b_0,v)\mapsto b_0+v$, gives the polynomial and
naturality assertions for $\Phi_{b_0}$. The formula for
$H_{b_0}$ gives the same assertions for $H_{b_0}$.
\end{proof}

In Lemma~\ref{lem:direct-newton-retraction}(ii)  the first map $\Phi_{b_0}(v)$ gives the polynomial correction to an exact solution from (i) of the displaced approximate solution $b_0+v$, while the second map $H_{b_0}(v)$ measures the discrepancy between its tangent projection and the original displacement. Equivalently,  correction
followed by tangent projection is expressed as $\mathrm{id}+H_{b_0}$.
The proof of Theorem~\ref{thm:relative-finite-factors} will use the estimate \eqref{eq:direct-newton-lipschitz} to invert
this map on a sufficiently small tangent lattice.

\subsection{Tangent lattices and algebraic finite factors}
The purpose of this subsection is to supply the linear parameters
for the relative finite-factor theorem. For an $A$-algebra $R$
and $b_0\in S_b(R)$, the tangent module of
Lemma~\ref{lem:direct-newton-retraction}(ii) is
\[
 M_{b_0}=\Gamma(\Spec E_R,b_0^*T_X)
        =\ker\bigl(df_{b_0}:E_R^N\longrightarrow E_R^c\bigr).
\]
For $v\in E_R^N$, the equality
\[
 f(b_0+\varepsilon v)=\varepsilon\,df_{b_0}(v)
       \quad\text{in }\bigl(E_R[\varepsilon]/(\varepsilon^2)\bigr)^c
\]
shows that $M_{b_0}$ parametrizes first-order deformations of
$b_0$ as a Laurent map into $X$, with the source and the
$A$-coefficients fixed. Such a deformation lies in the bounded
stage $S_b$ precisely when $v\in q^{-b}O_R^N$; in particular,
sufficiently high-order tangent displacements $v$ preserve this
bound. Since $X$ is smooth and $E_R$ is finite free over
$R(\!(q)\!)$, the module $M_{b_0}$ is finite projective over
$R(\!(q)\!)$, hence a Tate $R$-module by  Example \ref{exRTm} below.
A lattice in $M_{b_0}$, as defined below, thus consists of tangent displacements.
Its openness ensures that it contains every sufficiently
high-order tangent displacement, and multiplication by a
sufficiently large power of $q$ makes all its elements small
enough for the implicit-function argument.

In the following module-theoretic definitions, \(R\) denotes an arbitrary commutative ring.
A countable \emph{elementary Tate $R$-module} is a topological
module $P\oplus Q^\vee$, where $P,Q$ are countably generated
projective $R$-modules, $P$ is discrete, and
$Q^\vee=\operatorname{Hom}_R(Q,R)$ has the topology of pointwise
convergence, with $R$ discrete. A \emph{Tate module} is a
topological direct summand of an elementary Tate module, that is, a direct summand whose inclusion and projection are continuous.
These are the countable cases of the definitions in
\cite[Section~3.2.1]{Drinfeld}, where $P,Q$ may be arbitrary
projective modules.
A \emph{lattice} in a Tate $R$-module $M$ is an open
$R$-submodule $L\subset M$ such that $L/V$ is finitely generated
over $R$ for every open submodule $V\subset L$.
It is \emph{coprojective} if $M/L$ is a projective $R$-module.
A submodule of $M$ is \emph{bounded} if it is contained in a
lattice. These definitions are those of
\cite[Section~3.2.3]{Drinfeld}.

The \emph{continuous dual} of $M$ is
\[
 M^\vee:=\operatorname{Hom}^{\mathrm{cont}}_R(M,R),
\]
where $R$ has the discrete topology. Equip $M^\vee$ with
the topology whose basic neighborhoods of zero are
\[
 L^\perp
 :=\{\lambda\in M^\vee:\lambda(\ell)=0
                         \text{ for every }\ell\in L\},
\]
as $L$ ranges over the lattices in $M$.
This is the dual topology of \cite[Section~3.2.4]{Drinfeld},
where annihilators of open bounded submodules are used:
every such submodule $U\subset M$ lies in a lattice $L$, so
$L^\perp\subseteq U^\perp$, while every lattice is open and
bounded. Here $U^\perp$ denotes the continuous functionals
vanishing on $U$. The same reference gives that $M^\vee$ is
a Tate module and that evaluation induces an isomorphism
of topological $R$-modules $M\xrightarrow{\sim}(M^\vee)^\vee$.

\begin{lemma}[Coprojective lattices]
\label{lem:coprojective-lattice-criterion}
A Tate $R$-module $M$ admits a coprojective lattice if and
only if it is elementary Tate
\cite[Section~3.2.3, Remark~(ii), p.~268]{Drinfeld}.
For every coprojective lattice $L\subset M$, the modules
$L^\vee$ and $M/L$ are discrete and countably generated
projective over $R$. An $R$-linear section of $M\to M/L$
gives a topological decomposition
\[
 M\cong L\oplus(M/L).
\]
\end{lemma}

\begin{proof}
Put
\[
 H_R:=R^{(\mathbb N)}\oplus R^{\mathbb N},
\]
where the first summand is the direct sum with the discrete
topology and the second is the product with its product
topology. Every Tate module under consideration is a
continuous direct summand of $H_R$: a countably generated
projective module is a direct summand of $R^{(\mathbb N)}$,
and dualizing such a splitting gives the corresponding
continuous splitting of the product.

First let $D$ be a discrete Tate module, with continuous
maps $i:D\to H_R$ and $r:H_R\to D$ satisfying $ri=\mathrm{id}_D$.
Continuity and discreteness imply that $r$ factors through
the natural projection
\[
 \operatorname{pr}_F:H_R\longrightarrow
                 R^{(\mathbb N)}\oplus R^F
\]
for some finite $F\subset\mathbb N$. Writing
$r=\bar r\operatorname{pr}_F$, one has
$\bar r\operatorname{pr}_Fi=\mathrm{id}_D$.
Thus $D$ is a direct summand of a countably generated free
$R$-module, hence is countably generated projective.

If $M=P\oplus Q^\vee$ is elementary, then $L=Q^\vee$ is open.
Dualizing a split inclusion $Q\hookrightarrow R^{(\mathbb N)}$
gives a continuous surjection $R^{\mathbb N}\to Q^\vee$.
For every open submodule $V\subset Q^\vee$, its composite
with $Q^\vee\to Q^\vee/V$ is a continuous surjection to a
discrete module, so it factors through finitely many
coordinates of $R^{\mathbb N}$. Therefore $Q^\vee/V$ is
finitely generated. Thus $L$ is a lattice, and it is
coprojective because $M/L\cong P$ is projective.

Conversely, let $L\subset M$ be coprojective. The projective
quotient $M/L$ is discrete, so an $R$-linear section of the
quotient map is continuous. It gives the topological
decomposition $M\cong L\oplus(M/L)$, making both summands
Tate modules. By the duality recalled above, $L^\vee$ is Tate.
Since $L$ is a lattice in itself, its annihilator in
$L^\vee$ is the open zero submodule; hence $L^\vee$ is
discrete. Applying the discrete-module argument to
$L^\vee$ and $M/L$ shows that both are countably generated
projective. Applying the same biduality to $L$ identifies the
first summand with the dual of the projective module $L^\vee$;
hence this decomposition is elementary Tate.
\end{proof}

For a ring homomorphism $R\to R'$ and a Tate $R$-module $M$,
we now give an explicit construction of its
\emph{completed change of coefficients} $R'\widehat\otimes_R M$,
the completed tensor product used in
\cite[Theorem~3.4 and Section~3.5.2]{Drinfeld}.
For an elementary Tate module $E=P\oplus Q^\vee$, set
\[
 R'\widehat\otimes_R E
 :=
 (P\otimes_RR')\oplus
 \operatorname{Hom}_{R'}(Q\otimes_RR',R'),
\]
with the first summand discrete and the second endowed with
the pointwise-convergence topology.
For a general Tate module, choose a presentation
$M=\operatorname{im}(e:E\to E)$ by a continuous idempotent
on an elementary Tate module and set
\[
 R'\widehat\otimes_R M
 :=
 \operatorname{im}\!\left(
 e_{R'}:R'\widehat\otimes_R E
       \longrightarrow R'\widehat\otimes_R E
 \right),
\]
Here the induced map $e_{R'}$ is defined by the following
construction, which also proves independence of the chosen
presentation. Use the continuous direct-summand presentations
in $H_R$ from the proof of
Lemma~\ref{lem:coprojective-lattice-criterion}.
A continuous homomorphism
$H_R\to H_R$ has a block matrix
\[
 \begin{pmatrix} a&b\\ c&d\end{pmatrix},
\]
where each column of $a$ and each row of $d$ has finite
support, $b$ has only finitely many nonzero entries, and
$c$ is arbitrary. In particular, continuity of the
product-to-discrete block $b$ makes it factor through
finitely many input coordinates; the images of their basis
vectors have a finite union of supports in the direct sum.
Conversely these conditions define a continuous map.
They are preserved by composition, and each entry of a
composite is a finite sum. Applying $R\to R'$ entrywise
therefore preserves compositions, identities, and
idempotents. For maps between continuous direct summands,
extend by the chosen inclusion and projection, apply this
construction, and restrict to the images of the resulting
idempotents. Comparison maps between two presentations of
the same summand become mutually inverse comparison maps
after coefficient change, because their defining
composition identities are preserved. This proves
presentation independence and functoriality, including
compatibility with successive ring homomorphisms. On an
elementary summand it gives exactly the displayed formula
for $R'\widehat\otimes_R E$. We abbreviate the resulting
module $R'\widehat\otimes_R M$ to $M_{R'}$. 

This
construction preserves split decompositions and equalities
between composites of continuous maps, without any
flatness assumption on $R\to R'$. In particular, $M_{R'}$ is a Tate $R'$-module, and it is elementary if $M$ is elementary. Indeed, $P\otimes_RR'$ and $Q\otimes_RR'$ remain countably generated projective modules. For general $M$, the continuous idempotent $e_{R'}$ exhibits $M_{R'}$ as a topological direct summand of the elementary Tate module $R'\widehat\otimes_R E$.

\begin{example}\label{exRTm}
A finite projective $R(\!( q )\!)$-module, with the
topology induced by a finite free summand presentation, is a
Tate $R$-module. Indeed, separating negative and nonnegative powers gives an isomorphism
of topological $R$-modules
\begin{equation}\label{eqTm}
R(\!(q)\!)^d
=
\underbrace{(q^{-1}R[q^{-1}])^d}_{\text{finite polar parts}}
\oplus
\underbrace{R\llbracket q\rrbracket ^d}_{\text{power-series part}}
 \cong
 \underbrace{\bigoplus_{m\ge1}R^d}_{P}
 \oplus
 \underbrace{\prod_{m\ge0}R^d}_{Q^\vee},
 \qquad Q:=\bigoplus_{m\ge0}R^d.
\end{equation}
Here $P$ is discrete and the product carries the product topology
with $R$ discrete. Since $P$ and $Q$ are countably generated
free $R$-modules, $R(\!( q )\!)^d$ is an elementary Tate $R$-module, and $R\llbracket q\rrbracket ^d$ is a coprojective lattice.

If $M=eR(\!( q )\!)^d$ for an idempotent matrix $e$, its entries
have a common lower Laurent bound $-c$, so
$e(q^mR\llbracket q \rrbracket^d)\subset q^{m-c}R\llbracket q \rrbracket^d$. Thus $e$ is
continuous and $M$ is a topological direct summand. The induced
topology on $M$ is called its \emph{standard topology}; it is
independent of the presentation, since the comparison maps
extend to finite Laurent matrices and are continuous.
These are the examples of \cite[Section~3.2.2(1)--(2)]{Drinfeld}.

If $M$ is finite projective over $R(\!( q )\!)$,
then the construction with the block matrix above extends a presenting idempotent matrix
coefficientwise to $R'(\!( q )\!)$. In particular,
\[
 (R\llbracket q \rrbracket)_{R'}=R'\llbracket q \rrbracket,\qquad
 \left(\prod_{j\ge0}R\right)_{R'}=\prod_{j\ge0}R'.
\]
Completed change of coefficients are not necessarily usual tensor products,
 for example, $
\mathbb C[x]\otimes_{\mathbb C}\mathbb C\llbracket q\rrbracket
\cong \mathbb C\llbracket q\rrbracket[x]
\subsetneq \mathbb C[x] \llbracket q\rrbracket \cong\mathbb C[x]\widehat\otimes_{\mathbb C}\mathbb C\llbracket q\rrbracket.
$
\end{example}

The following is Drinfeld's coprojective-lattice theorem \cite[Theorem~3.4]{Drinfeld}.

\begin{theorem}
\label{thm:drinfeld-tate-lattice}
Let $R$ be a commutative ring and $M$ a Tate $R$-module.
There is a finite affine Nisnevich covering
\[
 \pi:\coprod_{a=1}^r\Spec R_a\longrightarrow\Spec R
\]
such that, for every $a$, the Tate $R_a$-module
$M_a:=(M)_{R_a}$ has a coprojective lattice $L_a$. Here the covering can be chosen \emph{split by a finite
filtration}: $\pi$ is \'etale and there are closed subschemes
\[
 \Spec R=F_0\supset F_1\supset\cdots\supset F_m=\varnothing,
\]
each defined by a finitely generated ideal, such that
$\pi^{-1}(F_{j-1}\setminus F_j)\to F_{j-1}\setminus F_j$
admits a section for every $j$. If $R$ is Noetherian, this is
equivalent to the Nisnevich covering condition stated in
Subsection~\ref{subsec:relative-finite-factor-theorem}.
\end{theorem}

As an application one obtains:

\begin{lemma}
\label{lem:free-tate-lattice}
Let $R$ be a finitely generated $\C$-algebra and let $M$ be a
finite projective $R(\!( q )\!)$-module of positive constant rank,
with its standard topology. There is a finite affine
Nisnevich covering $\coprod_a\Spec R_a\to\Spec R$ and, for
each $a$, a coprojective lattice $L_a\subset M_a$ with a
(noncanonical) isomorphism of topological $R_a$-modules
\begin{equation}\label{eq:direct-free-lattice}
                 L_a\cong\prod_{j\ge0}R_a.
\end{equation}
\end{lemma}

\begin{proof}
Apply Theorem~\ref{thm:drinfeld-tate-lattice}. Fix one member
of the covering and write $R,M,L$ for $R_a,M_a,L_a$.
The ring $R$ is still of finite type over $\C$, since the
affine \'etale covering maps are of finite presentation.
Lemma~\ref{lem:coprojective-lattice-criterion} gives that
$L^\vee$ is discrete and countably generated projective.
For every maximal ideal $\mathfrak m\subset R$, put
$F=R/\mathfrak m$. Completed change of coefficients in
the splitting supplied by that lemma, together with Tate
biduality, gives
\[
 M_F\cong L_F\oplus\bigl((M/L)\otimes_RF\bigr),
 \qquad
 L_F\cong(L^\vee\otimes_RF)^\vee.
\]
The second summand of the first decomposition is discrete,
so $L_F$ is open in $M_F$. The dual description of $L_F$
shows that its quotients by open subspaces are finite
dimensional: an open subspace contains the common kernel
of evaluation on finitely many vectors of
$L^\vee\otimes_RF$. Thus $L_F$ is a lattice. Moreover,
coefficientwise change of the finite Laurent idempotent
identifies $M_F$ with
$M\otimes_{R(\!( q )\!)}F(\!( q )\!)$, whose rank is the given
positive rank of $M$. Choose an $F(\!( q )\!)$-basis of $M_F$,
and write $r:=\operatorname{rank}M>0$. Openness and the
lattice condition give integers $u,v\ge0$ such that
\[
 q^uF\llbracket q \rrbracket^r\subset L_F\subset q^{-v}F\llbracket q \rrbracket^r:
\]
the second inclusion follows by bounding the poles of
representatives of a basis of
$L_F/(L_F\cap F\llbracket q \rrbracket^r)$. Hence $L_F$ is infinite
dimensional, while $L^\vee\otimes_RF$ is at most
countably generated. Its dimension is therefore
$\aleph_0$: a finite-dimensional vector space would have
a finite-dimensional dual. Since $L^\vee$ is projective,
its localization at $\mathfrak m$ is free
\cite[Tag~0593]{Stacks} and has this same infinite rank.
The ring $R$ is Noetherian, so its
quotient by its Jacobson radical is Noetherian as well.
These are precisely the hypotheses of Bass's freeness
theorem \cite[Tag~0GVK]{Stacks}, which gives
$L^\vee\cong\bigoplus_{j\ge0}R$.
Dualizing proves \eqref{eq:direct-free-lattice}.
\end{proof}

For each member $\Spec R_a$ of the covering in
Lemma~\ref{lem:free-tate-lattice}, the chosen product
decomposition identifies the functor on $R_a$-algebras
\[
 R'\longmapsto R'\widehat\otimes_{R_a}L_a
                 \cong\prod_{j\ge0}R'
\]
with the functor represented by
\[
 \A^{\N}_{\Spec R_a}
       =\Spec R_a[\lambda_j:j\ge0].
\] Here note that a coprojective lattice is itself an elementary Tate module, with its induced topology, hence the completed change of coefficients of $L_a$ is well-defined.
Applied to the tangent module along an exact Laurent map over
a finite approximation scheme, this supplies the linear
parameters used in Step~2 of
Subsection~\ref{subsec:relative-finite-factor-proof}.

\subsection{Proof of the relative finite-factor theorem}
\label{subsec:relative-finite-factor-proof}
We prove Theorem~\ref{thm:relative-finite-factors}.
Work first with $S_b$, $A$, and $D$ from
Subsection~\ref{subsec:relative-finite-factor-theorem}.
Throughout this proof, $R$ ranges over commutative $A$-algebras.
For an affine $A$-scheme $Z$, write $\mathcal O(Z)=\Gamma(Z,\mathcal O_Z)$.
For an $A$-morphism $\xi:\Spec R\to Z$, the notation $\xi^*$ on coefficient
series and matrices means coefficientwise application of
$\xi^\#:\mathcal O(Z)\to R$; it uses the source algebras and
coefficient notation fixed above. For Tate modules we use the
completed change of coefficients constructed in the preceding
subsection. Pullback notation is suppressed only after the
structural map has been specified.

\emph{Step 1: construct the finite approximation system and exact centers.}
We construct finite-type $A$-schemes $T_m$ with inverse limit
$S_b$ and exact solutions $b_m$ over $T_m$ approximating their
universal truncated series. These will provide the finite bases
and centers for the charts.

Fix the constant $C$ of Lemma~\ref{lem:direct-newton-retraction},
enlarging it once so that, for $u,v\in q^{-b}O_R^N$,
\[
 \nu(f(u)-f(v))\ge\nu(u-v)-C,\qquad
 \nu(dp_u)\ge-C,\qquad
 \nu(dp_u-dp_v)\ge\nu(u-v)-C.
\]
For a matrix $M=(M_{ab})$ over $E_R$, put
$\nu(M)=\min_{a,b}\nu(M_{ab})$. These polynomial bounds are uniform
in $R$, and $C$ remains fixed throughout the proof.
For an integer $m>2C$, take the affine space over $A$ with
coordinates $t_{\ell,i,r}$, where
$1\le\ell\le N$, $0\le i<k$, and $-b\le r<m$, and form
the truncated tuple with components
\[
 (t_m)_\ell=\sum_{i=0}^{k-1}\sum_{r=-b}^{m-1}
                  t_{\ell,i,r}z^iq^r.
\]
Define $T_m$ by the coefficient equations
$f(t_m)\equiv0\pmod {q^{m-C}}$, and denote the resulting
universal tuple in $E_{\mathcal O(T_m)}^N$ again by $t_m$.
The polynomial lower bounds make this a finite list of equations,
so $T_m$ is of finite type over $A$. We verify that truncation
defines morphisms
\[
 \tau_m:S_b\longrightarrow T_m,\qquad
 \tau_{m',m}:T_{m'}\longrightarrow T_m\quad(m'\ge m).
\]
For $\gamma\in S_b(R)$, let $u_m\in q^{-b}O_R^N$ be its
truncation to exponents less than $m$. Then
$\nu(\gamma-u_m)\ge m$, so $f(\gamma)=0$ and the first
polynomial estimate give $\nu(f(u_m))\ge m-C$.
Thus $u_m$ defines the point $\tau_m(\gamma)\in T_m(R)$.
Likewise, for an $R$-point of $T_{m'}$ represented by $u_{m'}$,
its truncation $u_m$ satisfies
\[
 \nu(u_{m'}-u_m)\ge m,\qquad
 \nu(f(u_{m'})-f(u_m))\ge m-C.
\]
Since $\nu(f(u_{m'}))\ge m'-C\ge m-C$, this gives a point
of $T_m(R)$. These coefficient projections are morphisms of
$A$-schemes and satisfy
$\tau_m=\tau_{m',m}\tau_{m'}$ and the corresponding composition
identities for successive truncations.

Moreover, $S_b\cong\varprojlim_m T_m$. Indeed, compatible
points $\xi_m\in T_m(R)$ give compatible tuples
$u_m=\xi_m^*t_m$, which determine a unique
$\gamma\in q^{-b}O_R^N$. The estimate
$\nu(f(\gamma)-f(u_m))\ge m-C$ gives
$f(\gamma)\in q^{m-C}O_R^c$ for every $m>2C$, hence
$f(\gamma)=0$ by separatedness. This construction is inverse
to truncation and commutes with $A$-algebra homomorphisms.

The scheme $W_{b,m-C}$ parametrizes full Laurent tuples,
whereas $T_m$ retains only their coefficients of $q$-degree
less than $m$. The preceding difference estimate gives morphisms
\[
 T_m
 \xrightarrow{\text{extension by zero}}
 W_{b,m-C}
 \xrightarrow{\text{truncation}}
 T_m
\]
whose composite is $\mathrm{id}_{T_m}$. Here extension by zero
sets all coefficients of $q$-degree at least $m$ equal to zero. The estimates also show that \(T_m\) is a finite factor of the approximate-solution scheme \(W_{b,m-C}\).  We will use this finite approximation to construct finite-factor charts for the exact-solution scheme \(S_b\).

To construct the exact center $b_m$, apply the morphism
$\rho_{b,m-C}$ of Lemma~\ref{lem:direct-newton-retraction}(i)
to the universal tuple $t_m\in W_{b,m-C}(\mathcal O(T_m))$.
This gives
\begin{equation}\label{eq:direct-centers}
 b_m=\rho(t_m)\in X(E_{\mathcal O(T_m)}),\qquad
 \nu(b_m-t_m)\ge m-2C.
\end{equation}
Its coefficients are regular functions on $T_m$ because
$\rho_{b,m-C}$ is a morphism of $A$-schemes.
For $\xi:\Spec R\to T_m$, coefficient-change compatibility gives
\[
 \xi^*b_m=\rho_R(\xi^*t_m).
\]
In particular, if $\gamma\in S_b(R)$ and $\xi=\tau_m(\gamma)$,
then $\nu(\gamma-\xi^*b_m)\ge m-2C$. 

\emph{Step 2: choose a product lattice over one finite approximation.}
We construct the linear parameters for the local charts by
choosing a finite Nisnevich cover at one fixed cutoff and,
on each chart, a coprojective product lattice with Laurent
bounds preserved by every change of coefficients. The cover
and this lattice will be used for all larger cutoffs.

Choose an integer $m_0\ge4C+1$. The universal center is a morphism
$b_{m_0}:\Spec E_{\mathcal O(T_{m_0})}\to X$. Put
\[
 \begin{gathered}
 M_0:=\Gamma\bigl(\Spec E_{\mathcal O(T_{m_0})},b_{m_0}^*T_X\bigr)
       =\operatorname{im}P_0,\\
 P_0:=dp_{b_{m_0}}\in\operatorname{Mat}_N(E_{\mathcal O(T_{m_0})}).
 \end{gathered}
\]
The module $M_0$ is finite projective of rank $n$ over
$E_{\mathcal O(T_{m_0})}$, hence of rank $nk$ over
$\mathcal O(T_{m_0})(\!(q)\!)$ by restriction of scalars;
the source algebra is finite free of rank $k$ over this ring.
Using the basis $1,z,\ldots,z^{k-1}$, we also write $P_0$
for its $kN\times kN$ restriction-of-scalars matrix.
We use this convention for the other $E_R$-linear operators
below. Their lower $q$-bounds persist in these matrix
representations, since the source multiplication table has
coefficients in $A[q]$.
Lemma~\ref{lem:free-tate-lattice} gives a finite affine Nisnevich
covering $\{V_\alpha\to T_{m_0}\}$ on which the completed
base changes of $M_0$ have coprojective product lattices.
Work on one chart $V=\Spec R_0$ and write
\[
 M_{0,V}:=R_0\widehat\otimes_{\mathcal O(T_{m_0})}M_0,
 \qquad L_0\subset M_{0,V},
\]
the latter denoting a coprojective lattice.
Put $J=\Z_{\ge0}$ and choose a topological isomorphism
\[
 \iota_0:\prod_{j\in J}R_0\xrightarrow{\sim}L_0.
\]
Here $V$ is affine of finite type over $\C$, and
$M_{0,V}/L_0$ is projective over $R_0$. We use $P_0$ also
for the coefficientwise base change of its restriction-of-scalars
matrix to $R_0(\!(q)\!)$. With
$E_{R_0}^N=R_0(\!(q)\!)^{kN}$ in the chosen basis,
there are integers $a,b'$ such that
\begin{equation}\label{eq:direct-lattice-bounds}
 q^aP_0(R_0\llbracket q \rrbracket^{kN})\subset L_0
                      \subset q^{-b'}R_0\llbracket q \rrbracket^{kN}.
\end{equation}
The first bound follows from openness of $L_0$ and a lower
bound on the entries of $P_0$. For the second, the quotient
$L_0/(L_0\cap R_0\llbracket q\rrbracket^{kN})$ is finitely
generated, and finitely many representatives of its generators
have a common lower Laurent bound. These bounds express
continuous factorizations
\[
 R_0\llbracket q \rrbracket^{kN}\xrightarrow{q^aP_0}L_0
       \hookrightarrow R_0(\!( q )\!)^{kN},
 \qquad
 L_0\longrightarrow q^{-b'}R_0\llbracket q \rrbracket^{kN}
       \hookrightarrow R_0(\!( q )\!)^{kN}.
\]
The inclusion of $L_0$ in the ambient Laurent module splits
continuously: $M_{0,V}/L_0$ is discrete projective, and
$M_{0,V}$ is the image of the continuous idempotent $P_0$.
Completed change of coefficients preserves this splitting and the two
factorizations. Thus \eqref{eq:direct-lattice-bounds} holds
after every homomorphism $R_0\to R$.

\emph{Step 3: transport and scale the lattice at the varying centers.}
We identify next the tangent modules at the longer-precision
centers with the one used in Step~2, then scale the transported
lattice so that the nonlinear correction preserves it and
is contracting. The centers need not agree under truncation;
the following estimate controls their difference.
For $m\ge m_0$, put
\[
 V_m=V\times_{T_{m_0}}T_m,\qquad
 B_m=\mathcal O(V_m)
     =R_0\otimes_{\mathcal O(T_{m_0})}\mathcal O(T_m).
\]
Pull $b_m$ back along $V_m\to T_m$ and denote this center
by $b_1\in S_b(B_m)$; set $P_1=dp_{b_1}$.
Pull $b_{m_0}$ back along $V_m\to V\to T_{m_0}$ and $P_0$
along $V_m\to V$, retaining their notation. Thus all these objects have coefficients
over the same ring $B_m$. Truncation and
\eqref{eq:direct-centers} give
$\nu(b_1-b_{m_0})\ge m_0-2C$, uniformly for $m\ge m_0$.
In the matrix representation of Step~2, set
\[
 U=P_1P_0+(1-P_1)(1-P_0)
   =1+(P_1-P_0)(2P_0-1).
\]
The bounds fixed in Step~1 give
$\nu(U-1)\ge m_0-4C\ge1$, so
$U\in1+q\operatorname{Mat}_{kN}(B_m\llbracket q\rrbracket)$.
The matrix $U$ is invertible by its geometric-series inverse,
whose coefficients are regular on $V_m$, and
$UP_0=P_1U$. 
Thus  $
U:\operatorname{im}P_0
\xrightarrow{\sim}\operatorname{im}P_1
$
transports tangent vectors at the old center $b_{m_0}$ to tangent vectors at the new center $b_1$.

For a point $\xi:\Spec R\to V_m$, the maps
$R_0\to B_m\xrightarrow{\xi^\#}R$ define
\[
 L_{0,R}:=R\widehat\otimes_{R_0}L_0,\qquad
 U_\xi:=\xi^*U,\qquad P_{i,\xi}:=\xi^*P_i\quad(i=0,1).
\]
Both $U_\xi$ and $U_\xi^{-1}$ preserve
$R\llbracket q\rrbracket^{kN}$, and
\eqref{eq:direct-lattice-bounds} gives
\[
 q^aP_{1,\xi}(R\llbracket q\rrbracket^{kN})
 \subset U_\xi L_{0,R}
 \subset q^{-b'}R\llbracket q\rrbracket^{kN}.
\]
Choose an integer
$s\ge\max\{b'+C+1,\,2b'+a+C\}$ and put $h=s-b'$. Then
\begin{equation}\label{eq:direct-small-lattice-choice}
                 h>C,\qquad 2h-C\ge s+a.
\end{equation}
For each $\xi$, the topological isomorphism
\[
 q^sU_\xi:\operatorname{im}P_{0,\xi}
       \xrightarrow{\sim}\operatorname{im}P_{1,\xi}
\]
carries $L_{0,R}$ onto the coprojective product lattice
\[
 L_\xi:=q^sU_\xi L_{0,R}\cong\prod_{j\in J}R,
 \qquad L_\xi\subset q^hO_R^N.
\]
Thus $L_\xi$ is complete for its induced topology, and every
element has order at least $h$. These constructions commute
with completed changes of coefficients along homomorphisms
of $B_m$-algebras.
For calculations at a fixed
$\xi$, we henceforth write $b_1,P_1,U$ for their coefficientwise
pullbacks to $R$ and put $L=L_\xi$.

\emph{Step 4: parametrize exact solutions by the lattice.}
We integrate the tangent displacements in $L$ of Step~3 to
obtain a parametrization $L\to S_b(R)$ and identify its image.
Step~5 will realize these parametrizations by a morphism of
$V_m$-schemes.

For $v\in L$, the map $H_{b_1}$ in
Lemma~\ref{lem:direct-newton-retraction} takes values in $L$:
its values are tangent and have order at least $2h-C\ge s+a$.
Indeed, a tangent vector $u\in\operatorname{im}P_1$ with
$\nu(u)\ge s+a$ satisfies
\[
 u=q^{s+a}P_1(q^{-s-a}u)
       \in q^{s+a}P_1(O_R^N)\subset L.
\]
Put $\delta=h-C>0$. Equation~\eqref{eq:direct-newton-lipschitz}
gives
\[
 \nu\bigl(H_{b_1}(v)-H_{b_1}(w)\bigr)
       \ge\nu(v-w)+\delta\qquad(v,w\in L).
\]
Thus $H_{b_1}:L\to L$ contracts for the topology induced by
the submodules $L\cap q^dO_R^N$, $d\in\Z$.
Define the map of sets
\[
 g:L\longrightarrow L,\qquad g(v)=v+H_{b_1}(v),
\]
the restriction of $\mathrm{id}+H_{b_1}$ to $L$.
For each $w\in L$, consider the sequence
\[
 v_0=w,\qquad
 v_{\ell+1}=w-H_{b_1}(v_\ell)\quad(\ell\ge0).
\]
Since $L$ is complete and $H_{b_1}$ satisfies the preceding
strict contraction estimate, this sequence converges in $L$
to the unique solution of $v+H_{b_1}(v)=w$. Thus $g$ is
bijective, with
\[
 g^{-1}(w)=\lim_{\ell\to\infty}v_\ell.
\]
For this $\xi$, define a map of sets $\psi_{b_1}:L\to S_b(R)$ by
\begin{equation}\label{eq:direct-lattice-parameterization}
 \psi_{b_1}(w)=\Phi_{b_1}(g^{-1}(w)),\qquad
 P_1(\psi_{b_1}(w)-b_1)=w.
\end{equation}
The ambient estimate in the proof of
Lemma~\ref{lem:direct-newton-retraction} also gives
\[
 \nu(\Phi_{b_1}(v)-b_1-v)\ge2\nu(v)-C
                         \qquad(v\in L).
\]
Since $h>C$ and $g^{-1}(w)\in L$, this implies
$\nu(\psi_{b_1}(w)-b_1)\ge h$ for every $w\in L$.
Its image is exactly
\[
 \psi_{b_1}(L)=
 \{x\in X(E_R):\nu(x-b_1)\ge h,\ P_1(x-b_1)\in L\}.
\]
Indeed, the forward inclusion follows from
\eqref{eq:direct-lattice-parameterization} and the preceding
estimate. If $x,y\in X(E_R)$ both differ from $b_1$ by
vectors of order at least $h$ and $P_1(x-y)=0$, Taylor
expansion of $p(x)-p(y)=x-y$ around the common center \(b_1\) gives
\[
             \nu(x-y)\ge\nu(x-y)+h-C,
\]
with the at-least-quadratic term giving the bound on the right-hand side,
and so $x=y$. Applying this uniqueness to $x$ in the displayed
set and $y=\psi_{b_1}(P_1(x-b_1))$ proves the reverse inclusion.
In particular, the image contains every exact $x$ with
$\nu(x-b_1)\ge\max\{h,s+a+C\}$, since then
$\nu(P_1(x-b_1))\ge s+a$ and hence $P_1(x-b_1)\in L$.

\emph{Step 5: prove algebraicity of the parametrizations.}
Next we realize the family of
parametrizations $\psi_{b_1}:L\to S_b(R)$ as a morphism
of $V_m$-schemes
\[
 \A^J_{V_m}\longrightarrow V_m\times_{\Spec A}S_b.
\]
Put $B_\lambda:=B_m[\lambda_j:j\in J]$, so that
$\A^J_{V_m}=\Spec B_\lambda$. For a $B_m$-algebra $R$, write
$b_{1,R},U_R$ for the coefficientwise images of $b_1,U$,
$L_R$ for the lattice of Step~3 at $\Spec R\to V_m$, and
$\iota_{0,R}$ for the completed change of coefficients of
$\iota_0$ from Step~2. Apply Step~4 with $R=B_\lambda$ and
the tuple of indeterminates
$\lambda=(\lambda_j)_{j\in J}\in\prod_{j\in J}B_\lambda$.
For
\[
 w^{\mathrm{univ}}:=q^sU_{B_\lambda}\iota_{0,B_\lambda}(\lambda)
               \in L_{B_\lambda},
\]
there is a unique $v^{\mathrm{univ}}\in L_{B_\lambda}$ satisfying
\[
 v^{\mathrm{univ}}+H_{b_{1,B_\lambda}}(v^{\mathrm{univ}})
       =w^{\mathrm{univ}}.
\]
The iteration of Step~4 converges with respect to the submodules
$q^dO_{B_\lambda}^N$, so each Laurent coefficient stabilizes in the
polynomial ring $B_\lambda$ itself. The resulting exact solution
\[
 x^{\mathrm{univ}}:=\Phi_{b_{1,B_\lambda}}(v^{\mathrm{univ}})
                \in S_b(B_\lambda)
\]
and the structure map $\Spec B_\lambda\to V_m$ together define a
$V_m$-morphism $\Psi_m:\A^J_{V_m}\to V_m\times_{\Spec A}S_b$.

For $\mu=(\mu_j)_{j\in J}\in\prod_{j\in J}R$, let
$\varphi_\mu:B_\lambda\to R$ be the $B_m$-algebra homomorphism
$\lambda_j\mapsto\mu_j$. Write \((\varphi_\mu)_*:E_{B_\lambda}\to E_R\) for the induced homomorphism given by
\[
  (\varphi_\mu)_*\!\left(\sum_{r\ge r_0}\sum_{i=0}^{k-1}
     a_{ir}(\lambda)z^iq^r\right)
  =\sum_{r\ge r_0}\sum_{i=0}^{k-1}
     a_{ir}(\mu)z^iq^r,
  \qquad a_{ir}(\lambda)\in B_\lambda.
\]
Use the same notation for its componentwise action on tuples.
By Step~3 it carries $w^{\mathrm{univ}}$
to $w_\mu:=q^sU_R\iota_{0,R}(\mu)\in L_R$.
By Lemma~\ref{lem:direct-newton-retraction}(ii),
$(\varphi_\mu)_*v^{\mathrm{univ}}$ satisfies
$v+H_{b_{1,R}}(v)=w_\mu$ in $L_R$. Uniqueness in Step~4 and
the same lemma for $\Phi$ give
\[
 (\varphi_\mu)_*x^{\mathrm{univ}}=\psi_{b_{1,R}}(w_\mu).
\]
Thus $\Psi_m$ realizes the parametrization of Step~4 over every $R$.

For $u\in E_R^N$ and $r\in\Z$, write
\[
 \mathbf c_r(u):=\bigl([z^iq^r]u_\ell\bigr)_
        {1\le\ell\le N,\,0\le i<k}\in R^{kN}.
\]
Let $F_r^\psi\subset J$ consist of the indices $j$ for which
$\lambda_j$ occurs in an entry of $\mathbf c_r(x^{\mathrm{univ}})$.
These entries are polynomials in $B_\lambda$, so $F_r^\psi$ is finite.

\emph{Step 6: construct the finite factors and prove the atlas assertion for $S_b$.}
We require the exact solution \(x^{\mathrm{univ}}\) parametrized by \(\Psi_m\) in Step~5 to agree modulo \(q^m\) with the pullback of \(t_m\) along \(\mathbb A^J_{V_m}\to V_m\to T_m\). We will show that this requires
only finitely many product coordinates and yields an algebraic
product chart for the entire pullback of $S_b$.

Choose an integer
\[
 m\ge\max\{m_0,\,2C+h,\,3C+s+a\}.
\]
Write $U_V=S_b\times_{T_{m_0}}V$. A point
$(\gamma,\zeta)\in U_V(R)$ consists of $\gamma\in S_b(R)$
and a lift $\zeta:\Spec R\to V$ of $\tau_{m_0}(\gamma)$.
It determines
\[
 \xi=(\zeta,\tau_m(\gamma)):\Spec R\longrightarrow V_m.
\]
Use this $\xi$ for $b_1,P_1,U,L$ and the parametrization of
Step~4. Then \eqref{eq:direct-centers} and $\nu(P_1)\ge-C$ give
\[
 \nu(\gamma-b_1)\ge m-2C\ge h,\qquad
 \nu(P_1(\gamma-b_1))\ge m-3C\ge s+a.
\]
Thus $\gamma$ belongs to the image described in Step~4.
The constants $C,m_0$ are fixed for the atlas; the choices
$a,b',s,h,m$ may depend on $V$ and $L_0$, but not on $R$
or $(\gamma,\zeta)$.

Pull $t_m$ back along $\A^J_{V_m}\to V_m\to T_m$ and impose
\begin{equation}\label{eq:direct-jet-matching}
                  x^{\mathrm{univ}}\equiv t_m\pmod {q^m}.
\end{equation}
Since both tuples lie in $q^{-b}O_{B_\lambda}^N$, this congruence
is equivalent to the finite system
\[
 \mathbf c_r(x^{\mathrm{univ}})=\mathbf c_r(t_m),\qquad -b\le r<m.
\]
Put
\[
 J_0:=\bigcup_{r=-b}^{m-1}F_r^\psi\subset J.
\]
The set $J_0$ is finite. The scalar equations above are defined
in $B_m[\lambda_j:j\in J_0]$ and cut out a closed subscheme
$Y\subset\A^{J_0}_{V_m}$, affine of finite type over $V_m$.
In $\A^J_{V_m}$, \eqref{eq:direct-jet-matching} cuts out the
scheme-theoretic inverse image of $Y$ under the coordinate
projection $\A^J_{V_m}\to\A^{J_0}_{V_m}$, namely
\[
 \A^J_{V_m}\times_{\A^{J_0}_{V_m}}Y\cong\A^{J\setminus J_0}_Y.
\]
Restrict $\Psi_m$ to the closed subscheme $\A^{J\setminus J_0}_Y$
defined by \eqref{eq:direct-jet-matching}. Let $R$ be a $B_m$-algebra
with structure point $\xi=(\zeta,t)\in V_m(R)$. For an $R$-point
$(\xi,\lambda)$ of $\A^{J\setminus J_0}_Y$, write
$\Psi_m(\xi,\lambda)=(\xi,\gamma)$. Then
$\tau_m(\gamma)=t$ by \eqref{eq:direct-jet-matching}, so
$\tau_{m_0}(\gamma)=\tau_{m,m_0}(t)$ is the image of $\zeta$
in $T_{m_0}$. Thus the second component of $\Psi_m$ and the
projection to $V$ induce a morphism over $V_m$
\[
 \A^{J\setminus J_0}_Y\longrightarrow U_V,\qquad
 (\xi,\lambda)\longmapsto(\gamma,\zeta).
\]

Conversely, let $(\gamma,\zeta)\in U_V(R)$ lie over $\xi$.
The estimates above place $\gamma$ in $\psi_{b_1}(L)$, and
\eqref{eq:direct-lattice-parameterization} gives its unique
preimage $w=P_1(\gamma-b_1)$. The product coordinates of $w$ are
\[
 \lambda=\iota_{0,R}^{-1}\bigl(q^{-s}U^{-1}P_1(\gamma-b_1)\bigr)
          \in\prod_{j\in J}R.
\]
Step~5 gives $\Psi_m(\xi,\lambda)=(\xi,\gamma)$. Since
$\xi=(\zeta,\tau_m(\gamma))$, the pair $(\xi,\lambda)$ satisfies
\eqref{eq:direct-jet-matching}. The displayed morphism is therefore
bijective on $R$-points for every $B_m$-algebra $R$. By Yoneda
this gives the isomorphism
\begin{equation}\label{eq:direct-chart-isomorphism}
 U_V\cong\A^{J\setminus J_0}_Y\qquad\text{over }V_m.
\end{equation}

As $V$ ranges over the finite Nisnevich cover of $T_{m_0}$ chosen
in Step~2, the maps $U_V\to S_b$ form its pullback along $\tau_{m_0}$.
They are finitely presented and \'etale. By
\eqref{eq:direct-chart-isomorphism} each chart has finite factor $Y$,
and its structure map to $\Spec A$ factors through $Y$. This proves
the atlas assertion of Theorem~\ref{thm:relative-finite-factors}
for $S_b$.

\emph{Step 7: impose the marked bounds and obtain the collision atlases.}
We now prove the atlas assertion for
$\mathcal Z_{I,\mathbf P}$ by including its marked pole
equations in the finite factors of Step~6.
Take $A=\C[a_i:i\in I]$ and $D(z)=\prod_i(z-a_i)$,
and use the affine generators $g_\ell$ fixed in
Subsection~\ref{conv:loop-volume-compatibility} as the coordinates
$x_\ell$ of the embedding. Let $r$ be the fixed positive integer
used there to define the pole bounds, and identify $q$ with
$q_0$ in Definition~\ref{def:kv-coefficient-presentation}.
Definition~\ref{def:kv-bounded-family} cuts
$\mathcal Z_{I,\mathbf P}$ out of $S_b$, for
$b=\max_i(rP_i)$, by finitely many coefficient equations.
Pull them back along $U_V\to S_b$. In
\eqref{eq:direct-chart-isomorphism}, they use only finitely
many additional coordinates, indexed by some
$F\subset J\setminus J_0$. Their zero scheme $Y'$ in
$\A^F_Y$ is affine of finite type over $A$, and
\[
 U_V\times_{S_b}\mathcal Z_{I,\mathbf P}
       \cong\A^{(J\setminus J_0)\setminus F}_{Y'}.
\]
The map to $\Spec A$ factors through $Y'$, giving the asserted
collision atlas.

\emph{Step 8: prove preservation under base changes and locally closed restrictions.}
We descend a specified restriction to one finite approximation,
then pull back the finite factors to it. This proves the
base-change and restriction assertions for the atlases.
Let $T\to\Spec A$ be of finite type, and let
$W\hookrightarrow S_b\times_{\Spec A}T$ be a finitely
presented locally closed immersion; the identity immersion is
allowed. By Step~1, its ambient scheme is the inverse limit
of $T_m\times_{\Spec A}T$, whose transition maps are affine
and whose terms are Noetherian. Finite-presentation
descent gives, for some $M\ge m_0$, a locally closed immersion
$W_M\hookrightarrow T_M\times_{\Spec A}T$ such that
\[
 W\cong
 (S_b\times_{\Spec A}T)
       \times_{T_M\times_{\Spec A}T}W_M
 \cong S_b\times_{T_M}W_M.
\]
In the last expression, $W_M\to T_M$ is the projection
from its indicated immersion. This fiber-product description
supplies the map $W\to W_M$ used below. To obtain this descent,
one may first factor the given immersion through a
quasi-compact open subscheme, descend that open subscheme,
and then descend the finitely presented closed immersion
inside it; see \cite[Tags~01Z4, 01ZM, and 01ZP]{Stacks}.
In particular, $W_M$ is Noetherian and of finite type over
$\C$.

Write the charts constructed in Step~6 as
\[
 u_\alpha:U_\alpha=S_b\times_{T_{m_0}}V_\alpha\longrightarrow S_b,
 \qquad U_\alpha\cong\A^{J_\alpha}_{Y_\alpha}.
\]
For each $\alpha$, the pullbacks of finitely many algebra
generators of $\mathcal O(T_M)$ along $\tau_Mu_\alpha$ use
only a finite subset $F_\alpha\subset J_\alpha$ of the tail
coordinates. Replace $Y_\alpha$ by $\A^{F_\alpha}_{Y_\alpha}$
and $J_\alpha$ by $J_\alpha\setminus F_\alpha$, then relabel
both. For the resulting projection $\pi_\alpha:U_\alpha\to Y_\alpha$,
the map $\tau_Mu_\alpha$ descends to a morphism
\[
 \varphi_\alpha:Y_\alpha\longrightarrow T_M,
 \qquad \varphi_\alpha\pi_\alpha=\tau_Mu_\alpha.
\]
Thus $U_\alpha\cong\A^{J_\alpha}_{Y_\alpha}$ over $T_M$, and
\[
 U_\alpha\times_{S_b}W
       \cong\A^{J_\alpha}_{Y_\alpha\times_{T_M}W_M}.
\]
The map to $T$ factors through the finite-type scheme
$Y_\alpha\times_{T_M}W_M$ and its projection to $W_M$.
Finite affine open covers of these factors give the required
finite affine Nisnevich atlas even when $T$ is not affine.

\emph{Step 9: construct the finite nice filtrations.}
For $W$ as in Step~8, we prove the remaining filtration
assertion by pulling back a splitting filtration of the
finite-type scheme $W_M$. We must show that each resulting
successive difference of $W$ is itself a finite factor times
an affine coordinate space.
Pull the Nisnevich cover $\coprod_\alpha V_\alpha\to T_{m_0}$
back to $W_M$. Choose a finite closed splitting filtration
of $W_M$, and refine it so that each successive difference
$Z$ is affine and has a splitting section into one component.
Such a refinement exists by Noetherian induction, using affine
open subsets of the splitting strata. For each $Z$, fix the
corresponding index $\alpha=\alpha(Z)$ and denote its section by
\[
 \sigma_Z:Z\longrightarrow V_\alpha\times_{T_{m_0}}Z.
\]
Put $W_Z=W\times_{W_M}Z$. The corresponding chart is
\[
 U_\alpha\times_{S_b}W_Z
 \cong U_\alpha\times_{T_M}Z
 \cong\A^{J_\alpha}_{Y_\alpha\times_{T_M}Z}.
\]
Its finite factor is affine, since $Y_\alpha$, $T_M$, and
$Z$ are affine. Pulling $\sigma_Z$ back along $W_Z\to Z$ gives
a section
\[
 \widetilde\sigma_Z:W_Z\longrightarrow U_\alpha\times_{T_M}Z.
\]
A section of an \'etale morphism is an open immersion, and
a section of a separated morphism is a closed immersion.
Thus $\widetilde\sigma_Z(W_Z)$ is open and closed. Here the chart morphism is a
base change of the separated \'etale map
$V_\alpha\to T_{m_0}$; this does not require $T$ to be separated.

To obtain a finite-factor presentation of \(W_Z\), it remains to show that the open and closed image of \(\widetilde\sigma_Z\) is the inverse image of an open and closed subscheme of the finite factor. We will refer later to this last step a few times.

\emph{Idempotents in polynomial algebras.} Write $B=\mathcal O(Y_\alpha\times_{T_M}Z)$ and choose an
idempotent $e\in B[\lambda_j:j\in J_\alpha]$ whose zero scheme is
$\widetilde\sigma_Z(W_Z)$. In fact $e=e(0)\in B$, where $e(0)$ sets all
$\lambda_j$ equal to zero. To see this for arbitrary $B$, reduce
modulo each prime of $B$: an idempotent polynomial over a
domain is constant. Thus $e-e(0)$ is nilpotent, since its
finitely many coefficients are nilpotent. The idempotents
$e(1-e(0))$ and $e(0)(1-e)$ are then nilpotent and hence zero.
Consequently
\[
 W_Z\cong\A^{J_\alpha}_{\Spec(B/(e))},
\]
and its map to $T$ factors through the finite factor by
\[
 \Spec(B/(e))\longrightarrow Y_\alpha\times_{T_M}Z
       \longrightarrow Z\longrightarrow W_M\longrightarrow T.
\]
Pulling back the refined closed filtration of $W_M$ therefore
gives the required finite nice filtration of $W$. Every
closed immersion in this filtration is finitely presented,
since it is a base change of a closed immersion between
Noetherian schemes. This includes $S_b$ itself by taking
$T=\Spec A$ and $W=S_b$.

\begin{remark}
The proof of Theorem~\ref{thm:relative-finite-factors} does not  give an algorithm to compute an atlas (\ref{eq:direct-relative-nice-atlas}) starting only from the equations of \(X\) and the collision algebra. The main missing ingredient is an effective version of Step 2, which chooses the Nisnevich cover and the product lattice.
\end{remark}

\subsection{Compatible finite-factor presentations of the collision diagram}

\begin{proposition}[Finite collision diagrams]
\label{prop:direct-collision-diagrams}
The atlases of Theorem~\ref{thm:relative-finite-factors} can be
chosen with the compatible finite-factor presentations required by
Definition~\ref{def:fusion-ready-data}(i)--(ii). Explicitly,
every finite subdiagram of the collision diagram of
Definition~\ref{def:geometric-collision-data}, together with a
finite list of its \'etale overlap maps, admits compatible
finite-factor presentations over its configuration bases, with factors
chosen for chart occurrences in finite composable strings as in
Definition~\ref{def:fusion-ready-data}(ii). The
finite-factor presentations include finite-factor models of the pole inclusions,
collision diagonals, open complements, product and relabeling
identifications, with their prescribed composition identities.
They may additionally
retain any finite list of regular functions and maps to
chosen finite factors.

Products and affine changes of source coordinate $z\mapsto uz+v$,
with $u$ invertible, preserve this assertion. Arbitrary finite-type
changes of configuration base preserve the algebraic atlases,
filtrations, presented morphisms, and common refinements constructed
here. Over separated new bases these are compatible finite-factor
presentations in the sense of
Subsection~\ref{subsec:relative-volume-preparations} and
Definition~\ref{def:fusion-ready-data}(i)--(ii). Any two finite
choices have common refinements. For an infinite diagram,
these assertions concern its finite subdiagrams and finite
lists of overlap maps; no single finite factor for the
whole infinite diagram is asserted.
\end{proposition}

\begin{proof}
Theorem~\ref{thm:relative-finite-factors} supplies the atlases
and filtrations of Definition~\ref{def:fusion-ready-data}(i).
We prove its finite-diagram requirement (ii).
For a finite list of pole bounds choose a common larger
bound. Lemma~\ref{lem:kv-pole-compatibility} makes every
smaller pole locus a finitely presented closed subscheme.
Present these closed inclusions and their open complements
by Subsection~\ref{subsec:volume-index-preparations}, over the
configuration base as in
Subsection~\ref{subsec:relative-volume-preparations}. A diagonal
restriction is a base change, so its pulled-back chart is
$Y\times_TT'\times\A^J$ with the same omitted coordinates.
Products of charts give the separation finite-factor presentations.

Transport charts along the geometric isomorphisms and take
their \'etale fiber products when two descriptions overlap.
After finite affine refinement where necessary, let
$W\to U=Y\times\A^J$ be one such overlap. Write
\[
 U=\varprojlim_{F\subset J,\ F\text{ finite}}
                         (Y\times\A^F).
\]
The terms are affine and the transition maps are affine.
Finite presentation gives an affine finite-type scheme
$W_F$ over $Y\times\A^F$ and an isomorphism
\[
 W\cong W_F\times_{Y\times\A^F}U
       \cong W_F\times\A^{J\setminus F}.
\]
After enlarging $F$, the descended map is \'etale
\cite[Tags~01ZM and 07RP]{Stacks}.

The same finite-presentation argument gives Cartesian models
for the induced geometric fiber products. For affine finitely
presented arrows $V_i\to U_0=Y_0\times\A^{J_0}$, $i=1,2$,
descend both algebra presentations to one finite subset
$F\subset J_0$. Write the resulting maps as
\[
 V_i\cong Y_i\times\A^{J_0\setminus F}
 \xrightarrow{\ h_i\times\mathrm{id}\ }
 (Y_0\times\A^F)\times\A^{J_0\setminus F},
 \qquad h_i:Y_i\to Y_0\times\A^F.
\]
Then
\[
 V_1\times_{U_0}V_2\cong
 (Y_1\times_{Y_0\times\A^F}Y_2)\times\A^{J_0\setminus F}.
\]
Enlarging $F$ preserves this identity by base change. Together
with the configuration pullbacks and product charts above,
this treats the finite sequence of geometric fiber products
in the selected diagram. It asserts Cartesianity for these
geometric finite-factor presentations; comparison maps between different
finite factors need not form Cartesian squares.

For the finite strings of chart maps stipulated in
Definition~\ref{def:fusion-ready-data}(ii), apply
Lemma~\ref{lem:finite-presentation-categories}(i)--(iii) with
the prescribed closed collections. Include the additionally
specified functions and the pullbacks of coordinate-ring
generators from the targets of the specified finite-factor
and affine configuration maps among the functions in
part~(iii). The lemma supplies compatible
finite-factor presentations, their smooth comparisons, and the cofinal
comparisons under omission and repetition of occurrences.
Equalities are preserved on common refinements: if $q:V\to Z$ is a finite-factor
projection on a common chart and $h_1,h_2:Z\to Y$ represent
two equal composites, then $h_1q=h_2q$. Composing with the
zero section of $q$ gives $h_1=h_2$.

For an affine source change put $z'=uz+v$ and $a_i'=ua_i+v$.
Then $z'-a_i'=u(z-a_i)$ and $q'=u^kq$, so it carries
the marked ideals to their counterparts multiplied by units.
In monic coefficient coordinates its maps and inverse use,
at every finite order, finitely many coefficients and the
regular functions $u,u^{-1},v$. It therefore fits the same procedure over
the enlarged finite-type parameter base. Products commute with
the constructions just described. After an arbitrary finite-type
base change $T'\to T$, each factor $Y\times_TT'$ is Noetherian
and has a finite affine open cover. The induced charts have
quasi-compact, quasi-separated overlaps, finitely presented
\'etale over either chart. Finite affine covers of these overlaps
and the same finite-presentation descent give nice overlap
finite-factor presentations. Retaining the maps to the finite factors of both
induced charts gives common refinements over $T'$. Thus these algebraic
assertions do not require $T'$ to be separated; the references
to compatible finite-factor presentations in the stated formalism
concern separated bases. Finally, the source-algebra identities of
Lemma~\ref{lem:kv-pole-compatibility} already commute before
choosing finite-factor presentations. Present their finite lists of maps;
the zero-section argument above then gives the same identities
for the descended maps on common refinements.
\end{proof}

\section{Fusion-ready families}
\label{sec:direct-fusion-ready}
\label{sec:fusion-ready-families}

In this section we show that our collision families can be made fusion-ready.

\subsection{Relative volume indices}
\label{subsec:direct-relative-volume-indices}

We extend the volume-index construction of
Section~\ref{sec:formal-loops} to the collision families
over their configuration bases.  As in the one-point case, we compute the index from the kernel of the coefficient-residue map, using the volume form to normalize the determinant. Here the marked points vary over the base, so we differentiate loop coefficients with the marked points held fixed and prove that the finite kernel quotients defining the index commute with the specified base changes, including restriction to collision diagonals. On each collision stratum, the residue splits as a sum over the distinct marked points, and the determinant comparison contributes zero to the index. Thus the index on a product finite factor is the sum of the one-disk indices, one for each distinct marked point. We also check that relabeling and affine changes of source coordinate preserve this normalization.

Fix the compatible finite-factor
presentations of Theorem~\ref{thm:relative-finite-factors} and
Proposition~\ref{prop:direct-collision-diagrams}.
Let $T\to\A^I$ be a specified base change, put $k=|I|$, and
take a family and a presented chart in this system:
\[
 \begin{gathered}
 f:S=T\times_{\A^I}\mathcal Z_{I,\mathbf P}\longrightarrow T,
 \qquad u:U=Y\times\A^J\longrightarrow S,\\
 \pi:U\longrightarrow Y,\qquad f_Y:Y\longrightarrow T,
 \qquad f\circ u=f_Y\circ\pi.
 \end{gathered}
\]
Here $u$ is a nice \'etale chart, $\pi$ is the finite-factor
projection, and $f_Y$ is the map through which $f u$ factors.
For an affine test scheme $V=\Spec R\to U$, let
$a=(a_i)_{i\in I}\in R^I$ be the configuration defined by
\[
 V\longrightarrow U\xrightarrow{u}S\xrightarrow{f}T\longrightarrow\A^I.
\]
Thus $a_i$ is the pullback of the $i$-th coordinate on $\A^I$.
Use the source algebras $\mathcal O_a$ and $\mathcal K_a$ of
Definitions~\ref{def:kv-multipoint-functor} and
\ref{def:kv-coefficient-presentation}, writing $q$ for $q_0$ in
\eqref{eq:kv-monic-source}. In this presentation
$q=D_a(z)=\prod_{i\in I}(z-a_i)$, and $1,z,\ldots,z^{k-1}$ form
an ordered $R\llbracket q\rrbracket$-basis of $\mathcal O_a$.
The universal loop restricts to $\gamma_V\in X(\mathcal K_a)$.
Define the $R(\!( q )\!)$-linear map
\[
 \lambda_R:\mathcal K_a\longrightarrow R(\!( q )\!),\qquad
 \lambda_R\!\left(\sum_{\nu=0}^{k-1}b_\nu z^\nu\right)
       =b_{k-1}.
\]
The pairing $(b,c)\mapsto\lambda_R(bc)$ is perfect: its
matrix in this basis satisfies
$\lambda_R(z^{r+s})=0$ for $r+s<k-1$ and
$\lambda_R(z^{r+s})=1$ for $r+s=k-1$, and has determinant
$(-1)^{k(k-1)/2}$.  This argument works over every $R$ and
does not assume that $\mathcal K_a$ is \'etale over $R(\!( q )\!)$.

Put $N_V=\gamma_V^*\Omega_X^1$. The Laurent cotangent module
in the parameter $q$ is
\[
 \Cot_q(\gamma_V):=
 \bigl(\operatorname{Res}_{\mathcal K_a/R(\!( q )\!)}N_V\bigr)
       \otimes_{R(\!( q )\!)}R(\!( q )\!)\,dq,
\]
where $\operatorname{Res}$ denotes restriction of scalars and
$dq$ is the formal generator. Thus $\Cot_q(\gamma_V)$ is finite
projective of rank $d=kn$ over $R(\!( q )\!)$, where $n=\dim X$.
To construct the integer volume indices, use the nowhere-vanishing
form $\Omega\in\Gamma(X,\Omega_X^n)$ of
Subsection~\ref{conv:loop-volume-compatibility}.
For a commutative \(B\)-algebra \(A\) finite free as a \(B\)-module, write
\[
\operatorname{Nm}_{A/B}(b)=\det_B(m_b),
\qquad m_b:A\longrightarrow A,\quad c\longmapsto bc.
\]
For \(G\in\operatorname{GL}_n(A)\), let \(\operatorname{Res}_{A/B}G\) denote its action on \(A^n\) as a \(B\)-module. Then
\[
\det_B(\operatorname{Res}_{A/B}G)
=\operatorname{Nm}_{A/B}(\det_A G).
\]
To check this identity, localize \(B\) at a prime. The resulting finite algebra is semilocal, so every invertible matrix over it is a product of elementary matrices and diagonal matrices. On an elementary matrix both sides are \(1\). On a diagonal matrix both sides are the product of the norms of its diagonal entries.
Apply this with \(B=R(\!(q)\!)\) and \(A=\mathcal K_a\). Locally on \(\operatorname{Spec}B\), choose an \(A\)-basis \(\eta_1,\ldots,\eta_n\) of \(N_V\) satisfying
\[
\eta_1\wedge\cdots\wedge\eta_n=\gamma_V^*\Omega.
\]
Such a basis exists near each prime \(\mathfrak p\) of \(B\). The constant-rank projective module \(N_V\otimes_B B_{\mathfrak p}\) is free over the semilocal ring \(A\otimes_B B_{\mathfrak p}\). A basis extends to a neighborhood of \(\mathfrak p\), where rescaling its first vector gives the displayed equality. On this neighborhood take as generator of \(\det_B\operatorname{Cot}_q(\gamma_V)\) the element
\[
\bigwedge_{i=1}^{n}
\left(\,\bigwedge_{\nu=0}^{k-1}(z^\nu\eta_i\otimes dq)\right).
\]
Two bases with exterior product \(\gamma_V^*\Omega\) differ by a matrix \(G\) with \(\det_A G=1\). The determinant identity above therefore makes their displayed generators equal. These generators glue and commute with coefficientwise change of \(R\), since the ordered basis \(1,z,\ldots,z^{k-1}\), multiplication matrices, and exterior products all commute with that change.

Define the relative coefficient-residue map $\kappa_{\pi,V}$ by
\begin{equation}\label{eq:direct-relative-residue-map}
 \begin{aligned}
 \kappa_{\pi,V}:\Cot_q(\gamma_V)&\longrightarrow
                   (\Omega^1_{U/Y})|_V=\bigoplus_{j\in J}R\,du_j,\\
 b\,dg\otimes dq&\longmapsto
 \operatorname{Res}_{q=0}
       \lambda_R\!\left(b\,(d_{U/Y}g(\gamma_U))|_V\right)dq.
 \end{aligned}
\end{equation}
Here the differential is taken coefficientwise on \(U\), before
restriction to \(V\). More precisely, put
\(C=\Gamma(U,\mathcal O_U)\) and
\(M=\Gamma(U,\Omega^1_{U/Y})\); the differential of
\(\sum_r c_rq^r\) is \(\sum_r(dc_r)q^r\), a Laurent series
with coefficients in \(M\) and a finite lower bound.
Write $a_i$ also for the marked coordinate on $U$ before
restriction to $V$. Since $f u=f_Y\pi$, it is pulled back from
$Y$, so $d_{U/Y}a_i=0$.
Fixing \(q,z\) therefore makes this derivation respect
\(D_a(z)-q=0\). Apply \(M\to M\otimes_C R\) coefficientwise
and extend \(\lambda_R\) by extracting the \(z^{k-1}\) component.
Only finitely many terms contribute to the residue.
The Leibniz rule and the equations of \(X\) show that the formula
respects the relations in \(\gamma_V^*\Omega_X^1\).
Thus $\kappa_{\pi,V}$ is an \(R\)-linear map to
\(M\otimes_C R=\bigoplus_{j\in J}R\,du_j\).
For $I=\{1\}$, one has $q=z-a_1$ and $\lambda_R=\mathrm{id}$,
so $\kappa_{\pi,V}$ has the same coefficient formula as
$q_{\pi,V}$ in \eqref{eq:elementary-residue-map}, written
in the source parameter $q$.

\begin{lemma}[Relative coefficient kernels]
\label{lem:direct-relative-coefficient-kernel}
The relative coefficient-residue map $\kappa_{\pi,V}$ is surjective.
It commutes with coefficientwise change of $R$.
Put $L_V=\ker\kappa_{\pi,V}$. If $\Cot_q(\gamma_V)$ is free over
$R(\!( q )\!)$, choose a basis. In these coordinates
\[
 q^NR\llbracket q \rrbracket^d\subset L_V\subset q^{-N}R\llbracket q \rrbracket^d
\]
for some $N\ge0$. Put $H_R=q^NR\llbracket q\rrbracket^d$ in this
basis. Then $L_V/H_R$ is finite projective over $R$, and for
every $R\to R'$ the induced basis gives the same bounds and
\[
 (L_V/H_R)\otimes_R R'\cong L_{V'}/H_{R'},
 \qquad V'=\Spec R'\to V,\quad H_{R'}=q^NR'\llbracket q\rrbracket^d.
\]
For a field-valued point
$x:\Spec F\to U$, where $F/\C$ is a field extension, define
the relative pointwise coefficient index $c_\pi(x)$, the analogue
of Definition~\ref{def:pointwise-coefficient-index}, by
\begin{equation}\label{eq:direct-relative-field-index}
 c_\pi(x)=\dim_F(L_x/q^NF\llbracket q \rrbracket^d)-dN
\end{equation}
using a determinant-normalized basis of $\Cot_q(\gamma_x)$,
where $\gamma_x$ is the loop induced by $x$ and $L_x=\ker\kappa_{\pi,x}$.
This integer is independent of the basis and of sufficiently
large $N$, and is unchanged by extension of the field $F$.
\end{lemma}
\begin{proof}
Write $\mathcal O(X)=\C[x_1,\ldots,x_m]/(f_1,\ldots,f_c)$, with
$g_i=x_i|_X$ the affine generators used in the coefficient presentation.
Choose $b\ge0$ so that on $U$
\[
 g_i(\gamma_U)=\sum_{\nu=0}^{k-1}\sum_{r\ge-b}
                      c_{i,\nu,r}z^\nu q^r .
\]
The determinant calculation for $\lambda_R$ holds over
$R\llbracket q\rrbracket$, so the basis $\beta_0,\ldots,\beta_{k-1}$
of $\mathcal O_a$ dual to $1,z,\ldots,z^{k-1}$ for this pairing
commutes with change of $R$. For $v\in\mathcal K_a$ it satisfies
\[
 \operatorname{Res}_{q=0}\lambda_R(q^{-r-1}\beta_\nu v)\,dq
       =[z^\nu q^r]v .
\]
Here $[z^\nu q^r]v$ denotes the coefficient of $q^r$ in the
$z^\nu$-component of $v$. Thus the free $R$-module on the symbols $dc_{i,\nu,r}$, $r\ge-b$,
identifies with $\mathcal K_a^m/q^b\mathcal O_a^m$ by
$dc_{i,\nu,r}\mapsto[q^{-r-1}\beta_\nu e_i]$, where $e_i$ is the
$i$-th standard basis vector.

Put $J_\alpha=((\partial f_\alpha/\partial x_i)(\gamma_V))_i$.
The cotangent presentation of $X$ identifies
$\Cot_q(\gamma_V)$ with
$\mathcal K_a^m/\sum_\alpha\mathcal K_a J_\alpha$ via the quotient map $$p:\mathcal K_a^m \to \Cot_q(\gamma_V),\quad e_i\mapsto \gamma_V^*(dg_i)\otimes dq.$$
By coefficient extraction, the differential of
$[z^\nu q^r]f_\alpha(\gamma_U)=0$, restricted to $V$, is represented
by $[q^{-r-1}\beta_\nu J_\alpha]$. These classes span the image of
$\sum_\alpha\mathcal K_a J_\alpha$ modulo $q^b\mathcal O_a^m$:
choose $h\ge0$ with all $J_\alpha\in q^{-h}\mathcal O_a^m$.
Truncating a multiplier in the $\beta_\nu$-basis modulo
$q^{b+h}\mathcal O_a$ gives a finite sum and changes its product
with $J_\alpha$ by an element of $q^b\mathcal O_a^m$.

Apply over $T$ the conormal sequence for a closed immersion,
as in the calculation for $S\hookrightarrow\Spec C_b(X)$
following \eqref{eq:elementary-coefficient-quotient} in the proof
of Lemma~\ref{lem:elementary-coefficient-differentials}.
Here the closed immersion is defined by the pullback of the ideal
$J_H$ of Definition~\ref{def:kv-bounded-family}. Its additional
cotangent relations are the differentials of the finitely many
coefficient generators of $J_H$, taken on the ambient coefficient
scheme over $T$ and then restricted to $S$.
Since $u$ is \'etale, $\Omega^1_{U/T}=u^*\Omega^1_{S/T}$.
Passing to $\Omega^1_{U/Y}$ further kills the differentials of
finitely many algebra generators of $\Gamma(Y,\mathcal O_Y)$.
Represent all these relations on $U$ by finite sums of coefficient
differentials. Right exactness of tensor product preserves these
cokernel presentations upon restriction to $V$. The resulting
quotient map is $\kappa_{\pi,V}$ by the coefficient extraction above.
It is therefore surjective, with $R$-module kernel
\[
 L_V=p(q^b\mathcal O_a^m)+\sum_{\ell=1}^{l}Rv_\ell ,
\]
where $v_1,\ldots,v_l$ are the images under $p$ of the chosen
representatives. These formulas commute with coefficientwise change
of $R$.

In a basis of $\Cot_q(\gamma_V)\cong R(\!( q )\!)^d$, choose a
section $\sigma$ of $p$. As in Lemma~\ref{lem:elementary-coefficient-differentials},
choose $a_0\ge0$ so that the entries of the matrices of $p,\sigma$
(using the basis $\beta_\nu e_i$ of $\mathcal K_a^m$) and of the
vectors $v_\ell$ have order at least $-a_0$. Then the kernel
formula gives the asserted bounds with $N=b+a_0$, also after any
change of $R$. Set $H_R=q^NR\llbracket q\rrbracket^d$. The sequence
\[
 0\longrightarrow L_V/H_R\longrightarrow\Cot_q(\gamma_V)/H_R
   \xrightarrow{\ \bar\kappa_{\pi,V}\ }\bigoplus_{j\in J}R\,du_j
   \longrightarrow0
\]
splits since its last term is free. Restricting a retraction onto
$L_V/H_R$ to the finite free submodule
$q^{-N}R\llbracket q\rrbracket^d/H_R$ shows that $L_V/H_R$ is
finite projective. For $R\to R'$, the classes $q^r\xi_i$ with $r<N$ and
$\xi_i$ in the chosen basis give an $R$-basis of $\Cot_q(\gamma_V)/H_R$ and an $R'$-basis after coefficient
change. Thus $\bar\kappa_{\pi,V}\otimes_R R'=\bar\kappa_{\pi,V'}$ in
these bases. Tensoring the split sequence proves the stated
base-change isomorphism. For a field extension this preserves the dimension in
\eqref{eq:direct-relative-field-index}. Independence of $N$ and
of the determinant-normalized basis follows from the calculation
after Definition~\ref{def:pointwise-coefficient-index}, with
parameter $q$ and rank $d$.
\end{proof}

For a finitely presented locally closed restriction in the
finite-factor presentation diagram, localizing and adding the
finitely many conormal relations of its defining equations gives
the same conclusions for its presented charts.

For $F$-subspaces $L,L_0$ with $L\cap L_0$ of finite codimension
in each, their \emph{relative dimension} is
\[
 [L:L_0]=\dim_F(L/H)-\dim_F(L_0/H),
 \qquad H=L\cap L_0.
\]
Any common subspace of finite codimension gives the same difference,
which is additive under direct sums. In a determinant-normalized
basis, \eqref{eq:direct-relative-field-index} becomes
$c_\pi(x)=[L_x:F\llbracket q\rrbracket^d]$.
The bounds of Lemma~\ref{lem:direct-relative-coefficient-kernel}
make $L_x$ an open linearly compact $F$-subspace. Thus $c_\pi(x)$
is the value on $L_x$ of the dimension theory normalized to vanish
on $F\llbracket q\rrbracket^d$, in the sense
of \cite[Section~3.5.1, p.~273, and Section~3.9, (3.9), p.~276]{Drinfeld}.

\begin{proposition}[Relative volume indices]
\label{prop:direct-relative-volume-indices}
The integers \eqref{eq:direct-relative-field-index} descend
to a locally constant function
$\delta_{Y/T}:Y\to\Z$.  These functions form a relative
volume-index normalization in the sense of
Definition~\ref{def:relative-volume-index}, on the whole
direct finite-factor presentation diagram.  It is uniquely characterized
by compatibility with $\Omega$ on the collision strata.
Explicitly, the following assertions hold.
\begin{enumerate}[(i),wide=0pt]
\item For a smooth comparison $f:Y'\to Y$ over $T$ in the
finite-factor presentation system, of relative dimension $e$, and a specified
base change
$b:T'\to T$,
\[
 \delta_{Y'/T}=f^*\delta_{Y/T}+e,
 \qquad
 \delta_{Y\times_TT'/T'}=\operatorname{pr}_Y^*\delta_{Y/T}.
\]
For a closed inclusion $i:Z\hookrightarrow Y$ in a finite-factor presentation with the same omitted coordinates,
$\delta_{Z/T}=i^*\delta_{Y/T}$.
\item Let $C_\pi\subset\A^I$ be the collision stratum for
a partition $I=\coprod_{j\in J}I_j$, and put
$P'_j=\sum_{i\in I_j}P_i$.  Under the canonical identification
\[
 \mathcal Z_{I,\mathbf P}|_{C_\pi}
       \cong C_\pi\times\prod_{j\in J}F_{P'_j}\mathscr LX,
\]
a product finite factor $C_\pi\times\prod_jY_j$ has index
\begin{equation}\label{eq:direct-stratum-volume-index}
 \delta_{C_\pi\times\prod_jY_j/C_\pi}
       =\sum_j\operatorname{pr}_j^*\delta_{Y_j}.
\end{equation}
Here $\delta_{Y_j}$ is the already defined pole-stage volume
index of Proposition~\ref{prop:finite-factor-normalization}.
The same formula holds after any specified base change of
$C_\pi$.
\item The assignments respect the separation, diagonal,
relabeling, and affine source-change identifications,
including their finite composition diagrams.  In the source
change $z\mapsto uz+v$, the parameters $u,u^{-1},v$ are
functions on the configuration base and contribute nothing
to the relative volume index.
\end{enumerate}
\end{proposition}
\begin{proof}
\emph{Step 1: local constancy and descent to the finite factor.}
We first show that for every algebraically closed field $F/\C$,
every affine $F$-scheme $D$ of finite type and every morphism
$h:D\to U$ over $\C$, the function $c_\pi\circ h$ is locally
constant. By the field-extension assertion of
Lemma~\ref{lem:direct-relative-coefficient-kernel}, its value
at a point of $D$ can be computed over that point's residue field.

On a reduced integral affine subscheme $\Spec R$ of $D$, let
$\gamma_R$ be the loop induced by $h$. Represent $\Cot_q(\gamma_R)$
as the image of an idempotent matrix over $R(\!( q )\!)$.
Since its rank is $d$, some $d\times d$ minor is nonzero over
$\operatorname{Frac}(R)(\!( q )\!)$. Invert its lowest nonzero
$q$-coefficient $c\in R$. The minor is a unit over $R_c(\!( q )\!)$,
so the selected rows split off the span of the selected columns
in the image of the idempotent. The complement has rank zero.
The columns therefore form a basis, which we normalize by the
determinant generator. For $V=\Spec R_c$, choose $N$ as in
Lemma~\ref{lem:direct-relative-coefficient-kernel}. On $V$,
\[
 c_\pi\circ h
   =\operatorname{rank}_{R_c}
       \bigl(L_V/q^NR_c\llbracket q\rrbracket^d\bigr)-dN
\]
by the base-change formula of that lemma. The quotient is finite
projective, so its rank is locally constant. Noetherian induction
now gives constructibility and finitely many values on $D$.

Put $A=F\llbracket s\rrbracket$ and let $\Spec A\to D$ be a
morphism. By the calculation in the proof of
Lemma~\ref{lem:coefficient-index-specialization}, $A(\!( q )\!)$ is a PID. Choose a
determinant-normalized basis of the Laurent cotangent module of
the induced loop and take $V=\Spec A$ in
Lemma~\ref{lem:direct-relative-coefficient-kernel}. For a common
bound $N$, the generic and special kernel quotients are the fibers
of the same finite free $A$-module
$L_V/q^NA\llbracket q\rrbracket^d$. Thus
\eqref{eq:direct-relative-field-index} has equal values at both
points. The curve argument of Step~4 of
Subsection~\ref{subsec:volume-index-proofs} applies over $F$:
a nonclosed constructible level set would give an integral curve
with generic point in that set and closed point outside it.
Normalizing the curve and completing at a point above the closed
point gives a morphism $\Spec F\llbracket s\rrbracket\to D$ with
the same generic value by field-extension invariance. Hence every level set is closed and, since there are
only finitely many, open. This proves the assertion about $h$.

Apply it to the zero section $s_0:Y\to U$ on affine opens of $Y$.
Then $\delta_{Y/T}:=c_\pi\circ s_0$ is locally constant. For a
geometric point $x=(y,(\alpha_j)_{j\in J})\in U(F)$, apply the
same assertion to the line
\[
 \A^1_F\longrightarrow U,\qquad
       \tau\longmapsto (y,(\tau\alpha_j)_{j\in J}).
\]
Its endpoints $s_0(y)$ and $x$ have the same index by connectedness
of $\A^1_F$. Field-extension invariance gives
$c_\pi=\pi^*\delta_{Y/T}$ on all scheme points.

\emph{Step 2: comparisons and part~(i).}
For a smooth comparison $f:Y'\to Y$ over $T$ with $\pi=f\pi'$,
fix a field-valued point $x$ of $U$. By
\eqref{eq:direct-relative-residue-map}, $\kappa_{\pi',x}$ is
$\kappa_{\pi,x}$ followed by the cotangent quotient
$(\Omega^1_{U/Y})|_x\to(\Omega^1_{U/Y'})|_x$. Hence, as in
Step~1 of Subsection~\ref{subsec:volume-index-proofs},
\[
 \ker\kappa_{\pi',x}/\ker\kappa_{\pi,x}
       \cong ((\pi')^*\Omega^1_{Y'/Y})|_x.
\]
Its dimension is the relative dimension $e$ of $f$ at $\pi'(x)$,
giving $c_{\pi'}(x)=c_\pi(x)+e$. For a specified base change or
a presented closed inclusion, corresponding points over the same
field have the same marked loop, determinant generator and
omitted-coordinate differentials. Formula~\eqref{eq:direct-relative-residue-map}
therefore gives the same kernel and index. Descent in Step~1
gives all three formulas in~(i).

For a finitely presented \'etale chart map $v:V\to U$ over $T$,
let $\rho:V\to Z$ be the finite-factor projection and $f:Z\to Y$
the smooth comparison satisfying $\pi v=f\rho$ in the
given presentation. Proposition~\ref{prop:direct-collision-diagrams}
provides the finite-factor refinements over $T$. On the refinements presenting $v$ with the same
omitted coordinates on both charts, \eqref{eq:direct-relative-residue-map} commutes
with \'etale pullback. The common-refinement calculation of Step~1
of Subsection~\ref{subsec:volume-index-proofs} thus applies over $T$:
$c_\rho=v^*c_\pi+\operatorname{reldim}(f)\circ\rho$.
Relative dimensions add under composition, giving the overlap
compatibilities of Definition~\ref{def:relative-volume-index}(i).

\emph{Step 3: the product formula on collision strata, part~(ii).}
Fix a field $F/\C$ and an $F$-point of $S$ whose configuration
$a=(a_i)_{i\in I}$ lies in $C_\pi(F)$. Write $c_j$
for the common value of the $a_i$ with $i\in I_j$, and put
$w_j=z-c_j$. The source factorizations of
Theorem~\ref{thm:kv-factorization} give
\[
 \mathcal O_a\cong\prod_{j\in J}F\llbracket w_j\rrbracket,\qquad
 \mathcal K_a\cong\prod_{j\in J}F(\!( w_j )\!).
\]
In the $j$-th factor,
$q=w_j^{|I_j|}\prod_{\ell\ne j}(w_j+c_j-c_\ell)^{|I_\ell|}$,
where the product over $\ell\ne j$ is a unit.
Write $\gamma_F\in X(\mathcal K_a)$ for the resulting loop and
$N_F=\gamma_F^*\Omega_X^1$. For $b\in \mathcal K_a$, let $b_j$ denote its component in
$F(\!( w_j )\!)$. The scalar source-residue functional
$b\mapsto\operatorname{Res}_q\lambda_F(b)\,dq$ satisfies
\begin{equation}\label{eq:direct-source-residues}
 \operatorname{Res}_q\lambda_F(b)\,dq
      =\sum_j\operatorname{Res}_{w_j}b_j\,dw_j.
\end{equation}
To check it, take $b=z^\nu q^m$, $0\le\nu<k$.
The right side is the sum of the finite residues of
$z^\nu D_a(z)^m dz$. For $m\ge0$ this form is polynomial.
For $m<0$, its Laurent expansion at infinity starts in
degree $\nu+km\le-1$ and has only lower degrees. Its
$z^{-1}dz$ coefficient is therefore zero except when
$\nu=k-1$ and $m=-1$, when it is $1$. The residue at
infinity is then $-1$, so the sum of the finite residues
is $1$, as is the left side. Only finitely many negative
$q$-powers can contribute, which extends the calculation
to Laurent series.

Let $\gamma_j$ be the loop obtained by projection to the
$j$th source factor. The source splitting defines an isomorphism
\[
 \Cot_q(\gamma_F)\longrightarrow\bigoplus_j\Cot_{w_j}(\gamma_j),\qquad
 b\,dg\otimes dq\longmapsto
       (b_j\,dg\otimes dw_j)_j.
\]
The target is an $F(\!( q )\!)$-module through the displayed
expressions for $q$. This map sends the formal symbol $dq$
to the formal disk symbols $dw_j$, including when $q$ has order
$|I_j|>1$ in the $j$-th disk.
Equation~\eqref{eq:direct-source-residues} shows that on a product
chart it intertwines the relative coefficient-residue map
\eqref{eq:direct-relative-residue-map} with the direct sum
of the disk coefficient-residue maps \eqref{eq:elementary-residue-map}.

Choose bases $e_{j,1},\ldots,e_{j,n}$ of
$\gamma_j^*\Omega_X^1$ whose wedges are
$\gamma_j^*\Omega$, and, under this isomorphism, put
\[
 \Lambda:=\bigoplus_j\bigoplus_{r=1}^n
                   F\llbracket w_j \rrbracket e_{j,r}\otimes dw_j.
\]
The disk bases combine to a $\mathcal K_a$-basis of $N_F$ whose
wedge is $\gamma_F^*\Omega$. The ordered basis $(1,z,\ldots,z^{k-1})$ and any
product integral basis of $\mathcal O_a$ are bases of the same
free $F\llbracket q \rrbracket$-module. Their determinants differ by a unit.
The determinant normalization of $\Cot_q(\gamma_F)$ therefore gives
order zero to the determinant of $\Lambda$.
More generally, for $b=(b_j)\in \mathcal K_a^\times$,
\begin{equation}\label{eq:kv-volume-norm-order}
 \operatorname{ord}_q\operatorname{Nm}_{\mathcal K_a/F(\!( q )\!)}(b)
                   =\sum_j\operatorname{ord}_{w_j}(b_j).
\end{equation}
Indeed, multiplication by $w_j$ on $F\llbracket w_j \rrbracket$ is
injective with cokernel $F$, of length one over
$F\llbracket q \rrbracket$; hence its determinant has $q$-order one.
Multiplication by a power-series unit has unit
determinant. Writing each $b_j$ as a power of $w_j$
times a unit and using multiplicativity proves the
formula, including negative powers. In particular, no
ramification multiplicity occurs.

After base change to the chosen point $\Spec F\to C_\pi$,
take a product chart with finite-factor projections
$\pi_j:U_j\to Y_j$. Let $x_j:\Spec F\to U_j$ be the
component points and put $x=(x_j)_j$.
The omitted cotangent quotient at $x$ is
$\bigoplus_j(\Omega^1_{U_j/Y_j})|_{x_j}$.
The preceding isomorphism therefore identifies
\[
 \ker\kappa_{\prod_j\pi_j,x}\cong\bigoplus_jL_j,
 \qquad L_j:=\ker q_{\pi_j,x_j}.
\]
Since the determinant of $\Lambda$ has order zero, taking relative
dimension against $\Lambda$ and using \eqref{eq:elementary-field-index}
gives $c_{\prod_j\pi_j}(x)=\sum_j c_{\pi_j}(x_j)$.
This proves \eqref{eq:direct-stratum-volume-index} on product
factors and, by (i), after the specified base changes.

\emph{Step 4: compatibility with the collision diagram, part~(iii).}
Under the identifications of Theorem~\ref{thm:kv-factorization} and
Lemma~\ref{lem:kv-pole-compatibility}, separation groups the disk
summands of Step~3, and both descriptions of a collision diagonal
give the same disk kernels and reference lattice $\Lambda$,
even when their monic degrees differ. Relabeling permutes the
summands and changes the determinant by at most a sign.

For an affine source change, let $V=\Spec R\to U$ be a test scheme
on the chart pulled back to the parameter base of the transformation.
Write $u\in R^\times$ and $v\in R$ for the values of its parameters, so
\[
 z'=uz+v,\qquad a_i'=ua_i+v,\qquad q'=u^kq.
\]
For the configuration $a'=(a_i')_{i\in I}$ in the coordinate $z'$,
these substitutions induce isomorphisms
\[
 \phi_q:R(\!( q' )\!)\xrightarrow{\sim}R(\!( q )\!),\qquad
 \phi_z:\mathcal K_{a'}\xrightarrow{\sim}\mathcal K_a.
\]
Let $\lambda'_R:\mathcal K_{a'}\to R(\!( q' )\!)$ extract the
coefficient of $(z')^{k-1}$. For $b'\in\mathcal K_{a'}$ and
$h'\in R(\!( q' )\!)$, coefficient extraction gives
\[
 \lambda_R(\phi_z(b'))
       =u^{k-1}\phi_q(\lambda'_R(b')),
 \qquad
 \operatorname{Res}_{q'}h'\,dq'
       =u^k\operatorname{Res}_q\phi_q(h')\,dq.
\]
Let $\gamma'_V\in X(\mathcal K_{a'})$ be the loop determined by
$\phi_z\circ(\gamma'_V)^*=\gamma_V^*$. Write $\kappa'$ and $\kappa$
for the relative coefficient-residue maps in the primed and
unprimed coordinates, respectively, with omitted-coordinate
cotangent quotients identified by the source change.
The $R$-linear, $\phi_q$-semilinear isomorphism
intertwining $\kappa'$ and $\kappa$ is
\[
 \Cot_{q'}(\gamma'_V)\longrightarrow\Cot_q(\gamma_V),\qquad
 b'\,dg\otimes dq'\longmapsto
             u\phi_z(b')\,dg\otimes dq.
\]
Indeed, $d_{U/Y}u=d_{U/Y}v=0$ because $u,v$ come from the base,
and the two coefficient identities give the asserted equality over $R$.

The image under $\phi_z$ of the ordered basis
$(1,z',\ldots,(z')^{k-1})$ has determinant $u^{k(k-1)/2}$
relative to $(1,z,\ldots,z^{k-1})$.
The displayed isomorphism also
multiplies each of its $kn$ cotangent basis vectors by
$u$. If $\varpi'$ and $\varpi$ denote the specified
determinant generators, their comparison is therefore
\[
                    \varpi'\longmapsto
                         u^{nk(k+1)/2}\varpi.
\]
This unit has order zero at every field point. The reference
lattices thus have relative dimension zero, and the identification
of the kernels gives equality of the indices. Under composition the
units $u$ multiply and the disk regroupings compose. With the
finite-factor comparisons of Proposition~\ref{prop:direct-collision-diagrams},
this proves (iii) on each finite composition diagram.

For uniqueness, every complex point lies on a collision stratum.
On a common smooth refinement with a product presentation, (ii)
fixes the index and (i) determines the index on the original factor.
These values determine a locally constant function on a finite-type
complex scheme.
\end{proof}

\begin{remark}[Relative volume indices on contact strata]
\label{rem:relative-contact-volume-indices}
Write $\operatorname{div}(\Omega)=-\sum_\alpha
a_\alpha D_\alpha$ and use the volume grading
$\mu_\Omega$ of \eqref{eq:volume-grading}.
For contact vectors $\mathbf v_j$ and sufficiently high jet
orders $\ell_j$, the factor
$Z=C_\pi\times\prod_jY_{\mathbf v_j,\ell_j}$ satisfies
\begin{equation}\label{eq:direct-contact-volume-index}
 \begin{aligned}
 \delta_{Z/C_\pi}
  &=\sum_j\left(n(\ell_j+1)-\sum_\alpha a_\alpha v_{j,\alpha}\right)\\
  &=\dim(Z/C_\pi)+\frac12\sum_j\mu_\Omega(\mathbf v_j).
 \end{aligned}
\end{equation}
This follows by summing the absolute contact formulas of
Proposition~\ref{prop:contact-volume-indices} in
\eqref{eq:direct-stratum-volume-index}.
For a general factor
$Y/T$, put $T_\pi=T\times_{\A^I}C_\pi$. Retain $e$
further coordinates and restrict to the scheme-theoretic
pullback to $T_\pi$ of the chosen product of absolute contact
loci, with the given scheme structure on $T_\pi$.
Write the restricted chart as $Z\times\A^{J'}$, with
$J'$ its omitted-coordinate set. Locally this chart is \'etale over
a product of smooth contact jet factors and affine tails pulled back
to $T_\pi$, so it is formally smooth over $T_\pi$.
The zero section of $Z\times\A^{J'}\to Z$ transfers the lifting
property to $Z/T_\pi$, so finite presentation gives smoothness. Under common smooth
refinements, (i) of Proposition~\ref{prop:direct-relative-volume-indices}
increases the index by the same amount as the relative dimension.
Applying \eqref{eq:direct-contact-volume-index} and subtracting the
$e$ retained coordinates therefore gives, at the corresponding point,
\[
 \delta_{Y/T}
   =\dim(Z/T_\pi)+\frac12\sum_j\mu_\Omega(\mathbf v_j)-e.
\]
The configuration-base dimension is not a summand.
The vectors $\mathbf v_j$ here are the actual contacts of the
loops on that stratum; only the pole bounds are added
on a collision diagonal.  
\end{remark}

\subsection{The fusion-ready package}
\label{subsec:all-multipoint-models}

We have now all the requirements to make the collision families fusion-ready.

\begin{proposition}
\label{lem:all-multipoint-models}
The families $\mathcal Z_{I,\mathbf P}\to\A^I$ admit
fusion-ready data in parity grading, without a volume
form. If a nowhere-vanishing volume form $\Omega$ is specified as in
Subsection~\ref{conv:loop-volume-compatibility}, the same finite-factor
presentations, equipped with the relative volume indices of
Proposition~\ref{prop:direct-relative-volume-indices}, give
fusion-ready data in integer grading. On a product finite factor
over a collision stratum, the relative index is the sum of the
pullbacks of the absolute pole-stage volume indices.
\end{proposition}
\begin{proof}
Theorem~\ref{thm:relative-finite-factors} and
Proposition~\ref{prop:direct-collision-diagrams} supply
Definition~\ref{def:fusion-ready-data}(i)--(ii).
When $\Omega$ is specified,
Proposition~\ref{prop:direct-relative-volume-indices} supplies (iii).
In parity grading, use the periodic purity comparisons of
Definition~\ref{def:relative-normalized-dual-complexes} and
Subsection~\ref{subsec:coherent-descent}. For a smooth
finite-factor comparison $h$ of relative dimension $e$ and
a two-periodic complex $K$ on its target, the comparison is
\[
 h^!K\simeq h^{-1}K[2e]\xrightarrow{\,u^{-e}\,}h^{-1}K,
\]
by \eqref{eq:oriented-smooth-purity}, with $u$ the periodicity
generator of degree two. Smooth purity respects base change and
composition, and $u^{-e}u^{-d}=u^{-(e+d)}$ for successive relative
dimensions $e,d$. This gives Definition~\ref{def:fusion-ready-data}(iii).
\end{proof}

\section{Normalized dualizing complexes of collision families}
\label{sec:relative-complexes}

We apply Subsection~\ref{subsec:relative-compact-support} to the
collision families. Their normalized dualizing complexes
are constructed on the total spaces; their direct images on the
configuration bases carry the factorizations used to define the
fusion product.

\subsection{Normalized complexes and finiteness}
\label{subsec:relative-pole-recovery}
\label{par:relative-pole-recovery}
In arguments common to the two gradings, use $B_P$ and $B$
for the derived pole stages and their homotopy colimit in parity
grading; in integer grading these symbols become
$B_{\Omega,P}$ and $B_\Omega$. Similarly, $\mathbf B$ and
$F_P\mathbf B\subseteq\mathbf B$ denote the chosen cochain model
and its filtration subcomplexes, with the subscript $\Omega$
understood in integer grading. Write $\omega_P$ for
$\omega_{P,\Omega}$ in integer grading and for $\omega_P^{\mathrm{per}}$
in parity grading. Set $B_{-1}=0$ and fix the
whole-diagram equivalence
$(F_P\mathbf B)_{P\ge0}\simeq(B_P)_{P\ge0}$ of
Theorem~\ref{thm:canonical-complex}; in particular, the inclusions
of subcomplexes represent the proper-trace pole transitions.

\begin{definition}[Normalized dualizing complexes for collision data]
\label{def:relative-fusion-complex}
Use the fixed pole-stage finite-factor presentations of
Section~\ref{sec:formal-loops} and the fusion-ready
data of Proposition~\ref{lem:all-multipoint-models}, including
their specified restrictions and induced finite-factor presentations on the
1-parameter binary slice.
In integer grading, use $\Omega$ and the relative volume indices
of Proposition~\ref{prop:direct-relative-volume-indices}; on a
pole-stage finite factor set
$\delta_{Y/\mathrm{pt}}:=\delta_Y$ as in
Proposition~\ref{prop:finite-factor-normalization}.
For each of these families $f:S\to T$, let
\[
 \omega^{\mathrm n}_{S/T}\in\mathcal D(S^{\mathrm{an}};\F),
 \qquad
 \mathcal B_{S/T},\ \mathcal E_{S/T}
       \in\mathcal D(T^{\mathrm{an}};\F)
\]
be the objects of
Definition~\ref{def:relative-normalized-dual-complexes}; thus
$\mathcal B_{S/T}=Rf_*\omega^{\mathrm n}_{S/T}$.
For a nonempty finite set $I$ and $\mathbf P\in\Z_{\ge0}^I$,
write
\[
 \begin{gathered}
 f_{I,\mathbf P}:\mathcal Z_{I,\mathbf P}\longrightarrow\A^I,\\
 \omega_{I,\mathbf P}:=
       \omega^{\mathrm n}_{\mathcal Z_{I,\mathbf P}/\A^I},\qquad
 \mathcal B_{I,\mathbf P}:=
       \mathcal B_{\mathcal Z_{I,\mathbf P}/\A^I},\qquad
 \mathcal E_{I,\mathbf P}:=
       \mathcal E_{\mathcal Z_{I,\mathbf P}/\A^I}.
 \end{gathered}
\]
Thus $\omega_{I,\mathbf P}$ is a normalized dualizing complex
on $\mathcal Z_{I,\mathbf P}^{\mathrm{an}}$, and
$\mathcal B_{I,\mathbf P}=Rf_{I,\mathbf P*}\omega_{I,\mathbf P}$
is its relative normalized Borel--Moore complex on
$(\A^I)^{\mathrm{an}}$.
In parity grading use the periodic construction of
Definition~\ref{def:relative-normalized-dual-complexes}, without
choosing $\Omega$ or relative volume indices. Use the same notation,
with the objects lying in $\mathcal D_{\mathrm{per}}(S^{\mathrm{an}};\F)$
and $\mathcal D_{\mathrm{per}}(T^{\mathrm{an}};\F)$, respectively.
\end{definition}

Proposition~\ref{lem:relative-compact-support} applies to these
finite-factor presentations and gives the following comparisons and finiteness
properties.

\begin{corollary}
\label{cor:relative-pole-recovery}
Use the objects and finite-factor presentation data of
Definition~\ref{def:relative-fusion-complex}.
In integer grading the following assertions hold.
\begin{enumerate}[(i)]
\item For the pole-stage family
$F_P\mathscr LX\to\mathrm{pt}$, $P\ge0$, there are equivalences
\[
 \omega^{\mathrm n}_{F_P\mathscr LX/\mathrm{pt}}
       \simeq\omega_{P,\Omega},\qquad
 \mathcal B_{F_P\mathscr LX/\mathrm{pt}}\simeq B_P
\]
in $\mathcal D((F_P\mathscr LX)^{\mathrm{an}};\F)$ and
$\mathcal D(\F)$, respectively. These equivalences identify
the closed pole traces and their coherent composition comparisons.

\item For each of the families
\[
 \begin{aligned}
 &F_P\mathscr LX\longrightarrow\mathrm{pt},\qquad P\ge0,\\
 &f_{I,\mathbf P}:\mathcal Z_{I,\mathbf P}\longrightarrow\A^I,
     \qquad \varnothing\ne I\text{ finite},\quad
                  \mathbf P\in\Z_{\ge0}^I,\\
 &\mathcal Z_{P,Q}\longrightarrow\A^1_s,\qquad P,Q\ge0,
 \end{aligned}
\]
write $f:S\to T$ for its structure map. Then
\[
 \mathcal E_{S/T},\ \mathcal B_{S/T}\in D^b_c(T^{\mathrm{an}};\F),
 \qquad
 \mathcal B_{S/T}\simeq\mathbb D_T\mathcal E_{S/T},\quad
 \mathcal E_{S/T}\simeq\mathbb D_T\mathcal B_{S/T}.
\]
Membership in $D^b_c$ refers to their images in the ordinary
derived category; the displayed equivalences are specified
in the derived enhancement.
\end{enumerate}
In parity grading, (i) holds with $\omega_P^{\mathrm{per}}$
in place of $\omega_{P,\Omega}$ and with the periodic categories.
In (ii), the relative images are constructible with
finite-dimensional stalk cohomology in each parity and satisfy
the same biduality.
\end{corollary}
\begin{proof}
For (i), the assignment $\delta_{Y/\mathrm{pt}}=\delta_Y$
makes the local sheaves, normalized purity equivalences, and
closed traces identical to those of
Proposition~\ref{prop:normalized-pole-sheaves}.
Descent of these coherent diagrams, including their arrows,
gives the sheaf equivalence. Applying $R\Gamma$ and
Definition~\ref{def:pole-stage-complexes} gives the equivalence
with $B_P$, compatibly with the pole maps.

The pole-stage finite nice filtration is supplied by
Theorem~\ref{thm:drinfeld-finite-factors}; the collision families
and their 1-parameter binary pullbacks have such filtrations by
Theorem~\ref{thm:relative-finite-factors} and
Proposition~\ref{lem:all-multipoint-models}.
To spell out the application of
Proposition~\ref{lem:relative-compact-support}(iv), write the
presented filtration as
$S=S_0\supset\cdots\supset S_m=\varnothing$ and its nice
differences as $V_a=S_a\setminus S_{a+1}
\simeq Y_a\times\A^{J_a}$, for $0\leq a<m$, with the map to
$T$ factoring through $g_a:Y_a\to T$. Their companions are
\[
 \mathcal E_{V_a/T}\simeq
 R(g_a)_!\bigl(\underline\F_{Y_a}[2\delta_{Y_a/T}]\bigr),
\]
which are bounded constructible because $g_a$ is algebraic
of finite type and $\delta_{Y_a/T}$ is locally constant
with finite image. The localization triangles
\[
 \mathcal E_{V_a/T}\longrightarrow\mathcal E_{S_a/T}
 \longrightarrow\mathcal E_{S_{a+1}/T}\xrightarrow{+1}
\]
prove the same assertion for $\mathcal E_{S/T}$ by finite
induction. The equivalence
$\mathcal B_{S/T}\simeq\mathbb D_T\mathcal E_{S/T}$ and
constructible Verdier duality give the assertion for
$\mathcal B_{S/T}$ and biduality.

In parity grading the companion of $V_a/T$ is
$\operatorname{Per}_2(R(g_a)_!\underline\F_{Y_a})$.
It has finite-dimensional stalk cohomology in each parity;
its periodic dual and biduality are obtained by periodizing
bounded constructible Verdier duality. Finite localization
triangles preserve constructibility, finite-dimensional
cohomology in each parity, and biduality. The same induction
therefore proves the parity assertions, without choosing
$\Omega$ or a global integer-graded lift of the periodic
system. No interchange of periodization with an infinite
\v Cech limit is used.
\end{proof}

\subsection{Factorization}
\label{subsec:relative-collision-complexes}

\begin{definition}[Coherent relative system]
\label{def:coherent-relative-system}
Use the families and complexes of
Definition~\ref{def:relative-fusion-complex}, including their
specified restrictions and refinements of finite-factor presentations. All sheaf
functors below act on analytic spaces, in the derived
$\infty$-categories $\mathcal D(-;\F)$ (or their parity versions).
A \emph{coherent relative system} is the following collection of
maps, equivalences, and compatibility homotopies.
\begin{enumerate}[(i),wide=0pt]
\item \emph{Pole traces and restrictions.}
For componentwise bounds $\mathbf P\le\mathbf Q$, let
$i_{\mathbf P,\mathbf Q}:\mathcal Z_{I,\mathbf P}\hookrightarrow
\mathcal Z_{I,\mathbf Q}$ be the closed pole inclusion. Specify
\[
 i_{\mathbf P,\mathbf Q}^!\omega_{I,\mathbf Q}
       \simeq\omega_{I,\mathbf P},\qquad
 \tau_{\mathbf P,\mathbf Q}:
 i_{\mathbf P,\mathbf Q*}\omega_{I,\mathbf P}
 \simeq i_{\mathbf P,\mathbf Q*}i_{\mathbf P,\mathbf Q}^!
                  \omega_{I,\mathbf Q}
 \xrightarrow{\ \epsilon\ }\omega_{I,\mathbf Q},
\]
where $\epsilon:i_*i^!\to\mathrm{id}$ is the adjunction counit.
For each open restriction $j:V\hookrightarrow \mathcal Z_{I,\mathbf P}$
in the collision diagram, specify
$j^{-1}\omega_{I,\mathbf P}\simeq\omega^{\mathrm n}_{V/\A^I}$,
with the induced finite-factor presentations and normalization.

For a surjection $\pi:I\twoheadrightarrow J$, put
$I_j=\pi^{-1}(j)$ and $P'_j=\sum_{i\in I_j}P_i$.
Let $\Delta_\pi:\A^J\hookrightarrow\A^I$ repeat the coordinate
$b_j$ at each $i\in I_j$, and write
\[
 \mathcal Z_\Delta:=\mathcal Z_{I,\mathbf P}\times_{\A^I,\Delta_\pi}\A^J,
 \qquad \widetilde\Delta_\pi:\mathcal Z_\Delta\hookrightarrow
                                      \mathcal Z_{I,\mathbf P}.
\]
Let $\kappa_\pi:\mathcal Z_\Delta\xrightarrow{\sim}\mathcal Z_{J,\mathbf P'}$
be the diagonal identification. Specify
\[
 \widetilde\Delta_\pi^!\omega_{I,\mathbf P}
                  \simeq\kappa_\pi^{-1}\omega_{J,\mathbf P'}.
\]
Here $i^!$ for a closed inclusion $i$ is the right adjoint of
$i_*$, so $i_*i^!$ is sections with support. The functor
$\kappa_\pi^{-1}$ is ordinary sheaf inverse image.

\item \emph{Products and source identifications.}
Let $U_\pi$ be the separation open of
Theorem~\ref{thm:kv-factorization}, and let
\[
 \phi_\pi:\mathcal Z_{I,\mathbf P}|_{U_\pi}
       \xrightarrow{\sim}
       \left.\prod_{j\in J}\mathcal Z_{I_j,\mathbf P|_{I_j}}\right|_{U_\pi}
\]
be the geometric separation isomorphism. Specify
\[
 \omega_{I,\mathbf P}|_{\mathcal Z_{I,\mathbf P}|_{U_\pi}}
 \simeq\phi_\pi^{-1}
       \left(\left.\boxtimes_{j\in J}
                    \omega_{I_j,\mathbf P|_{I_j}}\right|_{U_\pi}\right).
\]
We call this equivalence the \emph{normalized sheaf separation
comparison associated with $\pi$}; its underlying geometric
separation isomorphism is $\phi_\pi$. On $U_\pi$, the groups
$I_j$ are separated from one another, while points within each
group may still coincide. The comparison is equipped with the
coherence specified in part~(iii).
For a bijection $\sigma:I\xrightarrow{\sim}I'$, set
$(\sigma_*\mathbf P)_{\sigma(i)}=P_i$. The relabeling
isomorphism $r_\sigma:\mathcal Z_{I,\mathbf P}\xrightarrow{\sim}
\mathcal Z_{I',\sigma_*\mathbf P}$ is equipped with
\[
 r_\sigma^{-1}\omega_{I',\sigma_*\mathbf P}
                 \simeq\omega_{I,\mathbf P}.
\]
For affine source changes put $G=\A^1_b\times\Spec\C[u,u^{-1}]$ and let
\[
 A_I:G\times \mathcal Z_{I,\mathbf P}\xrightarrow{\sim}
                         G\times \mathcal Z_{I,\mathbf P},
 \qquad (b,u,a,\gamma)\longmapsto(b,u,b+ua,\gamma'),
\]
where $\gamma'(g)(z)=\gamma(g)((z-b)/u)$ for $g\in\mathcal O(X)$.
The source-algebra map uses $D_{b+ua}(z)=u^{|I|}D_a((z-b)/u)$.
Specify
\[
 A_I^{-1}(\omega_G\boxtimes\omega_{I,\mathbf P})
       \simeq\omega_G\boxtimes\omega_{I,\mathbf P},
\]
with the inverse and composition comparisons over powers of $G$.
The finite-factor presentations and relative volume-index invariance are
supplied by Propositions~\ref{prop:direct-collision-diagrams} and
\ref{prop:direct-relative-volume-indices}.

\item \emph{Coherence.}
Specify (i)--(ii) as natural transformations of the diagrams of
complexes on compatible finite-factor presentations, together
with their composition and
base-change homotopies, before descent. For example, for
$\mathbf P\le\mathbf Q\le\mathbf R$, specify a homotopy
\[
 H_{\mathbf P,\mathbf Q,\mathbf R}:
 \tau_{\mathbf P,\mathbf R}
 \simeq \tau_{\mathbf Q,\mathbf R}\circ
               i_{\mathbf Q,\mathbf R*}\tau_{\mathbf P,\mathbf Q}
\]
in the mapping space
\[
 \operatorname{Map}_{\mathcal D(\mathcal Z_{I,\mathbf R}^{\mathrm{an}};\F)}
 \bigl(i_{\mathbf P,\mathbf R*}\omega_{I,\mathbf P},
                                      \omega_{I,\mathbf R}\bigr),
\]
using $i_{\mathbf P,\mathbf R*}\simeq
i_{\mathbf Q,\mathbf R*}i_{\mathbf P,\mathbf Q*}$.
For four successive bounds, specify a homotopy between the two
paths from the direct trace to the threefold composite; for
longer strings specify the higher homotopies, compatible with
all faces and identities. Require the same coherence for
successive separations, diagonal restrictions, relabelings,
affine source comparisons, and their mixed compositions with
pole maps and open restrictions. These requirements mean
coherent diagrams in the sense of
Appendix~\ref{app:coherent-diagrams}, and include compatibility
with common refinements of finite-factor presentations.
\end{enumerate}
\end{definition}

The next lemma constructs this system and its induced systems on
the configuration bases. Their coherence is used to compare
successive collisions and source-rotation parameter families.

\begin{lemma}
\label{lem:relative-fusion-factorization}
Use the notation of Definition~\ref{def:coherent-relative-system}.
\begin{enumerate}[(i),wide=0pt]
\item The normalized dualizing complexes of
Definition~\ref{def:relative-fusion-complex} admit the maps,
equivalences, and homotopies required by
Definition~\ref{def:coherent-relative-system}. They are induced
by the normalized purity equivalences on the finite factors,
together with the adjunction and base-change transformations,
and descend with their sources and targets.

\item On the finite-dimensional bases, direct image gives the
covariant pole maps
\[
 u_{\mathbf P,\mathbf Q}:=
 Rf_{I,\mathbf Q*}\tau_{\mathbf P,\mathbf Q}:
 \mathcal B_{I,\mathbf P}\longrightarrow\mathcal B_{I,\mathbf Q},
\]
using $f_{I,\mathbf P}=f_{I,\mathbf Q}i_{\mathbf P,\mathbf Q}$.
Their Verdier duals are the contravariant pole maps
\[
 \mathbb D_{\A^I}u_{\mathbf P,\mathbf Q}:
 \mathcal E_{I,\mathbf Q}\longrightarrow\mathcal E_{I,\mathbf P}.
\]
For a base change $b:T'\to T$ between separated finite-type
complex bases, the induced system uses
$b^!\mathcal B_{S/T}\simeq\mathcal B_{S\times_TT'/T'}$
and $b^*\mathcal E_{S/T}\simeq\mathcal E_{S\times_TT'/T'}$
from \eqref{eq:relative-dual-base-change}, whenever the pulled-back
finite-factor presentations are used. In particular, separation and diagonal
restriction give
\begin{equation}\label{eq:relative-factorization-complexes}
 \begin{aligned}
 \Delta_\pi^!\mathcal B_{I,\mathbf P}
      &\simeq\mathcal B_{J,\mathbf P'},\\
 \left.\mathcal B_{I,\mathbf P}\right|_{U_\pi}
      &\simeq
       \left.\boxtimes_{j\in J}
       \mathcal B_{I_j,\mathbf P|_{I_j}}\right|_{U_\pi}.
 \end{aligned}
\end{equation}
The external product has base
$\prod_j\A^{I_j}=\A^I$. For a singleton,
\[
 \mathcal B_{\{1\},P}\simeq\omega_{\A^1}\otimes B_P,
 \qquad
 \mathcal E_{\{1\},P}\simeq
           \underline\F_{\A^1}\otimes B_P^\vee.
\]
At a point $a$ where precisely the groups $I_j$ coincide,
with inclusion $i_a:\{a\}\hookrightarrow\A^I$, one has
\begin{equation}\label{eq:relative-stalks}
 i_a^!\mathcal B_{I,\mathbf P}
       \simeq\bigotimes_{j\in J}B_{P'_j},
 \qquad
 i_a^*\mathcal E_{I,\mathbf P}
       \simeq\bigotimes_{j\in J}(B_{P'_j})^\vee.
\end{equation}
Both complexes are constructible for the diagonal
stratification. The dual form of
\eqref{eq:relative-factorization-complexes} for $\mathcal E$
uses ordinary inverse image. Its separation equivalence is
\[
 \left.\mathcal E_{I,\mathbf P}\right|_{U_\pi}
 \simeq
 \left.\boxtimes_{j\in J}
       \mathcal E_{I_j,\mathbf P|_{I_j}}\right|_{U_\pi}.
\]
We call this equivalence and the second equivalence in
\eqref{eq:relative-factorization-complexes} the
\emph{separation comparisons on the base $U_\pi$}.
They are induced by the normalized sheaf separation comparison
of Definition~\ref{def:coherent-relative-system}(ii).
We use the same terminology for their pullbacks under the
base-change hypotheses in this part, using $b^*$ for
$\mathcal E$ and $b^!$ for $\mathcal B$.
For topological base changes, we use the hypotheses and
comparisons of
Lemma~\ref{lem:topological-pullback-presentations}(i)--(ii).

\item Let $\sigma:\A^1_s\hookrightarrow\A^{\{1,2\}}$,
$\sigma(s)=(s,-s)$, and use the 1-parameter binary collision family
\[
 \mathcal Z_{P,Q}
   =\mathcal Z_{\{1,2\},(P,Q)}\times_{\A^{\{1,2\}},\sigma}\A^1_s.
\]
Its relative complexes are
\[
 \mathcal B_{P,Q}:=\mathcal B_{\mathcal Z_{P,Q}/\A^1_s}
       \simeq\sigma^!\mathcal B_{\{1,2\},(P,Q)},\qquad
 \mathcal E_{P,Q}:=\mathcal E_{\mathcal Z_{P,Q}/\A^1_s}
       \simeq\sigma^*\mathcal E_{\{1,2\},(P,Q)}.
\]
For $j:\Gm\hookrightarrow\A^1_s$ and
$i:\{0\}\hookrightarrow\A^1_s$, these satisfy
\begin{equation}\label{eq:binary-slice-factorization}
 \begin{aligned}
 j^*\mathcal B_{P,Q}
     &\simeq\omega_{\Gm}\otimes(B_P\otimes B_Q),
 &i^!\mathcal B_{P,Q}&\simeq B_{P+Q},\\
 j^*\mathcal E_{P,Q}
     &\simeq\underline{(B_P\otimes B_Q)^\vee}_{\Gm},
 &i^*\mathcal E_{P,Q}&\simeq B_{P+Q}^\vee.
 \end{aligned}
\end{equation}
More precisely, let
\[
 \begin{gathered}
 j_{P,Q}:\mathcal Z_{P,Q}|_{\Gm}\hookrightarrow\mathcal Z_{P,Q},
 \qquad
 i_{P,Q}:F_{P+Q}\mathscr LX\xrightarrow{\sim}
             (\mathcal Z_{P,Q})_0\hookrightarrow\mathcal Z_{P,Q},\\
 \phi_{P,Q}:\mathcal Z_{P,Q}|_{\Gm}
       \xrightarrow{\sim}\Gm\times F_P\mathscr LX\times F_Q\mathscr LX
 \end{gathered}
\]
be the open inclusion and the fiber identifications.
Then
\begin{equation}\label{eq:binary-slice-sheaf-factorization}
 \begin{aligned}
 j_{P,Q}^{-1}\omega^{\mathrm n}_{\mathcal Z_{P,Q}/\A^1_s}
  &\simeq\phi_{P,Q}^{-1}
      (\omega_{\Gm}\boxtimes\omega_P
                          \boxtimes\omega_Q),\\
 i_{P,Q}^!\omega^{\mathrm n}_{\mathcal Z_{P,Q}/\A^1_s}
  &\simeq\omega_{P+Q}.
 \end{aligned}
\end{equation}
These equivalences are compatible with the relative direct images
and Verdier duality giving \eqref{eq:binary-slice-factorization}.
As $P,Q$ vary, they are equivalences of the entire coherent
pole diagrams, including the pole traces, their dual restriction
maps, and comparisons of finite-factor presentations.
\end{enumerate}
The relative-image and compact-support constructions preserve
the homotopies of Definition~\ref{def:coherent-relative-system}(iii).
Thus the displayed maps and equivalences form coherent systems
for $\mathcal B$ and $\mathcal E$, respectively; their coherence
is a conclusion of the lemma, not additional data to be chosen.
\end{lemma}

\begin{proof}
The maps of Lemma~\ref{lem:kv-pole-compatibility},
including their finite composition diagrams and affine source
comparisons, have compatible finite-factor presentations and volume indices by
Propositions~\ref{prop:direct-collision-diagrams} and
\ref{lem:all-multipoint-models}.
On finite factors, for a closed inclusion $i:Y'\hookrightarrow Y$
and an open inclusion $j:V\hookrightarrow Y$, use
\[
 i^!\omega_Y\simeq\omega_{Y'},\qquad
 i_*\omega_{Y'}\simeq i_*i^!\omega_Y
                   \xrightarrow{\epsilon}\omega_Y,
 \qquad j^{-1}\omega_Y\simeq\omega_V.
\]
These give the closed restrictions, pole traces, and open
restrictions in Definition~\ref{def:coherent-relative-system}(i).
Exterior products give separation; transport along the
relabeling isomorphisms gives the relabeling comparisons.
Proposition~\ref{prop:direct-relative-volume-indices} gives
additivity under products, restriction under the specified closed
base changes, and invariance under affine source changes.
A smooth refinement of relative dimension $e$ increases the
volume index by $e$, canceling the purity shift $[2e]$. Pullback to the infinite
charts and descent of these entire diagrams, as in
Proposition~\ref{lem:relative-compact-support}, prove (i).

For each identity between factorization maps, use
the corresponding composition or base-change transformation
of the enhanced finite-factor functors. These transformations
already include the higher compatibilities. Transporting
them by exact inverse image and descent constructs condition
(iii) of Definition~\ref{def:coherent-relative-system}; equality
of the resulting maps in an ordinary derived category would
not suffice. For the pole maps, $H_{\mathbf P,\mathbf Q,\mathbf R}$
is the compatibility of the adjunction counit with composition;
the higher homotopies are those of the same enhanced adjunctions.
Diagram descent preserves them, including the mixed comparisons.
Explicitly, for a simplicial set $K$ indexing one of these
coherent diagrams on a total space $S$, and a nice \'etale
\v Cech nerve $W_\bullet\to S$, descent gives
\[
 \operatorname{Fun}\bigl(K,\mathcal D(S^{\mathrm{an}};\F)\bigr)
 \simeq\varprojlim_{[n]\in\Delta}
 \operatorname{Fun}\bigl(K,\mathcal D(W_n^{\mathrm{an}};\F)\bigr).
\]
Indeed, functor categories preserve limits. This descends
the compatible diagrams with their arrows and higher
simplices. The same argument applies in the periodic categories.

Applying $Rf_*$ to these diagrams gives $u_{\mathbf P,\mathbf Q}$
and its composition homotopies. On the compact-support finite
terms the pole map is restriction from the larger finite factor
to the smaller one; its Verdier dual is the displayed proper
trace. Descent and Corollary~\ref{cor:relative-pole-recovery}(ii)
therefore identify these maps with $\mathbb D_{\A^I}u_{\mathbf P,\mathbf Q}$.
For $b:T'\to T$, put $S'=S\times_TT'$ and use the
pulled-back finite-factor presentations. Proper-support base change on
each finite factor, followed by the realization defining
$\mathcal E$, gives
$b^*\mathcal E_{S/T}\simeq\mathcal E_{S'/T'}$; ordinary
inverse image preserves that realization. Consequently
\[
 \begin{aligned}
 b^!\mathcal B_{S/T}
   &\simeq b^!\mathbb D_T\mathcal E_{S/T}
    \simeq\mathbb D_{T'}b^*\mathcal E_{S/T}\\
   &\simeq\mathbb D_{T'}\mathcal E_{S'/T'}
    \simeq\mathcal B_{S'/T'}.
 \end{aligned}
\]
Here the duality comparison is
\eqref{eq:duality-inverse-image} on the finite-dimensional
bases. Together with exterior-product compatibility in
Proposition~\ref{lem:relative-compact-support}, these
comparisons give (ii). They are transformations of the whole
finite-factor diagrams, so realization and duality preserve
the composition homotopies constructed in (i).
For a singleton, $t=z-a_1$ identifies the collision family with
$\A^1\times F_P\mathscr LX$, and its normalized dualizing
complex is $\omega_{\A^1}\boxtimes\omega_P$. Direct image gives
the singleton formula, and duality gives its companion.

Let $V_J\subset\A^J$ be the distinct-point locus. Identify
the diagonal stratum $\Delta_\pi(V_J)$ with $V_J$ and write
$j_\pi:V_J\hookrightarrow\A^I$ for its locally closed
inclusion. Diagonal restriction, separation, and the singleton
calculation give
\[
 j_\pi^!\mathcal B_{I,\mathbf P}
       \simeq\omega_{V_J}\otimes\bigotimes_jB_{P'_j},
 \qquad
 j_\pi^*\mathcal E_{I,\mathbf P}
       \simeq\underline\F_{V_J}\otimes
                              \bigotimes_j(B_{P'_j})^\vee.
\]
Taking a point costalk or stalk proves
\eqref{eq:relative-stalks}. The second formula and
Corollary~\ref{cor:relative-pole-recovery}(ii) prove
constructibility of $\mathcal E$ for the diagonal stratification.
This is a Whitney stratification: each stratum is open in the
complex linear subspace $L\subset\C^I$ defined by its coordinate
equalities. For a stratum in its closure, the corresponding
subspace $L'$ satisfies $L'\subset L$. The real tangent spaces
of the first stratum are the underlying real vector space of
$L$, containing those of the second stratum; secant lines
between points of the two strata also lie in $L$. These
inclusions persist under limits and give Whitney conditions
(a) and (b). Constructible Verdier duality therefore gives
the same diagonal constructibility for $\mathcal B$, in integer
and parity grading. In particular, ordinary stalks of $\mathcal B$
are not being identified with the unshifted fiber complexes.

For (iii), apply the base-change equivalences of (ii) to
$\sigma$. On $\Gm$, separation and the singleton formulas give
\[
 (\sigma j)^!(\omega_{\A^2}\otimes B_P\otimes B_Q)
       \simeq\omega_{\Gm}\otimes B_P\otimes B_Q.
\]
At $s=0$, \eqref{eq:relative-stalks} gives $B_{P+Q}$.
Duality gives the two formulas for $\mathcal E_{P,Q}$.
On the total spaces the same argument uses the closed
base-change comparison
\[
 \widetilde\sigma^!\omega_{\{1,2\},(P,Q)}
       \simeq\omega^{\mathrm n}_{\mathcal Z_{P,Q}/\A^1_s},
 \qquad
 \widetilde\sigma:\mathcal Z_{P,Q}
           \hookrightarrow\mathcal Z_{\{1,2\},(P,Q)},
\]
and then the separation and special-fiber identifications in (i).
Restriction of the relative volume indices and their additivity
make these the normalized comparisons
\eqref{eq:binary-slice-sheaf-factorization}. Their direct images
are the formulas already obtained on the base.
All these operations are applied to the coherent systems of
(i)--(ii); hence they preserve the pole maps, duality, and
comparisons of finite-factor presentations, with their coherence homotopies.
\end{proof}

\subsection{Source coordinates and factorization}
\label{subsec:source-coordinates-factorization}
We make the factorization comparisons explicit in the source coordinates
$z=b_j+u_jt_j$. This will be help prove  properties of the fusion product and of the BV operator. We shall refer to the $u_j$ below as {\it multipliers}.

\begin{lemma}
\label{lem:source-coordinate-factorization}
Fix a nonempty finite set $I$ and $\mathbf P\in\Z_{\ge0}^I$.
Use the collision families of Definition~\ref{def:kv-bounded-family}
and the coherent relative system of
Lemma~\ref{lem:relative-fusion-factorization}, with the complexes,
gradings, and fixed finite-factor presentation data of
Definition~\ref{def:relative-fusion-complex}.
Use $\omega_P$ and $B_P$ as in
Subsection~\ref{subsec:relative-pole-recovery}.
Write $r_0$ for the integer $r$ in \eqref{eq:kv-marked-poles},
and $H_{a,\mathbf P}=\prod_i(z-a_i)^{r_0P_i}$. The source algebras $\mathcal O_a,\mathcal K_a$
and polynomial $D_a$ are those of
Definition~\ref{def:kv-multipoint-functor}; the affine generators
$g_\ell$ are those fixed in Definition~\ref{def:kv-bounded-family}.

Let $\pi:I\twoheadrightarrow J$ be a surjection onto a
nonempty finite set, and define
\[
 \begin{gathered}
 I_j=\pi^{-1}(j),\qquad P'_j=\sum_{i\in I_j}P_i,\qquad
 V_J=\operatorname{Conf}_J\subset\A^J,\qquad
 T_\pi=V_J\times(\Gm)^J,\\
 \Delta_\pi:\A^J\hookrightarrow\A^I,
       \quad(b_j)_j\longmapsto(b_{\pi(i)})_i,\qquad
 h_\pi=\Delta_\pi\circ\operatorname{pr}_{V_J}:T_\pi\to\A^I,\\
 S_\pi=\mathcal Z_{I,\mathbf P}\times_{\A^I,h_\pi}T_\pi,
 \qquad F_{\pi,\mathbf P}=\prod_{j\in J}F_{P'_j}\mathscr LX,
 \qquad B_{\pi,\mathbf P}=\bigotimes_{j\in J}B_{P'_j}.
 \end{gathered}
\]
Thus $V_J(R)=\{b\in R^J:b_j-b_{j'}\in R^\times\ (j\ne j')\}$;
write $(b,u)=((b_j),(u_j))$ for the universal coordinates on
$T_\pi$. For a complex $V$ write
$V^\vee=R\operatorname{Hom}_{\F}(V,\F)$ and
$\underline V_T=\underline\F_T\otimes V$; tensor products use
the fixed Koszul symmetry.
Then the following assertions hold in integer and parity grading.
\begin{enumerate}[(i),wide=0pt]
\item For every $\C$-algebra $R$ and $(b,u)\in T_\pi(R)$,
put $a_i=b_{\pi(i)}$. The substitutions $z=b_j+u_jt_j$ give
natural isomorphisms
\begin{equation}\label{eq:source-coordinate-crt}
 \begin{aligned}
 \mathcal O_a&\xrightarrow{\sim}\prod_{j\in J}R\llbracket t_j \rrbracket,&
 \mathcal K_a&\xrightarrow{\sim}\prod_{j\in J}R(\!( t_j )\!),\\[-2pt]
 f&\longmapsto\bigl(f(b_j+u_jt_j)\bigr)_j.
 \end{aligned}
\end{equation}
The completions have their adic topologies. On the punctured
algebras use the neighborhood bases of zero
$\{D_a^N\mathcal O_a\}_{N\ge0}$ and
$\{\prod_jt_j^NR\llbracket t_j \rrbracket\}_{N\ge0}$, respectively.
The displayed maps are isomorphisms of topological $R$-algebras. In each factor,
$t_j=(z-b_j)/u_j$. They identify the marked pole condition
$H_{a,\mathbf P}\gamma(g_\ell)\in\mathcal O_a$ with
$t_j^{r_0P'_j}\gamma_j(g_\ell)\in R\llbracket t_j \rrbracket$ for every
$j,\ell$, and induce an isomorphism of schemes over $T_\pi$:
\begin{equation}\label{eq:source-coordinate-identification}
 \begin{aligned}
 \Phi_{\pi,u}:S_\pi&\xrightarrow{\sim}
                         T_\pi\times F_{\pi,\mathbf P},\\
 (b,u,\gamma)&\longmapsto(b,u,(\gamma_j)_j),\qquad
 \gamma_j(g)(t_j)=\gamma(g)(b_j+u_jt_j)
                         \quad(g\in\mathcal O(X)).
 \end{aligned}
\end{equation}
The substitution means the $j$th component of
\eqref{eq:source-coordinate-crt}. In particular, $|I|=|J|=1$
gives the singleton identification, and $|J|=1$ gives
$S_\pi\simeq T_\pi\times F_{\sum_iP_i}\mathscr LX$.

\item On analytic spaces, the isomorphism in (i) gives
\begin{equation}\label{eq:source-coordinate-comparisons}
 \begin{aligned}
 \omega^{\mathrm n}_{S_\pi/T_\pi}
   &\simeq\Phi_{\pi,u}^{-1}
       \left(\omega_{T_\pi}\boxtimes
                          \boxtimes_j\omega_{P'_j}\right),\\
 e_{\pi,u}:=(\Phi_{\pi,u})_!:
 \mathcal E_{S_\pi/T_\pi}
   &\xrightarrow{\sim}\underline{B_{\pi,\mathbf P}^{\vee}}_{T_\pi},\\
 b_{\pi,u}:=(\Phi_{\pi,u})^!:
 \omega_{T_\pi}\otimes B_{\pi,\mathbf P}
   &\xrightarrow{\sim}\mathcal B_{S_\pi/T_\pi},
 \qquad b_{\pi,u}=\mathbb D_{T_\pi}(e_{\pi,u}).
 \end{aligned}
\end{equation}
We call $e_{\pi,u}$ and $b_{\pi,u}$ the
\emph{source-coordinate comparisons}. They combine the
diagonal identification for $\pi$, separation of the groups
$I_j$, and the singleton identifications in the coordinates
$z=b_j+u_jt_j$.
The first equivalence is in
$\mathcal D(S_\pi^{\mathrm{an}};\F)$ and the last two in
$\mathcal D(T_\pi^{\mathrm{an}};\F)$, or their periodic
versions. In the notation $(\Phi_{\pi,u})_!$, the operation is
transport of the normalized finite-factor trace diagram along
$\Phi_{\pi,u}$, followed by its product comparison and homotopy
realization. The notation $(\Phi_{\pi,u})^!$ denotes the
Verdier dual of this map on $T_\pi^{\mathrm{an}}$.
These symbols denote maps of complexes, not exceptional functors
on the infinite-dimensional total spaces. Their base-change
comparisons are those of
Lemma~\ref{lem:relative-fusion-factorization}(ii).

\item Put
\[
 \widehat T_\pi=V_J\times(\Gm)^J_u\times(\Gm)^J_\lambda,
 \quad p(b,u,\lambda)=(b,u),\quad
 m(b,u,\lambda)=(b,u\lambda).
\]
Since $h_\pi p=h_\pi m$, the families $p^*S_\pi$ and
$m^*S_\pi$ are canonically the same family
$\widehat S_\pi$ over $\widehat T_\pi$. Define the automorphism
of $\widehat T_\pi\times F_{\pi,\mathbf P}$ by
\[
 r_\lambda(b,u,\lambda,(\eta_j)_j)
     =(b,u,\lambda,(r_{\lambda_j}\eta_j)_j),\qquad
 (r_{\lambda_j}\eta_j)(g)(t)=\eta_j(g)(\lambda_jt).
\]
Denote the universal maps in (ii) by $e_\pi$ and $b_\pi$.
On $\widehat T_\pi$ write
$e_{\pi,u}=p^*e_\pi$, $e_{\pi,u\lambda}=m^*e_\pi$
and $b_{\pi,u}=p^!b_\pi$, $b_{\pi,u\lambda}=m^!b_\pi$,
using the base-change comparisons in (ii).
In particular, both $b$-maps have source
$\omega_{\widehat T_\pi}\otimes B_{\pi,\mathbf P}$.
The finite-factor transport of $r_\lambda$ and its Verdier dual
are denoted by $(r_\lambda)_!$ and $r_\lambda^!$, respectively.
Then
\begin{equation}\label{eq:source-coordinate-change}
 \begin{aligned}
 \Phi_{\pi,u\lambda}&=r_\lambda\circ\Phi_{\pi,u},\\
 e_{\pi,u\lambda}&=(r_\lambda)_!\circ e_{\pi,u},\\
 b_{\pi,u\lambda}&=b_{\pi,u}\circ r_\lambda^!.
 \end{aligned}
\end{equation}
These are identities of the geometric maps and coherent
comparisons of the sheaf maps over $\widehat T_\pi$.

For the affine source-coordinate map
$\varphi_{b,u}(t)=b+ut$, $u\in R^\times$, one has
\begin{equation}\label{eq:source-coordinate-composition}
 \varphi_{b,u}\circ\varphi_{c,v}=\varphi_{b+uc,uv}
                  \qquad(v\in R^\times).
\end{equation}
The affine source action of
Definition~\ref{def:coherent-relative-system}(ii) induces
$\widetilde A:G\times S_\pi\to S_\pi$, $G=\A^1\times\Gm$,
over
\[
 \alpha:G\times T_\pi\longrightarrow T_\pi,
 \qquad ((\beta,v),(b,u))\longmapsto(\beta+vb,vu).
\]
For the universal map $\Phi_\pi=\Phi_{\pi,u}$ in (i),
\[
 \Phi_\pi\circ\widetilde A
  =(\alpha\times\mathrm{id}_{F_{\pi,\mathbf P}})
                  \circ(\mathrm{id}_G\times\Phi_\pi).
\]
Both sides cover $\alpha$; on the left the source coordinates
are therefore evaluated at $(\beta+vb,vu)$.

The comparisons are natural for the closed pole maps
$i_{\mathbf P,\mathbf Q}$ and their traces
$\tau_{\mathbf P,\mathbf Q}$, relabelings $r_\sigma$,
separations $\phi_\rho$, and diagonal identifications
$\kappa_\rho$ of Definition~\ref{def:coherent-relative-system}.
For example, if $\rho:J\twoheadrightarrow K$ and
$P''_k=\sum_{j\in\rho^{-1}(k)}P'_j$, their diagonal
compatibility is the comparison
\[
 \kappa_{\rho\pi}
  =\kappa_\rho\circ\Delta_\rho^*\kappa_\pi:
 \mathcal Z_{I,\mathbf P}\times_{\A^I}\A^K
                  \xrightarrow{\sim}\mathcal Z_{K,\mathbf P''},
\]
where $\Delta_{\rho\pi}=\Delta_\pi\Delta_\rho$ and
$\Delta_\rho^*\kappa_\pi$ is first viewed as a map to
$\mathcal Z_{J,\mathbf P'}\times_{\A^J}\A^K$.
Form this comparison on the full diagonal families, then use
$\Phi$ on each collision stratum separately. Regrouping the
already separated factors on $T_\pi$ instead uses
\[
 \prod_{j\in J}R(\!( t_j )\!)\simeq
 \prod_{k\in K}\prod_{j\in\rho^{-1}(k)}R(\!( t_j )\!);
\]
it does not merge their centers. All these transported maps
preserve the composition and base-change homotopies of
Lemma~\ref{lem:relative-fusion-factorization}.

\item Let $v:T\to T_\pi^{\mathrm{an}}$ satisfy the topological
pullback hypotheses of
Lemma~\ref{lem:topological-pullback-presentations}(ii) for the
chosen finite-factor presentations. Put
\[
 S_T=S_\pi^{\mathrm{an}}\times_{T_\pi^{\mathrm{an}}}T,
 \qquad
 \Phi_T:S_T\xrightarrow{\sim}T\times F_{\pi,\mathbf P}^{\mathrm{an}}.
\]
Let $f_T:S_T\to T$ be the projection, and let
$\omega^{\mathrm n}_{S_T/T}$ be the normalized complex
constructed by that lemma. Set
\[
 \mathcal E_{S_T/T}:=v^{-1}\mathcal E_{S_\pi/T_\pi},
 \qquad
 \mathcal B_{S_T/T}:=R(f_T)_*\omega^{\mathrm n}_{S_T/T}.
\]
The comparisons in (ii) become
\[
 \begin{aligned}
 \omega^{\mathrm n}_{S_T/T}
   &\simeq\Phi_T^{-1}
       \left(\omega_T\boxtimes\boxtimes_j\omega_{P'_j}\right),\\
 e_T=v^{-1}e_\pi:
 \mathcal E_{S_T/T}
   &\xrightarrow{\sim}\underline{B_{\pi,\mathbf P}^{\vee}}_T,\\
 b_T=\mathbb D_T(e_T):
 \omega_T\otimes B_{\pi,\mathbf P}
   &\xrightarrow{\sim}\mathcal B_{S_T/T}
      \simeq\mathbb D_T\mathcal E_{S_T/T}
      \simeq v^!\mathcal B_{S_\pi/T_\pi}.
 \end{aligned}
\]
The maps preserve (iii) under further pullbacks satisfying the
same hypotheses. Here $\omega_T$ is the actual dualizing
complex, also when $T$ has boundary or corners.
\end{enumerate}
The coordinate maps require $b_j-b_{j'}$ and $u_j$ to be
invertible. No product identification across a further collision,
or coordinate identification at $u_j=0$, is asserted.
\end{lemma}
\begin{proof}
\emph{(i) Source coordinates.}
For $(b,u)\in T_\pi(R)$ and $m_j=|I_j|$, the ideals
$(z-b_j)$ are pairwise comaximal. The Chinese remainder maps
\[
 R[z]/(D_a^N)\xrightarrow{\sim}
       \prod_jR[z]/((z-b_j)^{m_jN})
\]
give the completed-source decomposition of
Lemma~\ref{lem:kv-pole-compatibility}. Taking inverse limits
and substituting $z=b_j+u_jt_j$ gives the first map in
\eqref{eq:source-coordinate-crt}. On its $j$th factor,
\[
 \begin{aligned}
 D_a(b_j+u_jt_j)
  &=(u_jt_j)^{m_j}
       \prod_{\ell\ne j}(b_j-b_\ell+u_jt_j)^{m_\ell}
        \in t_j^{m_j}R\llbracket t_j \rrbracket^\times,\\
 H_{a,\mathbf P}(b_j+u_jt_j)
  &=(u_jt_j)^{r_0P'_j}
       \prod_{\ell\ne j}(b_j-b_\ell+u_jt_j)^{r_0P'_\ell}
        \in t_j^{r_0P'_j}R\llbracket t_j \rrbracket^\times.
 \end{aligned}
\]
The first identity gives the punctured-source isomorphism;
the second identifies the marked pole conditions with the
bounds defining $F_{P'_j}\mathscr LX$.
Thus applying $\operatorname{Hom}_{\C\text{-alg}}(\mathcal O(X),-)$
to the punctured-source map gives
\eqref{eq:source-coordinate-identification}, by representability.
Put $m_{\max}=\max_jm_j$. The image of $D_a^N\mathcal O_a$ is
$\prod_jt_j^{m_jN}R\llbracket t_j \rrbracket$, and
\[
 \prod_jt_j^{m_{\max}N}R\llbracket t_j \rrbracket
 \subseteq\prod_jt_j^{m_jN}R\llbracket t_j \rrbracket
 \subseteq\prod_jt_j^NR\llbracket t_j \rrbracket.
\]
The neighborhood bases are therefore cofinal, proving continuity
in both directions. For $R\to R'$, the maps are induced on
finite quotients and then coefficientwise on inverse limits;
no commutation of ordinary tensor product with completion is used.
For the analytic identification in (ii), use instead the
analytic realization of Section~\ref{sec:formal-loops}: its
topology is initial for regular functions with values in $\C$
with its usual topology. The scheme isomorphism and its inverse
therefore induce continuous maps on complex points. This
analytic assertion is separate from the adic continuity above.

\emph{(ii) Transport of the normalized comparisons.}
Let $\Phi_{\pi,1}$ denote the pullback to $T_\pi$ of the
identification with all multipliers equal to $1$.
It is the diagonal identification
$\kappa_\pi$, followed by separation and the singleton
identifications of Lemma~\ref{lem:relative-fusion-factorization}.
For general $u$, one has
\[
 \Phi_{\pi,u}=r_u\circ\Phi_{\pi,1},\qquad
 (r_u\eta)_j(g)(t)=\eta_j(g)(u_jt).
\]
On each separated factor,
\[
 r_{u_j}=(A_{0,u_j})^{-1}=A_{0,u_j^{-1}},
\]
since $A_{0,u_j}$ sends $\eta(g)(t)$ to $\eta(g)(t/u_j)$.
Use its inverse and composition comparisons specified in
Definition~\ref{def:coherent-relative-system}(ii).
Apply the corresponding comparisons of
Lemma~\ref{lem:relative-fusion-factorization}(i)--(ii) to this
composite. In integer grading, on a product finite factor
$Y=T_\pi\times\prod_jY_j$, the normalization used there is
\[
 \delta_{Y/T_\pi}=\sum_j\operatorname{pr}_j^*\delta_{Y_j},
 \qquad
 \omega_Y[-2\delta_{Y/T_\pi}]
   \simeq\omega_{T_\pi}\boxtimes
                       \boxtimes_j\omega_{Y_j}[-2\delta_{Y_j}].
\]
Hence its transported sheaf comparison is the first line of
\eqref{eq:source-coordinate-comparisons}, including the base
factor $\omega_{T_\pi}$. The compact-support product comparison
of the same lemma gives
\[
 \mathcal E_{S_\pi/T_\pi}
     \xrightarrow{\ e_{\pi,u}\ }\underline{\bigotimes_j B_{P'_j}^\vee}_{T_\pi}
     \simeq\underline{B_{\pi,\mathbf P}^\vee}_{T_\pi}.
\]
The last equivalence uses dualizability of the pole-stage
complexes, supplied by Corollary~\ref{cor:relative-pole-recovery}.
Dualizing gives
$b_{\pi,u}:\omega_{T_\pi}\otimes B_{\pi,\mathbf P}
\xrightarrow{\sim}\mathcal B_{S_\pi/T_\pi}$.
These are applications of the already descended coherent
comparisons, so no new finite-factor presentations or descent choices are made.

\emph{(iii) Changes of coordinates.}
On $\widehat T_\pi$, for every $g\in\mathcal O(X)$,
\[
 (\Phi_{\pi,u\lambda}\gamma)_j(g)(t)
   =\gamma(g)(b_j+u_j\lambda_jt)
   =(r_\lambda\Phi_{\pi,u}\gamma)_j(g)(t).
\]
Transport of the specified trace diagrams through this identity
gives the equality for $e$ in \eqref{eq:source-coordinate-change}.
Verdier duality reverses composition and gives the equality
for $b$; its base-change comparison is
\[
 \mathbb D_{\widehat T_\pi}(p^*e_\pi)
       \simeq p^!\mathbb D_{T_\pi}(e_\pi),
\]
and similarly for $m$. This also specifies the sources of the
two $b$-maps as $\omega_{\widehat T_\pi}\otimes B_{\pi,\mathbf P}$.
Indeed, $p$ and $m$ are smooth of complex relative dimension
$|J|$, so $p^!\simeq p^*[2|J|]$ and $m^!\simeq m^*[2|J|]$;
transitivity gives
$p^!\omega_{T_\pi}\simeq\omega_{\widehat T_\pi}
 \simeq m^!\omega_{T_\pi}$.
The map $(r_\lambda)_!$ is an endomorphism of
$\underline{B_{\pi,\mathbf P}^{\vee}}_{\widehat T_\pi}$
in the derived sheaf category, and $r_\lambda^!$ is its
Verdier dual. Both are obtained from trace diagrams over
the whole parameter space; their parameter dependence is
part of the morphisms, rather than determined only by their maps
on the
cohomology of individual fibers.
The calculation $b+u(c+vt)=b+uc+uvt$ proves
\eqref{eq:source-coordinate-composition}.
For the affine action, the defining formula for $\widetilde A$
gives
\[
 \gamma'(g)(\beta+vb_j+vu_jt)=\gamma(g)(b_j+u_jt),
\]
which is precisely the stated equivariance of $\Phi_\pi$.

The full diagonal identity for $\kappa$ and the compatibility
of separation with regrouping are those of
Lemma~\ref{lem:kv-pole-compatibility}; their sheaf comparisons,
including pole maps and mixed composition homotopies, are
constructed by Lemma~\ref{lem:relative-fusion-factorization}.
Transport these diagrams through the coordinate maps in (i).
Functoriality preserves their homotopies and higher
compatibilities. Each collision stratum has its own map
$\Phi$, so this argument does not restrict a separated
product map across a vanishing difference $b_j-b_{j'}$.

\emph{(iv) Topological pullback.}
Apply Lemma~\ref{lem:topological-pullback-presentations}(ii)
to the comparisons in (ii)--(iii). It gives
$\mathcal E_{S_T/T}\simeq v^{-1}\mathcal E_{S_\pi/T_\pi}$
and
$\mathcal B_{S_T/T}\simeq v^!\mathcal B_{S_\pi/T_\pi}$,
compatibly with duality and composition. Consequently
\[
 e_T=v^{-1}e_\pi,\qquad
 b_T=\mathbb D_T(e_T)\simeq v^!b_\pi,
 \qquad
 v^!(\omega_{T_\pi}\otimes B_{\pi,\mathbf P})
      \simeq\omega_T\otimes B_{\pi,\mathbf P}.
\]
For the total-space sheaf formula, let $Y$ be a finite-type
complex factor and write
\[
 p_Y:T\times Y^{\mathrm{an}}\to T,\qquad
 q_Y:T\times Y^{\mathrm{an}}\to Y^{\mathrm{an}}
\]
for the projections. There is a natural equivalence
\[
 q_Y^{-1}\omega_Y\otimes p_Y^{-1}K
       \xrightarrow{\sim}p_Y^!K
       \qquad(K\in\mathcal D(T;\F)).
\]
To justify it, work locally on $Y^{\mathrm{an}}$ in a
Euclidean neighborhood retract, as in
Subsection~\ref{subsec:finite-presentations}.
Applied to the constant sheaf on this neighborhood,
\cite[Propositions~7.8 and~6.5(2),(4)]{ScholzeSixFunctors}
identify $p_Y^!\underline\F_T$ with $q_Y^{-1}\omega_Y$
after base change to $T$, and identify
$p_Y^!\underline\F_T\otimes p_Y^{-1}K$ with $p_Y^!K$.
These canonical comparisons agree on overlaps. Taking
$K=\omega_T$ and using transitivity of exceptional inverse image
gives
\[
 \omega_{T\times Y^{\mathrm{an}}}
       \simeq\omega_T\boxtimes\omega_Y.
\]
In integer grading, apply this to the product finite factors
$T\times\prod_jY_j^{\mathrm{an}}$ with index
$\sum_j\operatorname{pr}_j^*\delta_{Y_j}$, and shift by
minus twice this index. Ordinary inverse image along the
coefficient tails and the chart descent of
Proposition~\ref{lem:relative-compact-support} then give
the asserted formula for $\omega^{\mathrm n}_{S_T/T}$.
Naturality of these comparisons preserves the composition
and refinement diagrams. This uses the actual $\omega_T$;
no manifold assumption on $T$ is needed. The same comparisons
in the periodic coefficient formalism give the parity assertions.
\end{proof}

\section{Fusion product: construction}
\label{sec:fusion-construction}

We first apply the normalized specialization and restriction constructions of
Section~\ref{sec:topological-preparations} to the collision families.
The identifications of Lemma~\ref{lem:source-coordinate-factorization}
express their input and output complexes in terms of the pole stages.
For two inputs, specialization and dual restriction give the same map
$B_P\otimes B_Q\to B_{P+Q}$. Passing to the pole homotopy colimit
gives the fusion product on $B$ and on $LH^*(X)$.
The last subsection gives its formula when the pole stages admit
global finite factors.

We follow the grading convention of
Subsection~\ref{subsec:relative-pole-recovery}: in integer grading,
$B_P$, $B$, $\mathbf B$, and $F_P\mathbf B$ stand for
$B_{\Omega,P}$, $B_\Omega$, $\mathbf B_\Omega$, and
$F_P\mathbf B_\Omega$, respectively, with $\Omega$ the chosen
nowhere-zero algebraic volume form. Without $\Omega$, these
symbols denote the parity-graded objects.

\subsection{Normalized specialization and restriction}
\label{subsec:normalized-specialization}
Fix the fusion-ready collision families of
Proposition~\ref{lem:all-multipoint-models} and their relative
complexes from Definition~\ref{def:relative-fusion-complex}.
In integer grading their relative volume indices satisfy
Subsection~\ref{conv:loop-volume-compatibility}; in parity no
volume form is required. All arrows and their comparisons are
in the derived enhancements fixed in
Subsection~\ref{conv:section-five-sheaf-operations}.
Write $[k]=\{1,\ldots,k\}$. For $k\ge2$, $\mathbf P=(P_i)_{i=1}^k$,
and pairwise distinct $c_i\in\C$, put
\[
 h:\A^1_s\longrightarrow\A^k,
 \quad h(s)=(b+c_i s)_{i=1}^k,\qquad
 S=h^*\mathcal Z_{[k],\mathbf P},\qquad b\in\C.
\]
Here $h^*\mathcal Z$ denotes scheme-theoretic pullback. Set
\[
 S_0=S\times_{\A^1}\{0\},\quad
 S^\times=S\times_{\A^1}\Gm,\quad
 H=[0,\infty),\quad L=(0,\infty),\quad
 S_L=S^{\mathrm{an}}\times_{\C}L.
\]
Let $i:S_0^{\mathrm{an}}\hookrightarrow S^{\mathrm{an}}$,
$j:(S^\times)^{\mathrm{an}}\hookrightarrow S^{\mathrm{an}}$, and
$\ell:S_L\hookrightarrow(S^\times)^{\mathrm{an}}$ be the inclusions.
The ray $L$ is oriented toward zero.

\begin{definition}[Specialization and restriction arrows]
\label{def:geometric-specialization-arrows}
For the above pullback family, use the arrow of
Lemma~\ref{lem:specialization-arrow-descent}
\[
 \operatorname{Sp}^{\mathrm n}_S:
 \mathcal K^{\mathrm n}_{S,\rho}
 :=Rj_*\ell_*\omega^{\mathrm n}_{S_L/L}[-1]
 \longrightarrow i_*\omega^{\mathrm n}_{S_0/\{0\}}.
\]
Its derived-global-section map is
\begin{equation}\label{eq:geometric-specialization-sections}
 R\Gamma(S_L,\omega^{\mathrm n}_{S_L/L})[-1]
 \xrightarrow{R\Gamma(\operatorname{Sp}^{\mathrm n}_S)}
 R\Gamma(S_0^{\mathrm{an}},\omega^{\mathrm n}_{S_0/\{0\}}).
\end{equation}
For the compact-support companion put
\[
 E_{S,H}:=\mathcal E_{S/\A^1}|_H
       \simeq(h^*\mathcal E_{[k],\mathbf P})|_H.
\]
The corresponding restriction arrow is
\begin{equation}\label{eq:geometric-restriction-sections}
 \operatorname{res}_{S;H,L}:
 R\Gamma(H,E_{S,H})\longrightarrow R\Gamma(L,E_{S,H}|_L).
\end{equation}
The identifications in these formulas are those of
Lemma~\ref{lem:topological-pullback-presentations}, applied to
$H,L\hookrightarrow\C$. 
\end{definition}

\begin{lemma}[Comparison and naturality of the arrows]
\label{lem:source-specialization-comparison}
The complex $E_{S,H}$ is constructible for $H=L\sqcup\{0\}$,
with perfect stalks and the trivialization on $L$ defined by
Lemma~\ref{lem:source-coordinate-factorization}, with input
coordinates $t_i=z-b-c_i s$. Under the
oriented-ray duality comparisons
\[
 \begin{aligned}
 R\Gamma(S_L,\omega^{\mathrm n}_{S_L/L})[-1]
    &\simeq R\Gamma(L,E_{S,H}|_L)^\vee,\\
 R\Gamma(S_0^{\mathrm{an}},\omega^{\mathrm n}_{S_0/\{0\}})
    &\simeq R\Gamma(H,E_{S,H})^\vee,
 \end{aligned}
\]
\eqref{eq:geometric-specialization-sections} is the linear dual of
\eqref{eq:geometric-restriction-sections}. These comparisons and
arrows have the following compatibilities.
\begin{enumerate}[(i),wide=0pt]
\item An isomorphism of presented families preserving the relative
complexes and the oriented ray identifies the arrows. In a fixed
trivialization of the punctured compact-support companion, two
inward rays bounding an oriented sector give homotopic maps.
\item For a factorization into a colliding group and separated
spectators, the specialization arrow is the group's arrow
exterior-tensored with the identity of the spectator complexes.
The restriction arrow has the dual compatibility.
\item Proper traces for pole inclusions and normalized smooth
comparisons are morphisms of specialization arrows. Common
finite-factor refinements and cofinal coefficient cutoffs induce
equivalences of the realized arrow diagrams. If $v:V\hookrightarrow S$ is an open
neighborhood of $S_0$, then, with the same closed-fiber target,
\[
 R\Gamma(\operatorname{Sp}^{\mathrm n}_V)\circ\operatorname{res}_v
       \simeq R\Gamma(\operatorname{Sp}^{\mathrm n}_S).
\]
Here $\operatorname{res}_v$ is restriction on the source complex;
it need not be an equivalence.
\end{enumerate}
All statements are comparisons of coherent diagrams, including
the pole maps, rather than only identities on cohomology.
\end{lemma}
\begin{proof}
The finite-factor presentations supplied by
Proposition~\ref{lem:all-multipoint-models} satisfy the hypotheses
of Lemma~\ref{lem:specialization-arrow-descent}. That lemma gives
(iii) and the isomorphism assertion in (i), and
Remark~\ref{rem:restriction-arrow-descent} gives dual restriction.
Constructibility follows from
Lemma~\ref{lem:relative-fusion-factorization}(ii), and the input
trivialization is that of Lemma~\ref{lem:source-coordinate-factorization};
perfectness is
Corollary~\ref{cor:relative-pole-recovery}(ii). Evaluation at zero
is an equivalence by Lemma~\ref{lem:parameter-chain-properties}(ii),
which gives the displayed closed-target duality comparison.
For the ray comparison in (i), put
$\mathcal E:=\mathcal E_{S/\A^1}$ and write $V$ for the perfect
complex in the fixed trivialization
$\mathcal E|_{\C^*}\simeq\underline V_{\C^*}$.
Let $\mathcal E_0$ be the stalk at zero, and put
$C:=V^\vee$, $D:=\mathcal E_0^\vee$, and $I:=[0,1]$.
Choose a piecewise semialgebraic path
$q:I\to\{u\in\C:|u|=1\}$ through the directions of the
specified sector, with endpoints the directions of the two rays.
Define
\[
 a:[0,1]\times I\longrightarrow\C,
 \qquad a(r,t)=r q(t),\qquad E':=a^{-1}\mathcal E.
\]
Let $i_0:I\hookrightarrow[0,1]\times I$, $i_0(t)=(0,t)$,
and $j_+:(0,1]\times I\hookrightarrow[0,1]\times I$ be the
zero-end and positive-part inclusions. Their ordinary inverse
images satisfy
\[
 i_0^*E'\simeq\underline{D^\vee}_I,\qquad
 j_+^*E'\simeq\underline{C^\vee}_{(0,1]\times I}.
\]
Lemma~\ref{lem:parameter-chain-properties}(i) therefore gives
\[
 \gamma:\omega_I\otimes C\longrightarrow\omega_I\otimes D.
\]
Its extraordinary restrictions to the endpoints of $I$ are the
two ray-specialization maps, under the fixed trivialization:
the adjunction map defining the collar construction restricts
to the corresponding map on each ray, and
Remark~\ref{rem:restriction-arrow-descent} identifies its dual
with specialization. Truncating a ray at radius $1$ gives the
same map, since evaluation at zero and in the constant positive
stratum identifies both restriction constructions with the same
stalk-to-positive-stalk map.
By Lemma~\ref{lem:parameter-chain-operations}(i), the oriented
singular interval gives a degree-$-1$ homotopy $K:C\to D$;
for endpoint representatives $f_0,f_1:C\to D$ it satisfies
\[
 d_DK+Kd_C=f_1-f_0.
\]
Here $d_C,d_D$ are the differentials of cochain representatives
chosen for the whole parameter
and pole diagram. Naturality of the collar construction and of
the parameter-chain comparison makes this homotopy compatible
with the pole maps and their comparison homotopies.

For (ii), the factorization
isomorphism identifies the closed--open localization diagram
with the exterior product of the group's diagram and the spectator
complexes. Localization and its dual restriction preserve this
exterior product. These are the specified transformations of the
finite-factor diagrams, so their descent preserves all composition
and pole-map homotopies.
\end{proof}

\subsection{Fusion at pole stages}
\label{subsBiSp}
\label{subsec:ordered-collision-operations}
In statements about sheaves or topology below, $F_P\sL X$ denotes
the space of complex points of the pole stage, with the analytic topology
induced from $\sL X(\C)$ as in
Convention~\ref{conv:loop-analytic-topology} and
Definition~\ref{def:pole-stages}. With the notation of
Subsection~\ref{subsec:relative-pole-recovery},
\[
 B_P=R\Gamma(F_P\sL X,\omega_P),\qquad
 B_{\mathbf P}:=\bigotimes_{i=1}^kB_{P_i},\qquad
 N:=\sum_iP_i,
\]
by Definition~\ref{def:pole-stage-complexes} and
Corollary~\ref{cor:relative-pole-recovery}(i).
The notation $B_{\mathbf P}$ abbreviates $B_{\pi,\mathbf P}$
of Lemma~\ref{lem:source-coordinate-factorization} for the
identity surjection $\pi:[k]\to[k]$.
These are the existing pole-stage objects, not new cochain choices.
Each $B_P$ is perfect by
Proposition~\ref{lem:relative-compact-support}; hence linear duality
identifies $B_{\mathbf P}^\vee$ with
$\bigotimes_iB_{P_i}^\vee$.

\begin{definition}[Normalized specialization at pole stages]
\label{def:normalized-ray-specialization}
For $S=h^*\mathcal Z_{[k],\mathbf P}$, denote by
\[
 \vartheta_L:L\times\prod_i F_{P_i}\sL X\xrightarrow{\sim}S_L,
 \qquad
 \vartheta_0:F_N\sL X\xrightarrow{\sim}S_0^{\mathrm{an}}
\]
the inverses of the source-coordinate identifications of
Lemma~\ref{lem:source-coordinate-factorization}: the input
coordinates are $t_i=z-b-c_i s$, and the output coordinate
is $t=z-b$. Put
\[
 \omega_{\rm in}:=\omega_L\boxtimes\boxtimes_i\omega_{P_i}.
\]
The sheaf comparisons in that lemma identify
$\vartheta_L^{-1}\omega^{\mathrm n}_{S_L/L}$ with $\omega_{\rm in}$
and $\vartheta_0^{-1}\omega^{\mathrm n}_{S_0/\{0\}}$ with
$\omega_N$, as coherent pole diagrams.
Let $\rho_L:\F\xrightarrow{\sim}R\Gamma(L,\omega_L)[-1]$
send $1$ to the inward fundamental class. Define
\[
 \begin{aligned}
 \alpha_{\mathbf P}:B_{\mathbf P}
 &\xrightarrow{\rho_L\otimes\mathrm{id}}
 R\Gamma(L,\omega_L)[-1]\otimes\bigotimes_iR\Gamma(F_{P_i}\sL X,\omega_{P_i})\\
 &\xrightarrow{\boxtimes}
 R\Gamma(L\times\textstyle\prod_i F_{P_i}\sL X,\omega_{\rm in})[-1]
 \xrightarrow{\sim}
 R\Gamma(S_L,\omega^{\mathrm n}_{S_L/L})[-1],\\
 \beta_{\mathbf P}:B_N
 &\xrightarrow{\sim}
 R\Gamma(S_0^{\mathrm{an}},\omega^{\mathrm n}_{S_0/\{0\}}).
 \end{aligned}
\]
The last arrows transport sections through $\vartheta_L,\vartheta_0$
and the sheaf comparisons of Lemma~\ref{lem:source-coordinate-factorization}. The exterior product is the
natural map on derived global sections, with the ray factor first.
Set
\begin{equation}\label{eq:derived-normalized-ray-specialization}
 \operatorname{sp}_S
 :=\beta_{\mathbf P}^{-1}\circ
 R\Gamma(\operatorname{Sp}^{\mathrm n}_S)\circ\alpha_{\mathbf P}
 :B_{\mathbf P}\longrightarrow B_N.
\end{equation}
For an open neighborhood $V$ of $S_0$, use
$R\Gamma(\operatorname{Sp}^{\mathrm n}_V)\operatorname{res}_v$
in this formula; the preceding lemma identifies the resulting map
with $\operatorname{sp}_S$. Only $k=2,3$ will be needed below.
\end{definition}

\begin{definition}[Fusion product at pole stages]
\label{lp:def:binary-fusion-product}
\label{def:restriction-fusion-map}
For $P,Q\ge0$, take
\[
 \sigma:\A^1_s\to\A^2,\quad\sigma(s)=(s,-s),\qquad
 S=\mathcal Z_{P,Q}:=\sigma^*\mathcal Z_{[2],(P,Q)}.
\]
Write $\alpha_{P,Q},\beta_{P,Q}$ for the preceding input map
and output equivalence, which are explicitly supplied for this
slice by Lemma~\ref{lem:relative-fusion-factorization}(iii).
The \emph{fusion product at pole stages} is
\begin{equation}\label{lp:eq:derived-binary-fusion}
 \overline m_{2;P,Q}:=\operatorname{sp}_{\mathcal Z_{P,Q}}
   =\beta_{P,Q}^{-1}R\Gamma(\operatorname{Sp}^{\mathrm n}_{\mathcal Z_{P,Q}})
                         \alpha_{P,Q}
 :B_P\otimes B_Q\longrightarrow B_{P+Q}.
\end{equation}
Put $E_{P,Q}:=\mathcal E_{P,Q}|_H$. Its
\emph{restriction description} is the linear dual $c_{P,Q}$ of
\begin{equation}\label{eq:binary-restriction-map}
 B_{P+Q}^\vee\xrightarrow{\sim}R\Gamma(H,E_{P,Q})
 \xrightarrow{\operatorname{res}_{H,L}}R\Gamma(L,E_{P,Q}|_L)
 \xrightarrow{\sim}(B_P\otimes B_Q)^\vee.
\end{equation}
The first equivalence is inverse to evaluation at $0$, by
Lemma~\ref{lem:parameter-chain-properties}(ii)--(iii), for $r=1$.
The last uses the input trivialization with $t_1=z-s$, $t_2=z+s$
from Lemma~\ref{lem:source-coordinate-factorization}, equivalently
\eqref{eq:binary-slice-factorization}, followed by sections on the
contractible ray $L$. Equivalences are inverted
as equivalences of the whole pole diagrams.
\end{definition}

\begin{proposition}[Pole-stage construction and comparison]
\label{prop:relative-binary-identification}
The normalized specializations
\eqref{eq:derived-normalized-ray-specialization} are degree-zero
morphisms, independent of compatible finite-factor presentations and
coherently natural for the pole maps. In particular there is a
homotopy of natural transformations in $\mathcal D(\F)$
(or $\mathcal D_{\mathrm{per}}(\F)$)
\[
 \overline m_{2;\bullet,\bullet}\simeq c_{\bullet,\bullet}:
 \bigl((P,Q)\longmapsto B_P\otimes B_Q\bigr)
 \longrightarrow\bigl((P,Q)\longmapsto B_{P+Q}\bigr).
\]
Both sides are therefore descriptions of the same morphism of
coherent pole diagrams. Source, ray, spectator and excision
comparisons are those of
Lemma~\ref{lem:source-specialization-comparison} under factorization.
\end{proposition}
\begin{proof}
Work in $\mathcal C=\mathcal D(\F)$ in volume grading and in
$\mathcal C=\mathcal D_{\mathrm{per}}(\F)$ in parity grading.
Write $\mathbf1$ for its tensor unit and take tensor products
and linear duals in $\mathcal C$.
Lemma~\ref{lem:specialization-arrow-descent} gives the degree-zero map
\[
 R\Gamma(\operatorname{Sp}^{\mathrm n}_S):
 R\Gamma(S_L,\omega^{\mathrm n}_{S_L/L})[-1]
 \longrightarrow
 R\Gamma(S_0^{\mathrm{an}},\omega^{\mathrm n}_{S_0/\{0\}}).
\]
The inward orientation gives
$\rho_L:\mathbf1\xrightarrow{\sim}R\Gamma(L,\omega_L)[-1]$:
the dualizing complex on the oriented real one-manifold $L$
is the constant coefficient sheaf shifted by $[1]$.
Thus $\alpha_{\mathbf P}$, $\beta_{\mathbf P}$, and
$\beta_{\mathbf P}^{-1}R\Gamma(\operatorname{Sp}^{\mathrm n}_S)
\alpha_{\mathbf P}$ all have degree zero.
For a smooth finite-factor refinement $v:Y'\to Y$ of relative
complex dimension $e$, put
$\delta_Y=\delta_{Y/\A^1}$ and $\delta_{Y'}=\delta_{Y'/\A^1}$.
In volume grading, $\delta_{Y'}=v^*\delta_Y+e$ and purity gives
\[
 \omega_{Y'}[-2\delta_{Y'}]
 \simeq v^{-1}(\omega_Y[-2\delta_Y])[2e-2e]
 \simeq v^{-1}(\omega_Y[-2\delta_Y]).
\]
In parity grading the even purity shift is absorbed by the
periodic comparison. These are the normalized purity equivalences
used in Lemma~\ref{lem:relative-fusion-factorization}.

For $R\le R'$, denote the pole transition by
$u_{R,R'}:B_R\to B_{R'}$.
If $P_i\le Q_i$ for all $i$, put
$\mathbf Q=(Q_1,\ldots,Q_k)$, $N'=\sum_iQ_i$, and
$S'=h^*\mathcal Z_{[k],\mathbf Q}$, using the same parameter map
$h$ as for $S$. Lemma~\ref{lem:source-specialization-comparison}(iii),
together with the input and output comparisons, gives
\[
 u_{N,N'}\circ\operatorname{sp}_S
 \simeq
 \operatorname{sp}_{S'}\circ
       \bigotimes_{i=1}^k u_{P_i,Q_i}.
\]
These homotopies are compatible with successive pole
transitions. The same lemma identifies the entire specialization
arrow diagrams under common finite-factor refinements, proving
independence of the compatible finite-factor presentations.

For each $i\in\{1,\ldots,k\}$, let $C_i^\bullet$ denote
the cosimplicial section diagram
$C_{P_i}^\bullet$ in volume grading, or
$C_{P_i,\mathrm{per}}^\bullet$ in parity grading, constructed
in Remark~\ref{rem:coherent-pole-descent}.
Its totalization gives
\[
 B_{P_i}\simeq\operatorname{Tot}(C_i^\bullet).
\]
For every tuple $(p_1,\ldots,p_k)$ of nonnegative integers,
tensoring the limit projections gives
\[
 \bigotimes_{i=1}^k\operatorname{Tot}(C_i^\bullet)
 \longrightarrow
 C_1^{p_1}\otimes\cdots\otimes C_k^{p_k}.
\]
These maps form a coherent cone over the multicosimplicial
diagram on the right. Its limit universal property therefore
gives the natural map
\[
 \bigotimes_{i=1}^k\operatorname{Tot}(C_i^\bullet)
 \longrightarrow
 \operatorname{Tot}_{p_1,\ldots,p_k}
 \bigl(C_1^{p_1}\otimes\cdots\otimes C_k^{p_k}\bigr).
\]
Put
\[
 G_{\mathbf P}:=
 R\Gamma\left(\prod_i F_{P_i}\sL X,\boxtimes_i\omega_{P_i}\right).
\]
Compose the preceding map with the map to $G_{\mathbf P}$ induced
by the termwise exterior products on the product \v Cech diagram.
This gives the exterior-product morphism
$\varepsilon_{\mathbf P}:B_{\mathbf P}\to G_{\mathbf P}$.
For the section-diagram comparisons, applying
Lemma~\ref{lem:totalization-realization-representatives}
one \v Cech direction at a time gives simultaneous cochain
representatives of the induced maps on totalizations, including
the coherent pole and refinement comparisons.
The input map is then the composite
\[
 \begin{aligned}
 B_{\mathbf P}
 &\xrightarrow{\rho_L\otimes\varepsilon_{\mathbf P}}
 R\Gamma(L,\omega_L)[-1]\otimes G_{\mathbf P}\\
 &\xrightarrow{\boxtimes}
 R\Gamma\left(L\times\prod_i F_{P_i}\sL X,\omega_{\rm in}\right)[-1]
 \xrightarrow{\sim}
 R\Gamma(S_L,\omega^{\mathrm n}_{S_L/L})[-1],
 \end{aligned}
\]
where $\omega_{\rm in}$ and the last comparison are those of
Definition~\ref{def:normalized-ray-specialization}.
Only a natural map from the tensor product of totalizations to
the multicosimplicial totalization was used; no commutation of
tensor products with infinite limits is asserted.

To identify $\varepsilon_{\mathbf P}$ with the dual
compact-support product, put
\[
 \begin{gathered}
 A_i:=\mathcal E_{F_{P_i}\mathscr LX/\mathrm{pt}}
       \simeq B_{P_i}^\vee,\qquad
 A:=\bigotimes_iA_i,\\
 A_{\mathrm{prod}}:=
       \mathcal E_{\prod_iF_{P_i}\mathscr LX/\mathrm{pt}}.
 \end{gathered}
\]
Proposition~\ref{lem:relative-compact-support}(iv) gives
\[
 \kappa:A\xrightarrow{\sim}A_{\mathrm{prod}},\qquad
 G_{\mathbf P}\simeq A_{\mathrm{prod}}^\vee,
 \qquad B_{P_i}\simeq A_i^\vee.
\]
Let $\operatorname{ev}_i:B_{P_i}\otimes A_i\to\mathbf1$ and
$\operatorname{ev}_{\mathrm{prod}}:
G_{\mathbf P}\otimes A_{\mathrm{prod}}\to\mathbf1$
be the resulting evaluation pairings. Define
$\chi_{\mathbf P}:B_{\mathbf P}\to A^\vee$ to be adjoint to
\[
 B_{\mathbf P}\otimes A
 \xrightarrow{\tau}
 \bigotimes_i(B_{P_i}\otimes A_i)
 \xrightarrow{\bigotimes_i\operatorname{ev}_i}\mathbf1,
\]
where $\tau$ moves each $A_i$ next to $B_{P_i}$ using the
Koszul symmetry, preserving the order of the resulting pairs.
The map $\chi_{\mathbf P}$ is an equivalence because the $A_i$
are perfect.

The identity to be checked is the following homotopy of pairings
with source $B_{\mathbf P}\otimes A$ and target $\mathbf1$:
\[
 \operatorname{ev}_{\mathrm{prod}}\circ
       (\varepsilon_{\mathbf P}\otimes\kappa)
 \simeq
 \left(\bigotimes_i\operatorname{ev}_i\right)\circ\tau.
\]
On each tuple of finite-factor charts this is the compatibility
of Verdier evaluation with exterior products. For a transition
$t:T'\to T$ between compact-support trace terms, the compatibility
with its dual restriction is the evaluation identity
\[
 \operatorname{ev}_T\circ(\mathrm{id}_{T^\vee}\otimes t)
 \simeq
 \operatorname{ev}_{T'}\circ(t^\vee\otimes\mathrm{id}_{T'}):
 T^\vee\otimes T'\longrightarrow\mathbf1.
\]
Here $\operatorname{ev}_T$ and $\operatorname{ev}_{T'}$ are the
linear evaluation pairings. Naturality of these pairings and of
the exterior products preserves the composition homotopies of the
normalized traces, as in
Proposition~\ref{lem:relative-compact-support}. Thus the
chartwise identities identify the coherent cones over the full
multicosimplicial section diagram. They give the displayed
global pairing identity by the limit universal property.
For a simplicial compact-support trace diagram $T_\bullet$, the
duality comparison used here is
\[
 \left(\operatorname*{hocolim}_{[q]\in\Delta^{\mathrm{op}}}
                  T_q\right)^\vee
 \simeq
 \operatorname*{holim}_{[q]\in\Delta}T_q^\vee.
\]
It does not interchange a tensor product with a limit.
Adjunction now gives, under $G_{\mathbf P}\simeq A_{\mathrm{prod}}^\vee$,
\[
 \kappa^\vee\circ\varepsilon_{\mathbf P}
       \simeq\chi_{\mathbf P},\qquad
 \varepsilon_{\mathbf P}
       \simeq(\kappa^\vee)^{-1}\circ\chi_{\mathbf P}.
\]
Inserting $\rho_L$ as the first factor identifies
$\alpha_{\mathbf P}$ with the oriented dual input comparison,
with the same Koszul order. All these comparisons preserve the
pole-transition homotopies.

For the binary slice $S=\mathcal Z_{P,Q}$, let
$r_{P,Q}:B_{P+Q}^\vee\to(B_P\otimes B_Q)^\vee$
be the composite in \eqref{eq:binary-restriction-map}.
By definition, $c_{P,Q}=r_{P,Q}^\vee$, using perfect biduality
to identify its source with $B_P\otimes B_Q$ and its target
with $B_{P+Q}$. Under
\eqref{eq:binary-slice-factorization} and
\eqref{eq:binary-slice-sheaf-factorization},
Remark~\ref{rem:restriction-arrow-descent} gives the square
\[
 \begin{tikzcd}[column sep=large,row sep=large]
 B_P\otimes B_Q \arrow[r,"c_{P,Q}"] \arrow[d,"\alpha_{P,Q}"']
 & B_{P+Q} \arrow[d,"\beta_{P,Q}","\sim"']\\
 R\Gamma(S_L,\omega^{\mathrm n}_{S_L/L})[-1]
 \arrow[r,"R\Gamma(\operatorname{Sp}^{\mathrm n}_S)"']
 & R\Gamma(S_0^{\mathrm{an}},\omega^{\mathrm n}_{S_0/\{0\}}).
 \end{tikzcd}
\]
The left comparison is $\alpha_{P,Q}$ by the pairing calculation
above, with the inward ray first; the right comparison is
$\beta_{P,Q}$ by the closed-fiber sheaf identification.
The square commutes in the derived enhancement, so
\[
 \overline m_{2;P,Q}
 =\beta_{P,Q}^{-1}R\Gamma(\operatorname{Sp}^{\mathrm n}_S)
                         \alpha_{P,Q}
 \simeq r_{P,Q}^\vee=c_{P,Q}.
\]
The cited remark descends the duality comparison with its entire
arrow diagram. Hence, as $(P,Q)$ varies with the coordinatewise
pole order, the displayed square is a commutative square in
$\operatorname{Fun}(\Z_{\ge0}^2,\mathcal C)$.
It therefore gives the asserted homotopy
$\overline m_{2;\bullet,\bullet}\simeq c_{\bullet,\bullet}$
of natural transformations. Only the factorization of the
descended open family is used; no product decomposition of an
individual finite total factor is required.
\end{proof}

\subsection{The product on the loop complex and loop cohomology}
\label{subsec:loop-fusion-colimit}
The pole-stage product now passes to the filtered object
$B=\hocolim_PB_P$ of Definition~\ref{def:loop-complex}.

\begin{proposition}[Passage to the loop complex]
\label{lp:prop:binary-fusion}
The coherent pole-stage morphism of
Proposition~\ref{prop:relative-binary-identification} induces a
morphism $\overline m_2:B\otimes B\to B$ with the sum filtration,
whose underlying morphism is
\begin{equation}\label{eq:fusion-pole-colimit}
 B\otimes B\simeq\hocolim_{P,Q}(B_P\otimes B_Q)
 \xrightarrow{\hocolim\overline m_{2;P,Q}}
 \hocolim_{P,Q}B_{P+Q}\simeq B.
\end{equation}
At filtration level $N$ it is
\[
 \hocolim_{P+Q\le N}(B_P\otimes B_Q)\longrightarrow B_N,
\]
using the transitions $B_{P+Q}\to B_N$.
On the chosen telescope model of
Theorem~\ref{thm:canonical-complex} it admits a degree-zero
filtered cochain representative $m_2$, unique up to filtered
chain homotopy as a representative of this morphism.
\end{proposition}
\begin{proof}
Apply Lemma~\ref{lem:filtered-operations-colimit} to the coherent
natural transformation of
Proposition~\ref{prop:relative-binary-identification}. Tensor
product preserves homotopy colimits, and the addition functor
$(P,Q)\mapsto P+Q$ is cofinal, giving
\eqref{eq:fusion-pole-colimit}. The same lemma gives its sum
filtration and its representative on the chosen telescope.
\end{proof}

\begin{definition}[Fusion product on the loop complex and its cohomology]
\label{def:loop-fusion-product}
The \emph{fusion product on the loop complex} is the filtered
morphism
\[
 \overline m_2:B\otimes B\longrightarrow B
\]
constructed in Proposition~\ref{lp:prop:binary-fusion}, where
the source has the sum filtration.

On the fixed filtered cochain model $(\mathbf B,d_{\mathbf B},F_\bullet)$ of
Theorem~\ref{thm:canonical-complex}, choose a degree-zero
filtered cochain representative
\[
 m_2:\mathbf B\otimes\mathbf B\longrightarrow\mathbf B
\]
of $\overline m_2$, whose existence is supplied by
Proposition~\ref{lp:prop:binary-fusion}. For $P,Q\ge0$, write
\[
 m_{2;P,Q}:=
 m_2|_{F_P\mathbf B\otimes F_Q\mathbf B}:
 F_P\mathbf B\otimes F_Q\mathbf B\longrightarrow F_{P+Q}\mathbf B.
\]
Under the specified equivalence of entire filtered diagrams
$(F_P\mathbf B)_{P\ge0}\simeq(B_P)_{P\ge0}$ from
Theorem~\ref{thm:canonical-complex}, these restrictions represent
the coherent pole-stage morphism
$\overline m_{2;\bullet,\bullet}$. Thus the comparison includes
the transition morphisms and their higher coherence data, rather
than only the objects at each pole bound.

The \emph{fusion product on loop cohomology} is
\[
 \mu:LH^*(X)\otimes LH^*(X)\longrightarrow LH^*(X),
 \qquad \mu([a],[b]):=[m_2(a,b)],
\]
for cocycles $a,b\in\mathbf B$, using the specified identification
$H^*(\mathbf B)\cong H^*(B)=LH^*(X)$.
We write $ab:=\mu(a,b)$ for cohomology classes $a,b$ when
the fusion algebra is understood; the notation $\mu(a,b)$
is retained when the operation must be displayed explicitly.
Use volume grading and the subscript $\Omega$ when $\Omega$
is specified, and parity grading otherwise.
\end{definition}

Since $m_2$ is a degree-zero cochain map, for homogeneous
$a,b\in\mathbf B$ one has
\begin{equation}\label{lp:eq:binary-leibniz}
 d_{\mathbf B}m_2(a,b)
   =m_2(d_{\mathbf B}a,b)+(-1)^{|a|}m_2(a,d_{\mathbf B}b).
\end{equation}
Consequently the displayed formula for $\mu$ is well defined.

\begin{corollary}[Passage to cohomology]
\label{prop:loop-fusion-cohomology}
The induced map $\mu$ is independent of the compatible
finite-factor presentations and cochain representatives. It satisfies
\[
 \mu(LH^i,LH^j)\subset LH^{i+j},\qquad
 \mu(F_PLH^*,F_QLH^*)\subset F_{P+Q}LH^*.
\]
For $\iota_P:H^*(B_P)\to LH^*(X)$,
\begin{equation}\label{eq:fusion-cohomology-stages}
 \mu(\iota_P\alpha,\iota_Q\beta)
   =\iota_{P+Q}H^*(\overline m_{2;P,Q})(\alpha\otimes\beta).
\end{equation}
Consequently these stage maps determine $\mu$ under
$LH^*(X)\cong\varinjlim_PH^*(B_P)$; no injectivity
of $\iota_P$ is required.
\end{corollary}
\begin{proof}
Apply cohomology to the filtered morphism
\eqref{eq:fusion-pole-colimit} and use
Remark~\ref{rem:pole-cohomology-colimit}.
\end{proof}

\subsection{Finite-factor formulas and global finite models}
\label{subsec:global-binary-fusion}
The following formulas compute the product already constructed.
Local finite factors give maps of descent diagrams; compatible
global finite models permit computation without chart-overlap
terms.

\begin{lemma}[Finite-factor computation of fusion]
\label{lem:finite-factor-fusion-formula}
Fix $P,Q\ge0$, put $R_1=P$, $R_2=Q$, $R_0=P+Q$, and let
$S=\mathcal Z_{P,Q}\to\A^1_s$ be the binary collision family of
Definition~\ref{lp:def:binary-fusion-product}.
Choose compatible finite-factor presentations of its finite diagrams,
including the input and output identifications, pole maps, and
finite composable strings of chart maps, together with their
finite-factor refinements, as in
Definition~\ref{def:fusion-ready-data}(ii).
All section complexes below are taken on the associated analytic spaces.
\begin{enumerate}[(i),wide=0pt]
\item \emph{The finite-factor operation.}
Let $v:\mathcal U\to S$ be a nice chart in these finite-factor presentations, with finite-factor
projection $\pi:\mathcal U\to T$ over $\A^1_s$, and put
\[
 T^\times=T\times_{\A^1}\Gm,\qquad
 T_0=T\times_{\A^1}\{0\},\qquad
 T_L=T^{\mathrm{an}}\times_{\C}L,
 \quad L=(0,\infty).
\]
Assume that, after common refinements of finite factors of the
same nice charts, $T$ has specified identifications of its
\emph{whole} punctured and closed fibers
\begin{equation}\label{lp:eq:binary-finite-family}
 \begin{aligned}
 T^\times\cong G_{\rm in}
   &=\Gm\times Y_1\times Y_2\times\A^{e_{\rm in}},\\
 T_0\cong G_{\rm out}&=Y_0\times\A^{e_{\rm out}}.
 \end{aligned}
\end{equation}
Here $Y_i$ is a finite factor of a pole chart at bound $R_i$,
and $e_{\rm in},e_{\rm out}\ge0$ count further coordinates
retained from the existing affine tails. They do not count the
base coordinate $s$ and do not change the family.
The first identification is over $\Gm$; their compatibility
with the input and output morphisms is the condition in
\textup{(ii)}.
In integer grading, write $\delta_T=\delta_{T/\A^1}$ and
$\delta_i=\delta_{Y_i}$ for the volume indices. Their restrictions
satisfy Subsection~\ref{conv:loop-volume-compatibility}:
\begin{equation}\label{lp:eq:fusion-delta-balance}
 \delta_{G_{\rm in}/\Gm}=\operatorname{pr}_1^*\delta_1+
                    \operatorname{pr}_2^*\delta_2+e_{\rm in},
 \qquad
 \delta_{G_{\rm out}/\{0\}}=\operatorname{pr}_0^*\delta_0+e_{\rm out},
\end{equation}
where the left sides are the restrictions of $\delta_T$ and
$\operatorname{pr}_i$ denotes projection to $Y_i$.
Put
\[
 p_{\rm in}:Y_1\times Y_2\times\A^{e_{\rm in}}
                    \longrightarrow Y_1\times Y_2,
 \qquad p_{\rm out}:Y_0\times\A^{e_{\rm out}}\longrightarrow Y_0.
\]
Let $\operatorname{sp}_T$ be the degree-preserving Borel--Moore
specialization of Lemma~\ref{lp:lem:bm-specialization-formula},
using \eqref{lp:eq:binary-finite-family}. For $a,b\in\Z$, define
\[
 \begin{aligned}
 \mu_{\mathrm{fin}}:
 H_a^{\mathrm{BM}}(Y_1;\F)\otimes H_b^{\mathrm{BM}}(Y_2;\F)
 &\longrightarrow
 H_{a+b+2e_{\rm in}-2e_{\rm out}}^{\mathrm{BM}}(Y_0;\F),\\
 \mu_{\mathrm{fin}}
 &:=(p_{\rm out}^{!})^{-1}\circ\operatorname{sp}_T
                \circ p_{\rm in}^{!}\circ\boxtimes.
 \end{aligned}
\]
Here $\boxtimes$ is exterior product, the Gysin maps are the
complex-oriented affine Gysin isomorphisms, and
$\operatorname{sp}_T$ maps the homology of
$Y_1\times Y_2\times\A^{e_{\rm in}}$ to that of
$T_0\cong Y_0\times\A^{e_{\rm out}}$.
Extend $\mu_{\mathrm{fin}}$ to the graded groups by linearity.
If $\rho\in H_1^{\mathrm{BM}}(\Gm;\F)$ is the image of the
inward fundamental class under $L\hookrightarrow\Gm$, then
\begin{equation}\label{lp:eq:full-fusion-formula}
 \mu_{\mathrm{fin}}(\alpha,\beta)
 =(p_{\rm out}^{!})^{-1}\partial_{s=0}
     \bigl(\rho\boxtimes p_{\rm in}^{!}
                         (\alpha\boxtimes\beta)\bigr).
\end{equation}
The ray is placed first. Here $\partial_{s=0}$ is the localization
boundary for $T_0\hookrightarrow T\hookleftarrow T^\times$,
expressed using \eqref{lp:eq:binary-finite-family}.

For the derived operation, put
\[
 \begin{gathered}
 C_i:=R\Gamma(Y_i,\omega_{Y_i}[-2\delta_i]),\qquad i=0,1,2,\\
 A_{T,\mathrm{in}}:=
 R\Gamma(T_L,\omega_{T_L}[-2\delta_T|_{T_L}])[-1],\qquad
 A_{T,\mathrm{out}}:=
 R\Gamma(T_0,\omega_{T_0}[-2\delta_T|_{T_0}]).
 \end{gathered}
\]
In parity, use $\operatorname{Per}_2\omega_{Y_i}$,
$\operatorname{Per}_2\omega_{T_L}$, and
$\operatorname{Per}_2\omega_{T_0}$ in these definitions,
retaining the ray shift $[-1]$; no volume indices are chosen.
Using the endpoint identifications, define the following maps on
derived section complexes:
\[
 a_T:C_1\otimes C_2\longrightarrow A_{T,\mathrm{in}},\qquad
 b_T:C_0\xrightarrow{\sim}A_{T,\mathrm{out}}.
\]
The map $a_T$ is exterior product, followed by the normalized
input affine Gysin map and exterior product with the inward ray
class as the first factor; $b_T$ is the normalized output affine
Gysin map. The affine Gysin maps are induced by smooth purity
with the stated normalization shifts; the ray factor uses
$\F\xrightarrow{\sim}R\Gamma(L,\omega_L)[-1]$ from
Definition~\ref{def:normalized-ray-specialization}.
Writing
\[
 \Sigma_T:=R\Gamma(\operatorname{Sp}^{\mathrm n}_T):
 A_{T,\mathrm{in}}\longrightarrow A_{T,\mathrm{out}},\qquad
 m_T:=b_T^{-1}\Sigma_Ta_T:C_1\otimes C_2\longrightarrow C_0,
\]
for the arrow of Lemma~\ref{lem:specialization-arrow-descent}
and its indicated composite, $m_T$ has degree zero and induces
$\mu_{\mathrm{fin}}$ under the normalized Borel--Moore identifications.
All volume shifts and these identifications are read componentwise
when the indices vary.

\item \emph{Compatibility with pole-stage fusion.}
Let $u_i:U_i\to F_{R_i}\sL X$ be the pole charts in these
finite-factor presentations, with projections $\pi_i:U_i\to Y_i$.
On complex points the chart maps induce the restrictions
\[
 r_i:B_{R_i}=R\Gamma(F_{R_i}\sL X,\omega_{R_i})
 \longrightarrow R\Gamma(U_i,u_i^{-1}\omega_{R_i})
 \xrightarrow{\sim}C_i.
\]
On cohomology, $r_i$ sends a pole-stage class to its
finite-factor representative in the sense of
Definition~\ref{def:finite-factor-representatives}.
Likewise, restriction along $v$ and the finite-factor comparisons
for $\pi$ define
\[
 \begin{aligned}
 q_{\mathrm{in}}:
 R\Gamma(S_L,\omega^{\mathrm n}_{S_L/L})[-1]
     &\longrightarrow A_{T,\mathrm{in}},\\
 q_{\mathrm{out}}:
 R\Gamma(S_0,\omega^{\mathrm n}_{S_0/\{0\}})
     &\longrightarrow A_{T,\mathrm{out}},
 \end{aligned}
\]
where $S_L=S^{\mathrm{an}}\times_{\C}L$ and
$S_0=S\times_{\A^1}\{0\}$.
Neither these chart restrictions nor the $r_i$ are assumed to
be equivalences.

The requirement that \eqref{lp:eq:binary-finite-family} realize
the specified input and output comparisons of
Lemma~\ref{lem:relative-fusion-factorization}(iii) means that
the displayed arrows satisfy
\[
 q_{\mathrm{in}}\alpha_{P,Q}
       \simeq a_T(r_1\otimes r_2),\qquad
 q_{\mathrm{out}}\beta_{P,Q}\simeq b_Tr_0,
\]
with their overlap, refinement, and pole compatibilities.
Here $\alpha_{P,Q},\beta_{P,Q}$ are the input and output maps
of Definition~\ref{lp:def:binary-fusion-product}.
Then there is a homotopy of derived morphisms
\[
 r_0\overline m_{2;P,Q}
       \simeq m_T(r_1\otimes r_2):B_P\otimes B_Q\longrightarrow C_0.
\]
In particular, for $a\in H^*(B_P)$ and $b\in H^*(B_Q)$,
\[
 H^*(r_0)\bigl(H^*(\overline m_{2;P,Q})(a\otimes b)\bigr)
   =\mu_{\mathrm{fin}}\bigl(H^*(r_1)a,H^*(r_2)b\bigr),
\]
using the normalized Borel--Moore identifications in \textup{(i)}.
These homotopies satisfy the stated diagram compatibilities.
Global computation uses the underlying specialization arrows
on the full coherent descent diagrams of
Proposition~\ref{prop:relative-binary-identification}, followed
by Corollary~\ref{prop:loop-fusion-cohomology}; \textup{(i)}
gives their product-factor formula wherever its hypotheses hold.
This uses full \v Cech cocycles, including overlap components;
the identities on individual chart cohomology groups alone
need not determine the global class.

The whole-fiber identifications and their compatibility with
the endpoint morphisms are hypotheses of this formula; local
finite-factor presentations alone do not supply them. A restriction
to an input chart cannot be reversed to extend an arbitrary
class to $T^\times$: its map on cohomology need not be surjective.

\item \emph{The formula in loop cohomology for global finite-factor presentations.}
If the three $r_i$ are equivalences, as for compatible global
finite-factor presentations of the pole stages, let
\[
 \iota_R:H^*(B_R)\longrightarrow LH^*(X),\qquad
 j_i:=\iota_{R_i}\circ H^*(r_i)^{-1}:H^*(C_i)\longrightarrow LH^*(X)
\]
be the canonical pole-stage map and the indicated composites.
For $\alpha\in H^*(C_1)$ and $\beta\in H^*(C_2)$, one has
\[
 \mu(j_1\alpha,j_2\beta)
       =j_0\bigl(\mu_{\mathrm{fin}}(\alpha,\beta)\bigr).
\]
Again the right side uses normalized Borel--Moore homology.
No injectivity of $\iota_R$ is required. Compatible global
collision finite-factor presentations to which this applies are specified in
Proposition~\ref{prop:global-binary-fusion}.

\item \emph{Several inputs and excision.}
For $k\ge2$ and pole bounds $\mathbf P=(P_1,\ldots,P_k)$, put
$N=\sum_iP_i$ and take the ray family
$S=h^*\mathcal Z_{[k],\mathbf P}$ of
Definition~\ref{def:normalized-ray-specialization}.
Assume the analogous specified whole-fiber identifications
\[
 T^\times\cong\Gm\times\prod_iY_i\times\A^{e_{\rm in}},\qquad
 T_0\cong Y_{\rm out}\times\A^{e_{\rm out}},
\]
with $Y_i$ at bound $P_i$, $Y_{\rm out}$ at bound $N$, and
the endpoint morphism compatibilities of \textup{(ii)}.
Define their normalized complexes $C_i,C_{\rm out}$ and chart
restrictions $r_i:B_{P_i}\to C_i$, $r_{\rm out}:B_N\to C_{\rm out}$
as in \textup{(i)--(ii)}.
In integer grading the input and output relative indices are
$\sum_i\operatorname{pr}_i^*\delta_{Y_i}+e_{\rm in}$ and
$\operatorname{pr}_{\rm out}^*\delta_{Y_{\rm out}}+e_{\rm out}$,
where the $\delta$'s are the corresponding pole-factor indices
and the $\operatorname{pr}$'s are the projections to those factors.
Let $m_T^{(k)}:\bigotimes_iC_i\to C_{\rm out}$ be the composite
$b_T^{-1}\Sigma_Ta_T$ of \textup{(i)}, using the $k$-fold exterior
product to define $a_T$, and denote its map on normalized
Borel--Moore homology by $\operatorname{sp}_{\rm fin}$.
Then
\[
 r_{\rm out}\operatorname{sp}_S
       \simeq m_T^{(k)}\bigl(\bigotimes_i r_i\bigr),
\]
where $\operatorname{sp}_S:\bigotimes_iB_{P_i}\to B_N$ is
the normalized specialization of that definition.

For its Borel--Moore formula put
\[
 p_{\rm in}:\prod_iY_i\times\A^{e_{\rm in}}\to\prod_iY_i,
 \qquad p_{\rm out}:Y_{\rm out}\times\A^{e_{\rm out}}\to Y_{\rm out}.
\]
Choose an open neighborhood $V\subseteq T$ containing $T_0$,
put $V^\times=V\cap T^\times$, and let
$j_{\rm in}:V^\times\hookrightarrow T^\times$ be the inclusion.
For $\alpha_i\in H_*^{\mathrm{BM}}(Y_i;\F)$ one has
\begin{equation}\label{eq:normalized-ray-specialization}
 \operatorname{sp}_{\rm fin}(\alpha_1,\ldots,\alpha_k)
 =(p_{\rm out}^!)^{-1}\partial_{s=0}j_{\rm in}^*
 \bigl[\rho\boxtimes p_{\rm in}^!
                   (\alpha_1\boxtimes\cdots\boxtimes\alpha_k)\bigr],
\end{equation}
where $\partial_{s=0}$ is now the boundary for
$T_0\hookrightarrow V\hookleftarrow V^\times$.
For $V=T$ the restriction is omitted; in particular, the binary
case gives \eqref{lp:eq:full-fusion-formula}.
\end{enumerate}
\end{lemma}
\begin{proof}
\noindent\textup{(i)} The boundary on the finite factor is
\[
 \partial_T:H_{d+1}^{\mathrm{BM}}(T^\times;\F)
       \longrightarrow H_d^{\mathrm{BM}}(T_0;\F),\qquad d\in\Z.
\]
Lemma~\ref{lp:lem:bm-specialization-formula} gives
$\operatorname{sp}_T(\xi)=\partial_T(\rho\boxtimes\xi)$
under the prescribed punctured-fiber identification.
Composing with the two affine Gysin isomorphisms proves
\eqref{lp:eq:full-fusion-formula}.
The derived arrow $\Sigma_T$ uses the same inward-ray trace and
localization boundary. The maps $a_T,b_T$ use the same exterior
product and normalized affine Gysin maps just defined. Thus taking
cohomology of $b_T^{-1}\Sigma_Ta_T$ gives exactly $\mu_{\mathrm{fin}}$.

If $\alpha,\beta$ have Borel--Moore degrees $a,b$, input Gysin
adds $2e_{\rm in}$, the ray adds $1$, the boundary subtracts
$1$, and inverse output Gysin subtracts $2e_{\rm out}$.
On each connected component of $T$ meeting both fibers, the
same locally constant relative index gives
$\delta_0+e_{\rm out}=\delta_1+\delta_2+e_{\rm in}$. Hence
\[
 2\delta_0-(a+b+2e_{\rm in}-2e_{\rm out})
       =(2\delta_1-a)+(2\delta_2-b).
\]
This is also the cancellation of shifts in the defining
derived composite, so $m_T$ has degree zero.
Components with empty punctured fiber contribute zero.
The inward ray is first and has boundary $+[0]$; all affine
Gysin shifts are even. These conventions give the sign in the
formula and the same degree-zero assertion in parity.

\medskip
\noindent\textup{(ii)} The section equivalences underlying
Definition~\ref{def:finite-factor-representatives} and its
relative version in
Definition~\ref{def:relative-normalized-dual-complexes}
are natural under the finite-factor refinements in these diagrams by
Lemma~\ref{lem:infinite-product-descent}(i).
They define the comparisons in $r_i$ and the $q$'s;
no \'etale chart restriction is inverted.

Naturality of the specialization arrow in
Lemma~\ref{lem:specialization-arrow-descent}, including its
descent to the chart $\mathcal U$, gives
\[
 q_{\rm out}R\Gamma(\operatorname{Sp}^{\mathrm n}_S)
       \simeq\Sigma_Tq_{\rm in}.
\]
Using the two stipulated endpoint comparisons and the definition
of $\overline m_{2;P,Q}$, one obtains
\[
 \begin{aligned}
 b_Tr_0\overline m_{2;P,Q}
 &\simeq q_{\rm out}\beta_{P,Q}\overline m_{2;P,Q}\\
 &\simeq q_{\rm out}R\Gamma(\operatorname{Sp}^{\mathrm n}_S)\alpha_{P,Q}\\
 &\simeq\Sigma_Tq_{\rm in}\alpha_{P,Q}\\
 &\simeq\Sigma_Ta_T(r_1\otimes r_2)
   \simeq b_Tm_T(r_1\otimes r_2).
 \end{aligned}
\]
Since $b_T$ is an equivalence, this proves the asserted homotopy;
taking cohomology and using \textup{(i)} proves the class formula.
The specialization naturality and the specified endpoint comparisons
include their composition homotopies and pole compatibilities on
all overlaps and refinements. The displayed calculation
therefore takes place in their coherent diagrams.
Applying homotopy totalization to the underlying specialization
arrow diagrams gives the pole-stage maps by
Proposition~\ref{prop:relative-binary-identification}; passing
to the pole homotopy colimit gives fusion by
Proposition~\ref{lp:prop:binary-fusion}.
This applies to the full descent complexes, not to a replacement
by the cohomology groups of their individual terms.

\medskip
\noindent\textup{(iii)} Put $a=H^*(r_1)^{-1}\alpha$ and
$b=H^*(r_2)^{-1}\beta$. By \textup{(ii)},
\[
 H^*(\overline m_{2;P,Q})(a\otimes b)
       =H^*(r_0)^{-1}\mu_{\rm fin}(\alpha,\beta).
\]
Applying $\iota_{P+Q}$ and using
Corollary~\ref{prop:loop-fusion-cohomology} gives
\[
 \mu(j_1\alpha,j_2\beta)
 =\iota_{P+Q}H^*(\overline m_{2;P,Q})(a\otimes b)
 =j_0\mu_{\rm fin}(\alpha,\beta).
\]
Only the $r_i$ were inverted, so no injectivity of a pole-stage
map is needed.

\medskip
\noindent\textup{(iv)} Replace the binary exterior product and the maps
$\alpha_{P,Q},\beta_{P,Q}$ in \textup{(i)--(ii)} by the
$k$-fold exterior product and
$\alpha_{\mathbf P},\beta_{\mathbf P}$ of
Definition~\ref{def:normalized-ray-specialization}.
The same composite and naturality calculation give
$r_{\rm out}\operatorname{sp}_S\simeq
m_T^{(k)}(\bigotimes_i r_i)$ and the Borel--Moore formula for
$V=T$; the degree calculation uses $\sum_i a_i$ and
$\sum_i\delta_{Y_i}$ in place of $a+b$ and $\delta_1+\delta_2$,
where $a_i$ is the Borel--Moore degree of $\alpha_i$.
For general $V$, let $\partial_V$ be the localization boundary
for $T_0\hookrightarrow V\hookleftarrow V^\times$.
Naturality of localization gives
\[
 \partial_Vj_{\rm in}^*=\partial_T,
\]
because the map on the common closed fiber $T_0$ is the identity.
This proves \eqref{eq:normalized-ray-specialization} and uses
restriction only in its stated direction.
\end{proof}

The next criterion is the multiplicative counterpart of
Subsection~\ref{subsec:global-finite-presentations}. It is useful in computations.
See also Remark~\ref{rem:computing-fusion} for a further computational method for the fusion product.

\begin{proposition}
\label{prop:global-binary-fusion}
Use the 1-parameter binary collision families
$\mathcal Z_{P,Q}\to\A^1_s$ with the fusion-ready data of
Proposition~\ref{lem:all-multipoint-models}. Suppose their
pole-stage finite-factor presentations are global as in
Proposition~\ref{prop:global-finite-presentations}.
When $\Omega$ is given, put
\[
 B_{P,N}:=R\Gamma\bigl(Z_{P,N},
                   \omega_{Z_{P,N}}[-2\delta_{P,N}]\bigr).
\]
As for $B_P$, the dependence of $B_{P,N}$ on $\Omega$ is
suppressed in integer grading.
Assume the chosen collision finite-factor presentation diagrams, including
their pole maps and fiber identifications, admit simultaneous
cofinal towers of global finite factors
$\mathcal T_{P,Q,N}\to\A^1_s$, with affine-bundle transition
maps over $\A^1_s$ and specified isomorphisms
\begin{equation}\label{eq:global-binary-fibers}
 \begin{aligned}
 \mathcal T_{P,Q,N}^{\times}
  &\cong\Gm\times Z_{P,N_1}\times Z_{Q,N_2}
                                      \times\A^{e_{\rm in}},\\
 (\mathcal T_{P,Q,N})_0
  &\cong Z_{P+Q,N_0}\times\A^{e_{\rm out}}.
 \end{aligned}
\end{equation}
The global factors and these isomorphisms realize the chosen
families and their input and output product identifications, simultaneously on each
finite collection of pole bounds after common cofinal refinement.
Here $N_i=N_i(P,Q,N)$ are cofinal cutoffs in the respective
pole-stage towers; they need not be equal. Apply the notation of
Lemma~\ref{lem:finite-factor-fusion-formula} with
\[
 T=\mathcal T_{P,Q,N},\qquad
 (Y_1,Y_2,Y_0)=(Z_{P,N_1},Z_{Q,N_2},Z_{P+Q,N_0}).
\]
Realizing the chosen input and output product identifications includes the endpoint
comparison identities of that lemma's part \textup{(ii)},
coherently with the pole and refinement maps in these finite-factor presentation
diagrams.
The pole and refinement compatibilities are those of
Proposition~\ref{prop:direct-relative-volume-indices}. In integer
grading its relative normalization is required, so
\eqref{lp:eq:fusion-delta-balance} holds under these substitutions.
Then the following statements hold.
\begin{enumerate}[(i),wide=0pt]
\item The operations of
Lemma~\ref{lem:finite-factor-fusion-formula}(i) give degree-zero
morphisms
\[
 m_{2;P,Q,N}:=m_{\mathcal T_{P,Q,N}}:
 B_{P,N_1}\otimes B_{Q,N_2}\longrightarrow B_{P+Q,N_0}.
\]
On their cohomology, expressed as normalized Borel--Moore
homology, the formula is
\begin{equation}\label{eq:global-binary-product}
 \mu_{P,Q,N}(\alpha,\beta)
  =(p_{{\rm out},N}^{!})^{-1}
       \operatorname{sp}_{\mathcal T_{P,Q,N}}
       \bigl(p_{{\rm in},N}^{!}(\alpha\boxtimes\beta)\bigr),
\end{equation}
where $p_{{\rm in},N}$ and $p_{{\rm out},N}$ are the projections
$p_{\rm in}$ and $p_{\rm out}$ of
Lemma~\ref{lem:finite-factor-fusion-formula}(i) under these
substitutions, with its specialization and orientation conventions.
\item Under the normalized transition maps of
Proposition~\ref{prop:global-finite-presentations}, these finite
morphisms represent the coherent pole-stage fusion maps.
Lemma~\ref{lem:filtered-operations-colimit} identifies their
induced filtered derived morphism with
$\overline m_2:B\otimes B\to B$ of
Definition~\ref{def:loop-fusion-product}; on the fixed model it is
represented by $m_2:\mathbf B\otimes\mathbf B\to\mathbf B$. Its cohomology product
$\mu$ is computed by their direct-limit maps under
\eqref{eq:bm-limits}. For a fixed pair of pole-stage classes,
one sufficiently large compatible finite-factor presentation computes
the product; no chart-overlap descent is needed.
\end{enumerate}
Without $\Omega$, take $B_{P,N}=R\Gamma(Z_{P,N},
\operatorname{Per}_2\omega_{Z_{P,N}})$ and omit the integer
volume-index requirements. The same conclusions hold in parity.
\end{proposition}
\begin{proof}
By Proposition~\ref{prop:global-finite-presentations}, write
\[
 r_{R,n}:B_R\xrightarrow{\sim}B_{R,n}
\]
for the coherent comparisons at pole bound $R$ and chosen
cutoff $n$. Lemma~\ref{lem:finite-factor-fusion-formula}(i)--(ii),
applied with the specified factors, gives the degree-zero maps,
\eqref{eq:global-binary-product}, and
\[
 r_{P+Q,N_0}\overline m_{2;P,Q}
 \simeq m_{2;P,Q,N}(r_{P,N_1}\otimes r_{Q,N_2}).
\]
These homotopies are coherent on the simultaneous global
finite-factor presentation diagrams by that lemma, applied with the endpoint
comparison data in its part \textup{(ii)}. This proves \textup{(i)} and the pole-stage
comparison in \textup{(ii)}.

Lemma~\ref{lem:filtered-operations-colimit} therefore identifies
the induced filtered derived morphism with $\overline m_2$
in \eqref{eq:fusion-pole-colimit}.
Lemma~\ref{lem:finite-factor-fusion-formula}(iii) identifies
its cohomology product under \eqref{eq:bm-limits}.
Cofinality places any fixed pair of pole-stage classes in one
sufficiently large compatible finite diagram; normalized affine
refinements are quasi-isomorphisms, so further cutoffs preserve
the result. This proves \textup{(ii)}. The parity case is included
in the cited parts of the lemma.
\end{proof}

\subsection{Example -- The multiplicative group}
\label{ex:cstar-global-fusion}
\label{ex:cstar-fusion-ready-presentations}
Continue Example~\ref{ex:cstar-loop-complex}, with $X=\Gm$
and $\Omega=dx/x$. Put $C:=R\Gamma(\Gm(\C);\F)$. We identify
$
 H^1(C)=H^1(\Gm(\C);\bQ)=H^1(\Gm;\bQ),
$
where the last notation denotes the cohomology of the analytic
space $\Gm(\C)$. Let $\theta\in H^1(\Gm;\bQ)$ evaluate to
$1$ on the counterclockwise unit circle. In the winding-$r$
summand, that example identifies
$1,\theta$ with $e_r,\eta_r$, respectively, and gives
\[
 LH^*_{\Omega}(\Gm)=
       \bigoplus_{r\in\Z}(\F e_r\oplus\F\eta_r).
\]
For $\Omega'=c x^k\,dx/x$, with $c\in\C^*$ and $k\in\Z$,
use the same symbols $e_r,\eta_r$ for the corresponding
cohomology classes of the shifted winding complex $C[-2kr]$
of that example.

\begin{proposition}[The fusion algebra of the multiplicative group]
\label{prop:cstar-fusion-algebra}
For every $\Omega'=c x^k\,dx/x$, the fusion product on
$LH^*_{\Omega'}(\Gm)$ satisfies
\begin{equation}\label{eq:cstar-global-fusion-product}
 e_ve_w=e_{v+w},\qquad
 e_v\eta_w=\eta_ve_w=\eta_{v+w},\qquad
 \eta_v\eta_w=0 \quad(v,w\in\Z).
\end{equation}
It is associative, unital with unit $e_0$, and graded
commutative. With $T=e_1$ and $\eta=\eta_0$, there is an
isomorphism of graded algebras
\[
 LH^*_{\Omega'}(\Gm)
   \cong\F[T,T^{-1}]\otimes\Lambda_{\F}(\eta),
 \qquad |T|=2k,\quad |\eta|=1,
 \qquad T^v=e_v,\quad T^v\eta=\eta_v,
\]
where $\Lambda_{\F}(\eta)=\F\oplus\F\eta$ and $\eta^2=0$.
Thus the windingwise identification for a change from $\Omega$
to $\Omega'$ preserves multiplication and changes the degrees to
\[
 |e_v|=2kv,\qquad |\eta_v|=2kv+1.
\]
The same multiplication and algebra presentation hold in parity,
with $T$ even and $\eta$ odd.
\end{proposition}

We prove the proposition using three lemmas. In their statements
and proofs, use $\Omega=dx/x$, so $B_K=B_{\Omega,K}$.
For $K\in\Z_{\ge0}$ and $r\in\Z$ with $|r|\le K$, put
\[
 W_{K,r}:=D(x_r y_{-r})\subset F_K\mathscr L\Gm,
 \qquad y=x^{-1},\quad x_j=[t^j]x,\quad y_j=[t^j]y.
\]
Example~\ref{ex:cstar-loop-complex} gives
$F_K\mathscr L\Gm=\coprod_{r=-K}^{K}W_{K,r}$ and
$(W_{K,r})_{\mathrm{red}}\cong\Gm\times\A^{\Z_{>0}}$,
with leading-coefficient coordinate $a_0=[t^r]x$.
Its normalized comparisons give the summand inclusions and their
cohomology maps
\[
 B_K\simeq\bigoplus_{r=-K}^{K}C,\qquad
 j_{K,r}:C\longrightarrow B_K,\qquad
 \kappa_{K,r}:=H^*(j_{K,r}).
\]
Let $\iota_K:H^*(B_K)\to LH^*_{\Omega}(\Gm)$ be the
canonical pole-stage map. The pole transitions preserve the
existing winding summands, so
$\kappa_r:=\iota_K\kappa_{K,r}$ is independent of $K\ge|r|$;
it satisfies $\kappa_r(1)=e_r$ and $\kappa_r(\theta)=\eta_r$.

Fix $P,Q\ge0$ and let $\mathcal Z=\mathcal Z_{P,Q}\to\A^1_s$
be the collision family of
Definition~\ref{lp:def:binary-fusion-product}.
For a family $\mathcal Y\to\A^1_s$, write
$\mathcal Y^\times=\mathcal Y\times_{\A^1_s}\mathbb G_{m,s}$
and $\mathcal Y_0=\mathcal Y\times_{\A^1_s}\{0\}$,
where $\mathbb G_{m,s}=\A^1_s\setminus\{0\}$.
In the source coordinates $t_1=z-s$, $t_2=z+s$ and $t=z$
of that definition, the Chinese remainder decomposition and
Example~\ref{ex:cstar-loop-complex} give
$$
 \mathcal Z^\times\cong
   \coprod_{\substack{|v|\le P\\|w|\le Q}}\mathcal W_{v,w},\quad\quad
 \mathcal W_{v,w}:=\mathbb G_{m,s}\times W_{P,v}\times W_{Q,w},$$
$$ \mathcal Z_0\cong\coprod_{r=-P-Q}^{P+Q}W_{P+Q,r}.
 $$
Regard $\mathcal W_{v,w}$ as the corresponding clopen subscheme
of $\mathcal Z^\times$.

Fix integers $v,w$ with $|v|\le P$ and $|w|\le Q$.
For a complex algebra $R$ and $s\in R$, put
\[
 D_s=z^2-s^2,\qquad
 \mathcal O_s=\varprojlim_{n\ge1}R[z]/(D_s^n),\qquad
 \mathcal K_s=\mathcal O_s[D_s^{-1}],\qquad
 H_s=(z-s)^P(z+s)^Q.
\]
As in Example~\ref{ex:kv-binary-family}, $\mathcal Z$ parametrizes
$x,y\in\mathcal K_s$ with $xy=1$ and
$H_sx,H_sy\in\mathcal O_s$. Define
\[
 q=(z-s)^{P+v}(z+s)^{Q+w},\qquad
 p=(z-s)^{P-v}(z+s)^{Q-w}.
\]

\begin{lemma}[The closed family for a winding pair]
\label{lem:cstar-closed-winding-family}
The conditions $H_sx\in q\mathcal O_s$ and
$H_sy\in p\mathcal O_s$ define a finitely presented closed
inclusion $i_{v,w}:\mathcal E_{v,w}\hookrightarrow\mathcal Z$.
Here the notation suppresses the fixed bounds $P,Q$.
Over every $R$, its points have the unique expression
\begin{equation}\label{eq:cstar-collision-factorization}
 H_sx=qU,\qquad H_sy=pU^{-1},\qquad U\in\mathcal O_s^\times,
\end{equation}
equivalently $x=(z-s)^v(z+s)^wU$ and $y=x^{-1}$.
This expression commutes with coefficient-ring homomorphisms.
The endpoint maps identify the following closed subschemes:
\[
 \mathcal E_{v,w}^\times=(\mathcal W_{v,w})_{\mathrm{red}}
        \subset\mathcal Z^\times,\qquad
 (\mathcal E_{v,w})_0=(W_{P+Q,v+w})_{\mathrm{red}}
        \subset\mathcal Z_0.
\]
\end{lemma}

\begin{proof}
Put $M=\max\{P+|v|,Q+|w|\}$. Both $q$ and $p$ divide
$D_s^M$. For $h=q,p$, monic division gives
$\mathcal O_s/h\mathcal O_s\cong R[z]/(h)$: for $n\ge M$,
the remainder of a multiple of $D_s^n$ is zero and its quotient
by $h$ is divisible by $D_s^{n-M}$, so division extends to the
completion. The remainders of $H_sx$ modulo $q$ and of $H_sy$
modulo $p$ depend only on their classes modulo $D_s^M$.
Their vanishing gives finitely many coefficient equations,
proving the closed finite-presentation assertion.

Monic division also gives
$\mathcal O_s\cong R\llbracket u \rrbracket[z]/(z^2-s^2-u)$, with $u$ mapped
to $D_s$. This module is free over $R\llbracket u \rrbracket$ with basis $1,z$,
so $D_s$ is a nonzerodivisor. Since $q,p,H_s$ divide powers
of $D_s$, they too are nonzerodivisors. The unique quotients
$H_sx=qU$, $H_sy=pV$ satisfy
$H_s^2=qpUV=H_s^2UV$, hence $UV=1$.
Conversely, every unit $U$ gives the stated point, since $z-s$
and $z+s$ are invertible in $\mathcal K_s$. Uniqueness and
monic division give compatibility with coefficient-ring maps.
This includes $P=Q=0$, when $v=w=0$, $q=p=1$, and $U=x$.

For invertible $s$, the Chinese remainder theorem identifies
$\mathcal O_s$ with $R\llbracket t_1 \rrbracket\times R\llbracket t_2 \rrbracket$. On the two
branches the formula for $x$ becomes
\[
 x(s+t_1)=t_1^v(2s+t_1)^wU(s+t_1),\qquad
 x(-s+t_2)=t_2^w(-2s+t_2)^vU(-s+t_2).
\]
The additional factors are units. At $s=0$, the formula is
$x=z^{v+w}U$ with $U\in R\llbracket z \rrbracket^\times$.
Together with the finite-factor presentations of the reduced winding loci recalled above,
these formulas give inverse coefficient morphisms for both
endpoint identifications.
\end{proof}

For $N\ge1$, let
\[
 \mathcal V_N:=
 \{(s,a_0,b_0):(a_0+sb_0)(a_0-sb_0)\ne0\}\times\A^{2N-2},
\]
where the remaining coordinates are $a_j,b_j$ for $1\le j<N$.

\begin{lemma}[Global finite factors and their normalization]
\label{lem:cstar-winding-global-factors}
Retaining $U$ modulo $D_s^N$ gives a finite-factor presentation
\[
 \mathcal E_{v,w}\cong
 \mathcal V_N\times
 \A^{\{(j,\epsilon):j\ge N,\ \epsilon\in\{0,1\}\}}
 \quad\text{over }\A^1_s.
\]
The transition $\mathcal V_{N+1}\to\mathcal V_N$ is an
affine-plane projection, and
\begin{equation}\label{eq:cstar-global-collision-fibers}
 \mathcal V_N^\times\cong
 \mathbb G_{m,s}\times(\Gm\times\A^{N-1})^2,\qquad
 (\mathcal V_N)_0\cong\Gm\times\A^{2N-1}.
\end{equation}
The relative volume index inherited from the collision family
is $\delta_{\mathcal V_N/\A^1_s}=2N$.
\end{lemma}

\begin{proof}
Every element of $\mathcal O_s$ has a unique expansion
$\sum_{j\ge0}(a_j+zb_j)D_s^j$. Modulo $D_s$, multiplication
by $a_0+zb_0$ has determinant $a_0^2-s^2b_0^2$ in the basis
$1,z$. A series is a unit precisely when this determinant is
invertible; its inverse lifts by a convergent geometric series
in the $D_s$-adic completion. All remaining coefficients are
unrestricted. This proves the finite-factor presentation and the
transition assertion. For invertible $s$, reduction modulo
$D_s^N$ gives two unit jets of length $N$ at $z=s,-s$;
at zero it gives one unit jet of length $2N$ in $z$.
This proves the fiber identifications.

Include a finite-factor presentation of $i_{v,w}$ in the
chosen finite-factor presentation diagrams of $\mathcal Z$
and compare the restricted factors with $\mathcal V_N$ by
common smooth refinements. On $s\ne0$, the inherited index
is $N+N=2N$, by Example~\ref{ex:cstar-loop-complex}.
The transported index is locally constant and $\mathcal V_N$
is connected, so it is $2N$ everywhere. This agrees at zero
with the output index for a unit jet of length $2N$.
Each affine-plane refinement adds two to the index.
\end{proof}

Put $\mathcal V=\mathcal V_1$, with coordinates $(s,a,b)$,
and let $p_0:\mathcal V_0=\Gm\times\A^1\to\Gm$ be
$p_0(a,b)=a$. The actual input leading coefficients in
the parameters $t_1=z-s$, $t_2=z+s$ are
\[
 A_+=(2s)^w(a+sb),\qquad A_-=(-2s)^v(a-sb).
\]
They give an isomorphism
$\mathcal V^\times\to\mathbb G_{m,s}\times\Gm^2$,
$(s,a,b)\mapsto(s,A_+,A_-)$.
Let $L=(0,\infty)\subset\C^*$ have the inward orientation,
and put $\mathcal V_L=\mathcal V(\C)\times_{\C}L$.
Using this input isomorphism and placing the ray first,
Lemma~\ref{lp:lem:bm-specialization-formula} gives
\[
 \operatorname{sp}_{\mathcal V}^{v,w}:
 H_d^{\mathrm{BM}}(\Gm^2;\F)
       \longrightarrow H_d^{\mathrm{BM}}(\mathcal V_0;\F)
       \qquad(d\in\Z).
\]
For $\alpha\in H^i(C)$ and $\beta\in H^j(C)$, define
$m_{v,w}(\alpha,\beta)\in H^{i+j}(C)$ by
\[
 m_{v,w}(\alpha,\beta)\frown[\Gm]
 =(p_0^!)^{-1}\operatorname{sp}_{\mathcal V}^{v,w}
       \bigl((\alpha\frown[\Gm])\boxtimes
             (\beta\frown[\Gm])\bigr).
\]
Here brackets denote complex-oriented fundamental classes,
$\boxtimes$ is exterior product, and
$p_0^!:H_d^{\mathrm{BM}}(\Gm;\F)
\xrightarrow{\sim}H_{d+2}^{\mathrm{BM}}(\mathcal V_0;\F)$
is the affine-line Gysin isomorphism.
The right side lies in $H_{2-i-j}^{\mathrm{BM}}(\Gm;\F)$,
so smooth duality defines $m_{v,w}$ uniquely; extend it linearly.

\begin{lemma}[Comparison with pole-stage fusion]
\label{lem:cstar-winding-fusion-comparison}
With the preceding notation,
\[
 H^*(\overline m_{2;P,Q})
   \bigl(\kappa_{P,v}(\alpha)\otimes\kappa_{Q,w}(\beta)\bigr)
 =\kappa_{P+Q,v+w}\bigl(m_{v,w}(\alpha,\beta)\bigr).
\]
Here $\overline m_{2;P,Q}:B_P\otimes B_Q\to B_{P+Q}$ is
the morphism of Definition~\ref{lp:def:binary-fusion-product},
and K\"unneth identifies its input cohomology with
$H^*(B_P)\otimes H^*(B_Q)$. Consequently
\[
 \mu\bigl(\kappa_v(\alpha),\kappa_w(\beta)\bigr)
   =\kappa_{v+w}\bigl(m_{v,w}(\alpha,\beta)\bigr).
\]
\end{lemma}

\begin{proof}
Apply Lemma~\ref{lem:specialization-arrow-descent}(ii) to the
closed inclusion $i_{v,w}$ in these diagrams, with the normalization of
Lemma~\ref{lem:cstar-winding-global-factors}.
By Lemma~\ref{lem:cstar-closed-winding-family}, its endpoint
maps are reduction inclusions into the indicated clopen winding
loci. After the input and output comparisons of
Definition~\ref{lp:def:binary-fusion-product}, their normalized
proper traces are therefore the summand inclusions
\[
 j_{P,v}\otimes j_{Q,w}:C\otimes C\longrightarrow B_P\otimes B_Q,
 \qquad j_{P+Q,v+w}:C\longrightarrow B_{P+Q}.
\]
These comparisons use the actual leading coefficients
$A_+,A_-$ on the input and $a$ on the output.

The projection $\mathcal V_N\to\mathcal V$ has affine
relative dimension $2N-2$, equal to the difference of its
indices $2N$ and $2$. Its normalized affine comparison is an
equivalence commuting with specialization by the same descent
lemma. The calculation therefore reduces to $\mathcal V$.
The finite-family formula of
Lemma~\ref{lem:finite-factor-fusion-formula}(i) uses, in the
notation of that lemma,
\[
 T=\mathcal V,\qquad Y_1=Y_2=Y_0=\Gm,\qquad
 e_{\rm in}=0,\quad e_{\rm out}=1,
 \qquad p_{\rm in}=\operatorname{id},\quad p_{\rm out}=p_0.
\]
Its relative index is $2$ and each $Y_i$ has index $1$;
with the inward ray first, smooth duality identifies its
cohomology operation with $m_{v,w}$.
Definition~\ref{lp:def:binary-fusion-product} identifies the
ambient normalized specialization with $\overline m_{2;P,Q}$.
Proper-trace naturality gives the first identity. Applying
$\iota_{P+Q}$ and \eqref{eq:fusion-cohomology-stages} gives
the identity on $LH$.
\end{proof}

\begin{proof}[Proof of Proposition~\ref{prop:cstar-fusion-algebra}]
First take $\Omega=dx/x$. The invertible functions
$u_\pm=a\pm sb$ on $\mathcal V$ define
$c_\pm=u_\pm^*\theta\in H^1(\mathcal V(\C);\F)$.
On $L$, the factors $(2s)^w$ and $(-2s)^v$ have constant
arguments and contribute zero angular class. Under
$\mathcal V_L\cong L\times\Gm(\C)^2$ using $A_+,A_-$,
\[
 c_+|_{\mathcal V_L}=1\otimes\theta\otimes1,\qquad
 c_-|_{\mathcal V_L}=1\otimes1\otimes\theta,
 \qquad c_\pm|_{\mathcal V_0}=p_0^*\theta.
\]
The family $\mathcal V\to\A^1_s$ is smooth, so specialization
sends the complex-oriented fundamental class of $\Gm^2$ to
that of $\mathcal V_0$ with coefficient $+1$ for the inward
ray placed first. More explicitly, use the product interpretation
in Lemma~\ref{lem:specialization-module}(iii) with
$Y=\mathcal V$, $U=\Gm^2$, and the punctured trivialization
defined by $A_+,A_-$. Put
$\mathcal V_H:=\mathcal V(\C)\times_{\C}[0,\infty)$ and let
$i_0:\mathcal V_0(\C)\hookrightarrow\mathcal V(\C)$ be the
special-fiber inclusion. For a homogeneous class
$c\in H^d(\mathcal V(\C);\F)$, write $c_U$ for its restriction
to $\{s_0\}\times\Gm(\C)^2$ in this trivialization, where
$s_0>0$. Apply the cap-product identity of that lemma to
$c|_{\mathcal V_H}$ and $\xi=[\Gm^2]$. It gives
\[
 \begin{aligned}
 \operatorname{sp}_{\mathcal V}^{v,w}
       (c_U\frown[\Gm^2])
   &=(i_0^*c)\frown
       \operatorname{sp}_{\mathcal V}^{v,w}([\Gm^2])
     \\
   &=(i_0^*c)\frown[\mathcal V_0]
 \quad\text{in }H^{\mathrm{BM}}_{4-d}(\mathcal V_0;\F).
 \end{aligned}
\]
For $c=1,c_+,c_-,c_+\smile c_-$, respectively, the input
classes $c_U$ are $1,\theta\otimes1,1\otimes\theta,
\theta\otimes\theta$. Consider their special restrictions, that is,  the
ordinary cohomological pullbacks along $i_0$:
\[
 i_0^*1=1,\qquad i_0^*c_+=i_0^*c_-=p_0^*\theta,
 \qquad i_0^*(c_+\smile c_-)=p_0^*(\theta^2)=0.
\]
Here $\theta^2=0$ in $H^2(\Gm;\bQ)$.
Smooth duality and inverse output Gysin give
\[
 m_{v,w}(1,1)=1,\qquad
 m_{v,w}(\theta,1)=m_{v,w}(1,\theta)=\theta,
 \qquad m_{v,w}(\theta,\theta)=0.
\]
Lemma~\ref{lem:cstar-winding-fusion-comparison}, for arbitrary
$P\ge|v|$ and $Q\ge|w|$, proves
\eqref{eq:cstar-global-fusion-product} for $\Omega$.
These formulas and the additive basis from
Example~\ref{ex:cstar-loop-complex} give the stated algebra,
its unit, associativity, and graded commutativity.

For $\Omega'=c x^k\,dx/x$, Example~\ref{ex:cstar-loop-complex}
gives input indices $N+kv$ and $N+kw$ and output index
$2N+k(v+w)$ on the factors of $\mathcal E_{v,w}$.
Write $\delta^{\Omega'}_{\mathcal V_N/\A^1_s}$ for the
inherited relative volume index associated with $\Omega'$.
As in Lemma~\ref{lem:cstar-winding-global-factors}, connectedness
of $\mathcal V_N$ extends its separated-fiber index to the whole
family, yielding
\[
 \delta^{\Omega'}_{\mathcal V_N/\A^1_s}
    =2N+k(v+w)=(N+kv)+(N+kw).
\]
Thus the normalization $\omega_Y[-2\delta_Y]$ on a finite
factor $Y$ changes by compatible even shifts: the input shifts
$-2kv,-2kw$ add to the output shift $-2k(v+w)$.
The same specialization and endpoint traces give the same
multiplication, with no additional Koszul sign. The constant
$c$ has zero order and changes no index.
The shifted winding complex $C[-2kv]$ places $e_v,\eta_v$
in degrees $2kv,2kv+1$, proving the assertion for $\Omega'$.
Two-periodization gives the parity assertion.
\end{proof}

\begin{remark}[How the selected families meet]
\label{rem:cstar-selected-family-components}
The global factors also describe the reduced total collision
scheme, although this description is not needed to compute
the product of a fixed pair of winding classes.
As $(v,w)$ ranges over $|v|\le P$, $|w|\le Q$, one has
\[
 \mathcal E_{v,w}=\overline{\mathcal W_{v,w}}^{\,\mathrm{red}},
 \qquad
 \mathcal Z_{\mathrm{red}}=
       \bigcup_{\substack{|v|\le P\\|w|\le Q}}\mathcal E_{v,w}.
\]
The closure is taken in $\mathcal Z$ with reduced induced
structure, and the union is the scheme-theoretic union.
These families are precisely the irreducible components of
$\mathcal Z_{\mathrm{red}}$.
Indeed, the coordinate ring of each $\mathcal E_{v,w}$ is
\[
 \C[s,a_0,b_0,(a_0^2-s^2b_0^2)^{-1},a_1,b_1,a_2,b_2,\ldots],
\]
an integral domain in which $s\ne0$ is dense.
Every geometric point of $\mathcal Z$ lies in one such family:
away from zero use its winding pair, and at zero choose bounded
$v,w$ with $v+w$ equal to its output winding.
The finite union therefore has all of $|\mathcal Z|$; its ideal,
being an intersection of radical ideals, is the nilradical.
The punctured loci are nonempty and pairwise disjoint, so none
of these irreducible closed families contains another.

For $-P-Q\le r\le P+Q$, put
\[
 \mathcal Z^{(r)}:=
   \bigcup_{\substack{|v|\le P,\ |w|\le Q\\v+w=r}}\mathcal E_{v,w}.
\]
The endpoint descriptions show that these are disjoint closed
unions covering $\mathcal Z_{\mathrm{red}}$; their complements
are finite closed unions, so
$\mathcal Z_{\mathrm{red}}=\coprod_r\mathcal Z^{(r)}$ is clopen.
Within $\mathcal Z^{(r)}$, distinct families have disjoint
punctured loci and the same reduced special-fiber image
$(W_{P+Q,r})_{\mathrm{red}}$.
For $P=Q=1$, the families
$\mathcal E_{1,-1},\mathcal E_{0,0},\mathcal E_{-1,1}$
all have that image with $r=0$.
This does not assert that the scheme-theoretic fiber
$(\mathcal Z_{\mathrm{red}})_0$ is reduced.
\end{remark}

\begin{remark}[Pole bounds and the contact filtration]
\label{rem:cstar-fusion-contact-filtration}
Increasing $P$ multiplies $H_s,q,p$ by $z-s$, and increasing
$Q$ multiplies them by $z+s$. The pole embeddings therefore
identify the selected unit families and leave $U$ unchanged in
\eqref{eq:cstar-collision-factorization}.
By Example~\ref{ex:cstar-loop-complex}, the contact filtration,
for every $\Omega'=c x^k\,dx/x$ and also in parity, is
\[
 F_PLH^*_{\Omega'}(\Gm)
   =\bigoplus_{|v|\le P}\F\{e_v,\eta_v\}.
\]
The product preserves this filtration since $|v+w|\le|v|+|w|$.
The inequality can be strict: $e_1e_{-1}=e_0$ has contact
weight $0$, whereas both input classes have contact weight $1$.
\end{remark}

\section{Fusion product: properties}
\label{sec:incidence}
We have constructed the fusion product in Section~\ref{sec:fusion-construction}.
In this section we prove its ring-theoretic properties and compute
its associated-graded version. We will also compare the effect of different
compactifications and record additional valuation gradings. First goal is to prove:

\begin{theorem}[The loop cohomology ring]
\label{thm:geometric-fusion-ring}
\label{thm:general-binary-fusion}
\label{thm:general-associativity}
Let $(M,D,A)$ be a polarized SNC compactification of a smooth
affine complex variety $X$. The fusion product of
Definition~\ref{def:loop-fusion-product} makes
$LH^*(X)$ an associative unital graded-commutative
$\F$-algebra, graded by parity. Its unit is the image of
$1\in H^0(X)$ through the derived zero pole stage $B_0\to B$.
The product preserves the contact filtration:
\[
 \mu(F_PLH^i,F_QLH^j)\subseteq F_{P+Q}LH^{i+j}.
\]

The additive identification of
Theorem~\ref{thm:canonical-complex} gives an isomorphism
\[
 H^*(\gr_F\mathbf B)
       \xrightarrow{\;\sim\;} H^*_{\log}(M,D)
\]
of weight-graded, parity-graded unital $\F$-algebras.
Here the product on the left is induced by $\gr_Fm_2$,
and the ring on the right is that of
Definition~\ref{def:log-cohomology}, in parity grading. The summand
$H^*(\gr_P\mathbf B)$ corresponds to the logarithmic summands
of weight $w(\vct)=P$.

If $X$ admits a nowhere-zero algebraic volume form $\Omega$,
all these assertions hold with $B$, $\mathbf B$, and $LH$
replaced by $B_\Omega$, $\mathbf B_\Omega$, and $LH_\Omega$,
respectively, and with the integer-graded
logarithmic ring of Definition~\ref{def:log-cohomology}.
\end{theorem}

Associativity and graded-commutativity will come from moving and
colliding marked source points: the two orders of a
three-point collision give the same specialization proving associativity of the fusion product, while
exchanging two inputs gives graded commutativity.
At maximal contact, the input contact vectors add, and
specialization on the smooth finite models becomes a
diagonal Gysin map. This identifies the leading product
with cup product after restriction to the common boundary
face. The class $1\in H^0(X)$ acts as the identity on these
contact layers; associativity and the pole filtration then
establish the unit on loop cohomology.

\subsection{Associativity and graded commutativity}
\label{subsec:three-point-collision}
We follow the grading convention of
Subsection~\ref{subsec:relative-pole-recovery}: in integer grading,
$B_P$, $B$, $\mathbf B$, and $F_P\mathbf B$ stand for
$B_{\Omega,P}$, $B_\Omega$, $\mathbf B_\Omega$, and
$F_P\mathbf B_\Omega$, respectively, with $\Omega$ the chosen
nowhere-zero algebraic volume form. Without $\Omega$, these
symbols denote the parity-graded objects.

For $I=[k]=\{1,\ldots,k\}$, $k=2,3$, and pole bounds
$\mathbf P=(P_i)_{i\in I}$, recall that
$\mathcal E_{I,\mathbf P}$ is the normalized compact-support
companion of $f_{I,\mathbf P}:\mathcal Z_{I,\mathbf P}\to\A^I$
from Definition~\ref{def:relative-fusion-complex}, constructed
by the finite-factor realization \eqref{eq:relative-compact-factor}.
It is a complex on $(\A^I)^{\mathrm{an}}$, and
Corollary~\ref{cor:relative-pole-recovery}(ii) identifies it as
\[
 \mathcal E_{I,\mathbf P}
   \simeq\mathbb D_{\A^I}\mathcal B_{I,\mathbf P}
   =\mathbb D_{\A^I}
       \bigl(Rf_{I,\mathbf P*}\omega_{I,\mathbf P}\bigr).
\]
Here $\omega_{I,\mathbf P}$ is the normalized dualizing complex
on $\mathcal Z_{I,\mathbf P}^{\mathrm{an}}$, and
$\mathbb D_{\A^I}$ is Verdier duality on the analytic base.
The same notation denotes the periodic companions in parity
grading. 

We compare successive collisions on the two boundary axes of a
quadrant with simultaneous collision along an interior ray.
First we identify the restriction operation on the one-sided
1-parameter binary slice used on those axes.

\begin{lemma}[One-sided 1-parameter binary slice]
\label{lem:one-sided-binary-restriction}
Put $H=[0,\infty)$, $L=(0,\infty)$,
$q_2:H\to\A^2(\C)$, $q_2(t)=(0,t)$, and
$E=q_2^*\mathcal E_{[2],(P,Q)}$. The linear dual of
\[
 B_{P+Q}^\vee\simeq R\Gamma(H,E)
 \xrightarrow{\operatorname{res}}R\Gamma(L,E|_L)
 \simeq(B_P\otimes B_Q)^\vee
\]
agrees with $c_{P,Q}=\overline m_{2;P,Q}$ as a morphism of
coherent pole diagrams, where the input parameters are $z,z-t$
and the output parameter is $z$. Here $\operatorname{res}$ is
the open-star restriction $r_{\varnothing,\{1\}}$ of
Lemma~\ref{lem:parameter-chain-properties}(iii), since
$H=W_{\varnothing}$ and $L=W_{\{1\}}$.
\end{lemma}
\begin{proof}
Apply Remark~\ref{rem:restriction-arrow-descent} to the pullback
of $\mathcal Z_{[2],(P,Q)}$ along $t\mapsto(0,t)$. It identifies
the displayed dual restriction with normalized specialization
along $t>0$. In the source coordinate $z'=z-t/2$, the marked
points become $(-t/2,t/2)=(s,-s)$, where $s=-t/2$.
The identities
\[
 z'-s=z,\qquad z'+s=z-t,\qquad z'|_{t=0}=z
\]
preserve both labeled input parameters and the output parameter.
The source translation therefore identifies this specialization
with that of $\mathcal Z_{P,Q}$ on the inward negative $s$-ray.
Use the affine-source equivariance of
Lemma~\ref{lem:source-coordinate-factorization}(iii), pulled
back along the translation with parameters $(\beta,v)=(-t/2,1)$.
It transports the family and its normalized finite-factor
comparisons from the coordinate $z$ to $z'$. On $t>0$,
parts~(i)--(ii), with the two singleton groups and multipliers
$u_1=u_2=1$, identify the companion with
$\underline{(B_P\otimes B_Q)^\vee}$ through the local
parameters $z'-s=z$ and $z'+s=z-t$. At $t=0$, the same
parts with the single group $\{1,2\}$ identify the companion
with $B_{P+Q}^\vee$ through $z'=z$. Part~(iv) pulls these
comparisons to $L$ and $\{0\}$, respectively, and the
naturality in part~(iii)
identifies the restriction arrow together with these endpoint
comparisons. Thus the transport uses exactly the input and
output identifications in the statement.
The punctured input trivialization by $z-s,z+s$ is defined over
all $\Gm$; the sector comparison of
Lemma~\ref{lem:source-specialization-comparison}(i) identifies
the negative-ray and positive-ray maps. The latter is
$\overline m_{2;P,Q}=c_{P,Q}$ by
Proposition~\ref{prop:relative-binary-identification}.
All source, ray, and duality comparisons commute coherently
with the pole maps.
\end{proof}

Fix $\mathbf P=(P_1,P_2,P_3)$ and put $N=P_1+P_2+P_3$.
Set
\[
 \begin{gathered}
 b:\A^2_{r,s}\longrightarrow\A^3,\qquad b(r,s)=(0,r,r+s),\\
 \mathcal Q=[0,\infty)^2,\qquad
 \mathcal Q^\circ=(0,\infty)^2,\qquad q_3=b|_{\mathcal Q},
 \qquad E=q_3^*\mathcal E_{[3],\mathbf P}.
 \end{gathered}
\]
For a test $\C$-algebra $R$ and $r,s\in R$, the pullback
$b^*\mathcal Z_{[3],\mathbf P}$ has completed source and
marked pole conditions
\begin{equation}\label{eq:general-triple-poles}
 \widehat{R[z]}_{(z(z-r)(z-r-s))},\qquad
 z^{r_0P_1}(z-r)^{r_0P_2}(z-r-s)^{r_0P_3}g_j(\gamma)
 \in\widehat{R[z]}_{(z(z-r)(z-r-s))}\quad\text{for every }j.
\end{equation}
Here $g_j$ are the affine generators chosen using the fixed
very ample multiple $r_0A$ of the polarization divisor, and
$\gamma$ maps the punctured completed source to $X$.
The three input coordinates are $z,z-r,z-r-s$; the output
coordinate is $z$. 

Define
\[
 c_{P_1,P_2,P_3}:=
 \left[R\Gamma(\mathcal Q,E)\longrightarrow
       R\Gamma(\mathcal Q^\circ,E)\right]^\vee:
 B_{P_1}\otimes B_{P_2}\otimes B_{P_3}\longrightarrow B_N,
\]
the dual of the open-star restriction $r_{\varnothing,\{1,2\}}$
of Lemma~\ref{lem:parameter-chain-properties}(iii), using the
equivalences
\[
 \begin{aligned}
 R\Gamma(\mathcal Q,E)&\simeq E_{(0,0)}
                         \simeq B_N^\vee,\\
 R\Gamma(\mathcal Q^\circ,E)&\simeq
       (B_{P_1}\otimes B_{P_2}\otimes B_{P_3})^\vee.
 \end{aligned}
\]
The first uses open-star evaluation from
Lemma~\ref{lem:parameter-chain-properties}(ii) and the diagonal
identification of Lemma~\ref{lem:relative-fusion-factorization}(ii):
at the corner all three marked points coincide, and their pole
bounds add to $N$. Here $E_{(0,0)}$ denotes the stalk at the corner.
The second uses separation of the distinct points $0,r,r+s$
on $\mathcal Q^\circ$, with parameters $z,z-r,z-r-s$, which
identifies $E|_{\mathcal Q^\circ}$ with the constant complex
of the displayed value; $\mathcal Q^\circ$ is contractible.
The dual map uses these equivalences and canonical biduality.

\begin{lemma}[Composition of collisions]
\label{lem:collision-composition}

 Then
\begin{equation}\label{eq:ternary-collision-composition}
 \begin{aligned}
 c_{P_1+P_2,P_3}(c_{P_1,P_2}\otimes\mathrm{id})
 &=c_{P_1,P_2+P_3}(\mathrm{id}\otimes c_{P_2,P_3})\\
 &=c_{P_1,P_2,P_3}
 \end{aligned}
\end{equation}
 in the homotopy category of coherent pole-bound diagrams
$\operatorname{Fun}(\mathbb N^3,\mathcal D(\F))$
(and in its parity-graded version).
\end{lemma}
\begin{proof}
Lemma~\ref{lem:relative-fusion-factorization}(ii) gives
constructibility of $\mathcal E_{[3],\mathbf P}$ for the
equality strata of $\A^3$ and its diagonal and separation
comparisons. Under $q_3(r,s)=(0,r,r+s)$, the corner, the
open $r$-axis, the open $s$-axis, and the interior have equality
groups $\{1,2,3\}$, $\{1\},\{2,3\}$,
$\{1,2\},\{3\}$, and $\{1\},\{2\},\{3\}$,
respectively. Pulling back those comparisons therefore makes
$E=q_3^*\mathcal E_{[3],\mathbf P}$ constructible for the
four-face stratification of $\mathcal Q$.
For $J\subseteq\{1,2\}$ put $G(J)=R\Gamma(W_J,E)$, where
$W_J$ is the open star of the coordinate stratum indexed by $J$
as in Lemma~\ref{lem:parameter-chain-properties}(ii).
Part~(iii) of that lemma gives the coherent open-star restriction square
\begin{equation}\label{eq:ternary-star-square}
 \begin{array}{ccc}
 G(\varnothing)&\xrightarrow{r_{\varnothing,\{1\}}}&G(\{1\})\\
 \big\downarrow\rlap{$\scriptstyle r_{\varnothing,\{2\}}$}&&
 \big\downarrow\rlap{$\scriptstyle r_{\{1\},\{1,2\}}$}\\
 G(\{2\})&\xrightarrow{r_{\{2\},\{1,2\}}}&G(\{1,2\}).
 \end{array}
\end{equation}
At the corner the single group has parameter $z$. On the $s$-axis
$r=0$, the groups $\{1,2\},\{3\}$ have parameters $z,z-s$;
on the $r$-axis $s=0$, the groups $\{1\},\{2,3\}$ have
parameters $z,z-r$. In the interior the singleton parameters
are $z,z-r,z-r-s$. For each of these four groupings, let
$J_\pi$ be its set of groups and let $\pi:[3]\twoheadrightarrow J_\pi$
send each label to its group.
Lemma~\ref{lem:source-coordinate-factorization}(i) uses precisely
these parameters in $\Phi_{\pi,1}$, with all multipliers $u_j=1$; part~(ii) gives the
comparison $e_{\pi,1}$ of \eqref{eq:source-coordinate-comparisons}.
By part~(iv), its pullback identifies $E$ on the corresponding
real face with the constant complex of value
$B_{\pi,\mathbf P}^{\vee}$, where
$B_{\pi,\mathbf P}=\bigotimes_{j\in J_\pi}
B_{\sum_{i\in\pi^{-1}(j)}P_i}$.
The open-star evaluation equivalences $\epsilon_J$ at the
points $v_J$ of Lemma~\ref{lem:parameter-chain-properties}(ii)
consequently identify the vertices as
\[
 \begin{aligned}
 G(\varnothing)&\simeq B_N^\vee,&
 G(\{1\})&\simeq(B_{P_1}\otimes B_{P_2+P_3})^\vee,\\
 G(\{2\})&\simeq(B_{P_1+P_2}\otimes B_{P_3})^\vee,&
 G(\{1,2\})&\simeq(B_{P_1}\otimes B_{P_2}\otimes B_{P_3})^\vee.
 \end{aligned}
\]
To identify the arrows on the axes, put
\[
 \begin{aligned}
 d_s:H\to\mathcal Q,&\quad d_s(t)=(0,t),&
 \Delta_{12}:\A^2\to\A^3,&\quad\Delta_{12}(a,b)=(a,a,b),\\
 d_r:H\to\mathcal Q,&\quad d_r(t)=(t,0),&
 \Delta_{23}:\A^2\to\A^3,&\quad\Delta_{23}(a,b)=(a,b,b),
 \end{aligned}
\]
where $H=[0,\infty)$ and $q_2(t)=(0,t)$. The identities
$q_3d_s=\Delta_{12}q_2$ and $q_3d_r=\Delta_{23}q_2$,
together with the ordinary-inverse-image form of the diagonal
comparison in Lemma~\ref{lem:relative-fusion-factorization}(ii),
give
\[
 d_s^*E\simeq q_2^*\mathcal E_{[2],(P_1+P_2,P_3)},\qquad
 d_r^*E\simeq q_2^*\mathcal E_{[2],(P_1,P_2+P_3)}.
\]
For $d_s$, the coordinate-face restriction comparisons of
Lemma~\ref{lem:parameter-chain-properties}(iii), with $J=\{2\}$,
compare $r_{\varnothing,\{2\}}:G(\varnothing)\to G(\{2\})$
with the open-star restriction on this axis. The vertical
comparisons are equivalences by the open-star evaluations
of part~(ii).
Lemma~\ref{lem:one-sided-binary-restriction} therefore identifies
its dual with $c_{P_1+P_2,P_3}$. The same argument on $s=0$
identifies the dual of the upper arrow with $c_{P_1,P_2+P_3}$.

On $W_{\{2\}}=\{s>0\}$, the third point is separated from
the first two. Apply the separation comparison for
$\mathcal E_{[3],\mathbf P}$ of
Lemma~\ref{lem:relative-fusion-factorization}(ii) to the groups
$\{1,2\},\{3\}$. For the singleton, use
Lemma~\ref{lem:source-coordinate-factorization}(i)--(ii) with
center $r+s$ and multiplier $1$: its parameter is $z-(r+s)$
and its companion is the constant complex $B_{P_3}^\vee$.
Part~(iv) pulls this singleton comparison back to
$W_{\{2\}}$; the separation comparison for $\mathcal E$
is pulled back by $(q_3|_{W_{\{2\}}})^*$. Writing
$\operatorname{pr}_r:W_{\{2\}}\to H$ for $(r,s)\mapsto r$,
we obtain
\[
 E|_{W_{\{2\}}}\simeq
 \operatorname{pr}_r^*\bigl(q_2^*\mathcal E_{[2],(P_1,P_2)}\bigr)
       \otimes\underline{B_{P_3}^\vee}_{W_{\{2\}}}.
\]
Lemma~\ref{lem:source-specialization-comparison}(ii) identifies
the restriction $G(\{2\})\to G(\{1,2\})$ with the pair's
restriction in $r$, tensored with $\mathrm{id}_{B_{P_3}^\vee}$.
Taking sections in the contractible coordinate $s>0$ and
dualizing therefore gives $c_{P_1,P_2}\otimes\mathrm{id}$
for the dual of the lower arrow.

On $W_{\{1\}}=\{r>0\}$, pull back the separation comparison
for $\mathcal E_{[3],\mathbf P}$ of
Lemma~\ref{lem:relative-fusion-factorization}(ii) by
$(q_3|_{W_{\{1\}}})^*$, for the groups $\{1\},\{2,3\}$.
Set $z'=z-r$ for the second group.
Its centers are $0,s$ and its parameters are $z'=z-r$,
$z'-s=z-r-s$; at $s=0$ its output parameter is $z-r$.
Lemma~\ref{lem:source-coordinate-factorization}(iii), with
translation parameters $(\beta,v)=(-r,1)$, transports this
pair to the one-sided binary family. Parts~(i)--(ii) identify
its endpoint companions in these parameters and the singleton
companion at $0$ with $B_{P_1}^\vee$; part~(iv) restricts
the comparisons to the real faces. Thus, for
$\operatorname{pr}_s:W_{\{1\}}\to H$, $(r,s)\mapsto s$,
\[
 E|_{W_{\{1\}}}\simeq
 \underline{B_{P_1}^\vee}_{W_{\{1\}}}\otimes
 \operatorname{pr}_s^*\bigl(q_2^*\mathcal E_{[2],(P_2,P_3)}\bigr).
\]
The restriction-arrow assertion of
Lemma~\ref{lem:source-specialization-comparison}(ii), followed
by sections in $r>0$ and duality, gives
$\mathrm{id}\otimes c_{P_2,P_3}$ for the dual of the right arrow.

These are comparisons of the entire restriction square.
Lemma~\ref{lem:relative-fusion-factorization}(i)--(ii) supplies
the diagonal and separation maps together with their composition
homotopies. Lemma~\ref{lem:source-coordinate-factorization}(iii)
transports these homotopies under the affine substitutions and
the successive groupings; part~(iv) supplies their restrictions to the real
faces. In particular, the intermediate parameters $z$ for
$\{1,2\}$ and $z-r$ for $\{2,3\}$ agree in the axis
and separation comparisons above. Together with the
restriction-arrow comparison of
Lemma~\ref{lem:source-specialization-comparison}(ii), this
preserves the composition homotopy of the star-restriction square.
Thus both composites in \eqref{eq:ternary-star-square} agree with
the open-star restriction
$r_{\varnothing,\{1,2\}}:R\Gamma(\mathcal Q,E)\to
R\Gamma(\mathcal Q^\circ,E)$. Duality gives
\eqref{eq:ternary-collision-composition}.
The same transformations commute coherently with every pole map
and comparison of finite-factor presentations, proving the asserted identity
of pole-bound diagrams.
\end{proof}

For $r_1,s_1>0$, put
\[
 \ell_{r_1,s_1}:\A^1_t\longrightarrow\A^2_{r,s},\qquad
 t\longmapsto(tr_1,ts_1),
\]
and let $\operatorname{sp}_{r_1,s_1}$ be the normalized
specialization for $(b\ell_{r_1,s_1})^*\mathcal Z_{[3],\mathbf P}$
along $t>0$, as in Definition~\ref{def:normalized-ray-specialization}.
The following corollary will be used in
Theorem~\ref{mir:thm-low-trace-comparison} to identify a
three-input collision trace with the trace of an iterated product.

\begin{corollary}[Three-input specialization]
\label{lem:three-point-collision}
In the homotopy category of coherent pole-bound diagrams,
\[
 \operatorname{sp}_{r_1,s_1}=c_{P_1,P_2,P_3}:
 B_{P_1}\otimes B_{P_2}\otimes B_{P_3}\longrightarrow B_N.
\]
In particular the map is independent of $r_1,s_1>0$.
\end{corollary}
\begin{proof}
Restrict $\ell=\ell_{r_1,s_1}$ to $H=[0,\infty)$ and put
$L=(0,\infty)$. Section-to-pullback maps give the coherent square
\[
 \begin{array}{ccc}
 R\Gamma(\mathcal Q,E)&\longrightarrow&R\Gamma(\mathcal Q^\circ,E)\\
 \big\downarrow&&\big\downarrow\\
 R\Gamma(H,\ell^*E)&\longrightarrow&R\Gamma(L,\ell^*E).
 \end{array}
\]
The left vertical arrow is an equivalence by evaluation at the
corner; the right one is an equivalence by the constant interior
trivialization, using Lemma~\ref{lem:parameter-chain-properties}(ii).
The dual of the top arrow is $c_{P_1,P_2,P_3}$.
Along $\ell(t)=(r_1t,s_1t)$, the input coordinates above become
$z,z-r_1t,z-(r_1+s_1)t$ and the output coordinate remains $z$.
For $t>0$, apply Lemma~\ref{lem:source-coordinate-factorization}(i)
with the three singleton groups, centers $0,r_1t,(r_1+s_1)t$,
and multipliers $u_1=u_2=u_3=1$. At $t=0$, apply the same
part with the single group $[3]$, center $0$, and multiplier
$1$; its pole bound is $N$ and its parameter is $z$.
Part~(ii) supplies the corresponding companion comparisons,
and part~(iv) pulls them back to $L$ and $\{0\}$, respectively:
\[
 (\ell^*E)|_L\simeq
 \underline{(B_{P_1}\otimes B_{P_2}\otimes B_{P_3})^\vee}_L,
 \qquad
 (\ell^*E)_0\simeq B_N^\vee.
\]
Here $(\ell^*E)_0$ denotes the stalk at $t=0$. Further-pullback
compatibility in part~(iv) identifies these maps with the
pullbacks along $\ell$ of the comparisons on $\mathcal Q^\circ$
and at the corner. Taking derived global sections of their
Verdier-dual comparisons in parts~(ii) and~(iv), with the exterior
product and inward-ray class specified in
Definition~\ref{def:normalized-ray-specialization}, gives exactly
its input and output maps $\alpha_{\mathbf P},\beta_{\mathbf P}$.
The naturality in part~(iii) includes the pole maps in these
comparisons. Remark~\ref{rem:restriction-arrow-descent} identifies
the dual of the bottom restriction arrow with
$\operatorname{sp}_{r_1,s_1}$.
The square therefore gives the asserted equality, compatibly
with every pole map. Its interpretation as either iterated binary
product follows from Lemma~\ref{lem:collision-composition} and
Proposition~\ref{prop:relative-binary-identification}.
\end{proof}

\begin{proposition}[Associativity and graded commutativity]
\label{prop:binary-associativity}
The stage maps satisfy
\[
 \begin{gathered}
 \overline m_{2;P+Q,R}(\overline m_{2;P,Q}\otimes\mathrm{id})
 =\overline m_{2;P,Q+R}(\mathrm{id}\otimes\overline m_{2;Q,R}),\\
 \overline m_{2;P,Q}=\overline m_{2;Q,P}\circ\tau,
 \end{gathered}
\]
in the homotopy categories of the respective pole-bound diagrams,
where $\tau$ interchanges the factors with the Koszul sign.
Consequently the filtered derived product $\overline m_2$ is associative
and graded commutative; the comparisons are represented by
filtered chain homotopies on joint cochain representatives.
For homogeneous loop cohomology classes
\[
 \mu(\mu(a,b),c)=\mu(a,\mu(b,c)),\qquad
 \mu(a,b)=(-1)^{|a||b|}\mu(b,a).
\]
These statements hold in parity grading and in the integer grading
determined by $\Omega$.
\end{proposition}
\begin{proof}
Associativity is \eqref{eq:ternary-collision-composition}, using
Proposition~\ref{prop:relative-binary-identification}.
For commutativity, exchange the two marks in $\mathcal Z_{P,Q}$.
In the symmetric slice this sends $s$ to $-s$, fixes the output
coordinate $z$, and exchanges the input parameters $z-s,z+s$.
More explicitly, part~(i) of
Lemma~\ref{lem:source-coordinate-factorization}, with the two
singleton groups and $u_1=u_2=1$, identifies this relabeling
over $s\ne0$ with
\[
 \begin{aligned}
 \Gm\times F_P\mathscr LX\times F_Q\mathscr LX
   &\longrightarrow\Gm\times F_Q\mathscr LX\times F_P\mathscr LX,\\
 (s,\gamma_1,\gamma_2)&\longmapsto(-s,\gamma_2,\gamma_1).
 \end{aligned}
\]
Indeed, the new first parameter is $z-(-s)=z+s$ and the new
second parameter is $z+(-s)=z-s$, with multipliers still $1$.
Part~(ii) transports this interchange to the normalized
complexes; the relabeling naturality in part~(iii) identifies
the input map with the fixed Koszul symmetry
$\tau:B_P\otimes B_Q\to B_Q\otimes B_P$.
At $s=0$, parts~(i)--(ii) use the single group $\{1,2\}$
and the unchanged parameter $z$; by the diagonal compatibility
in part~(iii), the output map is $\mathrm{id}_{B_{P+Q}}$.
Part~(iv) pulls these comparisons to the positive and negative
rays, preserving the pole-map naturality of part~(iii).
Lemma~\ref{lem:source-specialization-comparison}(i) identifies
the inward negative-ray operation with the positive-ray operation:
their difference is the boundary of an oriented sector in $\C^*$.
Source relabeling and ray comparison preserve the complex
orientations and commute coherently with pole maps.
Thus only the Koszul sign of $\tau$ remains.
Lemma~\ref{lem:filtered-operations-colimit} carries both
pole-diagram identities to filtered derived morphisms on $B$;
joint cochain representatives realize the comparisons by filtered
chain homotopies. Taking cohomology gives the displayed formulas.
\end{proof}

\subsection{The leading product and the unit}
\label{lp:subsec:fusion-leading}
\label{subsec:general-ring-properties}
We identify the leading product with the common-face product of
Section~\ref{sec:setup}; this will also determine the action of the unit.

The next lemma extends to the boundary coordinates of an SNC
chart the factorization by marked-point powers and unit series
used for the multiplicative group in
Subsection~\ref{ex:cstar-global-fusion}.

\begin{lemma}[Highest-contact collision locus]
\label{lem:fusion-leading-contact}
Fix contact vectors $\vct,\wct$, and put
\[
 P=w(\vct),\quad Q=w(\wct),\quad N=P+Q,\qquad
 I=\supp\vct,\quad J=\supp\wct,\quad L=I\cup J.
\]
In the fixed-source family $\mathcal Z_{P,Q}\to\A^1_s$ of
Example~\ref{ex:kv-binary-family}, define
\[
 \begin{aligned}
 \mathcal C^\times_{\vct,\wct}
   &:=\Gm\times\mathscr P_{\vct}\times\mathscr P_{\wct},\\
 \mathcal H_{\vct,\wct}
   &:=\overline{\mathcal C^\times_{\vct,\wct}}
       \cap(\mathcal Z_{P,Q})_0
       \cap(F_N\mathscr LX\setminus F_{N-1}\mathscr LX).
 \end{aligned}
\]
The first locus is embedded by the inverse of
\eqref{eq:source-coordinate-identification} in
Lemma~\ref{lem:source-coordinate-factorization}(i), with the
two singleton groups, centers $s,-s$, and multipliers $1,1$.
Thus its exact input contacts are $(\vct,\wct)$ in the
parameters $z-s,z+s$. The closure is schematic; the displayed
intersection has its reduced induced structure, and the output
identification $(\mathcal Z_{P,Q})_0=F_N\mathscr LX$ uses $z$.
Then the following hold.
\begin{enumerate}[(i),wide=0pt]
\item One has
\[
                    \mathcal H_{\vct,\wct}
                         =\mathscr P_{\vct+\wct}.
\]
In particular, its output loops have center in $D_L^\circ$.

\item For every output point $\lambda\in\mathcal H_{\vct,\wct}$,
there exist an open neighborhood
\[
 \mathcal U\subset
 \overline{\mathcal C^\times_{\vct,\wct}}
 \subset\mathcal Z_{P,Q}
\]
of $\lambda$ and an affine open $V\subset M$ containing its
center, with the following properties. One can write
\[
 \mathcal U=\overline{\mathcal C^\times_{\vct,\wct}}\cap D(f),
 \qquad f\in\Gamma(\mathcal Z_{P,Q},\mathcal O_{\mathcal Z_{P,Q}}),
\]
where $D(f)$ is the principal open where $f$ is invertible,
and $f$ depends on $s$ and finitely many Laurent coefficients.
The open $V$ meets no boundary components other than those
indexed by $L$ and admits an \'etale coordinate system
containing local equations $x_i$ for $D_i\cap V$, $i\in L$.

For every complex algebra $R$ and every morphism
$\Spec R\to\mathcal U$, let $s\in R$ be the pullback of the
coordinate on $\A^1_s$, and put
\[
 \mathcal O_s:=R\llbracket q_0 \rrbracket[z]/(z^2-s^2-q_0).
\]
Here $q_0$ is the formal completion variable, satisfying
$q_0=(z-s)(z+s)$. The corresponding point of
$\mathcal Z_{P,Q}$ determines a punctured-source map
\[
 \gamma:\Spec\mathcal O_s[q_0^{-1}]\longrightarrow X.
\]
Viewed as a map to $M$, it extends uniquely to
\[
 \bar\gamma:\Spec\mathcal O_s\longrightarrow V,
 \qquad
 \bar\gamma|_{\Spec\mathcal O_s[q_0^{-1}]}=\gamma.
\]
Its boundary-coordinate functions have the factorization
\begin{equation}\label{lp:eq:fusion-snc-factors}
 X_i:=\bar\gamma^*x_i
     =(z-s)^{v_i}(z+s)^{w_i}U_i,
 \qquad U_i\in\mathcal O_s^\times,\quad i\in L.
\end{equation}
The extension and these unit series are compatible with
change of the test algebra $R$.

\item For the neighborhood $\mathcal U$ and chart $V$ in (ii),
put $T=\A^1_s$, $n=\dim M$, and $\ell=|L|$.
Choose an order $L=\{i_1,\ldots,i_\ell\}$ and write the
\'etale coordinate map as
\[
 \phi=(x_{i_1},\ldots,x_{i_\ell},y_{\ell+1},\ldots,y_n):
 V\longrightarrow\A^n.
\]
For the data of (ii), set
\[
 C_a=\begin{cases}
 U_{i_a},&1\le a\le\ell,\\
 \bar\gamma^*y_a,&\ell<a\le n.
 \end{cases}
\]
Since $1,z$ form an $R\llbracket q_0\rrbracket$-basis of
$\mathcal O_s$, there are unique coefficients $u_{a,r},\xi_{a,r}\in R$
such that
\[
 C_a=\sum_{r\ge0}(u_{a,r}+z\xi_{a,r})q_0^r.
\]
For $h\ge1$, put $C_a^{(h)}:=C_a\bmod q_0^h$, use
$u_{a,r},\xi_{a,r}$ also for the independent coordinate
functions on $\A^{2nh}$, and define
\[
 \mathcal D_h=
 D\!\left(\prod_{a=1}^{\ell}(u_{a,0}^2-s^2\xi_{a,0}^2)\right)
 \subset\A^1_s\times\A^{2nh},
 \qquad 1\le a\le n,\quad 0\le r<h.
\]
Thus an $R$-point of $\mathcal D_h$ records $s\in R$ and
$n$ truncated series $C_a^{(h)}\in\mathcal O_s/q_0^h$,
the first $\ell$ of which are units. Equivalently,
$u_{a,0}^2-s^2\xi_{a,0}^2\in R^\times$ for $a\le\ell$.
For $h'\ge h\ge1$, let
$d_{h',h}:\mathcal D_{h'}\to\mathcal D_h$ forget the
coefficients with $h\le r<h'$.

There are $h_0\ge1$, affine finite-type $T$-schemes $\mathcal W_h^\circ$
for $h\ge h_0$, \'etale maps
$c_h:\mathcal W_h^\circ\to\mathcal D_h$, and truncation maps
$\tau_{h',h}:\mathcal W_{h'}^\circ\to\mathcal W_h^\circ$
for $h'\ge h\ge h_0$, satisfying
\[
 c_h\tau_{h',h}=d_{h',h}c_{h'},\qquad
 \mathcal U\cong\varprojlim_{h\ge h_0}\mathcal W_h^\circ.
\]
The composite $\mathcal U\to\mathcal W_h^\circ\xrightarrow{c_h}
\mathcal D_h$ records the series $C_a^{(h)}$ defined above.
Each $\mathcal W_h^\circ\to T$ is smooth, and there are
product presentations
\[
 \mathcal W_{h'}^\circ\cong\mathcal W_h^\circ\times
 \Spec\C[u_{a,r},\xi_{a,r}:1\le a\le n,\ h\le r<h']
 \quad\text{over }\mathcal W_h^\circ,
\]
in which $\tau_{h',h}$ is projection onto the first factor.

\item For two choices $(\mathcal U,V)$ and $(\mathcal U',V')$
as in (ii)--(iii), use primed symbols for the second choice.
The overlap of the collision neighborhoods $\mathcal U\cap\mathcal U'$
has finite models as follows. For all sufficiently large $h$,
there are open subschemes and isomorphisms over $T$
\[
 O_h\subset\mathcal W_h^\circ,\qquad
 O_h'\subset(\mathcal W_h')^\circ,
 \qquad \chi_h:O_h\xrightarrow{\sim}O_h',
\]
whose inverse limits are identified with $\mathcal U\cap\mathcal U'$.
On the target-chart overlap $V\cap V'$, write $x_i'=a_ix_i$
with $a_i$ an invertible regular function. The maps $\chi_h$
are induced by $U_i'=(\bar\gamma^*a_i)U_i$ and by pullback
of the new tangential coordinates. They satisfy
\[
 \tau'_{h',h}\chi_{h'}=\chi_h\tau_{h',h}
 \quad\text{on }O_{h'}.
\]
Reversing the chart change gives $\chi_h^{-1}$; for three
choices, the composite from the first to the second and
then to the third equals the direct comparison on their
common open submodels.
\end{enumerate}
\end{lemma}
\begin{proof}
\emph{Step 1: construct an extension to an affine chart for (ii).}
Fix any point $\lambda\in\mathcal H_{\vct,\wct}$ and write
$\kappa=\kappa(\lambda)$ for its residue field.
Since $\mathcal Z_{P,Q}$ is affine, begin with an ambient
principal open $D(f_0)$ containing $\lambda$ and put
\[
 \Spec R=
 \overline{\mathcal C^\times_{\vct,\wct}}\cap D(f_0).
\]
Work with the $R$-point of $\mathcal Z_{P,Q}$ defined by
the composite
\[
 \Spec R\hookrightarrow
 \overline{\mathcal C^\times_{\vct,\wct}}
 \hookrightarrow\mathcal Z_{P,Q}.
\]
The composite $\Spec R\to\mathcal Z_{P,Q}$ specifies
$s\in R$ and $\gamma:\Spec\mathcal O_s[q_0^{-1}]\to X$
as in the statement. After each localization of this coefficient
algebra, use the same notation for the pulled-back data.

Choose $r>0$ and sections $s_0,\ldots,s_m$ as in
Proposition~\ref{prop:loop-sectors}, so that
$\operatorname{div}(s_0)=rA$ and $g_j=s_j/s_0$ are the
fixed pole generators. Put
\[
 H=(z-s)^{rP}(z+s)^{rQ},\qquad
 g_0=1,\qquad G_j=Hg_j(\gamma)\in\mathcal O_s
 \quad(0\le j\le m).
\]
The inclusions $G_j\in\mathcal O_s$ are the defining pole
conditions of $\mathcal Z_{P,Q}$.
At $\lambda$, the completed source is $\kappa\llbracket z \rrbracket$, with
$q_0=z^2$. Properness of $M$ gives the unique extension
$\bar\gamma_\lambda:\Spec\kappa\llbracket z \rrbracket\to M$; write
$p_\lambda=\bar\gamma_\lambda(0)$.
In a local trivialization of $\mathcal O_M(rA)$, exact output
weight $N$ gives
\[
 \bar\gamma_\lambda^*s_0=z^{rN}a(z),\qquad
 a(z)\in\kappa\llbracket z \rrbracket^\times.
\]
Choose $j$ with $s_j(p_\lambda)\ne0$. Then
$G_j|_\lambda=\bar\gamma_\lambda^*s_j/a(z)$ is a unit.

The algebra $\mathcal O_s/q_0\mathcal O_s=R[z]/(z^2-s^2)$
is free of rank two. Invert the determinant of multiplication
by $G_j$ on this algebra; it is a coefficient function nonzero
at $\lambda$. The element $G_j$ is then a unit in
$\mathcal O_s$, by $q_0$-adic completeness. Thus the ratios
\[
 (G_0/G_j,\ldots,G_m/G_j)
\]
define a map $\bar\gamma:\Spec\mathcal O_s\to
M\cap\{s_j\ne0\}$ extending $\gamma$ as a map to $M$.
Indeed, the equations of this affine chart of $M$ hold after
inverting $q_0$, and
\[
 \mathcal O_s\hookrightarrow\mathcal O_s[q_0^{-1}]
\]
is injective: $\mathcal O_s$ is free over $R\llbracket q_0 \rrbracket$, on
which multiplication by $q_0$ is injective. Every coefficient
of these ratios is a regular function of finitely many Laurent
coefficients and the inverse of the determinant.

\emph{Step 2: refine the extension to an SNC chart for (ii).}
Let $K$ index the boundary components through $p_\lambda$. Choose an affine open
$V\subset M\cap\{s_j\ne0\}$ containing $p_\lambda$,
meeting no boundary components outside $K$, with an \'etale
coordinate map $V\to\A^n$ containing local equations $x_i$
for $D_i\cap V$, $i\in K$.
We may choose $V$ to be principal in $M\cap\{s_j\ne0\}$.
Pull back its defining function by $\bar\gamma$ and invert
the determinant of multiplication by its reduction modulo
$q_0$. This determinant is nonzero at $\lambda$; after this
localization the pulled-back function is a unit, and hence
\[
                  \bar\gamma:\Spec\mathcal O_s\to V.
\]
The two marked-center maps are obtained by evaluation at
$(q_0,z)=(0,s)$ and $(0,-s)$. On the input contact locus their
images lie in $D_I$ and $D_J$, respectively. These are closed
conditions, so they hold on the input closure wherever the
extension is defined. At $s=0$ the two centers coincide.
Consequently $L\subseteq K$.

\emph{Step 3: prove $\mathcal H_{\vct,\wct}\subseteq\mathscr P_{\vct+\wct}$ and the center assertion in (i).}
For $i\in K$, put
\[
 X_i=\bar\gamma^*x_i,\qquad
 H_i=(z-s)^{v_i}(z+s)^{w_i},
\]
where contact entries outside the respective supports are zero.
For any integer $b\ge\max(v_i,w_i)$,
\[
 H_i(z-s)^{b-v_i}(z+s)^{b-w_i}=q_0^b.
\]
Thus multiplication by $H_i$ on $\mathcal O_s$ is injective.
Moreover, since $H_i$ is monic and $q_0^b=0$ modulo $H_i$,
\[
 \mathcal O_s/(H_i)\cong R[z]/(H_i),\qquad
 \operatorname{rank}_R R[z]/(H_i)=v_i+w_i.
\]
The condition $X_i\in H_i\mathcal O_s$ is therefore the
vanishing of finitely many remainder coefficients. It holds
on the separated input contact locus, by its exact contacts
in $z-s,z+s$, and hence on its schematic closure. These
equations are imposed over the coefficient algebra itself,
without passing to its reduction.

At $\lambda$, set $\nu_i=\ord_z\bar\gamma_\lambda^*x_i$.
Divisibility and exact output weight give
\[
 \nu_i\ge v_i+w_i,\qquad
 \sum_{i\in K}k_i\nu_i=N
       =\sum_{i\in K}k_i(v_i+w_i).
\]
Since every $k_i>0$, all the inequalities are equalities.
If $i\in K\setminus L$, this would give $\nu_i=0$, contrary
to $p_\lambda\in D_i$. Hence $K=L$ and the contact vector
of $\lambda$ is $\vct+\wct$. This proves
\[
 \mathcal H_{\vct,\wct}\subseteq\mathscr P_{\vct+\wct}
\]
on underlying points, together with the center assertion
$p_\lambda\in D_L^\circ$ in (i). The equality $K=L$ also
gives the boundary-support condition on $V$ in (ii).
If $D_L^\circ=\varnothing$, both
loci are empty and the local assertions are vacuous; this
already proves the lemma in that case.

\emph{Step 4: complete (ii) by obtaining unit factors, uniqueness, and base change.}
By Step~3 there is a unique quotient $U_i\in\mathcal O_s$
with $X_i=H_iU_i$, for every $i\in L$. At $\lambda$ this
quotient is a unit, because $\nu_i=v_i+w_i$.
Each Laurent coefficient of $U_i$ depends regularly on
finitely many coefficients of $X_i$ and on $s$: for $b\ge\max(v_i,w_i)$,
\[
 q_0^bU_i=(z-s)^{b-v_i}(z+s)^{b-w_i}X_i.
\]
Invert the determinant of multiplication by each $U_i$
modulo $q_0$. There are only finitely many such determinants,
all nonzero at $\lambda$. Together with the localizations in
Steps~1--2 they define an open neighborhood
\[
 \lambda\in\mathcal U\subset
       \overline{\mathcal C^\times_{\vct,\wct}}
\]
obtained by inverting finitely many regular functions, each
depending on $s$ and finitely many Laurent coefficients.
At each successive localization, write the next function
as $a/d^m$, where $a,d$ lie in the coordinate ring of the
closure, $m\ge0$, and $d$ is already invertible.
Its invertibility is equivalent to that of $a$.
Lift each numerator $a$ to the coordinate ring of
$\mathcal Z_{P,Q}$. Multiplying these finitely many lifts
and $f_0$ gives an ambient regular function $f$ with
\[
 \mathcal U=\overline{\mathcal C^\times_{\vct,\wct}}\cap D(f),
\]
as required in (ii). On $\mathcal U$ every
$U_i$ is a unit, giving
\eqref{lp:eq:fusion-snc-factors}. Since $K=L$, the open $V$
of Step~2 has exactly the properties required in the statement.

Apply these universal formulas to any morphism
$\Spec R\to\mathcal U$. They give the required extension
$\bar\gamma:\Spec\mathcal O_s\to V$ and its unit series.
The extension is unique: $V$ is affine, and two homomorphisms
from its coordinate ring to $\mathcal O_s$ that agree in
$\mathcal O_s[q_0^{-1}]$ agree in $\mathcal O_s$ by the
injectivity in Step~1.
For a homomorphism of complex algebras $\varphi:R\to R'$,
write $s'=\varphi(s)$. Applying $\varphi$ coefficientwise
gives the homomorphism
\[
 \mathcal O_s\longrightarrow
 \mathcal O_{s'}:=R'\llbracket q_0 \rrbracket[z]/(z^2-(s')^2-q_0).
\]
It preserves the displayed coordinate formulas and the
identities $X_i=H_iU_i$. The uniqueness of the extension and
of each quotient therefore identifies their images with the
constructions for the composite $\Spec R'\to\Spec R\to
\mathcal U$. This proves compatibility with change of the
test algebra, including nonreduced algebras; it uses the
coefficientwise map, not an identification of ordinary tensor
products with completed power-series algebras.
This completes (ii).

\emph{Step 5: verify the pole bounds and separated contacts of normal-form maps.}
To prepare the reverse inclusion in (i) and the construction
in (iii), use any affine SNC chart $V$ at
a point of $D_L^\circ$, meeting no boundary components outside
$L$ and with boundary equations $x_i$, $i\in L$.
For any complex algebra $R$ and $s\in R$, form
$\mathcal O_s$ as in (ii) and put
$H_i=(z-s)^{v_i}(z+s)^{w_i}$. Consider any map
$\bar\gamma:\Spec\mathcal O_s\to V$ satisfying
\eqref{lp:eq:fusion-snc-factors}, with each $U_i$ a unit.
After inverting $q_0$, both $z-s$ and $z+s$ are units, so
this map avoids all boundary components and restricts to a
map into $X$. Locally on $V$, the denominator of every
$g_j=s_j/s_0$ divides $\prod_{i\in L}x_i^{rk_i}$, and
\[
 \prod_{i\in L}(H_iU_i)^{rk_i}
       =H\prod_{i\in L}U_i^{rk_i}.
\]
Hence $Hg_j(\gamma)\in\mathcal O_s$ for every $j$,
which is precisely the pole condition for
$\gamma$ to define a point of $\mathcal Z_{P,Q}$.
When $s$ is invertible, the other marked factor is a unit
on each completed disk. The exact input contacts are therefore
$\vct,\wct$ in $z-s,z+s$, respectively. Thus the restriction
of every such family to the base locus where $s$ is invertible
maps to $\mathcal C^\times_{\vct,\wct}$.

\emph{Step 6: identify the coefficient models needed for (iii).}
For the series $C_a$ in (iii), multiplication by $u_{a,0}+z\xi_{a,0}$ on
$R[z]/(z^2-s^2)$ has matrix
\[
 \begin{pmatrix}
 u_{a,0}&s^2\xi_{a,0}\\
 \xi_{a,0}&u_{a,0}
 \end{pmatrix}.
\]
Its determinant is $u_{a,0}^2-s^2\xi_{a,0}^2$.
Thus $\mathcal D_h$ in the statement parametrizes precisely
these truncated series, with the first $\ell$ series invertible;
the tangential coefficients are unrestricted.
Set $E=(\Gm)^\ell\times\A^{n-\ell}$,
$u=(u_{a,0})_a$, and $\xi=(\xi_{a,0})_a$. The maps
\[
 e_\pm:\mathcal D_1\longrightarrow E,\quad
       (s,u,\xi)\longmapsto u\pm s\xi,
 \qquad
 p:(\mathcal D_1)_0\longrightarrow E,\quad (u,\xi)\longmapsto u
\]
give
\[
 \mathcal D_1|_{\Gm}\cong E\times E\times\Gm,
 \qquad (\mathcal D_1)_0\cong TE=E\times\A^n,
 \qquad \mathcal D_h=\mathcal D_1\times\A^{2n(h-1)}.
\]
Here $TE$ denotes the tangent bundle of the displayed open
subset $E\subset\A^n$. In particular, every
$\mathcal D_h\to\A^1_s$ is smooth.

\emph{Step 7: construct the finite models and comparisons for (iii)--(iv).}
Multiplying the first $\ell$ series $C_a$ by $H_{i_a}$ gives
the boundary-coordinate pullbacks; the remaining series are
the tangential-coordinate pullbacks. To lift the resulting
map to $\A^n$ through $\phi:V\to\A^n$,
we also record a lift modulo $q_0$; formal \'etaleness then
determines all higher lifts uniquely. The inverse limit of
the resulting finite models will be identified with
$\mathcal U$ in Step~8.
For $h\ge1$, put
\[
 Z_h=\Spec\C[s,z]/\bigl((z^2-s^2)^h\bigr).
\]
The map $Z_h\to T$ is finite locally free of rank $2h$.
For an affine complex scheme $W$, Weil restriction along
$Z_1\to T$ is characterized by
\[
 \operatorname{Res}_{Z_1/T}(W\times Z_1)(T')
       =\operatorname{Hom}_{\C}(Z_1\times_TT',W)
\]
for $T$-schemes $T'$. In the coordinate order fixed in (iii),
the expressions
\[
 \bigl(H_{i_a}(u_{a,0}+z\xi_{a,0})\bigr)_{a\le\ell},
 \qquad
 \bigl(u_{a,0}+z\xi_{a,0}\bigr)_{a>\ell}
 \quad\bmod (z^2-s^2)
\]
define the map from $\mathcal D_1$ in
\[
 \mathcal W_1=
 \mathcal D_1\mathop{\times}_{
   \operatorname{Res}_{Z_1/T}(\A^n\times Z_1)}
   \operatorname{Res}_{Z_1/T}(V\times Z_1),
 \qquad
 \mathcal W_h=\mathcal W_1\times_{\mathcal D_1}\mathcal D_h.
\]
The other map is induced by $\phi$; products with $Z_1$
are over $\C$. Expanding generators and relations in the basis
$1,z$ shows that these Weil restrictions are affine of finite
presentation over $T$. The map induced by $V\to\A^n$ is
\'etale: its nilpotent lifting problems, after base change by
$Z_1$, are the unique lifting problems for $V\to\A^n$.
Consequently $\mathcal W_h\to\mathcal D_h$ is \'etale and
$\mathcal W_h\to T$ is smooth.

For a complex algebra $R$ and $s\in R$, put
\[
 Z_{h,R}:=Z_h\times_T\Spec R
        =\Spec(\mathcal O_s/q_0^h).
\]
By the defining fiber product, an $R$-point of $\mathcal W_h$
over $s$ consists of the coefficient data of $\mathcal D_h(R)$,
written as $C_a^{(h)}\in\mathcal O_s/q_0^h$, together with a map
$\bar\gamma_1:Z_{1,R}\to V$ satisfying
\[
 \phi\circ\bar\gamma_1=\beta_h|_{Z_{1,R}},
\]
where
\[
 \beta_h:Z_{h,R}\longrightarrow\A^n,\qquad
 \beta_h=
 \bigl(H_{i_1}C_1^{(h)},\ldots,H_{i_\ell}C_\ell^{(h)},
       C_{\ell+1}^{(h)},\ldots,C_n^{(h)}\bigr).
\]
The immersion $Z_{1,R}\hookrightarrow Z_{h,R}$ is defined
by the nilpotent ideal $(q_0)/(q_0^h)$. Since $\phi$ is
formally \'etale, the diagram
\[
\begin{tikzcd}
 Z_{1,R}\arrow[r,"\bar\gamma_1"]\arrow[d,hook]
   &V\arrow[d,"\phi"]\\
 Z_{h,R}\arrow[r,"\beta_h"']
   \arrow[ur,dashed,"\bar\gamma_h"]
   &\A^n
\end{tikzcd}
\]
has a unique dashed lift, functorially in $R$.
Applied to the universal point of $\mathcal W_h$, it gives
the morphism
\[
 \bar\gamma_h:Z_h\times_T\mathcal W_h\longrightarrow V.
\]
For compatible $R$-points of $\mathcal W_h$ as $h$ varies,
the pulled-back maps satisfy
\[
 \bar\gamma_{h+1}|_{Z_{h,R}}=\bar\gamma_h.
\]
Because $V$ is affine and
$\mathcal O_s=\varprojlim_h\mathcal O_s/q_0^h$, these maps
determine a unique
\[
 \bar\gamma:\Spec\mathcal O_s\longrightarrow V,
 \qquad
 \bar\gamma|_{Z_{h,R}}=\bar\gamma_h.
\]
Its coordinate pullbacks are
$(H_{i_a}C_a)_{a\le\ell},(C_a)_{a>\ell}$, where
$C_a=\varprojlim_h C_a^{(h)}$ and the first $\ell$ series
are units. Conversely, such a normal-form map determines
all these truncations and the initial lift $\bar\gamma_1$.
This identifies the $R$-points of
$\varprojlim_h\mathcal W_h$ with the normal-form maps,
functorially in every complex algebra $R$.

Restriction to the punctured source, as in Step~5, therefore
defines a morphism
\[
 \Psi:\varprojlim_h\mathcal W_h\longrightarrow\mathcal Z_{P,Q}.
\]
We now express the open restriction $D(f)$ of Step~4 at
finite order. Since
\[
 \Psi^*f\in
 \varinjlim_h\Gamma(\mathcal W_h,\mathcal O_{\mathcal W_h}),
\]
choose an order $h_0$ and a function
$f_{h_0}\in\Gamma(\mathcal W_{h_0},\mathcal O_{\mathcal W_{h_0}})$
representing $\Psi^*f$. For $h\ge h_0$, define
\[
 \mathcal W_h^\circ
 :=\mathcal W_h\times_{\mathcal W_{h_0}}D(f_{h_0}),
 \qquad
 \varprojlim_{h\ge h_0}\mathcal W_h^\circ
 =\Psi^{-1}(D(f)).
\]
Restricting $\mathcal W_h\to\mathcal D_h$ gives $c_h$,
and the maps $d_{h',h}$ induce $\tau_{h',h}$. For
$h'\ge h\ge h_0$, the defining fiber products give
\[
 \mathcal W_{h'}^\circ
   =\mathcal W_h^\circ\times_{\mathcal D_h}\mathcal D_{h'},
 \qquad c_h\tau_{h',h}=d_{h',h}c_{h'}.
\]
The coefficients with $h\le r<h'$ are free on
$\mathcal D_{h'}\to\mathcal D_h$, giving exactly the product
presentations in (iii). The squares formed by truncation and
the open immersions $\mathcal W_h^\circ\hookrightarrow\mathcal W_h$
are Cartesian. Each $\mathcal W_h^\circ$ is smooth over $T$.

We first construct the comparisons on unrestricted overlap
models and then impose the open conditions defining
$\mathcal U$ and $\mathcal U'$. For two SNC charts $V,V'$,
write $\mathcal W_h(V)$ and $\mathcal W_h(V')$ for their
unrestricted models.
Let $T_1(V,V')\subset\mathcal W_1(V)$ be the open locus
where the initial lift $\bar\gamma_1$ takes values in $V\cap V'$.
It is open because its complement is the image of
$\bar\gamma_1^{-1}(V\setminus V')$ under the finite projection
$Z_1\times_T\mathcal W_1(V)\to\mathcal W_1(V)$.
Put
\[
 T_h(V,V'):=\mathcal W_h(V)\times_{\mathcal W_1(V)}T_1(V,V').
\]
On this open, every $\bar\gamma_h$ takes values in $V\cap V'$,
since $Z_{1,R}\hookrightarrow Z_{h,R}$ is a nilpotent
thickening. For a compatible family of points of the
$T_h(V,V')$, the associated full extension $\bar\gamma$ also
factors through $V\cap V'$, since $q_0$ belongs to the
Jacobson radical of $\mathcal O_s$.
Use the same order on $L$ for both charts, write $x_i'=a_ix_i$
with $a_i$ invertible on $V\cap V'$, and denote the new
tangential coordinates by $y_a'$ for $a>\ell$.
At order $h$, set
\[
 C_a^{\prime(h)}=(\bar\gamma_h^*a_{i_a})C_a^{(h)}
       \quad(a\le\ell),\qquad
 C_a^{\prime(h)}=\bar\gamma_h^*y_a'\quad(a>\ell),
\]
and use the same initial lift, now viewed as a map to $V'$.
These formulas define a morphism
\[
 \chi_h:T_h(V,V')\longrightarrow T_h(V',V).
\]
Indeed, they pull back regular functions along the universal
$\bar\gamma_h$ and take their coefficients in the finite free
algebra $\mathcal O_s/q_0^h$. For $a\le\ell$ one has
\[
 H_{i_a}C_a^{\prime(h)}
   =(\bar\gamma_h^*a_{i_a})H_{i_a}C_a^{(h)}
   =\bar\gamma_h^*x_{i_a}'.
\]
Thus formal \'etaleness identifies the lift reconstructed
from the primed data with $\bar\gamma_h$.
The inverse coordinate change recovers the original series
and initial lift, so $\chi_h$ is an isomorphism.
Reduction modulo $q_0^h$ commutes with these formulas;
successive coordinate changes pull back the same regular
functions along $\bar\gamma_h$. This proves respectively
the truncation and triple-composition identities in (iv),
without cancelling $H_i$ modulo $q_0^h$.

Let $h_0,h_0'$ be the cutoffs for $\mathcal U,\mathcal U'$.
For $h\ge\max(h_0,h_0')$, define
\[
 \begin{aligned}
 O_h&:=T_h(V,V')\cap\mathcal W_h^\circ(V)
       \cap\chi_h^{-1}\bigl(\mathcal W_h^\circ(V')\bigr),\\
 O_h'&:=\chi_h(O_h)\subset\mathcal W_h^\circ(V').
 \end{aligned}
\]
Here $\mathcal W_h^\circ(V)$ and $\mathcal W_h^\circ(V')$
are defined by the pullbacks of $f_{h_0}$ and $f'_{h_0'}$,
respectively. Use $\chi_h$ also for its restriction
$O_h\xrightarrow{\sim}O_h'$ in (iv). All these opens are
inverse images under truncation, so the truncation identity
holds on $O_{h'}$ for $h'\ge h$. For finitely many charts,
take the maximum of their cutoffs and intersect the
corresponding opens; the same formulas give the stated
inverse and triple-overlap identities on those restrictions.
Step~8 identifies their inverse limits with the overlaps
of the neighborhoods.

\emph{Step 8: identify the neighborhoods and finish (i), (iii), and (iv).}
Put
\[
 \mathcal W:=\varprojlim_h\mathcal W_h,\qquad
 \mathcal W^\times:=\mathcal W\times_T(T\setminus\{0\}).
\]
Thus $\mathcal W^\times$ is the open locus where the collision
parameter $s$ is invertible. Here the puncture is in the base
$T=\A^1_s$; the punctured formal source of Step~7 is obtained
by inverting $q_0$.
Since each $\mathcal W_h$ is affine and smooth over $T$,
multiplication by $s$ is injective on its coordinate ring.
Filtered direct limits preserve this injectivity, so
\[
 \begin{gathered}
 \Gamma(\mathcal W,\mathcal O_{\mathcal W})
   =\varinjlim_h\Gamma(\mathcal W_h,\mathcal O_{\mathcal W_h}),\\
 \Gamma(\mathcal W,\mathcal O_{\mathcal W})
   \hookrightarrow
 \Gamma(\mathcal W^\times,\mathcal O_{\mathcal W^\times})
   =\Gamma(\mathcal W,\mathcal O_{\mathcal W})[s^{-1}].
 \end{gathered}
\]
This injection, and its localizations, say precisely that
$\mathcal W^\times$ is schematically dense in $\mathcal W$.
By Step~5, the restriction of the morphism $\Psi$ of Step~7
factors as
\[
 \mathcal W^\times\longrightarrow
 \mathcal C^\times_{\vct,\wct}
 \hookrightarrow\mathcal Z_{P,Q}.
\]
Every local equation of
$\overline{\mathcal C^\times_{\vct,\wct}}$ therefore pulls
back to zero on $\mathcal W^\times$, and hence on
$\mathcal W$ by the displayed injectivity. Thus $\Psi$
factors through the schematic input closure.

For the open models of Step~7 this gives
\[
 \varprojlim_{h\ge h_0}\mathcal W_h^\circ
 =\Psi^{-1}(D(f))
 \xrightarrow[\sim]{\ \Psi\ }
 \overline{\mathcal C^\times_{\vct,\wct}}\cap D(f)
 =\mathcal U.
\]
Indeed, on every complex test algebra, the inverse sends a
point of $\mathcal U$ to the truncations of its extension
and unit quotients constructed in Step~4. Uniqueness of the
extension and of those quotients makes both compositions
the identity. This completes (iii).
For a test point of $\mathcal U\cap\mathcal U'$, its two
extensions to $V$ and $V'$ agree as maps to $M$: their
equalizer is closed because $M$ is separated and contains
the schematically dense open $\Spec\mathcal O_s[q_0^{-1}]$.
The density follows from the injectivity in Step~1, valid
also for nonreduced $R$. The common extension therefore
factors through $V\cap V'$. By Step~3, multiplication by
$H_i$ on the full ring $\mathcal O_s$ is injective; hence
$U_i'=(\bar\gamma^*a_i)U_i$. Thus its two sets of truncated
coefficients are related by $\chi_h$ and give a point of
every $O_h$. Conversely, a compatible point of the $O_h$
gives the same punctured-source map in both $\mathcal U$ and
$\mathcal U'$. These constructions are inverse on every
test algebra, identifying both limits in (iv) with
$\mathcal U\cap\mathcal U'$ and completing (iv).

It remains to prove the reverse inclusion in (i).
Let $\kappa$ be any extension field of $\C$ and let
$\alpha:\Spec\kappa\llbracket z \rrbracket\to M$ be an arc of contact
$\vct+\wct$. Choose an affine SNC chart $V$ at its center.
Since the closed point of $\Spec\kappa\llbracket z \rrbracket$ maps to $V$,
the whole arc $\alpha$ factors through $V$.
Construct the unrestricted models of Steps~6--7 in this
chart, using the coordinate notation of (iii). The series
associated to $\alpha$ are
\[
 \begin{aligned}
 C_a&:=z^{-(v_{i_a}+w_{i_a})}\alpha^*x_{i_a}
       \in\kappa\llbracket z \rrbracket^\times &&(1\le a\le\ell),\\
 C_a&:=\alpha^*y_a\in\kappa\llbracket z \rrbracket &&(\ell<a\le n).
 \end{aligned}
\]
At $s=0$, the identification $\mathcal O_0\cong\kappa\llbracket z \rrbracket$
sends $q_0$ to $z^2$. Hence
\[
 \begin{gathered}
 \mathcal O_0/q_0^h\mathcal O_0\cong\kappa[z]/(z^{2h}),\\
 \widetilde\alpha_h:=\bigl((C_a\bmod z^{2h})_{a=1}^n,
       \alpha|_{\Spec(\kappa[z]/(z^2))}\bigr)
       \in(\mathcal W_h)_0(\kappa).
 \end{gathered}
\]
These are points of the fiber products defining $\mathcal W_h$
because $H_i|_{s=0}=z^{v_i+w_i}$, so their coordinate maps
are the truncations of $\alpha$. They are compatible under
truncation and define
$\widetilde\alpha=(\widetilde\alpha_h)_h\in\mathcal W_0(\kappa)$,
where $\mathcal W_0=\mathcal W\times_T\{0\}$.
The factorization of $\Psi$ proved above gives
\[
 \Psi(\widetilde\alpha)=\alpha|_{\Spec\kappa(\!( z )\!)}
 \in\overline{\mathcal C^\times_{\vct,\wct}}
       \cap(\mathcal Z_{P,Q})_0.
\]
This loop has weight $w(\vct+\wct)=N$, so it belongs to
$\mathcal H_{\vct,\wct}$. This proves
$\mathscr P_{\vct+\wct}\subseteq\mathcal H_{\vct,\wct}$
on underlying points. Step~3 proves the opposite inclusion.
Both loci have their reduced induced structures in the output
stage, so
\[
             \mathcal H_{\vct,\wct}=\mathscr P_{\vct+\wct}
\]
as reduced loci. This completes (i) and the proof.
\end{proof}

\begin{lemma}[The leading product]
\label{lem:fusion-leading-product}
Let $\vct,\wct,P,Q,N,I,J,L$ be as in
Lemma~\ref{lem:fusion-leading-contact}.
Under the contact-layer identification \eqref{eq:loop-layer},
the map induced by $m_2$ on $H^*(\gr_F\mathbf B)$ is
\begin{equation}\label{lp:eq:fusion-leading-product}
 (\alpha t^{\vct})(\beta t^{\wct})
 =\bigl(r^*_{I,L}\alpha\smile
        r^*_{J,L}\beta\bigr)t^{\vct+\wct},
\end{equation}
for $\alpha\in H^*(S_I^\circ)$ and
$\beta\in H^*(S_J^\circ)$, with the grading shifts of
\eqref{eq:loop-layer}. The maps $r^*_{I,L},r^*_{J,L}$ are those
of Lemma~\ref{lem:face-restriction}; an empty common face gives
zero. In particular, the class in $H^0(B_0)$ represented by
$1\in H^0(X)$ acts as the identity on $H^*(\gr_F\mathbf B)$. On $B_0$ the
binary map is the ordinary cup-product morphism.
\end{lemma}
\begin{proof}
\emph{Step 1: describe the lower-pole union by a closed-cover diagram.}
We express the normalized complex of the lower-pole union as a
homotopy pushout of the complexes of its two closed pieces.
Use $\mathcal Z_{a,b}=\varnothing$ when either index is negative,
$F_{-1}\mathscr LX=\varnothing$, and $B_{-1}=0$. Put
\[
 \begin{gathered}
 U_a=F_a\mathscr LX\setminus F_{a-1}\mathscr LX\quad(a\ge0),\\
 \mathcal Z_-=\mathcal Z_{P-1,Q}\cup\mathcal Z_{P,Q-1},\qquad
 \mathcal Z_\cap=\mathcal Z_{P-1,Q}\cap\mathcal Z_{P,Q-1},\qquad
 \mathcal V=\mathcal Z_{P,Q}\setminus\mathcal Z_-.
 \end{gathered}
\]
Unions and intersections are schematic. Source factorization
and the pole-stage identifications give
\[
 \mathcal V^\times\simeq\Gm\times U_P\times U_Q,
 \qquad (\mathcal V_0)_{\mathrm{red}}\simeq(U_N)_{\mathrm{red}},
 \qquad ((\mathcal Z_-)_0)_{\mathrm{red}}
       \simeq(F_{N-1}\mathscr LX)_{\mathrm{red}}.
\]
All normalizing indices are restricted from $\mathcal Z_{P,Q}$;
reduction preserves the normalized analytic diagram by
\emph{Reduction of finite models} in
Subsection~\ref{subsec:volume-index-proofs}.
Include these closed--open diagrams in the compatible finite-factor
presentations of Lemma~\ref{lem:relative-fusion-factorization}.

On a chart in these finite-factor presentations, choose common finite factors
$Y_1,Y_2\subset Y$ for the restrictions of
$\mathcal Z_{P-1,Q},\mathcal Z_{P,Q-1}\subset\mathcal Z_{P,Q}$.
Put $Y_-=Y_1\cup Y_2$ and $Y_\cap=Y_1\cap Y_2$;
these present $\mathcal Z_-$ and $\mathcal Z_\cap$.
Write
\[
 i_a:Y_a\hookrightarrow Y_-\quad(a=1,2),\qquad
 i_\cap:Y_\cap\hookrightarrow Y_-.
\]
On $Y_-(\C)$, restriction gives an exact sequence of sheaves
\[
 0\longrightarrow\F_{Y_-}
 \longrightarrow(i_1)_*\F_{Y_1}\oplus(i_2)_*\F_{Y_2}
 \longrightarrow(i_\cap)_*\F_{Y_\cap}
 \longrightarrow0.
\]
The first map consists of the two restrictions; the second
is their difference on $Y_\cap$. Exactness can be checked on
stalks: on the intersection the maps are the diagonal
$\F\to\F^2$ and $(a,b)\mapsto a-b$, and away from the
intersection there is only one nonzero summand.
Verdier duality gives the cofiber sequence
\[
 (i_\cap)_*\omega_{Y_\cap}
 \xrightarrow{(\tau_1,-\tau_2)}
 (i_1)_*\omega_{Y_1}\oplus(i_2)_*\omega_{Y_2}
 \xrightarrow{\sigma_1+\sigma_2}\omega_{Y_-}.
\]
Here $\tau_a$ is the proper trace for
$Y_\cap\hookrightarrow Y_a$, pushed forward along $i_a$,
and $\sigma_a$ is the proper trace for $i_a$. Thus, in
$\mathcal D(Y_-(\C);\F)$,
\[
 \omega_{Y_-}\simeq
 (i_1)_*\omega_{Y_1}
 \mathop{\coprod}_{(i_\cap)_*\omega_{Y_\cap}}
 (i_2)_*\omega_{Y_2},
\]
where the maps defining the pushout are $\tau_1,\tau_2$.
On each open-and-closed piece where the inherited ambient
normalizing index is constant, say $d$, apply the same shift
$[-2d]$ to every term. Exact periodization gives the parity
version. Proposition~\ref{lem:relative-compact-support}(i)
and the normalized purity equivalences descend this entire
diagram, including its proper traces and cofiber structure.

\emph{Step 2: identify the associated-graded morphism by localization.}
We identify the associated-graded product with specialization
on the open complement $\mathcal V$ of the lower-pole union.
Over $s\ne0$, source factorization and inward-ray insertion
identify the corresponding input complex of $\mathcal Z_-$ with
\[
 M_-=(B_{P-1}\otimes B_Q)
       \mathop{\coprod}_{B_{P-1}\otimes B_{Q-1}}
             (B_P\otimes B_{Q-1}),
\]
where $\coprod$ denotes the pushout in the derived enhancement.
Exactness of tensor product gives
\[
 \begin{aligned}
 \operatorname{cofib}(M_-\longrightarrow B_P\otimes B_Q)
 &\simeq
 \operatorname{cofib}(B_{P-1}\longrightarrow B_P)
 \otimes\operatorname{cofib}(B_{Q-1}\longrightarrow B_Q)\\
 &\simeq\gr_P\mathbf B\otimes\gr_Q\mathbf B,
 \end{aligned}
\]
where the last equivalence is that of
Theorem~\ref{thm:canonical-complex}. Closed/open localization
and specialization therefore give the coherent diagram
\[
 \begin{tikzcd}[column sep=large]
 M_- \arrow[r] \arrow[d]
   & B_P\otimes B_Q \arrow[r]
       \arrow[d,"\overline m_{2;P,Q}"]
   & \gr_P\mathbf B\otimes\gr_Q\mathbf B
       \arrow[d,"\gr m_{2;P,Q}"]\\
 B_{N-1} \arrow[r]
   & B_N \arrow[r]
   & \gr_N\mathbf B .
 \end{tikzcd}
\]
Both rows are cofiber sequences. Here $\gr m_{2;P,Q}$ denotes
the associated-graded morphism of the filtered product.
By Definition~\ref{lp:def:binary-fusion-product} and
Proposition~\ref{prop:relative-binary-identification}, the middle
arrow is normalized specialization. Proper-trace naturality
in Lemma~\ref{lem:specialization-arrow-descent} identifies the
right arrow with specialization on $\mathcal V$.

This comparison includes the homotopies defining the cofiber
maps. Indeed, the closed inclusion of the doubly lowered family
in each of the two singly lowered families gives the following
specified diagram over $\A^1_s$:
\[
\begin{tikzcd}[column sep=large]
 \mathcal Z_{P-1,Q}\arrow[d,"\mathrm{id}"']
   &\mathcal Z_{P-1,Q-1}\arrow[l]\arrow[r]\arrow[d]
   &\mathcal Z_{P,Q-1}\arrow[d,"\mathrm{id}"]\\
 \mathcal Z_{P-1,Q}
   &\mathcal Z_\cap\arrow[l]\arrow[r]
   &\mathcal Z_{P,Q-1}.
\end{tikzcd}
\]
Every horizontal arrow is a closed inclusion. The middle
vertical arrow factors through the schematic intersection
$\mathcal Z_\cap$, and is an isomorphism over $s\ne0$.
Thus, over $s\ne0$, source factorization and inward-ray
insertion identify the proper-trace diagrams of both rows
with the diagram whose pushout is $M_-$. If $P,Q\ge1$, proper traces on the
closed fibers give the diagram
\[
\begin{tikzcd}[column sep=large]
 B_{N-1}\arrow[d,"\mathrm{id}"']
   &B_{N-2}\arrow[l]\arrow[r]\arrow[d]
   &B_{N-1}\arrow[d,"\mathrm{id}"]\\
 B_{N-1}
   &B_{N-1}\arrow[l,"\mathrm{id}"']\arrow[r,"\mathrm{id}"]
   &B_{N-1}.
\end{tikzcd}
\]
The three arrows with source $B_{N-2}$ are the pole-transition
maps. The pushout of the bottom row is canonically $B_{N-1}$;
the induced map from the pushout of the top row to $B_{N-1}$
restricts to the identity on each outer term. Composing this
map with specialization of the top row gives the left vertical
arrow $M_-\to B_{N-1}$ in the preceding cofiber diagram.
Specializing the bottom row gives the same arrow, with its
specified homotopy, by naturality of this coherent diagram.
If exactly one of $P,Q$ is zero, both middle families are
empty and at most one outer family is nonempty, so the
same comparison reduces to that outer family's specialization.
If $P=Q=0$, all the lower terms vanish. These comparisons
identify the maps and their homotopies before taking cofibers;
no identification of $\mathcal Z_\cap$ with
$\mathcal Z_{P-1,Q-1}$ over $s=0$ is used.

\emph{Step 3: restrict specialization to the chosen pair of contact summands.}
We isolate the $(\vct,\wct)$ input summand and prove the
empty-common-face assertion of the lemma.
Let $\mathcal T_{\vct,\wct}\subset\mathcal V$ be the reduced
schematic closure of $\mathcal C^\times_{\vct,\wct}$.
We check the finite-presentation hypotheses needed for its trace.
By Proposition~\ref{prop:loop-sectors},
$\mathcal C^\times_{\vct,\wct}$ is open and closed in
$(\mathcal V^\times)_{\mathrm{red}}$.
Choose a nice chart $\Spec A\to\mathcal V$, that is, a
finitely presented \'etale chart with a finite-factor presentation
\[
 A=R_0[y_\eta:\eta\in\Lambda],
\]
where $R_0$ is of finite type over $\C$, $s\in R_0$, and
$\Lambda$ is finite or countable. The pullback of $\mathcal C^\times_{\vct,\wct}$ to
$\Spec A[s^{-1}]$ is the reduced closed subscheme supported
on an open-and-closed subset of $\Spec A[s^{-1}]$.
Let $e\in A[s^{-1}]$ be its idempotent: $e^2=e$, and
$D(e)$ is precisely that subset, so $e$ restricts to $1$
on it and to $0$ on its complement. The
argument in \emph{Idempotents in polynomial algebras} in
Subsection~\ref{subsec:relative-finite-factor-proof},
applied to $A[s^{-1}]=R_0[s^{-1}][y_\eta:\eta\in\Lambda]$,
gives $e\in R_0[s^{-1}]$. Define
\[
 J_0=(1-e)+\sqrt{(0)}\subset R_0[s^{-1}],
\]
where $\sqrt{(0)}$ is the nilradical of $R_0[s^{-1}]$.
The quotient $R_0[s^{-1}]/J_0$ is the reduced ring of the
open-and-closed part selected by $e$; hence $J_0$ is
radical. The pullback of the chosen input locus
$\mathcal C^\times_{\vct,\wct}$ to this chart over $s\ne0$
is therefore
\[
 \Spec\bigl(A[s^{-1}]/J_0A[s^{-1}]\bigr).
\]
Let $\lambda_0:R_0\to R_0[s^{-1}]$ and
$\lambda:A\to A[s^{-1}]$ be the localization maps, and put
\[
 I_0=\lambda_0^{-1}(J_0)
     =\ker\bigl(R_0\longrightarrow R_0[s^{-1}]/J_0\bigr).
\]
Coefficientwise localization gives
\[
 \lambda^{-1}\bigl(J_0A[s^{-1}]\bigr)=I_0A.
\]
Indeed, a polynomial belongs to the left-hand ideal exactly
when the image under $\lambda_0$ of each of its coefficients
belongs to $J_0$, equivalently when each coefficient belongs
to $I_0$. Every polynomial has finitely many terms, so this
argument also applies when $\Lambda$ is countable.
This inverse-image formulation does not require
$s$ to be a non-zero-divisor.
It identifies the schematic closure on this chart with
\[
 \Spec(A/I_0A)=\Spec(R_0/I_0)\times\A^\Lambda.
\]
The ideal $I_0$ is radical, so this closure is already reduced,
and $I_0$ is finitely generated because $R_0$ is Noetherian.
For a flat affine comparison $A\to A'$, tensoring the
kernel sequence with $A'$ gives
\[
 I_0A'=\ker\bigl(A'\longrightarrow
                   A'[s^{-1}]/J_0A'[s^{-1}]\bigr).
\]
This is the closure ideal of the pulled-back chosen input
locus. Smooth comparisons also preserve reducedness.
Thus the reduced closures agree on common chart refinements
and give compatible finite-factor presentations of
$\mathcal T_{\vct,\wct}\subset\mathcal V$.
The open $\mathcal V$ is quasi-compact, being the complement
of a finitely presented closed subscheme in the quasi-compact
scheme $\mathcal Z_{P,Q}$. Thus $\mathcal T_{\vct,\wct}$ is
quasi-compact, separated, locally nice, and finitely presented
as a closed subscheme of $\mathcal V$.

By construction and Lemma~\ref{lem:fusion-leading-contact}(i),
\[
 \mathcal T_{\vct,\wct}^{\times}
    =\mathcal C^\times_{\vct,\wct},\qquad
 ((\mathcal T_{\vct,\wct})_0)_{\mathrm{red}}
    =\mathscr P_{\vct+\wct}.
\]
The proper traces for $\mathcal T_{\vct,\wct}\hookrightarrow\mathcal V$
are the inclusions of the chosen input and output contact summands.
Consequently specialization on $\mathcal T_{\vct,\wct}$ computes
the restriction of $\gr m_{2;P,Q}$ to the $(\vct,\wct)$ input
summand, and no other output summand of weight $N$ occurs.
If $D_L^\circ=\varnothing$, then
$\mathscr P_{\vct+\wct}=\varnothing$ and this specialization
has zero target. This proves the empty-common-face assertion.
For the remaining steps assume $D_L^\circ\ne\varnothing$.

\emph{Step 4: compute the local specialization as a diagonal Gysin map.}
We compute specialization on $\mathcal D_h$ and show that it
sends the input complex-oriented fundamental class to the output one.
Fix an output point and use the neighborhood $\mathcal U$,
SNC chart $V$, finite factors $\mathcal W_h^\circ$, and
coefficient notation of Lemma~\ref{lem:fusion-leading-contact}(ii)--(iii).
Their maps $c_h:\mathcal W_h^\circ\to\mathcal D_h$ are
\'etale, and $\mathcal U\simeq\varprojlim_{h\ge h_0}\mathcal W_h^\circ$.
Put $n=\dim M$ and $E=(\Gm)^{|L|}\times\A^{n-|L|}$.
The coefficient models satisfy
\[
 \mathcal D_h=\mathcal D_1\times\A^{2n(h-1)},\qquad
 \mathcal D_h|_{s\ne0}
 \simeq\Gm\times E\times E\times\A^{2n(h-1)}.
\]
Put $u=(u_{a,0})_{a=1}^n$ and $\xi=(\xi_{a,0})_{a=1}^n$.
The two $E$-coordinates are $e_+=u+s\xi$ and $e_-=u-s\xi$.
The input and output factors have the affine projections
\[
 \begin{aligned}
 p_{\mathrm{in},h}&:
 E\times E\times\A^{2n(h-1)}\longrightarrow E\times E,\\
 p_{\mathrm{out},h}&:
 (\mathcal D_h)_0=E\times\A^{n(2h-1)}\longrightarrow E.
 \end{aligned}
\]
The first projection refers to the factor displayed over
$s\ne0$; it forgets the coefficients of $q_0^r$ for $r\ge1$.
The second forgets these coefficients and $\xi$.
Their relative complex dimensions are $2n(h-1)$ and
$n(2h-1)$, respectively.

Let $\operatorname{sp}_h=\partial_{s=0}(\rho\boxtimes-)$
be the Borel--Moore specialization of
\eqref{lp:eq:bm-sp} for $\mathcal D_h$, using
the above-displayed identification over $s\ne0$ and placing the
inward ray $\rho$ first, as in the fusion formula
\eqref{lp:eq:full-fusion-formula}.
With the complex-oriented smooth
pullbacks, the common affine factors give
\[
 (p_{\mathrm{out},h}^{!})^{-1}\operatorname{sp}_h
       p_{\mathrm{in},h}^{!}
   =(p^!)^{-1}\operatorname{sp}_1
   =\Delta_E^!,
 \qquad \Delta_E:E\hookrightarrow E\times E,
\]
as maps
$H_k^{\mathrm{BM}}(E\times E;\F)\to
H_{k-2n}^{\mathrm{BM}}(E;\F)$.
Here $p:(\mathcal D_1)_0=TE\to E$ is the tangent-bundle
projection. Put $x=e_+=u+s\xi$ and $\eta=-2\xi$; then
$e_-=x+s\eta$. These coordinates identify $\mathcal D_1$
with the deformation to the normal bundle of $\Delta_E$.
Specialization followed by the inverse Thom isomorphism is
therefore $\Delta_E^!$. The change $\xi\mapsto-2\xi$ is
complex linear and preserves the underlying real orientation.
In particular, writing $[Z]$ for the complex-oriented
Borel--Moore fundamental class of a smooth complex variety $Z$,
\[
 \operatorname{sp}_1([E\times E])=[TE]=p^![E],
 \qquad
 (p^!)^{-1}\operatorname{sp}_1([E\times E])=[E].
\]
Thus $[E\times E]\in H_{4n}^{\mathrm{BM}}(E\times E;\F)$
maps to $[E]\in H_{2n}^{\mathrm{BM}}(E;\F)$ with coefficient $1$.
The additional affine factors give the same assertion for
$\mathcal D_h$.

Open restriction and smooth base change pull this local
orientation calculation back along
$c_h:\mathcal W_h^\circ\to\mathcal D_h$.
These pullbacks are compatible with the truncation
projections $\tau_{h',h}$ and, after open restriction, with
the $T$-isomorphisms $\chi_h:O_h\xrightarrow{\sim}O_h'$ of
Lemma~\ref{lem:fusion-leading-contact}(iii)--(iv).
Consequently, on every such finite chart, specialization
sends the input fundamental class to the output fundamental class;
under purity, the corresponding local orientation generator
$1\in\F$ maps to $1\in\F$.
The displayed affine Gysin inverses use only the projections
$p_{\mathrm{in},h}$ and $p_{\mathrm{out},h}$ of $\mathcal D_h$.

\emph{Step 5: identify the global output orientation class and its degree.}
We prove the leading-product formula for $\alpha=\beta=1$
in both gradings.
For each contact vector $\uct$, let $\epsilon_{\uct}$ be the
class corresponding under \eqref{eq:loop-layer}, or its
periodic counterpart, to $1\in H^0(S_{I(\uct)}^{\log};\F)$
in the $\uct$-summand and to zero in all other summands.
Thus $\epsilon_{\uct}=1t^{\uct}$ has degree
$\mu_\Omega(\uct)$ in integer grading and parity $\bar0$;
in parity grading we use the image of the ordinary class $1$.
On smooth contact factors, purity identifies this class with
the normalized complex-oriented fundamental class.

Take the neighborhoods $\mathcal U$ of Step~4 at all
closed-fiber points, restrict them to $\mathcal V$, and let
$\mathcal N$ be the union of their complex-point spaces over
$s\ge0$, shrinking about the closed fiber as needed.
It is an open neighborhood of the entire closed fiber in
$\mathcal T_{\vct,\wct}(\C)|_{s\ge0}$. Put
\[
 \begin{gathered}
 \mathcal N_0:=\mathcal N\cap\{s=0\}
       =\mathscr P_{\vct+\wct}(\C),\qquad
 \mathcal N^+:=\mathcal N\cap\{s>0\},\\
 i:\mathcal N_0\hookrightarrow\mathcal N,\quad
 j:\mathcal N^+\hookrightarrow\mathcal N.
 \end{gathered}
\]
The equality for $\mathcal N_0$ follows from
Lemma~\ref{lem:fusion-leading-contact}(i); no uniform interval
in $s$ is required.

On each $\mathcal U$, take $Y=\mathcal W_h^\circ$ for
$h\ge h_0$, or a specified smooth common refinement, and put
$d_Y=\dim(Y/\A^1_s)$, read componentwise. For
$Y=\mathcal W_h^\circ$, one has $d_Y=2nh$.
Write $Y_0$ for its fiber over zero and $Y_+$ for its
complex-point space over the positive real ray.

In integer grading, set
\[
 \mu=\mu_\Omega(\vct+\wct)
     =\mu_\Omega(\vct)+\mu_\Omega(\wct).
\]
The contact-index formula \eqref{eq:direct-contact-volume-index}
and Remark~\ref{rem:relative-contact-volume-indices} give
\[
 \delta_{Y/\A^1}=d_Y+\frac{\mu}{2}
\]
on $s\ne0$. Since $Y\to\A^1_s$ is smooth, multiplication
by $s$ is injective on its affine coordinate rings, so this
open is schematically dense. Both sides are locally constant;
the equality therefore holds also along $Y_0$. Oriented purity gives
\[
 \omega_{Y_0}[-2\delta_{Y/\A^1}]
       \simeq\F_{Y_0}[-\mu],\qquad
 \omega_{Y_+}[-2\delta_{Y/\A^1}][-1]
       \simeq\F_{Y_+}[-\mu].
\]
Removing the common shift $[-\mu]$ gives the ordinary
dimension-normalized endpoint complexes
\[
 \omega_{Y_0}[-2d_Y]\simeq\F_{Y_0},\qquad
 \omega_{Y_+}[-2d_Y][-1]\simeq\F_{Y_+}.
\]
The shift $[-1]$ accounts for the real ray. For parity grading,
where no volume form is required, use these same ordinary
complexes as an auxiliary model before periodization.

For a smooth comparison $f:Y'\to Y$ of relative complex
dimension $e$, one has $d_{Y'}=f^*d_Y+e$. Oriented smooth
purity cancels this additional shift $[-2e]$; in integer
grading the volume indices obey the same refinement formula.
The truncation maps and overlap isomorphisms of
Lemma~\ref{lem:fusion-leading-contact}(iii)--(iv), together
with the coherent purity and localization comparisons of
Lemma~\ref{lem:specialization-arrow-descent}, therefore
descend the dimension-normalized boundaries to an ordinary arrow
\[
 \sigma_{\mathcal N}^{\mathrm{dim}}:
 Rj_*\F_{\mathcal N^+}\longrightarrow i_*\F_{\mathcal N_0}
 \quad\text{in }\mathcal D(\mathcal N;\F).
\]
The source identification uses the inward orientation of the
real ray in $\omega_{Y_+}[-2d_Y][-1]\simeq\F_{Y_+}$.

We identify this arrow globally. At a point of $\mathcal N_0$,
the smooth finite models and their affine tails give a basis
of neighborhoods of the form
\[
 B\times[0,\varepsilon)\times W,
\]
where $B$ is a ball in finite-dimensional real coordinates
on a smooth complex fiber, $\varepsilon>0$, and $W$ is a convex
cylinder in a finite or countable product of copies of $\C$.
Their intersections with $\mathcal N^+$ replace
$[0,\varepsilon)$ by $(0,\varepsilon)$. These intersections
are weakly contractible and locally weakly contractible.
The product argument of Lemma~\ref{lem:infinite-product-descent}(i)
therefore proves that the adjunction units
\[
 \begin{aligned}
 \eta_{\F}:\F_{\mathcal N}
   &\xrightarrow{\sim}Rj_*\F_{\mathcal N^+},\\
 \eta_{\mathrm{per}}:\operatorname{Per}_2\F_{\mathcal N}
   &\xrightarrow{\sim}Rj_*\operatorname{Per}_2\F_{\mathcal N^+}
 \end{aligned}
\]
are equivalences; away from $\mathcal N_0$ this is immediate.
The set of morphisms between the following ordinary sheaf
complexes is
\[
 \pi_0\operatorname{Map}_{\mathcal D(\mathcal N;\F)}
       (\F_{\mathcal N},i_*\F_{\mathcal N_0})
       =H^0(\mathcal N_0;\F).
\]
Under this identification, $\sigma_{\mathcal N}^{\mathrm{dim}}
\eta_{\F}$ is a locally constant function. Step~4 shows that
its value at every point is $1$, so it is the restriction
morphism $\F_{\mathcal N}\to i_*\F_{\mathcal N_0}$.
Under exterior product, inward-ray insertion, and restriction
to $\mathcal N^+$, the input class
$\epsilon_{\vct}\otimes\epsilon_{\wct}$ is represented by
the constant section $1$ after removing the common shift in
integer grading, and by its periodic image in parity grading.
The ordinary dimension-normalized output is therefore
$1\in H^0(\mathcal N_0;\F)$.

For the periodic operation, the canonical coefficient-change map
\[
 \theta:\operatorname{Per}_2(Rj_*\F_{\mathcal N^+})
       \longrightarrow Rj_*\operatorname{Per}_2\F_{\mathcal N^+}
\]
satisfies $\theta\operatorname{Per}_2(\eta_{\F})=\eta_{\mathrm{per}}$.
The two units are equivalences, hence so is $\theta$.
The coefficient-change compatibility proved in
Proposition~\ref{lem:relative-compact-support}, including
the coherent trace, purity, and localization comparisons,
identifies the periodic specialization arrow with
$\operatorname{Per}_2(\sigma_{\mathcal N}^{\mathrm{dim}})$
through $\theta$ and the canonical equivalence
$\operatorname{Per}_2(i_*\F_{\mathcal N_0})
\simeq i_*\operatorname{Per}_2\F_{\mathcal N_0}$.
These identifications descend as arrows by
Lemma~\ref{lem:infinite-product-descent}(ii).
Thus its output on the input orientation classes is the
image of the ordinary class $1\in H^0(\mathcal N_0;\F)$;
it has no additional positive even-degree component.
This argument does not interchange periodization with an
infinite \v Cech limit or with derived global sections on
$\mathcal N$.

In integer grading, restoring $[-\mu]$ places this output in
degree $\mu$, the sum of the input degrees. The contact-layer
identification and open-restriction naturality of specialization
now give, in both gradings,
\[
 (\gr m_{2;P,Q})_*
       (\epsilon_{\vct}\otimes\epsilon_{\wct})
       =\epsilon_{\vct+\wct}.
\]
This proves the leading-product formula for the two constant classes.

\emph{Step 6: extend the input cohomology classes across the collision.}
We extend the input classes across $\mathcal N$ and compute
their restrictions to $\mathcal N_0$.
Write $\ell:\mathcal N_0\to S_L^\circ$
for the map recording the output center and the phases of
its nonzero leading normal coefficients. 

For an input arc $\gamma_+\in\mathscr P_{\vct}(\C)$,
let
\[
 \operatorname{ev}_{\vct}(\gamma_+)\in S_I^\circ
\]
record its center $\gamma_+(0)\in D_I^\circ$ and the
phases of its leading normal coefficients. Explicitly,
in local boundary coordinates, write
\[
 x_i(\gamma_+(t_+))
   =a_{i,+}t_+^{v_i}+O(t_+^{v_i+1}),
 \qquad a_{i,+}\ne0,\quad i\in I;
\]
the corresponding normal-circle coordinate is
$a_{i,+}/|a_{i,+}|$. Define
$\operatorname{ev}_{\wct}:\mathscr P_{\wct}(\C)\to S_J^\circ$
in the same way.

For $s>0$, source factorization associates to each point
of $\mathcal N$ two input arcs $(\gamma_+,\gamma_-)$,
in the parameters $t_+=z-s$ and $t_-=z+s$, respectively.
The maps
\[
 (s,\gamma_+,\gamma_-)\longmapsto
 j_I\operatorname{ev}_{\vct}(\gamma_+),\qquad
 (s,\gamma_+,\gamma_-)\longmapsto
 j_J\operatorname{ev}_{\wct}(\gamma_-)
\]
extend continuously across $s=0$ to maps
\[
 F_+:\mathcal N\longrightarrow S_I^{\log},\qquad
 F_-:\mathcal N\longrightarrow S_J^{\log}.
\]
We verify this extension in the normal-form coordinates of
Lemma~\ref{lem:fusion-leading-contact}(ii)--(iii). For $i=i_a\in L$,
abbreviate $u_i=u_{a,0}$ and $\xi_i=\xi_{a,0}$; thus
$U_i\bmod q_0=u_i+z\xi_i$. Evaluation modulo $q_0$ at
$z=s$, respectively $z=-s$, is evaluation at the center of
the corresponding formal input arc. The factorization
\eqref{lp:eq:fusion-snc-factors} therefore gives
\[
 \begin{aligned}
 a_{i,+}&:=[t_+^{v_i}]x_i(\gamma_+(t_+))
       =(2s)^{w_i}(u_i+s\xi_i),\\
 a_{i,-}&:=[t_-^{w_i}]x_i(\gamma_-(t_-))
       =(-2s)^{v_i}(u_i-s\xi_i).
 \end{aligned}
\]
This also defines $a_{i,+}$ when $v_i=0$ and $a_{i,-}$ when
$w_i=0$: in those cases they are the boundary coordinates
of $\gamma_+(0)$ and $\gamma_-(0)$, respectively.
For $s>0$, their phases are
\[
 \frac{a_{i,+}}{|a_{i,+}|}
     =\frac{u_i+s\xi_i}{|u_i+s\xi_i|},\qquad
 \frac{a_{i,-}}{|a_{i,-}|}
     =(-1)^{v_i}\frac{u_i-s\xi_i}{|u_i-s\xi_i|}.
\]
These are continuous functions of the finite coefficients
$s,u_i,\xi_i$, including at $s=0$, because there $u_i\ne0$.
For $F_+$, the $i$th radial coordinate is identically zero
if $v_i>0$. If $v_i=0$, then $w_i>0$ since $i\in L$,
and that radial coordinate is $|a_{i,+}|$, which tends to
zero at the closed fiber. The same argument for $F_-$
interchanges $v_i$ and $w_i$. The tangential center
coordinates are the values at $z=\pm s$ of the tangential
series modulo $q_0$; both tend to their value at $z=0$.
Thus these formulas define continuous maps to the real
oriented blowup near each point of $\mathcal N_0$.

We next check that these local maps agree on overlaps.
For $x_i'=a_i x_i$ with $a_i$ invertible on the overlap,
the overlap formula in Lemma~\ref{lem:fusion-leading-contact}(iv)
gives on the separated locus
\[
 a_{i,+}'=a_i(\gamma_+(0))a_{i,+},\qquad
 a_{i,-}'=a_i(\gamma_-(0))a_{i,-}.
\]
At $s=0$, both multiplying factors become $a_i(\gamma_0(0))$,
where $\gamma_0$ is the output arc. Taking phases gives
exactly the transition functions of the normal-circle
bundles; for nonzero radial coordinates the same identities
give the coordinate changes in $M^{\log}$. The tangential
coordinates change by pullback under the same extension
$\bar\gamma$. Consequently the local maps glue to the
stated maps $F_+$ and $F_-$ on $\mathcal N$.

The map $\ell$ is the output evaluation
$\operatorname{ev}_{\vct+\wct}$, defined as above using the
output parameter $z$. Since the leading coefficient of
$x_i(\gamma_0(z))$ is $u_i$, the preceding phase formulas give
\[
 F_+|_{\mathcal N_0}=\iota_{IL}j_L\ell,\qquad
 F_-|_{\mathcal N_0}=\iota_{JL}j_LR_{\vct}\ell.
\]
Here $j_L:S_L^\circ\hookrightarrow S_L^{\log}$ and
$\iota_{IL}:S_L^{\log}\hookrightarrow S_I^{\log}$,
$\iota_{JL}:S_L^{\log}\hookrightarrow S_J^{\log}$ are the
inclusions of Lemma~\ref{lem:face-restriction}, and
$R_{\vct}:S_L^\circ\to S_L^\circ$ multiplies the $i$th
normal-circle coordinate by $(-1)^{v_i}$ for every $i\in L$.
Multiplication by $\exp(\pi\sqrt{-1}\,\tau v_i)$,
$0\le\tau\le1$, defines an isotopy over $D_L^\circ$ from
the identity to $R_{\vct}$. Hence $R_{\vct}^*=\mathrm{id}$
on cohomology.

Now let $\alpha\in H^a(S_I^\circ;\F)$ and
$\beta\in H^b(S_J^\circ;\F)$, where $a,b$ are ordinary
cohomological degrees, and put $d=a+b$. Define
\[
 \widehat\alpha=(j_I^*)^{-1}\alpha,\qquad
 \widehat\beta=(j_J^*)^{-1}\beta,\qquad
 c=F_+^*\widehat\alpha\smile F_-^*\widehat\beta
       \in H^d(\mathcal N;\F).
\]
The inverse maps $ (j_I^*)^{-1}$ and $(j_J^*)^{-1}$ exist
by Lemma~\ref{lem:face-restriction}. Let
\[
 q:\mathcal N^+\longrightarrow
 \mathscr P_{\vct}(\C)\times\mathscr P_{\wct}(\C),
 \qquad (s,\gamma_+,\gamma_-)\longmapsto(\gamma_+,\gamma_-)
\]
be the restriction of the input projection under source
factorization. The definition of $F_\pm$ gives
\[
 c|_{\mathcal N^+}
   =q^*\bigl(\operatorname{ev}_{\vct}^*\alpha
                  \boxtimes\operatorname{ev}_{\wct}^*\beta\bigr),
\]
where $\boxtimes$ is the exterior product of cohomology
classes on the two contact spaces. On the closed fiber,
Lemma~\ref{lem:face-restriction} and $R_{\vct}^*=\mathrm{id}$ give
\[
 \begin{aligned}
 (F_+|_{\mathcal N_0})^*\widehat\alpha
   &=\ell^*j_L^*\iota_{IL}^*(j_I^*)^{-1}\alpha
     =\ell^*r_{I,L}^*\alpha,\\
 (F_-|_{\mathcal N_0})^*\widehat\beta
   &=\ell^*R_{\vct}^*j_L^*\iota_{JL}^*(j_J^*)^{-1}\beta
     =\ell^*r_{J,L}^*\beta.
 \end{aligned}
\]
Consequently
\begin{equation}\label{eq:leading-collar-class}
 c|_{\mathcal N_0}
   =\ell^*\bigl(r_{I,L}^*\alpha\smile r_{J,L}^*\beta\bigr).
\end{equation}

\emph{Step 7: use the cohomology action to prove the product formula.}
We combine the extending class of Step~6 with the cohomology
action on specialization to prove \eqref{lp:eq:fusion-leading-product}.
For a contact vector $\uct$, use the normalized pole-stage
complexes of Proposition~\ref{prop:normalized-pole-sheaves}
to define
\[
 \omega_{\uct}^{\mathrm n}:=
 \omega_{w(\uct)}|_{\mathscr P_{\uct}(\C)}.
\]
The restriction is ordinary sheaf inverse image along
\[
 \mathscr P_{\uct}(\C)
 \hookrightarrow F_{w(\uct)}\mathscr LX(\C).
\]
Thus $\omega_{\uct}^{\mathrm n}$ belongs to
$\mathcal D(\mathscr P_{\uct}(\C);\F)$ or
$\mathcal D_{\mathrm{per}}(\mathscr P_{\uct}(\C);\F)$,
in integer and parity grading, respectively.

The contact-layer comparison in the proof of
Theorem~\ref{thm:canonical-complex} identifies its derived
sections with the $\uct$-summand of \eqref{eq:loop-layer}.
Under this identification, we also write $\epsilon_{\uct}$
for the corresponding class in these derived sections.
The same contact-layer comparison identifies $\omega_{\uct}^{\mathrm n}$ with
$\underline\F_{\mathscr P_{\uct}}[-\mu_\Omega(\uct)]$
in integer grading, and with the constant periodic complex
in parity grading. Thus, for
$\eta\in H^*(S_{I(\uct)}^\circ;\F)$, the contact-layer notation means
\[
 \eta t^{\uct}
   =\operatorname{ev}_{\uct}^*\eta\frown\epsilon_{\uct}.
\]
Here $\operatorname{ev}_{\uct}$ is defined by the same
center-and-leading-phase construction as in Step~6, and
$\frown$ is the cohomology action of
Lemma~\ref{lem:specialization-module}.

Let $\omega_+^{\mathrm n}$ be the restriction to $\mathcal N^+$ of the
relative normalized complex over $(0,\infty)$, and let
$\omega_0^{\mathrm n}=\omega_{\vct+\wct}^{\mathrm n}$ on $\mathcal N_0$.
More explicitly, on each finite model $Y$ of Steps~4--5,
these complexes are represented in integer grading by
\[
 \omega_{Y_+}[-2\delta_{Y/\A^1}],\qquad
 \omega_{Y_0}[-2\delta_{Y/\A^1}],
\]
respectively, where $Y_+$ is the inverse image of the
positive real ray. In parity grading use
$\operatorname{Per}_2\omega_{Y_+}$ and
$\operatorname{Per}_2\omega_{Y_0}$, with the specified
purity comparisons. In particular, $\omega_+^{\mathrm n}$ includes the
dualizing factor of the real ray; its shift $[-1]$ is
applied separately below. Corollary~\ref{cor:normalized-specialization-module}
applies to $\mathcal T_{\vct,\wct}$ by Step~3 and gives
\[
 \mathcal K^{\mathrm n}_{\mathcal N,\rho}=Rj_*\omega_+^{\mathrm n}[-1],\qquad T_{\mathcal N}=i_*\omega_0^{\mathrm n},
 \qquad
 \sigma_{\mathcal N}^{\mathrm n}:\mathcal K^{\mathrm n}_{\mathcal N,\rho}\longrightarrow T_{\mathcal N}.
\]
The last arrow is the descended half-ray localization
boundary, shifted by $[-1]$; Step~5 identifies it with
$\sigma_{\mathcal N}^{\mathrm{dim}}[-\mu]$ in integer grading
and with $\operatorname{Per}_2(\sigma_{\mathcal N}^{\mathrm{dim}})$
in parity grading. These complexes and this
arrow lie in $\mathcal D(\mathcal N;\F)$, or in its periodic
counterpart. For a sheaf complex $A$ on $\mathcal N$, write
$m_{c,A}:A\to A[d]$ for the tensor-unit action of the class
$c$ constructed in Step~6, as in
Lemma~\ref{lem:specialization-module}; in parity grading use
$\operatorname{Per}_2(c)$. The same corollary gives
\begin{equation}\label{eq:leading-specialization-module}
 \sigma_{\mathcal N}^{\mathrm n}[d]\circ m_{c,\mathcal K^{\mathrm n}_{\mathcal N,\rho}}
   \simeq m_{c,T_{\mathcal N}}\circ\sigma_{\mathcal N}^{\mathrm n}
       :\mathcal K^{\mathrm n}_{\mathcal N,\rho}\longrightarrow T_{\mathcal N}[d].
\end{equation}
This identity is on $\mathcal N$ itself and respects the specified
comparisons; it requires no representative of $c$ on a finite factor.

We now define the input map
\[
 A_{\mathcal N}:
 R\Gamma(\mathscr P_{\vct},\omega_{\vct}^{\mathrm n})
 \otimes R\Gamma(\mathscr P_{\wct},\omega_{\wct}^{\mathrm n})
       \longrightarrow R\Gamma(\mathcal N,\mathcal K^{\mathrm n}_{\mathcal N,\rho})
\]
as follows. First take exterior product of derived sections
on $\mathscr P_{\vct}\times\mathscr P_{\wct}$. On the full
separated locus, source factorization gives
\[
 \mathcal T_{\vct,\wct}(\C)|_{s>0}
 \simeq(0,\infty)\times\mathscr P_{\vct}(\C)
                         \times\mathscr P_{\wct}(\C).
\]
Its normalized sheaf comparison identifies the relative
complex with
$\omega_{(0,\infty)}\boxtimes\omega_{\vct}^{\mathrm n}
\boxtimes\omega_{\wct}^{\mathrm n}$.
Insert the inward orientation class
$\rho=[(0,\infty)]\in H^{\mathrm{BM}}_1((0,\infty);\F)$
as the first factor, using the degree-zero morphism
\[
 \F\longrightarrow
 R\Gamma((0,\infty),\omega_{(0,\infty)})[-1]
\]
that it represents; in parity grading, periodize this
morphism. Finally restrict to $\mathcal N^+$
and use composition of derived direct images to identify
\[
 R\Gamma(\mathcal N^+,\omega_+^{\mathrm n}[-1])
 \simeq R\Gamma(\mathcal N,Rj_*\omega_+^{\mathrm n}[-1])
 =R\Gamma(\mathcal N,\mathcal K^{\mathrm n}_{\mathcal N,\rho}).
\]
These operations define $A_{\mathcal N}$.

Put $\xi=H^*(A_{\mathcal N})(\epsilon_{\vct}\otimes\epsilon_{\wct})$,
using the contact-summand interpretation of the two classes
specified above. Apply the ray-action calculation of
Lemma~\ref{lem:specialization-module}(iii) on the full
separated product, with the class
$\operatorname{ev}_{\vct}^*\alpha\boxtimes
\operatorname{ev}_{\wct}^*\beta$, and then restrict to
$\mathcal N^+$. Step~6 identifies its restriction with
$c|_{\mathcal N^+}$. Naturality of the action under restriction
and the projection formula for $Rj_*$ give the action on
$\mathcal K^{\mathrm n}_{\mathcal N,\rho}$, including its ray shift:
\[
 m_{c,\mathcal K^{\mathrm n}_{\mathcal N,\rho}}
   \simeq(-1)^d\bigl(Rj_*m_{c|_{\mathcal N^+},\omega_+^{\mathrm n}}\bigr)[-1].
\]
The sign $(-1)^d$ cancels the sign from moving the
degree-$d$ input class past $\rho$, whose cohomological
degree is $-1$. The integer normalization shifts are even,
so they introduce no further sign. Using the contact-layer
action identity above, we obtain
\[
 H^*(A_{\mathcal N})
     \bigl((\alpha t^{\vct})\otimes(\beta t^{\wct})\bigr)
       =c\frown\xi.
\]
Only the exterior class is used on the full separated
product; the class $c$ itself is used on $\mathcal N$.
Thus neither a product decomposition of $\mathcal N^+$
nor an extension of $c$ beyond $\mathcal N$ is required.

Set $s_{\mathcal N}=H^*R\Gamma(\mathcal N,\sigma_{\mathcal N}^{\mathrm n})$.
Since $\mathcal N$ contains the entire closed fiber,
\[
 R\Gamma(\mathcal N,T_{\mathcal N})
 \simeq R\Gamma(\mathcal N_0,\omega_0^{\mathrm n})
 \simeq R\Gamma(\mathscr P_{\vct+\wct},
                       \omega_{\vct+\wct}^{\mathrm n}).
\]
Thus restriction leaves the target unchanged; no equivalence
on source sections is needed. Specifically, open-restriction
naturality in Corollary~\ref{cor:normalized-specialization-module}
identifies specialization after restriction to $\mathcal N$
with specialization on $\mathcal T_{\vct,\wct}$, with this
common target. Steps~2--3 identify the latter with
$\gr m_{2;P,Q}$ on the $(\vct,\wct)$ input summand and
with output in the $(\vct+\wct)$ summand. Consequently
$s_{\mathcal N}\circ H^*(A_{\mathcal N})$ is exactly this
contact-summand map, and Step~5 gives
$s_{\mathcal N}(\xi)=\epsilon_{\vct+\wct}$.
Thus \eqref{eq:leading-specialization-module} and
\eqref{eq:leading-collar-class} yield
\[
 \begin{aligned}
 (\gr m_{2;P,Q})_*
   \bigl((\alpha t^{\vct})\otimes(\beta t^{\wct})\bigr)
   &=s_{\mathcal N}(c\frown\xi)\\
   &=(c|_{\mathcal N_0})\frown\epsilon_{\vct+\wct}\\
   &=\bigl(r_{I,L}^*\alpha\smile r_{J,L}^*\beta\bigr)t^{\vct+\wct}.
 \end{aligned}
\]
For the last equality, use \eqref{eq:leading-collar-class}
and the contact-layer action identity with
$\operatorname{ev}_{\vct+\wct}=\ell$.
In integer grading both sides have degree
$a+b+\mu_\Omega(\vct)+\mu_\Omega(\wct)
 =a+b+\mu_\Omega(\vct+\wct)$ by Step~5.
In parity grading the actions are their periodizations,
and the degree is $\overline{a+b}$. Additivity treats all
cohomology classes from the homogeneous ones considered
here. This proves \eqref{lp:eq:fusion-leading-product} in
both gradings; Step~3 treated the empty-common-face case.

\emph{Step 8: prove the zero-contact and unit assertions.}
We identify the binary morphism on $B_0$ with cup product and
prove the associated-graded unit assertion of the lemma.
When $\vct=\wct=\zero$, the support sets are empty,
$S_\varnothing^\circ=X$, and the evaluation maps record
only the arc centers. Since $B_{-1}=0$, the stage $B_0$
is its own associated-graded piece. Under
$B_0\simeq R\Gamma(X,\F_X)$, Step~7 therefore gives
\[
 H^*(\overline m_{2;0,0})(\alpha\otimes\beta)
       =\alpha\smile\beta.
\]
We check that this also identifies the binary morphism
itself in the derived category. Let
$\Delta_X:X\to X\times X$ be the diagonal. The ordinary
cup-product morphism is the composite
\[
 \begin{aligned}
 B_0\otimes B_0
 &\simeq R\Gamma(X,\F_X)\otimes R\Gamma(X,\F_X)\\
 &\xrightarrow{\ \boxtimes\ }
 R\Gamma(X\times X,\F_{X\times X})
 \xrightarrow{\ \Delta_X^*\ }R\Gamma(X,\F_X)
 \simeq B_0.
 \end{aligned}
\]
Here $\boxtimes$ is exterior product of derived sections.
By the K\"unneth isomorphism over $\F$, decomposable
classes span the cohomology of $B_0\otimes B_0$, so this
composite and $\overline m_{2;0,0}$ induce the same map on
all cohomology groups. Every complex of $\F$-vector spaces
splits as its cohomology with zero differential plus a
contractible complex: choose complements to the boundaries
inside the cycles and to the cycles in each cochain group.
Consequently two morphisms between such complexes with
the same maps on cohomology represent the same morphism
in the derived homotopy category. Hence the displayed
composite and $\overline m_{2;0,0}$ are homotopic in the
derived enhancement.

For parity coefficients, use
$B_0\simeq\operatorname{Per}_2R\Gamma(X,\F_X)$ and the
periodization of the displayed cup-product composite.
The same splitting argument in the two parities, together
with the periodic K\"unneth isomorphism, identifies this
morphism with $\overline m_{2;0,0}$. This proves the
zero-contact assertion of the lemma.

Finally, for every contact vector $\uct$ with support $J$,
the face-restriction map satisfies
$r_{\varnothing,J}^*(1)=1$. The product formula therefore gives
\[
 (1t^{\zero})(\alpha t^{\uct})
    =\alpha t^{\uct}
    =(\alpha t^{\uct})(1t^{\zero})
       \qquad\bigl(\alpha\in H^*(S_J^\circ;\F)\bigr).
\]
Thus the class represented by $1\in H^0(X)$ acts as the
identity on every contact summand, and hence on
$H^*(\gr_F\mathbf B)$, as asserted.
\end{proof}

\begin{proof}[Completion of the proof of Theorem~\ref{thm:geometric-fusion-ring}]
Proposition~\ref{prop:binary-associativity} gives associativity
and graded commutativity. It remains to prove the unit assertion
and identify the filtered and associated-graded products.
Work in $\mathcal C=\mathcal D(\F)$ in volume grading and in
$\mathcal C=\mathcal D_{\mathrm{per}}(\F)$ in parity grading.
Write $\mathbf1$ for the tensor unit of $\mathcal C$; it is
$\F$ in volume grading and $\operatorname{Per}_2(\F)$ in parity
grading.

\emph{Step 1: construct the endomorphisms given by multiplication by the candidate unit.}
Let $e_0:\mathbf1\to B_0$ represent $1\in H^0(X)$, and let
$e\in H^0(B)=LH^0(X)$ be its image under $B_0\to B$.
For every $P\ge0$, define endomorphisms of the derived pole
stage $B_P$ by the composites
\[
 \begin{aligned}
 L_{e,P}:B_P
 &\simeq\mathbf1\otimes B_P
 \xrightarrow{e_0\otimes\mathrm{id}_{B_P}}B_0\otimes B_P
 \xrightarrow{\overline m_{2;0,P}}B_P,\\
 R_{e,P}:B_P
 &\simeq B_P\otimes\mathbf1
 \xrightarrow{\mathrm{id}_{B_P}\otimes e_0}B_P\otimes B_0
 \xrightarrow{\overline m_{2;P,0}}B_P.
 \end{aligned}
\]
The coherent compatibility of the product with pole transitions
in Proposition~\ref{prop:relative-binary-identification} makes
these families endomorphisms $L_e,R_e$ of the entire diagram
$(B_P)_{P\ge0}$ in $\operatorname{Fun}(\Z_{\ge0},\mathcal C)$.
In the following identities of morphisms, we pass to the
homotopy category of $\mathcal C$ or of this diagram category,
as appropriate.

The zero-contact assertion of Lemma~\ref{lem:fusion-leading-product}
gives
\[
 \overline m_{2;0,0}\circ(e_0\otimes e_0)=e_0,
\]
using $\mathbf1\otimes\mathbf1\simeq\mathbf1$.
Associativity therefore gives, with the tensor-unit
identifications suppressed,
\[
 \begin{aligned}
 L_{e,P}^{\,2}
 &=\overline m_{2;0,P}\circ
   \bigl((\overline m_{2;0,0}\circ(e_0\otimes e_0))
          \otimes\mathrm{id}_{B_P}\bigr)\\
 &=L_{e,P}.
 \end{aligned}
\]
The same argument gives $R_{e,P}^{\,2}=R_{e,P}$.
Proposition~\ref{prop:binary-associativity} supplies these
associativity comparisons coherently in $P$, so
\[
 L_e^2=L_e,\qquad R_e^2=R_e
\]
in the homotopy category of the diagram category.

\emph{Step 2: prove by induction that these endomorphisms are equivalences at every stage.}
The unit assertion on the contact layers in
Lemma~\ref{lem:fusion-leading-product} gives
\[
 H^*(\gr_PL_e)=H^*(\gr_PR_e)
   =\mathrm{id}_{H^*(\gr_P\mathbf B)}\qquad(P\ge0).
\]
Here $\gr_PL_e$ and $\gr_PR_e$ are the induced endomorphisms
of $\operatorname{cofib}(B_{P-1}\to B_P)\simeq\gr_P\mathbf B$,
with $B_{-1}=0$. Cohomology detects equivalences in
$\mathcal C$ (in the two parities in the periodic case), so
both associated-graded endomorphisms are equivalences.

Extend $L_e$ and $R_e$ to $B_{-1}=0$ by the unique
endomorphism of the zero object, which is an equivalence.
The coherent pole compatibility gives, for every $P\ge0$,
a morphism of localization cofiber sequences
\[
 \begin{tikzcd}[column sep=large]
 B_{P-1}\arrow[r]\arrow[d,"{L_{e,P-1}}"']
   &B_P\arrow[r]\arrow[d,"{L_{e,P}}"]
   &\gr_P\mathbf B\arrow[r]\arrow[d,"{\gr_PL_e}"]
   &B_{P-1}[1]\arrow[d,"{L_{e,P-1}[1]}"]\\
 B_{P-1}\arrow[r]
   &B_P\arrow[r]
   &\gr_P\mathbf B\arrow[r]
   &B_{P-1}[1].
 \end{tikzcd}
\]
Assume inductively that $L_{e,P-1}$ is an equivalence.
The first, third, and fourth vertical maps are then equivalences.
Applying cohomology to the two cofiber sequences and using
the five lemma shows that $H^*(L_{e,P})$ is an isomorphism;
hence $L_{e,P}$ is an equivalence. This argument includes
$P=0$, starting from $B_{-1}=0$. The same diagram with $R_e$
in place of $L_e$ proves the assertion for $R_{e,P}$.
Thus both endomorphisms are equivalences at every pole stage.

\emph{Step 3: deduce the unit and finish the theorem.}
Equivalences of coherent diagrams are detected at each stage,
as recalled in Definition~\ref{def:filtered-derived-object}.
Consequently $L_e$ and $R_e$ are invertible in the homotopy
category of the diagram category. An invertible idempotent
is the identity: multiplying $L_e^2=L_e$ by $L_e^{-1}$ gives
$L_e=\mathrm{id}$, and likewise $R_e=\mathrm{id}$.
The homotopy-colimit functor carries these identities to
endomorphisms of $B\simeq\hocolim_PB_P$. Their induced maps
on $H^*(B)=LH^*(X)$ are left and right multiplication by $e$,
so
\[
 \mu(e,a)=a=\mu(a,e)\qquad(a\in LH^*(X)).
\]
This proves the unit assertion.

Corollary~\ref{prop:loop-fusion-cohomology} supplies the
filtration bound. Lemma~\ref{lem:fusion-leading-product}
identifies the associated-graded product and its unit with
the logarithmic product and its unit under
\eqref{eq:loop-layer}. This proves the remaining assertions
in both gradings.
\end{proof}

The image of $e_0$ may vanish in the cohomology of the pole
colimit; in that case the unital ring is the zero ring. 

\subsection{Compactification independence and valuation gradings}
\label{subsec:fusion-invariance}

Finite-factor presentation and cochain-choice independence for the fusion product are established in
Section~\ref{sec:fusion-construction}. We now compare different
compactifications and record  additional valuation gradings.

\begin{proposition}[Compactification independence]
\label{prop:fusion-invariance-gradings}
Fix $X$ and, in integer grading, $\Omega$.
For two polarized SNC compactifications of $X$, let
$(B^{(1)},\overline m_2^{(1)})$ and
$(B^{(2)},\overline m_2^{(2)})$ be their loop complexes and
fusion morphisms. The cofinal pole-stage comparisons induce
an equivalence $g:B^{(1)}\simeq B^{(2)}$ in
$\mathcal D_{\mathrm{per}}(\F)$ in parity grading, or in
$\mathcal D(\F)$ in integer grading, satisfying
\[
 g\circ\overline m_2^{(1)}
       =\overline m_2^{(2)}\circ(g\otimes g)
\]
in the corresponding homotopy category. The map $H^*(g)$ is
a unital graded algebra isomorphism. Thus the unfiltered fusion
algebra is independent of the polarized compactification;
the comparison need not identify the two contact filtrations.
\end{proposition}
\begin{proof}
Let $\mathbf g=(g_\ell)$ and $\mathbf h=(h_j)$ be the affine
generators from the two compactifications, with embedding
multiples $r_g,r_h$. Choose expressions
\[
 h_j=\Phi_j(\mathbf g),\qquad g_\ell=\Psi_\ell(\mathbf h),
 \qquad \deg\Phi_j\le d,\quad\deg\Psi_\ell\le e,
\]
where $d,e\ge1$, and integers $C,C'\ge1$ satisfying
$r_hC\ge dr_g$ and $r_gC'\ge er_h$.
Write $\mathcal Z^{\mathbf g}_{I,\mathbf P}$ and
$\mathcal Z^{\mathbf h}_{I,\mathbf P}$ for the corresponding
marked families. For a marked configuration $a=(a_i)_{i\in I}$
over a complex algebra $R$, use the source algebras
$\mathcal O_a,\mathcal K_a$ of \eqref{eq:kv-monic-source},
writing $q$ for $q_0$. Let $\gamma:\Spec\mathcal K_a\to X$, and put
$H_g=\prod_{i\in I}(z-a_i)^{r_gP_i}$.
Polynomial substitution gives
\[
 \bigl(H_gg_\ell(\gamma)\in\mathcal O_a\text{ for all }\ell\bigr)
 \ \Longrightarrow\ 
 H_g^dh_j(\gamma)\in\mathcal O_a\text{ for all }j.
\]
Since $H_g^d\mid\prod_i(z-a_i)^{r_hCP_i}$, this and the
reverse substitution define closed inclusions
\[
 \mathcal Z^{\mathbf g}_{I,\mathbf P}\hookrightarrow
       \mathcal Z^{\mathbf h}_{I,C\mathbf P},\qquad
 \mathcal Z^{\mathbf h}_{I,\mathbf P}\hookrightarrow
       \mathcal Z^{\mathbf g}_{I,C'\mathbf P}.
\]
They are finitely presented over arbitrary $R$, including
nonreduced $R$. For example, in the first target substitute
$g_\ell=\Psi_\ell(\mathbf h)$ and impose
\[
 [q^nz^k]\bigl(H_gg_\ell(\gamma)\bigr)=0
       \quad(n<0,\ 0\le k<|I|).
\]
Here $[q^nz^k]$ extracts the coefficient of $q^n$ in the
$z^k$-component of the expansion in the basis
$1,z,\ldots,z^{|I|-1}$ over $R(\!( q )\!)$.
The target bounds give a uniform lower bound on $n$; hence
only finitely many equations occur, each involving finitely
many coefficient variables. Their vanishing is precisely
the source condition \eqref{eq:kv-marked-poles}.

Every finite collection of stages from the two systems therefore
embeds as finitely presented closed loci in a common larger stage.
Include these loci, their products, and the specified base changes
in each finite collision diagram. Their finite equation lists
give closed finite-factor presentations with a common affine
tail on each ambient chart, by
Proposition~\ref{prop:direct-collision-diagrams} and
Subsection~\ref{subsec:volume-index-preparations}.
On each resulting closed factor, the omitted coordinates and
relative coefficient-residue map $\kappa_{\pi,V}$ of
\eqref{eq:direct-relative-residue-map} are inherited from the
ambient finite-factor projection $\pi$, for every affine test
scheme $V$ on the closed chart. The punctured-source map
and, in integer grading, $\Omega$ are unchanged.
Thus the index on such a closed factor is the restriction
of the ambient index by
Proposition~\ref{prop:direct-relative-volume-indices}(i).
Common refinements identify the normalized complexes and their
closed traces in these diagrams.

The inclusions preserve the punctured-source map and commute
with separation, enlargement of individual bounds, and collision:
for every colliding group $I_j\subseteq I$,
\[
 C\sum_{i\in I_j}P_i=\sum_{i\in I_j}CP_i.
\]
The descended maps satisfy these identities and their compositions: the
equalities can be checked after the finite-factor projections.
Lemma~\ref{lem:source-specialization-comparison}(iii)
then gives coherent product comparisons. Writing $B^{(1)}_P$
and $B^{(2)}_P$ for the respective pole stages, these include
\[
 \begin{tikzcd}[column sep=large]
 B^{(1)}_P\otimes B^{(1)}_Q
     \arrow[r,"\overline m^{(1)}_{2;P,Q}"] \arrow[d]
   & B^{(1)}_{P+Q} \arrow[d]\\
 B^{(2)}_{CP}\otimes B^{(2)}_{CQ}
     \arrow[r,"\overline m^{(2)}_{2;CP,CQ}"]
   & B^{(2)}_{C(P+Q)}.
 \end{tikzcd}
\]
The vertical maps are the closed traces of the displayed
inclusions. The reverse comparisons have the same property,
and their composites are the ordinary pole transitions
$B^{(i)}_P\to B^{(i)}_{CC'P}$ for $i=1,2$.
Cofinality and Lemma~\ref{lem:filtered-operations-colimit} give
$g$ and its product identity after taking the two input pole
homotopy colimits. At bound zero the comparison is the identity
on the intrinsic arc space, so it preserves $1\in H^0(X)$
and hence the fusion unit.
\end{proof}

\begin{proposition}[Valuation gradings from units]
\label{prop:fusion-unit-valuations}
For $u\in\mathcal O(X)^\times$ and $a\in\Z$, write
$LH^*(X)_{u=a}$ for the summand supported on
$\{\gamma:\ord_tu(\gamma)=a\}$.
This locus is open and closed by the local-constancy
argument in Remark~\ref{rem:volume-change}, which uses
only the invertibility of $u$ and applies in both gradings.
The fusion product satisfies
\[
 \mu\bigl(LH^i(X)_{u=a},LH^j(X)_{u=b}\bigr)
       \subseteq LH^{i+j}(X)_{u=a+b}
       \qquad(a,b\in\Z).
\]
This holds in parity grading and in the integer grading determined
by $\Omega$, and simultaneously for every finite list of units.
\end{proposition}
\begin{proof}
The one-point summands and their passage to the pole
homotopy colimit are those of Remark~\ref{rem:volume-change}.
To prove compatibility with fusion, we construct a locally
constant total valuation on the collision families.

Fix a nonempty finite set $I$ and pole bounds
$\mathbf P=(P_i)_{i\in I}\in\Z_{\ge0}^I$. Work on the
pole-bounded multipoint family
$\mathcal Z_{I,\mathbf P}\to\A^I$ of
Definition~\ref{def:kv-bounded-family}. An $R$-point, for a complex algebra
$R$, consists of a configuration $\mathbf a=(a_i)_{i\in I}$
and a punctured-source map $\gamma:\Spec\mathcal K_{\mathbf a}\to X$
satisfying \eqref{eq:kv-marked-poles}. Write $q$ for $q_0$ in
\eqref{eq:kv-monic-source}. Using the norm of
Subsection~\ref{subsec:direct-relative-volume-indices}, put
\[
 f_u:=\operatorname{Nm}_{\mathcal K_{\mathbf a}/R(\!( q )\!)}(\gamma^*u)
       \in R(\!( q )\!)^\times.
\]
Multiplicativity of the norm gives $f_{u^{-1}}=f_u^{-1}$.
These constructions commute with coefficientwise change
of the test algebra $R$.

Write $f_u=\sum_m f_{u,m}q^m$ and
$f_{u^{-1}}=\sum_m f_{u^{-1},m}q^m$. For
$\mathfrak p\in\Spec R$, define the norm valuation by
\[
 \nu_{f_u}(\mathfrak p)
   :=\ord_q(f_u\bmod\mathfrak p)
       \quad\text{in }\kappa(\mathfrak p)(\!( q )\!),
\]
where $\kappa(\mathfrak p)$ is the residue field.
The coefficient criterion in Remark~\ref{rem:volume-change},
applied in each residue field, gives
\[
 \nu_{f_u}^{-1}(c)=D(f_{u,c}f_{u^{-1},-c}),\qquad
 \Spec R=\coprod_{c\in\Z}D(f_{u,c}f_{u^{-1},-c}).
\]
Here $D(b)=\{\mathfrak p:b\notin\mathfrak p\}$ is the
principal open defined by $b\in R$. Only finitely many
terms are nonempty, by the Laurent lower bounds for
$f_u$ and $f_{u^{-1}}$; hence each is also closed.
The orders are always taken in residue fields, so this
argument does not require $R$ to be reduced.
On $\mathcal Z_{I,\mathbf P}$, polynomial expressions for
$u$ and $u^{-1}$ in the fixed affine generators give
uniform lower bounds for the powers of $q$ in these
determinants. Their relevant coefficients are regular
functions involving finitely many coefficient variables.
The residue-field valuations are compatible with change
of the test algebra and hence define a locally constant
function on $\mathcal Z_{I,\mathbf P}$,
again denoted $\nu_{f_u}$. Its finite open-and-closed
decomposition is algebraic. Its defining idempotents descend
to every finite factor by the
argument in \emph{Idempotents in polynomial algebras} in
Subsection~\ref{subsec:relative-finite-factor-proof}.

We identify this total valuation on each collision stratum.
Let $K/\C$ be a field extension and
$(\mathbf a,\gamma)\in\mathcal Z_{I,\mathbf P}(K)$.
Partition $I=\coprod_{j\in J}I_j$ by equality of marked
positions: $a_i=b_j$ for $i\in I_j$, where the $b_j\in K$
are pairwise distinct. Thus $\mathbf a$ lies in the collision
stratum $C_\pi$ of Definition~\ref{def:geometric-collision-data},
with $\pi:I\twoheadrightarrow J$ sending $I_j$ to $j$.
Put $w_j=z-b_j$. Completion and Chinese remainders give
\[
 \mathcal K_{\mathbf a}\simeq\prod_{j\in J}K(\!( w_j )\!).
\]
The loop $\gamma_j:\Spec K(\!( w_j )\!)\to X$ is defined by
the composite on coordinate rings
\[
 \mathcal O(X)\xrightarrow{\gamma^*}\mathcal K_{\mathbf a}
              \longrightarrow K(\!( w_j )\!).
\]
It belongs to $F_{P'_j}\mathscr LX(K)$, where
$P'_j=\sum_{i\in I_j}P_i$, by
Lemma~\ref{lem:kv-pole-compatibility}. Thus
\eqref{eq:kv-volume-norm-order}, applied to $b=\gamma^*u$, says
\[
 \nu_{f_u}(\mathbf a,\gamma)
      =\sum_{j\in J}\ord_{w_j}(\gamma_j^*u).
\]
Each distinct marked position contributes once, even when
$|I_j|>1$; the norm-order formula has no factor $|I_j|$.
It also shows that separation and collision identifications
preserve these valuation sectors: the sum is unchanged
when the same distinct local loops are regrouped. Pole
inclusions preserve them because they leave the
punctured-source map unchanged.

For the binary family $\mathcal Z_{P,Q}$, the marked
positions are $s,-s$. The preceding formula becomes
\[
 \nu_{f_u}|_{s\ne0}
       =\ord_{z-s}u(\gamma_1)+\ord_{z+s}u(\gamma_2),
 \qquad
 \nu_{f_u}|_{s=0}=\ord_z u(\gamma_0),
\]
where $\gamma_1,\gamma_2$ are the two input loops for
$s\ne0$ and $\gamma_0$ is the single output loop for $s=0$.
In particular, the input sector $(a,b)$ lies in the total
sector $c=a+b$, whose closed fiber has output valuation $c$.
The localization triangle decomposes over these
open-and-closed sectors. Its boundary, inward-ray insertion,
and normalized finite-factor comparisons preserve that
decomposition, so specialization has no component in an
output sector with valuation different from $a+b$.
Passing to the pole colimit proves the asserted grading
of the fusion product. For finitely many units, take the
common refinement of their finite sector decompositions;
the same argument preserves each component of the tuple
of valuations.
\end{proof}

\begin{proposition}[Change of volume form]
\label{prop:fusion-volume-change}
Let $\Omega$ be a nowhere-zero algebraic volume form on $X$,
let $u\in\mathcal O(X)^\times$, and put $\Omega'=u\Omega$.
The comparisons of Remark~\ref{rem:volume-change},
\[
 (B_{\Omega',P})_{u=a}\simeq(B_{\Omega,P})_{u=a}[-2a],
       \qquad P\ge0,\quad a\in\Z,
\]
preserve fusion after regrading. In particular, a nonzero
homogeneous class $\alpha$ in the sector $u=a$ satisfies
\begin{equation}\label{lp:eq:volume-change-new}
 \deg_{u\Omega}\alpha-\deg_\Omega\alpha=2a.
\end{equation}
Direct-sum periodization identifies these products with the
parity-graded fusion product.
\end{proposition}
\begin{proof}
Remark~\ref{rem:volume-change} already constructs the
displayed pole-stage comparisons and their degree shifts.
It remains to prove that they intertwine fusion.
Use the norm-valuation sectors and the local loops
$\gamma_j$ in parameters $w_j$ constructed in the proof
of Proposition~\ref{prop:fusion-unit-valuations}.
For a relative finite factor $Y\to T$ in the chosen finite-factor presentation, write
$\delta_{\Omega,Y/T}$ and $\delta_{\Omega',Y/T}$ for its two
relative volume indices. Restrict to the open-and-closed sector
$\nu_{f_u}=c$. We claim that
\[
 \delta_{\Omega',Y/T}-\delta_{\Omega,Y/T}=c.
\]
On a product factor over $C_\pi$, formula
\eqref{eq:direct-stratum-volume-index} expresses each
relative index as the sum of the absolute indices on the
individual loop factors. The absolute index difference
on the $j$th factor is $\ord_{w_j}(\gamma_j^*u)$ by
Remark~\ref{rem:volume-change}. Consequently
\[
 \delta_{\Omega',Y/T}-\delta_{\Omega,Y/T}
   =\sum_j\ord_{w_j}(\gamma_j^*u)=c.
\]
The last equality is the total-valuation calculation in
the preceding proof. Smooth refinements add the same
relative dimension to both indices, and the specified base
changes pull both indices back, by
Proposition~\ref{prop:direct-relative-volume-indices}(i).
The image in the configuration base of every complex
point lies on some $C_\pi$, so these comparisons prove
the claimed equality on all factors in the finite-factor presentation.

On this sector the two normalized finite-factor complexes
are therefore related by
\[
 \omega_Y[-2\delta_{\Omega',Y/T}]
   =\bigl(\omega_Y[-2\delta_{\Omega,Y/T}]\bigr)[-2c].
\]
Both sides use the same dualizing complex and the same
purity, proper-trace, and localization morphisms, with this
common additional shift. The comparisons consequently
intertwine the descended specialization arrows themselves.
On a binary collision, the preceding proposition gives
$c=a+b$ for input valuations $a,b$. Under the canonical
tensor-shift identification
\[
 \begin{aligned}
 &(B_{\Omega,P})_{u=a}[-2a]\otimes
             (B_{\Omega,Q})_{u=b}[-2b]\\
 &\hspace{12mm}\simeq
 \bigl((B_{\Omega,P})_{u=a}\otimes
             (B_{\Omega,Q})_{u=b}\bigr)[-2(a+b)],
 \end{aligned}
\]
the input and output shifts thus agree. This proves that
the pole-stage comparisons preserve fusion, and passage
to the pole homotopy colimit proves the asserted algebra
comparison after regrading. The degree formula is the one
already established in Remark~\ref{rem:volume-change}.
Finally, all shifts are even, and their identifications
under direct-sum periodization commute with tensor products.
The coefficient-change comparisons therefore identify
these products with the parity-graded fusion product.
\end{proof}

\section{The BV operator}\label{sec:bv}
\label{subsec:loop-bv}

BV operators enrich graded-commutative algebras with Gerstenhaber
brackets, obtained by measuring their failure to satisfy the Leibniz
rule. These are structures which occur in topological field theory and string
topology \cite{GetzlerBV,AbouzaidBV}. In symplectic cohomology, the BV operator
$\Delta:SH^k(X)\to SH^{k-1}(X)$ arises from rotating asymptotic
markers and, together with the pair-of-pants product, defines a
BV algebra. It also induces the first differential in the spectral
sequence for $S^1$-equivariant symplectic cohomology
\cite[Section~8(b)]{SeidelBiased}. For cotangent bundles, the Viterbo
isomorphism identifies it with rotation of free loops
\cite{AbouzaidBV}. These roles motivated our construction of the corresponding
operator on $LH^*(X)$. 

In this section we construct a degree-$-1$ source-rotation operator on the chosen cochain model $\mathbf B$ of the loop complex $B$ and
prove the BV identities for its action on the fusion algebra $LH$. The main result of this section is Theorem \ref{thm:loop-bv}. Together with the previous properties of the fusion product, it allows us to finish the proof of Theorem \ref{thm:summary-fusion} from the Introduction.

We follow the grading convention of
Subsection~\ref{subsec:relative-pole-recovery}: in integer grading,
$B_P$, $B$, $\mathbf B$, and $F_P\mathbf B$ stand for
$B_{\Omega,P}$, $B_\Omega$, $\mathbf B_\Omega$, and
$F_P\mathbf B_\Omega$, respectively, with $\Omega$ the chosen
nowhere-zero algebraic volume form. Without $\Omega$, these
symbols denote the parity-graded objects.

\subsection{The BV identities}

\begin{definition}[BV algebra]\label{def:loop-bv}
A \emph{Batalin--Vilkovisky algebra}, or \emph{BV algebra},
is a unital associative graded-commutative $\F$-algebra $H$
with a linear map $\Delta:H^j\to H^{j-1}$ such that
$\Delta^2=0$, $\Delta(1)=0$, and, for homogeneous $a,b,c$
with $p=|a|$, $q=|b|$,
\begin{equation}\label{eq:bv-seven-term}
\begin{split}
 \Delta(abc)={}&\Delta(ab)c+(-1)^p a\Delta(bc)
             +(-1)^{(p-1)q}b\Delta(ac)\\
 &-\Delta(a)bc-(-1)^p a\Delta(b)c
             -(-1)^{p+q}ab\Delta(c).
\end{split}
\end{equation}
Our degree-$-1$ bracket convention is
\begin{equation}\label{eq:general-rotation-bracket}
 \{a,b\}:=(-1)^{|a|-1}
       \bigl(\Delta(ab)-\Delta(a)b-(-1)^{|a|}a\Delta(b)\bigr).
\end{equation}
\end{definition}

\begin{lemma}[The bracket of a BV algebra]
\label{lem:algebraic-bv-bracket}
For the bracket \eqref{eq:general-rotation-bracket} on a BV algebra,
\[
 \begin{aligned}
 \{a,b\}&=-(-1)^{(|a|-1)(|b|-1)}\{b,a\},\\
 \{a,bc\}&=\{a,b\}c+(-1)^{(|a|-1)|b|}b\{a,c\},\\
 \sum_{\mathrm{cyclic}(a,b,c)}
 (-1)^{(|a|-1)(|c|-1)}\{a,\{b,c\}\}&=0.
 \end{aligned}
\]
The bracket has degree $-1$, that is,
$\{H^p,H^q\}\subseteq H^{p+q-1}$. It  makes $H$ a  \emph{Gerstenhaber algebra}, that is, an associative
graded-commutative algebra equipped with a bilinear bracket
of degree $-1$ satisfying the three identities above:
graded skew-symmetry, the graded Leibniz rule, and the
graded Jacobi identity, respectively.
\end{lemma}
\begin{proof}
Graded commutativity gives the first identity. Substitution of
\eqref{eq:general-rotation-bracket} into the second gives exactly
\eqref{eq:bv-seven-term}. The algebraic BV identity
\cite[Proposition~1.2]{GetzlerBV} then deduces Jacobi from
$\Delta^2=0$. Reversing the homological degrees in that reference
makes $\Delta$ have degree $-1$; our bracket is the negative of
its bracket. Neither change alters the displayed parity signs
or the Jacobi identity.
\end{proof}

\begin{example}[Volume divergence on polyvector fields]
\label{ex:volume-divergence-bv} Here is a classical example.
Let $Y$ be a smooth complex $n$-fold with a nowhere-zero
algebraic $n$-form $\Omega$. The graded algebra
\[
 \operatorname{PV}^p(Y):=\Gamma(Y,\textstyle\bigwedge^pT_Y),
 \qquad \operatorname{PV}^*(Y)=\bigoplus_{p=0}^n
                  \operatorname{PV}^p(Y),
\]
has wedge product and degree $p$ on $p$-vector fields. Define
\[
 I_p:\bigwedge^pT_Y\xrightarrow{\sim}\Omega_Y^{n-p},
 \qquad
 I_p(V_1\wedge\cdots\wedge V_p)
       =\iota_{V_p}\cdots\iota_{V_1}\Omega,
\]
and, with $d$ the algebraic de Rham differential, put
\[
 \left.\Delta_\Omega\right|_{\operatorname{PV}^p}
   =(-1)^{p-1}I_{p-1}^{-1}dI_p\quad(p\ge1),
 \qquad \left.\Delta_\Omega\right|_{\operatorname{PV}^0}=0.
\]
In \'etale coordinates $x_1,\ldots,x_n$ with
$\Omega=g\,dx_1\wedge\cdots\wedge dx_n$, write
$\partial_i=\partial/\partial x_i$. Then
\[
 \Delta_\Omega(f\partial_{i_1}\wedge\cdots\wedge\partial_{i_p})
  =\sum_{j=1}^p(-1)^{j-1}g^{-1}\partial_{i_j}(gf)
       \partial_{i_1}\wedge\cdots\wedge
       \widehat{\partial_{i_j}}\wedge\cdots\wedge\partial_{i_p}.
\]
Thus, if $\theta_i$ denotes the odd, that is, degree-one, generator $\partial_i$,
\[
 \Delta_\Omega=g^{-1}\sum_i\partial_{x_i}\,g\,
                         \partial_{\theta_i},
\]
where $\partial_{\theta_i}$ is left odd differentiation, that is, contraction with \(dx_i\).
The operators $g^{-1}\partial_{x_i}g$ commute and the odd
operators anticommute, proving $\Delta_\Omega^2=0$.
Each summand is a composite of two graded derivations plus a
first-order term. Expanding their Leibniz rules gives
\eqref{eq:bv-seven-term}; also $\Delta_\Omega(1)=0$.
Consequently $\operatorname{PV}^*(Y)$ is a BV algebra.

For an algebraic vector field \(V\), the algebraic Cartan identity \(\mathcal L_V=d\iota_V+\iota_Vd\) gives $\mathcal L_V\Omega=(\Delta_\Omega V)\Omega$; thus
$\Delta_\Omega V$ is its divergence. For a regular function $f$,
\[
 \Delta_\Omega(fV)=f\Delta_\Omega V+V(f),
 \qquad \{V,f\}=V(f)
\]
with precisely the convention \eqref{eq:general-rotation-bracket}.
\end{example}

\subsection{Source rotation on the pole-stage diagram}
Write $B_\bullet=(B_P)_{P\ge0}$ for the original derived
pole-stage diagram and $B=\hocolim_{P\ge0}B_P$ for its colimit.
Let $\mathbf B$ be the fixed filtered cochain model of
Theorem~\ref{thm:canonical-complex}, with its specified equivalence
\[
 (F_P\mathbf B)_{P\ge0}\simeq (B_P)_{P\ge0}
\]
of entire diagrams. As explained there, this is an equivalence
in the functor $\infty$-category: it includes the comparisons of
all transition morphisms and their higher compatibilities.
For $P\le Q$, write
\[
 \jmath_{P,Q}:F_P\mathbf B\hookrightarrow F_Q\mathbf B
\]
for the actual subcomplex inclusion; its derived image corresponds
to the original proper-trace morphism $B_P\to B_Q$. We write
$d_{F_P\mathbf B}:=d_{\mathbf B}|_{F_P\mathbf B}$ for the
restricted differential.
In integer grading, $B$, $B_P$, and $\mathbf B$ abbreviate
$B_\Omega$, $B_{\Omega,P}$, and $\mathbf B_\Omega$, respectively,
for the chosen form $\Omega$.

A \emph{coherent circle action} on $B_\bullet$ is a
functor
\[
 BS^1\longrightarrow\operatorname{Fun}(N(\N),\mathcal D(\F))
\]
with value $B_\bullet$ at the base point, cf.
\cite[Definition~3.1]{MathewTorus}.  In parity grading use the category of two-periodic complexes. Here $BS^1$ is the classifying homotopy type of the
topological circle. The adjective \emph{coherent} emphasizes that the action
is an $\infty$-functor, including the composition and unit
homotopies and all their higher compatibilities.
Explicitly,  let
$\operatorname{Aut}_{\mathcal D(\F)}(B_\bullet)$ be the space of
self-equivalences of the entire pole-stage diagram.
Such a coherent  action functor specifies a map of spaces
\[
 \rho:S^1\longrightarrow
       \operatorname{Aut}_{\mathcal D(\F)}(B_\bullet),
 \qquad \lambda\longmapsto\rho_\lambda,
\]
together with homotopies
\[
 \rho_1\simeq\mathrm{id}_{B_\bullet},
 \qquad
 \rho_\lambda\circ\rho_\mu\simeq\rho_{\lambda\mu},
\]
and all higher unit and associativity compatibilities.
For example, for three circle elements, the two successive
uses of the multiplication homotopy are themselves related
by a specified homotopy. These data are parametrized by
the spaces $(S^1)^r$ for all $r\ge0$, so they include the dependence on paths
and higher homotopies in the circle parameters.
All equivalences and homotopies take place in $\mathcal D(\F)$,
and hence respect the pole transitions coherently.

\begin{proposition}[Coherent source rotation]
\label{prop:coherent-source-rotation}
For $P\ge0$, let $\gamma:\Spec\C(\!( t )\!)\to X$ be a formal
loop belonging to $F_P\mathscr LX(\C)$. The source rotations
$r_\lambda\gamma(t)=\gamma(\lambda t)$, $\lambda\in G:=\C^*$,
preserve $F_P\mathscr LX$ and induce a coherent circle action on
$B_\bullet$. This action is independent of the finite-factor
presentations of the pole stages $F_P\mathscr LX$, under their
specified refinement and filtered comparisons, with
the prescribed normalizations. It admits parameter-chain
representatives in the sense of
Definition~\ref{def:parameter-chain-representative}, on the fixed
filtered cochain model:
\[
 \begin{aligned}
 \Theta_P^G&:C_{-*}(G;\F)\otimes F_P\mathbf B
                  \longrightarrow F_P\mathbf B,\\
 \Theta_P&:=\Theta_P^G(\iota_*\otimes\mathrm{id}):
       C_{-*}(S^1;\F)\otimes F_P\mathbf B
                  \longrightarrow F_P\mathbf B,
       \qquad \iota:S^1\hookrightarrow G.
 \end{aligned}
\]
They may be chosen as the restrictions of filtered cochain maps
$\Theta^G:C_{-*}(G;\F)\otimes\mathbf B\to\mathbf B$ and
$\Theta:=\Theta^G(\iota_*\otimes\mathrm{id})$, together with
representatives of the entire action diagram and its comparisons.
\end{proposition}
\begin{proof}
Source substitution defines an automorphism over $\Gm$:
\[
 A_P:\Gm\times F_P\mathscr LX\xrightarrow{\sim}
                  \Gm\times F_P\mathscr LX,
 \qquad(\lambda,\gamma)\longmapsto(\lambda,r_\lambda\gamma).
\]
Indeed, for every $\C$-algebra $R$, $\lambda\in R^\times$,
and $\gamma\in F_P\mathscr LX(R)$, it sends
\[
 \gamma(g)(t)=\sum_n a_nt^n
 \quad\longmapsto\quad
 (r_\lambda\gamma)(g)(t)=\sum_n\lambda^na_nt^n
 \qquad(g\in\mathcal O(X)).
\]
Thus it preserves the pole conditions, has inverse given by
$\lambda^{-1}$, and commutes with every pole inclusion.
We also write $A_P$ for the induced map over $G=\C^*$ on
complex-point spaces.

Apply Lemma~\ref{lem:source-coordinate-factorization}(ii)--(iii)
in the one-factor case $I=J=\{1\}$, with pole bound $P$,
center $b=0$, initial coordinate multiplier $u=1$, and variable
multiplier $\lambda$. Its automorphism $r_\lambda$ is
precisely $A_P$. The finite-factor construction in these parts
transports charts along $A_P$ and compares them with the
original charts on common refinements over $G$: the maps from
all charts to $G$ factor through their finite factors, and
the comparison maps commute with those maps to $G$. Source rotation preserves the
relative volume indices, and the normalization cancels the
smooth-purity shifts of the refinements, by
Proposition~\ref{prop:direct-relative-volume-indices}.
The resulting compact-support comparison is
\[
 (r_\lambda)_!:
 \underline{B_P^\vee}_G\longrightarrow
 \underline{B_P^\vee}_G,
\]
with the notation and variance of
Lemma~\ref{lem:source-coordinate-factorization}(iii).
Define its Verdier dual on $G$ by
\[
 \mathcal R_P:=\mathbb D_G\bigl((r_\lambda)_!\bigr):
 \omega_G\otimes B_P\longrightarrow\omega_G\otimes B_P.
\]
This is a morphism in $\mathcal D(G;\F)$, or its periodic
version. No action on an individual chosen finite factor
is required.

For a fixed $\lambda\in G$, let
$i_\lambda:\{\lambda\}\hookrightarrow G$ be the inclusion.
Under the product comparisons,
\[
 i_\lambda^!\mathcal R_P
   =r_\lambda^!:B_P\longrightarrow B_P.
\]
Its compact-support dual is pullback by $r_{\lambda^{-1}}$.
On smooth finite factors, Poincar\'e duality identifies
$r_\lambda^!$ with ordinary cohomological pullback
$r_\lambda^*$, with the prescribed normalization shifts.

The geometric identities $r_\lambda r_\eta=r_{\lambda\eta}$
and $r_1=\mathrm{id}$ are transported with their composition
and unit comparisons by part~(iii) of that lemma. Explicitly,
for multiplication $m:G^2\to G$, $m(\lambda,\eta)=\lambda\eta$,
and projections $\operatorname{pr}_1,\operatorname{pr}_2:G^2\to G$,
\[
 \operatorname{pr}_1^!\mathcal R_P\circ
 \operatorname{pr}_2^!\mathcal R_P\simeq m^!\mathcal R_P
\]
as endomorphisms of $\omega_{G^2}\otimes B_P$.
Part~(iii) supplies the unit comparison, all higher comparisons
over $G^r$, and their compatibility with pole transitions and
refinements of finite-factor presentations. Part~(iv) restricts these
diagrams to $(S^1)^r\hookrightarrow G^r$ by extraordinary
pullback. Its hypotheses hold: these parameter maps are
semialgebraic, and the pulled-back smooth finite-factor
comparisons are locally oriented product projections, as in
Lemma~\ref{lem:topological-pullback-presentations}(ii).
This gives the coherent circle action on the entire diagram
$B_\bullet$ and its independence under the specified comparisons.

Each $B_P$ is perfect by
Corollary~\ref{cor:relative-pole-recovery}(ii). Apply
Definition~\ref{def:parameter-chain-representative} to
$\mathcal R_P$ and to the whole diagram of its compositions
and pole comparisons. The adjunction
\eqref{eq:parameter-chain-adjunction} gives
\[
 C_{-*}(G;\F)\otimes B_P\longrightarrow B_P.
\]
By Lemma~\ref{lem:parameter-chain-operations}(i), restriction
to $S^1$ is precomposition with $\iota_*\otimes\mathrm{id}$;
parts~(ii)--(iii) preserve composition and the higher comparisons.
Use Lemma~\ref{lem:filtered-operations-colimit} with the fixed
factor $C_{-*}(G;\F)$ and its circle restriction, and transport
the entire diagram through
$(F_P\mathbf B)_P\simeq(B_P)_P$.
It supplies the filtered cochain maps $\Theta^G$, $\Theta$,
their restrictions $\Theta_P^G$, $\Theta_P$, and representatives
of all action comparisons on the fixed model $\mathbf B$.
The perfectness used here concerns the derived objects $B_P$;
it imposes no termwise finiteness on $F_P\mathbf B$ or
$\mathbf B$, and no perfectness on the colimit $B$.
\end{proof}

\begin{definition}[The BV operator]
\label{def:source-rotation-operator}
For a counterclockwise fundamental cycle
$\sigma\in Z_1(S^1;\F)$, put
\[
 \Delta:=\Theta(\sigma,-):\mathbf B\longrightarrow\mathbf B[-1],
 \qquad
 \Delta_P:=\Delta|_{F_P\mathbf B}
       =\Theta_P(\sigma,-):F_P\mathbf B\longrightarrow(F_P\mathbf B)[-1].
\]
Their derived images, transported through the specified comparisons,
are denoted
\[
 \overline\Delta_P:B_P\longrightarrow B_P[-1],\qquad
 \overline\Delta:B\longrightarrow B[-1].
\]
We call $\Delta$ the \emph{BV operator} and use the same notation
for its induced cohomology operator. 
Theorem~\ref{thm:loop-bv} will prove that its induced operator on
$LH$ satisfies Definition~\ref{def:loop-bv}.
\end{definition}

\begin{lemma}[Basic properties of source rotation]
\label{lem:source-dilation}
The operator $\Delta$ preserves every pole stage and satisfies
$d_{\mathbf B}\Delta+\Delta d_{\mathbf B}=0$. Its filtered homotopy class is
independent of the fundamental cycle and the parameter-chain
representatives, under the specified comparisons of models.
Let $e\in(F_0\mathbf B)^0$ be a closed cochain whose class
corresponds to $1\in H^0(X;\F)$ under the specified arc-stage
comparison. There are filtered nullhomotopies of $\Delta^2$
and of $\Delta(e)$; in fact $\Delta_0\simeq0$ on all of
$F_0\mathbf B$. Consequently $\Delta^2=0$ and $\Delta(1)=0$
on $LH$.
\end{lemma}
\begin{proof}
For a homogeneous graded map $f$, use the Hom differential
$d_{\operatorname{Hom}}f=d_{\mathrm{target}}f
 -(-1)^{|f|}f d_{\mathrm{source}}$.
Lemma~\ref{lem:parameter-chain-operations}(i), applied to
$\Theta_P$, gives
\[
 d_{\operatorname{Hom}}\Theta_P(c,-)=\Theta_P(\partial c,-).
\]
Since $\partial\sigma=0$, this proves
$d_{F_P\mathbf B}\Delta_P+\Delta_Pd_{F_P\mathbf B}=0$;
filtration preservation follows from
Proposition~\ref{prop:coherent-source-rotation}. If
$\sigma'-\sigma=\partial b$ with $b\in C_2(S^1;\F)$,
then the degree-$-2$ graded endomorphism $\Theta_P(b,-)$ of
$F_P\mathbf B$ satisfies
\[
 d_{\operatorname{Hom}}\Theta_P(b,-)
       =\Theta_P(\sigma',-)-\Delta_P.
\]
Using the same $b$ for all $P$ gives a filtered homotopy.
For changes of parameter-chain representatives, the final
assertion of Lemma~\ref{lem:parameter-chain-operations}(iii)
includes the comparison homotopies of the whole action diagram;
Lemma~\ref{lem:filtered-operations-colimit} represents them
by filtered homotopies on the fixed cochain model.

Let $m:S^1\times S^1\to S^1$, $m(\lambda,\eta)=\lambda\eta$,
be circle multiplication. The Pontryagin chain product is
\[
 *=m_*\circ\operatorname{EZ}:
 C_{-*}(S^1;\F)\otimes C_{-*}(S^1;\F)
       \longrightarrow C_{-*}(S^1;\F),
\]
where $\operatorname{EZ}$ is the Eilenberg--Zilber chain cross
product; in particular, $C_p(S^1;\F)*C_q(S^1;\F)
\subseteq C_{p+q}(S^1;\F)$.
Apply Lemma~\ref{lem:parameter-chain-operations}(i) to
restriction along $m$, part~(ii) to composition of the action
morphisms, and part~(iii) to the two independent circle
parameters. The multiplication comparison of the coherent
action then gives a degree-$-1$ homotopy between the maps
\[
 \begin{aligned}
 C_{-*}(S^1;\F)^{\otimes2}\otimes F_P\mathbf B
     &\longrightarrow F_P\mathbf B,\\
 c\otimes d\otimes x
     &\longmapsto\Theta_P(c,\Theta_P(d,x)),\\
 c\otimes d\otimes x
     &\longmapsto\Theta_P(c*d,x).
 \end{aligned}
\]
Evaluate this homotopy on $\sigma\otimes\sigma$, with sign
chosen so that its Hom differential is the first map minus
the second. This gives a degree-$-3$ graded endomorphism
$K_P$ of $F_P\mathbf B$, satisfying
\[
 d_{\operatorname{Hom}}K_P
      =\Delta_P^2-\Theta_P(\sigma*\sigma,-).
\]
Since $\sigma*\sigma$ is a $2$-cycle and $H_2(S^1;\F)=0$,
choose $b\in C_3(S^1;\F)$ with
$\partial b=\sigma*\sigma$. Part~(i) now gives
\[
 H_P:=K_P+\Theta_P(b,-),\qquad |H_P|=-3,\qquad
 d_{F_P\mathbf B}H_P+H_Pd_{F_P\mathbf B}=\Delta_P^2.
\]
By the joint choice of comparison homotopies in
Proposition~\ref{prop:coherent-source-rotation}, the maps $K_P$
may be chosen as restrictions of one filtration-preserving
graded endomorphism of $\mathbf B$. The same choice of $b$
for all $P$ makes the $H_P$ restrictions of a filtered
nullhomotopy of $\Delta^2$.

On $F_0\mathscr LX=\mathscr L^+X$, center evaluation to
$X$, $\gamma\mapsto\gamma(0)$, and the inclusion of constant
arcs are equivariant for
the trivial action on $X$. The contraction
$\gamma(t)\mapsto\gamma(st)$, $s\in[0,1]$, of
Lemma~\ref{lp:lem:arc-contraction} commutes with rotation.
At each finite jet stage, evaluation and the normalized
trace to $X$ also commute with rotation. Hence the
comparison $B_0\simeq C^*(X;\F)$ of
Theorem~\ref{thm:canonical-complex}, with two-periodization
in parity grading, intertwines the coherent rotation action
with the trivial action, including the parameter comparisons.
Transport this comparison and the
trivial action jointly through $F_0\mathbf B\simeq B_0$.
On the chosen cochain representatives it gives a graded map
\[
 U:C_{-*}(S^1;\F)\otimes F_0\mathbf B\longrightarrow
       F_0\mathbf B,\qquad |U|=-1,
\]
satisfying
\[
 d_{\operatorname{Hom}}U
   =\Theta_0-\varepsilon\otimes\mathrm{id}_{F_0\mathbf B},
 \qquad \varepsilon:C_{-*}(S^1;\F)\to\F,
\]
where $\varepsilon$ is the singular-chain augmentation, sending
each point to $1$ and every positive-dimensional chain to $0$.
Since $\partial\sigma=0$ and $\varepsilon(\sigma)=0$, evaluation
on $\sigma$ gives the degree-$-2$ graded endomorphism
$L_0:=U(\sigma,-)$ of $F_0\mathbf B$ and the identity
\[
 d_{F_0\mathbf B}L_0-L_0d_{F_0\mathbf B}=\Delta_0.
\]
This proves $\Delta_0\simeq0$. Define the cochain map
$e_0:\F\to F_0\mathbf B$ by $e_0(1)=e$. Then $L_0e_0$
is a nullhomotopy of $\Delta_0e_0$, and
$\jmath_{0,P}L_0e_0$ is a nullhomotopy of
$\Delta_P\jmath_{0,P}e_0$ for every $P$. These maps are
compatible under the pole inclusions and give the asserted
filtered nullhomotopy of $\Delta(e)$.
\end{proof}

\begin{remark}[Fixed source dilations]
\label{rem:fixed-source-dilation}
For any $\lambda\in\C^*$, a path from $1$ to $\lambda$
and $\Theta_P^G$ give a degree-$-1$ homotopy from the
identity on $F_P\mathbf B$ to the chosen cochain representative
of $r_\lambda^!:B_P\to B_P$, compatible with the inclusions
$\jmath_{P,Q}$ after the specified endpoint and unit comparisons. Thus each fixed
source dilation is filtered homotopic to the identity.
\end{remark}

\subsection{Framed collision operations}
\label{subsec:framed-collision-operations}

\begin{definition}[Framed disc parameters]
\label{def:framed-disc-parameters}
For $k\ge1$, let
\[
 \mathcal D_k=\left\{(a_i,r_i)_{i=1}^k\in(\C\times\C^*)^k:
 \begin{array}{l}
 |a_i|+|r_i|<1,\\
 |a_i-a_j|>|r_i|+|r_j|\quad(i\ne j)
 \end{array}\right\}.
\]
A point specifies the ordered affine embeddings
$t\mapsto a_i+r_it$ of closed discs into the open unit disc. For each of these maps we call \(r_i\in\mathbb C^\times\) the {\it  multiplier}, its modulus \(|r_i|\) is the disc radius,  its argument records the framing, and, as usual, $a_i$ is called the {\it center}.
Insertion in slot $i$ is the map
\[
 \circ_i:\mathcal D_k\times\mathcal D_l\longrightarrow
                       \mathcal D_{k+l-1}
\]
which replaces $(a_i,r_i)$ by
$(a_i+r_ib_j,r_is_j)_{j=1}^l$ and leaves the other entries in
their original order. 
Permutations act by relabeling the entries.

For parameter chains, $c\circ_i d$ means
$(\circ_i)_*\operatorname{EZ}(c\otimes d)$, where
$\operatorname{EZ}:C_{-*}(\mathcal D_k;\F)\otimes
C_{-*}(\mathcal D_l;\F)\to C_{-*}(\mathcal D_k\times\mathcal D_l;\F)$
is the Eilenberg--Zilber cross product: triangulate
$\Delta^p\times\Delta^q$ by its $(p,q)$-shuffles and sum the
resulting singular $(p+q)$-simplices with their orientation signs.
It satisfies
$\partial\operatorname{EZ}(c\otimes d)
=\operatorname{EZ}(\partial c\otimes d)
 +(-1)^p\operatorname{EZ}(c\otimes\partial d)$ for $c\in C_p$, see \cite[pp.~277--278]{Hatcher}.
Only operations with at most three inputs will be used.
\end{definition}

For $k$ inputs with pole bounds $P_1,\ldots,P_k\ge0$,
permutations reorder the factors of
$\bigotimes_{i=1}^kF_{P_i}\mathbf B$ by the Koszul symmetry, permuting the pole bounds together with their factors.

\begin{proposition}[Framed collision operations]
\label{lem:framed-collision-operations}
The coherent relative system of
Lemma~\ref{lem:relative-fusion-factorization}(i)--(ii) and the source-coordinate
comparisons of Lemma~\ref{lem:source-coordinate-factorization}
give the following maps and compatibilities.
\begin{enumerate}[(i),wide=0pt]
\item For $1\le k\le3$, $\mathbf P=(P_1,\ldots,P_k)$, and
$N=\sum_iP_i$, there are degree-zero parameter-chain maps
\begin{equation}\label{eq:framed-collision-operations}
 \Theta_{k,\mathbf P}:C_{-*}(\mathcal D_k;\F)\otimes
                  \bigotimes_iF_{P_i}\mathbf B
                       \longrightarrow F_N\mathbf B,
 \qquad d_{\operatorname{Hom}}\Theta_{k,\mathbf P}(c)
                         =\Theta_{k,\mathbf P}(\partial c).
\end{equation}
They are equivariant for relabeling and compatible with pole maps:
if $\mathbf P\le\mathbf P'$ and $N'=\sum_iP'_i$, then
\[
 \jmath_{N,N'}\Theta_{k,\mathbf P}
 \simeq\Theta_{k,\mathbf P'}
       (\mathrm{id}\otimes\textstyle\bigotimes_i
                                      \jmath_{P_i,P'_i}).
\]
These are comparisons of the whole parameter-chain diagrams on
the fixed filtered model. The maps $\jmath_{P,Q}$ are actual
inclusions; the specified equivalence identifies their derived
images with the original proper-trace transitions.
\item Suppose $k+l-1\le3$, $\mathbf Q=(Q_1,\ldots,Q_l)$,
and $P_i=\sum_jQ_j$. Let $\mathbf R$ be obtained from
$\mathbf P$ by replacing $P_i$ by $\mathbf Q$.
The two cochain maps obtained by insertion of parameters and by
composition of operations are homotopic, coherently with repeated
insertions and pole maps. In particular, for parameter cycles $c,d$,
\[
 \Theta_{k+l-1,\mathbf R}(c\circ_i d)
 \simeq\Theta_{k,\mathbf P}(c)\circ_i\Theta_{l,\mathbf Q}(d),
\]
where insertion of homogeneous maps uses the Koszul rule.
\end{enumerate}
For each $k$, the maps in \textup{(i)} may be chosen as the
restrictions of one filtered cochain map
\[
 \Theta_k:C_{-*}(\mathcal D_k;\F)\otimes
                   \mathbf B^{\otimes k}\longrightarrow\mathbf B.
\]
Here the parameter complex has filtration concentrated in index
zero, and the tensor power has the sum filtration
\[
 F_N(\mathbf B^{\otimes k})
   :=\sum_{\sum_iP_i\le N}\bigotimes_iF_{P_i}\mathbf B.
\]
The comparison and insertion homotopies are likewise chosen as
restrictions of filtered homotopies. Thus their output on inputs
of bounds $P_1,\ldots,P_k$ lies in $F_{\sum_iP_i}\mathbf B$.
Homologous parameter cycles induce homotopic operations.
\end{proposition}

\begin{proof}
\emph{Step 1: construct the derived operations in (i).} The geometry behind (i) is the construction of a 1-parameter collision family associated to a point $(a,r)\in \mathcal D_k$; for $k=1$ is there is no collision, just a marked point moving along the source curve. A map from the source curve $\bA^1_z\to X$ and the centers $a=(a_i)_i$ on this curve give $k$ formal loops into $X$, and we can reparametrize the formal neighborhood of each center using the multiplier $r_i$ by setting locally $z=a_i+r_it_i$, with $t_i$ the new formal coordinate. Now we scale by $\varepsilon$, so $z=\varepsilon a_i+\varepsilon r_it_i$ and let $\varepsilon$ go to $0$. We give next the details how this collision family produces the maps from which (\ref{eq:framed-collision-operations}) are extracted.

Put $T=\mathcal D_k$ and $N=\sum_iP_i$. The objects $B_N$
and $\bigotimes_iB_{P_i}$ are perfect by
Corollary~\ref{cor:relative-pole-recovery}(i)--(ii).
Let $i_0:T\hookrightarrow[0,1]\times T$, $i_0(a,r)=(0,(a,r))$,
and $j_+:(0,1]\times T\hookrightarrow[0,1]\times T$ be the
zero-end and positive-locus inclusions. Put
\[
 h:[0,1]\times T\longrightarrow(\A^k)^{\mathrm{an}},
 \qquad h(\varepsilon,(a,r))=(\varepsilon a_i)_i,
 \qquad E=h^{-1}\mathcal E_{[k],\mathbf P}.
\]

For $\varepsilon>0$, the collision family pulled back along
$h$ has distinct marked centers $\varepsilon a_i$.
Apply Lemma~\ref{lem:source-coordinate-factorization}(i)--(ii)
with the singleton partition $\pi=\mathrm{id}_{[k]}$ and
the substitutions
\[
 b_i=\varepsilon a_i,\qquad u_i=\varepsilon r_i.
\]
The corresponding formal source coordinates are therefore
\[
 t_i=\frac{z-\varepsilon a_i}{\varepsilon r_i},
 \qquad z=\varepsilon a_i+\varepsilon r_it_i.
\]
More precisely, part~(iv) pulls back the source identification
$\Phi_{\pi,u}$ and compact-support comparison $e_{\pi,u}$
along
\[
 (0,1]\times T\longrightarrow T_\pi^{\mathrm{an}},
 \qquad
 (\varepsilon,(a,r))
 \longmapsto
 \bigl((\varepsilon a_i)_i,(\varepsilon r_i)_i\bigr),
\]
where $T_\pi$ is the parameter space of that lemma.
We denote the resulting compact-support comparison by
$e_{\pi,\varepsilon r}$, as in
\eqref{eq:source-coordinate-comparisons}.
The multipliers $\varepsilon r_i$ are nonzero on this
positive locus.

At $\varepsilon=0$, all marked centers coincide at $z=0$,
with total pole bound $N$. Apply the same parts of the lemma
to the one-part partition $[k]\twoheadrightarrow\{1\}$,
with center $b=0$ and multiplier $u=1$.
The resulting output coordinate is $z$.
Part~(iv) pulls this comparison back to the zero face along
the constant parameter map
\[
 T\longrightarrow(\A^1\times\Gm)^{\mathrm{an}},
 \qquad (a,r)\longmapsto(0,1).
\]

The parameter maps just displayed, as well as $h$, are
semialgebraic. Their pulled-back finite factors satisfy
Lemma~\ref{lem:topological-pullback-presentations}(ii):
they are semialgebraic, and the smooth comparison maps
remain locally oriented product projections.
Thus the topological pullbacks used above satisfy the
hypotheses of
Lemma~\ref{lem:source-coordinate-factorization}(iv).

Consequently, with
$\operatorname{pr}_T:(0,1]\times T\to T$ denoting projection,
these comparisons give
\begin{equation}\label{eq:framed-collar-identifications}
 i_0^{-1}E\simeq\underline{B_N^\vee}_T,
 \qquad
 j_+^{-1}E\simeq\operatorname{pr}_T^{-1}
       \underline{\left(\bigotimes_iB_{P_i}\right)^\vee}_T.
\end{equation}
Parts~(i)--(ii) of
Lemma~\ref{lem:topological-pullback-presentations} also identify
\[
 \mathbb D_{[0,1]\times T}E
 \simeq h^!\mathcal B_{[k],\mathbf P}
\]
with the normalized relative Borel--Moore complex on the
base of the collision family pulled back along $h$.

The adjunction unit $\eta_E:E\to Rj_{+*}j_+^{-1}E$ gives
\begin{equation}\label{eq:framed-collar-restriction}
 r_{k,\mathbf P}:\underline{B_N^\vee}_T
 \simeq i_0^{-1}E\xrightarrow{i_0^{-1}\eta_E}
 i_0^{-1}Rj_{+*}j_+^{-1}E
 \simeq\underline{\left(\bigotimes_iB_{P_i}\right)^\vee}_T.
\end{equation}
The last equivalence is the one induced by the positive-locus
trivialization in \eqref{eq:framed-collar-identifications}:
Lemma~\ref{lem:parameter-chain-properties}(i) identifies
$i_0^{-1}Rj_{+*}\operatorname{pr}_T^{-1}
\underline{(\bigotimes_iB_{P_i})^\vee}_T$ with
$\underline{(\bigotimes_iB_{P_i})^\vee}_T$.
Verdier duality on $T$ therefore defines
\[
 \alpha_{k,\mathbf P}:=\mathbb D_T(r_{k,\mathbf P}):
 \omega_T\otimes\bigotimes_iB_{P_i}
       \longrightarrow\omega_T\otimes B_N.
\]
Here perfectness identifies
$\mathbb D_T\underline{V^\vee}_T$ with $\omega_T\otimes V$
for each coefficient complex $V$ in this formula.

For $\mathbf P\le\mathbf P'$, the proper-trace transitions
$B_{P_i}\to B_{P'_i}$ and $B_N\to B_{N'}$, where
$N'=\sum_iP'_i$, give the coherent square
\[
 \begin{array}{ccc}
 \omega_T\otimes\bigotimes_iB_{P_i}
   &\xrightarrow{\alpha_{k,\mathbf P}}&\omega_T\otimes B_N\\
 \big\downarrow&&\big\downarrow\\
 \omega_T\otimes\bigotimes_iB_{P'_i}
   &\xrightarrow{\alpha_{k,\mathbf P'}}&\omega_T\otimes B_{N'}.
 \end{array}
\]
Indeed, Lemma~\ref{lem:relative-fusion-factorization}(i)--(ii)
gives the pole comparisons for the relative complexes, and
Lemma~\ref{lem:source-coordinate-factorization}(iii) makes the
source-coordinate comparisons compatible with those maps.
Naturality of $\eta_E$ then gives the square before duality
and hence the displayed square. The same argument applies to
relabeling, with the Koszul permutation of the input factors.
It includes the comparison homotopies of the whole diagrams.
This constructs the derived families required in (i).
Their cochain representatives will be chosen together in
Step~6, after constructing the insertion comparisons.
All steps also apply in parity grading, using the periodic
versions of the cited results.

\emph{Step 2: construct the binary insertion square for (ii).}
If $k,l>1$ and $k+l-1\le3$, then $k=l=2$. By relabeling, it suffices to insert into the first slot.
Put $T=\mathcal D_2\times\mathcal D_2$, and let
$p_1,p_2:T\to\mathcal D_2$ be the two projections. Let $((a_1,r_1),(a_2,r_2);(b_1,s_1),(b_2,s_2))\in T$.
 The geometry behind (ii), with this particular insertion of parameters, is the insertion of the 2-loops collision family from Step 1 for $(b_1,s_1),(b_2,s_2)\in \mathcal D_2$ into the the first loop of the 2-loops collision family from Step 1 for  $(a_1,r_1),(a_2,r_2)\in\mathcal D_2$. This gives us a collision of three loops. 
 
With this in mind, for $(\varepsilon,\eta)\in[0,1]^2$, we can then define directly the three centers by
the first column below. On $\varepsilon,\eta>0$, the second
column specifies the formal source coordinate $t=(z-c)/u$
at the center $c$, where $u$ is the multiplier in the same row:
\[
 \begin{array}{c|c}
 \text{center }c&\text{multiplier }u\\ \hline
 \varepsilon a_1+\varepsilon\eta r_1b_1
      &\varepsilon\eta r_1s_1\\
 \varepsilon a_1+\varepsilon\eta r_1b_2
      &\varepsilon\eta r_1s_2\\
 \varepsilon a_2&\varepsilon r_2.
 \end{array}
\]
Let $\chi:[0,1]^2\times T\to(\A^3)^{\mathrm{an}}$
be this center map, 

Write the inserted bounds as $(P,Q)$ and the bounds of the
configuration receiving the insertion as $(P+Q,R)$.
Put $N=P+Q+R$.  Put
$E_\chi=\chi^{-1}\mathcal E_{[3],(P,Q,R)}$.
When $\varepsilon,\eta>0$, the centers are distinct:
the two inserted discs are disjoint, and
\[
 |\eta r_1b_j|+|\eta r_1s_j|
   =\eta|r_1|(|b_j|+|s_j|)<|r_1|\qquad(j=1,2)
\]
places them inside the first disc of the receiving configuration,
separated from its second disc.

We now specify the source-coordinate comparison on each scale
stratum. Apply Lemma~\ref{lem:source-coordinate-factorization}(i)--(ii)
with the following partition and pairs (center, multiplier):
\[
 \begin{array}{c|c|c}
 \text{stratum}&\text{groups of labels}&
                  \text{center--multiplier pairs}\\ \hline
 \varepsilon=0&\{1,2,3\}&(0,1)\\
 \eta=0<\varepsilon&\{1,2\},\{3\}&
       (\varepsilon a_1,\varepsilon r_1),
       (\varepsilon a_2,\varepsilon r_2)\\
 \varepsilon,\eta>0&\{1\},\{2\},\{3\}&
       \text{the three rows above}.
 \end{array}
\]
The bound of each group is the sum of its input bounds.
Thus on $\eta=0<\varepsilon$ the formal coordinates are
$(z-\varepsilon a_1)/(\varepsilon r_1)$ and
$(z-\varepsilon a_2)/(\varepsilon r_2)$, with bounds $P+Q,R$.
On all of $\varepsilon=0$ the coordinate is $z$, with bound $N$.
The pairs in the table define the maps to the spaces $T_\pi$
of that lemma. Its part~(iv) pulls the comparisons $e_{\pi,u}$
back along these maps. The topological hypotheses hold by the
semialgebraic pullback argument of Step~1.
Consequently the ordinary inverse-image restrictions of
$E_\chi$ to the four scale strata are constant with values
\[
 \begin{array}{c|c}
 \varepsilon=\eta=0&B_N^\vee\\
 \varepsilon=0<\eta&B_N^\vee\\
 \eta=0<\varepsilon&(B_{P+Q}\otimes B_R)^\vee\\
 \varepsilon,\eta>0&(B_P\otimes B_Q\otimes B_R)^\vee.
 \end{array}
\]
The morphisms between these objects will be obtained from
restriction of sections of $E_\chi$, as follows.

\emph{Step 3: identify the four arrows and prove binary insertion.}
For $A\subseteq\{1,2\}$, let $v_A\in[0,1]^2$ have $j$th
coordinate $1$ for $j\in A$ and $0$ otherwise. Put
\[
 W_A=\{(\varepsilon,\eta)\in[0,1]^2:
       \varepsilon>0\text{ if }1\in A,\quad
       \eta>0\text{ if }2\in A\}.
\]
These are the open stars of
Lemma~\ref{lem:parameter-chain-properties}(ii), intersected
with $[0,1]^2$; the coordinates are ordered as
$(\varepsilon,\eta)$. Writing $\operatorname{pr}_T$ for
each projection to $T$, open-star evaluation gives
\[
 R(\operatorname{pr}_T)_*(E_\chi|_{W_A\times T})
       \xrightarrow{\sim}E_\chi|_{\{v_A\}\times T}.
\]
The right side is inverse-image restriction, viewed as a
complex on $T$. To see this equivalence over $T$, take a
contractible open $U\subset T$. The restrictions in Step~2
are locally constant on each product of a scale stratum
with $U$. The calculation of
Lemma~\ref{lem:parameter-chain-properties}(ii) therefore
computes $R\Gamma(W_A\times U,E_\chi)$ as the limit of the
stratum diagram indexed by
$\{A'\subseteq\{1,2\}:A\subseteq A'\}$.
This poset has initial object $A$, so evaluation induces
\[
 R\Gamma(W_A\times U,E_\chi)
       \xrightarrow{\sim}
 R\Gamma(\{v_A\}\times U,E_\chi).
\]
These evaluation maps commute with restriction in $U$ and
hence give the stated sheaf equivalences.

For $A\subseteq A'$, denote by $r_{A,A'}$ the restriction
\[
 R(\operatorname{pr}_T)_*(E_\chi|_{W_A\times T})
 \longrightarrow
 R(\operatorname{pr}_T)_*(E_\chi|_{W_{A'}\times T}),
\]
induced by $W_{A'}\times T\subseteq W_A\times T$, and also
its transport through the evaluation equivalences.
Part~(iii) of the same lemma gives the coherent square
\[
 \begin{array}{ccc}
 \underline{B_N^\vee}_T
   &\xrightarrow{r_{\varnothing,\{1\}}}&
       \underline{(B_{P+Q}\otimes B_R)^\vee}_T\\
 \big\downarrow\rlap{$\scriptstyle r_{\varnothing,\{2\}}$}&&
 \big\downarrow\rlap{$\scriptstyle r_{\{1\},\{1,2\}}$}\\
 \underline{B_N^\vee}_T
   &\xrightarrow{r_{\{2\},\{1,2\}}}&
       \underline{(B_P\otimes B_Q\otimes B_R)^\vee}_T
 \end{array}
\]
in $\mathcal D(T;\F)$, or its periodic version. Its homotopy
between the two composites is supplied by the functor of
open-star restrictions.

We identify its arrows with collar restrictions on the four
edges of the square. Restriction to either zero-coordinate
face gives the comparison of diagrams in
Lemma~\ref{lem:parameter-chain-properties}(iii).
For a face with one coordinate equal to $1$, restrict the
open stars for which that coordinate is required to be positive.
For example, restriction to $\varepsilon=1$ compares the
arrow for $\{1\}\subseteq\{1,2\}$ with the collar in $\eta$.
Both endpoints of this comparison are equivalences: evaluation
on the square and on the face gives the same restrictions
at $(1,0)$ and $(1,1)$. Naturality of evaluation gives a
comparison of the arrows themselves. The same argument on
$\eta=1$ applies to $\{2\}\subseteq\{1,2\}$.
Thus the comparison with positive-coordinate faces follows
from evaluation as well as the comparison with zero-coordinate
faces.

On $\eta=0$, the diagonal identification for the group
$\{1,2\}$ gives the binary collar with centers
$\varepsilon a_1,\varepsilon a_2$, multipliers
$\varepsilon r_1,\varepsilon r_2$, and bounds $P+Q,R$.
It is the pullback by $p_1$ of the collar in Step~1. Hence
\[
 r_{\varnothing,\{1\}}
       \simeq p_1^{-1}r_{2,(P+Q,R)}.
\]
This uses the diagonal comparison of
Lemma~\ref{lem:relative-fusion-factorization}(ii), its
coherence in part~(i), and the source-coordinate comparison
of Lemma~\ref{lem:source-coordinate-factorization}(ii)--(iii).

On $\varepsilon=1$, the first two points remain separated
as a group from the point $a_2$, including at $\eta=0$.
The separation comparison of
Lemma~\ref{lem:relative-fusion-factorization}(ii) splits off
the singleton of bound $R$. Its coordinate $(z-a_2)/r_2$
is independent of $\eta$, so its source-coordinate comparison
identifies its companion complex with the constant
$B_R^\vee$. For the first group, the affine substitution $z=a_1+r_1t_1$
identifies the whole $\eta$-collar with the pullback by $p_2$
of the binary collar of bounds $P,Q$, since $r_1\ne0$.
Its output coordinate at $\eta=0$ is $t_1$. For $\eta>0$,
the formal coordinate $t'_j$ at its $j$th marked point satisfies
\[
 t_1=\eta b_j+\eta s_jt'_j\qquad(j=1,2).
\]
The composition comparison in
Lemma~\ref{lem:source-coordinate-factorization}(iii)
identifies the composed positive-locus substitution with
$z=a_1+\eta r_1b_j+\eta r_1s_jt'_j$.
The compatibility of restriction arrows in
Lemma~\ref{lem:source-specialization-comparison}(ii), together
with collar naturality and the product comparison in
Lemma~\ref{lem:parameter-chain-properties}(i),(iii), now gives
\[
 r_{\{1\},\{1,2\}}
       \simeq(p_2^{-1}r_{2,(P,Q)})\otimes
                      \mathrm{id}_{B_R^\vee}.
\]
Tensor products here and below use the perfect-dual
identifications and the fixed Koszul symmetry.

On $\varepsilon=0$, the pulled-back family and its comparison
with $B_N^\vee$ both have center $0$ and coordinate $z$,
independently of $\eta$. Restriction of this constant diagram
therefore gives
$r_{\varnothing,\{2\}}\simeq\mathrm{id}_{B_N^\vee}$.
On $\eta=1$, the centers and multipliers are
\[
 \varepsilon(a_1+r_1b_1),\quad
 \varepsilon(a_1+r_1b_2),\quad \varepsilon a_2;
 \qquad
 \varepsilon r_1s_1,\quad\varepsilon r_1s_2,\quad
 \varepsilon r_2.
\]
These are precisely the collar data of Step~1 pulled back
along the insertion map $\circ_1:T\to\mathcal D_3$.
The positive-face comparison above identifies the entire
restriction arrow, giving
\[
 r_{\{2\},\{1,2\}}
       \simeq(\circ_1)^{-1}r_{3,(P,Q,R)}.
\]
Dualizing the coherent square and using
$\mathbb D_T f^{-1}\simeq f^!\mathbb D$ for these maps of
parameter spaces proves
\[
 (\circ_1)^!\alpha_{3,(P,Q,R)}
 \simeq
 (p_1^!\alpha_{2,(P+Q,R)})\circ
 \bigl((p_2^!\alpha_{2,(P,Q)})\otimes\mathrm{id}_{B_R}\bigr).
\]
The right side is the composite
\[
 \omega_T\otimes B_P\otimes B_Q\otimes B_R
 \longrightarrow\omega_T\otimes B_{P+Q}\otimes B_R
 \longrightarrow\omega_T\otimes B_N;
\]
in particular, there is one factor $\omega_T$ throughout.
This proves the derived comparison in (ii) for $k=l=2$;
relabeling gives the other insertion slot.

\emph{Step 4: prove the unary insertion cases of (ii).}
The remaining cases have $k=1$ or $l=1$. In each case let
$T=\mathcal D_k\times\mathcal D_l$, with projections
$p_1:T\to\mathcal D_k$ and $p_2:T\to\mathcal D_l$. The geometrical strategy is clear by now, each parameter insertion is modelled by loops insertion. This  dictates the definition of the centers and multipliers below.
The symbols $W_A,v_A,r_{A,A'}$ will refer to the square and its
relative open-star restrictions as in Step~3, for the
pulled-back sheaf specified below.

If $k=1$, put $N=P_1=\sum_jQ_j$ and define
\[
 c_j=\varepsilon a_1+\varepsilon\eta r_1b_j,
 \qquad u_j=\varepsilon\eta r_1s_j\quad(1\le j\le l).
\]
Let $\chi:[0,1]^2\times T\to(\A^l)^{\mathrm{an}}$ be
$(c_j)_j$ and put $E_\chi=\chi^{-1}\mathcal E_{[l],\mathbf Q}$.
On $\varepsilon,\eta>0$ apply the singleton-partition
comparison with pairs $(c_j,u_j)$. On $\eta=0<\varepsilon$
apply the one-part comparison with center $\varepsilon a_1$,
multiplier $\varepsilon r_1$, and bound $N$; on
$\varepsilon=0$ apply it with center $0$ and multiplier $1$.
At $(0,0),(0,1),(1,0),(1,1)$ the resulting coefficient
objects are respectively
\[
 B_N^\vee,\quad B_N^\vee,\quad B_N^\vee,\quad
       \left(\bigotimes_jB_{Q_j}\right)^\vee.
\]
On $\eta=0$, the $\varepsilon$-collar is the unary collar
pulled back by $p_1$. On $\varepsilon=1$, substitution
$z=a_1+r_1t_1$ turns the centers and multipliers into
$\eta b_j,\eta s_j$, giving the $l$-input collar pulled
back by $p_2$. Therefore the route through $(1,0)$ is
\[
 \underline{B_N^\vee}_T
 \xrightarrow{p_1^{-1}r_{1,(N)}}\underline{B_N^\vee}_T
 \xrightarrow{p_2^{-1}r_{l,\mathbf Q}}
       \underline{(\bigotimes_jB_{Q_j})^\vee}_T.
\]
The route through $(0,1)$ is the identity on $B_N^\vee$
followed by $(\circ_1)^{-1}r_{l,\mathbf Q}$, where
$\circ_1:T\to\mathcal D_l$ inserts the second configuration
into the unary first configuration. Duality thus gives
\[
 (\circ_1)^!\alpha_{l,\mathbf Q}
 \simeq(p_1^!\alpha_{1,(N)})\circ(p_2^!\alpha_{l,\mathbf Q}),
\]
as maps from $\omega_T\otimes\bigotimes_jB_{Q_j}$ to
$\omega_T\otimes B_N$. This includes $k=l=1$.

If $l=1<k$, then $Q_1=P_i$. Put $N=\sum_hP_h$ and define
\[
 \begin{aligned}
 c_h&=\varepsilon a_h,&u_h&=\varepsilon r_h&&(h\ne i),\\
 c_i&=\varepsilon a_i+\varepsilon\eta r_ib_1,
 &u_i&=\varepsilon\eta r_is_1.
 \end{aligned}
\]
Let $\chi=(c_h)_h:[0,1]^2\times T\to(\A^k)^{\mathrm{an}}$
and $E_\chi=\chi^{-1}\mathcal E_{[k],\mathbf P}$.
On $\varepsilon,\eta>0$ the singleton-partition comparison
has pairs $(c_h,u_h)$. On $\eta=0<\varepsilon$ its pairs
are $(\varepsilon a_h,\varepsilon r_h)$ for every $h$,
including $h=i$. On $\varepsilon=0$ the one-part comparison
has center $0$, multiplier $1$, and bound $N$.
The four vertex objects, in the same order as above, are
\[
 B_N^\vee,\quad B_N^\vee,\quad
 \left(\bigotimes_hB_{P_h}\right)^\vee,\quad
 \left(\bigotimes_hB_{P_h}\right)^\vee.
\]
The first arrow of the route through $(1,0)$ is
$p_1^{-1}r_{k,\mathbf P}$. On $\varepsilon=1$, the substitution
$z=a_i+r_it_i$ identifies the $i$th family with the unary
collar of center $\eta b_1$ and multiplier $\eta s_1$.
The other centers $a_h$ and coordinates $(z-a_h)/r_h$ are
fixed in $\eta$. Separation therefore identifies the second
arrow with $p_2^{-1}r_{1,(P_i)}$ on the $i$th dual factor and
the identity on every other factor. After duality this route
first applies $p_2^!\alpha_{1,(P_i)}$ to the $i$th coefficient
factor of $\omega_T\otimes\bigotimes_hB_{P_h}$, with
identities on the other factors, and then applies
$p_1^!\alpha_{k,\mathbf P}$. The fixed Koszul symmetry places
the single factor $\omega_T$ next to that $i$th factor when
forming the first map.

In both unary cases, the comparisons just stated are obtained
from Lemma~\ref{lem:source-coordinate-factorization}(i)--(ii)
with the displayed substitutions, pulled back by its part~(iv).
Its part~(iii) compares the composed coordinates; the diagonal
and separation comparisons are those of
Lemma~\ref{lem:relative-fusion-factorization}(i)--(ii).
Collar naturality and the face comparisons proved in Step~3
identify the morphisms, not just their endpoint objects.
On $\varepsilon=0$ the restriction is the identity, and on
$\eta=1$ the center and multiplier formulas are those for
$\circ_i:T\to\mathcal D_{k+l-1}$.
The coherent square thus identifies the two routes and proves
all the remaining derived insertion comparisons in (ii).
The repeated coefficient objects at the vertices do not make
the unary arrows identities: the arrows retain the dependence
on the framed configurations. The identity asserted here is
the restriction on $\varepsilon=0$, where both the family and
its output-coordinate comparison are constant.

\emph{Step 5: verify coherence for repeated insertions and pole maps.}
We describe one cube for each finite pattern of nested
insertions. Represent the pattern by a finite ordered rooted
tree with $n\le3$ ordered leaves. Let $V$ be its set of
internal vertices. An internal vertex $v$ has $k_v$ ordered
children, a configuration
$((a_{v,j},r_{v,j}))_{j=1}^{k_v}\in\mathcal D_{k_v}$,
and a scale $t_v\in[0,1]$. Put
$T=\prod_{v\in V}\mathcal D_{k_v}$ and assign the ordered
leaves pole bounds $R_1,\ldots,R_n$, with $N=\sum_\ell R_\ell$.
For every vertex define its incoming center and multiplier
recursively: at the root they are $(0,1)$, and for the
$j$th child $w$ of $v$ set
\[
 c_w=c_v+u_vt_va_{v,j},\qquad
 u_w=u_vt_vr_{v,j}.
\]
These formulas define the leaf centers for all scales, hence
the map and its pulled-back sheaf
\[
 \chi:[0,1]^V\times T\longrightarrow(\A^n)^{\mathrm{an}},
 \qquad\chi=(c_\ell)_{\ell=1}^n,
 \qquad E_\chi=\chi^{-1}\mathcal E_{[n],(R_1,\ldots,R_n)}.
\]
If all scales are positive, the nesting and disjointness of
the discs make the leaf centers distinct, with nonzero
coordinate multipliers $u_\ell$.

On a stratum specified by which scales vanish, traverse the
tree from the root. At the first vertex $v$ with $t_v=0$,
replace all its descendant leaves by one group with center
$c_v$, multiplier $u_v$, and bound
$\sum_{\ell\text{ below }v}R_\ell$, and stop the traversal
below $v$. If no zero scale is encountered before a leaf,
retain that leaf with its pair $(c_\ell,u_\ell)$ and bound
$R_\ell$. The retained multipliers are nonzero: their
defining products involve only scales at strict ancestors
already known to be positive and nonzero disc multipliers.
In particular, $t_{\mathrm{root}}=0$ gives the single pair
$(0,1)$ and bound $N$, independently of all other scales.
The retained centers are pairwise distinct: at the first
branching vertex of any two retained groups, the groups lie
in disjoint discs. Thus these groups and pairs specify the
partition $\pi$ and a map to $T_\pi$ on each stratum. Applying
Lemma~\ref{lem:source-coordinate-factorization}(i)--(ii),(iv)
therefore identifies the restriction of $E_\chi$ there with
the constant dual tensor product of the corresponding $B$
complexes. A unary vertex is treated by the same rule even
though collapsing it does not change the number of groups.

For $A\subseteq V$, let
\[
 W_A=\{t\in[0,1]^V:t_v>0\text{ for all }v\in A\},
 \qquad
 G(A)=R(\operatorname{pr}_T)_*(E_\chi|_{W_A\times T}).
\]
Restriction along $W_{A'}\subseteq W_A$ for $A\subseteq A'$
defines a functor
$G:N\mathcal P(V)\to\mathcal D(T;\F)$.
The relative evaluation argument of Step~3 identifies $G(A)$
with the restriction at the vertex having coordinates $1$
on $A$ and $0$ elsewhere. By
Lemma~\ref{lem:parameter-chain-properties}(iii), a chain
$A_0\subseteq\cdots\subseteq A_q$ supplies a $q$-simplex of
this functor, with its composition homotopies, faces, and
degeneracies.

We must also identify the faces with the insertion diagrams
they compare. A face $t_v=0$ collapses the groups specified
above. If $v$ is not the root, the face $t_v=1$ composes the
configuration at $v$ with its slot in the parent configuration:
if $v$ is the $j$th child of $p$, and $w$ is the $h$th
child of $v$, the recurrence on $t_v=1$ gives
\[
 \begin{aligned}
 c_w&=c_p+u_pt_p(a_{p,j}+r_{p,j}a_{v,h}),\\
 u_w&=u_pt_pr_{p,j}r_{v,h}.
 \end{aligned}
\]
These are the center and multiplier formulas for the insertion of
Definition~\ref{def:framed-disc-parameters}, with the same
leaf centers and multipliers as the tree in which that edge
has been contracted, pulled back to $T$ along the insertion
map on the corresponding configuration factors. The zero-face
comparisons are those of
Lemma~\ref{lem:parameter-chain-properties}(iii); the
positive-face comparisons follow from evaluation at the
common vertices, as proved in Step~3. Restricting successively
to any such faces gives the same inverse-image sheaf, with
the canonical composition comparisons for those restrictions.
Lemma~\ref{lem:source-coordinate-factorization}(iii) makes
the coordinate comparisons agree under these affine
compositions, and
Lemma~\ref{lem:relative-fusion-factorization}(i)--(ii)
includes the diagonal and separation comparisons.
Consequently the simplices of $G$, and their duals, compare
the iterated operations with the operations for the inserted
configurations, compatibly for different orders of insertion.
There is no bound on $|V|$: this includes arbitrarily many
unary insertions with at most three final inputs.

Increasing the leaf pole bounds gives natural transformations
of these same diagrams $G$, by the proper-trace comparisons
in Lemma~\ref{lem:relative-fusion-factorization}(i)--(ii)
and their source-coordinate compatibility in
Lemma~\ref{lem:source-coordinate-factorization}(iii).
Naturality of restriction and duality preserves their
compatibility with all the simplices and face comparisons.
This proves the derived coherence in (ii) and its compatibility
with the pole comparisons of (i).

\emph{Step 6: realize the operations and homotopies on the fixed filtered model.}
Apply Definition~\ref{def:parameter-chain-representative}
to the entire diagrams of derived maps and comparisons in
Steps~1--5, including the pole and relabeling comparisons.
Lemma~\ref{lem:parameter-chain-operations}(i)
represents a family $\alpha_{k,\mathbf P}$ by a degree-zero
map with parameter complex $C_{-*}(\mathcal D_k;\F)$ and
gives its boundary formula. It represents pullback along
$\circ_i$ by precomposition with $(\circ_i)_*$.
Part~(ii) of that lemma represents composition over one
parameter space using its Alexander--Whitney diagonal;
part~(iii) compares exterior products with Eilenberg--Zilber
and the Koszul interchange. Together they turn the comparisons
of Steps~3--5 into comparisons with
$c\circ_i d=(\circ_i)_*\operatorname{EZ}(c\otimes d)$.

Transport this entire diagram through the specified equivalence
$(F_P\mathbf B)_P\simeq(B_P)_P$ of filtered diagrams.
The fixed-cochain-model assertion of
Lemma~\ref{lem:filtered-operations-colimit}, including its
assertion for a fixed parameter complex, realizes it on the
given exhaustive filtered complex $\mathbf B$.
For an operation the parameter complex is
$C_{-*}(\mathcal D_k;\F)$; for an insertion comparison it is
$C_{-*}(\mathcal D_k;\F)\otimes C_{-*}(\mathcal D_l;\F)$,
and for repeated comparisons it is the corresponding finite
tensor product. Each parameter complex has filtration in
index zero. The resulting filtered maps $\Theta_k$ restrict
to the maps $\Theta_{k,\mathbf P}$ of (i), and all comparison
homotopies have image in $F_{\sum_iP_i}\mathbf B$ on inputs
of bounds $\mathbf P$.

To obtain literal relabeling equivariance, replace each
$\Theta_k$ by
\[
 \frac{1}{k!}\sum_{\sigma\in\mathfrak S_k}\Theta_k\circ\sigma,
\]
where $\mathfrak S_k$ is the permutation group and $\sigma$
acts diagonally on the domain, by relabeling the parameter
configuration and by the Koszul permutation of the input
factors. Since the coefficients are rational, this averaging
is defined and preserves the sum filtration. Averaging the specified
relabeling homotopies identifies the averaged map with the
original map. Transport the insertion and higher comparisons
along these filtered homotopies, retaining the notation
$\Theta_k$ for the resulting equivariant maps.

For clarity, the two maps in (ii) and their homotopy have
the following meaning on these representatives. Put
$m=k+l-1$, $N=\sum_hP_h=\sum_hR_h$, with $\mathbf Q$ and
$\mathbf R$ as in (ii). Define $\mathsf L_i$ and $\mathsf R_i$
to be the degree-zero cochain maps with common domain and target
\[
 C_{-*}(\mathcal D_k;\F)\otimes C_{-*}(\mathcal D_l;\F)
       \otimes\bigotimes_{h=1}^mF_{R_h}\mathbf B
       \longrightarrow F_N\mathbf B
\]
given, for $c\in C_p(\mathcal D_k;\F)$,
$d\in C_q(\mathcal D_l;\F)$, and homogeneous inputs
$x_h\in F_{R_h}\mathbf B$, by
\[
 \begin{aligned}
 \mathsf L_i(c,d;x_1,\ldots,x_m)
   &=\Theta_{m,\mathbf R}
       ((\circ_i)_*\operatorname{EZ}(c\otimes d))
                         (x_1,\ldots,x_m),\\
 \mathsf R_i(c,d;x_1,\ldots,x_m)
   &=(-1)^{q\sum_{h<i}|x_h|}\Theta_{k,\mathbf P}(c)
       \bigl(x_1,\ldots,x_{i-1},
          \Theta_{l,\mathbf Q}(d)(x_i,\ldots,x_{i+l-1}),
          x_{i+l},\ldots,x_m\bigr).
 \end{aligned}
\]
Here $|x_h|$ is cohomological degree, read modulo two in
parity grading; the sign moves the parameter chain $d$, of
cohomological degree $-q$, past the preceding inputs.
The specified insertion comparison is a degree-$-1$ map $K_i$
from this same domain to $F_N\mathbf B$ such that
\[
 d_{\operatorname{Hom}}K_i=\mathsf L_i-\mathsf R_i.
\]
Thus $d_{\operatorname{Hom}}K_i=d_{F_N\mathbf B}K_i+
K_i d_{\mathrm{domain}}$, with the tensor-product differential
on the displayed domain. For cycles $c,d$, evaluation of
$K_i$ gives a homotopy of degree $-p-q-1$ between the
operations of degree $-p-q$. The higher comparisons of
Step~5 are realized jointly with these maps; none is chosen
independently of the others.

The actual inclusions $\jmath_{P,Q}$ are the filtration maps
of this fixed model, and the specified diagram equivalence
identifies them with the proper-trace transitions used in
Steps~1--5. This proves the pole comparisons, completes (i)
and (ii), and proves the filtered assertions following them.
Comparisons with the previously chosen fusion and rotation
maps are included in the same diagram, as described in
Corollary~\ref{cor:framed-basic-operations}.
Finally, if parameter cycles satisfy $c'-c=\partial b$,
Lemma~\ref{lem:parameter-chain-operations}(i) gives
\[
 d_{\operatorname{Hom}}\Theta_{k,\mathbf P}(b)
       =\Theta_{k,\mathbf P}(c')-\Theta_{k,\mathbf P}(c),
\]
with output bound $\sum_iP_i$. This is the required filtered
homotopy for homologous parameter cycles.
\end{proof}

For $k=1$, put
$
 f:\mathcal D_1\longrightarrow G=\C^*$, $f(a,r)=r$.
The unary family is related to the source-rotation family
$\Theta_P^G$ of
Proposition~\ref{prop:coherent-source-rotation} by
\[
 \Theta_{1,(P)}\simeq
 \Theta_P^G\circ(f_*\otimes\mathrm{id}_{F_P\mathbf B}).
\]
Here $f_*$ is the induced map on singular chains, and
$\simeq$ denotes a homotopy of cochain maps from
$C_{-*}(\mathcal D_1;\F)\otimes F_P\mathbf B$ to
$F_P\mathbf B$, compatible with the pole comparisons.
This comparison of the entire parameter-chain families is
one part of the Corollary~\ref{cor:framed-basic-operations} below. The other part of next Corollary is based on the observation that the family in Step 1 of the previous proof for $k=2$ can also be used to define the fusion product.

\begin{corollary}[Fusion and rotation among the framed operations]
\label{cor:framed-basic-operations}
For $m\in\mathcal D_2$ and a counterclockwise framing cycle
$\tau$ in $\mathcal D_1$ given by $a=0$, $r=\rho e^{i\phi}$,
$0<\rho<1$, the pole-stage operations satisfy
\[
 \Theta_{2,(P,Q)}([m])\simeq m_{2;P,Q},\qquad
 \Theta_{1,(P)}(\tau)\simeq\Delta_P.
\]
The homotopies preserve the pole comparisons, so the same
identifications hold on $\mathbf B$ and on $LH$.
\end{corollary}
\begin{proof}
\emph{Step 1: compare the varying input multipliers with fixed ones.}
For $k=1,2$ and bounds $\mathbf P$, put
$T=\mathcal D_k$, $N=\sum_iP_i$, $H=[0,1]$, $L=(0,1]$,
and $I=[0,1]$. The coordinates on $H$ and $I$ are
$\varepsilon$ and $v$, respectively. Define
\[
 h:H\times T\longrightarrow(\A^k)^{\mathrm{an}},
 \quad h(\varepsilon,(a,r))=(\varepsilon a_i)_i,
 \qquad E=h^{-1}\mathcal E_{[k],\mathbf P}.
\]
Let $\widetilde E$ be the inverse image of $E$ under the
projection $H\times T\times I\to H\times T$. On
$L\times T\times I$, the singleton-partition comparison
of Lemma~\ref{lem:source-coordinate-factorization}(i)--(ii)
is pulled back along
\[
 (\varepsilon,(a,r),v)\longmapsto
 \bigl((\varepsilon a_i)_i,(u_i(\varepsilon,v))_i\bigr),
 \qquad u_i(\varepsilon,v)=((1-v)\varepsilon+v)r_i.
\]
The centers are distinct and each $u_i(\varepsilon,v)$ is
nonzero there. These substitutions give input coordinates
$t_i=(z-\varepsilon a_i)/u_i(\varepsilon,v)$.
They are formal source coordinates; the intermediate
multipliers need not define disjoint discs. Part~(iv) of
that lemma applies by the semialgebraic pullback and locally
oriented product comparisons verified in Step~1 of
Proposition~\ref{lem:framed-collision-operations}.
On $\varepsilon=0$, the one-part comparison has center $0$,
multiplier $1$, and bound $N$, for every $v$.

Let $i_0:T\times I\hookrightarrow H\times T\times I$ and
$j_+:L\times T\times I\hookrightarrow H\times T\times I$
be the zero-end and positive-locus inclusions. The specified
comparisons and Lemma~\ref{lem:parameter-chain-properties}(i)
give
\[
 \underline{B_N^\vee}_{T\times I}
 \simeq i_0^{-1}\widetilde E
 \xrightarrow{i_0^{-1}\eta_{\widetilde E}}
 i_0^{-1}Rj_{+*}j_+^{-1}\widetilde E
 \simeq\underline{(\bigotimes_iB_{P_i})^\vee}_{T\times I},
\]
where $\eta_{\widetilde E}$ is the adjunction unit.
Denote its Verdier dual by
\[
 \beta:\omega_{T\times I}\otimes\bigotimes_iB_{P_i}
                \longrightarrow\omega_{T\times I}\otimes B_N.
\]
For $v=0,1$, let $\iota_v:T\to T\times I$ be
$x\mapsto(x,v)$ and write $\beta_v=\iota_v^!\beta$.
Then $\beta_0\simeq\alpha_{k,\mathbf P}$ of the preceding
proposition, while $\beta_1$ is the collar operation with
input multipliers $r_i$, independent of $\varepsilon$.
The interval parameter compares these entire families,
including their pole morphisms, by
Lemma~\ref{lem:parameter-chain-operations}(i),(iii).
The output coordinate remains $z$ throughout; no input
coordinate is asserted at $\varepsilon=0$.

\emph{Step 2: identify the binary operation with fusion.}
Take
$m_0=((1/2,1/4),(-1/2,1/4))\in\mathcal D_2$ and let
$i_{m_0}:\{m_0\}\hookrightarrow\mathcal D_2$ be its inclusion.
The center map at $m_0$ is
$(\varepsilon/2,-\varepsilon/2)$, so $s=\varepsilon/2$
identifies this collar with the binary collision slice
$\mathcal Z_{P,Q}$ on $0\le s\le1/2$.
The open-star evaluation equivalences of
Lemma~\ref{lem:parameter-chain-properties}(ii)--(iii)
identify its restriction arrow with
\eqref{eq:binary-restriction-map}, restricted to this interval.
With input coordinates $z-s,z+s$, its dual is therefore
$\overline m_{2;P,Q}$ by
Proposition~\ref{prop:relative-binary-identification}.

At $v=1$ in Step~1, the actual input coordinates are
$4(z-s),4(z+s)$. Changing the two multipliers from $1$ to
$1/4$ changes the compact-support comparison by
$(r_{1/4})_!\otimes(r_{1/4})_!$, by
Lemma~\ref{lem:source-coordinate-factorization}(iii).
Duality consequently gives the precise comparison
\[
 i_{m_0}^!\alpha_{2,(P,Q)}
 \simeq i_{m_0}^!\beta_1
 \simeq\overline m_{2;P,Q}\circ
                  (r_{1/4}^!\otimes r_{1/4}^!).
\]
All three maps have source $B_P\otimes B_Q$ and target
$B_{P+Q}$. Here $r_{1/4}^!$ is the source-rotation map on
the indicated input complex. The path
$v\mapsto(1-v)/4+v$ in $\C^*$ joins $1/4$ to $1$.
The action comparison of
Proposition~\ref{prop:coherent-source-rotation}, equivalently
Remark~\ref{rem:fixed-source-dilation}, therefore gives
$r_{1/4}^!\simeq\mathrm{id}$ compatibly with pole maps.
This identifies the derived operation at $m_0$ with fusion.

Transporting this comparison through the specified filtered
equivalences, using Lemma~\ref{lem:filtered-operations-colimit},
gives a degree-$-1$ homotopy
\[
 K_{0;P,Q}:F_P\mathbf B\otimes F_Q\mathbf B
                       \longrightarrow F_{P+Q}\mathbf B,
 \qquad
 d_{\operatorname{Hom}}K_{0;P,Q}
       =\Theta_{2,(P,Q)}([m_0])-m_{2;P,Q}.
\]
For any $m\in\mathcal D_2$, choose a path $c_m$ from $m_0$
to $m$, viewed as a singular $1$-chain; such a path exists
because $\mathcal D_2$ is path connected. Since
$\partial c_m=[m]-[m_0]$, the boundary formula in
Proposition~\ref{lem:framed-collision-operations}(i) gives
\[
 d_{\operatorname{Hom}}
       \bigl(K_{0;P,Q}+\Theta_{2,(P,Q)}(c_m)\bigr)
       =\Theta_{2,(P,Q)}([m])-m_{2;P,Q}.
\]
This proves the first assertion, on the stated filtration
subcomplexes and compatibly with their pole comparisons.

\emph{Step 3: identify the unary family and its framing circle with source rotation.}
For $k=1$, the singleton comparison with center
$\varepsilon a$ and multiplier $1$ extends over all of
$H\times\mathcal D_1$. Its zero-end value is the output
comparison of Step~1. Under this trivialization, the collar
restriction is the identity on $B_P^\vee$.
For the input multiplier $r$, let $e_1,e_r$ denote the
compact-support comparisons with multipliers $1,r$ over
$L\times\mathcal D_1$. Part~(iii) of
Lemma~\ref{lem:source-coordinate-factorization} gives
\[
 e_r=(r_r)_!\circ e_1.
\]
Let $f$ be the multiplier map defined after
Proposition~\ref{lem:framed-collision-operations}.
Recall that $\mathcal R_P=\mathbb D_G((r_\lambda)_!)$ is
the rotation morphism of
Proposition~\ref{prop:coherent-source-rotation}.
Collar naturality and duality, together with Step~1, give
\[
 \alpha_{1,(P)}\simeq\beta_1\simeq f^!\mathcal R_P
 \quad\text{on }\omega_{\mathcal D_1}\otimes B_P.
\]
This is a comparison of the families over $\mathcal D_1$;
its extraordinary restriction at $(a,r)$ is $r_r^!$.
By Lemma~\ref{lem:parameter-chain-operations}(i), pullback
along $f$ corresponds on parameter-chain maps to
precomposition with $f_*\otimes\mathrm{id}$.
Joint realization of this family comparison and its pole
comparisons on the fixed filtered cochain model therefore
proves the unary comparison stated after
Proposition~\ref{lem:framed-collision-operations}.

Let $\sigma\in Z_1(S^1;\F)$ be the counterclockwise
fundamental cycle used in
Definition~\ref{def:source-rotation-operator}. Define
\[
 j_\rho:S^1\longrightarrow\mathcal D_1,
       \quad j_\rho(\lambda)=(0,\rho\lambda),
 \qquad q_\rho:S^1\longrightarrow G,
       \quad q_\rho(\lambda)=\rho\lambda.
\]
Then $f\circ j_\rho=q_\rho$, and the preceding family
comparison gives
$j_\rho^!\alpha_{1,(P)}\simeq q_\rho^!\mathcal R_P$.
Restricting the resulting cochain homotopy along
$(j_\rho)_*\otimes\mathrm{id}$ compares
$\Theta_{1,(P)}((j_\rho)_*\sigma)$ with
$\Theta_P^G((q_\rho)_*\sigma,-)$.

The radial homotopy
\[
 h_\rho:I\times S^1\longrightarrow G,
       \qquad h_\rho(v,\lambda)=((1-v)\rho+v)\lambda
\]
joins $q_\rho$ to the inclusion $q_1:S^1\hookrightarrow G$.
Let $[I]$ be the singular $1$-simplex $v\mapsto v$, oriented
from $0$ to $1$, and put
$b=(h_\rho)_*\operatorname{EZ}([I]\otimes\sigma)
\in C_2(G;\F)$.
The chain boundary identity gives
\[
 \partial b=(q_1)_*\sigma-(q_\rho)_*\sigma.
\]
Thus Lemma~\ref{lem:parameter-chain-operations}(i) gives the
degree-$-2$ homotopy
\[
 d_{\operatorname{Hom}}\Theta_P^G(b,-)
       =\Delta_P-\Theta_P^G((q_\rho)_*\sigma,-).
\]
The specified framing cycle $\tau$ is homologous to
$(j_\rho)_*\sigma$ in the same framing circle. The last
assertion of Proposition~\ref{lem:framed-collision-operations}
therefore identifies their operations by a filtered homotopy.
Combining these comparisons proves
$\Theta_{1,(P)}(\tau)\simeq\Delta_P$.
The variance is the oriented pullback $r_\lambda^!$.
For $g\in\mathcal O(X)$, writing
$\gamma(g)(t)=\sum_{\nu\in\Z}x_\nu t^\nu$ gives
$(r_\lambda\gamma)(g)(t)=\sum_\nu\lambda^\nu x_\nu t^\nu$.
Thus a nonzero coefficient $x_\nu$ winds $+\nu$ times
under the counterclockwise circle parameter; no inversion
of that parameter occurs.

All comparisons above are formed with the entire pole
diagrams, by Lemma~\ref{lem:source-coordinate-factorization}(iii).
Lemma~\ref{lem:filtered-operations-colimit} realizes them
jointly, with the fixed endpoints $m_{2;P,Q}$ and
$\Delta_P$. The displayed homotopies therefore restrict
from homotopies on $\mathbf B$ and induce the asserted
equalities on $LH$.
\end{proof}

\subsection{The BV identity and invariance}
\label{subsec:bv-identity-invariance}

The following is well-known and it is a chain-level formulation of the seven-term BV relation (\ref{eq:bv-seven-term}) in the homology operad of framed little disks, cf. \cite[Theorem 4.5]{GetzlerBV}.

\begin{lemma}[The seven-term parameter cycle]
\label{lem:seven-term-disc-cycle}
Let $m\in Z_0(\mathcal D_2;\F)$ be a point, let
$\tau\in Z_1(\mathcal D_1;\F)$ be the counterclockwise circle
$a=0$, $r=\rho e^{i\phi}$, where $0<\rho<1$, and put
$m^{(3)}=m\circ_1m$. With composition of parameter chains
given by Eilenberg--Zilber and disc insertion, the cycle
\begin{equation}\label{eq:framed-bv-relations}
\begin{split}
 Z_7={}&\tau\circ_1m^{(3)}-m\circ_1(\tau\circ_1m)
       -m\circ_2(\tau\circ_1m)\\
 &-[m\circ_2(\tau\circ_1m)]\cdot(12)
       +\sum_{i=1}^3m^{(3)}\circ_i\tau
\end{split}
\end{equation}
is a boundary in $C_*(\mathcal D_3;\F)$. Here $(12)\in\mathfrak S_3$ exchanges the first two input
labels. For a singular chain $c$ on $\mathcal D_3$, the
notation $c\cdot(12)$ means its image under the chain map
induced by
\[
 ((a_1,r_1),(a_2,r_2),(a_3,r_3))
 \longmapsto
 ((a_2,r_2),(a_1,r_1),(a_3,r_3)).
\]
\end{lemma}
\begin{proof}
All terms are cycles because insertion is a chain map.
Write $\mathring D=\{z\in\C:|z|<1\}$ and
$\operatorname{Conf}_3(U)=\{(a_1,a_2,a_3)\in U^3:a_i\ne a_j
\text{ for }i\ne j\}$.
Forgetting the positive radii identifies the homotopy type of
$\mathcal D_3$ with
$\operatorname{Conf}_3(\mathring D)\times(S^1)^3$:
the radius inequalities define convex fibers with a continuous
section given by choosing all radii sufficiently small.
The inclusion $\operatorname{Conf}_3(\mathring D)
\hookrightarrow\operatorname{Conf}_3(\C)$ is a homotopy
equivalence, by simultaneous positive rescaling of the centers.
Moreover
\[
 \operatorname{Conf}_3(\C)\cong
 \C\times\C^*\times(\C\setminus\{0,1\}),\qquad
 (a_1,a_2,a_3)\longmapsto
 \left(a_1,a_2-a_1,\frac{a_3-a_1}{a_2-a_1}\right).
\]
Put $z=a_2-a_1$ and $w=(a_3-a_1)/(a_2-a_1)$ in this
description. The three differences are $z$, $zw$, and $z(w-1)$.
Their windings recover those of $z$, $w$, and $w-1$ by
subtracting the first winding from the second and third.
Consequently $H_1(\mathcal D_3;\Z)\cong\Z^6$ is detected
by the three windings of $a_i-a_j$ and the three windings
of $r_i$.

In $\tau\circ_1m^{(3)}$, every center difference and every
frame winds once counterclockwise. The next three terms rotate
respectively the clusters $\{1,2\}$, $\{2,3\}$, and
$\{1,3\}$. Such a rotation gives winding one to the
corresponding center difference and to both frames in its
cluster, and winding zero to all other coordinates.
Indeed a moving center stays in its cluster disc, which is
disjoint from the spectator disc, so its difference from the
spectator has winding zero. Finally $m^{(3)}\circ_i\tau$
has only the $i$th frame winding, equal to one.
Thus each center-difference winding of $Z_7$ is $1-1=0$,
and each frame winding is $1-2+1=0$.
This proves $[Z_7]=0$, already with integral coefficients.
\end{proof}

For $\mathbf P=(P,Q,R)$, put $N=P+Q+R$.
Define the degree-$-1$ \emph{seven-term defect}
$\mathcal R_{7,\mathbf P}:F_P\mathbf B\otimes F_Q\mathbf B
\otimes F_R\mathbf B\to F_N\mathbf B$
on homogeneous $a\in(F_P\mathbf B)^p$, $b\in(F_Q\mathbf B)^q$,
$c\in F_R\mathbf B$ by
\begin{equation}\label{eq:seven-term-cochain-defect}
 \begin{split}
 \mathcal R_{7,\mathbf P}(a,b,c)={}&
 \Delta((ab)c)-\Delta(ab)c-(-1)^p a\Delta(bc)\\
 &-(-1)^{pq+q}b\Delta(ac)
       +(\Delta(a)b)c+(-1)^p(a\Delta(b))c\\
 &+(-1)^{p+q}(ab)\Delta(c).
 \end{split}
\end{equation}
Every product denotes the appropriate pole-stage map $m_{2;A,B}$,
and $\Delta$ on $F_A\mathbf B$ denotes $\Delta_A$. For example,
$(ab)c=m_{2;P+Q,R}(m_{2;P,Q}(a,b),c)$.
Thus every term has target $F_N\mathbf B$, with the displayed
parenthesization. Strict associativity would require the equality
\[
 m_2\circ(m_2\otimes\mathrm{id}_{\mathbf B})
 =m_2\circ(\mathrm{id}_{\mathbf B}\otimes m_2)
 :\mathbf B^{\otimes3}\longrightarrow\mathbf B
\]
of cochain maps, using the canonical identification of the
tensor-product parenthesizations. This equality is not assumed;
strictness of the pole-stage diagram does not imply it.

\begin{theorem}[The loop BV algebra]\label{thm:loop-bv}
For the fusion and rotation maps already constructed, there are
pole-compatible homotopies
\[
 H_{7,\mathbf P}\in
 \operatorname{Hom}^{-2}
 (F_P\mathbf B\otimes F_Q\mathbf B\otimes F_R\mathbf B,
                                      F_{P+Q+R}\mathbf B),
 \qquad d_{\operatorname{Hom}}H_{7,\mathbf P}
                       =\mathcal R_{7,\mathbf P}.
\]
Consequently the fusion algebra $LH^*(X)$ with $\Delta$ is
a filtered BV algebra in parity grading. A chosen nowhere-zero
algebraic volume form $\Omega$ gives the integer-graded version
$LH^*_\Omega(X)$. Its bracket satisfies
\[
 \{F_PLH^a,F_QLH^b\}\subset F_{P+Q}LH^{a+b-1}.
\]
\end{theorem}
\begin{proof}
Choose $c_7\in C_2(\mathcal D_3;\F)$ with $\partial c_7=Z_7$
as in Lemma~\ref{lem:seven-term-disc-cycle}. Since a parameter
chain of dimension two gives an operation of degree $-2$,
Proposition~\ref{lem:framed-collision-operations}(i) gives
$|\Theta_{3,\mathbf P}(c_7)|=-2$ and
\[
 d_{\operatorname{Hom}}\Theta_{3,\mathbf P}(c_7)
                         =\Theta_{3,\mathbf P}(Z_7).
\]
Let $z_1,\ldots,z_7$ be the seven cycles in
\eqref{eq:framed-bv-relations}, in displayed order after
expanding the final sum and removing the outer signs. Put
\[
 (\epsilon_1,\ldots,\epsilon_7)=(1,-1,-1,-1,1,1,1),
 \qquad Z_7=\sum_{j=1}^7\epsilon_jz_j.
\]
For each $j$, let $R_{j,\mathbf P}$ be the corresponding
degree-$-1$ composite in \eqref{eq:seven-term-cochain-defect},
with its outer sign removed but its Koszul sign retained.
Thus $\mathcal R_{7,\mathbf P}=\sum_j\epsilon_jR_{j,\mathbf P}$;
for example,
$R_{1,\mathbf P}(a,b,c)=\Delta((ab)c)$ and
$R_{3,\mathbf P}(a,b,c)=(-1)^p a\Delta(bc)$.

Proposition~\ref{lem:framed-collision-operations}(ii), applied
successively to the insertions defining $z_j$, compares
$\Theta_{3,\mathbf P}(z_j)$ with the corresponding composite
of binary operations $\Theta_2(m)$ and the unary operation
$\Theta_1(\tau)$, with the required pole bounds at every vertex.
For $z_4$, relabeling equivariance in part~(i) identifies the
operation with the unpermuted operation at bounds $(Q,P,R)$,
precomposed with the Koszul transposition of the first two inputs.
The required replacements use the homotopies constructed in
the proof of Corollary~\ref{cor:framed-basic-operations}.
For input bounds $A,B$, Step~2 constructs a degree-$-1$
homotopy
\[
 \Theta_{2,(A,B)}([m])\simeq m_{2;A,B}
\]
by deforming the input multipliers, removing the resulting
source dilations, and joining the standard binary configuration
to $m$. Step~3 constructs a degree-$-2$ homotopy
\[
 \Theta_{1,(A)}(\tau)\simeq\Delta_A
\]
by evaluating the unary-family comparison with source rotation
on the framing cycle, then applying the radial homotopy from
the framing circle to the unit circle. These homotopies preserve
the pole comparisons and the specified cochain representatives
$m_{2;A,B}$ and $\Delta_A$.

Any change of parenthesization
uses the filtered associativity comparisons of
Section~\ref{sec:incidence}. These comparisons are realized
jointly on $\mathbf B$, with the specified endpoints
$m_{2;A,B}$ and $\Delta_A$ and their pole compatibility.

To combine the comparisons, use the graded Leibniz rule for
the Hom differential. If $f,f'$ are closed homogeneous maps
of the same degree, $g,g'$ are composable closed homogeneous
maps of the same degree, and $h,k$ have degrees
$|f|-1,|g|-1$ and satisfy
$d_{\operatorname{Hom}}h=f-f'$,
$d_{\operatorname{Hom}}k=g-g'$, then
\[
 d_{\operatorname{Hom}}
       \bigl(hg+(-1)^{|f|}f'k\bigr)=fg-f'g'.
\]
The same rule applies to insertion into tensor factors, with
the Koszul signs. Successive replacements add their homotopies,
since the intermediate differences cancel. They therefore give
\[
 K_{j,\mathbf P}\in\operatorname{Hom}^{-2}
 \left(F_P\mathbf B\otimes F_Q\mathbf B\otimes F_R\mathbf B,
       F_{P+Q+R}\mathbf B\right),
 \qquad
 d_{\operatorname{Hom}}K_{j,\mathbf P}
       =\Theta_{3,\mathbf P}(z_j)-R_{j,\mathbf P}.
\]
The degree is $-2$ because the compared operations have degree
$-1$. Set $K_{\mathbf P}:=\sum_j\epsilon_jK_{j,\mathbf P}$.
Linearity now gives
\[
 d_{\operatorname{Hom}}K_{\mathbf P}
       =\Theta_{3,\mathbf P}(Z_7)-\mathcal R_{7,\mathbf P}.
\]
For the transposed pair term, the Koszul permutation contributes
$(-1)^{pq}$ and insertion of $\Delta$ after $b$ contributes
$(-1)^q$. Hence its sign is $-(-1)^{pq+q}$, as displayed;
$pq+q\equiv(p-1)q\pmod2$ gives the convention of
\eqref{eq:bv-seven-term}. The remaining signs are the usual
tensor signs for inserting a degree-$-1$ map in each slot.
It follows that
\[
 H_{7,\mathbf P}:=\Theta_{3,\mathbf P}(c_7)-K_{\mathbf P}
\]
has the asserted differential. All the comparisons preserve pole
maps and assign to each intermediate output the sum of its input
bounds. Thus these homotopies have target $F_{P+Q+R}\mathbf B$
and are restrictions of filtered homotopies on $\mathbf B$ by
the filtered assertions following parts~(i)--(ii) of
Proposition~\ref{lem:framed-collision-operations}.

For cocycles $a,b,c$, the asserted Hom-differential identity gives
\[
 \mathcal R_{7,\mathbf P}(a,b,c)
       =d\bigl(H_{7,\mathbf P}(a,b,c)\bigr).
\]
Hence the seven-term defect vanishes on cohomology.
Theorem~\ref{thm:geometric-fusion-ring} and
Lemma~\ref{lem:source-dilation} supply the other BV axioms and
their required filtered homotopies. Lemma~\ref{lem:algebraic-bv-bracket}
gives the Gerstenhaber identities. The bracket bound follows
from its definition because fusion adds pole bounds and rotation
preserves them. All constructions and signs apply in both gradings.
\end{proof}

\begin{corollary}[Compatibility of rotation with the comparison maps]
\label{cor:loop-bv-invariance}
The compactification and volume-regrading comparisons of
Propositions~\ref{prop:fusion-invariance-gradings}
and~\ref{prop:fusion-volume-change} intertwine
the rotation operators. More precisely:
\begin{enumerate}[(i)]
\item For the unfiltered compactification comparison $g:B^{(1)}\simeq B^{(2)}$,
\[
 g[-1]\circ\overline\Delta^{(1)}
          =\overline\Delta^{(2)}\circ g.
\]
This identity lies in $D_{\mathrm{per}}(\F)$ in parity grading,
and in $D(\F)$ in the integer grading determined by $\Omega$.
\item Fix $u\in\mathcal O(X)^*$ and denote the summand of
valuation $\ord_tu(\gamma)=a$ by $B_{u=a}$. Rotation restricts to
$\overline\Delta:B_{u=a}\to B_{u=a}[-1]$. If
\[
 v_{u,a}:(B_{u\Omega})_{u=a}\xrightarrow{\sim}
                         (B_\Omega)_{u=a}[-2a]
\]
is the volume-regrading comparison, then
\[
 v_{u,a}[-1]\circ\overline\Delta_{u\Omega}
             =\overline\Delta_\Omega[-2a]\circ v_{u,a}
          \quad\text{in }D(\F).
\]
The periodization comparison likewise intertwines the two
rotation operators.
\end{enumerate}
Together with Proposition~\ref{prop:coherent-source-rotation} and
Lemma~\ref{lem:source-dilation}, these statements upgrade the
existing algebra comparisons to BV-algebra comparisons.
The bracket adds unit valuations.
\end{corollary}
\begin{proof}
The cofinal pole maps defining $g$ commute with
$r_\lambda\gamma(t)=\gamma(\lambda t)$ over the entire parameter
space $\C^*$. Naturality of the relative morphisms
$\mathcal R_P$ of Proposition~\ref{prop:coherent-source-rotation}
therefore intertwines their parameter-chain maps, and evaluation
on $\sigma$ proves (i).

For (ii), substitution gives
\[
 \ord_tu(r_\lambda\gamma)=\ord_tu(\gamma)
                        \qquad(\lambda\in\C^*).
\]
The rotation diagram thus preserves each open-and-closed
valuation locus. On that locus the volume-regrading comparison
is the even shift $[-2a]$ of its normalized dualizing complexes;
it commutes with the same rotation diagram, proving the displayed
identity. Periodization commutes with these even shifts and the
parameter-chain construction. Finally, bracket compatibility
follows from \eqref{eq:general-rotation-bracket}; its valuation
bound follows from (ii) and the product bound in
Proposition~\ref{prop:fusion-unit-valuations}.
\end{proof}

\begin{proof}[Proof of Theorem~\ref{thm:summary-fusion}]
The natural map is
\[
 \iota_X:H^*(X;\F)\xrightarrow{\sim}H^*(B_0)
              \longrightarrow H^*(B)=LH^*(X).
\]
Its image is $F_0LH$ by Definition~\ref{def:loop-cohomology}.
The remaining assertions follow from
Theorem~\ref{thm:canonical-complex},
Remark~\ref{rem:pole-cohomology-colimit},
Theorem~\ref{thm:geometric-fusion-ring},
Theorem~\ref{thm:loop-bv}, and
Corollary~\ref{cor:loop-bv-invariance}(i)--(ii).
\end{proof}

\subsection{The BV operator on global finite models}
\label{subsec:global-finite-bv}
The global finite-factor presentations of
Proposition~\ref{prop:global-finite-presentations}(i) also compute
source rotation. The geometric formula in part~(ii) below additionally requires
algebraic actions on the finite factors and equivariant
transition maps.

\begin{proposition}[Global finite-model formula for the BV operator]
\label{prop:global-finite-bv}
Suppose that $F_P\mathscr LX=\varprojlim_N Z_{P,N}$ has the
compatible global finite-factor presentations of
Proposition~\ref{prop:global-finite-presentations}(i). Use the notation
$B_{P,N}$ preceding part~(i) of Proposition~\ref{prop:global-binary-fusion}, so that
in integer grading
\[
 B_{P,N}=R\Gamma(Z_{P,N},
                \omega_{Z_{P,N}}[-2\delta_{P,N}])
       \simeq B_P,
 \qquad G=\C^*.
\]
Without $\Omega$, replace these complexes by their parity
versions as in that proposition. Let $\theta\in H^1(G;\F)$
evaluate to $1$ on the counterclockwise unit circle.
\begin{enumerate}[(i),wide=0pt]
\item Transport the whole rotation morphism of
Proposition~\ref{prop:coherent-source-rotation} through the
displayed equivalence, obtaining
\[
 \mathcal R_{P,N}:\omega_G\otimes B_{P,N}
                  \longrightarrow\omega_G\otimes B_{P,N}.
\]
The complex orientation $\omega_G\simeq\F_G[2]$ and the
constant-section map define
\[
 \Psi_{P,N}:B_{P,N}
 \xrightarrow{\,1\otimes-\,}
 R\Gamma(G;\F)\otimes B_{P,N}
 \xrightarrow{\,R\Gamma(\mathcal R_{P,N}[-2])\,}
 R\Gamma(G;\F)\otimes B_{P,N}.
\]
For $\alpha\in H^q(B_{P,N})$, there is a unique class
$\Delta_{P,N}\alpha\in H^{q-1}(B_{P,N})$ such that
\begin{equation}\label{eq:global-bv-coaction}
 (\Psi_{P,N})_*\alpha
     =1\otimes\alpha+\theta\otimes\Delta_{P,N}\alpha.
\end{equation}
These operators commute with the normalized refinement maps
and the proper-trace pole maps. Under \eqref{eq:bm-limits},
\[
 \Delta[\alpha]=[\Delta_{P,N}\alpha],
 \qquad
 \Delta_{P,N}:
 \bigoplus_C H^{\mathrm{BM}}_{2\delta_{P,N}(C)-q}(C;\F)
 \longrightarrow
 \bigoplus_C H^{\mathrm{BM}}_{2\delta_{P,N}(C)-q+1}(C;\F),
\]
where $C$ ranges over the components of $Z_{P,N}$.

\item For a geometric computation, suppose additionally that
the source rotations descend to algebraic actions
\[
 a_{P,N}:\Gm\times Z_{P,N}\longrightarrow Z_{P,N},
 \qquad
 \pi_{P,N}(r_\lambda\gamma)
       =a_{P,N}(\lambda,\pi_{P,N}(\gamma)),
\]
where $\pi_{P,N}:F_P\mathscr LX\to Z_{P,N}$ is the
finite-factor projection. Require the affine-bundle refinement
maps and the closed pole inclusions to be equivariant.
Since $\Gm$ is connected, these actions preserve each connected
component and hence the locally constant volume indices:
$a_{P,N}^*\delta_{P,N}=\operatorname{pr}_{Z_{P,N}}^*\delta_{P,N}$.
Then $\Psi_{P,N}$ is the normalized oriented pullback along
$a_{P,N}$.

More explicitly, on a component whose analytic space is a
smooth complex manifold $Y$ of dimension $d$, with volume
index $\delta$, put $C_Y:=R\Gamma(Y,\omega_Y[-2\delta])$
and identify
\[
 H^q(C_Y)=H^{q+2d-2\delta}(Y;\F).
\]
Writing $a:G\times Y\to Y$ for the induced action, its
ordinary cohomological pullback satisfies
\begin{equation}\label{eq:global-bv-smooth-formula}
 a^*\alpha=1\otimes\alpha+\theta\otimes\Delta_{P,N}\alpha.
\end{equation}
Thus the coefficient of $\theta$ computes the BV operator.
\end{enumerate}
The additional equivariance in \textup{(ii)} is not part of
the existence of global finite factors. If source rotation does not descend to an algebraic action on a
chosen finite factor, \textup{(i)} still applies, but its
evaluation cannot use the geometric formula of \textup{(ii)}
on that factor.
\end{proposition}
\begin{proof}
The equivalences $B_{P,N}\simeq B_P$ identify the entire
pole-stage diagrams by
Proposition~\ref{prop:global-finite-presentations}(i).
Tensoring them with $\omega_G$ and transporting
$\mathcal R_P$ therefore defines $\mathcal R_{P,N}$ with
its pole and refinement comparisons. The unit of the rotation
action says that the restriction of $\mathcal R_{P,N}$ at
$1\in G$ is the identity. Since
$H^*(G;\F)=\F\oplus\F\theta$, the K\"unneth decomposition
of $(\Psi_{P,N})_*\alpha$ consequently has the form
\eqref{eq:global-bv-coaction}.
Evaluation of the corresponding parameter-chain map on the
counterclockwise fundamental cycle extracts its $\theta$
coefficient, by
Definition~\ref{def:parameter-chain-representative}.
Under $B_{P,N}\simeq B_P\simeq F_P\mathbf B$, the resulting
cohomology operator is the one induced by $\Delta_P$ of
Definition~\ref{def:source-rotation-operator}: equivalently,
the corresponding derived morphism is $\overline\Delta_P$.
This is a comparison of derived morphisms and their cohomology
operators, not a literal equality between independently chosen
finite-model and telescope cochain maps.
The sign is positive with the chain convention of
Definition~\ref{def:parameter-chain-representative}.
Indeed, put $S_G=C_{-*}(G;\F)$, with differential $+\partial$.
The comparison from usual scalar singular cochains, whose
differential is $f\mapsto f\circ\partial$, to
$\operatorname{Hom}_{\F}(S_G,\F)$ sends a degree-$j$ cochain
$f$ to $(-1)^{j(j+1)/2}f$. In degree $1$ this contributes
$-1$; tensor--Hom evaluation of the degree-$-1$ coefficient
map on a chain of degree $-1$ contributes another $-1$.
Their product is $+1$, so the positive angular coefficient
is $\Delta_{P,N}$.
Its compatibility with the pole and refinement maps gives
the direct-limit formula and the stated Borel--Moore degrees.

Under \textup{(ii)}, the algebraic isomorphism
\[
 A_{P,N}(\lambda,y)=(\lambda,a_{P,N}(\lambda,y)),
 \qquad
 A_{P,N}^{-1}(\lambda,y)
             =(\lambda,a_{P,N}(\lambda^{-1},y))
\]
presents the whole source-rotation family. Let
$p:\Gm\times Z_{P,N}\to Z_{P,N}$ be projection. Since
$a_{P,N}=p\circ A_{P,N}$, this is a smooth morphism of
complex relative dimension one. Using the analytifications of
these maps, its normalized oriented pullback is the composite
\[
 \begin{aligned}
 R\Gamma(Z_{P,N},\omega_{Z_{P,N}}[-2\delta_{P,N}])
 &\longrightarrow
 R\Gamma(G\times Z_{P,N},
     a_{P,N}^{-1}(\omega_{Z_{P,N}}[-2\delta_{P,N}]))\\
 &\simeq
 R\Gamma(G\times Z_{P,N},
     \omega_{G\times Z_{P,N}}[-2p^*\delta_{P,N}-2])\\
 &\simeq R\Gamma(G;\F)\otimes B_{P,N}.
 \end{aligned}
\]
The first arrow is ordinary inverse image on global sections;
the second is induced by the smooth-purity transformation for
$a_{P,N}$, with the displayed normalization shifts. The comparison for
$A_{P,N}$ is precisely the finite-factor comparison defining
$\mathcal R_P$, so this composite equals $\Psi_{P,N}$.
On $Y$, substitute $\omega_Y\simeq\F_Y[2d]$; purity
then identifies the composite with ordinary pullback $a^*$,
proving \eqref{eq:global-bv-smooth-formula}. All shifts used
here are even, so the same argument gives the parity formulas.
\end{proof}

\subsection{Example -- The multiplicative group}

Continue Subsection~\ref{ex:cstar-global-fusion} with $X=\Gm$ and
$\Omega=dx/x$.

\begin{proposition}[The BV operator of the multiplicative group]
\label{ex:cstar-global-bv}
\leavevmode
\begin{enumerate}[(i),wide=0pt]
\item For the winding classes of
Example~\ref{ex:cstar-loop-complex},
\[
 \Delta e_v=0,\qquad \Delta\eta_v=v e_v
 \qquad(v\in\Z).
\]
In the algebra presentation of
Proposition~\ref{prop:cstar-fusion-algebra}, with
$T=e_1$ and $\eta=\eta_0$, this means
\[
 \Delta\bigl(f(T)+g(T)\eta\bigr)=Tg'(T)
 \qquad\bigl(f,g\in\F[T,T^{-1}]\bigr),
\]
where $g'$ is ordinary Laurent-polynomial differentiation.
Equivalently, $\Delta=T\partial_T\partial_\eta$, where
$\partial_T$ is ordinary differentiation in $T$ and
$\partial_\eta(f(T)+g(T)\eta)=g(T)$.
It preserves the contact filtration of
Remark~\ref{rem:cstar-fusion-contact-filtration}.

\item For $P,N\ge0$ and $|v|\le P$, put
$Y_{v,N}=\Gm\times\A^N$, with coordinates
$a_0\in\Gm$ and $a_1,\ldots,a_N$.
Under the normalized finite-factor equivalence
\[
 B_{P,N}\simeq\bigoplus_{v=-P}^{P}
                   R\Gamma(Y_{v,N}(\C);\F)
\]
of Example~\ref{ex:cstar-loop-complex}, the coaction
$\Psi_{P,N}$ of Proposition~\ref{prop:global-finite-bv}(i)
is the direct sum of ordinary pullbacks by the actions
\begin{equation}\label{eq:cstar-global-rotation}
 a_{v,N}:\Gm\times Y_{v,N}\longrightarrow Y_{v,N},
 \qquad (\lambda,(a_j)_{j=0}^N)
              \longmapsto(\lambda^{v+j}a_j)_{j=0}^N.
\end{equation}
These identifications commute with normalized refinement maps
and proper-trace pole maps.
\end{enumerate}
\end{proposition}

\begin{proof}
\emph{Step 1: realize source rotation on the finite factors and prove (ii).}
Use the actual factors
$Z_{P,N}=\coprod_{v=-P}^{P}T_{P,v}\times Y_{v,N}$
of Example~\ref{ex:cstar-loop-complex}.
In the winding-$v$ coordinates of
\eqref{eq:cstar-full-sector-factor}, the monic polynomials
$q,p$ have degrees $P+v,P-v$, respectively, and
$a(t)=\sum_{j\ge0}a_jt^j$.
The substitution $t\mapsto\lambda t$ acts by
\[
 q(t)\longmapsto\lambda^{-(P+v)}q(\lambda t),\qquad
 p(t)\longmapsto\lambda^{-(P-v)}p(\lambda t),\qquad
 a(t)\longmapsto\lambda^v a(\lambda t).
\]
The transformed polynomials remain monic and their product
is $t^{2P}$; substituting them in
\eqref{eq:cstar-full-sector-factor} gives
$x(\lambda t),y(\lambda t)$.
Thus these formulas define algebraic actions on $Z_{P,N}$
for every $N$, including their nonreduced structure.
They commute with forgetting $a_{N+1}$ and with the pole
inclusions $(q,p,a)\mapsto(tq,tp,a)$.
Their restrictions to the reductions
$Y_{v,N}=(T_{P,v}\times Y_{v,N})_{\mathrm{red}}$
are precisely \eqref{eq:cstar-global-rotation}.

Example~\ref{ex:cstar-loop-complex} identifies the normalized
complexes with $R\Gamma(Y_{v,N}(\C);\F)$:
the reduction inclusions induce the identity on analytic
spaces, and $\delta_{P,N}=N+1=\dim Y_{v,N}$.
Proposition~\ref{prop:global-finite-bv}(ii) therefore identifies
the transported coaction with $a_{v,N}^*$.
On each forgotten affine line, multiplication by
$\lambda^{v+N+1}$ preserves the complex orientation, and
the pole maps retain the coordinates of every existing
winding component. Hence the normalized refinement and pole
comparisons commute with these actions over the whole
parameter scheme $\Gm$, proving (ii).

\emph{Step 2: extract the angular coefficient and prove (i).}
Use the class $\theta\in H^1(\Gm;\F)$ fixed in
Subsection~\ref{ex:cstar-global-fusion}; it evaluates to $1$
on the counterclockwise unit circle.
The classes $1$ and $a_0^*\theta$ on $Y_{v,N}(\C)$
represent $e_v$ and $\eta_v$, respectively, where
$a_0:Y_{v,N}\to\Gm$ is the leading-coefficient projection.
Since $a_0\circ a_{v,N}=\lambda^v a_0$,
\[
 a_{v,N}^*1=1\otimes1,\qquad
 a_{v,N}^*(a_0^*\theta)
       =1\otimes a_0^*\theta+v\theta\otimes1.
\]
The coefficient $v$ is the degree of
$S^1\to S^1$, $\lambda\mapsto\lambda^v$, including
$v\le0$. Part~(ii) and
Proposition~\ref{prop:global-finite-bv}(i), specifically
\eqref{eq:global-bv-coaction}, identify the positive
$\theta$ coefficient with $\Delta$.
This gives $\Delta e_v=0$ and $\Delta\eta_v=v e_v$.
Finally, Proposition~\ref{prop:cstar-fusion-algebra} gives
$e_v=T^v$ and $\eta_v=T^v\eta$; linearity yields the
formula for $f(T)+g(T)\eta$.
Both images remain in the winding-$v$ summand, so the
contact filtration is preserved. This proves (i).
\end{proof}

\subsection{The logarithmic BV operator}
\label{subsec:logarithmic-bv}

Fix the polarized SNC compactification $(M,D,A)$ and use the
logarithmic ring of Definition~\ref{def:log-cohomology}, in parity
grading or in the integer grading determined by $\Omega$.
For a realized contact vector $\vct$, put $J=I(\vct)$.
The normal-torus bundle $S_J^\circ\to D_J^\circ$ carries the
fiberwise action of $(S^1)^J$. Write
\[
 a_j:S^1\times S_J^\circ\longrightarrow S_J^\circ
 \qquad(j\in J)
\]
for rotation of its $j$th normal-circle factor. Let
$\theta\in H^1(S^1;\F)$ evaluate to $1$ on the
counterclockwise fundamental class. Define the degree-$-1$ map
$\iota_j:H^d(S_J^\circ;\F)\to H^{d-1}(S_J^\circ;\F)$ by
the K\"unneth decomposition
\[
 a_j^*\alpha=1\otimes\alpha+\theta\otimes\iota_j\alpha.
\]
The acting circle is the first factor. These maps are defined
on the circle bundles themselves; no trivialization or angular
differential form is chosen. Define
\begin{equation}\label{eq:bv-first-page}
 \Delta_{\log}(\alpha t^{\vct})
   :=\left(\sum_{j\in J}v_j\iota_j\alpha\right)t^{\vct}.
\end{equation}
The integers $v_j$ are viewed in $\F$. For $\vct=\mathbf0$
the sum is empty, so $\Delta_{\log}|_{H^*(X;\F)}=0$.

\begin{proposition}[Source rotation on the logarithmic ring]
\label{prop:logarithmic-bv}
The multiplicative identification of
Theorem~\ref{thm:geometric-fusion-ring},
\[
 \psi:H^*(\operatorname{gr}_F\mathbf B)
                   \xrightarrow{\sim}H^*_{\log}(M,D),
\]
intertwines the associated-graded rotation operator with
\eqref{eq:bv-first-page}:
\[
 \psi\circ H^*(\operatorname{gr}_F\Delta)
                         =\Delta_{\log}\circ\psi.
\]
Consequently $(H^*_{\log}(M,D),\Delta_{\log})$ is a BV algebra.
The operator has degree $-1$ and preserves each contact-vector
summand, hence has bidegree $(0,-1)$ for the bigrading of
Definition~\ref{def:log-cohomology}.
\end{proposition}
\begin{proof}
The contact-layer comparison \eqref{eq:loop-layer} is induced
by leading-normal-coefficient evaluation on contact-jet finite
factors, followed by normalized affine-fiber comparisons and
Lemma~\ref{lem:face-restriction}. In an SNC chart, write
\[
 x_j(\gamma(t))=t^{v_j}(c_j+O(t)),\qquad
 c_j\ne0,\qquad j\in J.
\]
Under $t\mapsto e^{i\phi}t$, the leading coefficient becomes
$e^{iv_j\phi}c_j$, and the center in $D_J^\circ$ is fixed.
At every chosen jet order, substitution multiplies the
coefficient of $t^\nu$ by $e^{i\nu\phi}$, and hence commutes
with truncation and leading evaluation. These are identities
over the entire parameter circle. The normalized purity
equivalences commute with these automorphisms; their relative
dimensions and volume indices are unchanged. On each nonzero
normal line, radial retraction to its unit circle, for a chosen
Hermitian metric, also commutes with multiplication by
$e^{iv_j\phi}$. The face-restriction comparison is induced by
the inclusion $S_J^\circ\hookrightarrow S_J^{\log}$, which
commutes with rotation of these normal rays: on $S_J^{\log}$,
rotation of an index $j\in J$ rotates that ray and fixes
the base point and all other rays.
Thus the contact-layer comparison identifies the action with
the normal-torus action restricted along
\[
 \lambda_{\vct}:S^1\longrightarrow(S^1)^J,
 \qquad e^{i\phi}\longmapsto(e^{iv_j\phi})_{j\in J}.
\]
In homology,
\[
 (\lambda_{\vct})_*[S^1]
            =\sum_{j\in J}v_j[S^1_j]
                 \quad\text{in }H_1((S^1)^J;\F).
\]
Slant product with this class, with the acting circle first,
is $\sum_jv_j\iota_j$. This proves the displayed intertwining
identity on every positive-contact summand. On contact zero
it follows from $\Delta_0\simeq0$ in
Lemma~\ref{lem:source-dilation}.

The filtered square-zero and unit homotopies of that lemma,
and the filtered seven-term homotopy of
Theorem~\ref{thm:loop-bv}, give the BV identities on
$H^*(\operatorname{gr}_F\mathbf B)$. The multiplicative map $\psi$
transports them to the logarithmic ring. Rotation preserves
the contact vector and its volume shift, while $\iota_j$
lowers coefficient degree by one; this proves the degree and
bidegree assertions in both gradings.
\end{proof}

Formula~\eqref{eq:bv-first-page} determines the associated-graded
operator. It is the first-page operator in
Theorem~\ref{thm:summary}; it does not by itself determine
lower-filtration terms of $\Delta$ on $LH$.

\section{The contact spectral sequence}\label{sec:spectral-sequence}
This section proves Theorems~\ref{thm:summary}
and~\ref{thm:intro-proportional}. Localization constructs the
contact spectral sequence, and the filtered fusion and rotation
maps make it converge strongly as a spectral sequence of unital
BV algebras. In the case of clopen
contact loci the spectral sequence degenerates at $E_1$, but we also identify the additive splitting of $B$, and give a criterion using units for this hypothesis.
The final subsection uses the integral contact class to constrain
the differentials and prove degeneration for proportional boundary
classes.

\subsection{Spectral sequences of BV algebras and convergence}
\label{subsec:bv-spectral-sequence-conventions}
\begin{definition}[Spectral sequences and multiplicative structures]
\label{def:multiplicative-spectral-sequence}
Let $\Gamma=\Z$ in integer grading and $\Gamma=\Z/2$ in
parity grading, and let $(p,q)\in\Z\times\Gamma$.
All Koszul signs are computed from total cohomological degree:
for $a\in E_r^{p,q}$ homogeneous, put $|a|=p+q\in\Gamma$
and read sign exponents modulo two.
A cohomological spectral sequence from $E_1$ consists of
bigraded $\F$-vector spaces $E_r$, differentials, and specified
bigraded isomorphisms
\[
 d_r:E_r^{p,q}\longrightarrow E_r^{p+r,q-r+1},\qquad
 d_r^2=0,\qquad
 \kappa_r:H(E_r,d_r)\xrightarrow{\sim}E_{r+1}
 \quad(r\ge1),
\]
where
\[
 H(E_r,d_r)^{p,q}
 =\frac{\ker(d_r:E_r^{p,q}\to E_r^{p+r,q-r+1})}
 {\operatorname{im}(d_r:E_r^{p-r,q+r-1}\to E_r^{p,q})}.
\]
It is \emph{unital and multiplicative} if every $E_r$ is a
unital associative graded-commutative algebra, with product
\[
 \mu_r:E_r^{p,q}\otimes E_r^{p',q'}
          \longrightarrow E_r^{p+p',q+q'},\qquad
 \mu_r(a\otimes b)=ab,
\]
and unit $1_r\in E_r^{0,0}$, such that
\[
 d_r(ab)=(d_ra)b+(-1)^{|a|}a(d_rb),\qquad
 \kappa_r([a][b])=\kappa_r([a])\kappa_r([b]),\qquad
 \kappa_r([1_r])=1_{r+1}.
\]
The last two identities concern $d_r$-cycles. The zero algebra, with unit zero, is allowed.
\end{definition}

\begin{definition}[Spectral sequences of unital BV algebras]
\label{def:bv-spectral-sequence}
A \emph{spectral sequence of unital BV algebras}, with
filtration-preserving BV operator, is a unital multiplicative
spectral sequence with page products $\mu_r$ as in
Definition~\ref{def:multiplicative-spectral-sequence}, using its total-degree
sign convention, together with maps
\[
 \Delta_r:E_r^{p,q}\longrightarrow E_r^{p,q-1}
 \quad(r\ge1)
\]
such that $(E_r,\mu_r,\Delta_r)$ is a BV algebra in the sense
of Definition~\ref{def:loop-bv}, and
\begin{equation}\label{eq:bv-spectral-page-compatibility}
 d_r\Delta_r+\Delta_rd_r=0,\qquad
 \kappa_r\circ H(\Delta_r)=\Delta_{r+1}\circ\kappa_r.
\end{equation}
Thus each $(E_r,d_r,\mu_r,\Delta_r)$ is a differential BV
algebra and each $\kappa_r$ is an isomorphism of BV algebras.
A morphism consists of bidegree-zero maps $f_r:E_r\to E'_r$
preserving product, unit, and BV operator, with
\[
 f_rd_r=d'_rf_r,\qquad
 f_{r+1}\kappa_r=\kappa'_rH(f_r).
\]
\end{definition}

The terminology refers to spectral sequences of vector spaces
with these compatible operations: $d_r$ is a derivation, not
an algebra homomorphism. For this pagewise BV compatibility,
with a different BV bidegree, see
\cite[Theorem~4.3]{KuribayashiBVSS}.
The cohomology operations used in the definition are well defined:
$\ker d_r$ is a unital subalgebra, $\operatorname{im}d_r$ is
an ideal in it, and
\[
 \Delta_r(\ker d_r)\subset\ker d_r,\qquad
 \Delta_r(d_rx)=-d_r(\Delta_rx).
\]
Consequently the BV identities descend to $H(E_r,d_r)$.
For the bracket of \eqref{eq:general-rotation-bracket}, expanding
with the Leibniz rule and $d_r\Delta_r=-\Delta_rd_r$ gives
\begin{equation}\label{eq:bv-spectral-bracket-differential}
 d_r\{a,b\}_r
   =\{d_ra,b\}_r+(-1)^{|a|-1}\{a,d_rb\}_r.
\end{equation}
Thus the brackets descend to cohomology and are preserved by
$\kappa_r$ and by morphisms of BV spectral sequences; these are
consequences of the preceding definition, not additional axioms.

\begin{definition}[Strong convergence in the bounded-below setting]
\label{def:strong-bv-convergence}
Suppose that, for each $(p,q)$, both arrows incident to
$E_r^{p,q}$,
\[
 E_r^{p-r,q+r-1}\xrightarrow{d_r}E_r^{p,q}
          \xrightarrow{d_r}E_r^{p+r,q-r+1},
\]
vanish for all sufficiently large $r$. The maps $\kappa_r$
then canonically identify these terms for large $r$; write
$E_\infty^{p,q}$ for their stable value.

For a graded vector space $H$ with an increasing filtration
$F_\bullet H$, $F_PH=0$ for $P<0$, we say that the spectral
sequence \emph{converges strongly to $(H,F)$} if
$H=\bigcup_{P\ge0}F_PH$ and it is equipped with isomorphisms
\[
 E_\infty^{p,q}\xrightarrow{\sim}
       F_{-p}H^{p+q}/F_{-p-1}H^{p+q}.
\]
The pole filtration index $P\ge0$ corresponds to the spectral-sequence
column $p=-P$. 

If $H$ is a filtered unital BV algebra, convergence \emph{as
BV algebras} additionally requires that the displayed
$E_\infty\cong\gr_FH$ preserve product, unit, and BV operator.
The stable operations on $E_\infty$ are obtained by choosing a
page on which the relevant input and output bidegrees have
stabilized. Convergence as algebras is defined by omitting the
BV operators.
\end{definition}

\begin{definition}[Degeneration]
\label{def:bv-spectral-degeneration}
For $r_0\ge1$, a spectral sequence \emph{degenerates at
$E_{r_0}$} if
\[
 d_r=0\qquad\text{for every }r\ge r_0.
\]
The transition isomorphisms $\kappa_r$ then canonically identify
$E_{r_0}$ with every later page and with $E_\infty$.
For a spectral sequence of unital BV algebras, these are
isomorphisms of bigraded unital BV algebras. If it converges
strongly to $(H,F)$ as BV algebras in the sense of
Definition~\ref{def:strong-bv-convergence}, they give
\[
 E_{r_0}\xrightarrow{\sim}E_\infty
             \xrightarrow{\sim}\gr_FH
\]
as bigraded unital BV algebras.
\end{definition}


For the loop complex, the BV
identities on $E_1$ have already been established in
Proposition~\ref{prop:logarithmic-bv}; the proof below constructs
the compatible differentials and page transitions.

\subsection{The contact BV spectral sequence}
\label{subsec:contact-bv-spectral-sequence}
We follow the grading convention of
Subsection~\ref{subsec:relative-pole-recovery}: in integer grading,
$B_P$, $B$, $\mathbf B$, and $F_P\mathbf B$ stand for
$B_{\Omega,P}$, $B_\Omega$, $\mathbf B_\Omega$, and
$F_P\mathbf B_\Omega$, respectively, with $\Omega$ the chosen
nowhere-zero algebraic volume form. Without $\Omega$, these
symbols denote the parity-graded objects.

Fix a polarized SNC compactification $(M,D,A)$ of $X$. Use the
filtered subcomplex model $(\mathbf B,d_{\mathbf B},F_\bullet)$ of
Theorem~\ref{thm:canonical-complex}, with
\[
 F_P\mathbf B=0\quad(P<0),\qquad \mathbf B=\bigcup_{P\ge0}F_P\mathbf B.
\]
Fix the filtered cochain maps $m_2$ and $\Delta$ of
Definition~\ref{def:loop-fusion-product} and
Definition~\ref{def:source-rotation-operator}, together with
their comparisons, in the corresponding grading.
The abutment is equipped with the filtration of
Definition~\ref{def:loop-cohomology}:
\[
 F_PLH^n=\operatorname{im}\bigl(H^n(F_P\mathbf B)\longrightarrow H^n(\mathbf B)\bigr).
\]

\begin{proof}[Proof of Theorem~\ref{thm:summary}]
\emph{Step 1: construction and the first page.}
Put $L_P=\operatorname{cofib}(B_{P-1}\to B_P)$. The specified
comparison of pole-stage diagrams in
Theorem~\ref{thm:canonical-complex} identifies
$L_P\simeq F_P\mathbf B/F_{P-1}\mathbf B$ and the localization
triangles
\begin{equation}\label{eq:contact-localization-triangles}
 B_{P-1}\longrightarrow B_P\longrightarrow L_P
       \xrightarrow{\partial_P}B_{P-1}[1]
\end{equation}
with those of the filtered subcomplex model. The associated
spectral sequence therefore has
\[
 E_1^{p,q}=H^{p+q}(L_{-p}),\qquad
 d_r:E_r^{p,q}\longrightarrow E_r^{p+r,q-r+1},\qquad
 \kappa_r:H(E_r,d_r)\xrightarrow{\sim}E_{r+1};
\]
see \cite[Tag~012M]{Stacks}. The contact-layer equivalence
\eqref{eq:loop-layer} and Definition~\ref{def:log-cohomology} give
\[
 E_1^{-P,n+P}=H^n(L_P)\cong H_{\log}^{-P,n+P}(M,D).
\]

\emph{Step 2: the BV operations on the pages.}
With the sum filtration of
Lemma~\ref{lem:filtered-operations-colimit}, the filtered
pairing $m_2$ and the filtered cochain map
$\Delta:\mathbf B\to\mathbf B[-1]$ induce
\[
 \mu_r:E_r^{p,q}\otimes E_r^{p',q'}
       \longrightarrow E_r^{p+p',q+q'},\qquad
 \Delta_r:E_r^{p,q}\longrightarrow E_r^{p,q-1}.
\]
The filtered-pairing construction and functoriality give
\[
 d_r\mu_r(a,b)=\mu_r(d_ra,b)+(-1)^{|a|}\mu_r(a,d_rb),
 \qquad d_r\Delta_r=-\Delta_rd_r,
\]
where $a,b$ are homogeneous and $|a|$ is total degree.
The maps $\kappa_r$ identify the operations induced on
$H(E_r,d_r)$ with $\mu_{r+1}$ and $\Delta_{r+1}$;
see \cite[Section~4]{MaySpectral} and
\cite[Tag~012O]{Stacks}. The minus sign comes from
$d_{\mathbf B[-1]}=-d_{\mathbf B}$. These constructions apply
in parity grading with $q$ and total degree read modulo two
and the same Koszul signs.

Theorem~\ref{thm:geometric-fusion-ring} and
Proposition~\ref{prop:logarithmic-bv} identify
$(E_1,\mu_1,\Delta_1)$ with the logarithmic BV algebra.
Its algebra and BV identities pass inductively to $E_{r+1}$ through
$\kappa_r$. Choose a closed $e\in(F_0\mathbf B)^0$
representing $1\in H^0(X;\F)$. Its page classes give the
units, with $1_{r+1}=\kappa_r[1_r]$; the zero algebra is
allowed as in Definition~\ref{def:multiplicative-spectral-sequence}.
This proves the asserted unital BV structure on every page.

\emph{Step 3: stabilization and the abutment.}
Write $d=d_{\mathbf B}$ and put
$Z_P^n=\{x\in(F_P\mathbf B)^n:dx=0\}$.
A representative $x\in(F_P\mathbf B)^n$ on page $E_r$
satisfies $dx\in(F_{P-r}\mathbf B)^{n+1}$. For $r>P\ge0$
this last space is zero. The cycle and boundary formulas
\cite[Tag~012K]{Stacks} therefore give
\begin{equation}\label{eq:contact-stable-cycle-quotient}
 E_r^{-P,n+P}\cong
 \frac{Z_P^n}
 {Z_{P-1}^n+
  \bigl(Z_P^n\cap d((F_{P+r-1}\mathbf B)^{n-1})\bigr)}.
\end{equation}
By Step~1, each weight-$P$ column of $E_1$ is a finite sum
of finite-dimensional cohomology groups. For $r>P$ the
numerator in \eqref{eq:contact-stable-cycle-quotient} is fixed,
and the denominators have increasing images in the finite-dimensional
space $E_{P+1}^{-P,n+P}$. Hence the quotient stabilizes.
Since $\mathbf B=\bigcup_R F_R\mathbf B$, every primitive
of a boundary belongs to some finite pole stage. The stable
quotient is therefore identified by
\[
 E_\infty^{-P,n+P}\xrightarrow{\sim}F_PLH^n/F_{P-1}LH^n,
 \qquad x\longmapsto[x]\bmod F_{P-1}LH^n.
\]
Exhaustiveness also gives $LH^n=\bigcup_PF_PLH^n$, proving
strong convergence in both gradings. Under the displayed
isomorphism, the product, BV operator, and unit on both sides
are represented by $m_2(x,y)$, $\Delta x$, and the unit
cocycle $e$ of Step~2, respectively.
\end{proof}

\begin{proposition}[Independence of finite-factor presentations and cochain choices]
\label{prop:contact-bv-spectral-independence}
Fix a polarized SNC compactification $(M,D,A)$ and a grading,
including the volume form $\Omega$ in integer grading.
For two chosen models $\mathbf B,\mathbf B'$ of the same pole-stage diagram,
common refinements of finite-factor presentations and the specified cochain
comparisons induce canonical isomorphisms
\[
 g_r:E_r(\mathbf B)\xrightarrow{\sim}E_r(\mathbf B')\qquad(r\ge1)
\]
of the BV spectral sequences of Theorem~\ref{thm:summary}.
They respect composition and are unchanged by further common
refinements or changes of cochain representatives of the specified
comparisons. No comparison of filtrations from different
compactifications is asserted.
\end{proposition}
\begin{proof}
The comparisons from Sections~\ref{sec:loop-complex},
\ref{sec:fusion-construction}, and~\ref{sec:bv} are zigzags of
stagewise quasi-isomorphisms, with specified filtered homotopies
comparing product, unit, and rotation. For a comparison
$g:\mathbf B\to \mathbf B'$, the localization triangles
\eqref{eq:contact-localization-triangles} give a morphism of
exact couples and hence maps $g_r$ commuting with $d_r$ and
$\kappa_r$. Since every $g_P:F_P\mathbf B\to F_P\mathbf B'$ is a
quasi-isomorphism, so is the induced map
$\gr_P\mathbf B\to\gr_P\mathbf B'$. The specified comparisons
identify these quotients with the localization cofibers, so $g_1$ is an
isomorphism; the operation comparisons make it a unital
BV-algebra isomorphism. The recursion
\[
 g_{r+1}=\kappa'_r\circ H(g_r)\circ\kappa_r^{-1}
\]
proves the same on all later pages. Inverting these isomorphisms
treats zigzags. The specified composition and refinement
homotopies give equal maps on $E_1$ for equivalent comparisons;
the recursion then gives equality on every page. This proves
both the claimed independence and compatibility with composition.
\end{proof}

\subsection{Separated contact loci and logarithmic models}
\label{sec:split-models}
\label{subsec:separated-sectors}

Write
\[
 \mathcal V:=\{\mathbf0\}\cup
 \{\vct\in\N^I\setminus\{\mathbf0\}:
                     D_{I(\vct)}^\circ\ne\varnothing\}
\]
for the realized contact vectors, where
$I(\vct)=\{i\in I:v_i>0\}$.

\begin{corollary}[Separated contact loci]\label{thm:separated-sectors}
Suppose every $\mathscr P_{\vct}$ is open and closed in
$\mathscr LX(\C)$ in the analytic topology of
Convention~\ref{conv:loop-analytic-topology}.
Given a nowhere-zero algebraic volume form $\Omega$ on $X$:
\begin{enumerate}
\item There are canonical equivalences
\[
 B_{\Omega,P}\simeq
       \bigoplus_{\substack{\vct\in\mathcal V\\w(\vct)\le P}}
       C^{*-\mu_\Omega(\vct)}(S_{I(\vct)}^{\log};\F),
\]
compatible with the pole-stage maps and the summand inclusions
for $P\le Q$. Thus they give an equivalence of filtered objects
in $\operatorname{Fun}(N(\Z_{\ge0}),\mathcal D(\F))$.
\item There is a canonical isomorphism of filtered graded
vector spaces
\[
 LH^*_\Omega(X)\cong H^*_{\log}(M,D),
\]
where logarithmic cohomology has the integer grading of
Definition~\ref{def:log-cohomology} and its weight filtration.
\end{enumerate}
Without $\Omega$, the first assertion uses the summands
$\operatorname{Per}_2C^*(S_{I(\vct)}^{\log};\F)$ in the
periodic category, and the second gives
$LH^*(X)\cong H^*_{\log}(M,D)$ with the parity grading of
Definition~\ref{def:log-cohomology}. In either grading, the
contact BV spectral sequence degenerates at $E_1$ in the sense
of Definition~\ref{def:bv-spectral-degeneration}.
\end{corollary}
\begin{proof}
Put $S_P=F_P\mathscr LX(\C)$ and let
$j_{\vct,P}:\mathscr P_{\vct}\hookrightarrow S_P$ denote
the inclusion. The finite clopen decomposition
$S_P=\coprod_{w(\vct)\le P}\mathscr P_{\vct}$ gives
\[
 \omega_{P,\Omega}
 \simeq\bigoplus_{w(\vct)\le P}
      (j_{\vct,P})_*j_{\vct,P}^{-1}\omega_{P,\Omega},
 \qquad
 B_{\Omega,P}
 \simeq\bigoplus_{w(\vct)\le P}
      R\Gamma(\mathscr P_{\vct},
                    j_{\vct,P}^{-1}\omega_{P,\Omega}).
\]
For $P\le Q$, the proper-trace arrow of
Proposition~\ref{prop:normalized-pole-sheaves} restricts on
each old clopen summand to the trace for its identity map,
whose adjunction counit is an equivalence. These comparisons
preserve the composition compatibilities of that proposition.
The contact calculation in Theorem~\ref{thm:canonical-complex}
identifies each summand with
$C^{*-\mu_\Omega(\vct)}(S_{I(\vct)}^{\log};\F)$,
independently of $P\ge w(\vct)$. Thus the pole-stage diagram
is the finite-sum diagram with summand inclusions. Its homotopy
colimit is the direct sum over $\mathcal V$, giving the asserted
filtered cohomology isomorphism and degeneration. These are
comparisons of derived diagrams; they do not assert a literal
direct-sum decomposition of the chosen subcomplexes
$F_P\mathbf B_\Omega$. The periodic argument is identical.
\end{proof}

\begin{remark}[Leading products and lower-weight corrections]
In the volume grading write
\[
 H_{\vct}^k:=
 H^{k-\mu_\Omega(\vct)}(S_{I(\vct)}^{\log};\F)t^{\vct},
\]
and denote the maps of the preceding splitting by
\[
 \iota_{\vct}:H_{\vct}^*\longrightarrow LH^*_\Omega(X),
 \qquad
 \operatorname{pr}_{\uct}:LH^*_\Omega(X)\longrightarrow H_{\uct}^*.
\]
Fusion of the $\vct$ and $\wct$ summands has weight at most
$w(\vct)+w(\wct)$; Lemma~\ref{lem:fusion-leading-product}
identifies its component of that weight with the logarithmic
product \eqref{eq:end-product}. Consequently the additive
identification in Corollary~\ref{thm:separated-sectors}(2) is a
ring identification exactly when, for every
$a\in H_{\vct}^*$ and $b\in H_{\wct}^*$,
\begin{equation}\label{eq:separated-no-lower-products}
 \operatorname{pr}_{\uct}\mu(\iota_{\vct}a,\iota_{\wct}b)=0
 \quad\text{whenever}\quad
 w(\uct)<w(\vct)+w(\wct).
\end{equation}
The same criterion holds in parity, using the contact summands
and grading of Definition~\ref{def:log-cohomology}.
Clopenness alone does not imply this vanishing:
in Subsection~\ref{ex:cstar-global-fusion}, the winding loci
of $\Gm$ are open and closed, but
$e_1e_{-1}=e_0=1$ has contact weight $0$, whereas
the sum of the input contact weights is $2$.
\end{remark}

\begin{remark}[Finite contact support and ordinary cohomology]
\label{sec:normal-loop-models}
\label{lp:prop:clopen-log-product}
Give $\mathscr LX(\C)$ the analytic topology of
Convention~\ref{conv:loop-analytic-topology}.
Under the clopenness hypothesis, this space is the topological
coproduct of its contact loci. For $R=\Z$ or $R=\F$,
ordinary singular cohomology therefore satisfies
\[
 H^q(\mathscr LX(\C);R)
 \cong\prod_{\vct\in\mathcal V}H^q(\mathscr P_{\vct};R)
 \cong\prod_{\vct\in\mathcal V}
               H^q(S_{I(\vct)}^{\log};R).
\]
Indeed each singular simplex lies in one summand, so cochains
form a product; products of $R$-modules are exact. The second
identification uses Proposition~\ref{prop:loop-sectors} and
Lemma~\ref{lem:face-restriction}. The subgroup of finite contact
support is
\[
 H^q_{\fin}(\mathscr LX;R)
 :=\bigoplus_{\vct\in\mathcal V}H^q(\mathscr P_{\vct};R)
 \cong\bigoplus_{\vct\in\mathcal V}
               H^q(S_{I(\vct)}^{\log};R)
 \subseteq H^q(\mathscr LX(\C);R).
\]
Here $q$ is the ordinary cohomological degree. For $R=\F$,
Corollary~\ref{thm:separated-sectors}(2) identifies this
finite-support group with loop cohomology after regrading each
$\vct$ summand from degree $q$ to $q+\mu_\Omega(\vct)$,
or to parity $\bar q$ without $\Omega$. The common-face product
\eqref{eq:end-product}, with coefficients $R$, transports to
$H^*_{\fin}(\mathscr LX;R)$ and gives its logarithmic product.
This differs from ordinary cup product, which is componentwise
in the displayed product. For example, for $X=\Gm$ in
Subsection~\ref{ex:cstar-global-fusion}, the classes $e_v$
correspond additively to characteristic functions of the winding
components, so $e_1\smile e_{-1}=0$, whereas
$\mu(e_1,e_{-1})=e_0$.
\end{remark}

\begin{proposition}[Valuations of units]
\label{lp:prop:unit-valuation}
\label{subsec:unit-separation}
Set $U(X)=\mathcal O(X)^*/\C^*$ and, for a unit $u$, define
$\nu_u(\gamma)=\ord_t\gamma^*u$. Give
$\operatorname{Hom}(U(X),\Z)$ the subspace topology from
$\Z^{U(X)}$, with $\Z$ discrete. Then
\[
 \nu:\mathscr LX(\C)\longrightarrow\operatorname{Hom}(U(X),\Z),
 \qquad \nu_\gamma([u])=\nu_u(\gamma),
\]
is continuous. For $q\ge1$, source substitution
$c_q(\gamma)(t)=\gamma(t^q)$ is continuous and satisfies
$\nu_{c_q(\gamma)}=q\nu_\gamma$. On $\mathscr P_{\vct}$,
\begin{equation}\label{lp:eq:unit-matrix}
 \nu_u(\gamma)=\sum_{i\in I}(\ord_{D_i}u)v_i.
\end{equation}
In particular, each $\nu_u$ is constant on connected components
of $\mathscr LX(\C)$ in the analytic topology.
\end{proposition}
\begin{proof}
A unit $u$ defines a morphism $X\to\Gm$ and hence a continuous
map of formal loop spaces. The winding loci of $\mathscr L\Gm$
are open and closed by Example~\ref{ex:cstar-loop-complex}, so
$\nu_u$ is continuous. Coordinatewise continuity gives continuity
of $\nu$ in the stated topology; additivity of order under
multiplication makes its values homomorphisms on $U(X)$.
Substitution $t\mapsto t^q$ is coefficientwise continuous,
sends a Laurent pole bound $b$ to $qb$, and multiplies orders
by $q$. Finally, $\operatorname{div}(u)$ is supported on $D$.
In an SNC chart about the center of the extended arc, write
$u=\varepsilon\prod_i x_i^{\ord_{D_i}u}$ with $\varepsilon$
a regular unit. Taking orders gives \eqref{lp:eq:unit-matrix}.
\end{proof}

\begin{lemma}[Separation by units]\label{lem:unit-separation}
Recall the set $\mathcal V$ of realized contact vectors,
including $\mathbf0$, defined above. For units
$u_1,\ldots,u_r\in\mathcal O(X)^*$, define the integer matrix
\[
 W=(W_{ai})\in\operatorname{Mat}_{r\times |I|}(\Z),
 \qquad W_{ai}:=\ord_{D_i}(u_a).
\]
If $W|_{\mathcal V}:\mathcal V\to\Z^r$ is injective, every
contact locus is open and closed in $\mathscr LX(\C)$, and
Corollary~\ref{thm:separated-sectors}(1)--(2) applies.
\end{lemma}
\begin{proof}
The map $\nu_W=(\nu_{u_1},\ldots,\nu_{u_r})$ is continuous
to the discrete space $\Z^r$. Properness of $M$ assigns a unique
realized contact vector to every loop, and
\eqref{lp:eq:unit-matrix} and injectivity give
$\mathscr P_{\vct}=\nu_W^{-1}(W\vct)$. Thus each contact locus
is clopen; no connectedness assertion is needed.
\end{proof}

\subsection{Degeneration for proportional boundary classes}
\label{sec:proportional-degeneration}

Fix a polarized SNC compactification $(M,D,A)$, with
$D=\bigcup_{i=1}^rD_i$, $r\ge1$, and
$A=\sum_i k_iD_i$, $k_i>0$. Suppose that a line bundle $L$
and positive integers $n_i$ satisfy
\begin{equation}\label{eq:proportional-boundary}
 \mathcal O_M(D_i)\cong L^{\otimes n_i},\qquad
 \mathbf n=(n_1,\ldots,n_r).
\end{equation}
This is the hypothesis $D_i\sim n_iH$ of
Theorem~\ref{thm:intro-proportional}, with $L=\mathcal O_M(H)$.
Fix these isomorphisms and let $s_i\in H^0(M,L^{\otimes n_i})$
be the section with divisor $D_i$. Put
\[
 \Lambda:=\Z^r/\Z\mathbf n,\qquad
 N_A:=\sum_i k_in_i>0.
\]
Since $L^{\otimes N_A}\cong\mathcal O_M(A)$, $L$ is ample.
The contact-label group $\Lambda$ is integral, with the discrete
topology; it is not tensored with $\F$. Its torsion subgroup has
order $\gcd(n_1,\ldots,n_r)$.

For $\gamma\in X(\C(\!( t )\!))$, properness of $M$ gives its unique
extension $\bar\gamma\in M(\C\llbracket t \rrbracket)$, hence a contact vector
$\vct(\gamma)
=\bigl(\operatorname{ord}_{\bar\gamma}D_i\bigr)_{i=1}^r$. Define
\begin{equation}\label{eq:integral-contact-class}
 \lambda:\mathscr LX(\C)\longrightarrow\Lambda,\qquad
 \lambda(\gamma)=[\vct(\gamma)].
\end{equation}

\begin{lemma}[Integral contact classes on pole stages]
\label{lem:integral-contact-class}
For $P\ge0$, let $S_P=F_P\mathscr LX(\C)$, with our analytic
topology.
The map $\lambda$ is locally constant on $S_P$ and on
$\mathscr LX(\C)$. Each $S_P$ has only finitely many nonempty
label fibers $S_{P,\xi}:=\lambda^{-1}(\xi)\cap S_P$.

In parity grading, recall the pole-stage complex and define its
label terms by
\[
 B_P=R\Gamma(S_P,\omega_P^{\mathrm{per}}),\qquad
 B_{P,\xi}:=R\Gamma(S_{P,\xi},
                    \omega_P^{\mathrm{per}}|_{S_{P,\xi}})
       \quad\text{in }\mathcal D_{\mathrm{per}}(\F).
\]
Given a nowhere-zero algebraic volume form $\Omega$ on $X$, the
corresponding integer-graded objects are
\[
 B_{\Omega,P}=R\Gamma(S_P,\omega_{P,\Omega}),
 \qquad
 B_{\Omega,P,\xi}:=R\Gamma(S_{P,\xi},
                         \omega_{P,\Omega}|_{S_{P,\xi}})
       \quad\text{in }\mathcal D(\F).
\]
Here the sheaves are the normalized dualizing complexes of
Proposition~\ref{prop:normalized-pole-sheaves}. There are canonical
finite-sum equivalences
\[
 B_P\simeq\bigoplus_{\xi\in\Lambda}B_{P,\xi},\qquad
 B_{\Omega,P}\simeq\bigoplus_{\xi\in\Lambda}B_{\Omega,P,\xi}.
\]
They commute coherently with pole-stage maps, localization
triangles, and rotation. In particular they give equivalences of
filtered diagrams in
$\operatorname{Fun}(N(\Z_{\ge0}),\mathcal D_{\mathrm{per}}(\F))$
and, when $\Omega$ is chosen,
$\operatorname{Fun}(N(\Z_{\ge0}),\mathcal D(\F))$, respectively, and
induce
\[
 B\simeq\bigoplus_{\xi\in\Lambda}\hocolim_P B_{P,\xi},
 \qquad
 B_\Omega\simeq
     \bigoplus_{\xi\in\Lambda}\hocolim_P B_{\Omega,P,\xi}.
\]
Consequently all contact spectral-sequence differentials
preserve the integral label $\lambda$.
\end{lemma}
\begin{proof}
\emph{Specialization of a one-parameter family of loops.}
Let $R=\C\llbracket s \rrbracket$, $K=\C(\!( s )\!)$, and consider a 1-parameter family of loops
$\gamma\in X(R(\!( t )\!))$. The ring
\[
 R(\!( t )\!)=\C\llbracket s,t \rrbracket[t^{-1}]
\]
is factorial, being a localization of the regular local
factorial ring $\C\llbracket s,t \rrbracket$, so its Picard group is zero.
Choose a basis $e$ of $\gamma^*L$ and write
\[
 \gamma^*s_i=a_i e^{\otimes n_i},\qquad a_i\in R(\!( t )\!)^\times.
\]
The coefficients are units because $\gamma$ lands in $X$.
Every unit of $\C\llbracket s,t \rrbracket[t^{-1}]$ has the form $t^m u$,
with $m\in\Z$ and $u\in\C\llbracket s,t \rrbracket^\times$: its divisor in
the factorial ring can be supported only on the prime $t$.
Thus $a_i=t^{m_i}u_i$, and its $t$-order is $m_i$ both over
$K$ and after specializing $s=0$.

For $F=K$ or $\C$, let $\gamma_F\in X(F(\!( t )\!))$
be induced by the coefficient homomorphism $R\to F$. Hence $\gamma_K$ and $\gamma_\C$ are, respectively, the generic loop and the special loop in the family $\gamma$.
Properness of $M$ gives a unique extension
$\bar\gamma_F:\Spec F\llbracket t \rrbracket\to M$ of $\gamma_F$.
The line bundle $\bar\gamma_F^*L$ on $\Spec F\llbracket t \rrbracket$
is free of rank one. Choose a generator $\widetilde e_F$.
On the punctured disk, write
$e_F=c_F\widetilde e_F$ with $c_F\in F(\!( t )\!)^\times$,
and put $b_F=\ord_t(c_F)$.
 Taking orders gives
\begin{equation}\label{eq:contact-class-specialization}
 v_i(\gamma_F)=m_i+n_i b_F,\qquad
 \vct(\gamma_K)-\vct(\gamma_\C)
                  =(b_K-b_\C)\mathbf n.
\end{equation}
The two fibers therefore have the same class in $\Lambda$.
This argument does not require an extension of the whole
family over $\operatorname{Spec}\C\llbracket s,t \rrbracket$.

\emph{Finite constructible contact conditions.}
On $S_P$ there are finitely many possible vectors, since
$0\le v_i\le\lfloor P/k_i\rfloor$.
We verify that their loci pull back to finite constructible
coefficient cylinders on every nice chart
$U=Y\times\A^J\to F_P\mathscr LX$.
Choose $d>0$ with $dA$ very ample and projective coordinates
as in Proposition~\ref{prop:loop-sectors}, writing
$g_j=s'_j/s'_0$ on $X$.
Stratify the Laurent stage by the minimum order of
$1,g_1,\ldots,g_m$ and a coordinate attaining that order.
The order lies in $[-dP,0]$; its value and the chosen nonzero
leading coefficient impose finitely many coefficient conditions.
On each piece, multiplying the projective coordinates by the
corresponding power of $t$ and dividing by the chosen unit
power series constructs the extended arc. Each coefficient of
that arc is a regular function of finitely many Laurent
coefficients with the chosen leading coefficient inverted.

Choose a finite affine refinement of these projective charts on
which all $D_i$ are principal. Membership of the center in these
opens and the coefficients through order $\lfloor P/k_i\rfloor$
of the local equation of $D_i$ are again finite coefficient
conditions, with finitely many denominators inverted. They
determine every $v_i$. These constructions work over every
residue field, so they define constructible contact strata,
not only partitions of complex points. On $U$, retain in its
finite factor the finitely many equations defining these pieces
and the numerators and denominators of their coefficient
functions, as in
Subsection~\ref{subsec:volume-index-preparations}. After enlarging
$Y$ and relabeling the remaining tail, the entire finite
partition is pulled back from a constructible partition of $Y$.

\emph{Grouping by the integral class.}
Let $C_\xi\subset Y$ be the union of the constructible pieces
whose contact vectors represent $\xi$. We show that each
$C_\xi$ is closed. Otherwise its closure difference contains a
closed complex point $y$. Algebraic curve selection for
constructible subsets of the finite-type scheme $Y$ gives a
curve through $y$ with generic point in $C_\xi$. Normalizing
the curve and completing at a point over $y$ gives
$\operatorname{Spec}\C\llbracket s \rrbracket\to Y$, with generic point in
$C_\xi$ and special point outside it. The zero section lifts this morphism to $Y\times\A^J$.
The chart map then gives an element of $X(R(\!( t )\!))$ whose
generic and special contact classes differ. The residue-field
description above identifies these classes, contradicting
\eqref{eq:contact-class-specialization}. Thus every $C_\xi$ is closed;
the partition is finite, so every $C_\xi$ is also open.

The nice \'etale charts form a jointly surjective family of
local homeomorphisms on complex points, by \emph{Analytic realization}
in Subsection~\ref{subsec:coherent-descent}. Their clopen
label partitions therefore descend to $S_P$. They agree under
the pole-stage inclusions. The final analytic topology of
$\mathscr LX(\C)$ makes $\lambda$ locally constant there as well.

Restriction to the finite clopen decomposition of $S_P$ gives
the asserted equivalences of complexes in each coefficient
category. The proper-trace maps and localization triangles
restrict to the corresponding label fibers and preserve their
composition coherences. Rotation preserves each contact vector,
so the rotation action and its induced operator preserve the
label decomposition. Taking homotopy colimits gives the
displayed decompositions of $B$ and $B_\Omega$, since direct
sums commute with these colimits. The decomposed localization
exact couples then give the assertion about every differential.
\end{proof}

\begin{remark}[Restrictions on differentials]
\label{rem:integral-contact-differentials}
Lemma~\ref{lem:integral-contact-class} gives, on every page
$E_\ell$ of the contact spectral sequence, a decomposition
$E_\ell=\bigoplus_{\xi\in\Lambda}E_{\ell,\xi}$
preserved by $d_\ell$.
Its $\xi$ summand is supported only
at weights of realized vectors $\vct$ with $[\vct]=\xi$,
because it is a subquotient of the corresponding first-page
summand. This does not supply a decomposition by individual
contact vectors on arbitrary later pages.
If $d_1,\ldots,d_{\ell-1}$ vanish, the specified page
isomorphisms identify $E_\ell$ canonically with $E_1$.
A nonzero component of $d_\ell$ from contact $\vct$ to
contact $\vct'$ must then satisfy
\[
 \vct-\vct'=a\mathbf n,\qquad
 \ell=w(\vct)-w(\vct')=aN_A,\qquad a\in\Z_{>0}.
\]
The first equality follows from label preservation, and the
second from the bidegree of $d_\ell$. Since $N_A>0$ and
$\ell>0$, the integer $a$ is positive.
\end{remark}

\begin{proof}[Proof of Theorem~\ref{thm:intro-proportional}]
\phantomsection\label{proof:intro-proportional}
Let $\mathbf e_1,\ldots,\mathbf e_r$ be the standard basis
of $\Z^r$. For $I\subseteq\{1,\ldots,r\}$, put
$\mathbf1_I:=\sum_{i\in I}\mathbf e_i$. The first-page algebra is
generated by $H^*(X)$ and the classes
$\alpha t^{\mathbf1_I}$ for nonempty realized supports $I$
and $\alpha\in H^*(S_I^\circ;\F)$. Indeed, for any realized
$\vct\ne\mathbf0$ with $I=I(\vct)$, the common-face product
\eqref{eq:end-product} gives
\begin{equation}\label{eq:primitive-contact-generation}
 \alpha t^{\vct}
   =(\alpha t^{\mathbf1_I})
        \prod_{i\in I}
          (1_{S_{\{i\}}^\circ}t^{\mathbf e_i})^{v_i-1}.
\end{equation}
Every constant coefficient on the right restricts to the unit
of $H^*(S_I^\circ;\F)$, also when that face is disconnected.
The generators may be chosen homogeneous: with $\Omega$ the
degree is $|\alpha|+\mu_\Omega(\vct)$, and additivity and
evenness of $\mu_\Omega$ make the identity homogeneous;
without $\Omega$ use cohomological parity.

Suppose inductively that $d_1,\ldots,d_{\ell-1}$ vanish.
Under the resulting canonical identification $E_\ell=E_1$,
Remark~\ref{rem:integral-contact-differentials} shows that any
nonzero component of $d_\ell(\alpha t^{\mathbf1_I})$ has
contact vector
\[
 \vct'=\mathbf1_I-a\mathbf n\in\N^r,
 \qquad a\in\Z_{>0},\qquad \ell=aN_A.
\]
If $n_j>1$ for some $j$, its $j$th coordinate is at most
$1-an_j<0$, a contradiction. If all $n_i=1$ and
$D_1\cap\cdots\cap D_r=\varnothing$, every realized $I$
omits some $j$, giving $v'_j=-a<0$, again a contradiction.
Thus $d_\ell$ vanishes on all these generators. It also
vanishes on $H^*(X)$, which has weight zero and no target
of smaller weight.

By Theorem~\ref{thm:summary}, for homogeneous classes $x,y$,
\[
 d_\ell(xy)=d_\ell(x)y+(-1)^{|x|}x\,d_\ell(y),
\]
where $|x|$ is the total cohomological degree, read modulo two
in parity grading. Hence, for homogeneous factors $g_1,\ldots,g_m$,
\[
 d_\ell(g_1\cdots g_m)
   =\sum_{j=1}^m(-1)^{\sum_{i<j}|g_i|}
       g_1\cdots g_{j-1}\,d_\ell(g_j)\,g_{j+1}\cdots g_m.
\]
Every factor on the right-hand side of
\eqref{eq:primitive-contact-generation} is one of the
primitive-contact generators on which $d_\ell$ has just been
shown to vanish. The displayed Leibniz formula therefore gives
$d_\ell(\alpha t^{\vct})=0$ for every realized contact vector.
Together with $d_\ell|_{H^*(X)}=0$ and linearity, this proves
$d_\ell=0$ on the entire page $E_\ell$. Induction starting with
$\ell=1$ proves that every differential vanishes. Together with the BV
convergence of Theorem~\ref{thm:summary}, this is degeneration
at $E_1$ in the sense of
Definition~\ref{def:bv-spectral-degeneration}, in parity and,
for every chosen $\Omega$, in integer grading.
\end{proof}

\begin{remark}[Computing fusion from leading terms]
\label{rem:computing-fusion}
For a polarized SNC compactification $(M,D,A)$ of $X$, the
additive computation and the multiplicative spectral sequence
give a procedure for computing fusion. Write
$\sigma_p:F_pLH^*(X)\to\operatorname{gr}_F^pLH^*(X)$ for the
quotient map. Theorem~\ref{thm:summary} identifies
$E_\infty\cong\operatorname{gr}_FLH^*(X)$ multiplicatively;
one may replace $E_\infty$ by $E_1=H^*_{\log}(M,D)$ only after
proving degeneration. If chosen homogeneous classes have symbols
generating this associated-graded algebra, they generate
$LH^*(X)$: subtract a polynomial with the same leading symbol
and induct on the nonnegative filtration degree.

For $a\in F_pLH^i(X)$ and $b\in F_qLH^j(X)$, choose
$L(a,b)\in F_{p+q}LH^{i+j}(X)$ lifting
$\sigma_p(a)\sigma_q(b)$. In a filtration-compatible additive
basis the undetermined terms have the form
\[
 \mu(a,b)=L(a,b)+\sum_\lambda c_\lambda x_\lambda,
 \qquad x_\lambda\in F_{p+q-1}LH^{i+j}(X),
\]
with finite support. For $u\in\mathcal O(X)^\times$, if the
inputs have valuation labels $\alpha_u,\beta_u$, choose the
lift in their sum label; then every remaining $x_\lambda$
lies in $LH^*(X)_{u=\alpha_u+\beta_u}$ by
Proposition~\ref{prop:fusion-unit-valuations}.
The grading, unit, symmetries, and associativity further
restrict the coefficients. A zero leading symbol asserts only
$\mu(a,b)\in F_{p+q-1}LH$, not $\mu(a,b)=0$.

Compute any remaining coefficients using the cochain descent
diagram of Lemma~\ref{lem:finite-factor-fusion-formula}(ii), including
the overlap components of the input cocycles. The individual-factor
formula \eqref{lp:eq:full-fusion-formula} applies when its stated
endpoint identifications are available; in particular,
Proposition~\ref{prop:global-binary-fusion}(i)--(ii) applies to the global
finite models specified there. Algebraic cycle representatives
permit specialization with multiplicities; arbitrary classes
require the localization morphism on the cochain descent diagram
before passage to cohomology. A proposed algebra
presentation is complete only after its normal monomials are
proved to give the previously computed additive basis.
These statements apply in parity grading and in the grading
defined by $\Omega$.
\end{remark}

\section{K\"unneth and deformation invariance of loop cohomology}
\label{sec:conclusion}

We prove in this section a K\"unneth theorem for the filtered unital BV algebra
$LH^*(X)$ and deformation invariance for complements of relative
SNC divisors in smooth projective families. 
We follow the grading convention of
Subsection~\ref{subsec:relative-pole-recovery}: in integer grading,
$B_P$, $B$, $\mathbf B$, and $F_P\mathbf B$ stand for
$B_{\Omega,P}$, $B_\Omega$, $\mathbf B_\Omega$, and
$F_P\mathbf B_\Omega$, respectively, with $\Omega$ the chosen
nowhere-zero algebraic volume form. Without $\Omega$, these
symbols denote the parity-graded objects.

\subsection{A K\"unneth theorem for loop cohomology}
\label{subsec:external-products-functoriality}

Let $X,Y$ be smooth affine complex varieties, with polarized SNC
compactifications $(M_X,D_X,A_X)$ and $(M_Y,D_Y,A_Y)$.
Give $X\times Y$ the compactification
\[
 M=M_X\times M_Y,\qquad
 D=(D_X\times M_Y)\cup(M_X\times D_Y),\qquad
 A=\operatorname{pr}_X^*A_X+\operatorname{pr}_Y^*A_Y.
\]
Write $B_X$ and $B_{X,p}$ for the total derived loop object
and its original derived pole-stage object for $X$, respectively;
use the analogous notation for $Y$ and $X\times Y$.
Without volume forms these objects belong
to $\mathcal D_{\mathrm{per}}(\F)$. If nowhere-vanishing forms
$\Omega_X,\Omega_Y$ are specified, use their integer-graded
versions in $\mathcal D(\F)$ and the product form
\[
 \Omega_{X\times Y}
 =\operatorname{pr}_X^*\Omega_X\wedge
                       \operatorname{pr}_Y^*\Omega_Y.
\]
All tensor products are derived over $\F$. For $p,q\ge0$, write
\[
 R_{p,q}:=F_p\sL X\times F_q\sL Y
\]
for the rectangular pole stage. Put
\[
 B_{R_{p,q}}:=R\Gamma(R_{p,q},\omega^{\mathrm n}_{R_{p,q}/\mathrm{pt}}),
 \qquad E_{R_{p,q}}:=\mathcal E_{R_{p,q}/\mathrm{pt}},
\]
using Definition~\ref{def:relative-normalized-dual-complexes}
and the product relative volume-index normalization. In parity
grading use its periodized versions.

\begin{theorem}[K\"unneth theorem for loop cohomology]
\label{prop:fusion-external-products}
\label{prop:bv-external-products}
There is a canonical isomorphism of filtered unital BV algebras
\[
 \kappa_{X,Y}^{H}:LH^*(X)\otimes LH^*(Y)
       \xrightarrow{\sim}LH^*(X\times Y),
 \qquad a\otimes b\longmapsto a\boxtimes b.
\]
The source has the sum filtration
\[
 F_N\bigl(LH^*(X)\otimes LH^*(Y)\bigr)
   =\sum_{p+q\le N}F_pLH^*(X)\otimes F_qLH^*(Y).
\]
For homogeneous classes its product and BV operator are
\[
 (a\boxtimes b)(c\boxtimes d)
    =(-1)^{|b||c|}(ac)\boxtimes(bd),
 \qquad 1_X\boxtimes1_Y=1_{X\times Y},
\]
\begin{equation}\label{eq:bv-external-products}
 \Delta_{X\times Y}(a\boxtimes b)
   =\Delta_X(a)\boxtimes b
       +(-1)^{|a|}a\boxtimes\Delta_Y(b).
\end{equation}
The isomorphism is induced by a canonical equivalence
\begin{equation}\label{eq:loop-external-kunneth}
 \kappa_{X,Y}:B_X\otimes B_Y\xrightarrow{\sim}B_{X\times Y}.
\end{equation}
Here $a\boxtimes b$ is defined as the image of $a\otimes b$
under
\[
 LH^i(X)\otimes LH^j(Y)
   \longrightarrow H^{i+j}(B_X\otimes B_Y)
   \xrightarrow{H^{i+j}(\kappa_{X,Y})}
        LH^{i+j}(X\times Y).
\]
The filtered stages of \eqref{eq:loop-external-kunneth} are the equivalences
\[
 \operatorname*{hocolim}_{p+q\le N}
            (B_{X,p}\otimes B_{Y,q})
       \xrightarrow{\sim}B_{X\times Y,N}.
\]
These comparisons are associative and compatible with the
Koszul permutations of factors. Their filtered
naturality is with respect to isomorphisms identifying the
polarized compactifications and, in integer grading, the
volume-index normalizations. More precisely, for another
polarized compactification $(M_{X'},D_{X'},A_{X'})$, an
isomorphism $\bar f_X:M_X\to M_{X'}$ must satisfy
\[
 \bar f_X^{-1}(D_{X'})=D_X,\qquad
 \bar f_X^*A_{X'}=A_X,
\]
and similarly for the $Y$ factor. In integer grading, writing
$f_X=\bar f_X|_X$, the condition
$f_X^*\Omega_{X'}=c_X\Omega_X$ with $c_X\in\C^\times$
is sufficient to identify the volume indices; analogously,
$f_Y^*\Omega_{Y'}=c_Y\Omega_Y$ with $c_Y\in\C^\times$,
where $f_Y$ is the restriction of the isomorphism of the
$Y$ compactifications. With the displayed volume forms they
give the same statements for
$LH^*_{\Omega_X}(X)\otimes LH^*_{\Omega_Y}(Y)$ and
$LH^*_{\Omega_{X\times Y}}(X\times Y)$ in integer grading.
\end{theorem}

\begin{proof}
\emph{Step 1: Rectangular pole bounds.}
By the product assertion of Lemma~\ref{lp:lem:topological-products},
the projections induce a natural isomorphism of loop functors
\[
 \mathscr L(X\times Y)\xrightarrow{\sim}
       \mathscr LX\times\mathscr LY,\qquad
 \gamma\longmapsto
   (\operatorname{pr}_X\circ\gamma,\operatorname{pr}_Y\circ\gamma).
\]
The subfunctor defined by requiring these two components to
belong to $F_p\mathscr LX$ and $F_q\mathscr LY$, respectively,
is represented by $F_p\mathscr LX\times F_q\mathscr LY$.
Its space of complex points is $R_{p,q}$ with the compactly
generated product topology of that lemma.

\emph{Step 2: Product volume indices and finite-factor complexes.}
In integer grading, choose nice charts
$\widetilde U\to F_p\mathscr LX$ and
$\widetilde V\to F_q\mathscr LY$, with finite-factor projections
$\pi_X:\widetilde U\to U$ and $\pi_Y:\widetilde V\to V$.
For points $x\in\widetilde U(F)$ and $y\in\widetilde V(F)$
over a common field extension $F/\C$, let $\gamma_x,\gamma_y$
be the corresponding loops. The lattices used to define the
volume indices are the coefficient-residue kernels $L_X,L_Y$ of
Lemma~\ref{lem:elementary-coefficient-differentials}:
\[
 \begin{aligned}
 L_X&:=\ker\!\left(
 q_{\pi_X,x}:\Cot_t(\gamma_x)
       \longrightarrow(\Omega^1_{\widetilde U/U})_x\right),\\
 L_Y&:=\ker\!\left(
 q_{\pi_Y,y}:\Cot_t(\gamma_y)
       \longrightarrow(\Omega^1_{\widetilde V/V})_y\right).
 \end{aligned}
\]
Under the identification
$\Cot_t((\gamma_x,\gamma_y))\cong\Cot_t(\gamma_x)\oplus\Cot_t(\gamma_y)$,
formula \eqref{eq:elementary-residue-map} gives
\[
 q_{\pi_X\times\pi_Y,(x,y)}
       =q_{\pi_X,x}\oplus q_{\pi_Y,y},
 \qquad
 \ker q_{\pi_X\times\pi_Y,(x,y)}=L_X\oplus L_Y.
\]
Concatenating bases normalized by $\Omega_X$ and $\Omega_Y$
gives a basis normalized by $\Omega_{X\times Y}$.
Put $n_X=\dim X$ and $n_Y=\dim Y$. In these bases choose
$h\gg0$ with
$t^hF\llbracket t \rrbracket^{n_X}\subset L_X$ and
$t^hF\llbracket t \rrbracket^{n_Y}\subset L_Y$. The formula
\eqref{eq:elementary-field-index} gives
\[
 \begin{aligned}
 \delta_{U\times V}(\pi_X(x),\pi_Y(y))
 &=\dim_F\frac{L_X\oplus L_Y}
                  {t^hF\llbracket t \rrbracket^{n_X}\oplus t^hF\llbracket t \rrbracket^{n_Y}}
                          -(n_X+n_Y)h\\
 &=\left(\dim_F\frac{L_X}{t^hF\llbracket t \rrbracket^{n_X}}-n_Xh\right)
   +\left(\dim_F\frac{L_Y}{t^hF\llbracket t \rrbracket^{n_Y}}-n_Yh\right).
 \end{aligned}
\]
Thus the locally constant volume indices satisfy
\[
 \delta_{U\times V}
     =\operatorname{pr}_U^*\delta_U
        +\operatorname{pr}_V^*\delta_V.
\]
The finite-factor exterior-product comparison is therefore
\[
 \omega_U[-2\delta_U]\boxtimes\omega_V[-2\delta_V]
 \simeq\omega_{U\times V}[-2(\delta_U+\delta_V)].
\]
In parity grading use the periodized comparison. The normalized
purity equivalences and proper trace commute with this comparison, so
Proposition~\ref{lem:relative-compact-support}(iii) descends it
to the normalized dualizing complexes on $R_{p,q}$.

\emph{Step 3: K\"unneth on a rectangular stage.}
For the two pole schemes, choose finite nice \'etale atlases and
compatible finite-factor presentations of finite composable strings
of their overlap maps, with a factor for each chart occurrence
as in \emph{Compatible finite-factor presentations} in
Subsection~\ref{subsec:coherent-descent}.
Write $U^X_{r,\alpha}=V^X_{r,\alpha}\times\A^{J^X_{r,\alpha}}$
and $U^Y_{s,\beta}=V^Y_{s,\beta}\times\A^{J^Y_{s,\beta}}$
for the overlap charts, and put
\[
 \begin{aligned}
 C^X_r&:=\bigoplus_\alpha
 R\Gamma_c(V^X_{r,\alpha};\F)[2\delta^X_{r,\alpha}],&
 E_{X,p}&:=\operatorname*{hocolim}_{[r]\in\Delta^{\mathrm{op}}}C^X_r,\\
 C^Y_s&:=\bigoplus_\beta
 R\Gamma_c(V^Y_{s,\beta};\F)[2\delta^Y_{s,\beta}],&
 E_{Y,q}&:=\operatorname*{hocolim}_{[s]\in\Delta^{\mathrm{op}}}C^Y_s.
 \end{aligned}
\]
Here $\delta^X_{r,\alpha}$ and $\delta^Y_{s,\beta}$ are the
finite-factor volume indices, with shifts taken on their
open-and-closed level sets. The arrows are the normalized traces
of Definition~\ref{def:relative-normalized-dual-complexes}.
All displayed tensor products are in $\mathcal D(\F)$;
$R\Gamma_c$ is applied only to the finite-type factors $V$.
In parity grading replace each displayed term by its periodization.

For every $(r,s,\alpha,\beta)$, finite-type compact-support
K\"unneth and Step~2 give
\[
 \begin{aligned}
 &R\Gamma_c(V^X_{r,\alpha};\F)[2\delta^X_{r,\alpha}]
 \otimes R\Gamma_c(V^Y_{s,\beta};\F)[2\delta^Y_{s,\beta}]\\
 &\hspace{1cm}\xrightarrow{\sim}
 R\Gamma_c(V^X_{r,\alpha}\times V^Y_{s,\beta};\F)
             [2\delta^X_{r,\alpha}+2\delta^Y_{s,\beta}].
 \end{aligned}
\]
The product of the two atlases has these product overlap charts
in equal simplicial degrees. The exterior-product comparison
commutes with normalized trace and its composition homotopies;
hence it identifies its entire trace diagram with
$r\mapsto C^X_r\otimes C^Y_r$. The diagonal
$\Delta^{\mathrm{op}}\to
 \Delta^{\mathrm{op}}\times\Delta^{\mathrm{op}}$
is cofinal for geometric realizations, equivalently by the
Eilenberg--Zilber comparison for a bisimplicial complex. Thus
\begin{equation}\label{eq:rectangular-companion-kunneth}
 \begin{aligned}
 E_{R_{p,q}}
 &\simeq\operatorname*{hocolim}_{r}
                        (C^X_r\otimes C^Y_r)\\
 &\simeq\operatorname*{hocolim}_{r,s}
                        (C^X_r\otimes C^Y_s)
 \simeq E_{X,p}\otimes E_{Y,q}.
 \end{aligned}
\end{equation}
The last equivalence uses preservation of colimits separately by
tensor product. It does not exchange a tensor product with a
\v Cech homotopy limit.

By Theorem~\ref{thm:drinfeld-finite-factors}, each pole scheme
has a finite nice filtration. Proposition~\ref{lem:relative-compact-support}(iv)
therefore makes $E_{X,p}$ and $E_{Y,q}$ perfect. Its duality
comparison and \eqref{eq:rectangular-companion-kunneth} give
\begin{equation}\label{eq:rectangular-loop-kunneth}
 \begin{aligned}
 B_{X,p}\otimes B_{Y,q}
 &\simeq E_{X,p}^{\vee}\otimes E_{Y,q}^{\vee}
 \xrightarrow{\sim}(E_{X,p}\otimes E_{Y,q})^{\vee}\\
 &\simeq E_{R_{p,q}}^{\vee}\simeq B_{R_{p,q}}.
 \end{aligned}
\end{equation}
The middle arrow is an equivalence because both factors are
perfect. Indeed, on a nice difference the companion is the
shifted compact-support cohomology of a finite-type factor,
which has finite-dimensional total cohomology. The finite
filtration and its localization triangles preserve this
property. Thus the duality comparison is made at each fixed
pair $(p,q)$, before taking the pole colimit.
No assertion about dualizing an arbitrary tensor
product of infinite-dimensional complexes is required.
All comparisons are compatible with the pole inclusions: contravariant pole
maps on the compact-support companions become proper-trace maps
on their duals. Thus \eqref{eq:rectangular-loop-kunneth} is an
equivalence of the rectangular pole diagrams.

\emph{Step 4: The sum filtration and the pole colimit.}
Weights for the product polarization satisfy
$w(\mathbf v,\mathbf w)=w_X(\mathbf v)+w_Y(\mathbf w)$.
Thus on complex points the rectangles give the finite closed cover
\[
 F_N\sL(X\times Y)=\bigcup_{p+q\le N}R_{p,q},\qquad
 R_{p,q}\cap R_{p',q'}
     =R_{\min(p,p'),\min(q,q')}.
\]
The corresponding closed inclusions of pole schemes are finitely
presented. On each chart of the pole scheme $F_N\sL(X\times Y)$,
present their pullbacks and intersections with a common affine tail
as in Subsection~\ref{subsec:volume-index-preparations}.
For $p+q\le N$, let $i_{p,q}^{N}:R_{p,q}\hookrightarrow F_N\sL(X\times Y)$
denote the inclusion. The closed-inclusion
compatibility of the volume indices and
Proposition~\ref{lem:relative-compact-support}(i) give
\[
 \omega^{\mathrm n}_{R_{p,q}/\mathrm{pt}}
       \simeq (i_{p,q}^{N})^!\omega^{\mathrm n}_{F_N\sL(X\times Y)/\mathrm{pt}}.
\]
These are the normalized dualizing complexes in the chosen grading. Nilpotents in the
scheme-theoretic union do not affect its Betti sheaf complexes.

For two members $Z,Z'$ of a finite closed diagram, localization
gives the closed-cover triangle
\[
 B_{Z\cap Z'}\xrightarrow{(i_*,-i'_*)}
 B_Z\oplus B_{Z'}\longrightarrow B_{Z\cup Z'}
 \xrightarrow{+1}.
\]
On finite factors this is the dual of the constant-sheaf
closed-cover sequence, with the inherited volume shifts;
descent preserves this cofiber sequence. Applying the triangle
successively to the rectangles, including their displayed
intersections, gives
\[
 \operatorname*{hocolim}_{p+q\le N}B_{R_{p,q}}
    \xrightarrow{\sim}B_{X\times Y,N}.
\]
One may check the indexing directly: for each point of the union,
the rectangles containing it have a least pair of bounds, and
their nerve is contractible. Thus the rectangular diagram is
the augmented closed-cover diagram with its repeated
intersections, not a direct sum of overlapping rectangles.
The equivalence on its $(p,q)$ term is the explicit arrow
\[
 \kappa^N_{p,q}:B_{X,p}\otimes B_{Y,q}
     \xrightarrow{\eqref{eq:rectangular-loop-kunneth}}B_{R_{p,q}}
     \xrightarrow{(i_{p,q}^{N})_*}B_{X\times Y,N},
 \qquad p+q\le N.
\]
Here $(i_{p,q}^{N})_*$ is the normalized proper-trace map.
These arrows commute coherently with increases of $p,q,N$.
Taking their pole homotopy colimit defines $\kappa_{X,Y}$,
and gives
\[
 \begin{aligned}
 B_{X\times Y}
 &\simeq\operatorname*{hocolim}_{p,q}
                  (B_{X,p}\otimes B_{Y,q})\\
 &\simeq
  \left(\operatorname*{hocolim}_pB_{X,p}\right)\otimes
  \left(\operatorname*{hocolim}_qB_{Y,q}\right).
 \end{aligned}
\]
The second equivalence uses preservation of colimits separately
by tensor product. For the isomorphisms specified in the theorem,
transporting projective frames identifies the pole stages with
the same bounds. The cotangent identifications preserve the
coefficient-residue maps $q_{\pi_X,x},q_{\pi_Y,y}$ and their
pointwise coefficient indices $c_{\pi_X}(x),c_{\pi_Y}(y)$;
in the sufficient volume-form condition, $\ord_t(c_X)=\ord_t(c_Y)=0$. Naturality in this sense,
associativity and permutation compatibility therefore follow
from the corresponding coherent comparisons for the finite-factor
exterior products and traces before descent.

It remains to identify the filtration after taking cohomology;
this is not an assertion that cohomology commutes with an
arbitrary finite homotopy colimit. Choose exhaustive models
$F_pC_X\subset C_X$ and $F_qC_Y\subset C_Y$ by subcomplexes,
as in Theorem~\ref{thm:canonical-complex}, together with specified equivalences
of entire diagrams
$(F_pC_X)_p\simeq(B_{X,p})_p$ and
$(F_qC_Y)_q\simeq(B_{Y,q})_q$. Their tensor model
has stages $\sum_{p+q\le N}F_pC_X\otimes F_qC_Y$ and represents
the displayed convolution of filtered objects: the intersection
of two rectangular subcomplexes is the rectangle with the two
minimum bounds, so the latching maps in this finite diagram
are inclusions of subcomplexes. For either of these exhaustively
filtered complexes $C$, with $F_pC=0$ for $p<0$, filtered
linear algebra supplies splittings
\[
 Z^j(C)=dC^{j-1}\oplus H_C^j,\qquad
 C^j=Z^j(C)\oplus T_C^j
\]
compatible with every filtration step. Here the subspaces of
$C^j$ carry their intersection filtrations; in particular, the
filtration on $dC^{j-1}$ is $dC^{j-1}\cap F_pC^j$.
The complements are obtained by extending adapted bases
successively along $F_0C^j\subset F_1C^j\subset\cdots$;
exhaustiveness makes their union a basis of $C^j$.
The differential vanishes
on $H_C^j$; the induced filtration on this subspace is identified
with $\operatorname{im}(H^j(F_pC)\to H^j(C))$. Inclusion of
$H_C$ and projection $C\to H_C$ are filtered chain maps and
induce inverse isomorphisms on cohomology. Their tensor products
are filtered chain maps inducing inverse isomorphisms on total
cohomology, since tensor product over $\F$ preserves
quasi-isomorphisms. Projection sends a cycle in the $N$th tensor
stage to the corresponding sum of filtered subspaces of
$H_{C_X}\otimes H_{C_Y}$; inclusion gives the reverse containment.
This proves
\[
 F_NH^*(C_X\otimes C_Y)
   =\sum_{p+q\le N}
       F_pH^*(C_X)\otimes F_qH^*(C_Y).
\]
Together with the filtered equivalences already proved, this is
the filtration asserted in the theorem. The splittings are used
only to prove the equality; the map $\kappa_{X,Y}$ was defined
without them.

\emph{Step 5: Fusion and the unit.}
For rectangular input bounds $(p,q),(p',q')$, mapping the common
two-point source into $X\times Y$ is the same as mapping it into
$X$ and into $Y$. Thus the rectangular collision family is
\[
 \mathcal Z^X_{p,p'}\times_{\A^1_s}\mathcal Z^Y_{q,q'},
\]
with output $R_{p+p',q+q'}$ and its closed inclusion into the
pole stage of total weight $p+p'+q+q'$. These diagrams map to
the product-polarization collision diagram: in projective
product coordinates, pole bounds add because the coordinate
functions are products of the functions from the two factors.
More explicitly, use a common very ample multiple of the two
polarizations and its Segre product coordinates. These include
each individual factor coordinate as its product with $1$.
Within the total collision bound, imposing the two separate
marked bounds sets finitely many Laurent coefficients equal
to zero. It therefore defines a finitely presented closed
inclusion of the rectangular collision family into the total
collision family. Retain these equations in the finite factors.
The direct-sum calculation of Step~2 applies to the relative
coefficient-residue maps $\kappa_{\pi,V}$ of
\eqref{eq:direct-relative-residue-map} for the common two-point
source. It identifies the product volume index with the
restriction of the ambient index. The normalized
proper-trace maps for this inclusion consequently commute
with specialization by
Lemma~\ref{lem:specialization-arrow-descent}(ii), as in the proof
of Lemma~\ref{lem:source-specialization-comparison}(iii).
It suffices, therefore, to compute specialization on the
rectangular family and then apply its output trace.

To compare their specialization maps, first give the two
families independent parameters $s_X,s_Y$. Their exterior product
is a family over $\A^2_{s_X,s_Y}$; pullback by
$\delta_{\C}:\A^1_s\to\A^2$, $s\mapsto(s,s)$, is the
displayed common-parameter family. Put $H=[0,\infty)$ and
$L=(0,\infty)$, and let $E_X,E_Y$ be the restrictions to $H$
of the two relative compact-support companions. Factorization
identifies their stalks by
\[
 \begin{array}{ll}
 (E_X)_0\simeq E_{X,p+p'},&
 E_X|_L\simeq\underline{E_{X,p}\otimes E_{X,p'}}_L,\\
 (E_Y)_0\simeq E_{Y,q+q'},&
 E_Y|_L\simeq\underline{E_{Y,q}\otimes E_{Y,q'}}_L.
 \end{array}
\]
In particular, the hypotheses of
Lemma~\ref{lem:parameter-chain-properties}(ii)--(iii) hold.
By Proposition~\ref{lem:relative-compact-support}(ii),(iv),
applied before restriction to these real parameter spaces,
the companion of the common-parameter family is
$\delta^*(E_X\boxtimes E_Y)$, where
$\delta=\delta_{\C}|_H:H\to H^2$.
The natural section-to-pullback arrows give the commutative square
\[
 \begin{array}{ccc}
 R\Gamma(H^2,E_X\boxtimes E_Y)&\longrightarrow&
       R\Gamma(L^2,E_X\boxtimes E_Y)\\
 \big\downarrow&&\big\downarrow\\
 R\Gamma(H,\delta^*(E_X\boxtimes E_Y))&\longrightarrow&
       R\Gamma(L,\delta^*(E_X\boxtimes E_Y)),
 \end{array}
\]
whose horizontal arrows are restrictions. Evaluation at $(0,0)$
and $0$ identifies the left vertical arrow with the identity of
$E_{X,p+p'}\otimes E_{Y,q+q'}$. Evaluation at $(1,1)$ and $1$
identifies the right vertical arrow with the identity of
$(E_{X,p}\otimes E_{X,p'})\otimes(E_{Y,q}\otimes E_{Y,q'})$.
The evaluation arrows are equivalences by part (ii) of that lemma. Its
exterior-product assertion in part (iii) identifies the upper arrow with
$r_X\otimes r_Y$, where
\[
 \begin{aligned}
 r_X &: E_{X,p+p'}\simeq R\Gamma(H,E_X)
       \longrightarrow R\Gamma(L,E_X)
       \simeq E_{X,p}\otimes E_{X,p'},\\
 r_Y &: E_{Y,q+q'}\simeq R\Gamma(H,E_Y)
       \longrightarrow R\Gamma(L,E_Y)
       \simeq E_{Y,q}\otimes E_{Y,q'}.
 \end{aligned}
\]
Consequently the restriction arrow for the common-parameter
family is $r_X\otimes r_Y$ under these identifications.
The specialization--restriction comparison of
Proposition~\ref{prop:relative-binary-identification}, followed
by linear duality, proves
\begin{equation}\label{eq:loop-external-fusion}
 \overline m_{2,X\times Y}(\kappa_{X,Y}\otimes\kappa_{X,Y})
 \simeq\kappa_{X,Y}(\overline m_{2,X}\otimes\overline m_{2,Y})
                (\mathrm{id}\otimes\tau\otimes\mathrm{id}),
\end{equation}
first on rectangular pole stages and then on their colimit.
Here $\overline m_{2,X}$, $\overline m_{2,Y}$, and
$\overline m_{2,X\times Y}$ are the derived fusion morphisms,
and $\tau:B_Y\otimes B_X\to B_X\otimes B_Y$ is the Koszul
interchange. The inward-ray conventions agree under this
restriction comparison; thus its sole interchange sign is
$(-1)^{|b||c|}$. The constant arc classes give
$1_X\boxtimes1_Y=1_{X\times Y}$.

\emph{Step 6: Rotation.}
Source rotation on $\mathscr L(X\times Y)$ rotates both
components with the same circle parameter. The exterior-product
comparisons above preserve this diagonal action. Under the circle
diagonal $d:S^1\to S^1\times S^1$,
\[
 d_*[S^1]=[S^1]\times[\mathrm{pt}]
                  +[\mathrm{pt}]\times[S^1].
\]
Apply the parameter-chain construction of
Proposition~\ref{prop:coherent-source-rotation} to this equality.
Moving its degree-$-1$ operator past $a$ gives
\eqref{eq:bv-external-products}. The bracket is determined by
the product and $\Delta$, so the cohomological K\"unneth
isomorphism preserves the full unital BV structure.
\end{proof}

\begin{corollary}[Product projections]\label{cor:projection-functoriality}
In the setup of the theorem, the projection
$\operatorname{pr}_X:X\times Y\to X$ induces the filtered
unital BV-algebra homomorphism
\[
 \operatorname{pr}_X^*:LH^*(X)\longrightarrow LH^*(X\times Y),
 \qquad a\longmapsto a\boxtimes1_Y.
\]
It has degree zero in the stated volume gradings. Iterated
product projections give $a\mapsto a\boxtimes1_Y\boxtimes1_Z$.
\end{corollary}
\begin{proof}
The unit satisfies $1_Y\in F_0LH^0(Y)$ and $\Delta_Y1_Y=0$.
The product and BV formulas in the theorem prove preservation
of the filtration and the operations. Associativity gives the
iterated-product identity.
\end{proof}

\subsection{Deformation invariance of loop cohomology}
\label{subsec:loop-deformation-invariance}

Let $\pi:\mathcal M\to T$ be a smooth projective morphism of
relative dimension $n>0$ between complex varieties, with $T$ smooth, and let
$\mathcal D=\bigcup_{i=1}^r\mathcal D_i$ be a relative SNC divisor
with labeled components. Fix
\[
 \mathcal A=\sum_{i=1}^r k_i\mathcal D_i,
 \qquad k_i\in\Z_{>0},\qquad
 \mathcal A\text{ relatively ample},\qquad
 \mathcal X=\mathcal M\setminus\mathcal D.
\]
Write $(M_t,D_t,A_t)$ and $X_t$ for the fibers. The weights
$w(\vct)=\sum_i k_iv_i$ are independent of $t$.
All sheaves below are on analytic spaces. A complex is
\emph{locally constant} if it is locally equivalent, in the
derived $\infty$-category of sheaves, to a constant complex.
This is a condition on the complex, not only on the dimensions
of its stalk cohomology. Denote their category by
$\operatorname{Loc}(T;\mathcal D(\F))$; use
$\mathcal D_{\mathrm{per}}(\F)$ in parity grading.

\begin{theorem}[Deformation invariance for relative SNC pairs]
\label{thm:loop-deformation-invariance}
\label{prop:conclusion-relative-loop}
In the above setup, let $B_{X_t,P}$ denote the original derived
pole-stage object for $X_t$. The diagrams
$P\mapsto B_{X_t,P}$ are the fibers of a
coherent diagram of locally constant complexes on $T$.
Consequently a path $\gamma:[0,1]\to T$ induces an equivalence
\[
 T^B_\gamma:
 (B_{X_{\gamma(0)},P})_{P\ge0}
       \xrightarrow{\sim}
 (B_{X_{\gamma(1)},P})_{P\ge0}
 \quad\text{in }
 \operatorname{Fun}(N(\Z_{\ge0}),\mathcal D_{\mathrm{per}}(\F)).
\]
These equivalences preserve fusion, the unit, and source rotation.
They induce a local system of filtered unital BV algebras
$t\mapsto LH^*(X_t)$. In particular, their cohomological transports
satisfy
\[
 T^H_{\gamma_2*\gamma_1}
     =T^H_{\gamma_2}T^H_{\gamma_1},\qquad
 T^H_{\mathrm{const}_t}=\mathrm{id},\qquad
 T^H_\gamma=T^H_{\gamma'}\quad
   (\gamma\simeq\gamma'\text{ relative to their endpoints}).
\]
If a nowhere-vanishing relative volume form
$\Omega\in\Gamma(\mathcal X,\Omega^n_{\mathcal X/T})$ is given,
the same assertions hold in $\mathcal D(\F)$ and for the
integer-graded algebras $LH^*_{\Omega_t}(X_t)$.
The constructions are independent of finite-factor presentations
under their common-refinement comparisons.
\end{theorem}
 
We need to extend the finite factors discussion from  Section~\ref{sec:relative-finite-factors} for a varying smooth affine target.

\begin{lemma}
\label{lem:varying-target-finite-factors}
Let $A_0$ be a finitely generated $\C$-algebra and let
$\mathcal X=\Spec C\to T=\Spec A_0$ be smooth affine of pure
relative dimension $n>0$. Choose an $A_0$-algebra presentation
$C=A_0[x_1,\ldots,x_N]/(f_1,\ldots,f_c)$.
Let $A$ be a finite-type $A_0$-algebra, let $D(z)\in A[z]$
be monic of degree $k\ge1$. For $b\ge0$ and an $A$-algebra $R$, put
\[
 O_R=R\llbracket q \rrbracket[z]/(D(z)-q),\qquad E_R=O_R[q^{-1}],
 \qquad
 S_b(R)=\{u\in q^{-b}O_R^N:f_1(u)=\cdots=f_c(u)=0\}.
\]
Thus $S_b$ is the Laurent coefficient scheme of maps
$\Spec E_R\to\mathcal X$ over $T$ with the displayed bound.
The following assertions hold.
\begin{enumerate}[(i),wide=0pt]
\item Every $S_b\to\Spec A$ has a finite Nisnevich atlas
\[
 U_\alpha\longrightarrow S_b,\qquad
 U_\alpha\cong Y_\alpha\times_{\Spec A}\A_A^{J_\alpha},
\]
where $Y_\alpha$ is affine of finite type over $A$ and
$J_\alpha$ is finite or countable. It has a finite filtration
by finitely presented closed subschemes with nice successive
differences. These assertions persist under finite-type base
change and finitely presented locally closed restriction.
\item For the marked sources
$A=A_0[a_i:i\in I]$ and $D(z)=\prod_i(z-a_i)$, the pole and
collision diagrams obtained from these coefficient schemes
have compatible finite-factor presentations over $T\times\A^I$.
These finite-factor presentations include the pole inclusions, their open complements,
collision diagonals, separation and relabeling isomorphisms,
affine source changes, and any prescribed finite list of their
compositions and \'etale overlap maps. Two such finite choices
have common refinements. On a collision stratum, the product
is the fiber product of the input families over $T$.
\item Suppose
$\Omega\in\Gamma(\mathcal X,\det\Omega^1_{\mathcal X/T})$
is nowhere vanishing. The finite-factor presentations in (i)--(ii) carry
relative volume indices over their specified bases, compatible
with smooth comparisons, base changes, closed restrictions,
and factorization. On a fiber over $t\in T$ they are the
volume indices specified by $\Omega_t$ in the opening normalization
assertion of Proposition~\ref{prop:direct-relative-volume-indices}.
In parity grading the same finite-factor presentations and their
two-periodic purity comparisons require no volume form.
\end{enumerate}
\end{lemma}

\begin{proof}  \emph{Step 1: Relative coefficient estimates.}
Write $I=(f_1,\ldots,f_c)\subset A_0[x_1,\ldots,x_N]$.
Formal smoothness of $C/A_0$ gives an $A_0$-algebra lift of the
identity,
\[
 C\longrightarrow A_0[x_1,\ldots,x_N]/I^2
       \longrightarrow C.
\]
Polynomial representatives define $p:\A_{A_0}^N\to\A_{A_0}^N$ with
\[
 p(x)-x=G(x)f(x),\qquad
 f_\alpha(p(x))=
   \sum_{\beta,\gamma}H_{\alpha\beta\gamma}(x)
                         f_\beta(x)f_\gamma(x),
 \qquad G,H\text{ defined over }A_0.
\]
In the basis $1,z,\ldots,z^{k-1}$, multiplication in $O_R$
has coefficients in $A[q]$. Hence its $q$-order satisfies
$\nu(ab)\ge\nu(a)+\nu(b)$, and the polynomial estimates of
Lemma~\ref{lem:direct-newton-retraction}(i)--(ii) hold with a single
constant depending on $b,p,f,D$, not on the $A$-algebra $R$.
In particular, the correction
$\rho(u)=\lim_jp^{\circ j}(u)$ and its coefficients are
defined over $A$ and commute with homomorphisms of test
algebras. Differentiation here and below is relative to $A_0$.

\emph{Step 2: Finite approximations and lattices.}
Form the finite-type $A$-schemes $T_m$ with their universal
truncated tuples $t_m$ and exact centers $b_m=\rho(t_m)$ as
in Step~1 of Subsection~\ref{subsec:relative-finite-factor-proof},
using its conventions for coefficientwise pullback and completed
change of coefficients. For $R_m=\mathcal O(T_m)$, the universal
morphism $b_m:\Spec E_{R_m}\to\mathcal X$ gives
\[
 M_m=\Gamma(\Spec E_{R_m},b_m^*T_{\mathcal X/T})
     =\operatorname{im}(d p_{b_m})\subset E_{R_m}^N.
\]
This is finite projective of rank $n$ over $E_{R_m}$ and
therefore finite projective of rank $nk$ over $R_m(\!( q )\!)$.
The second assertion uses only that $E_{R_m}$ is finite
free of rank $k$ over $R_m(\!( q )\!)$; it does not require this
source algebra to be \'etale. Since $R_m$ is of finite type
over $\C$, Lemma~\ref{lem:free-tate-lattice} supplies the
Nisnevich cover and coprojective product lattices.
Apply Steps~2--6 of
Subsection~\ref{subsec:relative-finite-factor-proof} with
$T_{\mathcal X/T}$ in place of $T_X$ and the relative polynomial
bounds from Step~1: choose $m_0$ first, then a lattice on each
member $V_\alpha\to T_{m_0}$ of the resulting Nisnevich cover,
and finally its scale and the longer precision $m$.
The universal solution of Step~5 and the jet equations
of Step~6 then give
\[
 S_b\times_{T_{m_0}}V_\alpha
       \cong Y_\alpha\times_{\Spec A}\A_A^{J_\alpha}.
\]
This is an isomorphism on all test algebras. Steps~8--9 of that
proof give the finite-type base changes, finitely presented
locally closed restrictions, and finite nice filtrations,
proving (i).

\emph{Step 3: Diagrams.}
The marked pole equations are the equations of
Definition~\ref{def:kv-bounded-family}, now with coefficients
in $A_0$. Their source-algebra identities are unchanged over
$A_0$. For each finite diagram, repeat the construction in the
proof of Proposition~\ref{prop:direct-collision-diagrams}, using
the relative atlases of part~(i) and the equations over $A_0$.
Lemma~\ref{lem:finite-presentation-categories} supplies the
finite-factor presentations for fixed strings of chart maps, with a factor
for each occurrence. Products over the deformation base and
the specified base changes give the separation and diagonal
diagrams. Descended composites agree after the faithfully
flat finite-factor projections, hence agree. Different atlas
choices are compared as in \emph{Comparing different atlases}
in Subsection~\ref{subsec:coherent-descent}, with the maps to
the deformation base retained. This proves (ii) over that
base itself.

\emph{Step 4: Relative volume indices.}
On a chart $\pi:U=Y\times\A^J\to Y$ over $\Spec A$ and
a test scheme $V=\Spec R\to U$, let
$\gamma_V:\Spec E_R\to\mathcal X$ be the corresponding Laurent map
over $T$. Using relative differentials, put
\[
 N_V=\gamma_V^*\Omega^1_{\mathcal X/T},\qquad
 \Cot_q(\gamma_V)=
 \bigl(\operatorname{Res}_{E_R/R(\!( q )\!)}N_V\bigr)\otimes dq.
\]
For the marked sources in (ii), the images of the marks in $R$
give a configuration $a\in R^I$. With $q=q_0$,
\eqref{eq:kv-monic-source} identifies $O_R$ with $\mathcal O_a$
and $E_R$ with $\mathcal K_a$.
The relative volume form and the ordered basis $(1,z,\ldots,z^{k-1})$
trivialize $\det \Cot_q(\gamma_V)$ by the determinant--norm formula of
Subsection~\ref{subsec:direct-relative-volume-indices}.
Use \eqref{eq:direct-relative-residue-map} with $N_V$ as above
and $\lambda_R:E_R\to R(\!( q )\!)$ given by extraction of the
$z^{k-1}$ coefficient. The elements of $A$ have zero
$d_{U/Y}$-differential, so this formula respects $D(z)-q=0$ and
the coefficient calculation of Lemma~\ref{lem:direct-relative-coefficient-kernel}
applies using the Jacobian of $\mathcal X$ over $T$.
The finite quotients of $\ker\kappa_{\pi,V}$ in that lemma
therefore commute with base change here as well. At a field point
$x$ take $c_\pi(x)$ as in \eqref{eq:direct-relative-field-index},
using $L_x=\ker\kappa_{\pi,x}$. The local constancy and comparison
arguments in Steps~1--2 of
Proposition~\ref{prop:direct-relative-volume-indices} now show
that $c_\pi=\pi^*\delta_{Y/\Spec A}$ for a locally constant
function $\delta_{Y/\Spec A}$ on $Y$ and that
\[
 \delta_{Y'/\Spec A}=g^*\delta_{Y/\Spec A}+e,
 \qquad
 \delta_{Y\times_AA'/\Spec A'}
        =\operatorname{pr}_Y^*\delta_{Y/\Spec A}
\]
for a smooth comparison $g:Y'\to Y$ of relative dimension $e$
and a specified base change $A\to A'$ in the finite-factor presentation
diagram. The same arguments apply to finitely presented
locally closed restrictions presented on the finite factors, as explained after
Lemma~\ref{lem:direct-relative-coefficient-kernel}.
On separated sources, the relative coefficient-residue map
$\kappa_{\pi,V}$ is identified with the direct sum of the
corresponding maps on the disk factors; the pointwise
coefficient index $c_\pi$ is therefore the sum of their indices.
These are the factorization and source-change
comparisons of parts (ii)--(iii) of that proposition, proved in
Steps~3--4 of its proof, with
deformation parameters held fixed. Restricting to $t\in T$
recovers the map and determinant specified by $\Omega_t$,
and hence exactly its fiber volume index. No dimension
of $T$ is added to this relative index. Two-periodic
purity gives the last assertion of (iii).
\end{proof}

\begin{proof}[Proof of Theorem~\ref{thm:loop-deformation-invariance}]
We work on a connected component on which
$d=\dim_\C T$ and $n=\dim_\C X_t$ are constant.
In the integer-graded proof $B_{X_t,P}$ denotes the derived
stage $B_{\Omega_t,P}$ for the fiber $X_t$ with volume form
$\Omega_t$. Periodizing the finite-factor
construction gives the parity proof without choosing $\Omega$.

\emph{Step 1: the relative pole diagrams.}
Choose $m>0$ so that $\mathcal O_{\mathcal M}(m\mathcal A)$
is relatively very ample and its direct image is locally free,
with formation commuting with base change. Locally on $T$,
extend the section $s_0$ with divisor $m\mathcal A$ to a frame
$s_0,s_1,\ldots,s_N$. The functions $g_j=s_j/s_0$ embed
$\mathcal X$ as a closed subscheme of $\A_T^N$.
Let $f_P:S_P\to T$ be the Laurent coefficient scheme of maps
to $\mathcal X$ over $T$ satisfying, for the formal source
parameter $\tau$,
\[
 g_j(\gamma(\tau))\in \tau^{-mP}R\llbracket \tau \rrbracket\quad(1\le j\le N)
 \qquad\text{for a test algebra }R\text{ over }T.
\]
Changes of frame fixing $s_0$ preserve these bounds, so the local
schemes glue. Their fibers have the pole-stage coefficient presentations
of Proposition~\ref{prop:loop-sectors}; in particular,
$(S_P)_t(\C)=F_P\mathscr LX_t$. The inclusions
$S_{P-1}\hookrightarrow S_P$ are finitely presented and closed.

Apply Lemma~\ref{lem:varying-target-finite-factors}(i)--(iii) to these
families and to their relative collision diagrams. It supplies
compatible finite factors over $T$ and over $T\times\A^I$,
including the relative volume indices in integer grading.
Definition~\ref{def:relative-normalized-dual-complexes} and
the duality comparison and parts (i)--(ii) of
Proposition~\ref{lem:relative-compact-support} therefore give
coherent diagrams
\[
 \mathcal E_P:=\mathcal E_{S_P/T},\qquad
 \mathcal B_P:=\mathcal B_{S_P/T}
       \simeq\mathbb D_T\mathcal E_P,
 \qquad i_t^!\mathcal B_P\simeq B_{X_t,P},
 \quad i_t:\{t\}\hookrightarrow T.
\]
All identifications here include the closed pole maps.
Put $U_P=S_P\setminus S_{P-1}$ for $P\ge1$ and
$\mathcal G_P=\mathcal B_{U_P/T}$. Localization gives
\begin{equation}\label{eq:conclusion-relative-triangle}
 \mathcal B_{P-1}\longrightarrow\mathcal B_P
    \longrightarrow\mathcal G_P\xrightarrow{+1},\qquad
 \mathcal E_{U_P/T}\longrightarrow\mathcal E_P
    \longrightarrow\mathcal E_{P-1}\xrightarrow{+1}.
\end{equation}
These triangles are constructed on the relative finite factors
and descended as diagrams. Their costalks are the localization
triangles defining the contact filtration of each fiber.

\emph{Step 2: identification of the actual contact layers.}
For $\vct\in\N^r$, let
$\mathscr P_{\vct}(\mathcal M/T)$ denote the relative arc locus
with exact contacts $v_i$ with $\mathcal D_i$. In a relative SNC
coordinate chart its equations are
\[
 z_i(\tau)=\tau^{v_i}u_i(\tau),\quad u_i(0)\ne0\quad(i\in I(\vct)),
 \qquad \gamma(0)\in\mathcal D_{I(\vct)}^\circ.
\]
The projective-coordinate calculation in
Proposition~\ref{prop:loop-sectors} applies with coefficients
in $\mathcal O(T)$ and gives
\begin{equation}\label{eq:relative-transport-contact-strata}
 (U_P)_{\mathrm{red}}
   \cong\coprod_{w(\vct)=P}\mathscr P_{\vct}(\mathcal M/T).
\end{equation}
Indeed, on the open set where $[\tau^{-mP}]g_j$ is invertible,
write $g_i=\tau^{-mP}u_i$. The projective coordinates
$[\tau^{mP}:u_1:\cdots:u_N]$ define the extended arc and give
the inverse to restriction to the punctured disk. The inverse
series use only finitely many coefficients in each order.
The identity $\sum_i k_iv_i=P$ and the positivity of the $k_i$
then make the finitely many exact-contact loci open and closed
in the reduced layer, by the same coefficient argument as in
that proposition. Thus \eqref{eq:relative-transport-contact-strata}
is a relative scheme identification, not merely a fiberwise
cohomology comparison. The comparison in
\emph{Reduction of finite models} in
Subsection~\ref{subsec:volume-index-proofs}, with the relative
indices restricted as in Definition~\ref{def:relative-volume-index}(iii),
identifies the normalized Betti complexes and the
localization maps in the original finite-factor presentation diagrams.

Put $I=I(\vct)$ and, for $\ell\ge\max_i v_i$, set
\[
 Q_{\vct,\ell}=j_\ell\mathscr P_{\vct}(\mathcal M/T),\qquad
 N_I^\times=
 \underset{\mathcal D_I^\circ}{\prod_{i\in I}}
       (N_{\mathcal D_i/\mathcal M}|_{\mathcal D_I^\circ}\setminus0).
\]
For $I=\varnothing$, set $N_\varnothing^\times=\mathcal X$.
Leading normal coefficients define
$Q_{\vct,\ell}\to N_I^\times$. The relative SNC expansions
in Proposition~\ref{prop:stable-strata} exhibit this as an
iterated affine bundle. Its source is smooth over $T$ of
relative dimension
\[
 d_{\vct,\ell}=n(\ell+1)-\sum_i v_i,
\]
and further truncation adds an affine $n$-bundle. The omitted
coefficients on $\mathscr P_{\vct}$ are countable affine tails.

Write
\[
 \operatorname{div}_{\mathcal M/T}(\Omega)
       =-\sum_i a_i\mathcal D_i,\qquad
 \mu_\Omega(\vct)=2\sum_i(1-a_i)v_i.
\]
These integers are constant on the relevant labeled components.
The relative volume-index calculation of
Lemma~\ref{lem:varying-target-finite-factors}(iii), applied in the
displayed SNC coordinates, is
\[
 \delta_{Q_{\vct,\ell}/T}
    =n(\ell+1)-\sum_i a_iv_i
    =d_{\vct,\ell}+\tfrac12\mu_\Omega(\vct).
\]
Complex orientation consequently gives
\begin{equation}\label{eq:relative-transport-contact-sheaf}
 \omega_{Q_{\vct,\ell}}[-2\delta_{Q_{\vct,\ell}/T}]
       \simeq\underline\F_{Q_{\vct,\ell}}
                         [2d-\mu_\Omega(\vct)].
\end{equation}
The term $2d$ comes from the absolute dualizing complex of the
finite total factor; base coordinates do not contribute to its
relative volume index. The relative normalized purity equivalences
make these identifications compatible with truncations and with changes of relative SNC
coordinates. Ordinary inverse image transports them across the
affine tails, using Lemma~\ref{lem:infinite-product-descent}(i).

Let $a_I:N_I^\times\to T$ be the structural map.
The normalized traces of the finite affine-bundle projections,
followed by descent, now give
\[
 \mathcal E_{\mathscr P_{\vct}/T}
    \simeq R(a_I)_!\underline\F_{N_I^\times}
                                    [2n+\mu_\Omega(\vct)],
 \qquad
 \mathcal B_{\mathscr P_{\vct}/T}
    \simeq R(a_I)_*\underline\F_{N_I^\times}
                                    [2d-\mu_\Omega(\vct)].
\]
For example, an affine-bundle factor of relative rank
$d_{\vct,\ell}-n$ contributes $[-2(d_{\vct,\ell}-n)]$
to its compact-support image, which cancels that part of
$[2\delta_{Q_{\vct,\ell}/T}]$. This proves the first formula;
Verdier duality on $T$ proves the second.
The same calculation at $\vct=0$ yields
\[
 \mathcal E_0\simeq R(\pi|_{\mathcal X})_!\underline\F_{\mathcal X}[2n],
 \qquad
 \mathcal B_0\simeq R(\pi|_{\mathcal X})_*\underline\F_{\mathcal X}[2d].
\]

Let $p_I:\mathcal S_I^{\log}\to T$ denote the family of closed
normal faces in the fiberwise real oriented blowup of
$(\mathcal M,\mathcal D)$. Quotienting each nonzero normal vector
by positive real scaling gives the open face; its inclusion
into the closed face is a homotopy equivalence by a boundary
collar. Both constructions work over the base. Thus
\eqref{eq:relative-transport-contact-strata} and the preceding
calculations identify the localization cofiber as
\begin{equation}\label{eq:relative-transport-log-layer}
 \mathcal G_P\simeq
  \bigoplus_{w(\vct)=P}
     R(p_{I(\vct)})_*\underline\F_{\mathcal S_{I(\vct)}^{\log}}
                    [2d-\mu_\Omega(\vct)].
\end{equation}
In parity grading all the displayed even normalization shifts
are omitted.

\emph{Step 3: local constancy and the entire pole diagram.}
The fiberwise oriented blowup is a proper smooth family of
manifolds with corners, and its map to $T$ is a submersion on
every boundary stratum. On a coordinate ball in $T$, lift the
coordinate vector fields to vector fields tangent to all boundary
faces. Relative SNC charts give such lifts locally; a partition
of unity gives them on the total space. Over a smaller ball,
properness makes their flows exist for the parameter intervals
under consideration. Their compositions trivialize the family
and every labeled face. Restricting this trivialization to the
interior also trivializes $\mathcal X\to T$.
The normal line bundles are trivialized along the same parameter
directions, for example by parallel transport in chosen smooth
connections. Hence $a_I:N_I^\times\to T$ is locally a product
bundle as well. These trivializations prove that
$R(a_I)_*\F$, $R(a_I)_!\F$, and $R(p_I)_*\F$ are locally
constant complexes with finite-dimensional stalk cohomology.
They are used to prove local constancy; they are not additional
choices in the definition of transport.

There are finitely many $\vct$ of weight $P$.
The formulas in Step 2 make $\mathcal E_0$ and
$\mathcal E_{U_P/T}$ perfect locally constant complexes.
The compact-support triangle in
\eqref{eq:conclusion-relative-triangle} and stability of
$\operatorname{Loc}(T;\mathcal D(\F))$ then show inductively
that every $\mathcal E_P$ is perfect and locally constant.
Indeed, on a sufficiently small contractible ball the two
outer terms are constant complexes, and their connecting
morphism is a morphism between those constant complexes.
Its fiber is constant as well and is the middle term of the
triangle. This uses the actual localization morphism, including
its homotopies, rather than only the dimensions of stalk
cohomology.
Consequently the same holds for $\mathcal B_P=\mathbb D_T\mathcal E_P$.
Put
\[
 \mathcal B_P^{\mathrm{loc}}:=\mathcal B_P[-2d].
\]
For a locally constant complex $L$, complex orientation gives
$i_t^!L\simeq i_t^*L[-2d]$. Therefore
\begin{equation}\label{eq:relative-transport-ordinary-fiber}
 i_t^*\mathcal B_P^{\mathrm{loc}}\simeq i_t^!\mathcal B_P\simeq B_{X_t,P},
 \qquad
 \mathcal B_P^{\mathrm{loc}}\simeq
 R\mathcal Hom_T(\mathcal E_P,\underline\F_T).
\end{equation}
Equivalently, $i_t^*\mathcal B_P\simeq B_{X_t,P}[2d]$.
The same base shift occurs at every $P$, so these are
identifications of the entire pole diagrams.

On a contractible coordinate ball $V\subset T$, evaluation at
$t\in V$ is an equivalence
$\operatorname{Loc}(V;\mathcal D(\F))\simeq\mathcal D(\F)$.
Applying $\operatorname{Fun}(N(\Z_{\ge0}),-)$ trivializes
$(\mathcal B_P^{\mathrm{loc}})|_V$ together with every pole arrow and all its
composition homotopies. Thus one does not intersect infinitely
many stagewise trivializing neighborhoods. Its homotopy colimit
\[
 \mathcal B^{\mathrm{loc}}:=\hocolim_P\mathcal B_P^{\mathrm{loc}}
       \quad\text{in }\operatorname{Loc}(T;\mathcal D(\F))
\]
has fiber $B_{X_t}:=\hocolim_P B_{X_t,P}$, the total derived
loop object of $X_t$.
The corresponding functor from the fundamental homotopy type
of $T$ gives the asserted coherent path transport of filtered
complexes. Taking cohomology gives the contact filtration
$F_PLH^*(X_t)=\operatorname{im}(H^*B_{X_t,P}\to LH^*(X_t))$.

\emph{Step 4: fusion, rotation, and the unit over $T$.}
Let $\mathcal Z_{P,Q}^{T}\to T\times\A_s^1$ be the
1-parameter binary collision family for $\mathcal X/T$.
The source coordinates $z-s,z+s$ on $s\ne0$ and $z$ on $s=0$
give its separated-input and output identifications. On the
compact-support companions these are
\[
 \mathcal E^{\mathrm{coll}}_{P,Q}|_{T\times\Gm}
    \simeq\operatorname{pr}_T^{-1}(\mathcal E_P\otimes\mathcal E_Q),
 \qquad
 i_0^{-1}\mathcal E^{\mathrm{coll}}_{P,Q}
    \simeq\mathcal E_{P+Q},
 \quad i_0(t)=(t,0).
\]
The comparisons are the relative versions of
Lemma~\ref{lem:source-coordinate-factorization}(i)--(ii), supplied by
the source-algebra identities and the specified finite-factor
comparisons of Lemma~\ref{lem:varying-target-finite-factors}(ii)--(iii).

Put $H=[0,\infty)$ and $L=(0,\infty)$, and write
\[
 j:T\times L\hookrightarrow T\times H,\qquad
 i_0:T\hookrightarrow T\times H,\quad i_0(t)=(t,0),
 \qquad E_H=\mathcal E^{\mathrm{coll}}_{P,Q}|_{T\times H},
 \quad E_L=j^{-1}E_H.
\]
For $\eta_j:E_H\to Rj_*j^{-1}E_H$, the unit of
$j^{-1}\dashv Rj_*$, define
\[
 \mathfrak c^E_{P,Q}:\mathcal E_{P+Q}
 \xrightarrow{\ \sim\ }i_0^{-1}E_H
 \xrightarrow{i_0^{-1}\eta_j}i_0^{-1}Rj_*E_L
 \xrightarrow{\ \sim\ }\mathcal E_P\otimes\mathcal E_Q.
\]
The last equivalence is the collar comparison of
Lemma~\ref{lem:parameter-chain-properties}(i). Explicitly,
the input identification makes
$E_L=\operatorname{pr}_T^{-1}(\mathcal E_P\otimes\mathcal E_Q)$;
on a contractible coordinate ball in $T$, sections over a
short punctured half-interval $(0,\varepsilon)$ therefore
identify $i_0^{-1}Rj_*E_L$ with $\mathcal E_P\otimes\mathcal E_Q$.
The comparison is natural under restriction of the ball, so it
glues on $T$. At $t\in T$, this comparison and the corresponding
one for the family over $\{t\}$ both identify with
\[
 R\Gamma\!\left((0,\varepsilon),
   \underline{i_t^*\mathcal E_P\otimes i_t^*\mathcal E_Q}\right)
 \simeq i_t^*\mathcal E_P\otimes i_t^*\mathcal E_Q,
 \qquad \varepsilon>0,
\]
where the underlined complex is constant on $(0,\varepsilon)$.
Naturality of the adjunction unit identifies the restriction
maps under these equivalences. On each fiber the displayed adjunction arrow is
the collar restriction used in the comparison preceding part (i) of
Lemma~\ref{lem:source-specialization-comparison}; its dual is
the normalized specialization with the inward-ray orientation.

Dualizing on $T$ gives a map
\[
 \Delta_T^!(\mathcal B_P\boxtimes\mathcal B_Q)
   \simeq\mathbb D_T(\mathcal E_P\otimes\mathcal E_Q)
   \xrightarrow{\mathbb D_T\mathfrak c^E_{P,Q}}
          \mathcal B_{P+Q},
 \qquad \Delta_T:T\hookrightarrow T\times T.
\]
Since the $\mathcal B_P$ are perfect local systems,
\[
 \Delta_T^!(\mathcal B_P\boxtimes\mathcal B_Q)
    \simeq(\mathcal B_P\otimes\mathcal B_Q)[-2d].
\]
Shifting the last arrow by $[-2d]$ therefore produces
$\mathfrak m_{P,Q}:\mathcal B_P^{\mathrm{loc}}\otimes\mathcal B_Q^{\mathrm{loc}}
\to\mathcal B_{P+Q}^{\mathrm{loc}}$. By
\eqref{eq:relative-transport-ordinary-fiber}, its stalk is the
fusion map on the pole stages, with the inward-ray sign fixed
in Section~\ref{sec:fusion-construction}.

Similarly the relative source automorphisms
\[
 r_\lambda:S_P\longrightarrow S_P,
 \qquad (r_\lambda\gamma)(\tau)=\gamma(\lambda\tau),
 \qquad\lambda\in S^1,
\]
and their identities $r_\lambda r_\eta=r_{\lambda\eta}$ give
a coherent circle action on $(\mathcal B_P^{\mathrm{loc}})_P$ in local
systems on $T$. This is the construction of
Proposition~\ref{prop:coherent-source-rotation}, using finite-factor
presentations and comparison maps over $T$. Its parameter-chain map and the
counterclockwise fundamental class $\sigma\in H_1(S^1;\F)$
give $\mathfrak d_P:\mathcal B_P^{\mathrm{loc}}\to\mathcal B_P^{\mathrm{loc}}[-1]$.
Finally $\mathcal B_0^{\mathrm{loc}}\simeq
R(\pi|_{\mathcal X})_*\underline\F_{\mathcal X}$, so the
adjunction unit gives $\mathfrak e:\underline\F_T\to\mathcal B_0^{\mathrm{loc}}$.
We have thus constructed, rather than assumed, the morphisms
\begin{equation}\label{eq:relative-transport-operations}
 \begin{gathered}
 \mathfrak m_{P,Q}:\mathcal B_P^{\mathrm{loc}}\otimes\mathcal B_Q^{\mathrm{loc}}
       \longrightarrow\mathcal B_{P+Q}^{\mathrm{loc}},\\
 \mathfrak d_P:\mathcal B_P^{\mathrm{loc}}\longrightarrow\mathcal B_P^{\mathrm{loc}}[-1],\qquad
 \mathfrak e:\underline\F_T\longrightarrow\mathcal B_0^{\mathrm{loc}}.
 \end{gathered}
\end{equation}
They preserve the pole-map comparisons because all constructions
were made on the whole relative diagrams. Their stalks are the
previously constructed fusion, BV, and unit maps.

Morphisms of local systems commute with path transport.
Apply this to \eqref{eq:relative-transport-operations}, and
pass to the pole homotopy colimit by
Lemma~\ref{lem:filtered-operations-colimit}. On cohomology the
fiber identities of Theorem~\ref{thm:loop-bv} give the unital
BV-algebra identities, and transport preserves them and the
filtration. Common-refinement equivalences intertwine every
map used above, giving the stated finite-factor presentation independence.
\end{proof}

In particular, for $t_0\in T$ the construction gives a monodromy
representation
\[
 \pi_1(T,t_0)\longrightarrow
 \operatorname{Aut}_{\mathrm{filtered\ unital\ BV}}
                       (LH^*(X_{t_0})),
\]
and its integer-graded counterpart when $\Omega$ is specified.

\clearpage
\part{Mirror symmetry via loop spaces}\label{part:mirror}
\label{part:gs-comparison}

\section{Fundamental classes and the contact spectral sequence}
\label{sec:log-cy-zero}

This section proves Theorems~\ref{thm:intro-lcy-zero} and
\ref{thm:intro-theta-basis}. Let $(M,D)$ be a positive log
Calabi--Yau pair in the sense of
Subsection~\ref{subsec:intro-positive-lcy}, and put $X=M\setminus D$.
Choose an ample boundary divisor
$A=\sum_{i=1}^r k_iD_i$, $k_i\in\Z_{>0}$, so that $(M,D,A)$
is polarized. Choose separately a trivialization of $K_M(D)$,
represented by a meromorphic volume form $\Omega$ satisfying
\begin{equation}\label{eq:lcy-volume}
 \operatorname{div}_M(\Omega)=-D.
\end{equation}
We use its integer grading and coefficients $\F$. The contact
shift of \eqref{eq:volume-grading} is then
$\mu_\Omega(\vct)=0$ for every realized contact vector $\vct$.

\begin{definition}[Constant contact classes]
\label{def:lcy-constant-contact-classes}
Use the logarithmic ring of Definition~\ref{def:log-cohomology}.
For $\vct\ne0$, put $I=I(\vct)$. For
$C\in\pi_0(D_I^\circ)$, let
$\pi_I:S_{I(\vct)}^\circ\to D_{I(\vct)}^\circ$ be the normal
circle-torus bundle and define
\[
 e_{\vct,C}:=\pi_I^*(1_C)t^{\vct}
       \in H^0_{\log}(M,D),
 \qquad e_0:=1_Xt^0.
\]
Here $1_C$ is the characteristic function of $C$.
For a connected component $\overline C\subset D_I$, put
$\overline C^\circ=\overline C\cap D_I^\circ$.
The assignment $\overline C\mapsto\overline C^\circ$ identifies
$\pi_0(D_I)$ with $\pi_0(D_I^\circ)$: each smooth connected
component of $D_I$ is irreducible, and its deeper boundary is a
finite union of proper divisors. We may therefore index by either
component, with $e_{\vct,\overline C}:=e_{\vct,\overline C^\circ}$.
For $I=\varnothing$ this convention reads
$e_{0,M}=e_{0,X}=e_0$.
The fibers of $\pi_I$ are connected, so these classes form the basis
\begin{equation}\label{eq:lcy-log-zero-basis}
 H^0_{\log}(M,D)=\bigoplus_{\vct,C}\F e_{\vct,C},
 \qquad w_A(e_{\vct,C})=\sum_i k_iv_i,
\end{equation}
where the zero-contact label is $(0,X)$ and every element has
finite support. For the label $p=(\vct,C)$ we also write
$e_p:=e_{\vct,C}$; the distinguished label $0$ denotes $e_0$.
\end{definition}

There are no negative total degrees on the logarithmic first
page. This excludes differentials entering degree zero, but not
those leaving it. We construct a canonical lift of each
$e_{\vct,C}$ at its own pole weight.

\subsection{Fundamental classes in a finite closed filtration}

We first construct fundamental classes for a finite closed
filtration of a finite-type scheme under an equidimensionality assumption. 
The following two lemmas are basic and not new.

\begin{lemma}
\label{lem:top-dimensional-filtration}
Let $Y$ be a separated finite-type complex scheme with a finite
closed filtration
\[
 \varnothing=Y_{-1}\subset Y_0\subset\cdots\subset Y_P=Y.
\]
Suppose that $L_p=(Y_p\setminus Y_{p-1})_{\mathrm{red}}$ is
smooth of pure complex dimension $\delta\ge0$, or empty, for every $p$.
All sheaf operations below are on analytic realizations, so
$j_p:L_p(\C)\hookrightarrow Y_p(\C)$ is an open inclusion.
Set
\[
 \begin{gathered}
 \omega_{Y_p}^{\mathrm n}:=\omega_{Y_p}[-2\delta],\qquad
 C_p:=R\Gamma(Y_p,\omega_{Y_p}^{\mathrm n}),\\
 \omega_{Y_{-1}}^{\mathrm n}:=0,\qquad C_{-1}:=0.
 \end{gathered}
\]
With $i_p:Y_{p-1}\hookrightarrow Y_p$, localization gives
\[
 \begin{gathered}
 (i_p)_*\omega_{Y_{p-1}}^{\mathrm n}\longrightarrow \omega_{Y_p}^{\mathrm n}
 \longrightarrow Rj_{p*}\F_{L_p}\xrightarrow{+1},\\
 C_{p-1}\longrightarrow C_p
 \longrightarrow R\Gamma(L_p;\F)\xrightarrow{+1}.
 \end{gathered}
\]
Here $j_p^{-1}\omega_{Y_p}^{\mathrm n}\simeq\F_{L_p}$ is induced by complex-oriented
smooth purity and the shift $[-2\delta]$ defining $\omega_{Y_p}^{\mathrm n}$. One has $\omega_{Y_p}^{\mathrm n}\in\mathcal D^{\ge0}(Y_p;\F)$ and
$C_p\in\mathcal D^{\ge0}(\F)$.

For $V\in\pi_0(L_p)$, let
$Z_V=\overline V^{\,\mathrm{red}}\subset Y_p$ and let
$a_V:Z_V\hookrightarrow Y_p$ be its closed inclusion.
The subvarieties $Z_V$, for all $p,V$, are precisely the
irreducible components of $Y$, without repetition. The map
\[
 \ell_p:H^0(L_p;\F)\longrightarrow H^0(C_p)
           =H^{BM}_{2\delta}(Y_p;\F),
 \qquad 1_V\longmapsto(a_V)_*[Z_V]
\]
is a section of the open-restriction map
$r_p:H^0(C_p)\to H^0(L_p;\F)$. Thus localization gives a
split exact sequence
\[
 0\longrightarrow H^0(C_{p-1})\longrightarrow H^0(C_p)
 \xrightarrow{r_p}H^0(L_p;\F)\longrightarrow0,
 \qquad r_p\ell_p=\mathrm{id}.
\]
\end{lemma}

\begin{proof}
Reduction leaves the analytic spaces and their Borel--Moore and
sheaf complexes unchanged. A connected smooth complex variety is
irreducible, so $V$ is a dense open subset of its closure $Z_V$,
and $\dim Z_V=\delta$. The finitely many pieces $V$ cover $Y$.
Every irreducible component of $Y$ is therefore one of their
closures: an irreducible closed set contained in a finite union
of these closures is contained in one of them. All these closures
have dimension $\delta$, so no one is properly contained in
another. Two distinct pieces cannot have the same closure, since
they would be disjoint dense opens of an irreducible variety.
The same argument applies to $Y_p$. Closedness of $Y_p\subset Y$
identifies the closures taken in these two spaces.

Remove from a nonempty $Y_p$ the singular loci of its irreducible
components and their pairwise intersections, and denote the removed
closed subset by $T_p$. If $\delta=0$, then $T_p=\varnothing$;
otherwise $\dim_\C T_p\le\delta-1$. In either case,
\[
 H^{BM}_{2\delta}(T_p;\F)
 =H^{BM}_{2\delta-1}(T_p;\F)=0.
\]
The localization long exact sequence therefore identifies
$H^{BM}_{2\delta}(Y_p;\F)$ with the top Borel--Moore homology
of the remaining disjoint smooth dense opens. Their
complex orientations give a basis consisting of the component
classes $[Z_V]$. Restriction to $L_p$ sends $(a_V)_*[Z_V]$ to
$1_V$ and to zero on its other components. This proves
$r_p\ell_p=\mathrm{id}$.

Locally finite triangulations of $Y_p(\C)$ have real dimension
at most $2\delta$, so its local homology vanishes in degrees
greater than $2\delta$. Equivalently, the dimension bound for
the Verdier dualizing complex gives
$\mathcal H^j(\omega_{Y_p})=0$ for $j<-2\delta$, so
$\omega_{Y_p}^{\mathrm n}\in\mathcal D^{\ge0}$. Derived global sections preserve
connectivity, giving the assertion for $C_p$. For the closed-open
pair $(Y_{p-1},L_p)$, one has
$i_p^!\omega_{Y_p}^{\mathrm n}=\omega_{Y_{p-1}}^{\mathrm n}$ and $j_p^{-1}\omega_{Y_p}^{\mathrm n}\simeq\F_{L_p}$.
Their localization cofiber sequence is the first displayed
sequence, and $R\Gamma$ gives the second. Since
$H^{-1}(L_p;\F)=0$, its degree-zero long exact sequence and
$r_p\ell_p=\mathrm{id}$ give the asserted split short exact
sequence.
\end{proof}

\begin{lemma}
\label{lem:fundamental-smooth-comparison}
Let $Y$ and $Y'$ be separated finite-type complex schemes, let
\[
 \varnothing=Y_{-1}\subset Y_0\subset\cdots\subset Y_P=Y
\]
be a finite closed filtration whose reduced layers
$L_p=(Y_p\setminus Y_{p-1})_{\mathrm{red}}$ are smooth of pure
complex dimension $\delta\ge0$ or empty, and let $f:Y'\to Y$
be smooth of pure relative dimension $e\ge0$. Put
\[
 \begin{gathered}
 Y'_p=Y'\times_Y Y_p,\qquad f_p:Y'_p\to Y_p,\\
 L'_p=(Y'_p\setminus Y'_{p-1})_{\mathrm{red}},\qquad
 g_p:L'_p\to L_p,\\
 \omega_{Y_p}^{\mathrm n}=\omega_{Y_p}[-2\delta],\qquad
 \omega_{Y'_p}^{\mathrm n}=\omega_{Y'_p}[-2(\delta+e)],\\
 C_p=R\Gamma(Y_p,\omega_{Y_p}^{\mathrm n}),\qquad
 C'_p=R\Gamma(Y'_p,\omega_{Y'_p}^{\mathrm n}),\qquad C_{-1}=C'_{-1}=0.
 \end{gathered}
\]
All sheaf operations are on analytic realizations. Each $L'_p$ is
smooth of pure dimension $\delta+e$ or empty. Complex orientation
gives $f_p^{-1}\omega_{Y_p}^{\mathrm n}\simeq \omega_{Y'_p}^{\mathrm n}$ and hence a normalized
smooth-pullback map
\[
 f_p^\sharp:C_p\longrightarrow
 R\Gamma(Y'_p,f_p^{-1}\omega_{Y_p}^{\mathrm n})\xrightarrow{\sim}C'_p.
\]
On Borel--Moore homology this is the smooth Gysin map with its
$2e$ degree shift absorbed by the normalization. These maps give a
morphism of localization cofiber sequences
\[
 \begin{array}{ccccc}
 C_{p-1}&\longrightarrow&C_p&\longrightarrow&R\Gamma(L_p;\F)\\
 {\scriptstyle f_{p-1}^\sharp}\downarrow&&
 {\scriptstyle f_p^\sharp}\downarrow&&{\scriptstyle g_p^*}\downarrow\\
 C'_{p-1}&\longrightarrow&C'_p&\longrightarrow&R\Gamma(L'_p;\F),
 \end{array}
\]
including their connecting maps to $C_{p-1}[1]$ and $C'_{p-1}[1]$.

For $p\le q$, write $i_{p,q}:Y_p\hookrightarrow Y_q$ and
$i'_{p,q}:Y'_p\hookrightarrow Y'_q$ for the closed stage inclusions.
The square
\[
 \begin{array}{ccc}
 Y'_p&\xrightarrow{i'_{p,q}}&Y'_q\\
 {\scriptstyle f_p}\downarrow&&\downarrow{\scriptstyle f_q}\\
 Y_p&\xrightarrow{i_{p,q}}&Y_q
 \end{array}
\]
is Cartesian. Let $(i_{p,q})_*:C_p\to C_q$ be induced by the
proper trace $(i_{p,q})_*\omega_{Y_p}^{\mathrm n}\to \omega_{Y_q}^{\mathrm n}$, and define
$(i'_{p,q})_*:C'_p\to C'_q$ similarly. Then
\[
 f_q^\sharp\circ(i_{p,q})_*
   =(i'_{p,q})_*\circ f_p^\sharp
 \quad\text{in }D(\F).
\]

For a connected component $V\subset L_p$, set
$Z_V=\overline V^{\,\mathrm{red}}\subset Y_p$ and write
$a_V:Z_V\hookrightarrow Y_p$; define $Z_{V'}$ and $a_{V'}$
similarly for $V'\in\pi_0(L'_p)$. Define
\[
 \begin{aligned}
 \ell_p:H^0(L_p;\F)&\longrightarrow H^0(C_p),&
 \ell_p(1_V)&=(a_V)_*[Z_V],\\
 \ell'_p:H^0(L'_p;\F)&\longrightarrow H^0(C'_p),&
 \ell'_p(1_{V'})&=(a_{V'})_*[Z_{V'}].
 \end{aligned}
\]
These are the sections of Lemma~\ref{lem:top-dimensional-filtration},
and they satisfy
\[
 \begin{gathered}
 H^0(f_p^\sharp)\ell_p=\ell'_p H^0(g_p^*),\\
 H^0(f_p^\sharp)\bigl((a_V)_*[Z_V]\bigr)
   =\sum_{V'\in\pi_0(g_p^{-1}V)}(a_{V'})_*[Z_{V'}].
 \end{gathered}
\]
Closed pushforward sends $(a_V)_*[Z_V]$ to the class of the same
reduced closure in $Y_q$. If $Y''$ is a separated finite-type
complex scheme and $h:Y''\to Y'$ is smooth of pure relative
dimension $d$, put $Y''_p=Y''\times_{Y'}Y'_p$ and let
$h_p:Y''_p\to Y'_p$ be the projection. With
$\omega_{Y''_p}^{\mathrm n}=\omega_{Y''_p}[-2(\delta+e+d)]$, one has
\[
 (f_p\circ h_p)^\sharp=h_p^\sharp\circ f_p^\sharp.
\]
Thus the localization maps, closed pushforwards, and
fundamental-class sections are compatible with compositions of
smooth Cartesian changes of finite-factor presentations.
\end{lemma}

\begin{proof}
The maps $f_p$ and $g_p$ are smooth of relative dimension $e$
by base change. Indeed, smooth pullback of a reduced scheme is
reduced, so
$L'_p=Y'_p\times_{Y_p}L_p$. Consequently $L'_p$ is smooth
of dimension $\delta+e$, and
Lemma~\ref{lem:top-dimensional-filtration} applies to the
pulled-back filtration. The complex-oriented identifications
\[
 \omega_{Y'_p}\simeq f_p^!\omega_{Y_p}
       \simeq f_p^{-1}\omega_{Y_p}[2e]
\]
give $f_p^{-1}\omega_{Y_p}^{\mathrm n}\simeq \omega_{Y'_p}^{\mathrm n}$. In the Cartesian square of
closed stage inclusions, proper base change identifies
$f_q^{-1}(i_{p,q})_*\omega_{Y_p}^{\mathrm n}$ with
$(i'_{p,q})_*f_p^{-1}\omega_{Y_p}^{\mathrm n}$. Under the normalized purity
equivalences just obtained, the pullback of the counit $(i_{p,q})_*i_{p,q}^!\omega_{Y_q}^{\mathrm n}\to \omega_{Y_q}^{\mathrm n}$ is the
counit for $i'_{p,q}$. Taking global sections and using the
section-to-pullback maps defining $f_p^\sharp$ therefore gives
$f_q^\sharp(i_{p,q})_*=(i'_{p,q})_*f_p^\sharp$.
Applying this argument to $i_p:Y_{p-1}\hookrightarrow Y_p$
identifies the first arrows of the two localization cofiber
sequences. The induced map on cofibers is $g_p^*$ by smooth
base change for the open complements. Naturality of the
cofiber sequence also identifies the connecting maps.

The restriction $f_p^{-1}Z_V\to Z_V$ is smooth and therefore
open. Thus $g_p^{-1}V$ is dense in $f_p^{-1}Z_V$. The latter
scheme is reduced, since it is smooth over the reduced scheme
$Z_V$. Its irreducible components are exactly the reduced
closures $Z_{V'}$ of the connected components of $g_p^{-1}V$.
On these smooth dense open subsets, complex orientation sends
the fundamental class of $Z_V$ to that of its inverse image
with multiplicity one. The top Borel--Moore calculation in
Lemma~\ref{lem:top-dimensional-filtration} determines the class
on the whole inverse image and proves the displayed sum. Since
$g_p^*(1_V)=\sum_{V'\in\pi_0(g_p^{-1}V)}1_{V'}$, this also
proves $H^0(f_p^\sharp)\ell_p=\ell'_pH^0(g_p^*)$.
This is smooth pullback: if $f_p$ is finite \'etale, the
coefficient on each inverse-image component is still one.
The degree of a finite cover occurs in proper pushforward,
not in this pullback formula.

For a closed stage inclusion, $Z_V$ has the same reduced closure
in $Y_q$, and proper pushforward preserves its complex-oriented
class. Finally, complex orientations multiply under composites
of smooth maps, and the base-change transformations respect
these composites. This proves the composition formula with the
relative dimensions added, and all stated compatibilities.
\end{proof}

\subsection{Fundamental-class lifts and the theta basis}
\label{subsec:lcy-fundamental-lifts}

The construction below defines a theta class at its first pole stage
in \eqref{eq:lcy-theta-stage-class} and then in loop cohomology
in \eqref{eq:lcy-theta-colimit-class}. We first construct the
finite-factor sections needed in these definitions, the crucial observation in the log Calabi-Yau case being the  equidimensionality property of the next result. Recall from Section \ref{sec:formal-loops} that the volume-index functions $\delta$ are locally constant, that is, constant on each connected component.

\begin{proposition}[Equidimensional finite factors in the log Calabi--Yau case]
\label{prop:lcy-equidimensional-factors}
Fix the polarized positive log Calabi--Yau pair $(M,D,A)$ and
$\Omega$ satisfying \eqref{eq:lcy-volume}. Fix $P\ge0$.
The scheme $F_P\mathscr LX$ admits a finite nice \'etale atlas
whose charts have finite-factor presentations
\[
 \begin{gathered}
 h:U=Y\times\A^J\longrightarrow F_P\mathscr LX,\\
 U\times_{F_P\mathscr LX}F_p\mathscr LX
       =Y_p\times\A^J\quad(0\le p\le P),\\
 \varnothing=Y_{-1}\subset Y_0\subset\cdots\subset Y_P=Y,
 \end{gathered}
\]
with the following properties.
\begin{enumerate}[(i),wide=0pt]
\item The reduced layers admit decompositions
\[
 L_p:=(Y_p\setminus Y_{p-1})_{\mathrm{red}}
   =\coprod_{\substack{w_A(\vct)=p\\C,\beta}}
           V_{\vct,C,\beta}
\]
into finitely many connected smooth varieties. The inverse image
of $V_{\vct,C,\beta}$ in $U$ lies over the contact component
labeled by $(\vct,C)$. The index $\beta$ labels the connected components of the
finite contact factor on this chart. Since this factor is
smooth, these are also its irreducible components. For zero contact the label is $(0,X)$.
The volume indices satisfy
\[
 \delta_{Y_p}=\delta_Y|_{Y_p},\qquad
 \delta_Y|_{V_{\vct,C,\beta}}=\dim V_{\vct,C,\beta}.
\]
\item For a nonempty open-and-closed level set
$Y^{(d)}=\delta_Y^{-1}(d)$, put
$Y_p^{(d)}=Y_p\cap Y^{(d)}$. Then $d\ge0$, and every nonempty
$(Y_p^{(d)})_{\mathrm{red}}$ is pure of dimension $d$.
More precisely, the map
\[
 V\longmapsto\overline V^{\,\mathrm{red}}
\]
is a bijection from the connected contact pieces contained in
$\coprod_{j=0}^p(L_j\cap Y^{(d)})$ to
$\operatorname{Irr}((Y_p^{(d)})_{\mathrm{red}})$; the closure may
be taken in either $Y_p^{(d)}$ or $Y^{(d)}$.
\end{enumerate}
These finite-factor presentations can include any specified finite collection
of closed pole inclusions and finite composable strings of chart
overlap maps, with a factor for each chart occurrence. In particular, the
assertion applies to the compatible local finite factors used
for descent; it does not require a global finite factor for
$F_P\mathscr LX$.
\end{proposition}

\begin{proof}
Choose a finite nice \'etale atlas as in
Remark~\ref{rem:finite-nice-atlas}. On each chart retain the
finitely many coefficient functions cutting out the closed pole
stages $F_p\mathscr LX$, $0\le p\le P$. The construction in
Lemma~\ref{lem:finite-presentation-categories}(i) gives the displayed
closed filtration $Y_p\subset Y$ with a common affine tail and
compatible finite-factor presentations of each specified finite string of
overlap maps.

We may further retain the reduced contact pieces of each layer.
Indeed, if the coordinate ring of a presented chart is
$R_0[u_j:j\in J]$, with $R_0$ of finite type over $\C$, then
\[
 \sqrt{(0)}_{\,R_0[u_j:j\in J]}
       =\sqrt{(0)}_{\,R_0}\,R_0[u_j:j\in J].
\]
The same equality applies after quotienting by the retained pole
ideals. Each reduced layer has a finite principal-open cover,
and its clopen contact pieces are given on that cover by
idempotents. Retain the finitely many $u_j$ occurring in the
functions defining the cover. The idempotents then descend
to its finite factors by the
argument in \emph{Idempotents in polynomial algebras} in
Subsection~\ref{subsec:relative-finite-factor-proof},
so each contact piece is a product of a reduced locally
closed finite factor and the remaining affine tail,
as in the proof of
Proposition~\ref{prop:general-factor-contact-index}.
Proposition~\ref{prop:contact-volume-indices} makes those finite
factors smooth. Splitting them into connected components gives
the decomposition in (i); there are finitely many because each factor
is of finite type and only finitely many contact vectors have
weight at most $P$.

Closed-stage compatibility of the volume indices gives
$\delta_{Y_p}=\delta_Y|_{Y_p}$. On each connected contact piece,
\eqref{eq:contact-factor-volume-index} gives
\[
 \delta_Y|_{V_{\vct,C,\beta}}
  =\dim V_{\vct,C,\beta}+\tfrac12\mu_\Omega(\vct)
  =\dim V_{\vct,C,\beta},
\]
since \eqref{eq:lcy-volume} implies $\mu_\Omega(\vct)=0$.
Consequently every nonempty reduced layer of $Y_\bullet^{(d)}$
is smooth of pure dimension $d$. Applying
Lemma~\ref{lem:top-dimensional-filtration} to this closed
filtration proves (ii), including the description of every
irreducible component. All preceding retentions involve finite
lists of functions and maps and can therefore be made
simultaneously in each specified finite comparison diagram.
\end{proof}

\begin{proof}[Proof of Theorems~\ref{thm:intro-lcy-zero} and
\ref{thm:intro-theta-basis}]
\phantomsection\label{proof:lcy-zero-and-theta-basis}
\emph{Step 1: compatible equidimensional finite factors.}
Fix $P\ge0$ and put $S_p=F_p\mathscr LX(\C)$ for $0\le p\le P$.
Use the finite-factor presentations of
Proposition~\ref{prop:lcy-equidimensional-factors}(i)--(ii), restricting
each finite factor to a nonempty level set $\delta_Y=\delta$.
Write $Y_p$ for the resulting closed factors and
$V_{\vct,C,\beta}\subset(Y_p\setminus Y_{p-1})_{\mathrm{red}}$
for their connected contact pieces. Thus every nonempty $Y_p$
is pure of dimension $\delta$, and the reduced closures of
these pieces are its irreducible components. All comparisons
below present this finite closed filtration with a common affine tail.

\emph{Step 2: connectivity and degree-zero sheaf sections.}
For a finite factor $Y_p$ as above, the normalization is
$\omega_{Y_p}^{\mathrm n}=\omega_{Y_p}[-2\delta]$.
The dualizing complex on a complex variety of dimension $\delta$
has no cohomology sheaves below degree $-2\delta$. Hence
\[
 \omega_{Y_p}^{\mathrm n}\in\mathcal D^{\ge0}(Y_p;\F).
\]
For nonempty $Y_p$, let
$\nu_p:\widetilde{(Y_p)_{\mathrm{red}}}\to(Y_p)_{\mathrm{red}}$
be its normalization. On the common analytic space of $Y_p$
and its reduction, complex orientation gives
\[
 \mathcal H^0\omega_{Y_p}^{\mathrm n}
   \simeq(\nu_p)_*\underline\F_{\widetilde{(Y_p)_{\mathrm{red}}}}.
\]
Indeed, its stalk at $y\in Y_p(\C)$ is the local homology
$H_{2\delta}(Y_p(\C),Y_p(\C)\setminus\{y\};\F)$,
with one oriented generator for each local analytic irreducible
branch. These branches correspond to the points of the
normalization over $y$, and the identifications respect
restriction to smaller neighborhoods. Thus a global
irreducible-component class determines the section with
coefficient one on every local branch belonging to that
component and zero on the other branches.
Ordinary inverse image of sheaves is exact. The chart description
\eqref{eq:nice-chart-dualizing} and the jointly conservative nice
atlas therefore imply
$\omega_{p,\Omega}\in\mathcal D^{\ge0}(S_p;\F)$.
Derived global sections preserve connective objects, so
\[
 B_{\Omega,p}=R\Gamma(S_p,\omega_{p,\Omega})
       \in\mathcal D^{\ge0}(\F),\qquad
 \varepsilon_p:H^0(B_{\Omega,p})
       \xrightarrow{\sim}\Gamma(S_p,\mathcal H^0\omega_{p,\Omega}).
\]
The canonical edge isomorphism $\varepsilon_p$ follows by applying
$R\Gamma$ to the standard
truncation triangle
$\mathcal H^0\omega_{p,\Omega}\to\omega_{p,\Omega}
\to\tau_{\ge1}\omega_{p,\Omega}$.
Its last term has derived global sections in
$\mathcal D^{\ge1}(\F)$. The construction is natural for
morphisms of connective complexes and for the section-to-pullback
maps associated with ordinary inverse image.

\emph{Step 3: descent of the normalized fundamental classes.}
Fix a label $(\vct,C)$ of weight $p$. On a chart of $S_p$ define
\[
 Z_{\vct,C,\beta}
    :=\overline{V_{\vct,C,\beta}}^{\,\mathrm{red}}\subset Y_p,
 \qquad
 \xi_{Y_p,\vct,C}
    :=\sum_\beta(a_\beta)_*[Z_{\vct,C,\beta}]
           \in H_{2\delta}^{BM}(Y_p;\F)=H^0R\Gamma(Y_p,\omega_{Y_p}^{\mathrm n}).
\]
Here $a_\beta:Z_{\vct,C,\beta}\hookrightarrow Y_p$ is the
closed inclusion, and the sum is finite. Let
\[
 \varepsilon_{Y_p}:H^0R\Gamma(Y_p,\omega_{Y_p}^{\mathrm n})
       \xrightarrow{\sim}\Gamma(Y_p,\mathcal H^0\omega_{Y_p}^{\mathrm n})
\]
be the finite-factor edge isomorphism of Step 2. For a chart
$h_a:U_a\to F_p\mathscr LX$ with projection
$\pi_a:U_a\to Y_{a,p}$, define the section
\[
 \sigma_a:=\pi_a^{-1}
       \bigl(\varepsilon_{Y_{a,p}}(\xi_{Y_{a,p},\vct,C})\bigr)
 \in\Gamma\bigl(U_a^{\mathrm{an}},
           h_a^{-1}\mathcal H^0\omega_{p,\Omega}\bigr).
\]
Here $\pi_a^{-1}$ denotes pullback of the section, using
\eqref{eq:nice-chart-dualizing} to identify its target.
Lemma~\ref{lem:fundamental-smooth-comparison} gives
\[
 H^0(f_p^\sharp)(\xi_{Y_p,\vct,C})=\xi_{Y'_p,\vct,C}
\]
for every smooth Cartesian comparison of these closed diagrams.
Naturality of the edge isomorphisms makes the sections
$\sigma_a$ compatible with these finite-factor changes,
including \'etale comparisons; the sum over $\beta$ includes
all components with the fixed contact label, even when a
comparison disconnects an inverse image.

For two charts, put
$U_{ab}=U_a^{\mathrm{an}}\times_{S_p}U_b^{\mathrm{an}}$ and let
$q_a:U_{ab}\to U_a^{\mathrm{an}}$ and
$q_b:U_{ab}\to U_b^{\mathrm{an}}$ be the projections.
On nice charts of $U_{ab}$, pass to common finite-factor
presentations. The two sections pull back to the same displayed
sum of inverse-image component classes. Thus, under the descent
identifications of Proposition~\ref{prop:normalized-pole-sheaves},
\[
 q_a^{-1}\sigma_a=q_b^{-1}\sigma_b
 \quad\text{on }U_{ab}.
\]
Define $\sigma^{(p)}_{\vct,C}$ to be the unique section of
$\mathcal H^0\omega_{p,\Omega}$ whose pullback along every $h_a$
is $\sigma_a$. Ordinary sheaf descent gives its existence and
uniqueness: the nice \'etale maps are local homeomorphisms on
analytic realizations and jointly cover $S_p$, by \emph{Analytic
realization} in Subsection~\ref{subsec:coherent-descent}. On triple overlaps,
the pairwise equalities agree by transitivity of the normalized
purity equivalences and the descent identifications of the normalized sheaf.
We now \emph{define the pole-stage theta class} by
\begin{equation}\label{eq:lcy-theta-stage-class}
 \vartheta^{(p)}_{\vct,C}
   :=\varepsilon_p^{-1}(\sigma^{(p)}_{\vct,C})
       \in H^0(B_{\Omega,p}),\qquad p=w_A(\vct).
\end{equation}
Its finite-factor representative on $Y_p$ is
$\xi_{Y_p,\vct,C}$, in the sense of
Definition~\ref{def:finite-factor-representatives}.

Let
$\rho_p:H^0(B_{\Omega,p})\to H^0(\gr_p\mathbf B_\Omega)$ be the map
induced by the localization cofiber sequence
$B_{\Omega,p-1}\to B_{\Omega,p}\to\gr_p\mathbf B_\Omega$,
where the last object is identified with the localization
cofiber by the specified filtered comparison of
Theorem~\ref{thm:canonical-complex}.
The layer identification and the open restriction in
Lemma~\ref{lem:top-dimensional-filtration} give
\[
 H^0(\gr_p\mathbf B_\Omega)
   =\bigoplus_{w_A(\wct)=p,C'}\F e_{\wct,C'},
 \qquad
 \rho_p(\vartheta^{(p)}_{\vct,C})=e_{\vct,C}.
\]
The construction is independent of the chosen atlas and finite
factors: any two choices have the common comparisons used above.

\emph{Step 4: passage between pole stages.}
For $Q\ge p$, work on each nice chart of $S_Q$ and present
the pullback of $S_p\hookrightarrow S_Q$ by a closed inclusion
$i:Y_p\hookrightarrow Y_Q$ with a common affine tail on that
chart. Common smooth refinements compare these finite-factor presentations
on chart overlaps, as in Step~3.
Since $Y_p$ is closed,
$\overline{V_{\vct,C,\beta}}^{\,Y_p}
 =\overline{V_{\vct,C,\beta}}^{\,Y_Q}$.
Consequently $i_*[Z_{\vct,C,\beta}]$ is the same oriented
component class in $Y_Q$. Write $u_{p,Q}:B_{\Omega,p}\to B_{\Omega,Q}$ for the specified
proper-trace map. Define the class at any larger stage and its
colimit class by
\begin{equation}\label{eq:lcy-theta-colimit-class}
 \vartheta^{(Q)}_{\vct,C}
    :=H^0(u_{p,Q})(\vartheta^{(p)}_{\vct,C}),\qquad
 \vartheta_{\vct,C}
    :=\iota_Q(\vartheta^{(Q)}_{\vct,C})\in LH^0_\Omega(X),
\end{equation}
where $\iota_Q:H^0(B_{\Omega,Q})\to
\varinjlim_RH^0(B_{\Omega,R})=LH^0_\Omega(X)$ is the canonical
colimit map. The closed-pushforward formula just proved and
Lemma~\ref{lem:fundamental-smooth-comparison} show that
$\vartheta^{(Q)}_{\vct,C}$ is also the class obtained from the
same reduced closures directly in a finite-factor presentation of $S_Q$.
The identities
$H^0(u_{Q,R})\vartheta^{(Q)}_{\vct,C}
 =\vartheta^{(R)}_{\vct,C}$ for $R\ge Q\ge p$
make \eqref{eq:lcy-theta-colimit-class} independent of $Q$.

\emph{Step 5: the basis, degeneration, and the leading product.}
Since $H^{-1}(\gr_p\mathbf B_\Omega)=0$ and the displayed lifts make
$\rho_p$ surjective, the localization long exact sequence gives
\begin{equation}\label{eq:lcy-zero-stage-sequence}
 0\longrightarrow H^0(B_{\Omega,p-1})
 \longrightarrow H^0(B_{\Omega,p})
 \xrightarrow{\rho_p}
       \bigoplus_{w_A(\vct)=p,C}\F e_{\vct,C}
 \longrightarrow0.
\end{equation}
For $p=0$ the first term is zero. The section sending
$e_{\vct,C}$ to $\vartheta^{(p)}_{\vct,C}$ splits this sequence.
Induction on $p$, followed by the exact filtered colimit, yields
\begin{equation}\label{eq:lcy-zero-splitting}
 LH^0_\Omega(X)=\bigoplus_{\vct,C}\F\vartheta_{\vct,C},\qquad
 F_PLH^0_\Omega(X)
   =\bigoplus_{w_A(\vct)\le P,C}\F\vartheta_{\vct,C}.
\end{equation}
The transition maps
$H^0(B_{\Omega,p-1})\to H^0(B_{\Omega,p})$ are injective, so
these are also the images defining the filtration on $LH^0$.
Only finitely many labels occur at each finite pole bound.
Every class in the filtered colimit comes from such a bound,
which explains the direct sum in \eqref{eq:lcy-zero-splitting}.
Under the specified comparison
$H^0(F_p\mathbf B_\Omega)\simeq H^0(B_{\Omega,p})$, choose a
closed representative
$x\in(F_p\mathbf B_\Omega)^0$ of the pole-stage theta class
$\vartheta^{(p)}_{\vct,C}$. Thus
$d_{\mathbf B_\Omega}x=0$, and the image of $x$ in the
associated graded represents $e_{\vct,C}\in E_1^{-p,p}$.
Every outgoing differential of this class is computed from
$d_{\mathbf B_\Omega}x$, so all these differentials vanish.
An incoming differential would originate in total degree
$-1$, which is zero already on $E_1$. Thus
\begin{equation}\label{eq:lcy-zero-permanence}
 E_1^{p,-p}=E_\infty^{p,-p}\qquad(p\in\Z).
\end{equation}

We now identify the multiplication on $\gr_F LH^0_\Omega(X)$
and establish vanishing in negative degrees.
Theorem~\ref{thm:summary} identifies
$\bigoplus_{p\ge0}E_1^{-p,p}$ with $H^0_{\log}(M,D)$
as contact-weight graded algebras. More explicitly,
Lemma~\ref{lem:fusion-leading-product}, through
\eqref{lp:eq:fusion-leading-product}, identifies the product
with pullback to the common face followed by cup product.
Since $\mu_\Omega=0$, the degree-zero contact classes are
locally constant functions. For two component characteristic
functions, their product is therefore $1$ on precisely those
connected output components on which both pulled-back
functions equal $1$, and is $0$ on the others; an empty
common face gives zero.

By \eqref{eq:lcy-zero-permanence}, every degree-zero class
survives from $E_1$ to $E_\infty$. These identifications
preserve products because the spectral sequence is
multiplicative. Its multiplicative convergence consequently
gives
\[
 H^0_{\log}(M,D)
 \cong \bigoplus_{p\ge0}E_\infty^{-p,p}
 \cong \gr_F LH^0_\Omega(X)
\]
as contact-weight graded rings.

For negative degrees, Step~2 gives
$H^j(B_{\Omega,p})=0$ for every $p\ge0$ and $j<0$.
Exactness of filtered colimits gives
$LH^j_\Omega(X)\cong\varinjlim_p H^j(B_{\Omega,p})$,
so these groups vanish. We have therefore proved
\begin{equation}\label{eq:lcy-zero-graded-ring}
 LH^j_\Omega(X)=0\quad(j<0),\qquad
 \gr_F LH^0_\Omega(X)\cong H^0_{\log}(M,D),
\end{equation}
where the latter isomorphism preserves multiplication and
contact weight.

\emph{Step 6: the polarization and the unit.}
For another ample boundary divisor $A'=\sum_i k'_iD_i$ with
$k'_i>0$, choose integers $c,c'>0$ satisfying
\[
 w_{A'}(\vct)\le c\,w_A(\vct),\qquad
 w_A(\vct)\le c'\,w_{A'}(\vct)
 \quad(\vct\in\N^r).
\]
For the coefficient-scheme structures, enlarge these constants
as follows. Let $bA$ and $b'A'$ be the very ample multiples
used for the two coefficient presentations, with affine generators
$\mathbf g=(g_i)$ and $\mathbf h=(h_j)$ of $\mathcal O(X)$,
respectively. Choose polynomial expressions
\[
 h_j=\Phi_j(\mathbf g),\quad\deg\Phi_j\le d,\qquad
 g_i=\Psi_i(\mathbf h),\quad\deg\Psi_i\le e,
 \qquad d,e\ge1,
\]
and require in addition $b'c\ge db$ and $bc'\ge eb'$.
For $P\ge0$, any complex algebra $R$, and any loop
$\gamma\in X(R(\!( t )\!))$, the condition
$g_i(\gamma)\in t^{-bP}R\llbracket t \rrbracket$ for all $i$ then implies
$h_j(\gamma)\in t^{-b'cP}R\llbracket t \rrbracket$ for all $j$.
Similarly, $h_j(\gamma)\in t^{-b'P}R\llbracket t \rrbracket$ for all $j$
implies $g_i(\gamma)\in t^{-bc'P}R\llbracket t \rrbracket$ for all $i$, using
the $\Psi_i$. These polynomial estimates work also for
nonreduced $R$; numerical contact orders alone would only
control the residue-field loops. They give the cofinal
finitely presented closed inclusions
\[
 F_P^A\mathscr LX\hookrightarrow F_{cP}^{A'}\mathscr LX,
 \qquad
 F_P^{A'}\mathscr LX\hookrightarrow F_{c'P}^{A}\mathscr LX
\]
of Proposition~\ref{prop:fusion-invariance-gradings}.
Here the superscript records the chosen polarization and its
coefficient presentation. Closedness and finite presentation
follow by imposing the finitely many coefficients below the
source bound in the reverse polynomial expressions on the
target stage.

For a fixed contact label, place stages of both systems in a
larger stage. On each nice chart of that stage, present the
pulled-back closed inclusions and retain the finitely many layer equations,
principal opens, and contact idempotents from both systems.
They give finite factors with a common affine tail on this
chart. The reduced closure of each labeled contact piece is
the same in the two closed subfactors and the containing factor,
by closedness, as in Step~4. Proper trace and the normalized
purity equivalences identify its component classes; further smooth refinements
make these identifications agree on overlaps. Therefore the canonical unfiltered
comparison satisfies
\[
 \vartheta^A_{\vct,C}\longmapsto\vartheta^{A'}_{\vct,C}.
\]
It respects fusion by Proposition~\ref{prop:fusion-invariance-gradings}.
Hence the labeled basis and all its coefficient functionals are
independent of the polarization; \eqref{eq:lcy-zero-splitting}
records each polarization's filtration through its weight function.
At zero contact, the construction gives the constant class
$1_X\in H^0(X)$. Since $X$ is connected and its image is the
fusion unit by Theorem~\ref{thm:general-associativity},
$\vartheta_0=1$. This completes both proofs.
\end{proof}

The fundamental-class construction uses only the closed pole
diagrams, their normalized dualizing complexes, and the contact
volume indices. Fusion is used for the ring assertion and the
identification of the unit.

\begin{corollary}[Stanley--Reisner presentation for connected
boundary intersections]
\label{cor:lcy-stanley-reisner}
Assume every nonempty $D_I$ is connected.
Then
\begin{equation}\label{eq:lcy-stanley-reisner}
 \gr_FLH^0_\Omega(X)
 \cong
 \frac{\F[t_1,\ldots,t_r]}
 {\left(\prod_{i\in I}t_i:\ D_I=\varnothing\right)},
 \qquad w(t_i)=k_i .
\end{equation}
The ring on the right is the Stanley--Reisner ring
$\F[\Delta(D)]$ of the simplicial complex
\[
 \Delta(D):=\{I\subseteq\{1,\ldots,r\}:D_I\ne\varnothing\};
\]
see Definition~\ref{mir:def-stanley-reisner}.
Its generators have cohomological degree zero and
contact weights $w(t_i)=k_i$.\end{corollary}
\begin{proof}
Under the connectedness assumption there is one constant
coefficient class for each realized contact vector. Common-face
convolution multiplies its monomials by adding their vectors
when the union of their supports is a face, and gives zero
otherwise. This is the displayed presentation of
$H^0_{\log}(M,D)$; apply
Theorem~\ref{thm:intro-lcy-zero}.
\end{proof}

\begin{remark}[A smooth boundary]
\label{rem:lcy-smooth-boundary}
Suppose $D$ is nonempty, smooth, and connected. Then $A=kD$ for some
$k\in\Z_{>0}$, and Corollary~\ref{cor:lcy-stanley-reisner} gives
$\gr_FLH^0_\Omega(X)=\F[t]$ with $w(t)=k$.
Put $s=\vartheta_{1,D}$. Its symbol is $t$, so $s^m$ has
symbol $t^m$. A nonzero polynomial relation in $s$ would have a
nonzero highest-weight symbol; hence no such relation exists.
Conversely, subtracting a multiple of $s^m$ cancels the highest
symbol of any element of weight $km$. Induction on the weight
therefore gives
\[
 LH^0_\Omega(X)=\F[s],\qquad
 F_PLH^0_\Omega(X)
    =\operatorname{span}_{\F}\{1,s,\ldots,s^{\lfloor P/k\rfloor}\}.
\]
For $m\ge1$, this proves only the triangular relation
\[
 \vartheta_{m,D}
       =s^m+\sum_{0\le j<m}a_{mj}s^j,\qquad a_{mj}\in\F;
\]
it does not assert $\vartheta_{m,D}=s^m$.
\end{remark}

\begin{remark}
\label{rem:lcy-zero-scope}
Degree-zero degeneration means that all differentials entering
or leaving the diagonal $p+q=0$ vanish. The proof does not give
vanishing on any other diagonal: a positive-degree class on a
contact locus need not extend to its closure. Also,
\eqref{eq:lcy-zero-splitting} is an additive splitting. Products
of its fundamental classes can have terms of smaller contact
weight, so it does not identify the full fusion ring with
$H^0_{\log}(M,D)$.
\end{remark}

\section{The affine loop mirror and its Calabi--Yau property}
\label{mir:sec-affine-scheme}
\label{mir:sec:loop-cy}

Our main goal in this section is to prove Theorem~\ref{mir:thm-intro-cy}.
We first construct the affine loop mirror and its boundary
degeneration. Under maximal boundary, a refinement makes the special
fiber a Gorenstein face scheme with canonical module of shift zero
and semi-log-canonical singularities. A general Rees argument then
lifts both properties to the affine loop mirror.

\subsection{The affine loop mirror and its boundary degeneration}

Retain the positive log Calabi--Yau setting of
Section~\ref{sec:log-cy-zero}, and fix the
coefficient embedding $\F\hookrightarrow\C$.  Thus
$D=\bigcup_{i=1}^rD_i$ is reduced SNC, $K_M+D\sim0$, and the
contact weights $k_i$ are positive integers.  We write
\[
 R_X=LH^0_\Omega(X;\F),\qquad
 \check X=\Spec(R_X\otimes_{\F}\C).
\]
We call $\check X$ the \emph{affine loop mirror} of $(X,\Omega)$. 
If a geometric mirror $Y$ is given, we expect that $R_X\otimes_{\F}\C\simeq\Gamma(Y,\mathcal O_Y)$, and thus that
$\check X$ is its affinization
$Y^{\mathrm{aff}}:=\Spec\Gamma(Y,\mathcal O_Y)$.
Theorem~\ref{thm:intro-theta-basis} provides the normalized
fundamental classes $\vartheta_{\vct,C}$ and the filtered basis
\eqref{eq:lcy-zero-splitting}. Put $A_D:=H^0_{\log}(M,D)$.
Here $F_{-1}R_X=0$, $F_0R_X=\F1$, the filtration is exhaustive and multiplicative,
and each $F_pR_X$ is finite-dimensional, by positivity of the $k_i$.
For the algebraic calculations below use closed-stratum component
labels as in Definition~\ref{def:lcy-constant-contact-classes}.
At zero contact take $D_\varnothing=M$ and
$\vartheta_{\zero,M}=1$.

Use the constant classes $e_{\vct,C}$ of
Definition~\ref{def:lcy-constant-contact-classes} as a basis of
$A_D$. If $I=I(\vct)$ and $J=I(\wct)$, its product is
\begin{equation}\label{mir:eq:disconnected-face-product}
 e_{\vct,C}e_{\wct,C'}=
 \sum_{\substack{E\in\pi_0(D_{I\cup J})\\E\subset C\cap C'}}
             e_{\vct+\wct,E}.
\end{equation}
This follows by restricting constant classes to the components of
$D_{I\cup J}$ in the logarithmic product.  The empty sum is zero,
and $e_{\zero,M}$ is the unit.
Give this algebra contact grading $w(e_{\vct,C})=\sum_i k_iv_i$.

The \emph{dual complex} $\Delta(D)$ has a cell $(I,C)$ of
dimension $|I|-1$ for each nonempty $I$ and each
$C\in\pi_0(D_I)$. Its faces are the pairs $(J,C')$ with
$J\subseteq I$ and $C\subseteq C'$. We include
$(\varnothing,M)$ as the empty face. Thus distinct cells may
have the same vertices. The faces of a cell $(I,C)$ are
indexed by the subsets of $I$; a face poset with this property
is called a \emph{simplicial poset}.

\begin{definition}[Stanley--Reisner algebra]
\label{mir:def-stanley-reisner}
For a field $k$, the \emph{Stanley--Reisner algebra}
$k[\Delta(D)]$ consists of the collections
\[
 (f_{I,C})_{(I,C)},\qquad f_{I,C}\in k[x_i:i\in I],
\]
indexed by all faces, satisfying
\[
 f_{J,C'}=
 f_{I,C}\big|_{x_i=0\ (i\in I\setminus J)}
 \quad\text{if }J\subseteq I\text{ and }C\subseteq C'.
\]
Addition and multiplication are componentwise, and
$\deg x_i=k_i$. 
For $D=\varnothing$, the only face is empty and
$k[\Delta(D)]=k$. This is the compatible-polynomial description
of the face algebra of a finite simplicial poset
\cite[\S5.2, Theorem~1]{DaliDelucchi}.
\end{definition}

\begin{remark}
For an ordinary finite simplicial complex $\Delta$ on a vertex
set $V$, the same definition gives the usual presentation
\[
 k[\Delta]=k[x_v:v\in V]/I_\Delta,\qquad
 I_\Delta=\left(\prod_{v\in J}x_v:
                          J\subseteq V,\ J\notin\Delta\right).
\]
Positive vertex weights give its grading. In particular, if
every nonempty $D_I$ is connected, then
$\Delta(D)=\{I:D_I\ne\varnothing\}$ is an ordinary simplicial
complex and the presentation is that of
Corollary~\ref{cor:lcy-stanley-reisner} and \cite[(3.18)]{GP}.
For disconnected intersections the component labels must be
retained; replacing $\Delta(D)$ by this set of subsets can
change its algebra.
\end{remark}

Compatibility identifies the coefficient of a monomial of
support $I$ on a containing face with its coefficient on the
corresponding face $(I,C)$. Comparing componentwise
multiplication with \eqref{mir:eq:disconnected-face-product}
therefore identifies $A_D\cong\F[\Delta(D)],$
cf. \cite[\S3.2, (3.18)]{GP}, \cite[(1.1)--(1.2)]{PomerleanoIntrinsic}. Together with Theorem~\ref{thm:intro-lcy-zero}, this gives
\begin{equation}\label{mir:eq:loop-algebra-input}
 \gr_F LH^0_\Omega(X)\cong A_D\cong\F[\Delta(D)]
\end{equation}
as contact-weight graded algebras.

Call a pair $(K,E)$, with $E\in\pi_0(D_K)$, a \emph{maximal
incidence pair} if there is no pair $(K',E')$ with
$K\subsetneq K'$ and $E'\subset E$. This definition allows strata
of different dimensions to be maximal.
For each maximal incidence pair define
\begin{equation}\label{mir:eq:face-restriction}
 \rho_{K,E}:A_D\longrightarrow\F[x_i:i\in K],\qquad
 e_{\vct,C}\longmapsto
 \begin{cases}
 \prod_{i\in K}x_i^{v_i},&I(\vct)\subset K\text{ and }E\subset C,\\
 0,&\text{otherwise}.
 \end{cases}
\end{equation}

The following facts about finite simplicial-poset face algebras
are recorded in \cite[\S5.2, especially Lemma~5.6,
Proposition~5.7, Theorem~1, and Remark~5.9]{DaliDelucchi};
we state them in boundary-component notation.

\begin{lemma}\label{mir:lem:face-algebra-facts}
The algebra $A_D$ is finitely generated and reduced. Each map
$\rho_{K,E}$ is a surjective algebra homomorphism, and their product
is injective. Its Krull dimension is
\begin{equation}\label{mir:eq:boundary-dimension}
 d_D=\max\{|I|:D_I\ne\varnothing\},
\end{equation}
where the empty set is allowed, so $d_D=0$ when $D=\varnothing$.
The kernels of the maps $\rho_{K,E}$ are precisely its minimal
prime ideals.  One finite generating set is
\[
 \{e_{\mathbf1_I,C}:\varnothing\ne I\subset\{1,\ldots,r\},
                         \ C\in\pi_0(D_I)\},
 \qquad (\mathbf1_I)_i=\begin{cases}1&i\in I,\\0&i\notin I.\end{cases}
\]
\end{lemma}

\begin{proof}
Under $A_D\cong\F[\Delta(D)]$, the maps $\rho_{K,E}$ are
restriction to maximal faces. The cited face-algebra results
give their surjectivity, the injection into their product,
reducedness, and the stated minimal primes. Hence
\[
 \dim A_D=\max_{(K,E)}\dim\F[x_i:i\in K]=d_D.
\]
The explicit generators follow from
\eqref{mir:eq:disconnected-face-product}: putting
$t_i=e_{\mathbf e_i,D_i}$, where $\mathbf e_i$ is the $i$th
standard basis vector, gives
\[
 e_{\vct,C}=e_{\mathbf1_I,C}
       \prod_{i\in I}t_i^{v_i-1},\qquad I=I(\vct)\ne\varnothing.
\]
There are finitely many pairs $(I,C)$, and the zero-contact
class is the unit.
\end{proof}

\begin{theorem}\label{mir:thm:loop-affine-finitetype}
The algebra $R_X$ is finitely generated and geometrically reduced,
and its spectrum is geometrically connected.
The scheme $\check X$ is therefore connected, reduced, and of finite type, of dimension
$d_D$. In particular, $\dim\check X=\dim X=n$ if and only if
the chosen boundary has a zero-dimensional stratum, equivalently
$D_I\ne\varnothing$ for some $|I|=n$.
A finite generating set of $R_X$ is supplied by the
classes $\vartheta_{\mathbf1_I,C}$.
Moreover the Rees algebra
\begin{equation}\label{mir:eq:loop-Rees}
 \mathcal R_F(R_X)=\bigoplus_{p\ge0}F_pR_X\,s^p\subset R_X[s]
\end{equation}
is finitely generated and flat over $\F[s]$, with fibers
\[
 \mathcal R_F(R_X)/(s)\simeq A_D,
 \qquad
 \mathcal R_F(R_X)/(s-a)\simeq R_X\quad(a\in\F^*).
\]
After extending scalars to $\C$, it gives a flat degeneration
of $\check X$ to $\Spec(A_D\otimes\C)$.
\end{theorem}

\begin{proof}
If $D=\varnothing$, then $R_X=\F$ and the assertions are immediate.
Otherwise let $g_1,\ldots,g_N$ be the generators of $A_D$ in
Lemma~\ref{mir:lem:face-algebra-facts}, and let $b_j$ be their
fundamental-class lifts, of contact weights $w_j>0$.
Subtracting a polynomial in the $b_j$ with the prescribed
highest symbol, and inducting on the filtration degree, gives
\begin{equation}\label{mir:eq:weighted-filtered-generation}
 F_pR_X=\operatorname{span}_{\F}
 \{b_1^{m_1}\cdots b_N^{m_N}:
                   \textstyle\sum_jw_jm_j\le p\}.
\end{equation}
Thus $R_X$ is generated by the $b_j$, and
$\mathcal R_F(R_X)$ is generated by $s$ and $b_js^{w_j}$.
The filtered theta basis gives
\[
 \mathcal R_F(R_X)
 =\bigoplus_{\vct,C}\F[s]\,
       \vartheta_{\vct,C}s^{w(\vct)},
\]
so the Rees algebra is free over $\F[s]$.
Its fiber and localization identities are the standard Rees
identities \cite[Lemma~B.3]{PomerleanoIntrinsic}:
\[
 \mathcal R_F(R_X)/(s)\simeq A_D,\qquad
 \mathcal R_F(R_X)[s^{-1}]\simeq R_X[s,s^{-1}].
\]
The latter identifies every fiber at $a\in\F^*$ with $R_X$.

For every field extension $L/\F$,
$\gr_F(R_X\otimes_\F L)\simeq A_D\otimes_\F L$ is reduced,
by the injective product of restriction maps
\eqref{mir:eq:face-restriction}. The highest nonzero symbol
of a nonzero nilpotent would therefore give a contradiction.
Likewise, an idempotent of positive filtration degree would
have a nonzero square-zero highest symbol. Consequently
$R_X\otimes_\F L$ is reduced and its only idempotents are
$0,1\in L=F_0(R_X\otimes_\F L)$.

For a positively graded $\F$-algebra $A$ with finite-dimensional
graded pieces, write
$H_A(t)=\sum_{p\ge0}\dim_\F(A_p)t^p$.
Multiplication by the regular element $s$ gives
\[
 H_{\mathcal R_F(R_X)}(t)=\frac{H_{A_D}(t)}{1-t}.
\]
The Hilbert-series dimension theorem
\cite[Proposition~4.4.1]{BrunsHerzog} gives
$\dim\mathcal R_F(R_X)=d_D+1$.
Since $s$ is a nonzerodivisor, the open set $s\ne0$
meets every irreducible component; localization therefore
preserves dimension for this finite-type algebra over a field.
The localization identity gives $\dim R_X=d_D$.
Finally, an SNC intersection $D_I\ne\varnothing$ has
dimension $n-|I|$, proving the stated dimension criterion.
All constructions commute with extension of the coefficient
field, in particular from $\F$ to $\C$.
\end{proof}

\begin{remark}
The Rees degeneration depends on the boundary and its positive
contact weights. Its general fiber is the spectrum of the
unfiltered algebra, which is intrinsic to $(X,\Omega)$ by
Proposition~\ref{prop:fusion-invariance-gradings}.
\end{remark}

\subsection{The refined boundary and the central fiber}

Assume henceforth that $D$ has maximal boundary. A reduced scheme
of pure dimension $n$ is \emph{Gorenstein} if it is Cohen--Macaulay
and its dualizing sheaf is invertible. We must prove that this
sheaf is globally trivial on $\check X$, and that $\check X$ is
semi-log-canonical. The analogous filtered argument for symplectic
cohomology is given in \cite[\S5.3 and Appendix B]{PomerleanoIntrinsic};
we will use the loop-algebra identification
\eqref{mir:eq:loop-algebra-input} and
Proposition~\ref{prop:fusion-invariance-gradings}.

A normal pair
$(Z,C)$ with $C$ an effective  $\bQ$-divisor is \emph{log canonical}
if $K_Z+C$ is $\bQ$-Cartier, the coefficients of $C$ are at most
one, and all discrepancies of exceptional divisors on a log resolution $\mu$
are at least $-1$.  Thus, writing
\[
 K_{\widetilde Z}+\widetilde C
 =\mu^*(K_Z+C)+\sum_E a_E E,
\]
one requires $a_E\ge-1$. A reduced pure-dimensional scheme $Z$ is
\emph{semi-log-canonical} (slc) if
\[
 \operatorname{depth}\mathcal O_{Z,z}
 \ge\min\{2,\dim\mathcal O_{Z,z}\}\quad(z\in Z),
\]
it is seminormal, its codimension-one singularities are ordinary
double crossings, $K_Z$ is $\bQ$-Cartier, and the normalization
pair $(Z^\nu,C)$ is log canonical. Here $C$ is the reduced divisor
supported on the conductor, and seminormality means that every finite
birational morphism $Z'\to Z$ from a reduced scheme which is a
homeomorphism and induces isomorphisms on residue fields is an
isomorphism. The displayed depth condition is Serre's condition $S_2$.

\begin{lemma}[A boundary refinement]
\label{mir:lem:cy-boundary-refinement}
There is a composition of blowups along smooth closed boundary strata
\[
 \pi:(M',D')\longrightarrow(M,D)
\]
which is the identity over $X$, such that $D'$ is SNC, all nonempty
intersections of its irreducible components are connected, and its dual
complex is the barycentric subdivision of the original dual complex.
Furthermore, $K_{M'}+D'\sim0$, maximal boundary is preserved, and $D'$
supports an ample divisor with positive coefficients on every component.
\end{lemma}
\begin{proof}
Blow up the connected components of the closed boundary strata
of codimension at least two, in increasing dimension, taking
strict transforms as centers. By \cite[Remark~2.1]{PayneBoundary},
this gives the asserted SNC boundary, connected intersections,
and barycentric subdivision; in particular maximal boundary
is preserved. All centers lie in $D$, so the morphism is the
identity over $X$.

For one such blowup $\pi$ with center of codimension $c$ and
exceptional divisor $E$, let $D_{\widetilde M}$ be the reduced
inverse image of the boundary. Then
\[
 K_{\widetilde M}=\pi^*K_M+(c-1)E,\qquad
 D_{\widetilde M}=\pi^*D-(c-1)E,
\]
so $K_{\widetilde M}+D_{\widetilde M}=\pi^*(K_M+D)$.
If $A$ is an ample $\bQ$-divisor supported on $D$ with all
coefficients positive, then $\pi^*A-\varepsilon E$ is ample
for sufficiently small $\varepsilon\in\bQ_{>0}$. Choose also
$\varepsilon<\operatorname{coeff}_E(\pi^*A)$, the positive sum
of the coefficients of the boundary components containing the
center. Every boundary coefficient then remains positive.
Iteration and clearing denominators give the required ample
divisor on $M'$.
\end{proof}

Choose an ample boundary divisor
$A'=\sum_i k_i'D_i'$ with positive integral coefficients, as supplied
by the refinement lemma. Proposition~\ref{prop:fusion-invariance-gradings}
identifies the unfiltered degree-zero algebra computed from $(M',D')$
with $R_X$. Denote that algebra by $R'_X$ and the identification by
$\phi:R'_X\xrightarrow{\sim}R_X$. Transport its contact filtration by
setting
\[
 F'_pR_X:=\phi(F^{A'}_pR'_X)\qquad(p\ge-1).
\]
Since $\phi$ is a unital algebra isomorphism, this filtration is
exhaustive, has finite-dimensional pieces and $F'_0R_X=\F1$, and
satisfies $F'_pR_X\cdot F'_qR_X\subseteq F'_{p+q}R_X$ for $p,q\ge0$.
Thus the original and refined filtrations are
filtrations of the same unfiltered algebra, with
\begin{equation}\label{mir:eq:original-refined-filtrations}
 \gr_F R_X\simeq A_D,\qquad
 \gr_{F'}R_X\simeq A_{D'}.
\end{equation}
The second identification is Theorem~\ref{thm:intro-lcy-zero}
applied to $(M',D')$. The filtrations and their associated graded
algebras can differ; the affine loop mirror is unchanged. After
extending scalars to $\C$, put
\[
 R=R_X\otimes_{\F}\C,\qquad
 S=\gr_{F'}R\simeq A_{D'}\otimes_{\F}\C
  \simeq\C[\Delta(D')],
 \qquad\deg x_i=k_i'.
\]
We use $F'_pR=(F'_pR_X)\otimes_\F\C\subset R$. Flatness of
$\C$ over $\F$ gives, for each $p\ge0$,
\[
 (F'_pR_X/F'_{p-1}R_X)\otimes_\F\C
       \xrightarrow{\sim}F'_pR/F'_{p-1}R.
\]
These isomorphisms respect multiplication; taking their direct sum
justifies the associated graded identification after scalar extension.
Here $\C[\Delta(D')]$ is the ring of
Definition~\ref{mir:def-stanley-reisner}, with one variable $x_i$
for each component $D_i'$. The last isomorphism follows from
connected boundary intersections and
Corollary~\ref{cor:lcy-stanley-reisner}. We will lift the
properties of this refined central fiber to $R$. 

For a face $\sigma$ of a simplicial complex $\Delta$, write
\[
 \operatorname{link}_\Delta\sigma
 =\{\tau\in\Delta:\tau\cap\sigma=\varnothing,
                           \ \tau\cup\sigma\in\Delta\}.
\]
The empty face is permitted and has dimension $-1$.  Reduced homology
is augmented, with
$\widetilde H_{-1}(\{\varnothing\};\F)=\F$.
The following theorem collects results of Koll\'ar--Xu and
Pomerleano on the dual complex
\cite[Theorem~2(1), Proposition~31, and Paragraphs~32--33]{KollarXu},
\cite[Theorem 5.34 and Corollary 5.35]{PomerleanoIntrinsic},
and the affine case of \cite[Lemma 17.1]{KWY}. We check their
hypotheses and make explicit the grading of the canonical module.

\begin{theorem}[Boundary topology and the central face scheme]
\label{mir:thm:cy-external-inputs}
Let $(M',D')$ be as in Lemma~\ref{mir:lem:cy-boundary-refinement},
and put $\Delta=\Delta(D')$.
\begin{enumerate}[label=(\roman*)]
\item The complex $\Delta$ is pure of dimension $n-1$. For every
face $\sigma$, including the empty face, and every integer $j$,
\[
 \widetilde H_j(\operatorname{link}_{\Delta}\sigma;\F)
 \cong
 \begin{cases}
 \F,&j=n-\dim\sigma-2,\\
 0,&j\ne n-\dim\sigma-2.
 \end{cases}
\]
\item The Stanley--Reisner ring $S=\C[\Delta]$ of
Definition~\ref{mir:def-stanley-reisner} is reduced,
Cohen--Macaulay, and Gorenstein of pure dimension $n$.  Its canonical
module is graded-isomorphic to the ring with shift zero, for every
assignment of positive integer weights to the vertices.
\item The normalization of its spectrum is the disjoint union of
one affine space $\A^n$ for each facet of $\Delta$. On each
component the conductor is the full reduced coordinate boundary.
The spectrum is seminormal, has ordinary double crossings in
codimension one, and is semi-log-canonical.
\end{enumerate}
\end{theorem}

\begin{proof}
The refinement lemma gives a smooth projective pair with
\[
 K_{M'}+D'\sim0,\qquad
 \dim\Delta(D')=n-1,
\]
and connected nonempty boundary intersections. Thus it satisfies
the maximal-degeneration and simplicial hypotheses of
\cite[Theorem 5.34]{PomerleanoIntrinsic}. That theorem, deduced
there from the cited results of Koll\'ar--Xu, says that
$|\Delta|$ has the reduced homology of $S^{n-1}$ and the local
homology of an $(n-1)$-manifold over $\F$. For a point $x$ in
the relative interior of a nonempty face $\sigma$ of dimension $d$,
the local star decomposition gives
\[
 H_j(|\Delta|,|\Delta|\setminus\{x\};\F)
 \cong\widetilde H_{j-d-1}
                  (\operatorname{link}_\Delta\sigma;\F).
\]
The local homology assertion therefore gives (i) for nonempty
faces, including facets with link $\{\varnothing\}$. It also
implies that every facet has dimension $n-1$. The global homology
assertion gives (i) for the empty face.

We now check the hypotheses of \cite[Lemma 17.1]{KWY}, with
its complex $\Gamma=\Delta$ of dimension $n-1$ and coefficient
field $\C$. A generator of the top reduced homology orients the
facets compatibly. The link of a codimension-one face consists
of two points by (i), so that face lies in exactly two facets,
with opposite induced orientations. Consequently condition (1)
of that lemma holds and its boundary subcomplex is empty.
Part (i), after scalar extension to $\C$, gives both its link
vanishing condition (2) and its global vanishing condition (3).
Its affine conclusion therefore applies with $E=\varnothing$:
$\Spec S$ is reduced and Cohen--Macaulay of pure dimension $n$,
has ordinary double crossings in codimension one, has trivial
dualizing sheaf, and is semi-log-canonical. In particular it is
Gorenstein and seminormal. This also gives the Gorenstein
conclusion of \cite[Corollary 5.35]{PomerleanoIntrinsic}.

In the normalization description in the proof of
\cite[Lemma 17.1]{KWY}, the components and their conductor
divisors in the present no-boundary case are
\[
 \nu:\coprod_{\sigma\text{ facet}}
       \Spec\C[x_i:i\in\sigma]\longrightarrow\Spec S,
 \qquad
 C_\sigma=\sum_{i\in\sigma}\{x_i=0\}.
\]
Indeed, a facet has $n$ vertices, and each of its coordinate
hyperplanes corresponds to a codimension-one face shared with
exactly one other facet. Thus $C_\sigma$ is the full reduced
coordinate boundary. This identifies all the objects in (iii).

It remains to record the grading of the trivialization, which
is not part of the statement of that lemma. For each ordered
facet $\sigma=\{i_1,\ldots,i_n\}$, its proof uses the form
\[
 \eta_\sigma=\varepsilon_\sigma
       \frac{dx_{i_1}}{x_{i_1}}\wedge\cdots\wedge
       \frac{dx_{i_n}}{x_{i_n}},
\]
where $\varepsilon_\sigma\in\{1,-1\}$ records the compatible
orientation. Opposite logarithmic residues of the $\eta_\sigma$
along the conductor glue these
forms to a generator of the dualizing sheaf in codimension one,
as in the cited proof. More explicitly, put $Z=\Spec S$, write
$\omega_Z$ for its dualizing sheaf, and choose an open immersion
$j:U\hookrightarrow Z$ whose complement has codimension at least two
and on which this gluing gives $\mathcal O_U\xrightarrow{\sim}\omega_Z|_U$.
Both $\mathcal O_Z$ and $\omega_Z$ satisfy $S_2$, so the natural maps
\[
 \mathcal O_Z\xrightarrow{\sim}j_*\mathcal O_U,\qquad
 \omega_Z\xrightarrow{\sim}j_*(\omega_Z|_U)
\]
are isomorphisms. Applying $j_*$ therefore extends the displayed
trivialization uniquely over $Z$. For arbitrary positive integers $w_i$,
the action $a_\lambda^*x_i=\lambda^{w_i}x_i$ satisfies
\[
 a_\lambda^*\left(\frac{dx_i}{x_i}\right)
     =\frac{d(\lambda^{w_i}x_i)}{\lambda^{w_i}x_i}
     =\frac{dx_i}{x_i}
 \qquad(\lambda\in\C^*).
\]
By uniqueness of the extension, the global generator is therefore
invariant, hence homogeneous
of degree zero. The canonical module with its induced grading
thus satisfies
\begin{equation}\label{mir:eq:cy-face-canonical}
 \omega_S\simeq S
\end{equation}
without a grading shift, completing (ii).
\end{proof}

\subsection{Rees lifting of the Calabi--Yau and slc properties}

We isolate the algebraic step, since local Gorensteinness alone would
not imply that the dualizing line bundle on an arbitrary affine
scheme is globally trivial.  For a graded module $N$, we use
$N(a)_j=N_{a+j}$.

The following proposition collects standard Rees-degeneration
consequences; see \cite[Appendix~B]{PomerleanoIntrinsic},
\cite[Proposition~3.6.11 and Corollary~3.6.14]{BrunsHerzog},
and \cite[Lemma~2.6 and Definition--Theorem~2.3]{KollarFamilies}.

\begin{proposition}[Filtered Calabi--Yau and slc lifting]
\label{mir:prop:cy-rees-lifting}
\label{mir:lem:cy-rees-lifting}
\label{mir:lem:cy-slc-lifting}
Let $R$ be a commutative complex algebra with an exhaustive
multiplicative filtration $F_{-1}R=0$, $F_0R=\C$, and
finite-dimensional $F_pR$. Suppose $S=\gr_F R$ is finitely generated,
reduced, and pure $n$-dimensional, and is graded Gorenstein with
$\omega_S\simeq S$ with shift zero. Put
$\mathcal R=\bigoplus_{p\ge0}F_pR\,s^p$.
\begin{enumerate}[label=(\roman*)]
\item The algebra $R$ is finitely generated, reduced, connected, and
Gorenstein of pure dimension $n$, with $\omega_R\simeq R$.
The Rees algebra $\mathcal R$ is finitely generated and flat over
$\C[s]$, is Gorenstein of pure dimension $n+1$, and satisfies
\begin{equation}\label{mir:eq:relative-canonical}
 \omega_{\mathcal R}\simeq\mathcal R(-1),\qquad
 \omega_{\mathcal R/\C[s]}\simeq\mathcal R,
\end{equation}
where $s$ has weight one and the relative generator has weight zero.
Here the relative canonical module means
$\omega_{\mathcal R}\otimes_{\C[s]}\omega_{\C[s]}^{-1}$.
\item If, in addition, $\Spec S$ is semi-log-canonical, then
$\Spec R$ is semi-log-canonical.
\end{enumerate}
\end{proposition}
\begin{proof}
The initial-form and Rees arguments in the proof of
Theorem~\ref{mir:thm:loop-affine-finitetype} give finite
generation, reducedness and connectedness of $R$, and
\[
 \mathcal R/s\mathcal R\simeq S,\qquad
 \mathcal R/(s-1)\mathcal R\simeq R,\qquad
 \mathcal R[s^{-1}]\simeq R[s,s^{-1}].
\]
Lifts of homogeneous bases of $S$, multiplied by the
corresponding powers of $s$, give a free $\C[s]$-basis of
$\mathcal R$.

Let $\mathfrak m=\mathcal R_{>0}$ be the homogeneous vertex.
The regular-element criterion \cite[Tag~0BJJ]{Stacks}
applied to $s$ makes $\mathcal R_{\mathfrak m}$ Gorenstein
of dimension $n+1$. The grading action
$\alpha_\lambda^*(bs^p)=\lambda^pbs^p$ contracts every
closed-point orbit to $\mathfrak m$. The complement of the
open Gorenstein locus is closed and invariant, hence empty.
All minimal primes are homogeneous and contained in
$\mathfrak m$; local equidimensionality at the vertex
\cite[Theorem~2.1.2(a)]{BrunsHerzog} thus makes $\mathcal R$
pure of dimension $n+1$.

Graded adjunction
\cite[Proposition~3.6.11 and Corollary~3.6.14]{BrunsHerzog}
gives $\omega_{\mathcal R}\simeq\mathcal R(a)$, $a\in\Z$, and
\[
 S\simeq\omega_S
 \simeq(\omega_{\mathcal R}/s\omega_{\mathcal R})(1)
 \simeq S(a+1).
\]
The least nonzero graded degree is zero in $S$, so $a=-1$:
\begin{equation}\label{mir:eq:cy-rees-canonical}
 \omega_{\mathcal R}\simeq\mathcal R(-1).
\end{equation}
Since $\omega_{\C[s]}=\C[s]\,ds$ with $\deg(ds)=1$,
the relative generator has weight zero. Adjunction for $s-1$
gives $\omega_R\simeq R$. Every irreducible component meets
$s=1$ by the grading action, and this regular equation cuts
it in codimension one. This proves \textup{(i)}.

For \textup{(ii)}, $f:\Spec\mathcal R\to\A^1_s$ is flat
and Gorenstein. By \cite[Lemma~2.6 and
Definition--Theorem~2.3]{KollarFamilies}, an open neighborhood
$U$ of $f^{-1}(0)$ has semi-log-canonical fibers; this local
statement does not require properness. The union
$\bigcup_{\lambda\in\C^*}\alpha_\lambda(U)$ has the same
property, since $f\circ\alpha_\lambda=\lambda f$.
Its complement is closed and invariant and misses the vertex,
so contraction makes it empty. In particular,
$f^{-1}(1)=\Spec R$ is semi-log-canonical.
\end{proof}

\subsection{The Calabi--Yau property}

\begin{proof}[Proof of Theorem~\ref{mir:thm-intro-cy}]
Put $R=R_X\otimes_\F\C$. Choose the refinement $(M',D')$ of
Lemma~\ref{mir:lem:cy-boundary-refinement} and a positive integral
ample boundary divisor $A'$. As in
\eqref{mir:eq:original-refined-filtrations}, transport its filtration
$F'=F^{A'}$ to the same unfiltered algebra $R$. Its associated
graded algebra is
\[
 S=\gr_{F'}R\simeq
 \C[x_i:D_i'\subset D']/I_{\Delta(D')},\qquad \deg x_i=k_i'.
\]
Theorem~\ref{mir:thm:cy-external-inputs}(ii) gives a reduced,
pure $n$-dimensional Gorenstein algebra with
$\omega_S\simeq S$ as graded modules, and part (iii) gives
semi-log-canonicity. Proposition~\ref{mir:prop:cy-rees-lifting}(i)
therefore makes $\check X=\Spec R$ connected and Gorenstein of
pure dimension $n$, with
$\omega_{\check X}\simeq\mathcal O_{\check X}$;
part (ii) makes it semi-log-canonical.

On $\check X_{\mathrm{sm}}$ the dualizing sheaf is
$\Omega^n_{\check X_{\mathrm{sm}}/\C}$, so restricting a global
trivialization gives a nowhere-vanishing algebraic $n$-form.
Let $\nu:\check X^\nu\to\check X$ be the normalization and
$C$ the reduced conductor divisor. The normalization is finite,
so $\check X^\nu$ is affine. Semi-log-canonicity says that
$(\check X^\nu,C)$ is log canonical. The normalization formula
for the dualizing sheaf is
\[
 \bigl(\nu^*\omega_{\check X}\bigr)^{**}
 \simeq\omega_{\check X^\nu}(C)
 \simeq\mathcal O_{\check X^\nu},
\]
where the middle term is interpreted as a rank-one reflexive
sheaf. Hence $K_{\check X^\nu}+C$ is Cartier and linearly
equivalent to zero. If $\check X$ is normal, then $C=0$ and
semi-log-canonicity reduces to log canonicity. The irreducible
components of a normal Noetherian scheme are disjoint open-and-closed
subschemes; connectedness thus makes $\check X$ integral.
\end{proof}

The proof produces a trivialization of the dualizing sheaf, without
choosing a distinguished scalar normalization. Normality and
smoothness require additional information about the fusion algebra.
The normality question is the following.

\begin{conjecture}[Normality of the affine loop mirror]
\label{mir:conj:loop-normality}
Under the hypotheses of Theorem~\ref{mir:thm-intro-cy}, the algebra
$R_X\otimes_\F\C$ is normal.
\end{conjecture}

\section{The polyvector algebra of the affine loop mirror}
\label{globpv:comparison}

This section proves Theorem~\ref{thm:intro-global-polyvectors}. Let $(M,D)$ be a positive log Calabi--Yau pair of dimension
$n\ge1$ with maximal boundary, and put $X=M\setminus D$.
Fix a logarithmic volume form $\Omega$ with
$\operatorname{div}_M\Omega=-D$, and set
\[
 L=LH^*_{\Omega}(X;\F)\otimes_{\F}\C,\qquad
 R=L^0=R_X\otimes_{\F}\C,\qquad \check X=\Spec R.
\]
Assume throughout this section that $\check X$ is smooth. The definition of a
positive log Calabi--Yau pair includes connectedness of $M$.
Since $M$ is smooth, it is irreducible; hence $X$ is connected
and $H^0(X;\C)=\C$.

\subsection{The canonical homomorphism}

 Let $\Delta:L^j\to L^{j-1}$ be the
complexification of the loop BV operator induced by rotation of the
formal source. Products, brackets, and $\Delta$ are those of Theorem \ref{thm:summary-fusion} and
Theorem~\ref{thm:loop-bv}.

Polyvectors have degree $k$, and we write
\[
 \operatorname{PV}^k(\check X)
 =\Gamma(\check X,\wedge^kT_{\check X})
 =\bigwedge_R^k\operatorname{Der}_{\C}(R).
\]
Affineness identifies this with
$\bigoplus_{p+q=k}H^q(\check X,\wedge^pT_{\check X})$.
For a nowhere-vanishing algebraic $n$-form $\eta$, let
$\Delta_\eta$ be the divergence BV operator with the contraction
and sign conventions of Example~\ref{ex:volume-divergence-bv}.
Its Schouten bracket satisfies $\{V,f\}=V(f)$, in agreement
with Definition~\ref{def:loop-bv}. Here $\eta$ is arbitrary; the
form used in Theorem~\ref{thm:intro-global-polyvectors} is fixed
in Definition~\ref{globpv:def-mirror-volume} below.

By Theorem~\ref{thm:intro-lcy-zero}, $L^j=0$ for $j<0$. Therefore
$\{R,R\}=0$. For $a\in L^k$ and $f_1,\ldots,f_k\in R$, define
\begin{equation}\label{globpv:eq-sigma}
 \sigma_k(a)(f_1,\ldots,f_k)
 =(-1)^{k(k-1)/2}
 \{\cdots\{\{a,f_1\},f_2\},\ldots,f_k\},
 \qquad \sigma_0=\mathrm{id}_R.
\end{equation}
An alternating multiderivation of $R$ is an alternating multilinear
map which is a derivation in every argument. When $R$ is smooth,
such maps are precisely the sections of $\wedge^kT_{\check X}$.
The iterated-bracket construction \eqref{globpv:eq-sigma} has a
precedent in \cite[(3.14)--(3.15)]{LianZuckerman} in a different setup.

\begin{lemma}\label{globpv:lemma-sigma}
Formula~\eqref{globpv:eq-sigma} defines a homomorphism of graded
Gerstenhaber algebras $\sigma:L\to\operatorname{PV}^*(\check X)$.
It is compatible with localization in $R$.
\end{lemma}
\begin{proof}
For $f\in R$ and homogeneous $a\in L^k$, put
\[
 \iota_f(a)=(-1)^{k-1}\{a,f\}=-\{f,a\}.
\]
The product rule and the graded Jacobi identity, together with
$\{R,R\}=0$, give
\begin{align*}
 \iota_f(ab)&=\iota_f(a)b+(-1)^{|a|}a\iota_f(b),\\
 \iota_f\iota_g&=-\iota_g\iota_f,
 &\iota_{fg}&=f\iota_g+g\iota_f,\\
 \iota_f\{a,b\}&=\{\iota_f(a),b\}
                  +(-1)^{|a|-1}\{a,\iota_f(b)\}.
\end{align*}
In particular, $\iota_f$ is $R$-linear, and
\[
 \iota_{f_k}\cdots\iota_{f_1}a
 =(-1)^{k(k-1)/2}
       \{\cdots\{\{a,f_1\},f_2\},\ldots,f_k\}.
\]
Anticommutation gives alternation, and the identity for
$\iota_{fg}$ gives the derivation rule in each argument. Thus
\eqref{globpv:eq-sigma} defines an alternating multiderivation.
More intrinsically it is an element of
$\operatorname{Hom}_R(\wedge^k_R\Omega^1_{R/\C},R)$, evaluated on
$df_1\wedge\cdots\wedge df_k$. Smoothness makes
$\Omega^1_{R/\C}$ finite projective and identifies this dual with
$\wedge^k_R\operatorname{Der}_{\C}(R)$.

For $a\in L^k$ and $b\in L^\ell$, the first displayed identity
expands the evaluation of $\sigma_{k+\ell}(ab)$ as
\[
 \sum_{\tau\in\operatorname{Sh}(k,\ell)}\operatorname{sgn}(\tau)
 \sigma_k(a)(f_{\tau(1)},\ldots,f_{\tau(k)})\,
 \sigma_\ell(b)(f_{\tau(k+1)},\ldots,f_{\tau(k+\ell)}).
\]
This is exactly the evaluation of $\sigma_k(a)\wedge\sigma_\ell(b)$.
Contraction of polyvectors by $df$ satisfies the last displayed
identity with the Schouten bracket. Also, by definition,
\[
 \iota_{df}\sigma(a)=\sigma(\iota_fa).
\]
To compare the brackets, apply these two identities to both
$\sigma\{a,b\}$ and $[\sigma(a),\sigma(b)]$. Induction on
$|a|+|b|$, followed by evaluation on exact differentials, reduces
the assertion to a bracket with a function; that case is precisely
the displayed contraction identity. Exact differentials generate
$\Omega^1_{R/\C}$, so all these evaluations determine the
polyvector and prove bracket compatibility.

For a multiplicative subset $T\subset R$, the product and bracket
extend to $T^{-1}L$, with
$\{a,f^{-1}\}=-f^{-2}\{a,f\}$ for $f\in T$.
The same formula extends derivations of $R$ to $T^{-1}R$.
Consequently the two constructions of $\sigma$ after localization
have the same evaluations on functions and coincide.
\end{proof}

\subsection{The refined contact filtration and the Rees family}

Apply Lemma~\ref{mir:lem:cy-boundary-refinement} to $(M,D)$.
Corollary~\ref{cor:loop-bv-invariance}(i) identifies the unfiltered
BV algebra before and after this refinement, preserving the
$\Omega$-grading. Thus $L$, $R$, and $\check X$ are unchanged.
Write $(M,D)$ again for the refined positive log Calabi--Yau
pair, whose nonempty boundary intersections are connected, and
choose an ample boundary divisor with positive integral
coefficients for its contact filtration. These auxiliary choices
are used only in the proof: \eqref{globpv:eq-sigma} is defined
by the unfiltered bracket.

The degree-zero associated graded ring is
\begin{equation}\label{globpv:eq-face}
 S:=\gr_F R
 =\C[t_1,\ldots,t_r]/
 (\prod_{i\in I}t_i:D_I=\varnothing).
\end{equation}
By Theorem~\ref{mir:thm:cy-external-inputs}(ii), $S$ is reduced and
Gorenstein of pure dimension $n$, with $\omega_S\simeq S$
of weight zero for the chosen positive contact weights.
Theorem~\ref{mir:thm-intro-cy} gives connectedness and pure
dimension $n$ for $\check X$. Its smoothness makes each local
ring a domain, so distinct irreducible components cannot meet.
There are finitely many components, each therefore open and
closed; connectedness leaves exactly one. Thus $R$ is a domain.

Use the Rees algebra and module
\[
 \mathcal R=\bigoplus_{p\geq0}F_pR\,s^p\subset R[s],
 \qquad
 \operatorname{Rees}_F(L)=\bigoplus_{p\geq0}F_pL\,s^p\subset L[s].
\]
Here cohomological degree is retained:
\[
 \operatorname{Rees}_F(L)^j
 =\bigoplus_{p\geq0}F_pL^j\,s^p.
\]
Theorem~\ref{mir:thm:loop-affine-finitetype} makes $R$ and
$\mathcal R$ finitely generated. Moreover $\mathcal R$ is a
domain, being a subring of $R[s]$.

With a similar argument the next lemma for symplectic
cohomology $SH^*(X)$ is due to 
\cite[Lemma~4.24(2)]{PomerleanoIntrinsic}; see also the proof of
\cite[Theorem~5.30]{GP}. The argument is purely algebraic, there are no new loop calculations needed.

\begin{lemma}\label{globpv:lemma-finite}
The graded loop algebra $L$ is a finite $R$-module, and $\operatorname{Rees}_F(L)$
is finite over $\mathcal R$.
\end{lemma}
\begin{proof}
Write $E_1$ for the direct sum of all contact-filtration and
cohomological degrees on the logarithmic first page. Its
cohomological-degree-zero subalgebra is $S=\gr_FR$.
For each nonempty support $I$, choose a basis of the
finite-dimensional space $H^*(S_I^\circ;\C)$ at contact vector
$\mathbf1_I=\sum_{i\in I}e_i$. If $\mathbf v$ has support $I$,
multiplication by
$\prod_{i\in I}t_i^{v_i-1}$ carries this summand to contact
$\mathbf v$ and is the identity on its coefficient cohomology.
There are finitely many supports. Together with a basis of the
zero-contact group $H^*(X;\C)$, these classes generate $E_1$
as an $S$-module. Thus $E_1$ is finite over the Noetherian ring $S$.

Every element of $S$ is permanent by degree-zero degeneration,
so its action descends to each page and every $d_r$ is $S$-linear.
For precision, form $S$-submodules
$\mathsf B_r\subset\mathsf Z_r\subset E_1$ with
$E_r=\mathsf Z_r/\mathsf B_r$, starting with
$\mathsf Z_1=E_1$ and $\mathsf B_1=0$. Given them on page $r$,
let $\mathsf Z_{r+1}$ and $\mathsf B_{r+1}$ be the inverse images
in $\mathsf Z_r$ of $\ker d_r$ and $\operatorname{im}d_r$,
respectively. The identity $d_r^2=0$ gives
\[
 \mathsf B_r\subset\mathsf B_{r+1}
 \subset\mathsf Z_{r+1}\subset\mathsf Z_r.
\]
Noetherianity makes the increasing sequence $\mathsf B_r$
stationary from some page onward. Since
$\operatorname{im}d_r=\mathsf B_{r+1}/\mathsf B_r$, all
differentials on these later pages vanish. The stable page is a
subquotient of the finite $S$-module $E_1$, hence finite over $S$.
Strong convergence, including compatibility with multiplication,
identifies it with $\gr_F L$.

Choose homogeneous $S$-module generators $\bar a_j$ of $\gr_F L$,
with contact weights $p_j\geq0$, and lifts $a_j\in F_{p_j}L$.
For $a\in F_pL$, express its symbol of weight $p$ as
$\sum_j\bar r_j\bar a_j$, with $\bar r_j$ of weight $p-p_j$;
terms with $p<p_j$ are absent. Lift
$\bar r_j$ to $r_j\in F_{p-p_j}R$. Subtracting
$\sum_jr_ja_j$ lowers the filtration weight. Induction ends at
$F_{-1}L=0$, proving that the $a_j$ generate $L$ over $R$.
The same induction after multiplication by $s^p$ expresses
$a s^p$ in the $\mathcal R$-span of $a_js^{p_j}$: the remainder
of weight at most $p-1$ is $s$ times its Rees representative.
These finitely many elements therefore generate $\operatorname{Rees}_F(L)$.
\end{proof}

Proposition~\ref{mir:prop:cy-rees-lifting}(i) supplies flatness over
$\C[s]$, graded Gorensteinness, and the canonical-module
formulas \eqref{mir:eq:relative-canonical}. Fix a homogeneous
weight-zero generator $\eta_{\mathrm{rel}}$ of
$\omega_{\mathcal R/\C[s]}$. The following extension property
is the additional consequence of smoothness needed here, and it is standard argument in commutative algebra.

\begin{lemma}\label{globpv:lemma-rees}
The contact Rees algebra $\mathcal R$ is normal. The relative
smooth locus $U$ of $Y=\Spec\mathcal R\to\A^1_s$ has
complement of codimension at least two.
\end{lemma}
\begin{proof}
The Gorenstein algebra $\mathcal R$ satisfies Serre's
condition $S_2$ and, as just proved, is a Noetherian domain.
We check every height-one point of $Y$.
If its prime does not contain $s$, its local ring is regular
because
\[
 \mathcal R[s^{-1}]=R[s,s^{-1}]
\]
is smooth over $\C$. If its prime $\mathfrak p$ contains $s$,
then $\mathfrak p$ is minimal over $(s)$: the principal ideal
theorem and the fact that $s$ is nonzero in a domain force every
minimal prime over $(s)$ to have height one. Its image
$\mathfrak q$ in $S$ is consequently a minimal prime of $S$.
Reducedness gives a field $S_{\mathfrak q}$, so
$\mathcal R_{\mathfrak p}/s\mathcal R_{\mathfrak p}$ is a field.
The maximal ideal of this one-dimensional local domain is thus
generated by $s$, and the ring is a DVR. This proves $R_1$.
Serre's $R_1+S_2$ criterion now gives normality of $\mathcal R$.

In fact every height-one point lies in the relative smooth locus
$U$. This is immediate off $s=0$, where the family is the product
of the smooth variety $\check X$ and $\Gm$. At a height-one point
on $s=0$, the fiber is at a generic point of one of the
polynomial-ring components of the reduced face ring. It is
smooth over $\C$ there. Flatness and the fiber criterion for
smoothness give relative smoothness at this point as well.
Thus the closed complement $Y\setminus U$ contains no point
of codimension zero or one, proving the final assertion.
\end{proof}

\subsection{A maximal face and the generic ranks}
\label{globpv:generic}

The goal now is to compute the generic $R$-rank of each $L^j$.
Lemmas~\ref{globpv:lemma-local-page}
and~\ref{globpv:lemma-local-degeneration} identify the first page
after localization at a maximal face and prove that all
differentials of the localized spectral sequence vanish.
Lemma~\ref{globpv:lemma-generic-ranks} then
compares the special and generic ranks over a DVR in the Rees
family. This uses the logarithmic BV first page already established
in Theorem~\ref{thm:summary} and standard facts of algebra;
no new loop calculation is needed.

Fix a zero-dimensional stratum and number its incident boundary
components $D_1,\ldots,D_n$. Set $h=t_1\cdots t_n\in S$.
Let $\theta_1,\ldots,\theta_n$ be the degree-one cohomology
classes of its normal torus $(S^1)^n$ dual to the ordered,
positively oriented normal circles. Write
\[
 \Lambda(\theta_1,\ldots,\theta_n)
 :=\bigwedge_{\C}^{\bullet}
       (\C\theta_1\oplus\cdots\oplus\C\theta_n)
\]
for their exterior algebra. Thus $\theta_i^2=0$ and
$\theta_i\theta_j=-\theta_j\theta_i$; each $\theta_i$ has
cohomological degree one. In the localized first page below,
these generators have contact weight zero.

\begin{lemma}[The localized first page]
\label{globpv:lemma-local-page}
Localization at $h$ gives the graded-algebra identifications
and bracket identities
\begin{equation}\label{globpv:eq-local-page}
\begin{gathered}
 S[h^{-1}]=\C[t_1^{\pm1},\ldots,t_n^{\pm1}],\qquad
 E_1[h^{-1}]=S[h^{-1}]\otimes_{\C}
                  \Lambda(\theta_1,\ldots,\theta_n),\\
 \{\theta_i,t_j\}=\delta_{ij}t_j,\qquad
 \{\theta_i,\theta_j\}=0.
\end{gathered}
\end{equation}
Here $\delta_{ij}$ is $1$ if $i=j$ and $0$ otherwise.
The displayed brackets, the vanishing bracket on $S[h^{-1}]$,
and the Leibniz rule determine the entire localized bracket.
\end{lemma}
\begin{proof}
A variable outside this maximal face annihilates $h$,
so localization removes all other face components and makes
$t_1,\ldots,t_n$ invertible. On coefficient cohomology,
multiplication by $h$ carries a class whose support is contained
in this face to its restriction on the normal torus of the chosen
stratum. If this restriction is zero, multiplication by $h$
annihilates the class, so its localization is zero.
A class with incompatible support is annihilated.
Every cohomology class of this torus already occurs at contact
$e_1+\cdots+e_n$. Dividing that summand by $h$ gives all of
$\Lambda(\theta_1,\ldots,\theta_n)$ at weight zero. Thus these
restrictions identify the whole localized coefficient algebra,
including its cup product, with the displayed exterior algebra.
The source-rotation formula \eqref{eq:bv-first-page} and the
bracket convention give $\{\theta_i,t_j\}=\delta_{ij}t_j$.

Since the bracket has degree $-1$, its restriction to
$S[h^{-1}]$ vanishes. For fixed $i,j$ there are unique
$b_k\in S[h^{-1}]$ such that
$\{\theta_i,\theta_j\}=\sum_{k=1}^n b_k\theta_k$.
For every $\ell\in\{1,\ldots,n\}$, the Leibniz rule and the
graded Jacobi identity give
\[
 b_\ell t_\ell
 =\{\{\theta_i,\theta_j\},t_\ell\}
 =\{\theta_i,\{\theta_j,t_\ell\}\}
   -\{\theta_j,\{\theta_i,t_\ell\}\}
 =0.
\]
The two terms are both $\delta_{i\ell}\delta_{j\ell}t_\ell$.
Invertibility of $t_\ell$ gives $b_\ell=0$. This proves
$\{\theta_i,\theta_j\}=0$.
\end{proof}

\begin{lemma}[Degeneration after localization]
\label{globpv:lemma-local-degeneration}
Every differential $d_r$ on $E_r[h^{-1}]$ vanishes. Consequently,
the page maps and strong convergence identify
\begin{equation}\label{globpv:eq-local-stable-page}
 (\gr_F L)[h^{-1}]=E_\infty[h^{-1}]
 =S[h^{-1}]\otimes_{\C}\Lambda(\theta_1,\ldots,\theta_n),
\end{equation}
with the product and bracket of \eqref{globpv:eq-local-page}.
\end{lemma}
\begin{proof}
The element $h$ is permanent by degree-zero degeneration.
Localization is exact, so it commutes with taking the cohomology
of each page differential and hence with passage to the next page.
The product and bracket extend to the localized pages.

We prove the vanishing by induction on $r$. If all earlier
localized differentials vanish,
Lemma~\ref{globpv:lemma-local-page} identifies the current page with $S[h^{-1}]\otimes_{\C}\Lambda(\theta_1,\ldots,\theta_n)$.
The elements $t_j^{\pm1}$ are permanent. Bracket compatibility
therefore gives
\[
 \{d_r\theta_i,t_j\}
 =d_r\{\theta_i,t_j\}=d_r(\delta_{ij}t_j)=0.
\]
The element $d_r\theta_i$ has cohomological degree two.
On this exterior-degree-two component, bracketing with $t_j$
is $-t_j$ times contraction with the basis vector dual to
$\theta_j$. Since $t_j$ is invertible, every such contraction
of $d_r\theta_i$ vanishes, and hence $d_r\theta_i=0$.
The product derivation rule then gives $d_r=0$ on the whole
localized page. This proves the inductive assertion, including
$n=1$, when the exterior-degree-two component is zero.

By the proof of Lemma~\ref{globpv:lemma-finite}, the unlocalized
spectral sequence is stationary from a finite page onward.
Localizing that page therefore gives $E_\infty[h^{-1}]$.
Strong convergence identifies $E_\infty$ with $\gr_F L$,
proving \eqref{globpv:eq-local-stable-page}.
These identifications come from the page maps, so a surviving
class represented by a filtered cocycle is identified with the
corresponding symbol of its cohomology class in $L$.
\end{proof}

\begin{lemma}[Generic ranks]\label{globpv:lemma-generic-ranks}
Put $K=\operatorname{Frac}(R)$. Then
\[
 \dim_K(L^j\otimes_RK)=\binom nj\quad(0\leq j\leq n),
 \qquad L^j\otimes_RK=0\quad(j\notin[0,n]).
\]
\end{lemma}
\begin{proof}
Let $\mathfrak q$ be the minimal prime of $S$ corresponding to
the chosen maximal face, and let $\mathfrak p$ be its inverse
image under $\mathcal R\twoheadrightarrow S$. Set
\[
 T=\mathcal R_{\mathfrak p},\qquad
 \kappa=S_{\mathfrak q}=\C(t_1,\ldots,t_n),\qquad
 Q=\operatorname{Frac}(\mathcal R)=K(s).
\]
The equality for $\kappa$ follows from \eqref{globpv:eq-local-page}.
The prime $\mathfrak p$ has height one, and the proof of
Lemma~\ref{globpv:lemma-rees} shows that $T$ is a DVR with
uniformizer $s$ and residue field $\kappa$.

The modules $M^j=\bigl(\operatorname{Rees}_F(L)^j\bigr)_{\mathfrak p}$ are finite over
$T$ by Lemma~\ref{globpv:lemma-finite}. Multiplication by $s$
is injective on $\operatorname{Rees}_F(L)\subset L[s]$, and remains
injective after localization. Since every nonzero element of
$T$ is a unit times a power of $s$, $M^j$ is torsion-free,
and hence free. Its reduction is
\[
 M^j/sM^j=(\gr_F L^j)_{\mathfrak q}
 =\wedge_\kappa^j\langle\theta_1,\ldots,\theta_n\rangle
\]
by \eqref{globpv:eq-local-stable-page}. Thus its rank is
$\binom nj$ for $0\leq j\leq n$, and zero otherwise.
Finally, $\operatorname{Rees}_F(L)[s^{-1}]=L[s,s^{-1}]$ gives
\[
 M^j\otimes_TQ
 \simeq(L^j\otimes_RK)\otimes_KQ.
\]
Dimensions are unchanged by the field extension $Q/K$, proving
the assertions.
\end{proof}

\subsection{The arc volume class and its Rees extension}
\label{globpv:ree-volume}

The goal now is to obtain a nowhere-vanishing top polyvector on
$\check X$ from the logarithmic volume form $\Omega$, thereby
determining the mirror volume form. We use the established BV
spectral sequence and Rees-family properties, together with
standard extension and logarithmic residue arguments; no new
loop calculation is needed.

The zero-contact inclusion is a map of complexes
\[
 \jmath_0:F_0\mathbf B_\Omega\longrightarrow\mathbf B_\Omega,
 \qquad F_0\mathbf B_\Omega\simeq C^*(X;\C),
\]
after extension of scalars. The algebraic form $\Omega|_X$ is
closed because it has top degree. Algebraic de Rham comparison
therefore gives its class $[\Omega]\in H^n(X;\C)$, with the
usual integration normalization. Use the specified arc-stage
comparison $F_0\mathbf B_\Omega\simeq B_{\Omega,0}
\simeq C^*(X;\C)$ to regard this class as an element of
$H^n(F_0\mathbf B_\Omega)$. Define
\begin{equation}\label{globpv:eq-volume-class}
 u=H^n(\jmath_0)([\Omega])\in F_0L^n.
\end{equation}
This definition uses the actual subcomplex inclusion, but does
not assert that its induced map on cohomology is injective.
By \eqref{eq:bv-first-page}, $\Delta_{\log}$ vanishes on the
zero-contact summand $H^*(X;\C)$ of $E_1$. Its induced operator
on the zero-contact subquotient of $E_\infty$ therefore vanishes.
Strong BV convergence in Theorem~\ref{thm:summary} identifies
this operator with $\operatorname{gr}_0\Delta$ on
\[
 \operatorname{gr}_0L=F_0L/F_{-1}L=F_0L.
\]
Thus $\Delta|_{F_0L}=0$, and in particular
\begin{equation}\label{globpv:eq-volume-bv}
 \Delta u=0.
\end{equation}

On the relative smooth locus $U$ from Lemma \ref{globpv:lemma-rees}, the restriction
$\eta_{\mathrm{rel}}|_U$ generates
$\wedge^n\Omega^1_{U/\A^1}$. We denote by
$\eta_{\mathrm{rel}}^{-1}$ its dual generator:
\[
 \eta_{\mathrm{rel}}^{-1}\in
 \Gamma\!\left(U,
   \operatorname{Hom}_{\mathcal O_U}
   \bigl(\wedge^n\Omega^1_{U/\A^1},\mathcal O_U\bigr)\right),
 \qquad
 \eta_{\mathrm{rel}}^{-1}(\eta_{\mathrm{rel}}|_U)=1.
\]
Thus $\eta_{\mathrm{rel}}^{-1}$ is a nowhere-vanishing
relative top polyvector on $U$.

\begin{lemma}\label{globpv:lemma-volume}
The top polyvector $\sigma_n(u)$ is nowhere vanishing on $\check X$.
More precisely, for the chosen relative generator
$\eta_{\mathrm{rel}}$, the Rees extension of $\sigma_n(u)$ is
$c\eta_{\mathrm{rel}}^{-1}$ on $U$, for a constant $c\in\C^*$.
\end{lemma}
\begin{proof}
Put $Y=\operatorname{Spec}\mathcal R$. Extend the filtered loop
bracket to the full graded algebra $\operatorname{Rees}_F(L)$ by
\[
 \{a s^p,b s^q\}=\{a,b\}s^{p+q};
\]
it is $\C[s]$-linear, and $s$ is central. Since $u\in F_0L^n$,
its top iterated bracket defines the homogeneous relative
multiderivation
\[
 P\in\operatorname{Hom}_{\mathcal R}
 \bigl(\wedge^n_{\mathcal R}\Omega^1_{\mathcal R/\C[s]},
       \mathcal R\bigr)_0.
\]
Concretely, for $f_i\in F_{p_i}R$ it satisfies
\[
 P\bigl(d(f_1s^{p_1}),\ldots,d(f_ns^{p_n})\bigr)
 =s^{p_1+\cdots+p_n}\sigma_n(u)(f_1,\ldots,f_n),
\]
where differentials are relative to $\C[s]$. This also verifies
its weight zero. Its restriction to $s=1$ is $\sigma_n(u)$.

Using the dual generator $\eta_{\mathrm{rel}}^{-1}$ defined
above, write
\[
 P|_U=g\eta_{\mathrm{rel}}^{-1},\qquad
 g=P|_U(\eta_{\mathrm{rel}}|_U)\in\Gamma(U,\mathcal O_U).
\]
Normality of $Y$ and
$\operatorname{codim}(Y\setminus U)\geq2$ give
\[
 \Gamma(U,\mathcal O_U)=\mathcal R.
\]
Thus $g$ extends uniquely over $Y$. The open set $U$ is
invariant under the grading action, and both $P$ and
$\eta_{\mathrm{rel}}^{-1}$ have weight zero. Hence $g$ is
invariant; uniqueness makes its extension invariant too.
Consequently
\[
 g\in\mathcal R_0=F_0R=\C.
\]
This argument extends a scalar regular function. It does not
require $\operatorname{Rees}_F(L)$ or any of its graded pieces to be reflexive.

We prove that this scalar is nonzero by tracing the actual
zero-contact class through the localized spectral sequence.
At the chosen zero-dimensional boundary stratum take analytic
coordinates with $D=\{z_1\cdots z_n=0\}$. The divisor equation
$\operatorname{div}\Omega=-D$ gives
\[
 \Omega=b(z)\frac{dz_1}{z_1}\wedge\cdots\wedge\frac{dz_n}{z_n},
 \qquad b(0)\ne0.
\]
Orient a small normal torus by the ordered positive angular
coordinates. If $\theta_i$ integrates to $1$ on its $i$th
circle, then
\[
 [\Omega]|_{(S^1)^n}
 =\lambda\,\theta_1\cdots\theta_n,
 \qquad
 \lambda=\int_{(S^1)^n}\Omega=(2\pi i)^n b(0)\ne0.
\]
The integral equality is the iterated Cauchy integral formula.
The logarithmic first-page product of a zero-contact class with
a class supported at this stratum is restriction to this normal
torus followed by cup product, by
Lemma~\ref{lem:fusion-leading-product}. Hence the image of the
zero-contact class $[\Omega]$ in $E_1[h^{-1}]$ is precisely
$\lambda\theta_1\cdots\theta_n$.

The class $[\Omega]$ has a closed representative in $F_0\mathbf B_\Omega$,
so it is carried by the page maps to the symbol of its actual
image $u$ whenever that symbol survives.
Lemma~\ref{globpv:lemma-local-degeneration} shows that every
localized page differential is zero. Its displayed nonzero localized
class therefore survives to $E_\infty[h^{-1}]$. In particular
$u\ne0$, and its symbol in $(\gr_F L)[h^{-1}]$ equals
$\lambda\theta_1\cdots\theta_n$. This rules out its becoming
a boundary from a higher contact level; such a boundary would
also make this localized class zero on a later page.

Now reduce $P$ modulo $s$ at the generic point of this central
component, which belongs to $U$. On $U$, the module
$\Omega^1_{U/\A^1}$ is locally free of rank $n$, so its top
exterior power and its dual commute with restriction to $s=0$.
Thus reduction of $P|_U$ is the top multiderivation on the
smooth locus of the central fiber obtained by reducing its
displayed evaluation formula. Compatibility of the bracket
with the associated graded identifies the reduction with the
top iterated bracket of
$\lambda\theta_1\cdots\theta_n$. Formula
\eqref{globpv:eq-local-page}, with the sign in
\eqref{globpv:eq-sigma}, evaluates it on
$t_1,\ldots,t_n$ as
\[
 \lambda t_1\cdots t_n\ne0.
\]
Thus $P|_U$ is nonzero, and its constant ratio to
$\eta_{\mathrm{rel}}^{-1}$ is $c\in\C^*$. Since the smooth fiber $s=1$ is contained in
$U$, restriction to that fiber proves that $\sigma_n(u)$ is
nowhere vanishing on $\check X$.
\end{proof}

\subsection{The global isomorphism and BV compatibility}

The following differential-simplicity statement is classical;
see \cite[Theorem~1.3]{BavulaDifferentialOperators}.
We recall a proof for use in subsequent arguments.

\begin{lemma}[Differential simplicity]\label{globpv:lemma-differential-ideals}
Let $B$ be a smooth finitely generated connected $\C$-algebra.
An ideal $I\subset B$ stable under every $\C$-derivation of $B$
is either $0$ or $B$.
\end{lemma}
\begin{proof}
For $B=0$ the assertion is immediate. Otherwise smoothness
makes all local rings domains, so the finitely many
irreducible components of $\Spec B$ are disjoint and open.
Connectedness therefore makes $B$ a domain. Suppose
$0\ne I\ne B$, choose a maximal ideal $\mathfrak m\supset I$,
and put $T_0=B_{\mathfrak m}$ and $J=IT_0$.
Then $0\ne J\ne T_0$. Smoothness makes
$\Omega^1_{B/\C}$ finite projective, so localization gives
\[
 \operatorname{Der}_{\C}(T_0)
 =T_0\otimes_B\operatorname{Der}_{\C}(B).
\]
The quotient rule shows that every derivation of $B$ preserving
$I$ extends to a derivation preserving $J$. The displayed
equality therefore makes $J$ stable under all derivations of
$T_0$.

The regular local ring $T_0$ has residue field $\C$.
If its dimension is zero it is a field, already contradicting
$0\ne J\ne T_0$. Otherwise choose a regular system of
parameters $x_1,\ldots,x_d$. Their differentials form a basis
of $\Omega^1_{T_0/\C}$: their classes modulo the maximal ideal are a basis of the
cotangent space, and this module is free of rank $d$.
The dual derivations $\partial_i$ therefore satisfy
$\partial_i(x_j)=\delta_{ij}$.

By the Krull intersection theorem each nonzero element of $T_0$
has finite maximal-ideal order. Choose $a\in J\setminus\{0\}$
of least such order $e$. Properness gives $e\geq1$, and the
initial form of $a$ is a nonzero homogeneous polynomial of
degree $e$ in
\[
 \operatorname{gr}_{\mathfrak mT_0}T_0
 \simeq\C[\bar x_1,\ldots,\bar x_d].
\]
Each $\partial_i$ lowers the maximal-ideal filtration by at
most one and acts on this associated graded ring as ordinary
partial differentiation. In characteristic zero some partial
derivative of that positive-degree polynomial is nonzero.
The corresponding $\partial_i(a)\in J$ has order exactly
$e-1$, contradicting the choice of $a$. Thus $I$ cannot be a
nonzero proper ideal.
\end{proof}

\begin{definition}[The mirror volume form]
\label{globpv:def-mirror-volume}
Let $u$ be the arc volume class of \eqref{globpv:eq-volume-class}.
By Lemma~\ref{globpv:lemma-volume}, $\sigma_n(u)$ is a
nowhere-vanishing section of $\wedge^nT_{\check X}$. Define
$\eta\in\Gamma(\check X,\Omega^n_{\check X/\C})$ to be its
unique dual section:
\[
 \langle\eta,\sigma_n(u)\rangle=1.
\]
Thus $\eta$ is nowhere vanishing, and
$\eta^{-1}:=\sigma_n(u)$ denotes the dual top polyvector.
Equivalently, if
$\eta_1:=\eta_{\mathrm{rel}}|_{s=1}$ and
$\sigma_n(u)=c\eta_1^{-1}$ as in
Lemma~\ref{globpv:lemma-volume}, then
\[
 \eta=c^{-1}\eta_1,\qquad c\in\C^*.
\]
The defining pairing makes $\eta$ independent of the choice of
$\eta_{\mathrm{rel}}$.
\end{definition}

\begin{proof}[Proof of Theorem~\ref{thm:intro-global-polyvectors}]
Use the volume form $\eta$ of
Definition~\ref{globpv:def-mirror-volume}, so
$\sigma_n(u)=\eta^{-1}$.
We first prove surjectivity. Contraction with $\eta^{-1}$ gives
an isomorphism of finite projective $R$-modules
\[
 \wedge_R^{n-1}\Omega^1_{R/\C}
      \xrightarrow{\ \sim\ }\operatorname{Der}_{\C}(R),\qquad
 df_1\wedge\cdots\wedge df_{n-1}
      \longmapsto
 \iota_{df_{n-1}}\cdots\iota_{df_1}\eta^{-1}.
\]
The isomorphism is checked in a local basis in which
$\eta^{-1}$ is a unit times the wedge of the dual basis.
Exact differentials generate $\Omega^1_{R/\C}$; their
$(n-1)$-fold wedges therefore generate its exterior power.
For every such wedge the displayed image is
\[
 \sigma_1\bigl(\iota_{f_{n-1}}\cdots\iota_{f_1}u\bigr),
\]
by the contraction identity proved in
Lemma~\ref{globpv:lemma-sigma}. For $n=1$ this formula has no
contractions and uses $u\in L^1$. Taking $R$-linear combinations
proves that $\sigma_1$ is surjective in all dimensions under
consideration. The target exterior algebra is generated by
functions and vector fields, so multiplicativity of $\sigma$
proves its surjectivity in every degree.

For injectivity put $N=\ker\sigma$. It is a finite $R$-module
by Lemma~\ref{globpv:lemma-finite} and Noetherianity of $R$.
Surjectivity and Lemma~\ref{globpv:lemma-generic-ranks} imply
that $\sigma\otimes_RK$ is an isomorphism in every degree:
source and target have the same finite dimension $\binom nj$
for $0\leq j\leq n$, and both vanish in every other degree.
Consequently $N\otimes_RK=0$.
Each member of a finite generating set of $N$ is therefore
killed by some nonzero element of the domain $R$; the product
of these elements is a nonzero element of
$I=\operatorname{Ann}_R(N)$.
For $D\in\operatorname{Der}_{\C}(R)$ choose $a\in L^1$ with
$\sigma_1(a)=D$. Since $\sigma$ preserves brackets,
$\{a,N\}\subset N$. If $f\in I$ and $m\in N$, the product
rule, with $|a|-1=0$, gives
\[
 0=\{a,fm\}=D(f)m+f\{a,m\}=D(f)m.
\]
It follows that $D(f)\in I$. Thus $I$ is stable under every
$\C$-derivation of $R$.

The nonzero ideal $I$ is therefore $R$ by
Lemma~\ref{globpv:lemma-differential-ideals}. Hence $N=0$,
proving the isomorphism
\begin{equation}\label{globpv:eq-main}
 \sigma:L\xrightarrow{\ \sim\ }
 \operatorname{PV}^*(\check X)
 =\bigwedge_R^*\operatorname{Der}_{\C}(R).
\end{equation}
Since $\operatorname{Der}_{\C}(R)$ is projective of rank $n$,
its exterior powers give all the asserted generation, rank,
projectivity, and degree-vanishing conclusions for $L$.

Finally transport the source-rotation operator to polyvectors:
\[
 \widetilde\Delta=\sigma\Delta\sigma^{-1},\qquad
 \delta=\widetilde\Delta-\Delta_\eta.
\]
Both terms generate the same Gerstenhaber bracket. Subtracting
the bracket-generation identities gives
\[
 \delta(P\wedge Q)=\delta(P)\wedge Q
                       +(-1)^{|P|}P\wedge\delta(Q).
\]
Thus $\delta$ is a graded derivation of degree $-1$.
Theorem~\ref{thm:intro-lcy-zero} gives
$\Delta(R)\subseteq L^{-1}=0$, and divergence also vanishes on
functions, so $\delta|_R=0$.
Its restriction to vector fields is therefore an $R$-linear functional
$\operatorname{Der}_{\C}(R)\to R$. Since the tangent module
is finite projective, this functional is evaluation against a
unique $\alpha\in\Omega^1_{R/\C}$. The exterior algebra is
generated in degrees zero and one, so $\delta=\iota_\alpha$
in every degree. Equation~\eqref{globpv:eq-volume-bv} gives
$\widetilde\Delta(\eta^{-1})=0$. Divergence also kills
$\eta^{-1}$: in Example~\ref{ex:volume-divergence-bv}, contraction
with $\eta$ sends this top polyvector
to the constant function $1$, whose exterior derivative is
zero. Hence $\iota_\alpha(\eta^{-1})=0$. Contraction with a
nonvanishing top polyvector identifies one-forms with
$(n-1)$-vectors, so $\alpha=0$. Thus
\[
 \sigma\Delta\sigma^{-1}=\Delta_\eta,
\]
which proves the BV comparison.
\end{proof}

\begin{remark}
The final identification of the BV operator with divergence is
classical; see \cite[Section~2]{KoszulBV} and
\cite[Section~2, Corollary following Theorem~3]{HuebschmannBV}.
A related strategy using brackets with functions and a
distinguished top class to recover polyvectors appears in
\cite{GromanIndicatorV1}.
\end{remark}

\begin{remark}[Rational coefficients]\label{globpv:remark-rational}
The iterated-bracket formula \eqref{globpv:eq-sigma} is defined
over $\F$ and gives
\[
 \sigma_{\F}:LH^*_{\Omega}(X)
 \longrightarrow
 \bigwedge_{R_X}^{\bullet}\operatorname{Der}_{\F}(R_X).
\]
Its complexification is $\sigma$, so faithful flatness of
$\C/\F$ makes $\sigma_{\F}$ an isomorphism. Here $R_X$ is finitely
presented, and smoothness descends from its complexification
along the faithfully flat extension $\C/\F$. Thus its module
of differentials is finite projective. Base change of that
module, its dual, and exterior powers identifies the
complexification of the rational target with
$\operatorname{PV}^*(\check X)$. This is the base-change
identification under which $\sigma$ becomes the isomorphism
proved above. The logarithmic de Rham class
was used over $\C$ to prove bijectivity; it need not itself be a
rational cohomology class.
\end{remark}

\section{The loop trace}
\label{mir:sec-trace-potentials}

Let $(M,D)$ be a positive log Calabi--Yau pair, put
$X=M\setminus D$, and fix
$\Omega\in\Gamma(X,\Omega_X^n)$ with
$\operatorname{div}_M\Omega=-D$. Choose an ample boundary divisor
$A=\sum_i k_iD_i$, $k_i\in\Z_{>0}$, to define the pole filtration.
Write
\[
 R_X=LH^0_\Omega(X),\qquad
 R_X=\bigoplus_{p}\F\vartheta_p,
 \qquad p=(\vct,C),\quad\vartheta_0=1,
\]
for the basis of Theorem~\ref{thm:intro-theta-basis}, constructed in
Section~\ref{sec:log-cy-zero}. Connectedness of $M$ gives one arc
label $0$. In this section we define the arc-coefficient functional $\epsilon:R_X\to \F$ and
prove that the pairing $(a,b)\mapsto\epsilon(ab)$ is
nondegenerate under the maximal-boundary and smooth-mirror
hypotheses stated below. This will reduce comparison of the degree-zero loop cohomology
algebra with the Gross--Siebert algebra to comparison of their
two-input and three-input traces. The basis and $\epsilon$ are
independent of $A$; the ample divisor is used only for finite pole
bounds and, in subsequent sections, bounds on curve degrees.

\subsection{The arc coefficient and its moments}

\begin{definition}[Loop trace]\label{mir:def-arc-trace}
For a finite expansion $a=\sum_p a_p\vartheta_p\in R_X$, define
the \emph{loop trace} by
\[
 \epsilon(a)=a_0.
\]
It extracts the coefficient of $\vartheta_0$, whose label
corresponds to the arc component.
For $d\geq1$ and $a_1,\ldots,a_d\in R_X$, define the
\emph{$d$-input loop trace} by
\begin{equation}\label{mir:eq-trace-moments}
 T_d(a_1,\ldots,a_d):=\epsilon(a_1\cdots a_d).
\end{equation}
When the arguments are contact labels, write
$T_d(p_1,\ldots,p_d)=T_d(\vartheta_{p_1},\ldots,\vartheta_{p_d})$.
The empty product has trace $T_0=1$. In particular, the
\emph{two-input loop trace} and \emph{three-input loop trace} are
\[
 T_2:R_X^{\otimes2}\longrightarrow\F,\quad
       T_2(a,b)=\epsilon(ab),\qquad
 T_3:R_X^{\otimes3}\longrightarrow\F,\quad
       T_3(a,b,c)=\epsilon(abc).
\]
Thus $T_1=\epsilon$, and $T_d$ for $d>1$ is a multilinear
functional obtained from the loop trace and the fusion product.
\end{definition}

\begin{remark}[Finite-factor computation of the loop trace]
\label{rem:arc-trace-finite-factor}
Fix $P\ge0$ and identify
$H^0(B_{\Omega,P})$ with its image $F_PR_X\subset R_X$ by
Theorem~\ref{thm:intro-theta-basis} and its proof in
Section~\ref{sec:log-cy-zero}.
Choose a nice chart
\[
 h:U=Y\times\A^J\longrightarrow F_P\mathscr LX
\]
as in Proposition~\ref{prop:lcy-equidimensional-factors},
on a constant volume-index level $\delta_Y=\delta$ meeting
the arc stage. Choose a closed subscheme $Y_0\subset Y$ presenting
the pulled-back arc inclusion, so that
\[
 U\times_{F_P\mathscr LX}F_0\mathscr LX=Y_0\times\A^J.
\]
By parts~(i)--(ii) of that proposition, $Y_{\mathrm{red}}$
is pure of dimension $\delta$. Its irreducible components
$Z_{q,\beta}$ are the reduced closures of the local contact
components with labels $q$ of weight at most $P$; $\beta$
distinguishes components with the same global label.

For $a\in F_PR_X$, denote by $\rho_{P,Y}(a)$ its
finite-factor representative in
Definition~\ref{def:finite-factor-representatives}, using
its unique lift to $H^0(B_{\Omega,P})$.
The theta construction in
\eqref{eq:lcy-theta-stage-class}--\eqref{eq:lcy-theta-colimit-class}
gives
\[
 \begin{aligned}
 \rho_{P,Y}:F_PR_X&\longrightarrow H^{BM}_{2\delta}(Y;\F),\\
 \rho_{P,Y}\!\left(\sum_q a_q\vartheta_q\right)
   &=\sum_q a_q\sum_\beta[Z_{q,\beta}].
 \end{aligned}
\]
Here each component class is pushed forward along
$Z_{q,\beta}\hookrightarrow Y$.

Choose one component $Z_0=Z_{0,\beta_0}$ of
$(Y_0)_{\mathrm{red}}$. It is also an irreducible component of
$Y_{\mathrm{red}}$, by
Proposition~\ref{prop:lcy-equidimensional-factors}(ii). Hence
\[
 V:=(Z_0)_{\mathrm{sm}}\setminus
       \bigcup_{\substack{Z\in\operatorname{Irr}(Y_{\mathrm{red}})\\
                         Z\ne Z_0}}Z
\]
is a nonempty smooth connected open subset of $Y_{\mathrm{red}}$
of dimension $\delta$. Open restriction and complex orientation
give
\[
 \begin{gathered}
 \operatorname{res}_V\rho_{P,Y}(a)=a_0[V],\\
 H^{BM}_{2\delta}(V;\F)\xrightarrow{\sim}H^0(V;\F)
       \xrightarrow{\sim}\F,\qquad [V]\longmapsto1,\\
 \epsilon|_{F_PR_X}
   =\operatorname{coeff}_{[V]}\circ
        \operatorname{res}_V\circ\rho_{P,Y}.
 \end{gathered}
\]
Here $\operatorname{coeff}_{[V]}$ is the composite of the two
isomorphisms in the second line.
To compute $T_d(a_1,\ldots,a_d)$, choose a pole stage
containing $a_1\cdots a_d$ and apply this formula to
$\rho_{P,Y}(a_1\cdots a_d)$.

Since the formula recovers the coefficient $a_0$ in the
global theta expansion, it is independent of the chart,
its finite factor, the chosen arc component and open $V$,
and the pole bound. The calculation uses only one nice chart
and one open subset of its finite factor; no global
finite-factor presentation is required.
\end{remark}

Associativity and commutativity make $T_d$ a symmetric multilinear
functional. The word trace refers to these scalar-valued products:
$\epsilon$ is generally not an algebra homomorphism. For instance,
in the torus calculation $\epsilon(z)=\epsilon(z^{-1})=0$ whereas
$\epsilon(zz^{-1})=1$. Nor is the nondegeneracy of the bilinear form
$T_2$ part of the definition.

\subsection{The BV operator and nondegeneracy}
\label{mir:sec-trace-nondegeneracy}

The next two geometric lemmas do not require maximal boundary
or smoothness of the loop mirror. Put
\[
 L:=LH^*_{\Omega}(X;\F)\otimes_{\F}\C,
 \qquad R:=L^0=R_X\otimes_{\F}\C,
\]
and extend $\epsilon$ by $\C$-linearity. Write
$\Delta:L^j\to L^{j-1}$ for the complex-linear extension of
the BV operator of Definition~\ref{def:source-rotation-operator}.
Its pole-stage realizations are denoted by $\Delta_P$. We first prove
\[
                  \epsilon\circ\Delta|_{L^1}=0.
\]
The geometric input is that arcs on the positive log Calabi-Yau $X$ with sufficiently general centers
have neighborhoods in a fixed pole stage containing only arcs. This is not true in general, e.g. for $\bA^1$. We begin with the classical bounded-degree obstruction to
$\A^1$-uniruledness supplied by the log Calabi--Yau condition.
For the smooth-boundary version, compare
\cite[Proposition~1.1]{TakahashiAffineLines}; the following proof
gives the SNC version needed here.

\begin{lemma}
\label{mir:lem-no-bounded-affine-lines}
Fix a very ample line bundle $H$ on $M$ and an integer $d\ge1$.
There is a dense Zariski-open subset $V_d\subset X$ such that
no nonconstant morphism $f:\Pj^1\to M$ with
$\deg f^*H\le d$ and $f^{-1}(D)\subset\{\infty\}$ has
$f(0)\in V_d$.
\end{lemma}
\begin{proof}
Maps of bounded $H$-degree form a finite-type parameter space:
in a projective embedding they are given by tuples of homogeneous
polynomials of bounded degree, with no common zero, satisfying
the equations of $M$. The condition that the inverse image of
$D$ be supported at $\infty$ is constructible. Indeed, choose
$c_i$ with $c_iH-D_i$ ample; a map not contained in $D_i$
has $0\le\deg f^*D_i\le c_i d$. After a finite
stratification, its contact orders at $\infty$ are fixed and
the condition is the prescribed vanishing of finitely many
coefficients, with nonvanishing leading coefficients.

Suppose that evaluation at $0$ from one irreducible family of
nonconstant such maps dominates $X$. In characteristic zero,
at a general smooth point of the reduced parameter family its
evaluation differential is surjective. For its map $f$ put
$E=f^*T_M(-\log D)$. A first-order variation preserving the
contact orders is a section of $E$: in a boundary coordinate
with pullback $t^m u(t)$, its variation is divisible by $t^m$.
Thus $H^0(\Pj^1,E)\to E|_0$ is surjective. Write
$E=\bigoplus_j\mathcal O(b_j)$. Surjectivity forces every
$b_j\ge0$, whereas $K_M+D\sim0$ gives $\sum_jb_j=0$.
Hence $E\cong\mathcal O^n$.

The logarithmic differential of $f$ gives
\[
 T_{\Pj^1}(-\infty)=\mathcal O(1)
                   \longrightarrow E\cong\mathcal O^n.
\]
It is zero because $\operatorname{Hom}(\mathcal O(1),
\mathcal O)=0$, contradicting nonconstancy in characteristic
zero. No irreducible parameter family therefore dominates $X$.
Remove the closures of their finitely many evaluation images.
Their complement is the required $V_d$. More precisely, evaluation
from the bounded-degree parameter scheme has empty inverse image
over $V_d$. The same exclusion therefore holds after any extension
of the ground field.
\end{proof}

Recall from Proposition~\ref{prop:coherent-source-rotation} that
\[
 r_\lambda:F_P\mathscr LX\longrightarrow F_P\mathscr LX,
 \qquad (r_\lambda\gamma)(t)=\gamma(\lambda t),
 \qquad\lambda\in\C^*.
\]
An open subset $U\subset(F_P\mathscr LX)_{\mathrm{red}}$ is
\emph{source-rotation-invariant} if $r_\lambda(U)=U$ for every
$\lambda\in\C^*$. On the arc stage, center evaluation satisfies
\[
 \operatorname{ev}_0:\cL^+(X)\longrightarrow X,
 \qquad \operatorname{ev}_0(\gamma)=\gamma(0),\qquad
 \operatorname{ev}_0\circ r_\lambda=\operatorname{ev}_0.
\]
Thus $\operatorname{ev}_0^{-1}(V)$ is invariant for every
$V\subset X$; the content of the next lemma is its openness
in a larger pole stage.

\begin{lemma}[An open arc locus at general centers]
\label{mir:lem-generic-arc-open}
For every $P\ge0$, there is a dense Zariski-open subset
$V_P\subset X$ such that
\[
 U_P:=\operatorname{ev}_0^{-1}(V_P)
       \subset\cL^+(X)\subset(F_P\mathscr LX)_{\mathrm{red}}
\]
is open and source-rotation-invariant. The restriction of the
normalized dualizing complex satisfies
\[
 R\Gamma(U_P^{\mathrm{an}},\omega_{P,\Omega}|_{U_P^{\mathrm{an}}})
       \simeq C^*(V_P;\F)
       \quad\text{in }\mathcal D(\F).
\]
This equivalence is compatible with the coherent circle actions,
where the action on the right is trivial. Consequently the
BV operator on this restricted complex is nullhomotopic.
\end{lemma}
\begin{proof}
\emph{Step 1: the open set of centers.}
Choose $r>0$, sections $s_0,\ldots,s_m$ and affine coordinates
$g_j=s_j/s_0$ as in Proposition~\ref{prop:loop-sectors}, so that
$\operatorname{div}(s_0)=rA$. If $P=0$, take $V_0=X$.
For $P>0$, apply Lemma~\ref{mir:lem-no-bounded-affine-lines} with
\[
                H=\mathcal O_M(rA),\qquad d=rP,
\]
and write $V_P:=V_d$. This change of subscript records the pole
bound determining $d$; it is not a second choice of open subset.
Every coordinate of a point of $F_P\mathscr LX$ has Laurent
expansion starting at exponent $-d$.

\emph{Step 2: reduction to a one-parameter formal family.}
Put $S_P=(F_P\mathscr LX)_{\mathrm{red}}$ and
$S_0=\cL^+(X)$. Suppose that a geometric point
$\gamma_0\in U_P$ belongs to the Zariski closure of
$S_P\setminus S_0$. Choose a nice \'etale chart through
that point and a finite-factor presentation of the pulled-back
arc inclusion, as in
Subsection~\ref{subsec:volume-index-preparations}:
\[
 \begin{gathered}
 h:U=Y\times\A^J\longrightarrow S_P,
 \qquad h^{-1}(S_0)=Z\times\A^J,\qquad Z\hookrightarrow Y,\\
                    Y,Z\text{ of finite type over }\C.
 \end{gathered}
\]
Passing to reductions does not change this topological argument.
Retain in $Y$ all coefficients $g_{j,k}$ with $-d\le k\le0$;
there are finitely many of them. In particular the underlying closed subset of $Z$ is cut out
by the coefficients with $k<0$, and the constant coefficients
on $Z$ define a morphism $Z\to X$ through which arc-center
evaluation on $h^{-1}(S_0)$ factors.

After extending the residue field, choose a lift
$\widetilde\gamma_0=(y_0,(u_j^0)_{j\in J})$ of $\gamma_0$.
An \'etale morphism is open, so $\widetilde\gamma_0$ lies in
$\overline{(Y\setminus Z)\times\A^J}$. The product projection
is open: on an affine open of $Y$, a basic open in the product
is defined by a polynomial in finitely many affine coordinates,
and its image is the union of the principal opens defined by
the coefficients of that polynomial. Thus
$y_0\in\overline{Y\setminus Z}$. Over an
algebraically closed extension $K\supset\C$ making this a
$K$-rational point, choose an integral curve through $y_0$ in
$Y_K$ whose generic point lies outside $Z_K$. Such a curve
exists by intersecting an irreducible component of
$\overline{Y_K\setminus Z_K}$ through $y_0$ with sufficiently
general hypersurfaces through $y_0$, and taking a curve component
not contained in $Z_K$. Normalize the curve, choose a point over
$y_0$, and complete its local ring.
Enlarging $K$ if necessary,
this gives
\[
 c:\Spec K\llbracket s \rrbracket\longrightarrow Y,
 \qquad c(0)=y_0,\qquad c(\eta)\notin Z,
 \quad\eta=\Spec K(\!( s )\!).
\]
Send each omitted coordinate $u_j$ to the constant $u_j^0$.
The resulting map to $U$ and then to $S_P$ is a
\emph{one-parameter formal family of loops}: its coordinates are
\[
 g_j(s,t)=\sum_{k\ge-d}a_{jk}(s)t^k,
 \qquad a_{jk}(s)\in K\llbracket s \rrbracket,
 \qquad 1\le j\le m.
\]
This is the loop analog of a {\it wedge} used in Nash-type problems.  Since the
special loop is $\gamma_0$ and the generic loop is not an arc,
\[
 a_{j,-\ell}(0)=0\ (\ell>0),\qquad
 a_{j,-\ell}\ne0\text{ for at least one pair }(j,\ell),
 \qquad a_{j,0}(0)=g_j(x),
\]
where $x=\operatorname{ev}_0(\gamma_0)\in V_P(K)$.

\emph{Step 3: extracting a forbidden affine line.}
The nonzero negative coefficients are a finite set. Write its
minimum slope in lowest terms as
\[
 \frac ab=\min_{\substack{1\le j\le m,\ 1\le\ell\le d\\
                         a_{j,-\ell}\ne0}}
       \frac{\operatorname{ord}_s(a_{j,-\ell})}{\ell},
 \qquad a,b\in\Z_{>0},\quad\gcd(a,b)=1.
\]
Substitute $s=\tau^b$, $t=\tau^a w$. A term
$s^\nu t^k$ becomes $\tau^{b\nu+ak}w^k$. For $k<0$
the minimum-slope inequality makes this exponent nonnegative;
for $k\ge0$ it is nonnegative as well. For each $N$, only
finitely many pairs $(\nu,k)$ contribute modulo $\tau^{N+1}$,
because $k\ge-d$ and $a,b>0$. Thus
\[
 g_j(\tau^b,\tau^a w)\in K[w,w^{-1}]\llbracket \tau \rrbracket.
\]
Its reduction modulo $\tau$ has the form
\[
 p_j(w^{-1})=g_j(x)+
 \sum_{\substack{1\le\ell\le d\\
          b\operatorname{ord}_s(a_{j,-\ell})=a\ell}}
      c_{j\ell}w^{-\ell},
 \qquad c_{j\ell}\in K^*.
\]
At least one $p_j$ is nonconstant, and each has degree at most
$d$ as a polynomial in $u=w^{-1}$. Substitution and reduction
preserve every defining polynomial of $X\subset\A^m$, so
\[
 f:\A^1_K\longrightarrow X_K,
 \qquad u\longmapsto(p_1(u),\ldots,p_m(u)),
 \qquad f(0)=x,
\]
is a nonconstant morphism. Put $e=\max_j\deg p_j\le d$.
Homogenizing $(1,p_1,\ldots,p_m)$ to degree $e$ gives
\[
 \bar f([u:v])=
 [v^e:v^ep_1(u/v):\cdots:v^ep_m(u/v)]\in M_K\subset\Pj^m_K.
\]
There is no common zero: $v^e\ne0$ for $v\ne0$, and at
$v=0$ at least one leading coefficient of degree $e$ is nonzero.
The equations of $M$ hold on the dense open $v\ne0$, hence
everywhere. Therefore
\[
 \deg\bar f^*H=e\le d,
 \qquad\bar f^{-1}(D)\subset\{\infty\},\qquad\bar f(0)=x.
\]
This contradicts Step~1, including its base-change assertion.

We have proved $U_P\cap\overline{S_P\setminus S_0}=\varnothing$.
Since $S_0\setminus U_P$ is closed in $S_P$, the identity
\[
 U_P=S_P\setminus
       \bigl(\overline{S_P\setminus S_0}\cup(S_0\setminus U_P)\bigr)
\]
proves that $U_P$ is open. Its rotation invariance follows from
$\operatorname{ev}_0\circ r_\lambda=\operatorname{ev}_0$.

\emph{Step 4: the restricted complex and its BV operator.}
Let $i_{0,P}:\cL^+(X)^{\mathrm{an}}\hookrightarrow
(F_P\mathscr LX)^{\mathrm{an}}$ be the closed arc inclusion.
The proper-trace map of
Proposition~\ref{prop:normalized-pole-sheaves},
\[
 (i_{0,P})_*\omega_{0,\Omega}\longrightarrow\omega_{P,\Omega},
\]
restricts to an isomorphism on $U_P^{\mathrm{an}}$: the closed
inclusion obtained by base change to this open set is its
identity map. Thus the normalized dualizing complex on $U_P$
is the arc-stage complex restricted to arcs centered in $V_P$.
The normalized
jet comparisons of Proposition~\ref{prop:stable-strata} and
Theorem~\ref{thm:canonical-complex}, restricted over $V_P$, give
\[
 R\Gamma(U_P^{\mathrm{an}},\omega_{P,\Omega}|_{U_P^{\mathrm{an}}})
             \simeq C^*(V_P;\F).
\]
Center evaluation and inclusion of constant arcs are equivariant
for the trivial action on $V_P$. Their homotopy inverse comparison
is the contraction
\[
 [0,1]\times U_P^{\mathrm{an}}\longrightarrow U_P^{\mathrm{an}},
 \qquad(c,\gamma)\longmapsto(t\mapsto\gamma(ct)),
\]
which commutes with every $r_\lambda$. The normalized finite-jet
comparisons are equivariant as well. Thus the displayed derived
equivalence intertwines the coherent circle actions, exactly as in
the proof of Lemma~\ref{lem:source-dilation}. On the right the
parameter-chain action factors through
$C_{-*}(S^1;\F)\to\F$; this augmentation sends a fundamental
$1$-cycle to zero. Hence the induced BV operator is
nullhomotopic.
\end{proof}

\begin{proposition}
\label{mir:prop-trace-kills-bv}
For every $a\in L^1$, one has $\epsilon(\Delta a)=0$.
\end{proposition}
\begin{proof}
Set $B_{\Omega,P,\C}:=B_{\Omega,P}\otimes_{\F}\C$.
The filtered-colimit description of loop cohomology gives
\[
 L^1=\mathop{\operatorname{colim}}_P H^1(B_{\Omega,P,\C}).
\]
Choose $P$ and $a_P\in H^1(B_{\Omega,P,\C})$ representing $a$.
For the open inclusion $j_P:U_P^{\mathrm{an}}\hookrightarrow
(F_P\mathscr LX)^{\mathrm{an}}$, define the morphism
\[
 \rho_P:B_{\Omega,P,\C}
 \xrightarrow{\operatorname{res}_{j_P}}
 R\Gamma(U_P^{\mathrm{an}},j_P^{-1}\omega_{P,\Omega})
                      \otimes_{\F}\C
 \xrightarrow{\sim}C^*(V_P;\C)
\]
in $\mathcal D(\C)$, using Lemma~\ref{mir:lem-generic-arc-open}.
Equivariance of open restriction and of the final comparison gives
\[
 H^0(\rho_P)\bigl(\Delta_P a_P\bigr)
       =\Delta_{V_P}\bigl(H^1(\rho_P)(a_P)\bigr)=0.
\]
Here $\Delta_P$ denotes the cohomology operator induced by
$\overline\Delta_P\otimes_{\F}\C$ on
$B_{\Omega,P,\C}$; in degree one it is
$H^1(\overline\Delta_P\otimes_{\F}\C)$.
The specified pole-stage comparison identifies it with the operator induced
on $F_P\mathbf B_\Omega\otimes_{\F}\C$ by the chosen
cochain map; $\Delta_{V_P}=0$ by that lemma.

The pole-stage theta classes with nonzero contact restrict
to zero on $U_P$. Indeed, on a nice chart each local component
summand is supported on the product of its reduced finite-factor
closure with the affine tail. The corresponding non-arc contact piece is dense
in this support. If the support met the open inverse image of
$U_P$, that contact piece would also meet it, contradicting
$U_P\subset\cL^+(X)$. The arc component restricts,
under normalized complex orientation, to the constant class
$1\in H^0(V_P;\C)$. Therefore the component formula of
Remark~\ref{rem:arc-trace-finite-factor} gives, for
$b\in H^0(B_{\Omega,P,\C})$,
\[
 H^0(\rho_P)(b)=\epsilon(u_Pb)\,1,
 \qquad u_P:H^0(B_{\Omega,P,\C})\longrightarrow H^0(B_{\Omega,\C})=R.
\]
The BV maps commute with the pole maps, so
$u_P(\Delta_Pa_P)=\Delta a$. Substitution yields
$\epsilon(\Delta a)\,1=0$. Since $V_P\ne\varnothing$, its
constant class $1$ is nonzero, proving the assertion.
\end{proof}

\begin{theorem}[Nondegeneracy of the loop trace]
\label{mir:thm-trace-nondegeneracy}
Suppose that $D$ has a zero-dimensional stratum and the loop
mirror $\check X=\Spec R$ is smooth.
Then the bilinear pairing
\[
             R\times R\longrightarrow\C,
                  \qquad(f,g)\longmapsto\epsilon(fg)
\]
is nondegenerate, that is,
\[
 R\longrightarrow\operatorname{Hom}_{\C}(R,\C),\qquad
 f\longmapsto\bigl(g\longmapsto\epsilon(fg)\bigr)
\]
is injective.
The same holds over $\F$ for $R_X$.
\end{theorem}
\begin{proof}
Let $I=\{f\in R:\epsilon(fg)=0\text{ for all }g\in R\}$
be the radical of this pairing. Associativity makes $I$ an
ideal, and $\epsilon(1)=1$ makes it proper. By
Theorem~\ref{thm:intro-global-polyvectors}, every derivation $D$ of $R$ has
the form $D(f)=\{a,f\}$ for some $a\in L^1$.
By Theorem~\ref{thm:intro-lcy-zero}, $L^{-1}=0$, so
$\Delta f=0$ for every $f\in R=L^0$; hence the bracket convention
\eqref{eq:general-rotation-bracket} reads
\[
                  D(f)=\Delta(af)-(\Delta a)f.
\]
For $f\in I$ and $g\in R$, Proposition~\ref{mir:prop-trace-kills-bv}
and the derivation rule give
\[
 \begin{aligned}
 0=\epsilon(\Delta(afg))
  &=\epsilon\bigl(D(f)g+fD(g)+(\Delta a)fg\bigr)\\
  &=\epsilon(D(f)g).
 \end{aligned}
\]
Thus $I$ is stable under all derivations.

Lemma~\ref{globpv:lemma-differential-ideals} now gives $I=0$:
the only ideals of the smooth connected affine algebra $R$ stable
under all $\C$-derivations are $0$ and $R$, and $\epsilon(1)=1$
excludes $I=R$.
Finally, write $\epsilon_{\C}$ for the complexification of the
original functional $\epsilon:R_X\to\F$.
For $0\ne f\in R_X$, faithful extension of scalars gives
$f\otimes1\ne0$ in $R=R_X\otimes_{\F}\C$.
Complex nondegeneracy supplies $g\in R$ such that
$\epsilon_{\C}((f\otimes1)g)\ne0$.
Write $g=\sum_{\nu=1}^N g_\nu\otimes c_\nu$, where
$N\in\Z_{>0}$, $g_\nu\in R_X$, and $c_\nu\in\C$. Then
\[
 0\ne\epsilon_{\C}((f\otimes1)g)
       =\sum_{\nu=1}^N c_\nu\,\epsilon(fg_\nu).
\]
At least one $\epsilon(fg_\nu)$ is nonzero, proving
nondegeneracy over $\F$. This uses only finite sums in the
algebraic tensor product, not a base-change identification
for the full linear dual.
\end{proof}
\begin{proposition}[Reconstruction from two-input and three-input traces]
\label{mir:prop-trace-reconstruction}
Let $k$ be a field, let $V$ be a unital commutative associative
$k$-algebra, and let $\epsilon:V\to k$ be a $k$-linear functional. Put
\[
 T_2(a,b)=\epsilon(ab),\qquad T_3(a,b,c)=\epsilon(abc).
\]
Suppose that $T_2$ is nondegenerate, meaning that
\[
 V\longrightarrow\operatorname{Hom}_k(V,k),\qquad
 v\longmapsto T_2(v,-)
\]
is injective. Then the product is uniquely determined by
\[
             T_2(ab,c)=T_3(a,b,c)\qquad(a,b,c\in V).
\]
A $k$-linear isomorphism between two such algebras which preserves
$T_2$ and $T_3$ is a unital algebra isomorphism.
No finite-dimensionality assumption is required.
\end{proposition}
\begin{proof}
For fixed $a,b$, the element $ab$ satisfies the last displayed identity.
If $v\in V$ also satisfies $T_2(v,c)=T_3(a,b,c)$ for every $c$,
then $T_2(v-ab,c)=0$ for every $c$, hence $v=ab$.

Let $\Phi:V\xrightarrow{\sim}V'$ preserve the two traces, and
write $T'_2,T'_3$ for those on $V'$. For all $a,b,c\in V$,
\[
 \begin{aligned}
 T'_2\bigl(\Phi(ab)-\Phi(a)\Phi(b),\Phi(c)\bigr)
 &=T_2(ab,c)-T'_3\bigl(\Phi(a),\Phi(b),\Phi(c)\bigr)\\
 &=T_3(a,b,c)-T_3(a,b,c)=0.
 \end{aligned}
\]
Surjectivity of $\Phi$ and nondegeneracy of $T'_2$ give
$\Phi(ab)=\Phi(a)\Phi(b)$. Consequently $\Phi(1_V)$ acts as
identity on all of $V'$, so $\Phi(1_V)=1_{V'}$.
\end{proof}

The reconstruction from $T_2$ and $T_3$ in
Proposition~\ref{mir:prop-trace-reconstruction} has a direct
precedent in \cite[Proposition~3.5]{JohnstonComparison}.
There the same argument, applied to the two-input and three-input
coefficient traces, identifies the Keel--Yu and Gross--Siebert
mirror algebras \cite[proof of Corollary~1.2]{JohnstonComparison}.

\section{Comparison with the  Gross--Siebert algebra}
\label{mir:sec-enumerative}\label{mir:sec-gs-comparison}\label{mir:sec-potentials-periods}

The goal of this section is to finish the proof of 
Theorem~\ref{mir:thm-intro-gs}. Nondegeneracy of the loop trace was proved in the previous section.
Let $(M,D)$ be a positive log Calabi--Yau pair of dimension $n$, put
$X=M\setminus D$, and retain the volume form $\Omega$ with
$\operatorname{div}_M(\Omega)=-D$. Retain the ample boundary divisor
$A=\sum_{i=1}^r k_iD_i$, $k_i\in\Z_{>0}$, from
Section~\ref{mir:sec-trace-potentials},  write
\[
 R_X=LH^0_\Omega(X),\qquad
 \epsilon:R_X\longrightarrow\F
\]
for the degree-zero fusion algebra and its loop trace
from Definition~\ref{mir:def-arc-trace}, and write 
\[
 \check X=\Spec(R_X\otimes_{\F}\C)
\]
for  the affine loop mirror
of Section~\ref{mir:sec-affine-scheme}. 

\subsection{How a collision produces a global curve}
\label{mir:subsec-global-curve-geometry}
We give  the geometric intuition behind the proof 
of Theorem~\ref{mir:thm-intro-gs}. We defined the fusion product by colliding loops and  we need to show that it also counts global rational curves in some way, more precisely, stable logarithmic maps as in (\ref{mir:eq-log-descendant}). We explain now how collisions of loops produce global curves.

Let us focus on  the idea behind  showing that the two and three-input loop trace agrees with the Gross-Siebert trace.
Fix $k\in\{2,3\}$ and nonzero input labels
$p_1,\ldots,p_k\in\Lambda_D$. Choose distinct constants
$c_1=0,c_2=1,c_3$ (omit $c_3$ when $k=2$) and move the
marked points on the affine $z$-line by
\[
 z=sc_i,\qquad s\in\A^1.
\]
For $s\ne0$ their formal neighborhoods are disjoint; for $s=0$
the markings meet at $z=0$. Maps from the punctured formal
neighborhoods to $X$ are the input loops. Their local parameters
are $z-sc_i$; the output parameter is $z$. We assume that the output is an arc, since this is what the loop trace $\epsilon$ witnesses from collisions.

A \emph{collision branch} is a one-dimensional geometric branch
of the reduced horizontal closure of the prescribed input-contact
locus, at the generic point of the output arc component in a
finite factor. Here horizontal closure means closure from
$s\ne0$; components supported on $s=0$ are excluded. The geometric
points of the normalization above this generic output point,
and their discrete valuation rings, are specified in
Subsection~\ref{mir:subsec-collision-normal-forms}.
After lifting a branch to the loop family, extending its residue
field to an algebraically closed field $K$, and choosing a
uniformizer $\tau$, its collision parameter is $s=\tau^e$ for
some $e\ge1$. The resulting family is a morphism
\begin{equation}\label{mir:eq-geometric-collision-source}
 \gamma:\Spec K\llbracket \tau,z \rrbracket[D_{\tau^e}(z)^{-1}]
                  \longrightarrow X,
 \qquad D_{\tau^e}(z)=\prod_{i=1}^k(z-\tau^ec_i).
\end{equation}
Its restriction at $\tau=0$ extends to an arc
$\alpha:\Spec K\llbracket z \rrbracket\to X$; put $x=\alpha(0)$.
The variable $\tau$ parametrizes a branch in the space of maps,
whereas $z$ is a source coordinate. In particular,
$T=\Spec K\llbracket \tau,z \rrbracket$ is a two-dimensional parameter--source
surface.

Regard \eqref{mir:eq-geometric-collision-source} as a rational
map $T\dashrightarrow M$. 
Each curve $V(z-\tau^e c_i)\subset T$ maps isomorphically
to $\Spec K\llbracket \tau \rrbracket$ under the projection
$T\to\Spec K\llbracket \tau \rrbracket$; we call these curves the
\emph{marked sections}. Their strict transforms are
pairwise disjoint on the blowup
\begin{equation}\label{mir:eq-geometric-collision-blowup}
 B_e=\operatorname{Bl}_{(\tau^e,z)}T,
 \qquad z=\tau^ew\ \text{on one chart},\qquad
 zv=\tau^e\ \text{on the other},\qquad v=w^{-1}.
\end{equation}
The reduced exceptional curve is $\Pj^1_K$.
The strict transform of $V(z-\tau^e c_i)$ meets it
at $w=c_i$, whereas the strict transform of $V(\tau)$
meets it at $w=\infty$. This last strict
transform is the \emph{outer disk}; its target map is the arc $\alpha$,
so its attachment to the exceptional curve maps to $x\in X$.
Resolving the target map and making the source semistable can
insert further rational components. The reduced exceptional
locus is then a proper tree of rational curves carrying the
resolved map to $M$.

Restrict the target map to this proper tree, replace its
attachment to the outer disk by an output marking, and stabilize
the resulting logarithmic map. The input contacts are $p_i$;
the output has contact zero because its image is $x\in X$.
Only components contracted by the target map and unstable as
marked curves are removed. This is where the \emph{global curve}
enters: the exceptional components are proper, although the
original map \eqref{mir:eq-geometric-collision-source} was defined
on a punctured formal source.

After the target modification and genericity choices specified
in Subsection~\ref{mir:subsec-fixed-source-maps},
Lemma~\ref{mir:lem-collision-normal-form}(i) shows that the
stable logarithmic map extracted from a collision branch has
smooth source $\Pj^1_K$. Its proof first identifies this map
as a point of the logarithmic moduli space over the open set
of Theorem~\ref{mir:thm-generic-log-maps}(i), and then applies
that theorem. Lemma~\ref{mir:lem-collision-normal-form}(ii)
constructs the converse collision family. Finally, Lemma~\ref{mir:lem-collision-branches}
identifies collision branches with these maps and proves
that each branch has multiplicity one. Thus the loop trace
counts the extracted global curves.

\subsection{The Gross--Siebert algebra}
\label{mir:subsec-gs-specialization}
Recall the contact-label set $\Lambda_D$ and its contact weight
from Subsection~\ref{subsLa}. For $p=(\vct,C)\ne0$ and for the
zero label, respectively,
\[
 w_A(p)=\sum_{i=1}^r k_i v_i,\qquad w_A(0)=0.
\]
The label $p=(\vct,C)$ specifies the vector $\vct$ in the cone
associated with $C\in\pi_0(D_{I(\vct)}^\circ)$; $0$ labels the
arc component. These labels index both the loop classes
$\vartheta_p$ and the Gross--Siebert theta functions
$\vartheta_p^{\mathrm{GS}}$. There are finitely many labels of
any bounded weight, since $k_i>0$ and each boundary stratum has
finitely many connected components.

Let $\Lambda=N_1(M)_{\Z}/\mathrm{tors}$ be the lattice of numerical
curve classes, and put $\deg_A(\beta)=A\cdot\beta$. Write $\overline{\mathrm{NE}}(M)\subset N_1(M)_{\R}
=\Lambda_{\R}$ for the closure of the convex cone generated
by the numerical classes of effective curves on $M$.
Choose a rational polyhedral cone $\sigma\subset\Lambda_{\R}$ such that
\[
 \overline{\mathrm{NE}}(M)\subset\sigma,
 \qquad \deg_A(\beta)>0\quad(0\ne\beta\in\sigma).\]
Put $Q:=\sigma\cap\Lambda$, viewed as an additive monoid.
It contains all effective numerical curve classes and indexes
the curve monomials $q^\beta$ in the coefficient ring $\F[Q]$,
with $q^\beta q^\gamma=q^{\beta+\gamma}$.
Such a cone $\sigma$ exists: the section
$\overline{\mathrm{NE}}(M)\cap\{\deg_A=1\}$ is compact by
ampleness, and is contained in a rational polytope in that
hyperplane. The cone over this polytope has the required properties.
The monoid $Q$ is finitely generated and satisfies
\[
 m\beta\in Q,\quad \beta\in\Lambda,\quad m\in\Z_{>0}
       \ \Longrightarrow\ \beta\in Q,
 \qquad Q\cap(-Q)=\{0\};
\]
these are saturation and sharpness, respectively.

For $p,q,r\in\Lambda_D$ and $\beta\in Q$, let
$N_{p,q}^{r,\beta}\in\bQ$ be the genus-zero punctured logarithmic
Gromov--Witten invariant of Gross--Siebert
\cite[Definitions~3.11 and~3.22]{GSIntrinsic}, for basic stable
punctured logarithmic maps of class $\beta$ with three labelled
contacts $p,q,-r$ and the logarithmic output point constraint.
If $r=(\mathbf v,C_r)\ne0$, choose $z\in C_r$; if $r=0$,
choose $z\in X$. For $r\ne0$, the symbol $-r$ records the
negative output contact vector, not an element of $\Lambda_D$.
The constraint is the one over $z$ in
\cite[Proposition~3.9]{GSIntrinsic}. The invariant is the degree
of the zero-dimensional virtual fundamental class of the proper
constrained moduli stack
\cite[Proposition~3.20 and Lemma~3.21]{GSIntrinsic}.
It is primary: no cotangent-line class is inserted. Its
independence of $z$ and of the chosen logarithmic point constraint is
\cite[Remark~3.10 and Theorem~7.11]{GSIntrinsic}.
The Gross--Siebert multiplication rule is
\begin{equation}\label{mir:eq-gs-product}
 \vartheta_p^{\mathrm{GS}}\vartheta_q^{\mathrm{GS}}
 =\sum_{r\in \Lambda_D}\sum_{\beta\in Q}
       N_{p,q}^{r,\beta}q^\beta\vartheta_r^{\mathrm{GS}}.
\end{equation}
We distinguish the symbols $\vartheta_p^{\mathrm{GS}}$ in this externally constructed
algebra from the fundamental loop classes $\vartheta_p$ until the
comparison has been proved.

\begin{theorem}[Gross--Siebert,
{\cite[Theorem~1.15 and Remark~1.26]{GSIntrinsic}}]
\label{mir:lem-gs-polynomial}
For fixed $p,q\in\Lambda_D$, the sum
\eqref{mir:eq-gs-product} is finite. It defines a commutative,
associative $\F[Q]$-algebra
\[
 R^{\mathrm{GS}}_{M,D}
   =\bigoplus_{p\in\Lambda_D}\F[Q]\vartheta_p^{\mathrm{GS}},
 \qquad 1=\vartheta_0^{\mathrm{GS}}.
\]
\end{theorem}

\begin{remark}\label{rmkBbb}
We record a finiteness estimate used below. Balancing
\cite[Proposition~1.16]{GSIntrinsic}, applied to $A$, gives
\begin{equation}\label{mir:eq-gs-degree-balance}
 N_{p,q}^{r,\beta}\ne0
 \quad\Longrightarrow\quad
 \deg_A(\beta)+w_A(r)=w_A(p)+w_A(q).
\end{equation}
For $C\ge0$, the set $\sigma\cap\{0\le\deg_A\le C\}$
is compact, so its intersection with $\Lambda$ is finite.
Taking $C=w_A(p)+w_A(q)$ bounds both $\beta\in Q$ and
$r\in\Lambda_D$. The divisor $A$ bounds the sums; it does
not enter the definition of $N_{p,q}^{r,\beta}$.
\end{remark}

Define the specialization homomorphism and the specialized algebra by
\begin{equation}\label{mir:eq-gs-specialization}
 \operatorname{aug}:\F[Q]\longrightarrow\F,
 \qquad q^\beta\longmapsto1,
 \qquad
 R^{\mathrm{GS},1}_{M,D}
   =R^{\mathrm{GS}}_{M,D}\otimes_{\F[Q],\operatorname{aug}}\F.
\end{equation}
This is a base change of a polynomial multiplication law. In
particular it does not require evaluation of a completed series at
$q^\beta=1$. The specialized algebra still has the basis
$\{\vartheta_p^{\mathrm{GS}}\}$, and its arc coefficient is
\[
 \epsilon_{\mathrm{GS}}\Bigl(\sum_p b_p\vartheta_p^{\mathrm{GS}}\Bigr)=b_0.
\]

\subsection{Logarithmic descendants and the two-input and three-input comparison}\label{mir:subsec-low-trace-comparison}

We recall the logarithmic descendants defined in
\cite[\S9, before Remark~9.1]{GSIntrinsic}.
\begin{definition}[Stable logarithmic maps and contacts]
\label{mir:def-stable-log-map}
We use basic stable logarithmic maps in the sense of
\cite[Theorem~1.1 and Definitions~1.3, 1.6 and~1.20]{GSLogGW},
with fine saturated log structures. The target carries the
divisorial log structure
$\mathcal M_{(M,D)}=\mathcal O_M\cap j_*\mathcal O_X^\times$,
where $j:X\hookrightarrow M$. A map over a log scheme
$(T,\mathcal M_T)$ has the form
\[
 f:(C;x_1,\ldots,x_d,x_{\mathrm{out}})\longrightarrow(M,D),
 \qquad C\longrightarrow T,
\]
where $C\to T$ is a proper flat family of connected nodal
genus-zero curves with disjoint markings in its smooth locus. 
We call $x_1,\ldots,x_d$ the \emph{input markings}
and $x_{\mathrm{out}}$ the \emph{output marking}.
Every component contracted by the underlying map has at least
three markings or nodes. The map is enhanced to a morphism of
log schemes from the marked logarithmic curve to $(M,D)$.
The adjective \emph{basic} specifies the canonical minimal base
log structure; every stable log map is obtained from a basic one
by a unique change of log structure on the same underlying base,
up to unique isomorphism \cite[Proposition~1.24]{GSLogGW}.

Write $\overline{\mathcal M}_T=\mathcal M_T/\mathcal O_T^\times$
for the characteristic sheaf, and similarly for the source and
target. At a marking $x_i$ over a geometric point $t$ of $T$,
the contact homomorphism is
\[
 \overline{\mathcal M}_{(M,D),f(x_i)}
 \longrightarrow\overline{\mathcal M}_{C,x_i}
 \cong\overline{\mathcal M}_{T,t}\oplus\mathbb N
 \longrightarrow\mathbb N.
\]
Its values on the local boundary generators are the contact
orders. Together with the prescribed component of the boundary
stratum, these orders specify the label in $\Lambda_D$.
\end{definition}

Fix $d\ge2$, labels $p_j=(\vct_j,C_j)\in\Lambda_D\setminus\{0\}$,
a class $\beta\in Q$, and a general point $x\in X$.
Consider the proper moduli stack
\[
 \overline{\mathcal M}^{\log}_{0,d+1}
       ((M,D),\beta;p_1,\ldots,p_d,0)_x.
\]
Its objects are the basic maps of
Definition~\ref{mir:def-stable-log-map}, with curve class $\beta$
and contact labels $(p_1,\ldots,p_d,0)$. The compactified contact
condition includes specialization to deeper boundary strata;
it is not the open incidence condition
$f(x_j)\in C_j$ on the ordinary stable-map space.
The subscript $x$ denotes the fiber of
$\operatorname{ev}_{\mathrm{out}}([f])=f(x_{\mathrm{out}})$ at
$x$, with its point-constrained virtual class.
Let $\psi_{\mathrm{out}}$ be the first Chern class of the cotangent
line of the source curve at that marking, using the stable-map
cotangent line. Define
\begin{equation}\label{mir:eq-log-descendant}
 L_\beta(p_1,\ldots,p_d):=
 \int_{[\overline{\mathcal M}^{\log}_{0,d+1}
       ((M,D),\beta;p_1,\ldots,p_d,0)_x]^{\mathrm{vir}}}
       \psi_{\mathrm{out}}^{d-2}.
\end{equation}
We call $L_\beta(p_1,\ldots,p_d)$ a genus-zero
\emph{logarithmic descendant Gromov--Witten invariant with
one point insertion}, or a \emph{logarithmic point descendant}
for short. Here $\psi_{\mathrm{out}}^{d-2}$ is the descendant
insertion; for $d=2$, it is $1$, giving a primary invariant.
All marking contacts are nonnegative, so these invariants use
ordinary stable logarithmic maps, without punctures.
The point-constrained virtual class is the one of
\cite[Proposition~3.20]{GSIntrinsic}.

Equivalently, the integral on the unconstrained moduli stack has
the insertion
\[
 \psi_{\mathrm{out}}^{d-2}\operatorname{ev}_{\mathrm{out}}^*[x],
\]
as in \cite[Proposition~4.6]{JohnstonPeriod}.
The logarithmic virtual dimension before the point constraint is
$n-3+(d+1)+c_1(T_M(-\log D))\cdot\beta=n+d-2$.
Since $K_M+D=0$, the constraint leaves dimension $d-2$, explaining
the exponent. An empty moduli stack contributes zero. We use the
standard virtual fundamental class, so the definition does not
assume that the moduli stack is smooth or that a count is transverse.

Since the output contact is zero, logarithmic balancing gives
\[
 L_\beta(p_1,\ldots,p_d)\ne0
 \quad\Longrightarrow\quad
 \deg_A(\beta)=\sum_{j=1}^d w_A(p_j).
\]
Thus Remark~\ref{rmkBbb} ensures that only finitely many
$\beta\in Q$ contribute to
$\sum_{\beta\in Q}L_\beta(p_1,\ldots,p_d)$.

Let $[q^\beta\vartheta_0^{\mathrm{GS}}](a)$ denote the coefficient
of $q^\beta\vartheta_0^{\mathrm{GS}}$ in $a\in R^{\mathrm{GS}}_{M,D}$.
The weak Frobenius structure theorem of Johnston, in the form
\cite[Proposition~4.6]{JohnstonPeriod} following \cite{JohnstonRoot},
states that for $d\ge2$ and nonzero contact labels,
\[
 [q^\beta\vartheta_0^{\mathrm{GS}}]
 (\vartheta_{p_1}^{\mathrm{GS}}\cdots\vartheta_{p_d}^{\mathrm{GS}})
 =L_\beta(p_1,\ldots,p_d).
\]
Summing the finite set of curve
classes and applying \eqref{mir:eq-gs-specialization} gives
\begin{equation}\label{mir:eq-gs-weak-frobenius}
 \epsilon_{\mathrm{GS}}
 (\vartheta_{p_1}^{\mathrm{GS}}\cdots\vartheta_{p_d}^{\mathrm{GS}})
 =\sum_\beta L_\beta(p_1,\ldots,p_d).
\end{equation}
For $d=2,3$ the coefficient identity is already
\cite[Theorem~9.4]{GSIntrinsic}. These are identities in the Gross--Siebert algebra. The comparison
with the loop traces $T_2,T_3$ of
Definition~\ref{mir:def-arc-trace} is the following theorem,
proved in the next three subsections.

\begin{theorem}[Two-input and three-input trace comparison]
\label{mir:thm-low-trace-comparison}
For a positive log Calabi--Yau pair with maximal boundary,
and for $k=2,3$ and nonzero contact labels $p_1,\ldots,p_k$,
\begin{equation}\label{mir:eq-low-trace-comparison}
 \epsilon(\vartheta_{p_1}\cdots\vartheta_{p_k})
 =\epsilon_{\mathrm{GS}}
   (\vartheta_{p_1}^{\mathrm{GS}}\cdots\vartheta_{p_k}^{\mathrm{GS}}).
\end{equation}
The products on the two sides are taken in $R_X$ and
$R^{\mathrm{GS},1}_{M,D}$, respectively. By
\eqref{mir:eq-gs-weak-frobenius}, the right-hand side is
$\sum_{\beta\in Q}L_\beta(p_1,\ldots,p_k)$.
No smoothness assumption on $\check X$ is needed here.
\end{theorem}

\subsection{Generic fixed-source logarithmic maps}
\label{mir:subsec-fixed-source-maps}

Assume from now until the proof of
Theorem~\ref{mir:thm-low-trace-comparison} that $D$ has maximal
boundary. Fix $k\in\{2,3\}$ and
$\mathbf p=(p_1,\ldots,p_k)\in(\Lambda_D\setminus\{0\})^k$.
Write $\mathcal M_{0,k+1}$ for the moduli scheme of smooth genus-zero
curves with $k+1$ distinct labeled markings, and
$\overline{\mathcal M}_{0,k+1}$ for its compactification by stable
nodal marked curves. Here \emph{stable} means that the
normalization of every irreducible component carries at least
three special points, namely markings or preimages of nodes.
A \emph{fixed source}
is a point
$c=[(\Pj^1;c_1,\ldots,c_k,\infty)]\in\mathcal M_{0,k+1}$.

\emph{Divisorial labels and the modified target.}
Here a \emph{logarithmic modification} of the SNC pair is a
proper birational morphism
\[
 \pi:(M',D')\longrightarrow(M,D),\qquad
 D'=(\pi^{-1}D)_{\mathrm{red}},\qquad \pi|_X=\mathrm{id}_X,
\]
which, in the local toric charts of the SNC boundary, is induced
by a subdivision of the boundary cones. We choose projective
subdivisions, refined to regular ones: every cone of a regular
subdivision is generated by part of a basis of its integral
lattice. Thus $M'$ is smooth and $D'$ is SNC. Such modifications
can be obtained by blowups of smooth closed boundary strata and
further toric subdivisions. A nonzero label $p$
is \emph{divisorial on $(M',D')$} if
\[
 p=\ell e_E,\qquad \ell\in\Z_{>0},
\]
where $E$ is an irreducible component of $D'$ and $e_E$ is the
primitive integral generator of its boundary ray. Equivalently,
the valuation represented by $p$ is $\ell\operatorname{ord}_E$.
Subdivide along the finitely many rays $\R_{\ge0}p_i$, then take
a projective regular refinement. This makes
$p_i=\ell_i e_{E_i}$ for components $E_i\subset D'$;
only the selected labels are required to be divisorial.
This is the modification used in \cite[\S5, before Theorem~5.1
and after Lemma~5.7]{JohnstonComparison}; see also
\cite[Definition~1.1]{KeelYu}.

Fix such a $\pi$, followed by the refinement of
Lemma~\ref{mir:lem:cy-boundary-refinement}. Choose an ample divisor
$A'=\sum_\alpha k'_\alpha D'_\alpha$ with $k'_\alpha>0$.
Define 
\[
 \Gamma=\left\{\beta'\in N_1(M')_{\Z}/\mathrm{tors}:
 \begin{array}{l}
 \beta'\text{ is effective},\\[-2pt]
 D'_\alpha\cdot\beta'
   =\displaystyle\sum_{i:E_i=D'_\alpha}\ell_i
   \text{ for every }\alpha
 \end{array}\right\}.
\]
This is the set of possible numerical curve classes. Indeed,
apply the balancing identity used in Remark~\ref{rmkBbb},
in its SNC form \cite[Corollary~1.17]{GSIntrinsic}, to each
boundary component $D'_\alpha$ of the modified pair.
For a stable logarithmic map $f:C\to M'$ of class $\beta'$
with input contacts $p_i=\ell_i e_{E_i}$ and output contact
zero, it gives
\[
 D'_\alpha\cdot\beta'
 =\deg_C f^*\mathcal O_{M'}(D'_\alpha)
 =\sum_{i:E_i=D'_\alpha}\ell_i.
\]
Indeed, the $i$th input contributes $\ell_i$ along $E_i$
and zero along every other boundary component, while the
output contributes zero. These are precisely the equalities
imposed in the definition of $\Gamma$. Consequently,
\[
 A'\cdot\beta'
 =\sum_\alpha k'_\alpha
       \sum_{i:E_i=D'_\alpha}\ell_i
 =\sum_i k'_{E_i}\ell_i,
\]
where $k'_{E_i}$ is the coefficient of $E_i$ in $A'$.
The compactness argument of Remark~\ref{rmkBbb}, applied to
$M'$ and the ample divisor $A'$, shows that only finitely many
effective numerical curve classes have this degree.
Thus $\Gamma$ is finite.

\emph{The moduli spaces and their transverse locus.}
For $\beta'\in\Gamma$, put
\[
 \overline{\mathcal M}^{\log}_{\beta'}
 :=\overline{\mathcal M}^{\log}_{0,k+1}
       ((M',D'),\beta';p_1,\ldots,p_k,0).
\]
For a represented map
$f:(C;x_1,\ldots,x_k,x_{\mathrm{out}})\to M'$, let
$C^{\mathrm{st}}$ be the stable marked curve obtained by
forgetting $f$ and the log structures, and successively
contracting rational components with fewer than three special
points. A contracted component may have been mapped
nonconstantly by $f$. Let $\sigma:C\to C^{\mathrm{st}}$ denote the stabilization
morphism. The induced markings are
\[
 x_i^{\mathrm{st}}:=\sigma(x_i)\quad(1\le i\le k),
 \qquad
 x_{\mathrm{out}}^{\mathrm{st}}:=\sigma(x_{\mathrm{out}}).
\]
They remain distinct smooth points. In particular, when a
rational tail carrying a marking is contracted, that marking
is sent to the attachment point, which becomes smooth after
the contraction. We subsequently omit the superscript
$\mathrm{st}$ from these markings.
The output evaluation and source stabilization are the morphisms
\[
\begin{aligned}
 \operatorname{ev}_{\mathrm{out}}:
 \overline{\mathcal M}^{\log}_{\beta'}&\longrightarrow M',
 & [f]&\longmapsto f(x_{\mathrm{out}}),\\
 \operatorname{st}:
 \overline{\mathcal M}^{\log}_{\beta'}&\longrightarrow
       \overline{\mathcal M}_{0,k+1},
 & [f]&\longmapsto
       [C^{\mathrm{st}};x_1,\ldots,x_k,x_{\mathrm{out}}].
\end{aligned}
\]
Restricting their product gives
\[
 T:=X\times\mathcal M_{0,k+1},\qquad
 \mathcal M^{\log}_{\beta'}
 :=(\operatorname{ev}_{\mathrm{out}},\operatorname{st})^{-1}(T),
 \qquad E_{\beta'}:\mathcal M^{\log}_{\beta'}\longrightarrow T.
\]
In the ordinary stable-map stack of class $\beta'$ with $k+1$
markings and $(\operatorname{ev}_{\mathrm{out}},\operatorname{st})$
in $T$, let $\mathcal M^{\mathrm{tr}}_{\beta'}$ denote the locus
satisfying
\[
 f^{-1}(D')=\sum_{i=1}^k\ell_i x_i
 \quad\text{as effective Cartier divisors},\qquad
 f(x_i)\in E_i^\circ
 :=E_i\setminus\bigcup_{D'_\alpha\ne E_i}D'_\alpha.
\]
For families the divisor equality is imposed on the total source
curve as an equality of relative effective Cartier divisors.
Locally along $x_i$, if $t_i$ is an equation for the marking and
$y_i$ is an equation for $E_i$, it says
$f^*y_i=t_i^{\ell_i}a_i$ with $a_i$ invertible.
Thus no irreducible component of the source lies in $D'$, and
there are no unmarked boundary intersections. By
\cite[Lemma~5.2]{JohnstonComparison}, these maps have unique
basic logarithmic enhancements, identifying
$\mathcal M^{\mathrm{tr}}_{\beta'}$ with an open substack of
$\mathcal M^{\log}_{\beta'}$.
Finally, $T_{M'}(-\log D')$ denotes the locally free sheaf of
derivations preserving the ideal of every $D'_\alpha$.

\begin{theorem}[Keel--Yu and Johnston]
\label{mir:thm-generic-log-maps}
In the preceding setup there is a dense open $\mathcal V\subset T$
with the following properties.
\begin{enumerate}[(i)]
\item For every $\beta'\in\Gamma$,
\[
 E_{\beta'}^{-1}(\mathcal V)\subset\mathcal M^{\mathrm{tr}}_{\beta'},
 \qquad E_{\beta'}^{-1}(\mathcal V)\longrightarrow \mathcal V
       \text{ is representable and finite \'etale}.
\]
Moreover, every map represented by a geometric point of \(E_{\beta'}^{-1}(\mathcal V)\)
 has smooth source $C\cong\Pj^1$ and
$f^*T_{M'}(-\log D')\cong\mathcal O_{\Pj^1}^{\oplus n}$.
\item For $(x,c)\in \mathcal V(\C)$, let
$\mathcal H_{x,c}:=\coprod_{\beta'\in\Gamma}
E_{\beta'}^{-1}(x,c)$. Then
\begin{equation}\label{mir:eq-log-naive-count}
 \sum_{\beta\in Q}L_\beta(p_1,\ldots,p_k)
   =\#\mathcal H_{x,c}(\C).
\end{equation}
The right side is the \emph{total naive rational-curve count}
of \cite[Definition~1.1]{KeelYu}, summed over curve classes.
It counts isomorphism classes of maps to $(M',D')$ with
the prescribed labeled contacts, marked source representing
$c$, and output mapping to $x$, each with weight one.
The descendants on the left are those of the original pair.
\end{enumerate}
\end{theorem}

\begin{proof}
We verify the hypotheses of the relevant results from
\cite{GSCanonical,JohnstonComparison,GSLogGW,KeelYu}.
By Lemma~\ref{mir:lem:cy-boundary-refinement}, $(M',D')$ is a
smooth projective log Calabi--Yau pair with maximal boundary
and connected nonempty boundary intersections. Moreover
$\pi^*A$ is nef and has positive coefficient on every boundary
component. Thus
\cite[Proposition~1.6]{GSCanonical}, applied to its identity
resolution, verifies the assumptions of
\cite[Theorem~1.1]{JohnstonComparison}. The pullback of $\Omega$ has a simple pole along every component
of $D'$, so all the selected boundary rays lie in the
Kontsevich--Soibelman skeleton required there.

Fix $\beta'\in\Gamma$. The geometric argument of
\cite[Lemma~5.8 and the proof of Theorem~1.1, Steps~1 and~3--5]
{JohnstonComparison} supplies a geometric fiber of $E_{\beta'}$
consisting entirely of transverse maps. Here the generic-source
condition is $\operatorname{st}([f])=c$ for general
$c\in\mathcal M_{0,k+1}$. Lemma~5.8 excludes
boundary-contained components on the smallest connected subtree
of the source dual graph containing all markings, and Step~1
excludes such components in the attached unmarked trees.
The local charts of logarithmic curves then give exactly the
boundary-contact conditions defining $\mathcal M^{\mathrm{tr}}_{\beta'}$
\cite[Theorem~1.1]{GSLogGW}.

Since $D'$ is SNC, local equations of its components differ
by units and therefore define global sections generating the
characteristic sheaf of $(M',D')$.
The boundedness and properness results
\cite[Theorems~3.9, 3.12 and~0.2]{GSLogGW} therefore apply
to the fixed discrete data; in particular, $E_{\beta'}$ is proper.
Let
\[
 Z_{\beta'}:=\mathcal M^{\log}_{\beta'}
               \setminus\mathcal M^{\mathrm{tr}}_{\beta'}
\]
with its reduced closed substack structure. The preceding
geometric fiber is disjoint from $Z_{\beta'}$, so
\[
 \mathcal V^{\mathrm{tr}}_{\beta'}
 :=T\setminus E_{\beta'}(Z_{\beta'})
\]
is nonempty open, hence dense because $T$ is irreducible.
Its whole inverse image is transverse.

By \cite[Notation~3.3, Lemma~3.6 and Proposition~3.13]{KeelYu},
as formulated in \cite[Theorem~5.1]{JohnstonComparison},
there is a dense open $\mathcal V_{\beta'}\subset
\mathcal V^{\mathrm{tr}}_{\beta'}$ on which
$E_{\beta'}^{-1}(\mathcal V_{\beta'})\to\mathcal V_{\beta'}$
is representable and finite \'etale. Every represented map has
source stable as a marked curve and trivial logarithmic tangent
pullback. Its stabilized source lies in $\mathcal M_{0,k+1}$,
so its source itself is smooth and isomorphic to $\Pj^1$.

The prescribed contacts fix the degree
\[
 A'\cdot\beta'=\sum_{i=1}^k k'_{E_i}\ell_i
 \qquad(\beta'\in\Gamma).
\]
As shown before the theorem, ampleness of $A'$ therefore makes
$\Gamma$ finite. Consequently
\[
 \mathcal V:=\bigcap_{\beta'\in\Gamma}\mathcal V_{\beta'}
\]
is dense open in the irreducible scheme $T$ and satisfies all
assertions of (i). Empty moduli spaces impose no restriction;
if $\Gamma=\varnothing$, take $\mathcal V=T$.

For (ii), we spell out the source and multiplicity identifications
in Johnston's enumerativity theorem
\cite[Theorem~1.1]{JohnstonComparison}. On a general output
fiber, Step~2 of its proof gives
\[
 \psi_{\mathrm{out}}=\operatorname{st}^*\bar\psi_{\mathrm{out}},
 \qquad
 \bar\psi_{\mathrm{out}}^{\,k-2}=[c]
 \quad\text{in }A^{k-2}(\overline{\mathcal M}_{0,k+1}).
\]
Here $\bar\psi_{\mathrm{out}}$ is the cotangent-line class on
the moduli space of stable marked curves; the second identity
uses $\overline{\mathcal M}_{0,3}=\mathrm{pt}$ and
$\overline{\mathcal M}_{0,4}\cong\Pj^1$.
Each fiber point has multiplicity one: with source and output
fixed its deformation and obstruction groups are
\[
 H^i\bigl(\Pj^1,f^*T_{M'}(-\log D')(-\infty)\bigr)
 =H^i(\Pj^1,\mathcal O(-1)^{\oplus n})=0,
 \qquad i=0,1.
\]
Writing $L^{(M',D')}_{\beta'}$ for the descendants of the modified
pair, these identifications give
\[
 L^{(M',D')}_{\beta'}(p_1,\ldots,p_k)
 =\deg\bigl(E_{\beta'}^{-1}(\mathcal V)\longrightarrow \mathcal V\bigr).
\]
Since $\mathcal V$ is irreducible, this degree equals
$\#E_{\beta'}^{-1}(x,c)(\C)$ for every $(x,c)\in \mathcal V(\C)$.

Finally, we explain the descendant compatibility in logarithmic
modification invariance \cite[Lemma~5.7]{JohnstonComparison}.
Fix $\beta\in Q$ and $x\in X(\C)$. Denote by $\mathscr M_{\beta,x}$ and
$\mathscr M'_{\beta',x}$ the point-constrained stable-log-map
stacks in \eqref{mir:eq-log-descendant} for $(M,D)$ and
$(M',D')$, respectively, with the prescribed contacts
$(p_1,\ldots,p_k,0)$. For $\beta'\in\Gamma$ with
$\pi_*\beta'=\beta$, composition with $\pi$, followed by
stabilization, induces
\[
 \Pi_{\beta'}:\mathscr M'_{\beta',x}
                      \longrightarrow\mathscr M_{\beta,x}.
\]
Every source component contracted in this stabilization maps
into $D'$. Indeed, an originally nonconstant component which
becomes constant under $\pi$ lies in $D'$, since $\pi$ is an
isomorphism over $X$. An originally constant component is
initially stable; if it becomes unstable later, it meets a
previously contracted component and hence has its constant
image in $D'$ as well. Induction through the contractions
proves the assertion, including chains of constant components.
No contracted component therefore contains the output marking,
whose image is $x\in X$. The contraction is an isomorphism
near that marking. Writing $\psi_{\mathrm{out}}$ and
$\psi'_{\mathrm{out}}$ for the output cotangent-line classes
on the original and modified stacks, respectively, we obtain
\[
 \psi'_{\mathrm{out}}
       =\Pi_{\beta'}^*\psi_{\mathrm{out}}.
\]
The evaluation maps commute with these morphisms. Thus virtual
pushforward \cite[Theorem~1.1.1]{AbramovichWise}, together with
commutation of proper pushforward and refined pullback along
$\{x\}\hookrightarrow X$, gives
\[
 \sum_{\substack{\beta'\in\Gamma\\\pi_*\beta'=\beta}}
  (\Pi_{\beta'})_*[\mathscr M'_{\beta',x}]^{\mathrm{vir}}
       =[\mathscr M_{\beta,x}]^{\mathrm{vir}}.
\]
Here balancing restricts all nonempty modified moduli stacks
to the finite set $\Gamma$. Applying the projection formula
with $\psi_{\mathrm{out}}^{k-2}$ gives
\[
 L_\beta(p_1,\ldots,p_k)
 =\sum_{\substack{\beta'\in\Gamma\\\pi_*\beta'=\beta}}
       L^{(M',D')}_{\beta'}(p_1,\ldots,p_k).
\]
Summing this identity proves
\eqref{mir:eq-log-naive-count}.
\end{proof}

Proposition~\ref{prop:fusion-invariance-gradings} identifies the
unfiltered fusion algebras of the original and modified pairs
through the identity on punctured loops. We check now that this
comparison also preserves the loop theta classes and the loop trace.

\begin{lemma}
\label{mir:lem-modification-theta-trace}
For the modification $\pi:(M',D')\to(M,D)$ fixed above, the
subdivision identifies the integral contact labels of the two
pairs. Under this identification, the unfiltered fusion algebra
comparison of Proposition~\ref{prop:fusion-invariance-gradings},
induced by the identity on punctured loops, carries each theta
class $\vartheta_p$ to the theta class with the same label $p$
and preserves the loop trace $\epsilon$.
\end{lemma}
\begin{proof}
For every field $K/\C$, an arc
$\gamma:\Spec K\llbracket t \rrbracket\to M$ whose generic point lies in $X$
has a unique lift to $M'$ by properness of the modification.
In local toric charts, write
\[
 \gamma^*x_i=t^{v_i}a_i(t),\qquad
 a_i(t)\in K\llbracket t \rrbracket^\times,
\]
where the $x_i$ are boundary equations and the $v_i$ are their
contact orders. The unique cone of the subdivision whose relative
interior contains the contact label determines the modified
contact vector: its chart coordinates are Laurent monomials in
the $x_i$, multiplied by units, and their orders are nonnegative.
Those of order zero are units. Thus corresponding labels,
including their stratum-component labels, describe the same
field-valued punctured loops.

On each nice chart of a common containing pole stage, choose a
finite factor presenting the pullbacks of the two closed pole
subschemes and these locally closed contact pieces, as in
the proof of Proposition~\ref{prop:fusion-invariance-gradings}.
The resulting reduced contact pieces in the common finite
factor have the same geometric points, hence coincide. Their
reduced closures also coincide, since both pole subschemes are
closed. The volume form on $X$ is unchanged, so both
normalizations restrict from the same ambient volume index.
The two closed traces consequently give the same normalized
fundamental class. Lemma~\ref{lem:fundamental-smooth-comparison}
makes this equality compatible with refinements and chart
overlaps. The construction in Section~\ref{sec:log-cy-zero}
therefore identifies the entire labeled theta bases. In
particular, it preserves the selected $\vartheta_{p_i}$ and
the coefficient of the unchanged zero label, namely $\epsilon$.
\end{proof}

For the local calculation below, write $(M,D,A)$ for the already
chosen $(M',D',A')$. Thus each selected $p_i$ is divisorial on
this target. On the logarithmic side,
\eqref{mir:eq-log-naive-count} retains the original sum of
descendants. By this identity and
Lemma~\ref{mir:lem-modification-theta-trace}, it suffices to
prove the trace comparison for this refined pair.

Choose distinct real numbers
\[
 c_1=0<c_2=1<\cdots<c_k
\]
whose marked source is general for
Theorem~\ref{mir:thm-generic-log-maps}(i), and put
$\mathcal V_c=\{x\in X:(x,c)\in \mathcal V\}$. This is dense open; for $k=3$
the remaining cross-ratio can be chosen outside a finite set.
The base change of the finite disjoint union in that theorem is
a finite \'etale scheme
\[
 \operatorname{ev}:\mathcal H_c
 :=\coprod_{\beta'\in\Gamma}E_{\beta'}^{-1}(\mathcal V_c\times\{c\})
 \longrightarrow \mathcal V_c.
\]
The first two input markings and the output marking determine
the unique coordinate $w$ on each smooth source satisfying
\[
 w(x_1)=0,\qquad w(x_2)=1,\qquad
 w(x_{\mathrm{out}})=\infty.
\]
For $k=2$ these are all the markings; for $k=3$ the remaining
input has coordinate $c_3$. This identifies the universal
source over $\mathcal H_c$ with $\mathcal H_c\times\Pj^1_w$, and the universal map is
\begin{equation}\label{mir:eq-fixed-source-universal-map}
 F:\mathcal H_c\times\Pj^1_w\longrightarrow M,\qquad
 F(u,\infty)=\operatorname{ev}(u).
\end{equation}
For fixed $u$, the restriction $F(u,-)$ maps 
$(\Pj ^1_w; c_1,\ldots,c_k,\infty)$ into $M$. The output marking $\infty$ has local coordinate $w^{-1}$,
vanishing at that marking.

For every field extension $K/\C$ and $x\in\mathcal V_c(K)$, write
\[
 \mathcal H_{x,c}
 :=\mathcal H_c\times_{\mathcal V_c,x}\Spec K
 \cong\coprod_{\beta'\in\Gamma}E_{\beta'}^{-1}(x,c).
\]
This agrees with the notation of
Theorem~\ref{mir:thm-generic-log-maps}(ii) when $K=\C$.
If $K$ is algebraically closed, $\#\mathcal H_{x,c}(K)$ is the
degree of the finite \'etale map $\operatorname{ev}$, constant
on the irreducible open $\mathcal V_c$, and hence equals the
count in \eqref{mir:eq-log-naive-count}.

\begin{lemma}[The logarithmic differential of the universal map]
\label{mir:lem-fixed-source-log-differential}
The differential of $F$ in the $\mathcal H_c$ directions is an isomorphism
\begin{equation}\label{mir:eq-fixed-source-log-differential}
 d_{\mathcal H_c}F:\operatorname{pr}_{\mathcal H_c}^*T_{\mathcal H_c}
       \xrightarrow{\ \sim\ }F^*T_M(-\log D).
\end{equation}
\end{lemma}
\begin{proof}
Locally near $\mathcal H_c\times\{c_i\}$, if $x$ is an equation for its
boundary component, the fixed contact condition reads
$F^*x=(w-c_i)^{\ell_i}a$, with $a$ invertible. For a vector
field $\xi$ on $\mathcal H_c$,
\[
 \xi(F^*x)=(w-c_i)^{\ell_i}\xi(a)
          =F^*x\,\frac{\xi(a)}a.
\]
Thus $d_{\mathcal H_c}F$ takes values in $F^*T_M(-\log D)$; away from the
marked divisors this is the ordinary differential. On
$\{u\}\times\Pj^1$, both bundles in
\eqref{mir:eq-fixed-source-log-differential} are trivial of
rank $n$, by Theorem~\ref{mir:thm-generic-log-maps}(i). Their
bundle map is therefore a constant matrix. At $w=\infty$ it
is $d\operatorname{ev}_u$, which is invertible because
$\operatorname{ev}$ is \'etale. The determinant is nonzero at
every geometric point, so $d_{\mathcal H_c}F$ is an isomorphism.
\end{proof}

\subsection{Finite-factor collision branches and their normal forms}
\label{mir:subsec-collision-normal-forms}
Recall we have fixed  nonzero contact labels $p_1,\ldots,p_k$  and numbers $c_1=0<c_2=1<\cdots<c_k$ in
Subsection~\ref{mir:subsec-fixed-source-maps}, with $k\in\{2,3\}$.
As there, $(M,D,A)$ denotes the refined pair and polarization,
and $p_i=\ell_i e_{E_i}$.
Put $I=\{1,\ldots,k\}$, $P_i=w_A(p_i)$ and choose the affine generators $g_j$ from a
very ample multiple $rA$, as in
Lemma~\ref{lem:kv-pole-compatibility}.  Write
\[
 m_i=rP_i,\qquad b=\sum_i m_i,\qquad
 G_s(z)=\prod_i(z-sc_i)^{m_i},\qquad
 D_s(z)=\prod_i(z-sc_i).
\]
Here $\ell_i$ is the contact order along $E_i$, whereas $m_i$
is the pole bound for the chosen affine generators at $z=sc_i$.

\emph{The input locus and its finite-factor closure.}
Consider the pullback
\[
 \iota_{\mathbf c}:\A^1_s\longrightarrow\A^k,\quad
 s\longmapsto(sc_i)_i,\qquad
 S=\A^1_s\times_{\A^k}\mathcal Z_{I,\mathbf P}.
\]
Let $\mathcal C_i\subset F_{P_i}\mathscr LX$ be the
exact-contact component labeled by $p_i$. Separated-source
factorization, with local parameters $t_i=z-sc_i$, identifies
\[
 S|_{\Gm}\cong
 \left(\prod_{i=1}^kF_{P_i}\mathscr LX\right)\times\Gm,
 \qquad \Gm=\Spec\C[s,s^{-1}],
\]
and the prescribed exact-contact locus is
\[
 \mathcal C^\circ:=
 \left(\prod_{i=1}^k\mathcal C_i\right)\times\Gm
 \subset S|_{\Gm}.
\]
By Theorem~\ref{thm:relative-finite-factors} and
Proposition~\ref{prop:direct-collision-diagrams}, choose a nice
chart $h:U=Y\times\A^J\to S$ and a finite factor over
$\A^1_s$ presenting the pulled-back contact loci and output
arc component. Write
\[
 \operatorname{pr}_Y^{-1}(V^\circ)
   =(U\times_S\mathcal C^\circ)_{\mathrm{red}},
 \qquad V^\circ\subset Y_{s\ne0},
\]
where $V^\circ$ is reduced and locally closed. Require that the
output-center map factors through the finite factor of the
output arc locus and that, for every chosen generator $g_j$, the functions
\[
 f_{j,d}:=[z^d]\bigl(G_s(z)g_j(\gamma)\bigr),
 \qquad 0\le d\le b,
\]
descend to $Y$. Explicitly, require functions
$f^Y_{j,d}\in\Gamma(Y,\mathcal O_Y)$ such that
\[
 h^*f_{j,d}=\operatorname{pr}_Y^*f^Y_{j,d}.
\]
Here $\gamma$ is the universal punctured-source map
parametrized by $S$, and we recall now the coefficient functions. For a complex algebra $R$ and an $R$-point of $S$ with
parameter $s\in R$, the formal neighborhood of the marked
divisor $D_s(z)=0$ in $\A^1_R$ has completed coordinate ring
\[
 \mathcal O_s
 :=\varprojlim_N R[z]/(D_s(z)^{N+1})
 \cong R\llbracket q_0 \rrbracket[z]/(D_s(z)-q_0),
\]
as in 
\eqref{eq:kv-monic-source}. The pole conditions give
$G_s(z)g_j(\gamma)\in\mathcal O_s$. Since
$1,z,\ldots,z^{k-1}$ form an $R\llbracket q_0 \rrbracket$-basis, write uniquely
\[
 G_s(z)g_j(\gamma)=\sum_{\nu=0}^{k-1}\sum_{h\ge0}
          n_{j,\nu,h}z^\nu q_0^h,
 \qquad n_{j,\nu,h}\in R.
\]
Since $c_1=0$, one has $D_s(z)\in zR[z]$, hence
$(D_s(z)^N)\subseteq(z^N)$ for every $N\ge1$.
Thus the canonical inclusion
$R[z]\hookrightarrow R\llbracket z \rrbracket$ is continuous for the $D_s$-adic
topology on $R[z]$ and the $z$-adic topology on $R\llbracket z \rrbracket$.
It therefore extends uniquely to a continuous
$R$-algebra homomorphism
\begin{equation}\label{mir:eq-marked-source-completion-map}
 \Phi_s:\mathcal O_s\longrightarrow R\llbracket z \rrbracket,
 \qquad q_0\longmapsto D_s(z),\quad z\longmapsto z.
\end{equation}
The coefficient notation above means
\[
 f_{j,d}:=[z^d]\Phi_s(G_s(z)g_j(\gamma))
 =\sum_{\nu=0}^{\min\{k-1,d\}}
    \sum_{h=0}^{d-\nu}
       n_{j,\nu,h}\,[z^{d-\nu}]D_s(z)^h.
\]
Indeed, $z^\nu D_s(z)^h$ is divisible by $z^{\nu+h}$,
so terms with $\nu+h>d$ contribute zero.
Each displayed polynomial coefficient of $D_s(z)^h$ is
polynomial in $s$. Thus $f_{j,d}$ depends on only finitely
many numerator coefficients and the parameter $s$, and is
a regular coefficient function that can be retained in $Y$.

The contact conditions are finite jet conditions at the fixed
pole bounds. Thus these locally closed loci, the output-center
map, and the finitely many functions $f_{j,d}$ can be retained
simultaneously, by \emph{Retaining data in a finite factor} in
Subsection~\ref{subsec:volume-index-preparations}.
Fix such an enlarged factor $Y$, with $V^\circ$ defined by
the displayed inverse-image identity for this factor, and put
\[
 W=\overline{V^\circ}^{\,\mathrm{red}}\subset Y,
 \qquad W_{s\ne0}:=W\cap Y_{s\ne0},\qquad
 Z\subset Y_0:=Y\times_{\A^1_s}\{0\},
\]
where $Z$ is the irreducible component representing the output
arc component. Denote its generic point by $\eta_Z$, put
$\kappa:=\C(Z)$, and choose an algebraic closure $\bar\kappa$.
Use the resulting geometric point $\Spec\bar\kappa\to Z$
over $\eta_Z$, and restrict to the open-and-closed level
$\delta_{Y/\A^1_s}=\delta$ containing it.

\begin{lemma}[The finite-factor input cycle]
\label{mir:lem-input-cycle-closure}
In this finite-factor presentation the following assertions hold.
\begin{enumerate}[(i),wide=0pt]
\item One has $\dim Z=\delta$, and $V^\circ\to\Gm$ is smooth
of relative dimension $\delta$. Every irreducible component of
$W$ dominates $\A^1_s$ and has dimension $\delta+1$.
Write $\vartheta_{p_i}^{(P_i)}\in H^0(B_{\Omega,P_i})$ for
the canonical pole-stage class of
\eqref{eq:lcy-theta-stage-class} with label $p_i$.
The exterior product of these classes, together with the
complex fundamental class $[\Gm]$, has finite-factor
representative
$(\iota_W)_*[W_{s\ne0}]\in H^{BM}_{2\delta+2}(Y_{s\ne0};\F)$,
as in
Definition~\ref{def:relative-normalized-dual-complexes}, where
\[
 j:V^\circ\hookrightarrow W_{s\ne0},\qquad
 \iota_W:W_{s\ne0}\hookrightarrow Y_{s\ne0}
\]
are respectively open and closed inclusions. Open restriction
on this support gives
\begin{equation}\label{mir:eq-input-cycle-restriction}
 \begin{gathered}
 \operatorname{res}_j:
 H^{BM}_{2\delta+2}(W_{s\ne0};\F)
 \longrightarrow H^{BM}_{2\delta+2}(V^\circ;\F),\\
 \operatorname{res}_j[W_{s\ne0}]=[V^\circ]
       =\sum_{T\in\operatorname{Irr}(V^\circ)}[T].
 \end{gathered}
\end{equation}
Here $\operatorname{Irr}$ denotes the set of irreducible
components, and the brackets denote their complex-oriented
fundamental classes. Every component has coefficient one.
These normalized classes and restriction maps are compatible
with the smooth comparisons in these finite-factor presentations.
\item Let $\nu:W^\nu\to W$ be the normalization of the reduced
scheme $W$. This morphism is finite. If $\eta_Z\in W$, then
$\mathcal O_{W,\eta_Z}$ is reduced and one-dimensional, and
multiplication by $s$ is injective on it. For each
$y\in\nu^{-1}(\eta_Z)$, the ring $\mathcal O_{W^\nu,y}$ is a
discrete valuation ring, $s$ is a nonzero nonunit in this ring,
and its residue field $\kappa(y)$ is a finite extension of
$\kappa$.
\item Let $R$ be a complex algebra and
$\varphi:\Spec R\to Y$ a morphism. Put $s_R=\varphi^*s\in R$. If
$R\to R[s_R^{-1}]$ is injective and the restriction of
$\varphi$ to $\Spec R[s_R^{-1}]$ factors through $V^\circ$,
then $\varphi$ factors through $W$.
\end{enumerate}
\end{lemma}
\begin{proof}
\emph{(i)}
Proposition~\ref{prop:lcy-equidimensional-factors}(ii), applied
to the output pole factor, gives $\dim Z=\delta$.
Over $s\ne0$, separated-source factorization compares the
chosen contact factors with products of the smooth contact
factors in part (i) of that proposition. Smoothness descends
through the common smooth comparison charts. Additivity of
the relative volume indices therefore makes
$V^\circ\to\Gm$ smooth of relative dimension $\delta$.
Each of its irreducible components is open in $V^\circ$,
has dimension $\delta+1$, and has nonempty open image in
$\Gm$. Hence it dominates $\Gm$. Taking reduced closures
proves the assertion about the components of $W$; this is the
horizontal closure, with no components supported at $s=0$.

By \eqref{eq:lcy-theta-stage-class}, each
$\vartheta_{p_i}^{(P_i)}$ is constructed from closed
pushforwards of the fundamental classes of its reduced
finite-factor contact-piece closures. Exterior products and
the complex orientation of $\Gm$ give the fundamental class
of the corresponding product of closures. The chosen
finite-factor presentations compare this product with $W_{s\ne0}$ on
common smooth charts. Smooth maps are open, so inverse image
commutes with closure; smooth pullback also preserves
reducedness. The two pulled-back reduced closures therefore
coincide, since they are closures of the same pulled-back
contact locus.
Lemma~\ref{lem:fundamental-smooth-comparison} preserves these
classes, their closed pushforwards, and their multiplicity-one
coefficients under the smooth comparison maps. This gives
$(\iota_W)_*[W_{s\ne0}]$ with the prescribed normalization.
Since $V^\circ$ is locally closed in $Y_{s\ne0}$, it is open
in its reduced closure $W_{s\ne0}$. Open restriction of the
fundamental class of each closure is the fundamental class
of its contact piece, proving
\eqref{mir:eq-input-cycle-restriction} and its compatibility.

\emph{(ii)}
Part (i) gives $\dim W_\alpha=\delta+1=\dim Z+1$ for every
irreducible component $W_\alpha$ containing $\eta_Z$.
Thus $\mathcal O_{W,\eta_Z}$ is reduced and one-dimensional.
Horizontality makes $s$ belong to none of its minimal primes,
so multiplication by $s$ is injective.
Since $W$ is of finite type over $\C$, its normalization is
finite \cite[Tag~035S]{Stacks}. A point
$y\in\nu^{-1}(\eta_Z)$ has residue field finite over $\kappa$;
on the component of $W^\nu$ containing $y$, its closure has
dimension $\delta$, whereas that component has dimension
$\delta+1$. Hence $\mathcal O_{W^\nu,y}$ is a one-dimensional
Noetherian normal local domain, and therefore a discrete
valuation ring \cite[Tag~00PD]{Stacks}. Horizontality makes
$s$ nonzero there, and $s(y)=0$ makes it a nonunit.

\emph{(iii)}
Every local defining equation of $W$ pulls back to zero
after inverting $s_R$, since $V^\circ\subset W$.
The assumed injectivity, and its localizations on $\Spec R$,
make those pullbacks zero already on $\Spec R$.
This is precisely the factorization through $W$.
\end{proof}

Fix an algebraic closure $\kappa:=\C(Z)\subset\bar\kappa$.
If $\eta_Z\notin W$, put $\mathcal B_Z=\varnothing$. Otherwise define
the set of \emph{geometric collision branches} by
\[
 \mathcal B_Z:=
 \bigl(\nu^{-1}(\eta_Z)\times_{\Spec\kappa}
                         \Spec\bar\kappa\bigr)(\bar\kappa).
\]
Thus $\mathfrak b\in\mathcal B_Z$ consists of a point
$y\in\nu^{-1}(\eta_Z)$ and a $\kappa$-embedding
$\iota:\kappa(y)\hookrightarrow\bar\kappa$. Put
\[
 e_{\mathfrak b}:=\operatorname{ord}_{\mathcal O_{W^\nu,y}}(s)
       \in\Z_{>0}\qquad(\mathfrak b=(y,\iota)).
\]
These orders are defined by
Lemma~\ref{mir:lem-input-cycle-closure}(ii). Since the finite field extension $\kappa(y)/\kappa$ is
separable in characteristic zero, it has exactly
$[\kappa(y):\kappa]$ distinct $\kappa$-embeddings into
$\bar\kappa$. Thus exactly this many geometric branches
$\mathfrak b=(y,\iota)$ lie above $y$. Replacing the geometric normalization fiber
$\nu^{-1}(\eta_Z)\times_{\Spec\kappa}\Spec\bar\kappa$
by its reduction leaves its set $\mathcal B_Z$ of geometric
points unchanged. For $\mathfrak b=(y,\iota)$, however,
$e_{\mathfrak b}$ is the vanishing order of $s$ in the
ambient discrete valuation ring $\mathcal O_{W^\nu,y}$,
before restricting to the fiber. 

Put $K:=\bar\kappa$, with the fixed inclusion
$\kappa\hookrightarrow K$ defining the chosen algebraic
closure. Let $\bar\eta_Z$ denote the composite
\[
 \bar\eta_Z:\Spec K\longrightarrow\Spec\kappa
 \xrightarrow{\eta_Z}Z\hookrightarrow Y,
\]
where the first arrow is induced by this inclusion and
$\eta_Z$ is the generic-point morphism.
In the finite-factor presentation $h:U=Y\times\A^J\to S$ chosen above,
define
\[
 \xi:=(\bar\eta_Z,0)\in U(K),
 \qquad 0=(0)_{j\in J}\in\A^J(K).
\]
Thus the $Y$-coordinate of $\xi$ is $\bar\eta_Z$, and all
its affine-factor coordinates are zero. Since $Z\subset Y_0$,
the point $h(\xi)$ lies over $s=0$. The product $Z\times\A^J$
lies in the pulled-back output arc locus, so the output loop represented by
$h(\xi)$ extends uniquely to an arc
\[
 \alpha:\Spec K\llbracket z \rrbracket\longrightarrow X.
\]
Taking the center of an arc gives a dominant morphism $Z\to X$, so the
 center
\[
 x:=\alpha(0)\in\mathcal V_c(K),
\]
since $\mathcal V_c$ is open and dense in $X$, and $x$ is determined by $\bar\eta_Z$.

For $\mathfrak b=(y,\iota)$, choose a
$\C$-algebra embedding of $\kappa(y)$ into
$\widehat{\mathcal O_{W^\nu,y}}$ whose composite with the
residue-field quotient
$\widehat{\mathcal O_{W^\nu,y}}\twoheadrightarrow\kappa(y)$
is the identity. Together with a uniformizer it gives
\[
 \widehat{\mathcal O_{W^\nu,y}}\cong\kappa(y)\llbracket \tau \rrbracket,
 \qquad \kappa(y)\llbracket \tau \rrbracket\longrightarrow K\llbracket \tau \rrbracket,
\]
where the second map applies $\iota:\kappa(y)\hookrightarrow K$
coefficientwise. The image of the chosen embedding
is a coefficient field; such an embedding exists by the equicharacteristic Cohen
structure theorem
\cite[Tag~032A]{Stacks}. The induced morphism
$\Spec K\llbracket \tau \rrbracket\to W^\nu\to Y$, together with zero
$\A^J$-coordinates, gives a lift of the branch to $U$
with special point $\xi$. The image of $s$ is
a unit times $\tau^{e_{\mathfrak b}}$. Since $K$ is
algebraically closed of characteristic zero, extracting a root
of this unit and changing the uniformizer gives
\[
 s=\tau^{e_{\mathfrak b}}.
\]
Every chosen branch lift specializes to the fixed arc \(\alpha\).



\emph{The source algebra over a complete parameter ring.}
Let $R$ be a commutative $\C$-algebra and let $\tau\in R$
be a non-zero-divisor. Assume that the canonical map
\[
 R\longrightarrow\varprojlim_{n\ge1}R/(\tau^n)
\]
is an isomorphism, so that $R$ is separated and complete
for the $\tau$-adic topology. Fix $s\in\tau R$.
In the applications below, we take
\[
 R=K\llbracket \tau \rrbracket
 \quad\text{or}\quad
 R=K[\lambda]\llbracket \tau \rrbracket,
 \qquad s=\tau^e,\quad e\in\Z_{\ge1}.
\]
Under these hypotheses, the homomorphism $\Phi_s$ of
\eqref{mir:eq-marked-source-completion-map} is an isomorphism
of $R$-algebras:
\begin{equation}\label{mir:eq-complete-marked-source}
 \Phi_s:R\llbracket q_0 \rrbracket[z]/(D_s(z)-q_0)
          \xrightarrow{\ \sim\ }R\llbracket z \rrbracket.
\end{equation}
Indeed, its source and target are $\tau$-adically complete,
separated, and torsion-free; for the source use the basis
$1,z,\ldots,z^{k-1}$ over $R\llbracket q_0 \rrbracket$.
Modulo $\tau$ the map is the monic-division isomorphism
\[
 (R/\tau R)\llbracket q_0 \rrbracket[z]/(z^k-q_0)
 \xrightarrow{\ \sim\ }(R/\tau R)\llbracket z \rrbracket.
\]
Successive lifting modulo $\tau^N$ proves surjectivity of
\eqref{mir:eq-complete-marked-source}; if an element lies in
its kernel, this reduction and cancellation of $\tau$ put it
in every $\tau^N$ times the source, proving injectivity.
Thus every series in $R\llbracket z \rrbracket$ has a unique preimage under
$\Phi_s$ in the completed marked-source algebra $\mathcal O_s$.
The same calculation applies with $R=K\llbracket s \rrbracket$ and parameter $s$.

Recall from Subsection~\ref{mir:subsec-fixed-source-maps} and (\ref{mir:eq-fixed-source-universal-map})
that $\operatorname{ev}:\mathcal H_c\to \mathcal V_c$ is the finite \'etale
evaluation map at the output marking for the fixed source
configuration. 
Its fiber over the center $x=\alpha(0)\in\mathcal V_c(K)$ is
\[
 \mathcal H_{x,c}
 =\mathcal H_c\times_{\mathcal V_c,x}\Spec K.
\]
For $u\in\mathcal H_{x,c}(K)$, also write
$u:\Spec K\to\mathcal H_c$ for its composite with the
projection to $\mathcal H_c$; thus $\operatorname{ev}\circ u=x$.
Specializing the universal map $F$ along $u$ gives a
$K$-morphism
\[
 F_u:\Pj^1_K\longrightarrow M_K,
 \qquad M_K:=M\times_{\C}\Spec K,
 \qquad F_u(\infty)=x.
\]
We write $F(u,-)$ for this specialized map. Its labeled
markings have coordinates $c_1,\ldots,c_k,\infty$, viewed
as $K$-points of $\Pj^1_K$.

\emph{The wedge of a collision branch.}
For $\mathfrak b\in\mathcal B_Z$, choose a lift of its
normalized branch to $U$ with special point $\xi$.
Choose a uniformizer $\tau$ of its parametrizing discrete
valuation ring so that this lift has the form
\[
 \Spec R\longrightarrow U\xrightarrow{h}S,
 \qquad R=K\llbracket \tau \rrbracket,\qquad s=\tau^e,
 \qquad e=e_{\mathfrak b}.
\]
The semistable reduction used below may require a finite
ramified base change. After such a base change, we again
write $R=K\llbracket \tau \rrbracket$ and $s=\tau^e$; the new exponent $e$
is the original exponent multiplied by the ramification
index and need not equal $e_{\mathfrak b}$.
Pulling back the universal punctured-source map along
$\Spec R\to S$ gives a morphism
\[
 \Spec\bigl(\mathcal O_s[q_0^{-1}]\bigr)\longrightarrow X_K,
 \qquad
 \mathcal O_s=R\llbracket q_0 \rrbracket[z]/(D_s(z)-q_0),
 \qquad X_K:=X\times_{\C}\Spec K\hookrightarrow M_K.
\]
Localizing the
isomorphism $\Phi_s$ of
\eqref{mir:eq-complete-marked-source} gives
\[
 \mathcal O_s[q_0^{-1}]
 \xrightarrow{\ \sim\ }
 K\llbracket \tau,z \rrbracket[D_{\tau^e}(z)^{-1}].
\]
Define the \emph{wedge} $\gamma$ associated with this lift as the composite
\[
 \gamma:
 \Spec\bigl(K\llbracket \tau,z \rrbracket[D_{\tau^e}(z)^{-1}]\bigr)
 \xrightarrow{\ \sim\ }
 \Spec\bigl(\mathcal O_s[q_0^{-1}]\bigr)
 \longrightarrow X_K,
\]
where the first arrow is induced by this localized
isomorphism. Here a wedge means a wedge of loops, distinguished
from the wedges of arcs used in Nash-type problems, see for example
\cite{EmbeddedNash,BBLN,BobadillaPePereiraNash,deFernexDocampoNash}.
Its restriction at $\tau=0$ is
$\alpha|_{\Spec K(\!( z )\!)}$, where $\alpha$ is also regarded
as the corresponding $K$-morphism to $X_K$.
We retain the notation $\gamma$ after further ramified
parameter extensions. 

\begin{lemma}[The geometric branch of a collision family]
\label{mir:lem-collision-family-branch}
Keep the fixed finite-factor presentation $h:U=Y\times\A^J\to S$, the point
$\xi$, and the inclusion $\kappa\hookrightarrow K$. Let
\[
 \Lambda=K\quad\text{or}\quad K[\lambda],
 \qquad R=\Lambda\llbracket t \rrbracket,
\]
and let $\xi_\Lambda:\Spec\Lambda\to U$ be the composite of
$\Spec\Lambda\to\Spec K$ with $\xi$.
Suppose that a morphism $\beta:\Spec R\to S$ satisfies
\[
 \beta^*s=t^e,\qquad e\in\Z_{\ge1},\qquad
 \beta|_{\Spec\Lambda}=h\circ\xi_\Lambda,
\]
where $\Spec\Lambda\subset\Spec R$ is defined by $t=0$, and
that $\beta|_{\Spec R[t^{-1}]}$ factors through $\mathcal C^\circ$.
Then the following assertions hold.
\begin{enumerate}[(i),wide=0pt]
\item There is a unique lift
\[
 \widetilde\beta:\Spec R\longrightarrow U,
 \qquad h\circ\widetilde\beta=\beta,
 \qquad \widetilde\beta|_{\Spec\Lambda}=\xi_\Lambda.
\]
Its projection to $Y$ factors through a morphism
$\rho_\beta:\Spec R\to W$. In particular, $\eta_Z\in W$.
\item There is a unique morphism
\[
 \rho_\beta^\nu:\Spec R\longrightarrow W^\nu,
 \qquad \nu\circ\rho_\beta^\nu=\rho_\beta.
\]
Its restriction to $\Spec\Lambda$ is the composite of
$\Spec\Lambda\to\Spec K$ with a unique geometric point
\[
 \mathfrak b:=\operatorname{br}(\beta)=(y,\iota)\in\mathcal B_Z,
 \qquad \iota:\kappa(y)\hookrightarrow K.
\]
If $\Lambda=K$, this lift induces a local injection
\[
 \phi_\beta:\mathcal O_{W^\nu,y}\hookrightarrow K\llbracket t \rrbracket,
 \qquad \phi_\beta(s)=t^e,
\]
with residue-field map $\iota$. Let $\operatorname{ord}_y$ denote
the normalized valuation of the discrete valuation ring
$\mathcal O_{W^\nu,y}$. For any uniformizer $\pi_y$, put
\[
 m_\beta:=\operatorname{ord}_t\bigl(\phi_\beta(\pi_y)\bigr)
 \in\Z_{\ge1}.
\]
This integer is independent of the choice of $\pi_y$. Extending
$\phi_\beta$ to the fraction fields, for every
$g\in\operatorname{Frac}(\mathcal O_{W^\nu,y})^\times$ one has
\begin{equation}\label{mir:eq-collision-family-valuations}
 \operatorname{ord}_t\bigl(\phi_\beta(g)\bigr)
 =m_\beta\operatorname{ord}_y(g),
 \qquad e=e_{\mathfrak b}m_\beta.
\end{equation}
In particular, $e_{\mathfrak b}\mid e$, and $e=1$ implies
$e_{\mathfrak b}=1$.
\item The lifts in (i)--(ii) commute with the $K$-linear
substitutions $\lambda=c$, for $c\in K$, and $t=v^r$, for
$r\in\Z_{\ge1}$ and a new formal parameter $v$, whenever
applicable. These substitutions preserve $\operatorname{br}(\beta)$.
\end{enumerate}
\end{lemma}
\begin{proof}
\emph{Step 1: construct the lift to $U$ and prove (i).}
Formal \'etaleness of $h$ gives compatible unique lifts of
$\xi_\Lambda$ over $R/(t^N)$ for every $N\ge1$
\cite[Tag~04FD]{Stacks}. Affineness of $U$ and $S$, together with $t$-adic
completeness of $R$, gives the unique lift $\widetilde\beta$ over
$R$. The projection to $Y$ factors through $V^\circ$ after
inverting $t$: the source is reduced and the lifted family lies
in $U\times_S\mathcal C^\circ$. Since $R\hookrightarrow R[t^{-1}]$,
Lemma~\ref{mir:lem-input-cycle-closure}(iii) makes the projection
factor through $W$. Its reduction is the constant geometric
point $\bar\eta_Z$, so $\eta_Z\in W$. This proves (i).

\emph{Step 2: construct the normalization lift and the branch in (ii).}
Since $W=\overline{V^\circ}_{\mathrm{red}}$ is closed in the affine
scheme $Y$, it is affine. For each
$g\in\Gamma(W,\mathcal O_W)$ with $g(\eta_Z)\ne0$, the element
$\rho_\beta^*g$ lies in $c+tR$, where $c\in K^\times$ is the
image of $g(\eta_Z)$ under $\kappa\hookrightarrow K$.
Completeness makes this element a unit. Hence $\rho_\beta$
induces a homomorphism
\[
 \mathcal O_{W,\eta_Z}\longrightarrow R,\qquad s\longmapsto t^e.
\]
By Lemma~\ref{mir:lem-input-cycle-closure}(ii), its source is
a reduced one-dimensional local ring. Its prime kernel
$\mathfrak q$ avoids $s$, which belongs to the maximal ideal.
Thus $\mathfrak q$ is not maximal and is therefore minimal. Put
$A_{\mathfrak q}:=\mathcal O_{W,\eta_Z}/\mathfrak q$, and let
$\overline A_{\mathfrak q}$ be its finite integral closure in
$\operatorname{Frac}(A_{\mathfrak q})$.
The injection $A_{\mathfrak q}\hookrightarrow R$ induces an
embedding
\[
 \operatorname{Frac}(A_{\mathfrak q})\hookrightarrow
 \operatorname{Frac}(R).
\]
The image of every element of $\overline A_{\mathfrak q}$ is
integral over $R$, since its monic equation over $A_{\mathfrak q}$
maps to one over $R$. The ring $R$ is a regular domain, hence
integrally closed, so these images belong to $R$. This uniquely
extends $A_{\mathfrak q}\hookrightarrow R$ to
$\overline A_{\mathfrak q}\hookrightarrow R$ and gives the unique
normalization lift $\rho_\beta^\nu$.

Since $\Lambda$ is reduced, the restriction of $\rho_\beta^\nu$
to $t=0$ factors through
\[
 \Spec\Lambda\longrightarrow
 \left(\nu^{-1}(\eta_Z)\times_{\Spec\kappa}\Spec K\right)_{\mathrm{red}}
 \cong\coprod_{\mathfrak b\in\mathcal B_Z}\Spec K.
\]
Connectedness of $\Spec\Lambda$ makes this morphism select a
unique summand, defining $\operatorname{br}(\beta)=(y,\iota)$.
This proves the assertions in (ii) concerning the normalization
lift and its geometric special point.

\emph{Step 3: prove the local injection and valuation assertions in (ii).}
Suppose $\Lambda=K$. The closed point of $\Spec K\llbracket t \rrbracket$ maps
to $y$ under $\rho_\beta^\nu$, so pullback on local rings gives
a local homomorphism
\[
 \phi_\beta:\mathcal O_{W^\nu,y}\longrightarrow K\llbracket t \rrbracket,
 \qquad \phi_\beta(s)=t^e.
\]
Writing $\mathfrak m_y$ for the maximal ideal of
$\mathcal O_{W^\nu,y}$, locality means
\[
 \phi_\beta^{-1}((t))=\mathfrak m_y.
\]
Reduction of $\phi_\beta$ modulo $t$ is the composite
\[
 \mathcal O_{W^\nu,y}\longrightarrow\kappa(y)
 \xrightarrow{\iota}K,
\]
by the definition of the geometric special point in Step~2.
The kernel of $\phi_\beta$ is prime because $K\llbracket t \rrbracket$ is a domain.
Its source is a discrete valuation ring, whose only prime ideals
are $0$ and $\mathfrak m_y$. Since $s\in\mathfrak m_y$ but
$\phi_\beta(s)=t^e\ne0$, this kernel cannot be $\mathfrak m_y$.
Thus $\phi_\beta$ is injective.

Units map to units, whereas a uniformizer $\pi_y$ maps to a
nonzero element of $(t)$. Hence $m_\beta\ge1$, and replacing
$\pi_y$ by a unit multiple does not change $m_\beta$.
Every nonzero element $g$ of
$\operatorname{Frac}(\mathcal O_{W^\nu,y})$ can be written
$g=\pi_y^{\operatorname{ord}_y(g)}v$ with
$v\in\mathcal O_{W^\nu,y}^{\times}$.
Applying the induced embedding of fraction fields and taking
$t$-orders gives the first equality in
\eqref{mir:eq-collision-family-valuations}. Applying it to $s$
gives
\[
 e=\operatorname{ord}_t(t^e)
  =\operatorname{ord}_t\bigl(\phi_\beta(s)\bigr)
  =m_\beta\operatorname{ord}_y(s)
  =m_\beta e_{\mathfrak b}.
\]
Thus $e_{\mathfrak b}\mid e$; if $e=1$, both positive integral
factors are $1$. This completes (ii).

\emph{Step 4: prove the substitution compatibilities in (iii).}
After either substitution in (iii), the resulting family
satisfies the same hypotheses. The base changes of the two
lifts are lifts of this family, so uniqueness identifies them
with the lifts in (i)--(ii). The reduction still selects the
same geometric point $(y,\iota)$, proving (iii).
\end{proof}

For $R$ and $e$ as in the lemma, a wedge
\[
 \gamma:\Spec\bigl(R\llbracket z \rrbracket[D_{t^e}(z)^{-1}]\bigr)\longrightarrow X_K
\]
with the prescribed pole bounds
$G_{t^e}(z)g_j(\gamma)\in R\llbracket z \rrbracket$ for every chosen generator $g_j$
corresponds, via the localized isomorphism
\eqref{mir:eq-complete-marked-source}, to a unique family
$\beta\in S(R)$. Indeed, the inverse of $\Phi_{t^e}$ recovers
its numerator coefficients in the completed marked-source algebra.
When this family satisfies the hypotheses of
Lemma~\ref{mir:lem-collision-family-branch}, define the
\emph{geometric collision branch determined by $\gamma$} to be
\[
 \operatorname{br}(\gamma):=\operatorname{br}(\beta).
\]
For a wedge constructed from a chosen branch
$\mathfrak b=(y,\iota)$ above, uniqueness of the lifts recovers
$\operatorname{br}(\gamma)=\mathfrak b$, also after the ramified
parameter extensions in Lemma~\ref{mir:lem-collision-family-branch}(iii).

In the formulas below, if $a$ and $w$
are morphisms from the same $K$-scheme to $\mathcal H_c$ and
$\Pj^1$ over $\C$, respectively, $F(a,w)$ denotes the
$K$-morphism to $M_K$ whose composite with $M_K\to M$ is
$F\circ(a,w)$. Thus these formulas evaluate the base change
of $F$ to $K$; the functions $g_j$ are likewise pulled back
to $X_K$.

\begin{lemma}[Normal form along a collision branch]
\label{mir:lem-collision-normal-form}
Keep the finite-factor presentation of Lemma~\ref{mir:lem-input-cycle-closure}
and the fixed output point $\xi$, arc $\alpha$, and center $x$
chosen above. Then:
\begin{enumerate}[(i),wide=0pt]
\item For a lifted branch and its wedge $\gamma$ as above, after a
finite ramified extension with $s=\tau^e$, there are a uniquely
determined $u\in\mathcal H_{x,c}(K)$
and a morphism $a:\Spec K\llbracket \tau,z \rrbracket\to \mathcal H_c$, with $a(0,0)=u$,
such that
\begin{equation}\label{mir:eq-collision-normal-form}
 \gamma(\tau,z)=F\bigl(a(\tau,z),z/\tau^e\bigr).
\end{equation}
This equality is between morphisms from the punctured marked
source to $M_K$, with $\gamma$ composed with $X_K\hookrightarrow M_K$.
Here $z/\tau^e$ denotes the morphism to $\Pj^1_K$ given by
$[z:\tau^e]$. It is defined on that source because $c_1=0$
makes $z$ invertible there; at $\tau=0$ it equals $\infty$.
Consequently, the restriction of $a$ to $\tau=0$ satisfies
$a(0,z)=a_0(z)$, where
\[
 a_0:\Spec K\llbracket z \rrbracket\longrightarrow \mathcal H_c,
 \qquad \operatorname{ev}\circ a_0=\alpha,
 \qquad a_0(0)=u
\]
is the unique lift of $\alpha$ through the finite \'etale map
$\operatorname{ev}:\mathcal H_c\to\mathcal V_c$
whose value at $z=0$ is $u$.
\item Conversely, for each $u\in\mathcal H_{x,c}(K)$ this lift
$a_0$ defines a point $\gamma_u\in S(K\llbracket s \rrbracket)$. Its associated
wedge, also denoted by $\gamma_u$, is
\[
 \begin{gathered}
 \gamma_u:\Spec\bigl(K\llbracket s,z \rrbracket[D_s(z)^{-1}]\bigr)\longrightarrow X_K,
 \\
 \gamma_u(s,z)=F\bigl(a_0(z),z/s\bigr),
 \end{gathered}
\]
where the equality is taken in $M_K$.
Here $z/s$ likewise denotes $[z:s]\in\Pj^1_K$. The special
fiber of $\gamma_u$ is the loop underlying $\alpha$, and its
generic fiber has the
prescribed exact-contact labels $p_i$. The lift of $\gamma_u$
through $h$ with special point $\xi$ projects to $W\subset Y$.
Its parameter is $s$, with no ramified change of coefficients.
\end{enumerate}
\end{lemma}

\begin{proof}
\emph{Step 1: resolve the wedge.}
Fix $\mathfrak b\in\mathcal B_Z$ and a lifted branch as in (i),
with parameter $s=\tau^e$. Put $R=K\llbracket \tau \rrbracket$ and
$T=\Spec R\llbracket z \rrbracket$.
The pole conditions for the wedge $\gamma$ give
\[
 g_j(\gamma)=\frac{N_j(\tau,z)}{G_{\tau^e}(z)},
 \qquad N_j\in R\llbracket z \rrbracket.
\]
Composing $\gamma$ with $X_K\hookrightarrow M_K$ gives a
rational map $T\dashrightarrow M_K$. Choose a sequence of point blowups on which the rational map
extends to a morphism and the reduced total transform of
$V(\tau D_{\tau^e}(z))$ has simple normal crossings. We call this
extension a \emph{resolved wedge}. Resolution by
point blowups is available for this excellent regular surface
\cite[Tags 0C5H and 0BIC]{Stacks}. Refine the resolution to dominate the blowup
$\pi:B_e\to T$ of \eqref{mir:eq-geometric-collision-blowup}.
As described in Subsection~\ref{mir:subsec-global-curve-geometry},
its reduced exceptional line has markings $w=c_i$ and an
attachment at $w=\infty$ to the outer disk. The restriction of the resolved wedge to
that disk is $\alpha$, and the attachment maps to $x\in X(K)$.

\emph{Step 2: semistable reduction.}
The finitely many point blowups of Step~1 extend to blowups
of $\Pj^1_R$ at points of its special fiber $\tau=0$: inductively, each
center is a $K$-point in the formal neighborhood already retained.
Semistable reduction dominating this model after a finite
ramified base change is classical
\cite[Proposition~10.4.6]{LiuArithmeticCurves}.
We recall the  local calculation to make the domination
and preservation of the marked sections explicit.
We call a vertical component an irreducible component of the
special fiber $\tau=0$ of the proper source model over $\Spec R$
obtained by the point blowups above; its multiplicity in that fiber is the
order of vanishing of $\tau$ along it.
After the embedded resolution, the base equation is,
\'etale locally and up to a unit, $\tau=x^\mu$ or
$\tau=x^\mu y^\nu$, where $\mu,\nu\in\mathbb Z_{>0}$ are these
multiplicities. Choose an integer
$N_{\mathrm{ss}}>0$ divisible by all these multiplicities and
make the base change $\tau=t^{N_{\mathrm{ss}}}$. Roots of the
units may be taken \'etale locally in characteristic zero.
For a crossing put
\[
 g=\gcd(\mu,\nu),\qquad
 \mu'=\mu/g,\qquad \nu'=\nu/g,\qquad
 q=N_{\mathrm{ss}}/\operatorname{lcm}(\mu,\nu).
\]
The base-changed equation splits locally into factors
$x^{\mu'}y^{\nu'}=\zeta t^{\mu'\nu'q}$ with $\zeta^g=1$.
After absorbing $\zeta$ into a coordinate, the normalization
of each local factor is
\[
 uv=t^q,\qquad x=u^{\nu'},\qquad y=v^{\mu'}.
\]
Indeed, this map is finite, and coprimality of $\mu',\nu'$
recovers $u,v$ as rational monomials in $x,y,t$, so it is
birational; the displayed nodal surface is normal. The charts
coming from $\tau=x^\mu$ normalize to smooth charts.
To resolve $uv=t^q$ for $q>1$, blow up the ideal $(u,t)$.
Its charts $t=ut_1$ and $u=tu_1$, with new local coordinates
$t_1,u_1$, have equations
$v=u^{q-1}t_1^q$ and $u_1v=t^{q-1}$, respectively.
The first is smooth with base equation $t=ut_1$; the second
is resolved by induction on $q$, starting with $uv=t$.
The resulting exceptional curves form a rational chain,
and these base equations show that every vertical component
has multiplicity one. Thus normalization followed by these resolutions gives a
reduced nodal fiber and dominates the base-changed chosen model.

The marked sections are disjoint and meet the resolved special
divisor transversely at its smooth points. If a section meets
a vertical component of multiplicity $\mu$, restricting its
local equation to the section gives
$1=\mu\,\operatorname{ord}_\tau(x|_{\mathrm{section}})$;
hence $\mu=1$. Normalized base change is smooth along these
sections, so they remain disjoint sections. All the modifications
pull back the resolved wedge on its existing neighborhood of the
exceptional locus; no extension to the whole proper source
model is needed. The exceptional part of the resulting nodal
model is a proper tree of rational curves. After this base change,
continue to write $R=K\llbracket \tau \rrbracket$, $T=\Spec R\llbracket z \rrbracket$, and
$s=\tau^e$; the notation $\pi:B_e\to T$ now denotes the
blowup of $(\tau^e,z)$ for this parameter.

\emph{Step 3: extract the marked curve and identify its parameter $u$.}
We construct a logarithmic enhancement to apply
Theorem~\ref{mir:thm-generic-log-maps}(i).
Let $\mathcal Y\to\Pj^1_R$ be the proper marked semistable model
of Step~2. Its pullback along $T\to\Pj^1_R$ gives the resolved
source $\widetilde T$ carrying the wedge. Equip $\mathcal Y$
with the divisorial log structure of its reduced special fiber and
marked sections, and write $\widetilde T^\dagger$ for $\widetilde T$
with the pullback log structure.
The logarithmic base is
$(\Spec R,\{0\})$, where $\{0\}=V(\tau)$; here a pair
denotes the scheme with the divisorial log structure of the
displayed divisor. Write $\widetilde f:\widetilde T\to M_K$ for the
resolved wedge. Together with the structural map to $\Spec R$,
it gives a logarithmic morphism
\[
 \widetilde T^\dagger\longrightarrow
 (M_K,D\times_{\C}\Spec K)
       \times_{\Spec K}(\Spec R,\{0\})
\]
over $(\Spec R,\{0\})$, with $\Spec K$ carrying the trivial
log structure in this product. Indeed, away from the special
fiber and marked sections the map is $\gamma$ with values in
$X_K$, so local boundary equations pull back to monomials times
units in the normal-crossing boundary of $\widetilde T$.

Let $C_{\mathrm{exc}}$ be the reduced exceptional tree of $\mathcal Y$
and $p_{\mathrm{out}}$ its attachment to the outer component.
Since $C_{\mathrm{exc}}$ lies over the closed $K$-point $(\tau,z)=(0,0)$
of $\Pj^1_R$, pullback to $T$ identifies it with the same algebraic
$K$-curve in $\widetilde T$. On $C_{\mathrm{exc}}\setminus\{p_{\mathrm{out}}\}$
the nodal and marked log charts are the restrictions of those on
$\mathcal Y$. In each such chart the unit factors in the boundary
monomials restrict to regular units on the corresponding open of
$C_{\mathrm{exc}}$, so the resolved log morphism restricts there. At $p_{\mathrm{out}}$ use
the smooth marked chart over the standard log point of $\Spec K$.
On a punctured neighborhood of $p_{\mathrm{out}}$ both log structures
are pulled back from that log point. Since
$\widetilde f(p_{\mathrm{out}})=x\in X_K$, the target log structure
is trivial near its image. The two charts and the map thus glue to
give a proper prestable log map on $C_{\mathrm{exc}}$ with output
contact zero.

Recall from the beginning of this subsection that the fixed
input labels are $p_i=\ell_i e_{E_i}$, where
$E_i\subset D$ and $e_{E_i}$ is the primitive integral generator
of its boundary ray. The generic input loop at $z=sc_i$
belongs to the exact-contact locus $\mathcal C_i$; hence its
contact order is $\ell_i$ along $E_i$ and zero along every
other boundary component. Let $\Sigma_i\subset\widetilde T$
be the strict transform of $z=\tau^e c_i$. Resolution and
ramified changes of the base parameter preserve these
horizontal multiplicities. Consequently,
\[
 \widetilde f^*(D_\alpha\times_{\C}\Spec K)
   =\sum_{i:E_i=D_\alpha}\ell_i\Sigma_i+V_\alpha,
 \qquad \operatorname{Supp}(V_\alpha)\subseteq\{\tau=0\},
\]
where $V_\alpha$ is an effective vertical divisor. At each
marked point, the marked section meets the reduced special
fiber transversely at a smooth point. The vertical factor
therefore belongs to the base log structure, whereas the
coefficient of $\Sigma_i$ records the logarithmic contact
at that marking. Thus the input contact remains $p_i$, even
if the special marked point lies on additional boundary
components.

Stabilize this prestable log map by successively contracting rational components
on which it is constant and which at the same time have fewer than three special
points, namely markings or nodes. Direct image of the source
log structure along this contraction makes the induced map logarithmic
\cite[Lemma~3.6]{JohnstonBirational}.
Then \cite[Proposition~1.24]{GSLogGW} gives a basic stable
logarithmic map with the same underlying ordinary stable map:
the preceding log map is its pullback along a morphism of log
points whose underlying map is $\mathrm{id}_{\Spec K}$.
These operations preserve the marking contacts just computed,
including the zero output contact.
After forgetting the map and stabilizing the marked curve,
its source is $(\Pj^1_K;c_1,\ldots,c_k,\infty)$, since the
additional exceptional components contract to the principal
exceptional line of $B_e$.

Denote this stable logarithmic map by $f:C\to M_K$ and
its numerical curve class by $\beta$. The SNC balancing
identity recalled in Subsection~\ref{mir:subsec-fixed-source-maps},
applied to each boundary component $D_\alpha$, gives
\[
 D_\alpha\cdot\beta
 =\deg_C f^*\mathcal O_{M_K}(D_\alpha\times_{\C}\Spec K)
 =\sum_{i:E_i=D_\alpha}\ell_i.
\]
The inputs contribute their specified contacts, the output
contributes zero, and the two contributions at each node cancel
when the componentwise balancing identities are summed.
This line-bundle degree identity remains valid when components
of $C$ map into the boundary. Writing
$A=\sum_\alpha k_\alpha D_\alpha$ and weighting these equalities
by $k_\alpha$ yields
\[
 A\cdot\beta
 =\sum_\alpha k_\alpha\sum_{i:E_i=D_\alpha}\ell_i
 =\sum_i w_A(p_i).
\]
Since $\beta$ is effective, the boundary-degree equalities above
give $\beta\in\Gamma$. In the moduli notation of
Subsection~\ref{mir:subsec-fixed-source-maps}, the map $f$
defines a point
\[
 [f]\in\mathcal M^{\log}_\beta(K),
 \qquad E_\beta([f])=(x,c)\in\mathcal V(K):
\]
its output is $x\in\mathcal V_c(K)$ and its stabilization as
a marked curve represents $c$.
Theorem~\ref{mir:thm-generic-log-maps}(i) therefore gives a
smooth marked source for $f$ and identifies its underlying
ordinary map with $F(u,-)$ for a unique
$u\in\mathcal H_{x,c}(K)$.
This constructs the point $u$ in (i); subsequent steps use
this underlying ordinary map.

\emph{Step 4: descend the resolved wedge to $B_e$.}
By Step~3, the resolved wedge is constant on every exceptional
component other than the principal line: stabilization of the
map contracts only constant components, and its stabilized
source is smooth by Theorem~\ref{mir:thm-generic-log-maps}(i).
The resolved wedge descends to $B_e$. Indeed, its reduced graph
image over $B_e$ is proper with zero-dimensional fibers:
each fiber of the source resolution is a point or a connected
union of additional exceptional components, and the resolved wedge
is constant on such a union. The graph image is
hence finite and birational; both charts in \eqref{mir:eq-geometric-collision-blowup}
are normal, so that image is $B_e$.

\emph{Step 5: construct the morphism to $\mathcal H_c$ on the formal neighborhood.}
Write $E\subset B_e$ for the reduced exceptional line and
$\widehat B_e$ for the completion of $B_e$ along $E$.
Let $\gamma_e:B_e\to M_K$ be the descended wedge and
$w=[z:\tau^e]:B_e\to\Pj^1_K$; use the same notation for
their restrictions to $\widehat B_e$.
Put $D_K:=D\times_{\C}\Spec K$ and
$\mathcal H_{c,K}:=\mathcal H_c\times_{\C}\Spec K$, and
let $F_K$ be the base change of $F$ to $K$.
On $E$ we have $\gamma_e|_E=F_K(u,w)$.
Set
\[
 P^\circ:=\Pj^1_K\setminus\{c_1,\ldots,c_k\},\qquad
 E^\circ=E\setminus\bigcup_{i=1}^k\{w=c_i\},
\]
and let $\widehat B_e^\circ=\widehat B_e|_{E^\circ}$
be the open formal subscheme whose underlying space is $E^\circ$.
The prescribed contacts imply that $F_K$ maps
$\mathcal H_{c,K}\times_K P^\circ$ into $X_K$. Consider
\[
 G:=(F_K,\operatorname{pr}_2):
 \mathcal H_{c,K}\times_K P^\circ\longrightarrow X_K\times_K P^\circ,
\]
where $\operatorname{pr}_2$ is projection to $P^\circ$.
On $X_K$, the logarithmic tangent bundle is the ordinary tangent
bundle. Thus Lemma~\ref{mir:lem-fixed-source-log-differential},
after base change to $K$, makes the differential $d_{\mathcal H_{c,K}}F_K$ an
isomorphism. Under the product decompositions
of the tangent bundles,
\[
 dG=
 \begin{pmatrix}
 d_{\mathcal H_{c,K}}F_K & d_{P^\circ}F_K\\
 0 & \mathrm{id}
 \end{pmatrix}
\]
is therefore an isomorphism. Both schemes are smooth over $K$:
$\mathcal H_{c,K}$ is finite \'etale over
$\mathcal V_c\times_{\C}\Spec K$, a smooth open of $X_K$.
The Jacobian criterion shows
that $G$ is \'etale throughout its source.

The pair $(\gamma_e,w)$ maps $\widehat B_e^\circ$ to
$X_K\times_K P^\circ$, since its restriction to $E^\circ$ does
and these are open conditions. We spell out the formal lifting
step. Let $\mathcal I$ be the ideal sheaf of $E^\circ$ in
$\widehat B_e^\circ$, and define its infinitesimal neighborhoods
\[
 E_m^\circ:=
 \bigl(E^\circ,\mathcal O_{\widehat B_e^\circ}/\mathcal I^{m+1}\bigr),
 \qquad m\geq0.
\]
Then $E_0^\circ=E^\circ$, and $E_m^\circ\hookrightarrow E_{m+1}^\circ$
is defined by the square-zero ideal
$\mathcal I^{m+1}/\mathcal I^{m+2}$.
Since \'etale morphisms are formally \'etale
\cite[Tag~02HM]{Stacks}, the lift
$h_0=(u,w):E^\circ\to\mathcal H_{c,K}\times_K P^\circ$
extends uniquely, successively, to morphisms
\[
 h_m:E_m^\circ\longrightarrow\mathcal H_{c,K}\times_K P^\circ,
 \qquad G\circ h_m=(\gamma_e,w)|_{E_m^\circ},
 \qquad h_{m+1}|_{E_m^\circ}=h_m.
\]
Here $u$ denotes the constant morphism with value $u$.
The lifts on affine opens glue by uniqueness, and the compatible
family defines the unique formal morphism
\[
 \begin{gathered}
 (a^\circ,w):\widehat B_e^\circ
       \longrightarrow\mathcal H_{c,K}\times_K P^\circ,\\
 G\circ(a^\circ,w)=(\gamma_e,w),
 \qquad a^\circ|_{E^\circ}=u.
 \end{gathered}
\]
We next extend its first component
$a^\circ:\widehat B_e^\circ\to\mathcal H_{c,K}$ across the marked points.

Fix $i\in\{1,\ldots,k\}$. Choose an open neighborhood
$V_i\subset M_K$ of $F_K(u,c_i)$ and an \'etale coordinate map
\[
 \chi_i=(x_1,\ldots,x_n):V_i\longrightarrow\A^n_K,
 \qquad n=\dim X,
\]
such that $D_K\cap V_i=\{x_1=0\}$.
Choose an open neighborhood
$W_i\subset\mathcal H_{c,K}\times_K\A^1_{K,w}$ of $(u,c_i)$
with $F_K(W_i)\subset V_i$ on which the fixed contact condition reads
\[
 x_1(F_K(v,w))=(w-c_i)^{\ell_i}\eta_i(v,w),
 \qquad \eta_i\in\Gamma(W_i,\mathcal O_{W_i}^{\times}).
\]
Define
\[
 \begin{aligned}
 Y_i&:=\mathbb G_{m,K}\times_K\A^{n-1}_K\times_K\A^1_K,\\
 \Phi_i:W_i&\longrightarrow Y_i,\\
 (v,w)&\longmapsto
 \bigl(\eta_i(v,w),x_2(F_K(v,w)),\ldots,x_n(F_K(v,w)),w\bigr).
 \end{aligned}
\]
The differential of the first $n$ components of $\Phi_i$ in the
$\mathcal H_{c,K}$ directions is obtained from the logarithmic
differential of $F_K$ by multiplying the
first row by the unit $\eta_i$; its last coordinate is $w$.
The same Jacobian criterion therefore makes $\Phi_i$ \'etale at
$(u,c_i)$. Shrink $W_i$ so that $\Phi_i$ is \'etale throughout $W_i$.

The horizontal multiplicity calculation in Step~3 gives order
$\ell_i$ along the marked section $w=c_i$ and no other
horizontal boundary divisor near that section. There are no
vertical boundary components, since neither the exceptional
line nor the outer disk maps into $D$. Normality of $B_e$ therefore
shows that
\[
 q_i:=x_1(\gamma_e)/(w-c_i)^{\ell_i}
\]
extends as a unit near that section.
Choose an open neighborhood $U_i\subset E$ of $w=c_i$
containing no other marking, and put $\widehat U_i=\widehat B_e|_{U_i}$.
Shrink $U_i$ so that $(u,w)(U_i)\subset W_i$,
$\gamma_e|_{\widehat U_i}$ factors through $V_i$, and $q_i$ is a unit
on $\widehat U_i$. Restricting the defining equality for $q_i$ to
$U_i$ and using $\gamma_e|_E=F_K(u,w)$ gives
\[
 q_i|_{U_i}=\eta_i(u,w),
\]
because $(w-c_i)^{\ell_i}$ is a nonzero divisor on $U_i$.
Consequently the morphism
\[
 \psi_i:=\bigl(q_i,x_2(\gamma_e),\ldots,x_n(\gamma_e),w\bigr):
 \widehat U_i\longrightarrow Y_i
\]
satisfies $\psi_i|_{U_i}=\Phi_i\circ(u,w)$.
Applying the same successive lifting argument using formal \'etaleness to $\Phi_i$ gives
a unique formal morphism
\[
 (a_i,w):\widehat U_i\longrightarrow W_i,
 \qquad \Phi_i\circ(a_i,w)=\psi_i,
 \qquad a_i|_{U_i}=u.
\]
These equalities imply
$\chi_i\circ F_K(a_i,w)=\chi_i\circ\gamma_e$.
The two maps to $V_i$ agree on $U_i$; formal unramifiedness of
the \'etale morphism $\chi_i$ makes them equal on every
infinitesimal neighborhood and hence on $\widehat U_i$. Thus
\[
 a_i:\widehat U_i\longrightarrow\mathcal H_{c,K},
 \qquad F_K(a_i,w)=\gamma_e|_{\widehat U_i}.
\]
Uniqueness of lifting through $G$ makes $a_i$ agree with
$a^\circ$ on $\widehat U_i\cap\widehat B_e^\circ$.
Since $U_i\cap U_j\subset E^\circ$ for $i\ne j$, these
morphisms glue to
\[
 \widehat a:\widehat B_e\longrightarrow\mathcal H_{c,K},
 \qquad F_K(\widehat a,w)=\gamma_e,
 \qquad \widehat a|_E=u.
\]

\emph{Step 6: construct $a:T\to \mathcal H_c$ and complete part (i).}
Shrink $\mathcal H_c$ to an affine neighborhood of $u$. We descend
$\widehat a$ using the equalities
\[
 \pi_*\mathcal O_{B_e}=\mathcal O_T,
 \qquad
 \Gamma(\widehat B_e,\mathcal O_{\widehat B_e})=K\llbracket \tau,z \rrbracket.
\]
The first equality follows from proper birationality over the
normal scheme $T$; the second follows by completing along
$(\tau,z)$ and applying the theorem on formal functions
\cite[Tags 0AY8 and 02OC]{Stacks}. This is the completion
already used: the ideal of the reduced exceptional line is
the radical of $(\tau,z)\mathcal O_{B_e}$, and on this
Noetherian scheme these ideals define equivalent adic
topologies. Since $\mathcal H_{c,K}$ is affine after the shrinking of $\mathcal H_c$, $\widehat a$
is determined by its $K$-algebra homomorphism
\[
 \Gamma(\mathcal H_{c,K},\mathcal O_{\mathcal H_{c,K}})
       \longrightarrow K\llbracket \tau,z \rrbracket.
\]
This defines a $K$-morphism $a_K:T\to\mathcal H_{c,K}$ with
$(a_K\circ\pi)|_{\widehat B_e}=\widehat a$.
Composing $a_K$ with the projection
$\mathcal H_{c,K}\to\mathcal H_c$ gives the required morphism
$a:T\to\mathcal H_c$.
The descended wedge and $F(a\circ\pi,[z:\tau^e])$
agree on $\widehat B_e$. At a closed point of the exceptional
line, the local ring of $B_e$ injects into its completion along
that line by the Krull intersection theorem. Hence the two maps
agree at the generic point of the integral scheme $B_e$;
separatedness of $M_K$ gives equality on $B_e$.
Restricting to the original punctured marked source proves
\eqref{mir:eq-collision-normal-form}. Reduction on the outer
disk gives $\operatorname{ev}(a(0,z))=\alpha(z)$ and
$a(0,0)=u$. Formal \'etaleness of $\operatorname{ev}$ gives
the asserted uniqueness of $a_0$, hence $a(0,z)=a_0(z)$.
The restriction to the exceptional line is the marked map
$F(u,-)$; the unique point of $\mathcal H_{x,c}(K)$ representing
this map was identified in Step~3, giving the asserted uniqueness of $u$. Together with
\eqref{mir:eq-collision-normal-form}, this proves (i).

\emph{Step 7: construct the converse wedge and complete part (ii).}
For (ii), return to the original finite \'etale family
$\operatorname{ev}:\mathcal H_c\to \mathcal V_c$ and fix any $u\in \mathcal H_{x,c}(K)$.
Lift $\alpha$ uniquely to $a_0$ and shrink
to an affine neighborhood of $u$. Each $g_j\circ F$ has poles
bounded by $\sum_i m_i\{c_i\}$ and is regular at infinity.
Multiplication by $\prod_i(w-c_i)^{m_i}$ therefore gives a
section of $\mathcal O_{\Pj^1}(b)$, whose basis
$1,w,\ldots,w^b$ supplies regular coefficient functions $P_{jd}$
on this neighborhood, with
\[
 g_j(F(v,w))=
 \frac{\sum_{d=0}^{b}P_{jd}(v)w^d}
      {\prod_i(w-c_i)^{m_i}}.
\]
For a morphism $a:\Spec K\llbracket s,z \rrbracket\to \mathcal H_c$ whose restriction at
$s=0$ is $a_0$, substitution
gives
\begin{equation}\label{mir:eq-collision-normal-form-coordinates}
 g_j\bigl(F(a(s,z),z/s)\bigr)=
 \frac{\sum_{d=0}^{b}P_{jd}(a(s,z))z^d s^{b-d}}
      {G_s(z)}.
\end{equation}
Take $a(s,z)=a_0(z)$. By
\eqref{mir:eq-complete-marked-source} for $R=K\llbracket s \rrbracket$, the
numerator in $R\llbracket z \rrbracket$ has a unique preimage under $\Phi_s$ in
\[
 \mathcal O_s=R\llbracket q_0 \rrbracket[z]/(D_s(z)-q_0).
\]
The denominator is invertible in
$\mathcal K_s=\mathcal O_s[D_s^{-1}]$. The identities of
$X$ hold because they hold for $F$, and the factor
$(w-c_i)^{\ell_i}$ times a unit gives exactly contact $p_i$.
Setting $s=0$ in the punctured source gives
$g_j(\operatorname{ev}(a_0(z)))=g_j(\alpha(z))$.
These equations therefore define $\gamma_u\in S(K\llbracket s \rrbracket)$
with the generic and special fibers asserted in (ii).
The coefficient ring is $K\llbracket s \rrbracket$ itself, so this construction
requires no ramified change of the parameter $s$.

Lemma~\ref{mir:lem-collision-family-branch}(i), applied with
$R=K\llbracket s \rrbracket$ and $t=s$, gives the unique lift of $\gamma_u$
through $h$ with special point $\xi$ and makes its projection
factor through $W$. This completes (ii). The same lemma shows
that existence of $\gamma_u$ forces $\eta_Z\in W$; thus
$\eta_Z\notin W$ implies $\mathcal H_{x,c}=\varnothing$.
\end{proof}

\subsection{Branch multiplicities and the trace comparison}
\label{mir:subsec-collision-multiplicities}

\begin{lemma}[Branches and multiplicities of the fixed-source collision]
\label{mir:lem-collision-branches}
In the setup of Subsection~\ref{mir:subsec-collision-normal-forms},
extraction of the marked curve induces a bijection
\[
 \mathcal B_Z\xrightarrow{\ \sim\ }\mathcal H_{x,c}(K),
 \qquad\mathfrak b\longmapsto f_{\mathfrak b},
\]
and $e_{\mathfrak b}=1$ for every $\mathfrak b\in\mathcal B_Z$.
Consequently the normalized ray specialization of the exterior
product $\boxtimes_i\vartheta_{p_i}^{(P_i)}$, computed on the
support $W$ of Lemma~\ref{mir:lem-input-cycle-closure}(i), has coefficient
$\#\mathcal H_{x,c}(K)$ along $[Z]$. This coefficient
is unchanged by compatible finite-factor presentations and their smooth
comparisons.
\end{lemma}
\begin{proof}
If $\mathcal B_Z=\varnothing$, then $\eta_Z\notin W$, so the
specialization supported on $W$ has coefficient zero along $[Z]$.
Lemma~\ref{mir:lem-collision-normal-form}(ii) also gives
$\mathcal H_{x,c}=\varnothing$, and all assertions follow.
Assume henceforth that $\mathcal B_Z\ne\varnothing$.

\emph{Step 1: extraction map is well defined.} We must check that every lift of a fixed geometric branch
$\mathfrak b\in\mathcal B_Z$ yields the same extracted point
$u\in\mathcal H_{x,c}(K)$.
Fix $\mathfrak b=(y,\iota)\in\mathcal B_Z$ and choose a lifted
branch. Lemma~\ref{mir:lem-collision-normal-form}(i), after a
possible ramified parameter extension, writes its wedge as
\[
 \gamma(\tau,z)=F(a(\tau,z),z/\tau^e),
 \qquad a(0,0)=u\in\mathcal H_{x,c}(K),
 \qquad f=F(u,-).
\]
Let $\rho:\Spec K\llbracket \tau \rrbracket\to W\subset Y$ be its projection,
so $\rho^*s=\tau^e$.
The functions $f_{j,d}$ retained in
Subsection~\ref{mir:subsec-collision-normal-forms} descend to
$f^Y_{j,d}$ and satisfy
\[
 \rho^*f^Y_{j,d}
 =[z^d]\bigl(G_{\tau^e}(z)g_j(\gamma(\tau,z))\bigr),
 \qquad 0\le d\le b.
\]
Let $P_{jd}$ be the regular numerator coefficients on a neighborhood
of $u$ used in \eqref{mir:eq-collision-normal-form-coordinates}.
That formula gives
\[
 \frac{\rho^*f^Y_{j,d}}{\tau^{e(b-d)}}
 =\sum_{\ell=0}^{d}\tau^{e(d-\ell)}
       [z^{d-\ell}]P_{j\ell}(a(\tau,z)).
\]
View $f^Y_{j,d}|_W$ and $s$ as functions on the component of
$W^\nu$ containing $y$. By
Lemma~\ref{mir:lem-collision-family-branch}(ii)--(iii), $\rho$
induces a local injection
$\mathcal O_{W^\nu,y}\hookrightarrow K\llbracket \tau \rrbracket$.
The valuation comparison
\eqref{mir:eq-collision-family-valuations} and the displayed
equality imply $f^Y_{j,d}/s^{b-d}\in\mathcal O_{W^\nu,y}$.
Reducing that equality modulo $\tau$ shows that the image of
this quotient under
\[
 \mathcal O_{W^\nu,y}\longrightarrow\kappa(y)
 \xrightarrow{\iota}K
\]
is $P_{jd}(u)$. Consequently
\[
 g_j(f(w))=
 \frac{\sum_{d=0}^{b}\overline{f^Y_{j,d}/s^{b-d}}\,w^d}
      {\prod_i(w-c_i)^{m_i}},
\]
where the bar denotes the residue class at $y$ followed by $\iota$.
These residue classes depend only on $\mathfrak b=(y,\iota)$,
independently of the auxiliary uniformizer used to
parametrize the branch and of the chosen lift to
$U=Y\times\A^J$, including its $\A^J$-coordinates.
The affine generators $g_j$ determine $f$ away from the markings,
and separatedness determines its extension to $\Pj^1_K$.
Thus extraction defines the map
$\mathcal B_Z\to\mathcal H_{x,c}(K)$ asserted in the statement.

\emph{Step 2: interpolate to the converse family without changing the branch.}
This auxiliary step will prove the branch identification needed
for injectivity and for the multiplicity assertion.
Fix $\mathfrak b$ and retain $a,u,e$ from Step~1. Put
$a_0(z)=a(0,z)$, the unique lift of $\alpha$ through
$\operatorname{ev}$ with $a_0(0)=u$.
Choose an affine open neighborhood $H\subset\mathcal H_c$ of the image
of $u$, contained in the coefficient neighborhood of Step~1,
and an \'etale coordinate map
$q:H\to\A^n$, where $n=\dim X$.
For an indeterminate $\lambda$, set
\[
 q(a_\lambda)=q(a_0)+\lambda\bigl(q(a)-q(a_0)\bigr).
\]
Formal \'etaleness uniquely defines $a_\lambda$ over
$K[\lambda]\llbracket z \rrbracket\llbracket \tau \rrbracket$, with reduction $a_0(z)$ modulo $\tau$.
By uniqueness, its endpoint restrictions satisfy
\[
 a_\lambda|_{\lambda=0}=a_0,
 \qquad a_\lambda|_{\lambda=1}=a.
\]
Insert it in \eqref{mir:eq-collision-normal-form-coordinates}
with $s=\tau^e$. Put
\[
 R:=K[\lambda]\llbracket \tau \rrbracket,\qquad
 \mathcal O_R:=R\llbracket q_0 \rrbracket[z]/(D_{\tau^e}(z)-q_0).
\]
The inverse of $\Phi_s$ in \eqref{mir:eq-complete-marked-source}
sends each resulting numerator to $\mathcal O_R$; inversion of
$G_{\tau^e}$ gives an element of
$\mathcal O_R[D_{\tau^e}^{-1}]$ with the stated pole bounds.
The result is an $R$-point of $S$ with constant special arc
$\alpha$. Over $\tau\ne0$ it has exactly the prescribed input contacts,
by the local contact factorizations in Step~5 of the proof of
Lemma~\ref{mir:lem-collision-normal-form}: after substituting
$z=\tau^e w$, their unit factors reduce modulo $\tau$ to
$\eta_i(u,w)$.

Apply Lemma~\ref{mir:lem-collision-family-branch} to the resulting
$R$-point of $S$. Its two endpoint families determine the same
geometric branch. At $\lambda=1$, the identity
$a_\lambda|_{\lambda=1}=a$ and
Lemma~\ref{mir:lem-collision-normal-form}(i) identify the wedge
with the original $\gamma$, whose branch is $\mathfrak b$.
At $\lambda=0$, the wedge $F(a_0(z),z/\tau^e)$ is the base change
under $s=\tau^e$ of the family
$\gamma_u(s,z)=F(a_0(z),z/s)$ over $K\llbracket s \rrbracket$ supplied by
Lemma~\ref{mir:lem-collision-normal-form}(ii).
Compatibility with ramified substitution in
Lemma~\ref{mir:lem-collision-family-branch}(iii) therefore gives
\[
 \operatorname{br}(\gamma_u)
 =\operatorname{br}\bigl(F(a_0(z),z/\tau^e)\bigr)
 =\operatorname{br}(\gamma)=\mathfrak b.
\]

\emph{Step 3: prove the asserted bijection.}
Lemmas~\ref{mir:lem-collision-normal-form}(ii) and
\ref{mir:lem-collision-family-branch} define a map
\[
 \sigma:\mathcal H_{x,c}(K)\longrightarrow\mathcal B_Z,
 \qquad u\longmapsto\operatorname{br}(\gamma_u).
\]
If $\mathfrak b$ extracts to $u$, Step~2 gives
$\sigma(u)=\mathfrak b$. Conversely, fix $u$ and put
$\mathfrak b=(y,\iota):=\sigma(u)$. Let
$\rho_u:\Spec K\llbracket s \rrbracket\to W$ be the projection of the lift of
$\gamma_u$ supplied by Lemma~\ref{mir:lem-collision-family-branch}(i).
The local injection in part (ii) of that lemma identifies
residue classes at $(y,\iota)$ with reduction modulo $s$.
Equation~\eqref{mir:eq-collision-normal-form-coordinates} gives
\[
 \frac{\rho_u^*f^Y_{j,d}}{s^{b-d}}
 =\sum_{\ell=0}^{d}s^{d-\ell}
       [z^{d-\ell}]P_{j\ell}(a_0(z))
 \equiv P_{jd}(u)\pmod{s},
 \qquad 0\le d\le b.
\]
The intrinsic coefficient reconstruction in Step~1 therefore
identifies the curve extracted from $\mathfrak b$ with $F(u,-)$.
Thus $\sigma$ is inverse to extraction, proving the bijection.

\emph{Step 4: prove $e_{\mathfrak b}=1$ for every geometric branch.}
For $\mathfrak b=(y,\iota)$ with extracted point $u$, Step~2 gives
$\operatorname{br}(\gamma_u)=\mathfrak b$.
The family $\gamma_u$ is defined over $K\llbracket s \rrbracket$ with collision
parameter $s$. Thus the case $t=s$ and $e=1$ of
Lemma~\ref{mir:lem-collision-family-branch}(ii) gives
$e_{\mathfrak b}=1$, proving the multiplicity assertion.

\emph{Step 5: compute the asserted specialization coefficient.}
Lemma~\ref{mir:lem-input-cycle-closure}(i) gives coefficient one
on every input component. By Step~4, $s$ is a uniformizer in
$\mathcal O_{W^\nu,y}$ for every $y\in\nu^{-1}(\eta_Z)$.
At a general smooth point of the corresponding boundary component,
$s$ is therefore a transverse local parameter.
The inward ray on a transverse complex disk contributes $+1$
to its localization boundary, by
\eqref{eq:specialization-half-ray-boundary} of
Lemma~\ref{lp:lem:bm-collision-specialization}, with the ray
placed first as in Definition~\ref{def:normalized-ray-specialization}.
The closure of $y$ maps generically finitely to $Z$ with degree
$[\kappa(y):\kappa]$. Proper trace under the finite normalization,
equivalently the order formula \cite[Tag~02MJ]{Stacks}, gives
the coefficient
\[
 \sum_{y\in\nu^{-1}(\eta_Z)}[\kappa(y):\kappa]
 =\#\mathcal B_Z=\#\mathcal H_{x,c}(K)
\]
along $[Z]$. The first equality counts geometric branches
above each $y$, and the second is the bijection proved in Step~3.

\emph{Step 6: prove independence of the finite-factor presentation and its comparisons.}
The smooth comparison maps are over $\A^1_s$. Their normalized
pullbacks commute with specialization by
Lemma~\ref{lem:specialization-arrow-descent}(ii) and preserve
fundamental classes by
Lemma~\ref{lem:fundamental-smooth-comparison}.
In particular, for a finite \'etale comparison each component
of the inverse image has coefficient one under pullback;
no covering degree is inserted. Thus the coefficient computed
in Step~5 is unchanged by compatible finite-factor presentations and
their smooth comparisons, proving the final assertion.
\end{proof}

\begin{proof}[Proof of Theorem~\ref{mir:thm-low-trace-comparison}]
Lemma~\ref{mir:lem-modification-theta-trace} and
Theorem~\ref{mir:thm-generic-log-maps}(ii) allow us to use the
refined pair. Retain the fixed configuration $\mathbf c$
chosen in Subsection~\ref{mir:subsec-fixed-source-maps}.

Collapse the chosen configuration $z=sc_i$ along $s>0$,
with the output coordinate $z$, using the normalized ray
specialization of Definition~\ref{def:normalized-ray-specialization}.
For $k=2$, the configuration is $(0,s)$ and the input parameters
are $z,z-s$. Lemma~\ref{lem:one-sided-binary-restriction}
identifies this operation with $\overline m_{2;P_1,P_2}$,
including these endpoint identifications.
For $k=3$, put $r=s$ and $t=(c_3-1)s$, so that
\[
 (sc_1,sc_2,sc_3)=(0,r,r+t),\qquad c_3-1>0.
\]
This is the interior ray of the quadrant in
Corollary~\ref{lem:three-point-collision}, with its two positive
gap coefficients $1$ and $c_3-1$. That corollary identifies
its normalized specialization with $c_{P_1,P_2,P_3}$.
Lemma~\ref{lem:collision-composition} then identifies this
morphism with either iterated binary product, with the same
input parameters $z-sc_i$ and output parameter $z$.
Thus in both cases the induced operation on cohomology is the
product used in the loop trace.

By Remark~\ref{rem:arc-trace-finite-factor}, the left side of
\eqref{mir:eq-low-trace-comparison} is the coefficient along
the output arc component $[Z]$ of this specialization.
Lemma~\ref{lem:specialization-arrow-descent}(i)--(ii) identifies
the normalized localization morphism with its finite-factor
representative, compatibly with the smooth comparison maps in
the chosen finite-factor presentations.
Applying Lemma~\ref{mir:lem-collision-branches} to the canonical
pole-stage input classes $\vartheta_{p_i}^{(P_i)}$ therefore gives
\[
 \epsilon(\vartheta_{p_1}\cdots\vartheta_{p_k})=\#\mathcal H_{x,c}(K).
\]
Theorem~\ref{mir:thm-generic-log-maps}(ii) identifies this cardinality
with $\sum_{\beta\in Q}L_\beta(p_1,\ldots,p_k)$.
Equation~\eqref{mir:eq-gs-weak-frobenius} identifies the latter sum
with $\epsilon_{\mathrm{GS}}
(\vartheta_{p_1}^{\mathrm{GS}}\cdots\vartheta_{p_k}^{\mathrm{GS}})$,
which proves \eqref{mir:eq-low-trace-comparison}.
\end{proof}
\subsection{Proof of the Gross--Siebert comparison theorem}
\label{mir:subsec-proof-gs-comparison}

We return to the original pair $(M,D)$ and coefficient monoid $Q$.

\begin{proof}[Proof of Theorem~\ref{mir:thm-intro-gs}]
By hypothesis, $D$ has maximal boundary and
$\check X=\Spec(R_X\otimes_\F\C)$ is smooth.
By Theorem~\ref{mir:thm-trace-nondegeneracy}, the two-input trace
$T_2(a,b)=\epsilon(ab)$ of Definition~\ref{mir:def-arc-trace}
is nondegenerate. Both contact bases are bases of direct-sum vector
spaces. Hence
\begin{equation}\label{mir:eq-gs-ring-isomorphism}
 \Phi_{\mathrm{GS}}:R_X\xrightarrow{\ \sim\ }R^{\mathrm{GS},1}_{M,D},
 \qquad \Phi_{\mathrm{GS}}(\vartheta_p)=\vartheta_p^{\mathrm{GS}}
\end{equation}
is a linear isomorphism preserving the unit and arc coefficient:
$\Phi_{\mathrm{GS}}(1)=1$ and $\epsilon_{\mathrm{GS}}\Phi_{\mathrm{GS}}=\epsilon$.

Theorem~\ref{mir:thm-low-trace-comparison} gives
\[
 \epsilon(a_1\cdots a_d)
 =\epsilon_{\mathrm{GS}}\bigl(\Phi_{\mathrm{GS}}(a_1)\cdots\Phi_{\mathrm{GS}}(a_d)\bigr),
 \qquad d=2,3.
\]
First this holds for nonzero contact labels. If an input is
$\vartheta_0=1$, remove it; the remaining cases are the same
identities in lower arity or preservation of the arc coefficient.
Multilinearity then gives the displayed identities for all inputs. Since $\Phi_{\mathrm{GS}}$ is bijective and $T_2$ is
nondegenerate, the case $d=2$ also proves that the pairing
$(a,b)\mapsto\epsilon_{\mathrm{GS}}(ab)$ on
$R^{\mathrm{GS},1}_{M,D}$ is nondegenerate.
Transport the Gross--Siebert product by $\Phi_{\mathrm{GS}}$. It and loop fusion
have the same two-input and three-input traces. By
Proposition~\ref{mir:prop-trace-reconstruction} they coincide, so
$\Phi_{\mathrm{GS}}$ is an algebra isomorphism. This uses injectivity
$R_X\to R_X^\vee=\operatorname{Hom}_{\F}(R_X,\F)$, $a\mapsto T_2(a,-)$; it does not require
$R_X$ to be finite-dimensional.

For arbitrary $d\ge2$, multiplicativity and trace preservation now
give
\begin{equation}\label{mir:eq-log-trace}
 \epsilon(\vartheta_{p_1}\cdots\vartheta_{p_d})
 =\epsilon_{\mathrm{GS}}
   (\vartheta_{p_1}^{\mathrm{GS}}\cdots\vartheta_{p_d}^{\mathrm{GS}})
 =\sum_\beta L_\beta(p_1,\ldots,p_d),
\end{equation}
where the last equality is \eqref{mir:eq-gs-weak-frobenius} and
the finiteness was established after
\eqref{mir:eq-log-descendant} using Remark~\ref{rmkBbb}.

The loop theta classes and their trace are independent of the
boundary polarization by Theorem~\ref{thm:intro-theta-basis}.
The Gross--Siebert coefficients are likewise independent of the
ample divisor used to bound their finite sums
(Theorem~\ref{mir:lem-gs-polynomial}). The map $\Phi_{\mathrm{GS}}$ is the unique
linear map identifying these labeled bases; thus the unfiltered
algebra isomorphism with trace is independent of $A$.
\end{proof}

\section{Loop-Floer Conjecture for surfaces}
\label{sec:surface-loop-floer}
We prove Theorem~\ref{thm:surface-loop-floer} by first computing
degree-zero symplectic cohomology from the surface mirror and then
comparing loop and symplectic cohomology in all degrees.

\subsection{Degree-zero symplectic cohomology and the Gross--Siebert ring}
\label{subsec:surface-sh-zero-gs}
Let $S=R^{\mathrm{GS},1}_{M,D}\otimes_{\F}\C$, where the superscript
$1$ denotes the specialization $q^\beta\mapsto1$ of
\eqref{mir:eq-gs-specialization}. We show here that a collection of results from the mirror symmetry literature imply that if $\Spec S$ is smooth,
then $SH^0_\Omega(X;\C)\cong S$ as unital $\C$-algebras. 

\emph{The mirror at the identity character.}
Recall that a \emph{Looijenga pair} is a smooth projective surface $Y$
together with a connected, reduced divisor $B\in|-K_Y|$ whose
singular locus is nonempty and consists of nodes. Then $Y$ is
rational, and $B$ is either a rational curve with one node or a
cycle of smooth rational curves, cf.
\cite[Definition~1.1]{GHKModuli} and the discussion following it.
Here $D\in|-K_M|$ since $\operatorname{div}_M(\Omega)=-D$.
Choose an ample divisor $A=\sum_i k_iD_i$ with $k_i>0$ as in
Subsection~\ref{subsec:intro-positive-lcy}. Its support $D$ is
connected. Since $D$ is SNC, its components are smooth and its
singularities are nodes. Maximal boundary in dimension two means
that two components meet. Thus $(M,D)$ is a Looijenga pair with
$D$ a cycle of smooth rational curves. Since $A^2>0$, the pair
is positive in the sense of
\cite[Lemma~6.9(1.1) and Definition~6.10]{GHKSurfaceMirror}.

After blowups at nodes of $D$, the resulting pair admits a toric model. It is
obtained from a smooth toric surface by blowups at smooth points
of the successive boundaries \cite[Proposition~2.3]{LaiZhou}.
The node blowups preserve $X$, which therefore contains the dense
open torus of the toric model.

Use the numerical curve lattice $\Lambda$ and the monoid $Q$ of
Subsection~\ref{mir:subsec-gs-specialization}. Let $E_{\rm eff}\subset\Lambda$
be the monoid of classes of effective $1$-cycles on $M$.
Gross--Hacking--Keel use the lattice points $P$ of the real cone
generated by $E_{\rm eff}$ \cite[\S0.1]{GHKSurfaceMirror}.
Thus $E_{\rm eff}\subseteq P\subseteq Q\subseteq\Lambda$. Since classes
of curves generate $\Lambda$ as a group, all three monoids have the
same group completion. We compare the mirror algebras over the torus
\[
 T=\Spec\C[\Lambda]\cong\operatorname{Pic}(M)\otimes\Gm,
\]
where the identification uses the perfect intersection pairing on
the rational surface $M$. Since $X$ contains a dense torus and $D$
supports the nef divisor $A$, Johnston's comparison
\cite[Corollary~1.2]{JohnstonComparison} identifies the Gross--Siebert
and Keel--Yu theta multiplications after extension to $\C[\Lambda]$.
The surface comparison \cite[Proposition~7.2]{LaiZhou} identifies
the latter with the Gross--Hacking--Keel multiplication over the same ring.

Let $\mathcal Y\to T$ be the restriction of the Gross--Hacking--Keel
mirror family, and put
$B=R^{\mathrm{GS}}_{M,D}\otimes_{\F}\C$, with its $\C[Q]$-algebra
structure from \eqref{mir:eq-gs-product}. Thus $\mathcal Y$ is affine
with coordinate ring $B\otimes_{\C[Q]}\C[\Lambda]$.
At the identity $e:\C[\Lambda]\to\C$, $q^\beta\mapsto1$,
its fiber $\mathcal Y_e$ has coordinate ring
\[
 (B\otimes_{\C[Q]}\C[\Lambda])\otimes_{\C[\Lambda],e}\C
 \cong B\otimes_{\C[Q],\,e|_{\C[Q]}}\C=S
\]
by \eqref{mir:eq-gs-specialization}. Hence $\mathcal Y_e\cong\Spec S$.

By \cite[Theorem~1.8]{LaiZhou}, $\mathcal Y_e$ has a
compactification $(\overline{\mathcal Y}_e,D_e)$ in the family of
generalized Looijenga pairs of the deformation type of $(M,D)$.
The surface $\overline{\mathcal Y}_e$ is smooth along $D_e$ and
on its interior $\mathcal Y_e=\Spec S$, the latter by the smoothness
assumed in this subsection.
Write $D_e=\sum_j D_{e,j}$ and choose a cyclic orientation to
identify $\operatorname{Pic}^0(D_e)$ with $\C^*$. The ordinary
period of the compactified fiber is the restriction character
\[
 \begin{aligned}
 D_e^\perp&=\{L\in\operatorname{Pic}(\overline{\mathcal Y}_e):
                      L\cdot D_{e,j}=0\text{ for every }j\},\\
 \phi_e:D_e^\perp&\longrightarrow\operatorname{Pic}^0(D_e),
                   \qquad L\longmapsto L|_{D_e}.
 \end{aligned}
\]
By \cite[Proposition~6.23]{LaiZhou}, the ordinary period of a fiber
is the image of its parameter under a homomorphism of tori. Thus
$\phi_e(L)=1$ for every $L\in D_e^\perp$.
Torelli in this deformation type \cite[Theorem~2.10]{HackingKeating}
identifies $(\overline{\mathcal Y}_e,D_e)$ with the distinguished
pair of trivial ordinary period.

\emph{The wrapped category.}
Let $\mathcal W$ denote the wrapped Fukaya category over $\C$ of the
standard Weinstein completion of $X$, with graded spin branes.
We use the logarithmic grading determined by $\Omega|_X$.
Any two forms on $M$ with divisor $-D$ differ by a nonzero constant
and hence give homotopic trivializations of the canonical bundle of $X$.
Varolg\"une\c{s}'s theorem
\cite[Theorem~1.12, Case~1]{Varolgunes} identifies this completion,
by a grading-preserving strongly exact symplectomorphism, with the
almost toric model specified by the fan and boundary blowup
multiplicities of the toric model above. Its algebraic mirror uses
the same fan and multiplicities, with the first blowup on each
affected boundary component at $-1$ in its open torus orbit.
Each further blowup on that component is at the intersection of its
strict transform with the last exceptional curve
\cite[\S1.1.2]{Varolgunes}. These centers give the distinguished
pair in the deformation type of the corner refinement of $(M,D)$.
Varying the centers deforms both the refinement and its corner
blowdown \cite[Definition~2.2 and \S\S2.1.1, 2.2]{HackingKeating}.

To compare with $(\overline{\mathcal Y}_e,D_e)$, perform the corresponding
blowups at boundary nodes on this pair. For one such blowup
$\pi:(Y',D')\to(Y,D)$, let $C$ be the exceptional component of $D'$.
Every class in $\operatorname{Pic}(Y')$ has the form $\pi^*L+aC$,
with $L\in\operatorname{Pic}(Y)$ and $a\in\Z$. Since
$(\pi^*L+aC)\cdot C=-a$, the projection formula gives
$(D')^\perp=\pi^*(D^\perp)$. For $L\in D^\perp$ one has
$(\pi^*L)|_{D'}=(\pi|_{D'})^*(L|_D)$, so triviality of the
ordinary period is preserved. Torelli therefore identifies the
refinement of $(\overline{\mathcal Y}_e,D_e)$ with the distinguished
pair constructed at $-1$. Each corner blowup restricts to an
isomorphism $Y'\setminus D'\cong Y\setminus D$, so the interior
of this distinguished pair is isomorphic to $\mathcal Y_e$.
The wrapped equivalence of Hacking--Keating
\cite[Theorems~1.1 and~1.2]{HackingKeating}, in the form
\cite[Theorem~1.10]{Varolgunes}, now gives a $\C$-linear exact
equivalence
\[
 H^0(\operatorname{Perf}\mathcal W)
       \simeq D^b\operatorname{Coh}(\mathcal Y_e).
\]
Here $\operatorname{Perf}\mathcal W$ is the category of perfect
modules over $\mathcal W$ and the right side is the bounded derived
category of coherent sheaves.

\emph{The cohomology rings.}
Since $\mathcal Y_e=\Spec S$ is smooth affine of finite type,
$\mathcal O_{\mathcal Y_e}$ generates $D^b\operatorname{Coh}(\mathcal Y_e)$
under shifts, cones, and direct summands. In a differential graded
model of $\operatorname{Perf}\mathcal W$, choose the corresponding
object $E$ and put $A_E=\operatorname{End}(E)$. Then $E$ is a
split-generator, and the equivalence gives
\[
 H^j(A_E)\cong\operatorname{Ext}^j_S(S,S)=
 \begin{cases}S,&j=0,\\0,&j\ne0,\end{cases}
\]
with composition identified with multiplication in $S$.
Let $\tau_{\leq0}A_E$ be the differential graded subalgebra
consisting of its negative degrees and its degree-zero cocycles.
Inclusion and projection give quasi-isomorphisms of algebras
\[
 A_E\longleftarrow\tau_{\leq0}A_E\longrightarrow H^0(A_E)\cong S.
\]
Derived Yoneda for $E$
therefore gives a Morita equivalence between $\mathcal W$ and $S$.

For the finite-type Weinstein completion of $X$, the open--closed
map from the finite collection of Weinstein cocores hits the unit
\cite[Remark~1.3 and Corollary~1.4]{CDRGG}.
Thus Ganatra's nondegeneracy hypothesis holds, and the closed--open
map $\mathrm{CO}$ is an unshifted isomorphism of unital graded rings
\cite[Definition~1.1 and Theorem~1.1]{Ganatra}. Write $HH^*$ for
Hochschild cohomology over $\C$. Composing $\mathrm{CO}$ with the
cup-product isomorphism induced by the above Morita equivalence
\cite[\S\S3.8, 5.4]{KellerDG} gives
\[
 \Theta:SH^*_\Omega(X;\C)
     \xrightarrow[\sim]{\mathrm{CO}}HH^*(\mathcal W)
     \xrightarrow{\sim}HH^*(S).
\]
In degree zero, $HH^0(S)$ is the center of $S$, which equals $S$
since $S$ is commutative. Thus $\Theta^0:SH^0_\Omega(X;\C)
\xrightarrow{\sim}S$ identifies the pair-of-pants product with
multiplication in the specialized Gross--Siebert ring.

\subsection{Proof of Theorem~\ref{thm:surface-loop-floer}}
\label{subsec:surface-loop-sh-proof}
Put $R=R_X\otimes_{\F}\C$, so $\check X=\Spec R$. The
Gross--Siebert comparison of Theorem~\ref{mir:thm-intro-gs} gives
\[
 \alpha=\Phi_{\mathrm{GS}}\otimes_{\F}\C:R\xrightarrow{\sim}S.
\]
Since $\check X$ is smooth by hypothesis, so is $\Spec S$.
We may therefore apply the isomorphism $\Theta$ of
Subsection~\ref{subsec:surface-sh-zero-gs}.
The polyvector comparison of Theorem~\ref{thm:intro-global-polyvectors}
identifies $LH^*_\Omega(X)\otimes_{\F}\C$ with
$\Gamma(\check X,\bigwedge^\bullet T_{\check X})
=\bigwedge_R^\bullet\operatorname{Der}_{\C}(R)$.
Transport of derivations by $\alpha$ sends
$\partial\mapsto\alpha\circ\partial\circ\alpha^{-1}$ and induces
the graded algebra isomorphism $\alpha_*$ below. Since $S$ is smooth
over $\C$, the multiplicative Hochschild--Kostant--Rosenberg
isomorphism for $S$ \cite[Theorem~5.2]{HKR} gives the third arrow
in the composite
\[
 \begin{aligned}
 LH^*_\Omega(X)\otimes_{\F}\C
 &\xrightarrow[\sim]{\sigma}
     \bigwedge_R^\bullet\operatorname{Der}_{\C}(R)
  \xrightarrow[\sim]{\alpha_*}
     \bigwedge_S^\bullet\operatorname{Der}_{\C}(S)\\
 &\xrightarrow[\sim]{\mathrm{HKR}}HH^*(S)
  \xrightarrow[\sim]{\Theta^{-1}}SH^*_\Omega(X;\C).
 \end{aligned}
\]
Each arrow preserves the unit and the graded product. In degree zero
the composite is $(\Theta^0)^{-1}\circ\alpha$, since $\sigma^0=\mathrm{id}_R$
and the Hochschild--Kostant--Rosenberg map in degree zero is
$\mathrm{id}_S$.
$\hfill\Box$

\clearpage
\addtocontents{toc}{\protect\setcounter{tocdepth}{1}}

\part{Examples}\label{part:examples}


\section{Affine spaces}\label{sec:affine-spaces}\label{subsec:affine-spaces}\label{sec:line}

For $n\ge1$, write $\A^n=\prod_{i=1}^n\A^1_{x_i}$ and use the
product compactification
\[
 M=(\Pj^1)^n,\qquad
 D=\bigcup_{i=1}^nD_{i,\infty},\qquad
 A=\sum_{i=1}^n k_iD_{i,\infty},\quad k_i\in\Z_{>0},\qquad
 \Omega=dx_1\wedge\cdots\wedge dx_n.
\]
Here $D_{i,\infty}$ is the inverse image of $\infty$ under the
$i$th projection. The contact weight of a loop is
$\sum_i k_i\max(0,-\ord_t x_i)$, with the summand understood as
zero when $x_i=0$.

\begin{proposition}
\label{prop:affine-spaces}\label{thm:line}
In parity grading, and in the integer grading of any nowhere-zero
algebraic volume form,
\[
 LH^*(\A^n;\F)=0.
\]
This is the zero unital algebra, and its BV operator is zero.
\end{proposition}
\begin{proof}
The affine-line calculation \eqref{eqBA1}, with boundary
polarization $\{\infty\}$, gives $B_{dx}\simeq0$ for
$(\A^1,dx)$. Its proof uses the zero morphisms
$\bQ[2P]\to\bQ[2P+2]$ between successive pole stages, since
$\operatorname{Hom}_{\mathcal D(\bQ)}
(\bQ[2P],\bQ[2P+2])=0$. Extension of coefficients and
periodization preserve these zero morphisms; the vanishing is proved
before passing to the two-periodic category. These are morphisms
of the original derived pole-stage diagram. On a chosen filtered
cochain model they are represented by the inclusions
$F_P\mathbf B_{dx}\hookrightarrow F_{P+1}\mathbf B_{dx}$,
whose images in the derived category are zero; the inclusions
themselves need not be zero cochain maps.
For boundary polarization $k\{\infty\}$, with
$k\in\Z_{>0}$, the stages indexed by $kP$ form a cofinal
subdiagram and correspond to the stages just used. Thus the
same vanishing holds for each polarized affine-line factor above.
Theorem~\ref{prop:fusion-external-products}, iterated over the
$n$ affine-line factors, therefore gives, for $(\A^n,\Omega)$,
\[
 B_\Omega\simeq0,
 \qquad
 LH^*_{\Omega}(\A^n)
   \cong\bigotimes_{i=1}^n LH^*_{dx_i}(\A^1)=0
 \quad(n\ge1).
\]
The second isomorphism respects multiplication, unit, and BV
operator. Every nowhere-zero algebraic volume form on $\A^n$
is $c\Omega$, $c\in\C^*$, and this constant does not change
its volume indices. Proposition~\ref{prop:fusion-invariance-gradings}
shows that the unfiltered result also holds for any other polarized
SNC compactification, in particular $(\Pj^n,H_\infty)$.
\end{proof}

It is well-known that
$SH^*(\C^n;\F)=0$ for $n\ge1$, cf. \cite{Cieliebak}, \cite[Example 1.1]{GP}.

\section{Algebraic tori}\label{sec:torus}

\subsection{Loop cohomology and its logarithmic comparison}
\label{subsec:example-28-1}\label{subsec:example-28-3}
Let $T=(\C^*)^n$, $n\ge1$, and compactify it by $M=(\Pj^1)^n$.
Write $D_{i,0},D_{i,\infty}$ for the boundary divisors and put
\[
 A=\sum_i(k_{i,0}D_{i,0}+k_{i,\infty}D_{i,\infty}),
 \qquad k_{i,0},k_{i,\infty}\in\Z_{>0}.
\]
Every complex loop has the unique expression
$x_i(t)=t^{v_i}c_{i,0}(1+c_{i,1}t+\cdots)$, with
$v_i\in\Z$ and $c_{i,0}\ne0$. Its contact weight is
\[
 w(v)=\sum_i\bigl(k_{i,0}\max(v_i,0)
                     +k_{i,\infty}\max(-v_i,0)\bigr).
\]
Use the volume form
$\Omega=cx^b\,d\log x_1\wedge\cdots\wedge d\log x_n$,
where $c\in\C^*$ and $b\in\Z^n$; these are all the nowhere-zero
algebraic volume forms on $T$. Let $z_i,\theta_i$ be the classes
$T,\eta$ of Subsection~\ref{ex:cstar-global-fusion} in the $i$th
factor, and put $z^v=\prod_i z_i^{v_i}$. Write $\iota_i$ for
the degree-$-1$ derivation of $\Lambda_\F(\theta_1,\ldots,\theta_n)$
determined by $\iota_i\theta_j=\delta_{ij}$.

\begin{proposition}[Loop cohomology of a torus]
\label{prop:torus-fusion}\label{thm:torus}\label{lp:thm:torus-filtered-product}
\begin{enumerate}[(i),wide=0pt]
\item The filtered loop BV algebra is
\begin{equation}\label{eq:torus-fusion}
 \begin{aligned}
 LH^*_{\Omega}(T)&\cong\F[z_1^{\pm1},\ldots,z_n^{\pm1}]
                  \otimes\Lambda_\F(\theta_1,\ldots,\theta_n),\\
 (z^v\alpha)(z^w\beta)&=z^{v+w}(\alpha\wedge\beta),
 \end{aligned}
\end{equation}
with $|z^v|=2\langle b,v\rangle$, $|\theta_i|=1$, and
\[
 F_PLH^*_{\Omega}(T)
    =\bigoplus_{w(v)\le P}z^v\Lambda_\F(\theta_1,\ldots,\theta_n).
\]
Its operator is
\begin{equation}\label{eq:torus-bv}
 \Delta(z^v\alpha)=z^v\sum_i v_i\iota_i\alpha.
\end{equation}
In parity the same formulas hold with all $z_i$ even.

\item In either grading, for the displayed compactification let $\ell_v\alpha$
denote the logarithmic class of contact $v$ and coefficient
$\alpha\in H^*((S^1)^n;\F)$. The map
\[
 \Phi_{\log}:H^*_{\log}(M,D)\xrightarrow{\sim}LH^*(T),
 \qquad \ell_v\alpha\longmapsto z^v\alpha
\]
is a filtered graded-vector-space isomorphism satisfying
$\Phi_{\log}\Delta_{\log}=\Delta\Phi_{\log}$.
It is not an algebra homomorphism: for the standard basis vector
$e_i\in\Z^n$,
\[
 \ell_{e_i}\ell_{-e_i}=0,\qquad z_i z_i^{-1}=1.
\]
There is instead an isomorphism of associated-graded BV algebras
$\operatorname{gr}_F LH^*(T)\cong H^*_{\log}(M,D)$.
For $b=0$ and unit boundary weights, the Rees algebra is
\[
 \Lambda_\F(\theta_1,\ldots,\theta_n)\otimes
 \frac{\F[s,U_1,V_1,\ldots,U_n,V_n]}
      {(U_iV_i-s^2:1\le i\le n)},
 \qquad U_i=sz_i,\quad V_i=sz_i^{-1}.
\]
Its fibers at $s=0$ and $s=1$ are the logarithmic and Laurent
rings, respectively.

\item The graded ring in \eqref{eq:torus-fusion} is isomorphic
to $SH^*(T)$, in parity or with the corresponding volume grading.
\end{enumerate}
\end{proposition}
\begin{proof}
For (i), the one-dimensional calculations in
Example~\ref{ex:cstar-loop-complex},
Subsection~\ref{ex:cstar-global-fusion}, and
Proposition~\ref{ex:cstar-global-bv} give
\[
 LH^*_{x_i^{b_i}dx_i/x_i}(\Gm)
      =\F[z_i^{\pm1}]\otimes\Lambda_\F(\theta_i),\qquad
 \Delta(z_i^v\theta_i)=v z_i^v,\quad\Delta(z_i^v)=0.
\]
Changing the positive boundary weights only relabels the winding
filtration by $k_{i,0}\max(v,0)+k_{i,\infty}\max(-v,0)$.
Apply Theorem~\ref{prop:fusion-external-products} to the $n$
factors. Its tensor multiplication gives \eqref{eq:torus-fusion};
its sum filtration gives $w(v)$. Its BV sign on the $i$th
factor of $z^v\theta_1^{\epsilon_1}\cdots\theta_n^{\epsilon_n}$
is $(-1)^{\epsilon_1+\cdots+\epsilon_{i-1}}$, precisely the
sign of $\iota_i$. This proves \eqref{eq:torus-bv}.
The scalar $c$ does not change the volume indices.

For (ii), the winding loci are clopen and retract to $(S^1)^n$,
so Corollary~\ref{thm:separated-sectors}(2) gives $\Phi_{\log}$.
Use angular classes of the original coordinates $x_i$; the
normal angular class at infinity is then $-\theta_i$.
Equation~\eqref{eq:bv-first-page} gives
$\Delta_{\log}(\ell_v\alpha)=\ell_v\sum_i v_i\iota_i\alpha$,
which proves operator compatibility. The logarithmic product
vanishes when opposite coordinate windings occur, since the
two boundary points of that factor are disjoint. Otherwise it
is the product of \eqref{eq:torus-fusion}. More explicitly,
for $v,v'\in\Z^n$,
\[
 w(v)+w(v')-w(v+v')
 =\sum_{i:v_iv'_i<0}(k_{i,0}+k_{i,\infty})
                         \min(|v_i|,|v'_i|).
\]
The positive boundary weights make this difference positive
exactly when opposite nonzero coordinate windings occur.
Thus the products agree on the associated graded.
For unit weights the Rees algebra is generated by $s$, $sz_i$,
$sz_i^{-1}$, and the $\theta_i$, with precisely the displayed
relations: in each Laurent winding its coefficient module is
$\F[s]s^{\sum_i|v_i|}z^v$. Indeed, the relations reduce every
monomial to
\[
 s^a\prod_i U_i^{\max(v_i,0)}V_i^{\max(-v_i,0)},
 \qquad a\ge0,\quad v\in\Z^n,
\]
whose images $s^{a+\sum_i|v_i|}z^v$ are linearly independent.
This proves that there are no further relations and verifies both fibers.
All direct sums consist of finitely supported elements.

Part (iii) follows from the known computation for $SH^*(T)$, see \cite[Example 1.1]{GP}.
\end{proof}

\section{Smooth affine curves}\label{sec:curves}\label{sec:punctured-line}
For a smooth connected affine curve $C$ other than $\A^1$ and
$\Gm$, which were addressed in Sections~\ref{sec:affine-spaces}
and~\ref{sec:torus}, this section computes $LH^*(C)$ as a filtered BV algebra.
Writing $\overline C$ for its smooth projective completion and
$D=\overline C\setminus C$ for its nonempty finite boundary, we prove
\[
 LH^*(C)\cong H^*_{\log}(\overline C,D)
\]
as filtered BV algebras, give the contact summands and the
formulas for their products and BV operator, and prove degeneration
of the contact spectral sequence at $E_1$. Logarithmic PSS of \cite{GP} also
gives an isomorphism of graded rings $LH^*(C)\cong SH^*(C)$.
The assertions use parity grading, or the integer grading of a
chosen nowhere-vanishing algebraic volume form.

Write $A=\sum_{p\in D}k_pp$, $k_p\in\Z_{>0}$, and
$\ord_p\Omega=-a_p$ when a volume form is specified. The positive
contact labels are $(m,p)$, $m\ge1$, with
\[
 w(m,p)=k_pm,\qquad\mu_\Omega(m,p)=2(1-a_p)m.
\]
We recall the following classical result.

\begin{lemma}
\label{lem:curve-surface-extension}
Let $k$ be an algebraically closed extension of $\C$, put
$R=k\llbracket s,t \rrbracket$, and let $T$ be a smooth proper connected curve
of positive genus over $k$. Every rational $k$-map
$\Spec R\dashrightarrow T$ extends uniquely to a morphism.
\end{lemma}
\begin{proof}
Resolve the rational map by a finite sequence of blowups at closed
points, obtaining $\rho:S'\to\Spec R$ and a morphism $f:S'\to T$;
see \cite[Tag 0C5I]{Stacks}. Here the normalized blowups in that
result are ordinary point blowups, since all the surfaces are
regular. Every nontrivial fiber of $\rho$ is a connected union
of rational curves. A morphism $\Pj^1_k\to T$ is constant by
Riemann--Hurwitz. The constants agree on intersecting components,
so connectedness makes $f$ constant on each such fiber.
Let $\Gamma\subset\Spec R\times_k T$ be the reduced closure
of the graph of the original rational map. The proper map
$(\rho,f)$ has image $\Gamma$, and the constancy just proved
shows that every fiber of $\Gamma\to\Spec R$ has one point.
This projection is proper and quasi-finite, hence finite.
It is birational, and $R$ is normal, so it is an isomorphism.
Projection to $T$ gives the extension. Uniqueness follows from
separatedness of $T$.
\end{proof}

The following proposition recovers the component classification
of Bouaziz \cite[Theorem~4.10]{BouazizCurves}. We give
an independent proof.

\begin{proposition}
\label{prop:global-loop-components}
The arc locus $\mathscr L^+C$ and the contact
loci $\mathscr P_{m,p}$, with $p\in D$ and $m\ge1$, are open
and closed in each reduced pole stage in which they occur.
In our analytic topology of
Convention~\ref{conv:loop-analytic-topology}, these loci give
a topological coproduct
\[
 \mathscr LC=\mathscr L^+C\;\sqcup\!
 \coprod_{p\in D,\,m\ge1}\mathscr P_{m,p}.
\]
The zero-contact locus retracts to $C$. Choose an algebraic local
parameter $x_p\in\mathcal O_{\overline C,p}$ and the resulting
identification $\widehat{\mathcal O}_{\overline C,p}=\C\llbracket x_p \rrbracket$.
Then $\mathscr P_{m,p}$, with its subspace analytic topology,
is identified with the space of formal power series of the form
\begin{equation}\label{eq:local-global-homeomorphism}
 x_p(t)=t^m(c_0+c_1t+c_2t^2+\cdots),
                  \qquad c_0\ne0.
\end{equation}
It is homeomorphic to $\C^*\times\prod_{j\ge1}\C$ and retracts
to $S^1$. Consequently
\[
 \pi_0(\mathscr LC)=\{\text{arc component}\}
                \sqcup(D\times\Z_{\ge1}).
\]
\end{proposition}
\begin{proof}
\emph{Step 1: separate the contact labels in genus zero.}
If $\overline C=\Pj^1$, write
$D=\{\infty,a_1,\ldots,a_r\}$, where $r\ge2$. The functions $g_i=x-a_i$ are units on $C$.
By Proposition~\ref{lp:prop:unit-valuation}, their loop orders
are locally constant. Their order vectors are respectively
\[
 (0,\ldots,0),\qquad m e_i,\qquad(-m,\ldots,-m)
\]
on the arc locus, on $\mathscr P_{m,a_i}$, and on
$\mathscr P_{m,\infty}$, where $e_i$ is the $i$-th standard
basis vector. These vectors are distinct, so they separate all
contact labels. On a reduced coefficient chart the same
unit-order decomposition is algebraically open and closed,
as in the proof of that proposition. This proves the clopenness
assertion when $g(\overline C)=0$.

\emph{Step 2: show that specialization preserves contact labels in positive genus.}
Suppose $g(\overline C)>0$, and let $k$ be an algebraically closed
extension of $\C$. A family of loops over
$\Spec k\llbracket s \rrbracket$ in a fixed pole stage is a morphism
\[
 \gamma:\Spec\bigl(k\llbracket s \rrbracket(\!( t )\!)\bigr)
        =\Spec\bigl(k\llbracket s,t \rrbracket[1/t]\bigr)\longrightarrow C_k.
\]
Lemma~\ref{lem:curve-surface-extension} extends it to
$\bar\gamma:\Spec k\llbracket s,t \rrbracket\to\overline C_k$.
For $p\in D$, its inverse image is an effective Cartier
divisor, possibly empty, supported on $t=0$, since $\gamma$
takes values in $C_k$. The ring $k\llbracket s,t \rrbracket$ is factorial, so
\[
                    \bar\gamma^*(p)=m_p[t=0]
                    \qquad(m_p\ge0).
\]
If $m_p>0$, a local parameter at $p$ therefore pulls back to
$t^{m_p}u$ with $u\in k\llbracket s,t \rrbracket^\times$. The unit remains
invertible after setting $s=0$ and after passing to $k(\!( s )\!)\llbracket t \rrbracket$.
Thus the generic and special loops have the same contact order
at every $p$. At most one $m_p$ is positive, because distinct
points of $D$ have disjoint inverse images. If all are zero,
both loops are arcs in $C_k$.

\emph{Step 3: prove the clopenness and topological coproduct assertions.}
Continue with $g(\overline C)>0$. Fix $P$.
 By Theorem~\ref{thm:drinfeld-finite-factors}, the
pole stage admits nice charts $U=Y\times\A^J$. There are only
finitely many possible labels, namely contact zero and the pairs
$(m,p)$ with $k_pm\le P$. After retaining finitely many tail
coordinates in $Y$, their inverse images in $U_{\mathrm{red}}$
have the form $Z_\lambda\times\A^J$, where the $Z_\lambda$
form a finite constructible partition of $Y_{\mathrm{red}}$.
Indeed, the equations of the finitely many closed pole substages
use finitely many coordinates. The reduced differences have
the finite open-and-closed contact decomposition of
Proposition~\ref{prop:loop-sectors}; these pieces also descend
after finite retention. This is the geometric argument for finite-factor presentations in Step~3 of Subsection~\ref{subsec:volume-index-proofs},
which does not require a volume form.

For any map $\Spec\C\llbracket s \rrbracket\to Y_{\mathrm{red}}$, set all remaining
tail coordinates equal to zero to lift it to $U$. This does not
change either endpoint label, since every piece has the form
$Z_\lambda\times\A^J$. The resulting
coefficient family has the same generic and special contact label
by Step~2. Each $Z_\lambda$ is therefore closed:
otherwise a boundary point of this constructible subset can be
reached by an integral algebraic curve whose generic point lies
in it; normalization and completion at that point give just such
a map from $\Spec\C\llbracket s \rrbracket$, a contradiction. The same reasoning
applies to its constructible complement, so $Z_\lambda$ is also
open. These properties descend along the surjective \'etale atlas,
proving the clopenness assertion on the reduced pole stage.
In either genus, analytic realization and the final topology
over pole stages give the asserted topological coproduct.

\emph{Step 4: prove the homeomorphism, retraction, and connected-component assertions.}
Proper extension to $\overline C$ and the completed local
parameter give the bijection \eqref{eq:local-global-homeomorphism}.
For an affine generator $g$, substitution in its Laurent expansion
at $p$ shows that every coefficient of $g(\gamma)$ is a regular
function of finitely many $c_j$ and $c_0^{-1}$. Conversely write
$x_p=a/b$ with $a,b\in\mathcal O(C)$. On this contact locus the
leading coefficient of $b(\gamma)$ is nonzero, so coefficientwise
division gives continuous inverse maps, using finitely many
coefficients at each order. This proves the homeomorphism.
Contracting the positive tail and radially retracting $\C^*$
give $S^1$. The arc retraction is
Lemma~\ref{lp:lem:arc-contraction}. Hence all the displayed
open-and-closed loci are connected and are precisely the
connected components, proving the remaining assertions.
\end{proof}

The next lemma determines the support of fusion under the contact
decomposition. This requires a collision argument in addition
to the specialization argument for a single loop above.

\begin{lemma}[Contact support of fusion for curves]
\label{lem:curve-fusion-support}
Let
\[
 \mathcal V_C:=\{0\}\cup\{me_p:p\in D,\ m\ge1\}\subset\N^D,
\]
where $e_p$ is the standard basis vector indexed by $p$ and $0$
denotes the arc contact. For $v\in\mathcal V_C$, let $B_v$ denote
the corresponding contact summand of the loop complex in
Corollary~\ref{thm:separated-sectors}, which applies by
Proposition~\ref{prop:global-loop-components}. Use the integer
grading when a volume form is specified and the parity grading
otherwise. For $\uct,v,w\in\mathcal V_C$, the component
\[
 m_{2;v,w}^{\uct}:B_v\otimes_{\F}B_w\longrightarrow B_{\uct}
\]
of binary fusion, obtained by inclusion of the two input summands
and projection to the output summand, is zero unless $\uct=v+w$.
In particular, fusion has no lower-weight terms, as required in
\eqref{eq:separated-no-lower-products}.
\end{lemma}
\begin{proof}
\emph{Step 1: prove the support assertion in genus zero.}
If $\overline C=\Pj^1$, then $|D|\ge3$. 
The coefficient vectors of the principal divisors on $\overline C$
of global units $g\in\mathcal O(C)^\times$ span
$\{(a_p)_{p\in D}\in\Z^D:\sum_pa_p=0\}$.
Choose finitely many units $g_1,\ldots,g_r$ whose divisor vectors
generate this lattice. Their valuation matrix $ V:=\bigl(\operatorname{ord}_p(g_i)\bigr)_{
                 1\le i\le r,\ p\in D}$ has $\ker V=\Z\mathbf1_D$, where $\mathbf1_D=(1)_{p\in D}$.
The sector decomposition in the proof of
Proposition~\ref{prop:fusion-unit-valuations}
therefore implies that a nonzero output component of contacts
$v,w$ must have contact $\uct$ satisfying
\[
 \uct-v-w=c\mathbf1_D,\qquad c\in\Z.
\]
If $c<0$, a coordinate outside $\operatorname{supp}(v+w)$ is
negative in $\uct$. If $c>0$, every coordinate of $\uct$ is positive.
Both are impossible for a realized curve contact; hence $c=0$.
Thus $\uct=v+w$, which is realized precisely when the positive
inputs lie at the same puncture, or at least one input is zero.

\emph{Step 2: identify the generic input orders in a one-parameter collision in positive genus.}
Suppose $g(\overline C)>0$. Choose a nice chart
$U\cong Y\times\A^J$ of the one-parameter binary collision
family of Definition~\ref{def:geometric-collision-data},
described explicitly in Example~\ref{ex:kv-binary-family},
such that the collision parameter factors through $Y$ and the
pulled-back input and output contact pieces are inverse images
under $\operatorname{pr}_Y:U\to Y$ of locally closed subschemes
of $Y$. Test a specialization from a separated input piece with
contacts $v,w$ by a morphism $\Spec k\llbracket u \rrbracket\to Y$,
after extension to an algebraically closed field $k$. Lift it to
$U$ by the zero section $Y\to Y\times\A^J$; the
inverse-image descriptions preserve its input and output contact
pieces. Its composite with the
collision family is an actual family of maps from the marked
source. Put $R=k\llbracket u \rrbracket$ and $K=k(\!( u )\!)$. Its collision parameter
is $s=s(u)\in uR\setminus\{0\}$, and its completed
marked-source algebra is
\[
 \mathcal O_s=R\llbracket q \rrbracket[z]/(z^2-s^2-q).
\]
Every element of $\mathcal O_s$ has a unique expression
$A(q)+zB(q)$ with $A,B\in R\llbracket q \rrbracket$. Substitution gives an isomorphism
\[
 \begin{aligned}
 \phi:\mathcal O_s&\xrightarrow{\sim}R\llbracket z \rrbracket,\\
 A(q)+zB(q)&\longmapsto
              A(z^2-s^2)+zB(z^2-s^2).
 \end{aligned}
\]
Indeed, for $H(z)=\sum_{r\ge0}h_rz^r\in R\llbracket z \rrbracket$, put
$E(q)=\sum_{r\ge0}h_{2r}q^r$ and
$O(q)=\sum_{r\ge0}h_{2r+1}q^r$. Then
$H(z)=E(z^2)+zO(z^2)$, and the inverse sends $H$ to
$E(q+s^2)+zO(q+s^2)$. These substitutions are
continuous for the joint $(u,q)$-adic and $(u,z)$-adic
topologies; at each fixed coefficient, the required sums
converge $u$-adically because $s\in uR$.

We call the contact orders $v_p,w_p$ at $p\in D$ of the two
input loops over $K=k(\!( u )\!)$ the generic input orders.
They are measured in the local source parameters
$t_+=z-s$ and $t_-=z+s$, respectively.
To compare them with orders along the prime divisors
$(z-s)$ and $(z+s)$ of $R\llbracket z\rrbracket$, use the
coefficientwise homomorphism
\[
 R\llbracket q\rrbracket\longrightarrow K\llbracket q\rrbracket,
 \qquad
 \sum_{j\ge0}a_jq^j\longmapsto\sum_{j\ge0}(a_j/1)q^j,
 \qquad a_j\in R.
\]
Together with $z\mapsto z$, it induces
\[
 \mathcal O_s\longrightarrow
 \mathcal O_{s,K}:=K\llbracket q\rrbracket[z]/(z^2-s^2-q).
\]
The two formal factors of $\mathcal O_{s,K}$ are
$K\llbracket t_+\rrbracket$ and $K\llbracket t_-\rrbracket$,
with factor maps
\[
 \begin{aligned}
 z&\longmapsto s+t_+,&q&\longmapsto t_+(t_++2s),\\
 z&\longmapsto -s+t_-,&q&\longmapsto t_-(t_--2s).
 \end{aligned}
\]
Since $s\in uR$, the changes of variables $z=\pm s+t_\pm$
are defined over $R$. The factor maps above identify
\[
 K\llbracket t_\pm\rrbracket
 \cong
 \widehat{\bigl(R\llbracket z\rrbracket\bigr)_{(z\mp s)}},
\]
where the completion is taken at the maximal ideal of the
indicated discrete valuation ring. These maps send the
uniformizer $z\mp s$ to $t_\pm$ and therefore preserve orders
of vanishing. Thus the generic input contact orders are
exactly the orders along the prime divisors $(z-s)$ and
$(z+s)$ of $R\llbracket z\rrbracket$.

\emph{Step 3: prove that every such specialization has output contact $v+w$.}
Pulling back the collision family along the lift to $U$ in
Step~2 gives the wedge
\[
 \gamma:\Spec R\llbracket z \rrbracket[(z^2-s^2)^{-1}]\longrightarrow C_k.
\]
Since $g(\overline C)>0$, Lemma~\ref{lem:curve-surface-extension} extends it to
$\bar\gamma:\Spec R\llbracket z \rrbracket\to\overline C_k$.
For each $p\in D$, factoriality and the generic input orders give
\[
 \bar\gamma^*(p)
    =v_p[z-s(u)=0]+w_p[z+s(u)=0].
\]
There is no additional component: outside those two divisors the
map already takes values in $C_k$. In particular, the generic
point of $u=0$ lies in this open locus, since there $q=z^2\ne0$.
If $v_p+w_p>0$, the local parameter $x_p$ chosen in
Proposition~\ref{prop:global-loop-components} pulls back to
$(z-s)^{v_p}(z+s)^{w_p}$ times a unit of $R\llbracket z \rrbracket$. Hence
\[
 \operatorname{ord}_z\!\left(
       \left.\bar\gamma^*x_p\right|_{u=0}\right)=v_p+w_p,
\]
independently of $\operatorname{ord}_u s(u)$; this is the output
contact order. If all these integers are zero, the output is
an arc in $C_k$. Moreover, if $v_p>0$ and $w_q>0$
for $p\ne q$, the common source origin would map to both $p$
and $q$, which is impossible. Thus every specialization has
exactly the output contacts asserted above.

\emph{Step 4: deduce vanishing of the excluded fusion morphisms.}
The wedge calculation in Steps~2--3 controls the specialization
morphism on each finite factor in this finite-factor presentation. On such a
factor $Y$, let
$W\subset Y|_{s\ne0}$ be the chosen reduced input contact
piece and let $Z\subset Y$ be its reduced closure. As in the
closure construction in the proof of
Lemma~\ref{lem:fusion-leading-product}, first descend the
clopen input idempotent to the finite factor; the closure
ideal is then obtained by contraction from the localization
at $s$. Thus the closure is independent of the
remaining affine tail. Every closed point of
$Z_0:=Z\cap V(s)$ is reached by an integral curve in $Z$
whose generic point lies in $W$. Normalizing the curve and
completing at a point above the chosen closed point gives a morphism
$\Spec\C\llbracket u\rrbracket\to Z$. The zero-section lift of its
composite with $Z\hookrightarrow Y$ gives a wedge to which
Steps~2--3 apply.
Hence $Z_0$ misses every
output contact piece $E\subset Y_0:=Y\cap V(s)$ whose contact differs from $v+w$.

Proper-trace naturality in
Lemma~\ref{lp:lem:bm-collision-specialization}, applied to
$Z\hookrightarrow Y$, makes specialization from $W$ factor
through the normalized complex on $Z_0$ before its closed
trace to $Y_0$. Restriction of this factorization to $E$ is
the zero morphism: proper base change identifies the restriction
to $E$ of the pushed-forward normalized complex on $Z_0$ with
the corresponding pushforward from
$Z_0\times_{Y_0}E=\varnothing$, hence with zero. Equivalently, open
restriction to $Y\setminus Z$ annihilates this contribution
to the localization boundary. The comparison maps in these finite-factor presentations
on refinements and on intersections of charts preserve the
proper traces and localization diagrams, and hence this zero
factorization. It remains zero after totalization and descent.
Thus the fusion component with output contact $\uct\ne v+w$ is
the zero morphism. If $v+w$ is realized, its weight is the sum
of the input weights. Hence there are no lower-weight terms,
proving the last assertion as well.
\end{proof}

\begin{proposition}[Loop cohomology of the remaining affine curves]
\label{prop:curve-fusion}\label{prop:curve-bv}
\label{thm:curves}\label{thm:global-loop-curves}\label{lem:weighted-punctured-torus}
For the curve $C=\overline C\setminus D$ fixed above, the following
hold in parity, or in the integer grading of a specified volume form.
\begin{enumerate}[(i),wide=0pt]
\item There is an equivalence of filtered objects in the sense of
Definition~\ref{def:filtered-derived-object},
\begin{equation}\label{eq:curve-complex}
 B_\Omega\simeq C^*(C;\F)\oplus
  \bigoplus_{p\in D,\,m\ge1}
                  C^{*-2(1-a_p)m}(S^1;\F),
\end{equation}
with weights $0$ and $k_pm$, respectively. In parity, replace
these summands by their two-periodizations. The induced
contact-preserving additive identification is
\begin{equation}\label{eq:global-loop-sh}
 LH^j_\Omega(C)\cong H^j(C;\F)\oplus
   \bigoplus_{p\in D,\,m\ge1}H^{j-2(1-a_p)m}(S^1;\F)
       =H^j_{\log}(\overline C,D).
\end{equation}
The contact spectral sequence degenerates at $E_1$.

\item The identification in (i) preserves multiplication.
Writing $r_p:H^*(C;\F)\to H^*(S^1;\F)$ for restriction to a
small positively oriented circle around $p$, its formulas are
\begin{align}
 ab&=a\smile b, &&a,b\in H^*(C;\F),\notag\\
 a(\alpha t_p^m)&=(r_pa\smile\alpha)t_p^m,
       &&m>0,\label{eq:curve-fusion}\\
 (\alpha t_p^m)(\beta t_q^l)
   &=\begin{cases}(\alpha\smile\beta)t_p^{m+l},&p=q,\\
                   0,&p\ne q,
      \end{cases} &&m,l>0.\notag
\end{align}
Here $t_p^m$ has degree $2(1-a_p)m$, or even parity, and
weight $k_pm$.

\item Let $\theta_p$ evaluate to $1$ on the positive normal
circle at $p$. Then
\begin{equation}\label{eq:curve-bv}
 \Delta(a)=0\quad(a\in H^*(C;\F)),\qquad
 \Delta(t_p^m)=0,\qquad
 \Delta(\theta_p t_p^m)=m\,t_p^m.
\end{equation}
Thus (i) identifies $LH^*(C)$ with
$H^*_{\log}(\overline C,D)$ as a filtered BV algebra.

\item There is an isomorphism of graded rings
$LH^*(C)\cong SH^*(C)$.
\end{enumerate}
\end{proposition}
\begin{proof}
\emph{(i)}
Apply Corollary~\ref{thm:separated-sectors}(1)--(2) to
Proposition~\ref{prop:global-loop-components}. Contact zero has
normal space $C$, and the normal circle at a boundary point is
$S^1$. Only finitely many of these contact summands occur at
each pole stage. For the contact locus $\mathscr P_{m,p}$, a smooth finite contact factor
of dimension $d$ has volume index $d+(1-a_p)m$, so its normalized
dualizing complex is $\underline\F[-2(1-a_p)m]$.
On every old clopen contact piece, a pole transition restricts
to the identity proper trace. These identifications respect the
finite-factor comparisons and their compositions, as in
Corollary~\ref{thm:separated-sectors}(1), so the coherent pole diagram
has summand inclusions as its maps. Its homotopy colimit is the
displayed direct sum. This proves \eqref{eq:curve-complex}, its
cohomology, and degeneration. In parity use the two-periodized
unshifted contact complexes; no numbers $a_p$ or volume form
are needed in that case. This proves (i).

\emph{(ii)}
Lemma~\ref{lem:curve-fusion-support} shows that the product of
contacts $v,w\in\mathcal V_C$ has only the component of contact
$v+w$. If both contacts are positive at distinct punctures,
$v+w\notin\mathcal V_C$, so the product is zero. Otherwise
$v+w$ is realized and has the sum of the input weights.
Lemma~\ref{lem:fusion-leading-product} therefore computes the
entire product and gives \eqref{eq:curve-fusion}: the relevant
face restrictions are the identity at a fixed puncture and
$r_p$ from contact zero to that puncture. On contact zero the
product is the ordinary cup product. These are exactly the
formulas of Definition~\ref{def:log-cohomology}, proving (ii).

\emph{(iii)}
Source rotation preserves the arc locus $\mathscr L^+C$ and
each contact locus $\mathscr P_{m,p}$. Consequently, under
\eqref{eq:curve-complex}, the composite of the inclusion of
any summand, the BV operator, and the projection to any
distinct summand is zero.
On contact zero the operator vanishes by
Lemma~\ref{lem:source-dilation}. On the contact locus $\mathscr P_{m,p}$, the leading
coefficient in \eqref{eq:local-global-homeomorphism} transforms
by $c_0\mapsto\lambda^m c_0$. The contact comparison in the proof
of Proposition~\ref{prop:logarithmic-bv} identifies this with
\[
 a_m:S^1\times S^1\longrightarrow S^1,\qquad
              (\lambda,c)\longmapsto\lambda^m c.
\]
With the acting circle first and generator $\theta$, one has
$a_m^*\theta_p=1\otimes\theta_p+m\theta\otimes1$ and
$a_m^*1=1\otimes1$. Slant product gives \eqref{eq:curve-bv},
which agrees with \eqref{eq:bv-first-page} on each summand.
Together with (ii), this proves the filtered BV algebra
identification in (iii).

\emph{(iv)}
Use the integrable complex structure on $\overline C$.
The multiplicatively topological condition of
\cite[Definition~5.4]{GP} excludes nonconstant spheres with
at most two boundary preimages in every boundary stratum,
including the ambient compactification. If $g(\overline C)>0$,
there is no nonconstant holomorphic map
$\Pj^1\to\overline C$. If $\overline C=\Pj^1$, then
$|D|\ge3$, and every nonconstant self-map is surjective;
the nonempty fibers over three distinct boundary points give
at least three distinct boundary preimages. Every nonempty
proper boundary stratum is a point and admits no nonconstant
sphere. Thus all the strata satisfy the condition. The
logarithmic PSS theorem \cite[Theorem~5.12]{GP} identifies
$H^*_{\log}(\overline C,D)$ with $SH^*(C)$.
Its volume grading is the shift
$2(1-a_p)m$ used here, and its parity version is covered by
\cite[Remark~4.32]{GP}. Composing this ring isomorphism with
the identification in (ii) proves (iv).
\end{proof}

\section{Projective hypersurface complements}
\label{sec:plane-complements}
\label{sec:hypersurface-complements}
\label{sec:hyperplanes}
Let $D=\bigcup_{i=1}^kD_i\subset\Pj^n$ be an SNC union of
smooth irreducible hypersurfaces of degrees $d_i>0$, where
$n\ge2$ and $k\ge n+1$. Fix
\[
 X=\Pj^n\setminus D,\qquad
 A=\sum_{i=1}^k k_iD_i\quad(k_i>0),\qquad
 \mathbf d=(d_1,\ldots,d_k),\quad
 d_\Sigma=\sum_i d_i,\quad g=\gcd(d_1,\ldots,d_k).
\]
We assume $d_\Sigma\ge n+2$: the excluded case consists of
$n+1$ hyperplanes in general position and has complement
$(\C^*)^n$, computed in Section~\ref{sec:torus}.
Write $\mathbf e_i$ for the standard basis of $\Z^k$ and
$w_A(\vct)=\sum_i k_iv_i$ for the contact weight.

The section determines the following loop invariants and comparisons.
\begin{enumerate}[(i)]
\item For every pair above, the contact decomposition identifies
$LH^*(X)$ additively with $H^*_{\log}(\Pj^n,D)$ and computes
the full loop BV operator. Logarithmic PSS of \cite{GP} gives a filtered
additive comparison with $SH^*(X)$.
\item The fusion product is the logarithmic product if
$k\ge2n+1$, if all $D_i$ are hyperplanes, or if $n=2$ and
$g\mid3$. Thus in these ranges the contact identification is
an isomorphism of filtered unital BV algebras. For hyperplane
complements we give a presentation of the entire algebra
and its BV operator.
\item The comparison with $SH^*(X)$ is proved to be
multiplicative when $k\ge2n+1$ or all $D_i$ are hyperplanes via \cite{GP}.
For plane complements with $k\in\{3,4\}$ and $g\mid3$,
other than four lines, the loop algebra is determined here
but multiplicativity of the symplectic comparison is not
verified. With $k\in\{3,4\}$ and $g\nmid3$, the plane
fusion ring remains undetermined here. In higher dimension,
the general nonlinear case $n+1\le k\le2n$ is likewise
restricted by the output and degree tests below, without a
general computation of its fusion ring or a multiplicative
comparison.
\end{enumerate}
All assertions use parity grading and every available volume
grading.

\subsection{Contact decomposition, volume grading, and the BV operator}
\label{subsec:example-30-2}
Choose homogeneous equations $G_i$ for $D_i$ and put
\[
 u_i=G_i^{d_1}/G_1^{d_i}\in\mathcal O(X)^*\quad(2\le i\le k),
 \qquad
 W:\Z^k\longrightarrow\Z^{k-1},\quad
 W\vct=(d_1v_i-d_iv_1)_{i=2}^k.
\]

\begin{proposition}[Separation by projective units]
\label{prop:hypersurface-separation}
Every realized contact locus is clopen, and there is a
canonical filtered additive identification
\[
 J_D:LH^*(X)\xrightarrow{\sim}H^*_{\log}(\Pj^n,D),
 \qquad J_D\Delta=\Delta_{\log}J_D.
\]
Here each contact monomial $t^{\vct}$ has weight $w_A(\vct)$,
and each element has finite contact support.
\end{proposition}
\begin{proof}
The valuation matrix satisfies
$\ker(W\otimes\bQ)=\bQ\mathbf d$.
If $W\vct=W\wct$ for realized contacts, then
$\vct-\wct=c\mathbf d$ for $c\in\bQ$.
For $c>0$, every entry of $\vct$ would be positive; for
$c<0$, every entry of $\wct$ would be positive. Both
contradict $|I(\vct)|,|I(\wct)|\le n<k$. Thus $c=0$.
Lemma~\ref{lem:unit-separation} and
Corollary~\ref{thm:separated-sectors}(2) give the clopen contact
decomposition and $J_D$.

Source rotation preserves each contact locus. The equivariant
contact comparison of Proposition~\ref{prop:logarithmic-bv}
therefore computes the operator on the actual summand of
$LH^*(X)$, rather than only its associated-graded quotient:
\begin{equation}\label{eq:hypersurface-bv}
 \Delta(\alpha t^{\vct})
   =\left(\sum_{i\in I(\vct)}v_i\iota_i\alpha\right)t^{\vct}
   =\Delta_{\log}(\alpha t^{\vct}).
\end{equation}
The operators $\iota_i$ are those of the normal-circle
actions defined in Subsection~\ref{subsec:logarithmic-bv}.
On contact zero the operator is zero by
Lemma~\ref{lem:source-dilation}.
\end{proof}

A nowhere-vanishing volume form on $X$ extends to a meromorphic
$n$-form on $\Pj^n$ with divisor $-\sum_i a_iD_i$, $a_i\in\Z$.
Conversely, such a form exists precisely when
\begin{equation}\label{eq:hypersurface-grading-condition}
 \sum_i a_id_i=n+1,
 \qquad\text{equivalently}\qquad g\mid n+1.
\end{equation}
For a chosen form $\Omega$,
the same contact decomposition reads
\[
 LH^j_\Omega(X)\cong
 \bigoplus_{\substack{\vct\text{ realized}}}
 H^{j-\mu_\Omega(\vct)}(S_{I(\vct)}^\circ;\F)t^{\vct},
 \qquad
 \mu_\Omega(\vct)=2\sum_i(1-a_i)v_i.
\]
The zero-contact coefficient is $H^*(X;\F)$; for nonzero
contact every connected component of $D_{I(\vct)}^\circ$
contributes separately. These are the gradings of
Definition~\ref{def:log-cohomology}.

\subsection{Possible fusion outputs and a range without lower terms}
For input contacts $\vct,\wct$, a nonzero output component
must have contact
\begin{equation}\label{eq:projective-fusion-output}
 \boldsymbol\nu=\vct+\wct-c\mathbf d,
 \qquad
 c=\min_i\frac{v_i+w_i}{d_i}\in\bQ_{\ge0}.
\end{equation}
Indeed, Proposition~\ref{prop:fusion-unit-valuations}
gives $W\boldsymbol\nu=W\vct+W\wct$.
The kernel of $W$ gives the displayed expression with some
$c\in\bQ$. If $c<0$, the output has full support, which is
impossible. Nonnegativity gives
$c\le\min_i(v_i+w_i)/d_i$, and a realized output has at
least one zero coordinate, giving equality. Thus there is
at most one output contact vector. It must be integral,
nonnegative, and realized; its coefficient can still have
several connected-component contributions. Integrality implies
\begin{equation}\label{eq:plane-lower-integrality}
 c\in\tfrac1g\Z_{\ge0}:
 \quad\sum_i b_id_i=g\ (b_i\in\Z)\quad\Longrightarrow\quad
 cg=\sum_i b_i(cd_i)\in\Z.
\end{equation}
For $c=0$ the output weight is the sum of the input weights,
so Lemma~\ref{lem:fusion-leading-product} identifies the
entire product with \eqref{eq:end-product}. Unit labels alone
do not determine an output with $c>0$.

If $\Omega$ is chosen and the input coefficients have ordinary
cohomological degrees $r,s$, the output coefficient has degree
\begin{equation}\label{eq:projective-fusion-degree}
 r+s+\mu_\Omega(\vct)+\mu_\Omega(\wct)
      -\mu_\Omega(\boldsymbol\nu)
 =r+s+2c(d_\Sigma-n-1).
\end{equation}
For each support $I$, the stratum $D_I^\circ$ is smooth affine
of complex dimension $n-|I|$ and has the homotopy type of a
CW complex of real dimension at most $n-|I|$. In the Serre
spectral sequence of its normal-torus bundle, the base and
fiber degrees satisfy $p\le n-|I|$ and $q\le |I|$.
Consequently $H^j(S_I^\circ;\F)=0$ for $j>n$.
Thus a component whose degree in
\eqref{eq:projective-fusion-degree} exceeds $n$ is zero.
When $\Omega$ exists, this vanishing holds for the underlying
operation and hence for its parity periodization.

\begin{proposition}[A range with no lower terms]
\label{prop:many-divisors-fusion}
If $k\ge2n+1$, the contact identification $J_D$ is an
isomorphism of filtered unital BV algebras
\[
 (LH^*(X),\mu,\Delta)\cong
 (H^*_{\log}(\Pj^n,D),\cdot,\Delta_{\log}).
\]
\end{proposition}
\begin{proof}
The support of $\vct+\wct$ has at most $2n$ elements. Some coordinate
is zero, so subtracting $c\mathbf d$ with $c>0$ would make it
negative. Thus \eqref{eq:projective-fusion-output} allows only $c=0$.
Lemma~\ref{lem:fusion-leading-product} identifies that term
with logarithmic multiplication. If its face is empty there
is no realized output. Equation~\eqref{eq:hypersurface-bv}
gives compatibility with the BV operators.
\end{proof}

\subsection{Hyperplane coefficients and the full fusion algebra}\label{subsec:projective-hyperplanes}
Suppose in this subsection that every $D_i$ is a hyperplane.
Then $k\ge n+2$, and the SNC condition means general position. Put $T^a=(S^1)^a$.
\begin{proposition}
\label{prop:hyperplane-topology}
For a contact vector with support $I$ of cardinality $s\leq n$,
there is an isomorphism of graded algebras
\begin{equation}\label{eq:hyperplane-loop-type}
 H^*(\mathscr P_{\vct};\F)\cong
 \Lambda^*(\F^s)\otimes_{\F}\Lambda^{\le n-s}(\F^{k-s-1}).
\end{equation}
Here the generators have degree one, and $\Lambda^{\le d}(V)$ denotes
$\Lambda^*(V)/\bigoplus_{j>d}\Lambda^j(V)$.
Every such vector occurs and its contact locus is connected. In particular,
$H^*(X;\F)\cong\Lambda^{\le n}(\F^{k-1})$.
The coefficient Poincare polynomial is
\begin{equation}\label{eq:hyperplane-coefficient-polynomial}
 (1+q)^s\sum_{j=0}^{n-s}\binom{k-s-1}{j}q^j.
\end{equation}
\end{proposition}

\begin{proof}
The stratum $D_I^\circ$ is $\Pj^{n-s}$ minus $k-s$ hyperplanes in
general position.  Choose one omitted hyperplane.  Its defining section
trivializes $\mathcal O(1)$ on this open stratum, hence trivializes each
of the $s$ normal lines.  The normal-torus bundle is therefore a product
$T^s\times D_I^\circ$ and has the same cohomology as the contact locus by Proposition~\ref{prop:loop-sectors}.
For $s<n$, removing the chosen hyperplane identifies $D_I^\circ$ with
the complement of $k-s-1$ affine hyperplanes in $\C^{n-s}$ in general position.
By \cite[Definition~3.4 and Corollary~3.7]{DimcaArrangements}, its cohomology
is the exterior algebra on the classes of the defining linear units,
truncated above degree $n-s$: intersections of at most $n-s$ hyperplanes
have the expected codimension, and larger intersections are empty.
When $s=n$, the stratum is a point. These strata are nonempty and connected,
as are the normal tori, so the contact loci are nonempty and connected.
The K\"unneth formula gives \eqref{eq:hyperplane-loop-type} and
\eqref{eq:hyperplane-coefficient-polynomial}.
\end{proof}

Put $R=H^*(X;\F)$. By
Proposition~\ref{prop:hyperplane-topology}, this is the exterior
algebra on $k-1$ degree-one generators, truncated above degree $n$.
The maps
\[
 r_I:R\longrightarrow H^*(S_I^\circ;\F)
\]
are surjective. Indeed, choose an omitted hyperplane $D_{j_0}$ as
infinity. The classes of the units $G_j/G_{j_0}$ restrict to the
normal-circle classes for $j\in I$ and to the usual generators
of the hyperplane-complement base for the remaining indices.
These generate the two algebra factors in
\eqref{eq:hyperplane-loop-type}.

\begin{theorem}[The full hyperplane fusion ring]
\label{thm:hyperplane-ring}
In the setup of Section~\ref{sec:hypersurface-complements},
assume $d_i=1$ for every $i$. Then the fusion product is the
logarithmic common-face product.
Set
\[
 \mathcal S_{n,k}=\F[t_1,\ldots,t_k]/
 (t_{i_1}\cdots t_{i_{n+1}}:i_1,\ldots,i_{n+1}\text{ distinct}),
 \qquad t_I=\prod_{i\in I}t_i.
\]
Then the fusion ring has the presentation
\begin{equation}\label{eq:hyperplane-ring-presentation}
 LH^*(X)=
 \frac{\mathcal S_{n,k}\otimes_{\F} R}
 {\bigl(t_I\alpha: |I|\le n,\ \alpha\in\ker r_I\bigr)}.
\end{equation}
The unit is the image of $1\in R$.
\end{theorem}
\begin{proof}
Let $t_i$ be the constant coefficient class of contact $\mathbf e_i$.
We use associativity of fusion from
Theorem~\ref{thm:general-associativity} to determine products
whose lower terms are not excluded by unit valuations alone.

If the union of two input supports is a proper subset of
$\{1,\ldots,k\}$, equation~\eqref{eq:projective-fusion-output}
forces $c=0$. Their entire product is therefore common-face
multiplication; it is zero when that face is empty.
It follows by induction inside a realized face that
$\prod_i t_i^{v_i}$ is its constant coefficient class.
Multiplication by $\alpha\in R$ is $r_I(\alpha)$ on that summand,
by the same output test and the leading-product theorem.
Surjectivity of $r_I$ therefore shows that $R$ and the $t_i$
generate the full fusion ring.

For $n+1$ distinct indices, first multiply $n$ of their $t_i$ and
then the remaining one. Their combined support still omits a
hyperplane because $k\ge n+2$. Its size is $n+1$, so no output
contact is realized and this product is zero. Associativity and
graded commutativity allow these even contact generators to be
grouped inside every larger monomial without additional signs. The relations $t_I\alpha=0$ for $\alpha\in\ker r_I$
follow from the preceding restriction calculation. We have thus
constructed a surjective algebra map from the right side of
\eqref{eq:hyperplane-ring-presentation} to $LH$.

In multidegree $\vct$ with support $I$, that quotient has precisely
the coefficient $R/\ker r_I=H^*(S_I^\circ;\F)$. Kernels from
smaller faces are contained in $\ker r_I$, so impose no further
relations. The map is the known additive isomorphism on each
summand and hence is an algebra isomorphism.
\end{proof}

\emph{The BV operator in the ring presentation.}
Fix $D_k$ as reference hyperplane, and let
$\beta_j\in H^1(X;\F)$ be the class of the phase of the unit
$G_j/G_k$, for $1\le j<k$. For a contact monomial $t^{\vct}$,
the operator in \eqref{eq:hyperplane-ring-presentation} is
\begin{equation}\label{eq:hyperplane-bv}
 \Delta\left(t^{\vct}\beta_{j_1}\cdots\beta_{j_r}\right)
 =\sum_{\ell=1}^r(-1)^{\ell-1}(v_{j_\ell}-v_k)\,
   t^{\vct}\beta_{j_1}\cdots\widehat{\beta_{j_\ell}}
                         \cdots\beta_{j_r}.
\end{equation}
In particular $\Delta(t^{\vct})=0$ and
$\Delta|_{H^*(X)}=0$. To verify the formula, restrict to the
coefficient space $S_{I(\vct)}^\circ$. Source rotation multiplies
the leading coefficient of $G_j/G_k$ by
$e^{i(v_j-v_k)\phi}$, so its circle operator sends
$r_I\beta_j$ to $(v_j-v_k)1$. The circle operator on coefficient
cohomology is a graded derivation of degree $-1$. Passing it
through the preceding $\ell-1$ degree-one classes produces
$(-1)^{\ell-1}$, giving \eqref{eq:hyperplane-bv}.
This geometric calculation proves that the formula respects the
relations in \eqref{eq:hyperplane-ring-presentation}; contraction
is applied after restriction to the contact coefficient, not to
the truncated exterior algebra before restriction.

\subsection{Plane coefficients and the additional fusion range}
\label{subsec:example-24-1}
Suppose now that $n=2$. Write $S_i=S_{\{i\}}^\circ$ and, for
$p\in D_i\cap D_j$, write
$T_{ij,p}=S_{\{i,j\}}^\circ|_p\cong(S^1)^2$.
The following theorem gives the contact decomposition and
BV operator explicitly. A volume grading exists precisely when
$g\mid3$, by \eqref{eq:hypersurface-grading-condition}; its
contact shifts are the $\mu_\Omega(\vct)$ defined above.

\begin{theorem}[Loop cohomology of a plane complement]
\label{thm:plane-complements}
In the setup of Section~\ref{sec:hypersurface-complements},
assume $n=2$.
\begin{enumerate}[(i),wide=0pt]
\item The contact decomposition of
Proposition~\ref{prop:hypersurface-separation} is
\begin{equation}\label{eq:plane-additive-decomposition}
 \begin{split}
 LH^*(X)\cong{}&H^*(X;\F)
 \oplus\bigoplus_{i=1}^k\bigoplus_{m\ge1}
           H^*(S_i;\F)t^{m\mathbf e_i}\\
 &\oplus\bigoplus_{i<j}\bigoplus_{p\in D_i\cap D_j}
       \bigoplus_{a,b\ge1}
           H^*(T_{ij,p};\F)t^{a\mathbf e_i+b\mathbf e_j}.
 \end{split}
\end{equation}

\item The full loop BV operator preserves every summand of
\eqref{eq:plane-additive-decomposition} and vanishes on $H^*(X;\F)$.
Choose a normal-circle trivialization of $S_i$ and let $\theta_i$
be its positive angular class, placed before classes from $D_i^\circ$.
For $\beta,\gamma\in H^*(D_i^\circ;\F)$ and
$\alpha\in\Lambda(\theta_i,\theta_j)=H^*(T_{ij,p};\F)$, it is
\begin{equation}\label{eq:plane-bv}
 \begin{aligned}
 \Delta((\beta+\theta_i\gamma)t^{m\mathbf e_i})
    &=m\gamma t^{m\mathbf e_i},\\
 \Delta(\alpha t^{a\mathbf e_i+b\mathbf e_j}|_p)
    &=(a\iota_i+b\iota_j)(\alpha)
                         t^{a\mathbf e_i+b\mathbf e_j}|_p.
 \end{aligned}
\end{equation}
Here $\iota_i,\iota_j$ are the degree-minus-one derivations
with $\iota_i(\theta_i)=\iota_j(\theta_j)=1$ and
$\iota_i(\theta_j)=\iota_j(\theta_i)=0$; in particular
$(a\iota_i+b\iota_j)(\theta_i\theta_j)=a\theta_j-b\theta_i$.
\end{enumerate}
\end{theorem}
\begin{proof}
Part (i) is the contact decomposition already proved, since
the nonzero realized contacts have support of size one or two.

For (ii), the smooth affine curve $D_i^\circ$ has the homotopy
type of a one-dimensional CW complex, so
$H^2(D_i^\circ;\Z)=0$. Its normal line bundle is therefore
topologically trivial, giving the normal-circle trivialization
used in the statement. Specialize \eqref{eq:hypersurface-bv}.
Rotation on
$S_i$ has weight
$m$ on its normal circle, and on $T_{ij,p}$ it is
$\lambda\mapsto(\lambda^a,\lambda^b)$. Slant product with
$[S^1]$, with the acting circle first, gives
\eqref{eq:plane-bv}; the zero-contact assertion is already
part of Proposition~\ref{prop:hypersurface-separation}.
\end{proof}

Among the complements of Theorem~\ref{thm:plane-complements}(i),
two lines and a conic form a borderline case for determining
the fusion product from the additive decomposition:
$(d_1,d_2,d_3)=(1,1,2)$ has $d_\Sigma=4$, and the volume
grading still allows lower terms of ordinary coefficient degree
$2$. The following lemma proves that these terms vanish.

\begin{lemma}[Two lines and a conic]\label{lem:plane-two-lines-conic}
Suppose $(d_1,d_2,d_3)=(1,1,2)$ and choose a volume grading.
Let $\diamond$ be a filtered unital graded-commutative associative
product on the graded vector space $H^*_{\log}(\Pj^2,D)$.
Assume that it adds the unit labels $W\vct$ and that its
associated-graded product is \eqref{eq:end-product}. Then
$\diamond$ equals the logarithmic product.
\end{lemma}
\begin{proof}
The output and degree tests
\eqref{eq:projective-fusion-output}--\eqref{eq:projective-fusion-degree}
use only the hypotheses imposed on $\diamond$, so apply to
this product. Here $g=1$ and $d_\Sigma=4$. A lower output
has $c\in\Z_{>0}$ and ordinary coefficient degree $r+s+2c$.
The degree bound $H^{>2}(S_I^\circ;\F)=0$ therefore leaves
only $r=s=0$ and $c=1$. If the input supports omit a component,
the output test forces $c=0$ and the product is logarithmic.

Let $t_i$ be the constant class of contact $\mathbf e_i$, and
write
\[
 a_p=1_p t^{\mathbf e_1+\mathbf e_3}
       \quad(p\in D_1\cap D_3),\qquad
 b_q=1_q t^{\mathbf e_2+\mathbf e_3}
       \quad(q\in D_2\cap D_3).
\]
By B\'ezout and the SNC hypothesis,
$|D_1\cap D_2|=1$ and
$|D_1\cap D_3|=|D_2\cap D_3|=2$.
Thus there are two $a_p$ and two $b_q$.
Multiplication within a single face expresses every
positive-contact constant as
\[
 t_i^m,\qquad t_1^a t_2^b,\qquad
 a_p t_1^{a-1}t_3^{b-1},\qquad
 b_q t_2^{a-1}t_3^{b-1}\qquad(a,b,m\ge1),
\]
where juxtaposition denotes $\diamond$. The unique
$D_1\cap D_2$ crossing accounts for the second expression.
A product $a_p\diamond t_2$ or $b_q\diamond t_1$ would need
$c=1/2$ to have a nonnegative realized output, contradicting
$c\in\Z$. Therefore
\[
 a_p\diamond t_2=b_q\diamond t_1=0,
 \qquad t_1t_2t_3=0.
\]
The remaining primitive pair has
$a_p\diamond b_q=\lambda_{pq}\in H^2(X;\F)$.
Associativity, and evenness of all these constant classes, give
\[
 \lambda_{pq}\diamond t_1=\lambda_{pq}\diamond t_2=0.
\]
Multiplication by $t_i$ on contact zero is exactly restriction
$\mathrm{res}_i:H^2(X;\F)\to H^2(S_i;\F)$ followed by multiplication
by $t_i$. Thus $\mathrm{res}_1\lambda_{pq}=\mathrm{res}_2\lambda_{pq}=0$.

We verify that $(\mathrm{res}_1,\mathrm{res}_2)$ is injective.
The dual graph $G$ has three vertices and five edges: one
joining $1,2$ and two each joining $1,3$ and $2,3$.
Since $D_1$ is a line, the map
$H_2(D;\F)\to H_2(\Pj^2;\F)$ is surjective.
The homology sequence of $(\Pj^2,D)$ and
Poincar\'e--Lefschetz duality therefore give the first two
isomorphisms below. Since every component of $D$ is rational,
the remaining identification is with the cycle space of $G$:
\[
 H^2(X;\F)\cong H_2(\Pj^2,D;\F)\cong H_1(D;\F)
 \cong H_1(G;\F)
 =\ker\bigl(\partial:\F^{E(G)}\to\F^{V(G)}\bigr).
\]
Under this identification, restriction to a crossing torus,
followed by integration over that torus, is its signed edge
coefficient. Indeed, in a small ball $B_p$ about a node,
$D\cap B_p=\{uv=0\}\cap B_p$ consists of two disks and
\[
 H_1(D\cap B_p,D\cap\partial B_p;\F)\cong\F.
\]
Excision sends a graph cycle to its coefficient on this local
edge; naturality of Poincar\'e--Lefschetz duality identifies
this map with restriction to the linking torus
$\{|u|=|v|=\varepsilon\}$. Since every edge meets vertex $1$
or $2$, vanishing of $\mathrm{res}_1\lambda$ and $\mathrm{res}_2\lambda$ makes every
edge coefficient zero. Hence $\lambda=0$, proving
$a_p\diamond b_q=0$.

Products of the displayed monomials now give logarithmic
multiplication: products within one face are already known;
incompatible crossing labels give zero; a union of three
supports contains one of the zero factors
$a_pb_q$, $a_pt_2$, $b_qt_1$, or $t_1t_2t_3$.
Associativity permits these factors to be grouped. Finally,
any input of positive ordinary coefficient degree excludes
lower terms by the first paragraph. This proves the assertion
for all classes.
\end{proof}

\begin{proposition}[The additional plane fusion range]
\label{prop:plane-fusion-range}
If $n=2$ and $g\mid3$, then $J_D$ is an isomorphism of
filtered unital BV algebras.
\end{proposition}
\begin{proof}
Choose a volume grading, which exists by
\eqref{eq:hypersurface-grading-condition}. The assumptions of
this section give $d_\Sigma\ge4$. A lower component has
$c\in\frac1g\Z_{>0}$ by \eqref{eq:plane-lower-integrality};
its ordinary coefficient degree is
$r+s+2c(d_\Sigma-3)$ by
\eqref{eq:projective-fusion-degree}, and must be at most two.
Since $g\mid3$, either $g=1$ or $g=3$. If $g=1$ and
$d_\Sigma\ge5$, then $2c(d_\Sigma-3)\ge4$. If $g=3$, then
$d_i\ge3$ for every $i$ and $k\ge3$, so
$2c(d_\Sigma-3)\ge2(9-3)/3=4$. Thus lower components vanish
in both cases.

If $d_\Sigma=4$, the degrees are either four ones or
$(1,1,2)$ up to permutation. The first case is
Theorem~\ref{thm:hyperplane-ring}. In the second case, transport
fusion through $J_D$. The resulting product satisfies the
hypotheses of Lemma~\ref{lem:plane-two-lines-conic} by
Theorem~\ref{thm:general-associativity},
Proposition~\ref{prop:fusion-unit-valuations}, and
Lemma~\ref{lem:fusion-leading-product}; hence it is logarithmic.
Finally $J_D$ commutes with the BV operators by
Proposition~\ref{prop:hypersurface-separation}. The integer-graded
calculation also determines its parity periodization.
\end{proof}

\subsection{Symplectic ring comparisons}
\label{subsec:plane-sh-rings}
The comparisons below use the loop calculations just established
and the logarithmic PSS theorems of \cite{GP}.

\begin{theorem}[Additive and multiplicative ranges]
\label{thm:hypersurface-comparison}
Logarithmic PSS gives a filtered additive isomorphism
\[
 \operatorname{PSS}_{\log}\circ J_D:
 LH^*(X)\xrightarrow{\sim}SH^*(X).
\]
It is multiplicative for $k\ge2n+1$ and, independently of
this bound, when all $D_i$ are hyperplanes. The assertions
respect every chosen volume grading.
\end{theorem}
\begin{proof}
Let a nonconstant relative holomorphic sphere be generically contained
in $D_I^\circ$, and put $s=|I|$.  Necessarily $s<n$.  Its projective
degree $e$ is positive, so it intersects each of the $k-s$ remaining
divisors with positive intersection number $ed_j$.  At one domain point
the image can lie on at most $n-s$ of those remaining divisors.  Thus
the number $r$ of distinct boundary preimages satisfies
\begin{equation}\label{eq:hypersurface-point-bound}
 r\geq\left\lceil\frac{k-s}{n-s}\right\rceil.
\end{equation}
The argument also applies to an adapted almost complex structure:
positive symplectic area gives $e>0$, and positivity of intersections
gives the required meetings with every remaining divisor.
For $k\geq n+1$ the bound is at least two.  For $k\geq2n+1$ it is at
least three, since $k-s>2(n-s)$.  The logarithmic PSS theorems
\cite[Theorems~5.10 and~5.12]{GP} give the stated filtered additive
and ring identifications. Composing with $J_D$ gives the
additive comparison. For $k\ge2n+1$,
Proposition~\ref{prop:many-divisors-fusion} makes $J_D$
multiplicative, proving the ring comparison.

If all $D_i$ are hyperplanes, $k\ge n+2$ and
\cite[Lemma~5.18]{GP} proves multiplicativity of logarithmic
PSS. Theorem~\ref{thm:hyperplane-ring} computes the loop
product as the same common-face product. Hence
\begin{equation}\label{eq:hyperplane-sh-ring}
 (LH^*(X),\mu)\cong SH^*(X)
\end{equation}
as filtered graded rings, including the generalized pair of
pants $k=n+2$.
\end{proof}

Four lines form the borderline plane hyperplane case: the
bound \eqref{eq:hypersurface-point-bound} permits a two-pointed
relative sphere. Indeed, the line through $D_1\cap D_2$ and
$D_3\cap D_4$ meets $D$ at only these two points. Thus the
comparison for four lines uses \cite[Lemma~5.18]{GP}, rather
than the exclusion of two-pointed relative spheres.

\section{The hyperbola complement}
\label{sec:hyperbola-comparison}

Let $X=\C^2\setminus\{xy=1\}$ and let
$f:X\to\C^*$ be the regular function $f(x,y)=xy-1$.
It is a unit in $\mathcal O(X)$, with regular reciprocal
$f^{-1}$. Introducing an invertible coordinate for this
function gives the presentation
\begin{equation}\label{eq:hyperbola-X}
 X=\Spec\C[x,y,f,f^{-1}]/(xy-1-f),\qquad
 \Omega=\frac{dx\wedge dy}{f}.
\end{equation}
Theorem~\ref{thm:hyperbola-ring} determines the integer-graded
ring $LH^*_\Omega(X)$ and its BV operator.
Proposition~\ref{prop:hyperbola-degeneration} identifies the
associated graded rings of the loop and Floer contact filtrations
with $H^*_{\log}(M,D)$ for the line--conic compactification
$(M,D)$ below. The full loop ring includes lower-weight fusion terms.

\begin{theorem}
\label{thm:hyperbola-ring}
The degree-zero ring has the presentation
$R_X=LH^0_\Omega(X)\cong\F[u,v,w,w^{-1}]/(uv-1-w)$.
There is an isomorphism of
integer-graded unital BV algebras
\begin{equation}\label{eq:hyperbola-full-ring}
 \bigl(LH^*_\Omega(X),\mu,\Delta\bigr)
 \cong \bigl(R_X\otimes\Lambda(\xi_u,\xi_v),\wedge,
                \Delta_\eta\bigr),
\end{equation}
Here $|\xi_u|=|\xi_v|=1$,  \(\Lambda\) denotes the exterior algebra over \(\F\), and $\eta=du\wedge dv/w$.
Identify $R_X\otimes\Lambda(\xi_u,\xi_v)$ with
$\bigwedge_{R_X}^{\bullet}\operatorname{Der}_{\F}(R_X)$
by $\xi_u\mapsto\partial_u$ and $\xi_v\mapsto\partial_v$.
Under this identification, the isomorphism
\eqref{eq:hyperbola-full-ring} is the canonical map
$\sigma_{\F}$ of Remark~\ref{globpv:remark-rational}, defined by
the iterated-bracket formula \eqref{globpv:eq-sigma}.
The operator $\Delta_\eta$ is the divergence BV operator
of Example~\ref{ex:volume-divergence-bv}.
\end{theorem}

The line--conic compactification of Subsection~\ref{subsec:example-35-1}
is positive with maximal boundary, and the presentation of $R_X$ above
identifies $\check X$ with the smooth open set $\{uv\ne1\}\subset\A^2_{\C}$.
Thus Theorem~\ref{thm:surface-loop-floer} gives an isomorphism of
unital graded $\C$-algebras
\[
 LH^*_\Omega(X)\otimes_{\F}\C\xrightarrow{\sim}SH^*_\Omega(X;\C).
\]
In degree zero, Pascaleff's computation
\cite[\S7.2, (139) and (144)]{Pascaleff} gives
$SH^0_\Omega(X;\F)\cong\F[u,v,(uv-1)^{-1}]$,
which is the presentation of $LH^0_\Omega(X)$ in
Theorem~\ref{thm:hyperbola-ring} after substituting $w=uv-1$.

\subsection{Geometry, contact loci, and degeneration}\label{subsec:example-35-1}

Let $M=\Pj^2$ have homogeneous coordinates $[Z_0:Z_1:Z_2]$.
The embedding $(x,y)\mapsto[1:x:y]$ identifies $X$ with
$M\setminus D$, where
\[
 L=\{Z_0=0\},\qquad
 Q=\{Z_1Z_2=Z_0^2\},\qquad D=L\cup Q.
\]
The line and conic meet transversely at two points.
Since $L+Q$ has degree three, $D=-K_M$ and
$\operatorname{div}_M\Omega=-D$. Take $A=L+Q$.
For contact $(\ell,q)$ the weight is $\ell+q$, and the
valuation of the unit $f$ is
\begin{equation}\label{eq:hyperbola-unit-weight}
 k=q-2\ell.
\end{equation}
Here $k$ is a grading supplied by a global algebraic unit,
as in Proposition~\ref{lp:prop:unit-valuation}.

The Gysin sequence for the smooth hyperbola $\C^*\subset\C^2$
gives
\[
 H^j(X;\F)=
 \begin{cases}\F,&j=0,1,2,\\0,&\text{otherwise}.\end{cases}
\]
The degree-one generator is the pullback of the angular class
by $f:X\to\C^*$, so its square is zero. All other products
of positive-degree classes vanish by degree.
Both $L^\circ$ and $Q^\circ$ are $\C^*$, their normal line
bundles are trivial on these open curves, and their normal
circle bundles have homotopy type $T^2$. Each of the two
crossings also has coefficient space $T^2$. All contact
shifts are zero.

\begin{proposition}\label{prop:hyperbola-degeneration}
For this compactification both the loop and
Ganatra--Pomerleano spectral sequences degenerate at $E_1$
in every total degree. Their common associated graded ring
is $H^*_{\log}(M,D)$. In particular
$LH^*_\Omega(X)$ is noncanonically isomorphic to this vector space
as a filtered graded vector space, with its weight filtration. Its
underlying vector space is
\begin{equation}\label{eq:hyperbola-first-page}
 H^*(X)\ \oplus\
 \bigoplus_{\ell\ge1}H^*(T^2)t_L^\ell\
 \oplus\bigoplus_{q\ge1}H^*(T^2)t_Q^q\
 \oplus\bigoplus_{\ell,q\ge1}
       \bigl(H^*(T^2)\oplus H^*(T^2)\bigr)t_L^\ell t_Q^q .
\end{equation}
\end{proposition}
\begin{proof}
Let $H$ be the hyperplane class on $\Pj^2$. Since
$L\sim H$ and $Q\sim2H$, the loop assertion follows
from Theorem~\ref{thm:intro-proportional}; the
Floer assertion follows from Corollary~1.2 of \cite{GP}.
Theorem~\ref{thm:summary} and \cite[Theorem~1.1]{GP}
identify the two multiplicative first pages with
$H^*_{\log}(M,D)$, hence identify their associated graded rings. The cohomology and normal-torus calculations above
give \eqref{eq:hyperbola-first-page}. Over the field $\F$,
choose splittings of the successive filtration quotients to obtain
the asserted additive isomorphism. 
\end{proof}

For a field extension $\kappa/\C$ and a loop
$\gamma:\Spec\kappa(\!(t)\!)\to X$, write
$x(t)=\gamma^*x$, $y(t)=\gamma^*y$, and $f(t)=\gamma^*f$.
In coefficient formulas we abbreviate these Laurent series by
$x,y,f$; in particular, $\ord_t f$ means
$\ord_t(\gamma^*f)$, and similarly for the other functions.
Here $\ord_t$ is the least exponent with nonzero coefficient,
with $\ord_t0=+\infty$.
Since $f$ is invertible on $X$, put
$k=\ord_t f\in\Z$ and $h(k)=\min(0,k)$.

The coefficient models below are reduced and carry the inherited
volume indices described in
\emph{Reduction of finite models} in
Subsection~\ref{subsec:volume-index-proofs}.
For specialization, reduction is applied to the total coefficient
family over the parameter line, with the relative index restricted
as in Definition~\ref{def:relative-volume-index}(iii).
The formal source rings and their finite jet algebras are unchanged.

\begin{lemma}
\label{lem:hyperbola-coefficient-models}
\begin{enumerate}[(i),wide=0pt]
\item For $k\in\Z\setminus\{0\}$ and $n\in\Z$, let
$\mathscr L_{k,n}X$ be the reduced locus in $\mathscr LX$
defined by $\ord_t f=k$ and $\ord_t x=n$.
We call it the \emph{$(k,n)$-sector}. Its loops have the unique expression
\begin{equation}\label{eq:hyperbola-nonzero-sectors}
 f=t^k F(t),\qquad x=t^n A(t),\qquad
 y=t^{-n}(1+t^kF(t))A(t)^{-1},
 \qquad A(0)F(0)\ne0,
\end{equation}
where $A,F\in\kappa\llbracket t\rrbracket$ on
$\kappa$-valued loops.
Retaining their coefficients through degree $N$ gives the
finite coefficient models of $\mathscr L_{k,n}X$
\[
 Y_N=(\Gm)^2\times\A^{2N},\qquad
 \delta_{Y_N}=2N+2,\qquad N\ge0,
\]
with $k,n$ fixed and suppressed from the notation.
Increasing $N$ gives normalized affine pullback, and pole
enlargement preserves the reduced model and its index.
The $(k,n)$-sector has normalized cohomology
$\Lambda(\alpha,\beta)$, where
$\alpha$ and $\beta$ are the pullbacks of the angular class
in $H^1(\C^*;\F)$ under the leading-coefficient maps
$\gamma\mapsto A(0)$ and $\gamma\mapsto F(0)$, respectively.

\item The irreducible components of the reduced locus
$\ord_t f=0$ are indexed by $n\in\Z$. The component indexed
by $n$ is the reduced closure of the locus
$\ord_t x=n$, $\ord_t y=-n$, and is described by
\[
 x=t^n A(t),\qquad y=t^{-n}B(t),\qquad
 A,B\in\kappa\llbracket t\rrbracket,\qquad A(0)B(0)\ne1,
 \qquad n\in\Z.
\]
The component $n=0$ is the arc stage.
The contact of $\mathscr L_{k,n}X$ for $k\ne0$ and the generic
contact of the component indexed by $n$ for $k=0$ are given by
\begin{equation}\label{eq:hyperbola-contact-dictionary}
 \ell(k,n)=\max\{0,-n,n-h(k)\},\qquad
 q(k,n)=k+2\ell(k,n).
\end{equation}
The indices $(k,n)$ parametrize all contact components: for
a fixed contact with $\ell,q>0$, there are two choices of
$n$, whereas a single-divisor contact gives one.
\end{enumerate}
\end{lemma}
\begin{proof}
\emph{(i)}
When $k\ne0$, the series $1+f$ is nonzero, so the relation
$xy=1+f$ forces $x$ and $y$ to be nonzero.
Writing $f=t^kF$ and $x=t^nA$ with unit series $F,A$ gives the unique
expression \eqref{eq:hyperbola-nonzero-sectors}.
Put $A_j=[t^j]A(t)$ and $F_j=[t^j]F(t)$ for $j\ge0$.
These are regular coefficient coordinates: lower coefficients
vanish on every residue field and hence on the reduced scheme,
while $A_0,F_0$ are inverted. Moreover,
\[
 \ord_t(1+f)=h(k),\qquad
 \ord_t y=h(k)-n,\qquad \ord_t x^{-1}=-n.
\]
Thus bounds on $x,y,f,f^{-1}$ and on $x,x^{-1},f,f^{-1}$
are cofinal. Polynomial changes of affine generators compare
them with the original pole bounds. Each coefficient formula
in either direction uses finitely many coefficients, so the
models $Y_N$ in (i) retain $A_0,F_0\ne0$ and
$A_1,\ldots,A_N,F_1,\ldots,F_N$.
Proposition~\ref{prop:lcy-equidimensional-factors}(ii) gives
$\delta_{Y_N}=\dim Y_N=2N+2$.
Proposition~\ref{prop:global-finite-presentations}(i),(iii)
identifies the normalized cohomology and transition maps with
those of $Y_N$. Since its affine factor contributes no cohomology,
the result is $H^*((\C^*)^2;\F)=\Lambda(\alpha,\beta)$.
This proves (i).

\emph{(ii)}
At $k=0$, one has $\ord_t(xy)\ge0$, so one can choose
$n\in\Z$ with $-\ord_t y\le n\le\ord_t x$.
Then $A=t^{-n}x$ and $B=t^ny$ are power series and
$A(0)B(0)\ne1$, since $f(0)\ne0$. Thus every loop lies
in one of the models in (ii). On a stage where
$x,y\in t^{-P}\kappa\llbracket t\rrbracket$, with $P\ge0$,
one can choose $-P\le n\le P$, so only finitely many models
are needed. For fixed $n$, the model is
closed in the $k=0$ locus, being defined by the vanishing
of the coefficients of $x$ below degree $n$ and of $y$
below degree $-n$. Its coefficient ring is the domain
\[
 \C[A_j,B_j:j\ge0][(A_0B_0-1)^{-1}],
 \qquad A_j=[t^j]A,\qquad B_j=[t^j]B.
\]
The additional condition $A_0B_0\ne0$ defines a dense open
with orders $n,-n$, which belongs to no model with a
different index. Hence these models are precisely the
irreducible components. 
The extension to $\Pj^2$ is represented by the homogeneous arc
$[t^\ell:t^\ell x:t^\ell y]$, where
$\ell=\max\{0,-\ord_t x,-\ord_t y\}$. Its contact with
$L$ is $\ell$, and substitution in the equation of $Q$ gives
$t^{2\ell}(xy-1)=t^{2\ell}f$, of order $2\ell+k$.
For $k\ne0$, the orders on $\mathscr L_{k,n}X$ computed in (i) give
\eqref{eq:hyperbola-contact-dictionary}; for $k=0$, use the
orders $n,-n$ on the dense opens just described.
For a fixed contact $(\ell,q)$, one has $k=q-2\ell$, and
the possible indices are $n=-\ell$ and $n=h(k)+\ell$.
They are distinct exactly when $\ell,q>0$, and coincide
for a single-divisor contact or the arc contact. This proves
the indexing assertion and completes (ii).
\end{proof}

\emph{Global finite models.}
For integers $M\ge0$ and $N\ge M$, retain $x_j,y_j$ for
$-M\le j\le N$ and give
\[
 \begin{split}
 E_{M,N}:={}&
 \left(\bigcup_{r=-M}^{M}
   \{x_i=0\ (i<r),\quad y_j=0\ (j<-r)\}\right)\\
 &\quad\cap
 \left\{1-\sum_{i=-M}^{M}x_i y_{-i}\ne0\right\}
 \ \subset\ \A^{2(N+M+1)}
 \end{split}
\]
its reduced scheme structure. Every coefficient index here lies
between $-M$ and $N$. By Lemma~\ref{lem:hyperbola-coefficient-models}(ii),
the reduced locus $\ord_t f=0$ with
$x,y\in t^{-M}\C\llbracket t\rrbracket$ is
\[
 E_{M,N}\times\A^{\{x_j,y_j:j>N\}}.
\]
Indeed, $\ord_t(xy)\ge0$ gives the displayed union, and the
open condition is $f(0)\ne0$. This includes $x=0$ or $y=0$.
Each component of $E_{M,N}$ has dimension $2N+2$, since it
imposes $2M$ coordinate equations.

Fix $P\in\Z_{\ge0}$ and work on the reduced pole stage
$(F_P\mathscr LX)_{\mathrm{red}}$. Put
\[
 I_P=\{(k,n)\in\Z^2:k\ne0,\ \ell(k,n)+q(k,n)\le P\}.
\]
By \eqref{eq:hyperbola-contact-dictionary}, every pair in
$I_P$ satisfies $-2P\le k\le P$ and $-P\le n\le P$,
so $I_P$ is finite.

On this pole stage, the loci with fixed $\ord_t f=k$ are
the principal opens
\[
 \{\ord_t f=k\}
 =
 \bigl\{([t^k]f)([t^{-k}]f^{-1})\ne0\bigr\}.
\]
Indeed, nonvanishing of these two coefficients gives
$\ord_t f\le k$ and $\ord_t f^{-1}\le-k$, forcing equality.
Since $k=q-2\ell$ and $\ell+q\le P$, only
$-2P\le k\le P$ occur. These finitely many disjoint opens
cover the stage, so each is also closed.

For $(k,n)\in I_P$, the entire sector $\mathscr L_{k,n}X$
lies in this pole stage. Within $\{\ord_t f=k\}$, it is
the principal open
\[
 \mathscr L_{k,n}X
 =
 \bigl\{([t^n]x)([t^{h(k)-n}]y)\ne0\bigr\}.
\]
Here the nonvanishing conditions and
$\ord_t x+\ord_t y=h(k)$ force
$\ord_t x=n$ and $\ord_t y=h(k)-n$.
For fixed $k\ne0$, these sectors form a finite open
partition of $\{\ord_t f=k\}$. Thus every
$\mathscr L_{k,n}X$ with $(k,n)\in I_P$ is open and closed
in $(F_P\mathscr LX)_{\mathrm{red}}$.

On $\{\ord_t f=0\}$, the contact weight is
\[
 3\max\{0,-\ord_t x,-\ord_t y\}.
\]
Consequently its intersection with the pole stage is
defined by
$\ord_t x,\ord_t y\ge-\lfloor P/3\rfloor$.
The preceding coefficient models therefore give
\begin{equation}\label{eq:hyperbola-global-factors}
 Z_{P,N}
 =E_{\lfloor P/3\rfloor,N}\;\sqcup\!
   \coprod_{(k,n)\in I_P}Y_N,
 \qquad N\ge\lfloor P/3\rfloor,
 \qquad \delta_{Z_{P,N}}=2N+2,
\end{equation}
with
\[
 (F_P\mathscr LX)_{\mathrm{red}}
       \cong\varprojlim_N Z_{P,N}.
\]
Each copy of $Y_N$ represents the sector indexed by its
$(k,n)$. The volume index follows from
Proposition~\ref{prop:lcy-equidimensional-factors}(ii).

The projection $Z_{P,N+1}\to Z_{P,N}$ forgets the two
highest coefficients and is an affine bundle of relative
dimension two. For $P\le P'$ and a common sufficiently
large cutoff $N\ge\lfloor P'/3\rfloor$, the pole inclusion
is the identity on the copies of $Y_N$ indexed by $I_P$.
On the $k=0$ part it is
\[
 E_{\lfloor P/3\rfloor,N}
 \hookrightarrow E_{\lfloor P'/3\rfloor,N},
 \qquad
 x_j=y_j=0
 \quad
 \bigl(-\lfloor P'/3\rfloor\le j<-\lfloor P/3\rfloor\bigr).
\]
These closed inclusions form Cartesian squares with the
coefficient projections. Hence
Proposition~\ref{prop:global-finite-presentations} computes
the reduced pole-stage diagram with its inherited volume
indices from the models \eqref{eq:hyperbola-global-factors}.

Let $\theta_{k,n}$ be the theta class constructed in
\eqref{eq:lcy-theta-colimit-class} for the contact component
indexed by $(k,n)$, with contact orders given by
\eqref{eq:hyperbola-contact-dictionary}. In particular,
$\theta_{0,0}=1$.

\subsection{The degree-zero collision calculation}

Set
\[
 U=\theta_{1,0},\qquad V=\theta_{-1,0},\qquad
 W=\theta_{0,1},\qquad W_-=\theta_{0,-1}.
\]
Since a contact $(\ell,q)$ has weight $\ell+q$,
\eqref{eq:hyperbola-contact-dictionary} gives
\[
 w(U)=1,\qquad w(V)=2,\qquad w(W)=w(W_-)=3.
\]
We compute fusion by applying normalized specialization to
exterior products of the finite-factor representatives of the
canonical pole-stage lifts of these classes. By
\eqref{eq:lcy-theta-stage-class}, these representatives are the
normalized fundamental classes of the corresponding reduced
closures, in the sense of
Definition~\ref{def:finite-factor-representatives}.

\begin{lemma}
\label{lem:hyperbola-exchange}
The fusion product $\mu$ on $R_X=LH^0_\Omega(X)$ satisfies
\begin{equation}\label{eq:hyperbola-exchange}
 WW_-=1,\qquad UV=1+W,\qquad
 \theta_{k,n}=
 \begin{cases}U^kW^n,&k\ge0,\\V^{-k}W^n,&k<0.\end{cases}
\end{equation}
The homomorphism
$\F[u,v,w,w^{-1}]/(uv-1-w)\to R_X$,
$u\mapsto U$, $v\mapsto V$, $w\mapsto W$, is an isomorphism.
\end{lemma}
\begin{proof}
\emph{Step 1: construct the normalized Laurent shifts.}
For $r\in\Z$, define
\[
 T_r(x(t),y(t),f(t))=(t^rx(t),t^{-r}y(t),f(t)).
\]
This preserves $xy=1+f$ and has inverse $T_{-r}$.
In the models of
Lemma~\ref{lem:hyperbola-coefficient-models}(i), it leaves
$A,F$ unchanged and gives
\[
 T_r:\mathscr L_{k,n}X\xrightarrow{\sim}\mathscr L_{k,n+r}X
 \qquad(k\ne0,\ n\in\Z).
\]
For $k=0$, it acts on the component models of part~(ii) by
\[
 T_r(t^nA,t^{-n}B,AB-1)
   =(t^{n+r}A,t^{-(n+r)}B,AB-1),
 \qquad A(0)B(0)\ne1,
\]
identifying the component indexed by $n$ with that indexed
by $n+r$.

With the source parameter $t$ held fixed,
\[
 T_r^*\Omega
 =\frac{d(t^rx)\wedge d(t^{-r}y)}{f}
 =\frac{dx\wedge dy}{f}=\Omega.
\]
To realize $T_r$ on finite factors, use the identity on the
coefficients $A,F$ of each nonzero-$k$ sector. For the $k=0$
part, let $E_{M;N_x,N_y}$ denote the reduced model defining
$E_{M,N}$ in \emph{Global finite models} in
Subsection~\ref{subsec:example-35-1},
with $x_j$ retained for $-M\le j\le N_x$ and $y_j$ for
$-M\le j\le N_y$, where $N_x,N_y\ge M$.
For $N\ge M+2|r|$, put $M'=M+|r|$. Then
\[
 j_r:E_{M,N}\hookrightarrow E_{M';N+r,N-r},\qquad
 x'_j=x_{j-r},\quad y'_j=y_{j+r},
\]
where an input coefficient below $-M$ is zero. The upper
cutoffs ensure that every coefficient on the right is retained.
The image is the union of components with indices
$-M+r,\ldots,M+r$, so $j_r$ is a closed immersion.
Every component of either factor has dimension $2N+2$.
The asymmetric factor is compared with a symmetric one by
the affine projection
\[
 E_{M',N+|r|}\longrightarrow E_{M';N+r,N-r},
\]
which forgets $2|r|$ free coefficients. These finite-factor presentations
and their refinements therefore belong to the finite-factor
diagram of Proposition~\ref{prop:global-finite-presentations}.

Since $\mu_\Omega(\ell,q)=0$, the contact-index formula
\eqref{eq:contact-factor-volume-index} gives index $2N+2$
on the dense smooth contact loci of both factors. Local
constancy extends this equality to their closures, by
Proposition~\ref{prop:lcy-equidimensional-factors}(ii).
Thus $j_r$ gives the normalized proper trace
\[
 (j_r)_*\omega_{E_{M,N}}[-2(2N+2)]
 \longrightarrow\omega_{E_{M';N+r,N-r}}[-2(2N+2)].
\]
Together with the nonzero-sector identifications and normalized
affine pullbacks, these traces define the map induced by $T_r$.
The contact bound increases by at most $3|r|$, by
\eqref{eq:hyperbola-contact-dictionary}. The constructions
for $r$ and $-r$ compose to the pole inclusions after common
refinement, so $T_r$ induces a degree-zero linear automorphism
of $LH^*_\Omega(X)$.
It preserves the normalized fundamental classes of the
corresponding contact-component closures, so
\eqref{eq:lcy-theta-stage-class}--\eqref{eq:lcy-theta-colimit-class}
give
\[
 T_r\theta_{k,n}=\theta_{k,n+r}
 \qquad(k,n,r\in\Z).
\]
In particular, $T_r1=\theta_{0,r}$, since $\theta_{0,0}=1$.

\emph{Step 2: construct normalized marked-point maps and source-unit isotopies on finite factors.}
Let $\mathscr L_{\{1,2\}}X\to\A^{\{1,2\}}$ be the
two-point loop space of
Definition~\ref{def:kv-multipoint-functor}, constructed
from the affine line with coordinate $z$ and marked
sections $z=a_1$, $z=a_2$.
For $n,m\in\Z$, define an automorphism
\[
 S_{n,m}:\mathscr L_{\{1,2\}}X
 \xrightarrow{\sim}\mathscr L_{\{1,2\}}X
 \quad\text{over }\A^{\{1,2\}}
\]
as follows. For a $\C$-algebra $R$, an $R$-point consists
of marks $a=(a_1,a_2)\in R^2$ and a homomorphism
$\gamma:\mathcal O(X)\to\mathcal K_a$. Denote the images
of the coordinate functions by $x_\gamma,y_\gamma,f_\gamma$;
they satisfy
\[
 x_\gamma y_\gamma=1+f_\gamma,
 \qquad f_\gamma\in\mathcal K_a^\times.
\]
Since $D_a=(z-a_1)(z-a_2)$ is inverted in $\mathcal K_a$,
both factors are units. Hence
\[
 u_a:=(z-a_1)^n(z-a_2)^m\in\mathcal K_a^\times.
\]
Fix the marks and assign to  $\gamma$  a new homomorphism
$\gamma'$ determined by
\[
 \gamma'(x)=u_ax_\gamma,\qquad
 \gamma'(y)=u_a^{-1}y_\gamma,\qquad
 \gamma'(f)=f_\gamma.
\]
This preserves the defining relation, since
$(u_ax_\gamma)(u_a^{-1}y_\gamma)=1+f_\gamma$.
The construction commutes with homomorphisms of test
algebras $R\to R'$ and has inverse obtained by replacing
$(n,m)$ with $(-n,-m)$, so it defines the asserted
automorphism, including when the marks coincide. 
Writing $x,y,f$ for their images under $\gamma$, we abbreviate this as
\[
 S_{n,m}(x,y,f)=\bigl(u_a x,\,u_a^{-1}y,\,f\bigr).
\]
We use the same convention for similar maps below.
For a pole-bounded family $\mathcal V$ of
Definition~\ref{def:kv-bounded-family}, choose an enlarged
bounded family $\mathcal V'$ containing its image. Such bounds
exist coefficientwise: $u_a^{\pm1}$ are finite Laurent
polynomials in $q_0=D_a(z)$ in the basis $1,z$, and the
pullbacks of the fixed affine generators of $\mathcal O(X)$
therefore have uniform finite lower bounds. Define
\[
 \mathcal V^S
 =\{\gamma'\in\mathcal V':S_{-n,-m}(\gamma')\in\mathcal V\}.
\]
The pole equations of $\mathcal V$, pulled back through
$S_{-n,-m}$, impose finitely many negative-coefficient
vanishings at each fixed bound. They define a closed
subscheme $\mathcal V^S\subset\mathcal V'$ and a factorization
\[
 \mathcal V\xrightarrow[\sim]{S_{n,m}}\mathcal V^S
 \hookrightarrow\mathcal V'.
\]
Include this factorization in the finite-factor presentation diagrams of
Proposition~\ref{prop:direct-collision-diagrams}. On each chosen
chart of $\mathcal V'$, present the pulled-back inclusion of
$\mathcal V^S$ by a closed finite-factor inclusion with the same
affine tail. Transport the resulting charts of $\mathcal V^S$
through $S_{-n,-m}$ to charts of $\mathcal V$; these identified
charts present the first arrow as well.

We verify preservation of the relative normalization before
applying specialization. At a point over a field $K/\C$ with
marks $a=(a_1,a_2)\in K^2$, use the punctured-source algebra
$\mathcal K_a$ of \eqref{eq:kv-monic-source}.
The frame $(dx,dy/f)$ has exterior product $\Omega$.
With $q_0,z,a_1,a_2$ held fixed, the differential of $S_{n,m}$
in this frame is $G=\operatorname{diag}(u_a,u_a^{-1})$.
Multiplication by $z-a_i$ on $\mathcal K_a$, in the basis $1,z$,
has determinant $-q_0$. Hence, with $\operatorname{Nm}$ denoting the norm defined in
Subsection~\ref{subsec:direct-relative-volume-indices},
\[
 \operatorname{Nm}_{\mathcal K_a/K(\!( q_0 )\!)}(u_a)=(-q_0)^{n+m},\qquad
 \det_{K(\!( q_0 )\!)}(G)=
 \operatorname{Nm}_{\mathcal K_a/K(\!( q_0 )\!)}(u_a)
 \operatorname{Nm}_{\mathcal K_a/K(\!( q_0 )\!)}(u_a^{-1})=1.
\]
For two charts identified by $S_{n,m}$, use the same finite
factor and omitted coordinates. If $\kappa,\kappa'$ are their
relative coefficient-residue maps \eqref{eq:direct-relative-residue-map} and
$L=\ker\kappa$, $L'=\ker\kappa'$, coefficientwise
differentiation gives
\[
 \kappa'=\kappa\circ G,\qquad L'=G^{-1}L.
\]
Since $G$ has determinant one, it carries a determinant-normalized
basis to another such basis. The relative pointwise coefficient-index
formula \eqref{eq:direct-relative-field-index} therefore gives
equal values of $c_\pi$ for the corresponding finite-factor
projections of the transported charts. This calculation holds at
every field-valued point, including $a_1=a_2$, and hence
identifies their locally constant relative volume indices.
The isomorphism onto $\mathcal V^S$ and the closed trace into
$\mathcal V'$ consequently preserve the relative normalization;
their common refinements use the same normalized purity equivalences.

Pull back this family and $S_{n,m}$ along
\[
 \A^1_s\longrightarrow\A^{\{1,2\}},
 \qquad s\longmapsto(s,-s),
\]
as in Example~\ref{ex:kv-binary-family}, retaining the notation
$S_{n,m}$ for the pulled-back automorphism.
Fix $\varepsilon>0$ and restrict their analytic spaces
to $J_\varepsilon:=(0,\varepsilon)\subset\C_s$.
Writing $s=b\in J_\varepsilon$, factorization identifies the
restricted loop family with
\[
 J_\varepsilon\times\mathscr LX(\C)\times\mathscr LX(\C),
\]
using the local source coordinates $t_1=z-b$ and $t_2=z+b$.
The function multiplying $x$ in $S_{n,m}$
has the respective expansions
\[
 \begin{aligned}
 (z-b)^n(z+b)^m
   &=t_1^n u_{1,b}(t_1),
 &u_{1,b}(t)&=(2b+t)^m,\\
 (z-b)^n(z+b)^m
   &=t_2^m u_{2,b}(t_2),
 &u_{2,b}(t)&=(-2b+t)^n.
 \end{aligned}
\]
Thus on the first input the Laurent shift $T_n$ is
accompanied by multiplication of $x$ by $u_{1,b}$,
and on the second input $T_m$ is accompanied by
multiplication of $x$ by $u_{2,b}$.
In both cases $y$ is multiplied by the inverse unit
and $f$ is fixed.

For each $i\in\{1,2\}$, identify $t_i$ with $t$ and,
in the following construction, write $u_b=u_{i,b}$. Since
$u_b\in\C\llbracket t\rrbracket^\times$, write
\[
 u_b(t)=c_b(1+r_b(t)),\qquad
 c_b=u_b(0)\in\C^*,\qquad
 r_b\in t\C\llbracket t\rrbracket.
\]
Choose a continuous logarithm of $c_b$ on $J_\varepsilon$ and set
\[
 u_{b,\tau}(t)=\exp(\tau\log c_b)(1+\tau r_b(t)),
 \qquad 0\le\tau\le1.
\]
Then $u_{b,0}=1$, $u_{b,1}=u_b$, and the coefficients
of $u_{b,\tau}^{\pm1}$ are continuous in $(b,\tau)$.
These {\it source units} define a continuous family of automorphisms
\[
 \begin{aligned}
 \Psi^{(i)}_{b,\tau}:\mathscr LX(\C)&\longrightarrow\mathscr LX(\C),\\
 (x,y,f)&\longmapsto
       (u_{b,\tau}x,u_{b,\tau}^{-1}y,f).
 \end{aligned}
\]
For $i=1,2$, this acts on the $i$th input loop in its
coordinate $t_i$. Its endpoints are
\[
 \Psi^{(i)}_{b,0}=\mathrm{id}_{\mathscr LX(\C)},\qquad
 \Psi^{(i)}_{b,1}(x,y,f)
   =(u_{i,b}x,u_{i,b}^{-1}y,f).
\]
Each $\Psi^{(i)}_{b,\tau}$ preserves the orders of $x,y,f$,
hence the pole bounds.
The contact-index argument of Step~1 gives equal volume
indices on corresponding finite coefficient models.

The following formulas, with $c=u_{b,\tau}$ for the chosen
input $i$, describe the automorphisms of the finite coefficient
models $Y_N$ and $E_{M,N}$ induced by $\Psi^{(i)}_{b,\tau}$. For a source unit $c\in\C\llbracket t\rrbracket^\times$,
write
\[
 c_h=[t^h]c,\qquad d_h=[t^h]c^{-1}
 \quad(h\ge0),\qquad c_0\ne0.
\]
On the model $Y_N$ of $\mathscr L_{k,n}X$, multiplication by
$c$ on $x$ induces
\[
 Y_N\longrightarrow Y_N,\qquad
 (A_j,F_j)_{j=0}^N\longmapsto
 \left(\sum_{h=0}^j c_hA_{j-h},\,F_j\right)_{j=0}^N.
\]
On $E_{M,N}$ it induces
\[
 x_j\longmapsto\sum_{h=0}^{j+M}c_hx_{j-h},\qquad
 y_j\longmapsto\sum_{h=0}^{j+M}d_hy_{j-h}
 \qquad(-M\le j\le N).
\]
These maps preserve each component of $E_{M,N}$ and
$[t^0](xy)$; their inverses interchange $c$ and $c^{-1}$.
The coefficients $d_h$ are polynomials in
$c_0^{-1},c_1,\ldots,c_h$, so $c=u_{b,\tau}$ gives continuous
families of automorphisms with continuous inverses.
The formulas commute with coefficient truncation and with
$E_{M,N}\hookrightarrow E_{M',N}$ for $M\le M'\le N$,
where the latter map sets the new coefficients below $-M$
equal to zero.

\emph{Step 3: prove shift equivariance and the relation $WW_-=1$.}
For input pole bounds $P,Q$ and a sufficiently large common
cutoff $N$, use the factors $Z_{P,N},Z_{Q,N}$ of
\eqref{eq:hyperbola-global-factors}. Set
\[
 Y=J_\varepsilon\times(Z_{P,N}\times Z_{Q,N})(\C).
\]
For $b\in J_\varepsilon$ and $\tau\in[0,1]$, let
\[
 \begin{aligned}
 \psi^{(1)}_{b,\tau}&:Z_{P,N}(\C)\longrightarrow Z_{P,N}(\C),\\
 \psi^{(2)}_{b,\tau}&:Z_{Q,N}(\C)\longrightarrow Z_{Q,N}(\C)
 \end{aligned}
\]
be the finite coefficient automorphisms of Step~2 induced
by $\Psi^{(1)}_{b,\tau}$ and $\Psi^{(2)}_{b,\tau}$,
respectively. For each $\tau\in[0,1]$, define
\[
 \begin{aligned}
 \Psi_\tau:Y&\longrightarrow Y,\\
 (b,y_1,y_2)&\longmapsto
 \bigl(b,\psi^{(1)}_{b,\tau}(y_1),
          \psi^{(2)}_{b,\tau}(y_2)\bigr).
 \end{aligned}
\]
Let $p:Y\to J_\varepsilon$ be the projection and put $I=[0,1]$.
Then $p\circ\Psi_\tau=p$ and $\Psi_0=\mathrm{id}_Y$.
Also write $\Psi_\tau$ for the product on
$J_\varepsilon\times\mathscr LX(\C)^2$ given by
\[
 (b,\gamma_1,\gamma_2)\longmapsto
 \bigl(b,\Psi^{(1)}_{b,\tau}(\gamma_1),
          \Psi^{(2)}_{b,\tau}(\gamma_2)\bigr);
\]
the preceding coefficient maps are its finite representatives.
The map
\[
 \widetilde\Psi:I\times Y\xrightarrow{\sim}I\times Y,
 \qquad(\tau,y)\longmapsto(\tau,\Psi_\tau(y))
\]
is a homeomorphism over $I\times J_\varepsilon$.
Set
\[
 \begin{gathered}
 \delta=\delta_{Z_{P,N}}+\delta_{Z_{Q,N}}=4N+4,\\
 \mathcal B_{Y/J_\varepsilon}=Rp_*\omega_Y[-2\delta],\\
 \mathcal C_Y=R\Gamma(J_\varepsilon,\mathcal B_{Y/J_\varepsilon}).
 \end{gathered}
\]
These are the relative normalized complex and its derived
global sections from
Definition~\ref{def:relative-normalized-dual-complexes};
$\mathcal C_Y[-1]$ is the finite-factor source of
specialization in \eqref{eq:geometric-specialization-sections}.

To construct the homotopy on $\mathcal B_{Y/J_\varepsilon}$, let
$a:I\times Y\to I$, $\pi_Y:I\times Y\to Y$, and
$\pi_J:I\times J_\varepsilon\to J_\varepsilon$ be the
projections, and put $\widetilde p=\mathrm{id}_I\times p$.
The relative dualizing object for $a$ satisfies
\[
 \begin{gathered}
 a^!\underline\F_I\simeq\pi_Y^{-1}\omega_Y,\\
 R\widetilde p_*\bigl(a^!\underline\F_I\bigr)[-2\delta]
       \simeq\pi_J^{-1}\mathcal B_{Y/J_\varepsilon}.
 \end{gathered}
\]
For the second comparison, compute on products with contractible
interval neighborhoods in $I$ and apply the homotopy-invariance
argument of Lemma~\ref{lem:infinite-product-descent}(i).
Dualizing functoriality for the homeomorphism
$\widetilde\Psi$ induces an endomorphism of
$\pi_J^{-1}\mathcal B_{Y/J_\varepsilon}$ whose ordinary inverse-image
restriction at $\tau$ is $(\Psi_\tau)_*$.
By \cite[Theorem~0.6]{HainePortaTeyssier},
\[
 \pi_J^{-1}:\mathcal D(J_\varepsilon;\F)
       \longrightarrow\mathcal D(I\times J_\varepsilon;\F)
\]
is fully faithful. The restriction at $\tau=0$ is the
identity, so the restriction at $\tau=1$ is homotopic to it.
The coefficient truncations and pole inclusions of Step~2
are equivariant for the whole $I$-family. Relative dualizing
functoriality, the normalized purity equivalences, and closed proper
trace therefore give the same endomorphism of the corresponding diagram of
normalized complexes.
Applying the fully faithful functor above to that diagram
gives compatible homotopies, and
Proposition~\ref{prop:global-finite-presentations}(i) transports
them to the pole-stage diagram.

On global cochains this can also be represented by a proper
prism. Indeed,
$H:I\times Y\to Y$, $H(\tau,y)=\Psi_\tau(y)$, is proper,
since it is the composite of $\widetilde\Psi$ and the
projection with compact fiber $I$. In the model
$\mathcal C_Y^j=C^{\mathrm{lf}}_{2\delta-j}(Y;\F)$, let
$K(c)=H_*([I]\times c)$, using the oriented singular-chain
prism, and let $d$ be the singular boundary. Then $K$ has
degree $-1$ and
\[
 dK+Kd=(\Psi_1)_*-\mathrm{id}_{\mathcal C_Y}.
\]
The relative construction above supplies the compatibility
with refinements and restriction that is needed for
specialization; no restriction of raw singular chains to
arbitrary open subsets is used.

Under the input and output identifications,
\[
 S_{n,m}|_{J_\varepsilon}=(T_n\times T_m)\circ\Psi_1,\qquad
 S_{n,m}|_{s=0}=T_{n+m}.
\]
Use the finite-factor presentation of the closed-image
factorization and the
relative-normalization comparison of Step~2. On the input
normalized complexes, the preceding compatible homotopy gives
\[
 (S_{n,m}|_{J_\varepsilon})_*
   =(T_n\times T_m)_*\circ(\Psi_1)_*
   \simeq(T_n\times T_m)_*.
\]
Apply Lemma~\ref{lem:source-specialization-comparison}(i)
to the isomorphism onto $\mathcal V^S$, and
Lemma~\ref{lem:specialization-arrow-descent}(ii) to its
closed proper trace. Their naturality identifies specialization
after $(S_{n,m}|_{J_\varepsilon})_*$ with $(T_{n+m})_*$
after specialization.
With the input and output comparisons of
Definition~\ref{lp:def:binary-fusion-product}, these
specialization maps give the pole-stage fusion maps.
Corollary~\ref{prop:loop-fusion-cohomology} identifies
their induced product on the pole colimit with $\mu$.
The preceding homotopy and specialization compatibility
therefore give
\begin{equation}\label{eq:hyperbola-shift-equivariance}
 \mu(T_na,T_mb)=T_{n+m}\mu(a,b),
 \qquad a,b\in LH^*_\Omega(X).
\end{equation}
Only $S_{n,m}$ is used at $s=0$; the homotopy is over $J_\varepsilon$.
Since $T_n1=\theta_{0,n}$, the unit identity gives
\[
 \theta_{0,n}a=T_na,\qquad WW_-=T_0(1)=1.
\]
Thus multiplication by $W$ is $T_1$. On
$\mathscr L_{k,n}X$ it leaves the series $A,F$ unchanged,
so it preserves their angular classes $\alpha,\beta$.

\emph{Step 4: prove the formula for $\theta_{k,n}$ by same-sign fusion.}
Let $k_1,k_2\in\Z$ with $k_1k_2>0$. By
\eqref{eq:hyperbola-shift-equivariance}, it suffices to take
inputs in $\mathscr L_{k_1,0}X$ and $\mathscr L_{k_2,0}X$.
For a $\C$-algebra $R$ and marks $a=(a_1,a_2)\in R^2$,
use the source algebras $\mathcal O_a,\mathcal K_a$ of
\eqref{eq:kv-source-algebras}, and put
\[
 H_a=(z-a_1)^{k_1}(z-a_2)^{k_2}\in\mathcal K_a^\times.
\]
The two units $x,F\in\mathcal O_a^\times$ define a two-point
loop by
\[
 (x,F)\longmapsto
 \bigl(x,(1+H_aF)x^{-1},H_aF\bigr)\in X(\mathcal K_a).
\]
When $a_1-a_2\in R^\times$, write $t_i=z-a_i$ and let
$F_i$ be the expansion of $F$ on the $i$-th disc. Then
\[
 f=t_i^{k_i}(a_i-a_j+t_i)^{k_j}F_i,
 \qquad \{i,j\}=\{1,2\}.
\]
The last two factors are units. Factorization therefore identifies
this family over the locus of distinct marks with
$\mathscr L_{k_1,0}X\times\mathscr L_{k_2,0}X$ over that locus.

Pull this family back along $s\mapsto(s,-s)$ as in
Example~\ref{ex:kv-binary-family}, and put $q_0=z^2-s^2$.
For a $\C[s]$-algebra $R$, its source rings are
\[
 \mathcal O_{s,R}
   =R\llbracket q_0\rrbracket[z]/(z^2-s^2-q_0),
 \qquad
 \mathcal K_{s,R}=\mathcal O_{s,R}[q_0^{-1}].
\]
Write $H=(z-s)^{k_1}(z+s)^{k_2}$. We first show that the
family parametrized by $x,F\in\mathcal O_{s,R}^\times$
is the entire reduced collision closure of the chosen inputs,
in a sufficiently large pole-bounded family.
For a loop $(x,y,f)$ in that ambient family, set
\[
 F=H^{-1}f,\qquad G=Hf^{-1}.
\]
Impose the following conditions:
\[
 \begin{cases}
 x,y,F,G\in\mathcal O_{s,R},&k_1,k_2>0,\\
 x,H^{-1}y,F,G\in\mathcal O_{s,R},&k_1,k_2<0.
 \end{cases}
\]
These are closed coefficient conditions, over arbitrary $R$.
Indeed, $H^{\pm1}$ are finite Laurent polynomials in $q_0$
with coefficients in $\C[s,z]$, because
$(z-s)^{-1}=(z+s)/q_0$ and
$(z+s)^{-1}=(z-s)/q_0$.
At a fixed pole bound, each displayed expression has a fixed
lower $q_0$-degree and a unique expansion
$\sum_j(r_{0,j}+zr_{1,j})q_0^j$.
Membership in $\mathcal O_{s,R}$ is exactly the vanishing
of the finitely many coefficients $r_{0,j},r_{1,j}$ with $j<0$.

The free basis $1,z$ over $R\llbracket q_0\rrbracket$
shows that $q_0$ is a non-zero-divisor in $\mathcal O_{s,R}$,
so $\mathcal O_{s,R}\hookrightarrow\mathcal K_{s,R}$.
On this closed family, $FG=1$ therefore holds in
$\mathcal O_{s,R}$, so $F$ is a source unit.
If $k_1,k_2>0$, then $H\in q_0\mathcal O_{s,R}$ and
$xy=1+HF$ is a unit, forcing $x$ and $y$ to be units.
If $k_1,k_2<0$, then $H^{-1}\in q_0\mathcal O_{s,R}$ and
\[
 x(H^{-1}y)=F+H^{-1}=F(1+H^{-1}G)
\]
is a unit, again forcing $x$ to be a unit.
Here $1+q_0\mathcal O_{s,R}$ consists of units by completeness.
In either case $y=(1+HF)x^{-1}$. Conversely, every pair of
units $x,F$ satisfies the stated closed conditions.
Thus these conditions define precisely the two-unit family,
as a functor of $R$, and give its closed inclusion in the
bounded collision family.

For $h\ge1$, let $\mathcal D_h\to\A^1_s$ parametrize
the images of $x,F$ in $(R[z]/(q_0^h))^\times$.
Writing $x\equiv u_x+zv_x$ and
$F\equiv u_F+zv_F\pmod{q_0}$ gives
\[
 \mathcal D_1=
 \left\{(s,u_x,v_x,u_F,v_F)\in\A^5:
 (u_x^2-s^2v_x^2)(u_F^2-s^2v_F^2)\ne0\right\},
 \qquad
 \mathcal D_h\cong\mathcal D_1\times\A^{4(h-1)}.
\]
The last isomorphism records the two coefficients of $q_0^j$
in each unit for $1\le j<h$.
Each $\mathcal D_h$ is reduced and smooth, hence flat,
over $\A^1_s$, and its coordinate ring injects into its
localization at $s$. Adjoining all remaining coefficients
preserves reducedness and this injectivity.
Consequently the closed two-unit family is the reduced
schematic closure of its locus $s\ne0$, which proves the
claimed identification with the full collision closure.

Include this closed inclusion and its input and output maps in
the finite-factor presentation diagram. On each ambient nice chart,
enlarge its finite factor to present the pulled-back closed
inclusion with the same affine tail. For each finite list of
chosen input and output factors, enlarge the source factors so
that the maps descend. Compare the resulting finite-factor presentations on
the closed family with its explicit models $\mathcal D_h$ by
common finite-factor refinements, taking $h$ sufficiently large.
These choices require only finitely many coefficients: the
closed conditions above are finite, and each specified
coefficient of the parametrization or its inverse
$(x,y,f)\mapsto(x,H^{-1}f)$ depends on finitely many coefficients.
Further truncations of the units give the
affine-bundle comparisons
$\mathcal D_{h+1}\to\mathcal D_h$ of relative dimension $4$.
Lemma~\ref{lem:specialization-arrow-descent}(ii) therefore
computes the restriction of fusion to these inputs by
specialization on this closed family, followed by its closed
proper trace to the output pole stage. In particular, this
computes the actual fusion product, including all output
contacts.

Put $E=(\C^*)^2$. For $s\ne0$, the constant coefficients
of $x,F$ on the two input discs are
\[
 e_\pm=(u_x\pm sv_x,u_F\pm sv_F)\in E.
\]
At $s=0$, their output constant coefficients are
$(u_x,u_F)\in E$, and
\[
 \mathcal D_h|_{s\ne0}
   \cong\Gm\times E\times E\times\A^{4h-4},
 \qquad
 (\mathcal D_h)_0\cong E\times\A^{4h-2}.
\]
On a separated fiber, the two factors are the models
$Y_{h-1}$ of Lemma~\ref{lem:hyperbola-coefficient-models}(i),
so their volume indices sum to $4h$.
Proposition~\ref{prop:direct-relative-volume-indices} and
connectedness of $\mathcal D_h$ give
$\delta_{\mathcal D_h/\A^1_s}=4h$ everywhere.
Thus every fiber is normalized by $\omega[-8h]$;
the output factor is $Y_{2h-1}$, with the same index $4h$.

Apply Borel--Moore specialization along $J_\varepsilon$
toward $s=0$, followed by the normalized affine-bundle
comparisons with $E\times E$ on the inputs and $E$ on
the output. The calculation in Step~4 of the proof of
Lemma~\ref{lem:fusion-leading-product} applies to these
finite models: $e_--e_+=-2s(v_x,v_F)$ identifies
$\mathcal D_1$ with the deformation to the normal bundle
of the diagonal of $E$, and the factor $-2$ preserves
complex orientation. Consequently the resulting operation
in ordinary cohomology is
\[
 \begin{aligned}
 \operatorname{diag}_E^*&:
 H^*(E\times E;\F)\longrightarrow H^*(E;\F),\\
 \operatorname{diag}_E&:E\longrightarrow E\times E,
 \qquad e\longmapsto(e,e).
 \end{aligned}
\]
The higher truncations add the same free coefficients to
this specialization diagram, so normalized affine pullback
makes these operations compatible with refinement.

Let $\alpha\in H^1(E;\F)$ be the angular class of its
first coordinate, corresponding to the constant coefficient
of $x$. The other-mark factors in the displayed expression
for $f$ do not change $x$ or its angular class. Hence
\[
 \operatorname{diag}_E^*(1\boxtimes1)=1,\qquad
 \operatorname{diag}_E^*(\alpha\boxtimes1)
 =\operatorname{diag}_E^*(1\boxtimes\alpha)=\alpha.
\]
The central fiber maps isomorphically onto
$\mathscr L_{k_1+k_2,0}X$ by
$(x,F)\mapsto(x,(1+z^{k_1+k_2}F)x^{-1},z^{k_1+k_2}F)$.
Its normalized fundamental class is therefore
$\theta_{k_1+k_2,0}$ after the closed proper trace.
The preceding computation gives
\[
 \mu(\theta_{k_1,0},\theta_{k_2,0})
       =\theta_{k_1+k_2,0}
\]
with coefficient $+1$ and no other output term, and also
computes the products with the class $\alpha$ on either input.
Iteration, followed by
\eqref{eq:hyperbola-shift-equivariance}, proves the formula
for $\theta_{k,n}$ in \eqref{eq:hyperbola-exchange}.

\emph{Step 5: identify the reduced collision closure for $U$ and $V$.}
Use the binary slice $a_1=s$, $a_2=-s$ of
Example~\ref{ex:kv-binary-family}. Let $\mathcal W$ be the
reduced closure of the separated family
$\Gm\times\mathscr L_{1,0}X\times\mathscr L_{-1,0}X$,
with parameter $s\in\Gm$, in a sufficiently large pole-bounded
collision family, with the relative volume index restricted
from the ambient family.
Write its coordinate series as
\[
 x=a(z),\qquad y=\frac{b(z)}{z+s},\qquad
 f=\frac{a(z)b(z)-z-s}{z+s}.
\]
The condition $f(s)=0$ is
\begin{equation}\label{eq:hyperbola-node-family}
 a(s)b(s)=2s.
\end{equation}
Expand
\[
 a(z)=\sum_{j\ge0}(a_{2j}+a_{2j+1}z)(z^2-s^2)^j,
 \qquad
 b(z)=\sum_{j\ge0}(b_{2j}+b_{2j+1}z)(z^2-s^2)^j.
\]
Set $A=a(s)=a_0+sa_1$ and $B=b(s)=b_0+sb_1$, so that
$AB=2s$, and put
\[
 H=1-Ab_1-Ba_1+2sa_1b_1,\qquad J=(ab)'(s)-1,
\]
where the prime denotes differentiation with respect to $z$.
For $s\ne0$, one has $a(-s)b(-s)=2sH$; hence $H\ne0$
and $J\ne0$ impose respectively
$\ord_{z+s}f=-1$ and $\ord_{z-s}f=1$.
The equalities $AB=2s$ and $a(-s)b(-s)=2sH$ also give
$\ord_{z-s}x=\ord_{z+s}x=0$.
At $s=0$, one has $J=-H$, giving
\begin{equation}\label{eq:hyperbola-central-two-components}
 AB=0,\qquad 1-Ab_1-Ba_1\ne0.
\end{equation}

We claim that $\mathcal W$ has coefficient ring
\[
 \mathcal A=
 \bigl(\C[s,a_j,b_j:j\ge0]/(AB-2s)\bigr)[H^{-1},J^{-1}].
\]
Each equation uses only finitely many coefficients. In the
formal source ring
\[
 \mathcal O_s=
 \mathcal A\llbracket q_0\rrbracket[z]/(z^2-s^2-q_0),
\]
equation~\eqref{eq:hyperbola-node-family} gives
\[
 ab-z-s=(z-s)C,\qquad C\in\mathcal O_s.
\]
To obtain $C$, subtract the value at $s$ and divide each
coefficient polynomial by $z-s$, using $q_0=(z-s)(z+s)$.
This factorization is unique because $q_0$ is a
non-zero-divisor in the free rank-two
$\mathcal A\llbracket q_0\rrbracket$-module $\mathcal O_s$,
and hence multiplication by $z-s$ is injective.
The evaluations and the determinant of multiplication by
$C$ modulo $q_0$ are
\[
 C(s)=J,\qquad C(-s)=-H,\qquad
 \det\bigl(C\cdot(-):\mathcal O_s/(q_0)\to\mathcal O_s/(q_0)\bigr)
 =-JH.
\]
The determinant identity holds for $s\ne0$ and extends as
a polynomial identity to $s=0$.
Thus $C$ is a unit modulo $q_0$, and completeness gives
$C\in\mathcal O_s^\times$. Consequently
\[
 f=\frac{z-s}{z+s}C,\qquad
 f^{-1}=\frac{z+s}{z-s}C^{-1}
 \quad\text{in }\mathcal O_s[q_0^{-1}].
\]
These formulas define a family of maps to $X$ with the
required separated input orders.

Conversely, on that separated input,
\[
 a=x,\qquad b=(z+s)y,\qquad
 C=\frac{(z+s)^2}{q_0}f,\qquad
 C^{-1}=\frac{(z-s)^2}{q_0}f^{-1}
\]
are regular on both formal discs. These expressions use
$q_0^{-1}$ but not $s^{-1}$. At fixed pole and coefficient
bounds, their regularity is the vanishing of their negative
$q_0$-coefficients, a closed condition. It therefore holds
on $\mathcal W$, where $C(s)=J$ and $C(-s)=-H$ force
$J,H$ to be units. This gives $\mathcal W\subset\Spec\mathcal A$.
For the reverse inclusion, use the polynomial coordinates
$s,A,B,a_1,b_1,a_2,b_2,\ldots$. At every finite cutoff
retaining $H,J$, the equation $AB=2s$ defines a smooth
total space, flat over $\A^1_s$, with dense locus $s\ne0$.
Adjoining the remaining free coefficients and inverting
$H,J$ preserve these properties. Hence
$\Spec\mathcal A\subset\mathcal W$, proving the claim.

Include the closed inclusion of $\mathcal W$ and its input
and output maps in the finite-factor presentation diagram over
$\A^1_s$, using sufficiently large cutoffs for each finite
collection of chosen factors as in Step~4.
Lemma~\ref{lem:specialization-arrow-descent}(ii) computes
$\mu(U,V)$ by normalized specialization on $\mathcal W$
followed by the closed proper trace to the output pole stage.
Normalized affine pullback removes the free coefficients.

\emph{Step 6: compute the specialization multiplicities and prove $UV=1+W$.}
The two components of $\mathcal W_0$ in
\eqref{eq:hyperbola-central-two-components} are
$\mathcal C_A=\{A=0\}$ and $\mathcal C_B=\{B=0\}$.
Their generic output orders are
\[
 \begin{array}{c|ccc|c}
  &\ord_z x&\ord_z y&\ord_z f&\text{theta class}\\ \hline
  \mathcal C_B&0&0&0&\theta_{0,0}\\
  \mathcal C_A&1&-1&0&\theta_{0,1}.
 \end{array}
\]
Indeed, at $s=0$ the output map is
\[
 (a(z),b(z))\longmapsto(x(z),y(z))=(a(z),b(z)/z).
\]
On $\mathcal C_B$ one has $b(0)=0$, so this identifies
$\mathcal C_B$ with the component indexed by $n=0$ in
Lemma~\ref{lem:hyperbola-coefficient-models}(ii).
On $\mathcal C_A$ one has $a(0)=0$, so it identifies
$\mathcal C_A$ with the component indexed by $n=1$ there.
In both cases the inverse is $(a,b)=(x,zy)$, and the open
condition in \eqref{eq:hyperbola-central-two-components}
is precisely $f(0)\ne0$. Thus, after matching coefficient
cutoffs, the output maps have generic degree one on both
components. Since $s=AB/2$, its order along each component
is also one. Writing $\operatorname{sp}_{\mathcal W}$ for
normalized specialization, Lemma~\ref{lp:lem:bm-specialization-formula}
therefore gives
\[
 \operatorname{sp}_{\mathcal W}[\mathcal W_s]
   =[\mathcal C_A]+[\mathcal C_B]\qquad(s\ne0).
\]
The output proper trace sends these fundamental classes to
$\theta_{0,1}$ and $\theta_{0,0}$, respectively, with coefficient
$+1$. Hence $UV=\theta_{0,1}+\theta_{0,0}=1+W$, proving the
remaining relation in \eqref{eq:hyperbola-exchange}.

\emph{Step 7: prove the asserted ring isomorphism.}
Every irreducible component of $Z_{P,N}$ in
\eqref{eq:hyperbola-global-factors} has dimension $2N+2$,
equal to its volume index. Hence
\[
 H^0\!\left(R\Gamma\!\left(Z_{P,N},
       \omega_{Z_{P,N}}[-2(2N+2)]\right)\right)
 =H^{\mathrm{BM}}_{4N+4}(Z_{P,N};\F)
\]
has the fundamental classes of these components as a basis.
Coefficient refinement carries each such class to the
fundamental class of its affine-bundle inverse image, and
the pole inclusions include the old components unchanged. Therefore
Proposition~\ref{prop:global-finite-presentations} identifies
$R_X$ with the vector space having basis
$\{\theta_{k,n}:k,n\in\Z\}$, with the theta classes defined in
\eqref{eq:lcy-theta-stage-class}--\eqref{eq:lcy-theta-colimit-class}.

The proved relations define
\[
 \phi:\F[u,v,w,w^{-1}]/(uv-1-w)\longrightarrow R_X,
 \qquad (u,v,w)\longmapsto(U,V,W).
\]
Using $uv=1+w$, every element of the source is a linear
combination of
\[
 u^kw^n\quad(k\ge0,\ n\in\Z),\qquad
 v^{-k}w^n\quad(k<0,\ n\in\Z).
\]
Their images are the basis elements $\theta_{k,n}$.
Thus they span the source and have linearly independent
images spanning $R_X$, so $\phi$ is an isomorphism.
\end{proof}

The additive identification in Proposition~\ref{prop:hyperbola-degeneration}
is not multiplicative. Indeed, $W$ and $W_-$ have weight three,
whereas $WW_-=1$ has weight zero. Also $UV=1+W$ contains the
lower-weight term $1$. More intrinsically, let $p_+,p_-$ be the
two points of $L\cap Q$, and let $e_+,e_-$ be their component
units in $H^0(S_{\{L,Q\}}^{\log})$. Then
\[
 (e_+t_Lt_Q)(e_-t_Lt_Q)=0,
 \qquad e_+t_Lt_Q\ne0\ne e_-t_Lt_Q.
\]
Thus $H^0_{\log}(M,D)$ has zero divisors, whereas
$LH^0_\Omega(X)=\F[u,v,(uv-1)^{-1}]$ is a domain. The two
rings are not isomorphic; the collision computation resolves the
multiplicative extension left by the degenerate spectral sequence.

\subsection{Source rotation and the positive-degree classes}

\begin{lemma}[Rotation on the $(k,n)$-sectors]
\label{lem:rotation-bracket}
For $k\ne0$ and $n\in\Z$, source rotation on
$\mathscr L_{k,n}X$ gives
\begin{equation}\label{eq:hyperbola-rotation}
 \Delta(\theta_{k,n}\gamma)
   =\theta_{k,n}(n\iota_\alpha+k\iota_\beta)\gamma,
 \qquad \gamma\in\Lambda(\alpha,\beta).
\end{equation}
Here $\theta_{k,n}\gamma$ denotes the class in $LH^*_\Omega(X)$
represented by $\gamma$ in the normalized cohomology of
$\mathscr L_{k,n}X$; in particular, $\theta_{k,n}1=\theta_{k,n}$.
The classes $\alpha,\beta$ are those of
Lemma~\ref{lem:hyperbola-coefficient-models}(i), and
$\iota_\alpha,\iota_\beta$ are the degree-$-1$ derivations
satisfying $\iota_\alpha(\alpha)=\iota_\beta(\beta)=1$ and
$\iota_\alpha(\beta)=\iota_\beta(\alpha)=0$.
\end{lemma}
\begin{proof}
On the finite coefficient models $Y_N$ of $\mathscr L_{k,n}X$, write
$A(t)=\sum_{j\ge0}A_jt^j$ and
$F(t)=\sum_{j\ge0}F_jt^j$. With the positive convention
$t\mapsto\lambda t$, the whole source-rotation family is
\[
 a_N:\Gm\times Y_N\longrightarrow Y_N,\qquad
 (\lambda,(A_j,F_j)_{j=0}^N)
 \longmapsto(\lambda^{n+j}A_j,\lambda^{k+j}F_j)_{j=0}^N.
\]
These actions commute with coefficient truncation and pole
enlargement. Let $\vartheta\in H^1(\C^*;\F)$ be the angular
class of $\lambda$. Under the cohomology identification of
Lemma~\ref{lem:hyperbola-coefficient-models}(i), one has
\[
 a_N^*\alpha=1\otimes\alpha+n\vartheta\otimes1,\qquad
 a_N^*\beta=1\otimes\beta+k\vartheta\otimes1.
\]
Proposition~\ref{prop:global-finite-bv}(ii) computes $\Delta$
as the coefficient of $\vartheta$. Multiplicativity of
$a_N^*$ therefore gives $n\iota_\alpha+k\iota_\beta$, proving
\eqref{eq:hyperbola-rotation}.
\end{proof}

Put $H=LH^*_\Omega(X)$, so $H^0=R_X$.
The Fano SNC reduced anticanonical compactification of
Subsection~\ref{subsec:example-35-1} has two zero-dimensional
strata. Since $R_X$ is a localization of $\F[u,v]$,
Theorem~\ref{thm:intro-global-polyvectors} and
Remark~\ref{globpv:remark-rational} give the canonical graded
algebra isomorphism
\begin{equation}\label{eq:hyperbola-canonical-polyvectors}
 \sigma:H\xrightarrow{\ \sim\ }
 \bigwedge_{R_X}^*(R_X\partial_u\oplus R_X\partial_v).
\end{equation}
To identify $\sigma^{-1}(\partial_u)$ and
$\sigma^{-1}(\partial_v)$, put
\[
 P=\theta_{-1,0}\alpha,\qquad Q_1=\theta_{1,0}\alpha
 \quad\text{in }H^1,
\]
using the sector-class convention of Lemma~\ref{lem:rotation-bracket}.

\begin{proposition}
\label{prop:hyperbola-exterior}
Under \eqref{eq:hyperbola-canonical-polyvectors}, the inverse
images of $\partial_u$ and $\partial_v$ are, respectively,
\[
 \xi_u=W^{-1}P,\qquad \xi_v=W^{-1}Q_1.
\]
\end{proposition}
\begin{proof}
By \eqref{globpv:eq-sigma}, the degree-one part of $\sigma$ is
the anchor isomorphism
\[
 \rho:H^1\xrightarrow{\sim}\operatorname{Der}_{\F}(R_X),
 \qquad \rho(a)(h)=\{a,h\}.
\]
Formula~\eqref{eq:hyperbola-rotation} gives
$\Delta(P)=\Delta(Q_1)=0$. Applying
\eqref{eq:general-rotation-bracket} to the same-sign products
and Laurent shifts of Lemma~\ref{lem:hyperbola-exchange} gives
\[
 \begin{aligned}
 \{P,V\}&=\Delta(\theta_{-2,0}\alpha)=0,&
 \{P,W\}&=\Delta(\theta_{-1,1}\alpha)=VW,\\
 \{Q_1,U\}&=\Delta(\theta_{2,0}\alpha)=0,&
 \{Q_1,W\}&=\Delta(\theta_{1,1}\alpha)=UW.
 \end{aligned}
\]
The derivation identity and $UV=1+W$ then give
\[
 V\{P,U\}=VW,\qquad U\{Q_1,V\}=UW.
\]
Cancelling $V$ and $U$ in the domain $R_X$ proves
\begin{equation}\label{eq:hyperbola-derivations}
 \rho(P)=W\partial_u,\qquad
 \rho(Q_1)=W\partial_v.
\end{equation}
Since $W$ is invertible, these identities give
$\sigma(\xi_u)=\partial_u$ and $\sigma(\xi_v)=\partial_v$.
Together with \eqref{eq:hyperbola-canonical-polyvectors}, this
proves the loop-side ring presentation in
\eqref{eq:hyperbola-full-ring} with the stated generators.
\end{proof}

\begin{proposition}[The hyperbola BV operator]
\label{prop:hyperbola-bv-formula}
For the generators of Proposition~\ref{prop:hyperbola-exterior}
and every $h\in R_X$, source rotation satisfies
\begin{align*}
 \Delta(h)&=0,\\
 \Delta(h\xi_u)&=\partial_u h-\frac vw h,\qquad
 \Delta(h\xi_v)=\partial_v h-\frac uw h,\\
 \Delta(h\xi_u\xi_v)
   &=\left(\partial_u h-\frac vw h\right)\xi_v
     -\left(\partial_v h-\frac uw h\right)\xi_u.
\end{align*}
\end{proposition}
\begin{proof}
Using $\Delta(P)=\Delta(Q_1)=0$ and
\eqref{eq:hyperbola-derivations}, the bracket convention
\eqref{eq:general-rotation-bracket} gives
\[
 \Delta(W^{-1}P)=\{P,W^{-1}\}=-V/W,
 \qquad
 \Delta(W^{-1}Q_1)=\{Q_1,W^{-1}\}=-U/W.
\]
The Jacobi identity gives
\[
 \rho(\{\xi_u,\xi_v\})
   =[\rho(\xi_u),\rho(\xi_v)]
   =[\partial_u,\partial_v]=0.
\]
Injectivity of $\rho$ implies $\{\xi_u,\xi_v\}=0$.
Together with $\{\xi_u,h\}=\partial_u h$ and
$\{\xi_v,h\}=\partial_v h$, the bracket convention
\eqref{eq:general-rotation-bracket} now gives the displayed
formulas for every $h\in R_X$.
By Example~\ref{ex:volume-divergence-bv}, these are the
formulas for $\Delta_\eta$, with $\eta$ as in
Theorem~\ref{thm:hyperbola-ring}, proving its BV assertion.
\end{proof}

\begin{remark}[Generators from the arc volume class]
\label{rem:hyperbola-volume-generators}
The volume-class argument in
Theorem~\ref{thm:intro-global-polyvectors} gives a nonzero image
of $H^2(X;\F)$ in $\ker(\Delta:H^2\to H^1)$.
Proposition~\ref{prop:hyperbola-bv-formula} gives, for $h\in R_X$,
\[
 \Delta(h\xi_u\xi_v)=0
 \quad\Longleftrightarrow\quad d(h/w)=0
 \quad\Longleftrightarrow\quad h/w\in\F.
\]
Thus the image of an arc generator can be normalized to
$\tau=W\xi_u\xi_v$, and
\[
 \sigma_2(\tau)=w\partial_u\wedge\partial_v,\qquad
 \xi_u=W^{-1}\{\tau,V\},\qquad
 \xi_v=-W^{-1}\{\tau,U\}.
\]
Indeed, \eqref{globpv:eq-sigma} gives
$\rho(\{\tau,V\})=w\partial_u$ and
$\rho(\{\tau,U\})=-w\partial_v$; injectivity of $\rho$
proves the two bracket identities. Hence $H$ is generated from
$R_X$ and $\tau$ by products and brackets.
\end{remark}

\clearpage

\appendix
\addtocontents{toc}{\protect\addvspace{1.5em}}
\addtocontents{toc}{\protect\setcounter{tocdepth}{1}}
\section{Derived complexes, homotopy limits, and homotopy colimits}
\label{app:derived-colimits}

We fix conventions for complexes, duality, and parity grading,
then describe coherent diagrams and their cochain representatives.
We next define homotopy limits, colimits, and cofibers, including
concrete topological models. The remaining subsections treat
\v Cech totalization and filtered objects with their telescope models.
Throughout, the coefficient field is $\bQ$; tensor products
and linear duals are taken over this field.

\subsection{Complexes, duality, and gradings}
\label{app:complexes-duality}

The \emph{derived $\infty$-category} $\mathcal D(\bQ)$ is obtained
from cochain complexes by making quasi-isomorphisms equivalences.
Its morphisms form mapping spaces: cochain maps give points,
chain homotopies give paths between them, and higher homotopies
give homotopies between paths, and so on.

For two complexes $K,L$, the homotopy information in their
mapping space $\operatorname{Map}_{\mathcal D(\bQ)}(K,L)$
is encoded by the mapping complex
\[
 \underline{\Hom}^{r}(K,L)
   =\prod_{n\in\Z}\Hom_{\bQ}(K^n,L^{n+r}),
 \qquad d(f)=d_Lf-(-1)^r f d_K.
\]
Degree-zero cycles are cochain maps; an element $H$ of degree
$-1$ is a homotopy from $f$ to $g$ when $d(H)=g-f$.
Further negative degrees encode higher homotopies. For a cycle $f\in\underline{\Hom}^{r}(K,L)$, a
\emph{nullhomotopy} is an element
$H\in\underline{\Hom}^{r-1}(K,L)$ such that
\[
 d(H)=d_LH-(-1)^{r-1}Hd_K=f.
\]
Thus $f$ is homotopic to zero. Explicitly, if $d_Kx=0$,
then $f(x)=d_L(H(x))$, so the induced map
$H^n(K)\to H^{n+r}(L)$ is zero for every $n$.

A \emph{differential graded category}, or \emph{dg category},
has a complex of morphisms for each pair of objects, with
composition compatible with the differentials. Here the objects
are cochain complexes, the morphism complexes are the displayed
$\underline{\Hom}(K,L)$, and composition is composition of
graded linear maps. Its \emph{dg nerve} organizes these data
into a simplicial set, a system of vertices, edges, triangles,
and higher-dimensional simplices with compatible faces and
degeneracies. Vertices are complexes and edges are cochain maps.
A triangle consists of three such maps and a specified chain
homotopy between the composite of two edges and the third edge;
higher simplices record the higher compatibility homotopies.
Degenerate simplices encode identities. This gives a concrete
model of $\mathcal D(\bQ)$; see
\cite[Section~1.3.1, Definition~1.3.5.8, and
Proposition~1.3.5.15]{HA}.

Every acyclic
complex of vector spaces is contractible, by choosing splittings
of its cycles and boundaries. Thus quasi-isomorphisms are already
homotopy equivalences, and these mapping complexes compute derived
maps without further resolutions. The associated mapping space
satisfies
\[
 \pi_r\operatorname{Map}_{\mathcal D(\bQ)}(K,L)
       =H^{-r}\underline{\Hom}(K,L),\qquad r\ge0,
\]
with zero as basepoint for $r>0$. The ordinary derived category
$D(\bQ)$ is the homotopy category of $\mathcal D(\bQ)$: it has
the same objects, but its morphisms are only the sets
$\pi_0\operatorname{Map}_{\mathcal D(\bQ)}(K,L)
=H^0\underline{\Hom}(K,L)$.
Thus $\mathcal D(\bQ)$ also remembers homotopies between maps
and their higher compatibility data.

\emph{Parity grading and periodization.}
A two-periodic cochain complex is a diagram
\[
 K^{\bar0}\xrightarrow{d^{\bar0}}K^{\bar1}
 \xrightarrow{d^{\bar1}}K^{\bar0},
 \qquad d^{\bar1}d^{\bar0}=d^{\bar0}d^{\bar1}=0.
\]
Its quasi-isomorphisms are the maps inducing isomorphisms on both
cohomology groups. Inverting them gives the derived
$\infty$-category $\mathcal D_2(\F)$. Put
$\Lambda=\F[u,u^{-1}]$, with $|u|=2$ and $du=0$, and let
$\mathcal D(\Lambda)$ be the  derived $\infty$-category
of differential graded $\Lambda$-modules. 
Equivalently,
$\mathcal D_2(\F)\simeq\mathcal D(\Lambda)$: unfolding a
two-periodic complex $K$ puts $K^{\bar n}$ in every degree
$n\in\Z$, and
multiplication by $u$ identifies degrees $n$ and $n+2$.
A map is a quasi-isomorphism precisely when its underlying
unfolded complex of $\F$-vector spaces is. The mapping-complex
and splitting descriptions above apply with degrees read modulo $2$.
The ordinary homotopy category is denoted $D_2(\F)$.

For an integer-graded complex $K$, its direct-sum
\emph{periodization} is
\[
 (\operatorname{Per}_2K)^{\bar\epsilon}
   =\bigoplus_{m\in\Z}K^{\epsilon+2m},\qquad \epsilon=0,1,
\]
with differential induced by $d_K$.
It is left adjoint to the functor which unfolds a two-periodic
complex into an integer-graded complex. Both functors preserve
quasi-isomorphisms, so the adjunction descends to the derived
$\infty$-categories. For complexes of sheaves use the same
formula. Reindexing gives canonical equivalences
\[
 \operatorname{Per}_2(K[2e])\simeq\operatorname{Per}_2K
 \qquad(e\in\Z).
\]
All Koszul signs in parity grading use degrees modulo $2$.

\emph{Sheaf complexes and perfect duality.}
For a fixed topological space $T$, let
$\mathcal D(T;\F)$ be the derived $\infty$-category obtained
from cochain complexes of sheaves of $\F$-vector spaces on $T$
by inverting quasi-isomorphisms; write $\mathcal D_2(T;\F)$
for its two-periodic counterpart, denoted
$\mathcal D_{\mathrm{per}}(T;\bQ)$ in
Subsection~\ref{subsec:coherent-descent}. For $C\in\mathcal D(\F)$,
write $\underline C_T$ for the corresponding constant sheaf
complex; use the same notation in parity grading.

The unfolding construction applies to sheaves as well, identifying
two-periodic sheaf complexes with sheaves of differential graded
$\Lambda$-modules. The proof of
\cite[Proposition~7.1]{ScholzeSixFunctors} applies to these sheaf
modules: the free module $\Lambda$, extended by zero from each
open subset, replaces the rational free sheaf, and the augmented
hypercover resolutions remain exact after tensoring with $\Lambda$.
Thus $\mathcal D_2(T;\F)$ is the category of hypersheaves valued
in $\mathcal D(\Lambda)$. Hyperdescent of module objects is
detected after forgetting their action, since the forgetful
functor preserves limits and detects equivalences. In this model
direct-sum periodization is extension of scalars,
$K\mapsto K\otimes_{\F}\Lambda_T$, where
$\Lambda_T:=\underline{\Lambda}_T$ is the constant sheaf complex
with value $\Lambda$.

A complex $C\in\mathcal D(\F)$
is \emph{perfect} if it is equivalent to a bounded cochain complex
of finite-dimensional $\F$-vector spaces; equivalently,
\[
 \sum_{j\in\Z}\dim_\F H^j(C)<\infty.
\]
A two-periodic complex
$C=(C^{\bar0}\xrightarrow{d}C^{\bar1}\xrightarrow{d}C^{\bar0})$,
with $d^2=0$, is \emph{perfect} if it is equivalent to one with
both terms finite-dimensional; equivalently,
$\dim_\F H^{\bar0}(C)+\dim_\F H^{\bar1}(C)<\infty$.
In either grading put
$C^\vee:=R\operatorname{Hom}_\F(C,\F)$, with $\F$ in degree zero
or even parity, respectively.
For perfect $C,D$, the evaluation and tensor-duality maps are
equivalences:
\[
 C\xrightarrow{\sim}C^{\vee\vee},
 \qquad
 C^\vee\otimes D^\vee\xrightarrow{\sim}(C\otimes D)^\vee.
\]
Indeed, splitting cycles and boundaries over $\F$ reduces these
assertions to finite-dimensional graded vector spaces.
For a space $T$ with dualizing complex $\omega_T$, perfectness
also gives
\[
 \mathbb D_T(\underline C_T)
 :=R\mathcal Hom_T(\underline C_T,\omega_T)
 \simeq\omega_T\otimes\underline{C^\vee}_T.
\]
In the two-periodic category $\omega_T$ denotes its periodization.
Perfectness of stalk complexes does not assert that the sheaf
complex is locally constant on the whole space.

\subsection{Coherent diagrams}
\label{app:coherent-diagrams}


An ordinary category is \emph{small} if its objects and morphisms
form sets in the chosen Grothendieck universe. For a small
ordinary category $I$, its \emph{nerve} $N(I)$ is the simplicial
set whose $m$-simplices are strings of $m$ composable arrows
of $I$. Categories of sheaves and complexes may be large, but
we apply the same definition using larger Grothendieck universes;
see \cite[\S 1.2.15]{HTT}.

A \emph{coherent diagram} in a derived
$\infty$-category $\mathcal C$ is a functor
$N(I)\to\mathcal C$; a morphism of coherent diagrams is a
natural transformation. Here ``coherent'' concerns the diagram
and its compatibility data, not coherent algebraic sheaves.
We abbreviate their $\infty$-category
$\operatorname{Fun}(N(I),\mathcal C)$ to
$\operatorname{Fun}(I,\mathcal C)$.

For a functor $F:I\to J$ between small ordinary categories
and an object $j\in J$, the \emph{comma category}
$(j\downarrow F)$ has objects $(i,u)$ with $i\in I$ and
$u:j\to F(i)$; a morphism $(i,u)\to(i',u')$ is an arrow
$h:i\to i'$ satisfying $F(h)u=u'$.
We call $F$ \emph{cofinal} if the geometric realization of
$N(j\downarrow F)$ is contractible for every $j\in J$.
This is cofinality of the induced map $N(F)$ in the
homotopical sense of \cite[Theorem~4.1.3.1]{HTT}.

In the dg model of $\mathcal D(\F)$, a coherent diagram
assigns a complex $K_i$ to each object and a cochain map
$f_a:K_i\to K_j$ to each arrow $a:i\to j$, together with
specified homotopies expressing compatibility with composition
and all their higher compatibilities. For composable arrows
$a:i\to j$ and $b:j\to k$, this includes a degree-$-1$ map
$H_{b,a}$ satisfying
\[
 d(H_{b,a})=d_{K_k}H_{b,a}+H_{b,a}d_{K_i}
          =f_bf_a-f_{ba}.
\]
For three composable arrows, the two resulting homotopies
between the triple composite and its assigned map are themselves
related by a specified  homotopy of degree $-2$. Longer strings
require the corresponding higher homotopies; identity arrows
carry the identity compatibility data. A morphism of coherent
diagrams likewise includes compatibility homotopies at every
level. These data distinguish a coherent diagram from a diagram
specified only in the ordinary derived category.

\subsection{Representatives of diagrams}
\label{app:diagram-representatives}

For a topological space $T$, let $\operatorname{Sh}(T;\bQ)$
be the abelian category of sheaves of $\bQ$-vector spaces and put
\[
 \mathcal M_T:=\operatorname{Ch}(\operatorname{Sh}(T;\bQ)),
 \qquad
 \mathcal M_{T,\mathrm{per}}:=
       \operatorname{Ch}_{\mathrm{per}}(\operatorname{Sh}(T;\bQ)).
\]
These are ordinary categories of unbounded cochain complexes
and two-periodic cochain complexes, respectively, with cochain
maps as morphisms. Differentials have degree $1$, read modulo
$2$ in the periodic case.
For either choice $\mathcal M$ and a small ordinary category
$A$, write
\[
 \mathcal M^A:=\operatorname{Fun}(A,\mathcal M)
\]
for the ordinary category of \emph{strict $A$-indexed diagrams}.
An object $K$ assigns complexes to objects and cochain maps to
arrows, with
\[
 K(gf)=K(g)\circ K(f),\qquad
 K(\mathrm{id}_a)=\mathrm{id}_{K(a)}
\]
as equalities of cochain maps for every composable pair
$a\xrightarrow{f}b\xrightarrow{g}c$ in $A$.
A morphism $u:K\to K'$ is an ordinary natural transformation:
its components are cochain maps satisfying
\[
 K'(f)\circ u_a=u_b\circ K(f)
 \qquad(f:a\to b).
\]
Let $W_{\mathrm{qis}}$ be the quasi-isomorphisms in $\mathcal M$
and $W_{\mathrm{obj}}$ the objectwise quasi-isomorphisms in
$\mathcal M^A$. Explicitly, $u\in W_{\mathrm{obj}}$ means that
\[
 \mathcal H^n(u_a):\mathcal H^n(K(a))
       \xrightarrow{\sim}\mathcal H^n(K'(a))
\]
is an isomorphism of sheaves for every $a\in A$ and every
$n\in\Z$, or $n\in\Z/2$ in the periodic case.

For an ordinary category $\mathcal M$ and a class $W$ of its
arrows, $N(\mathcal M)[W^{-1}]$ denotes its
\emph{$\infty$-categorical localization}: it is universal for
functors from $N(\mathcal M)$ to $\infty$-categories that send
the arrows of $W$ to equivalences. This construction preserves
higher homotopy information rather than passing to an ordinary
category. Here $N$ is the nerve defined in
Subsection~\ref{app:coherent-diagrams}.

\begin{proposition}[Diagram-model comparison]
\label{prop:diagram-model-comparison}
For every topological space $T$ and small ordinary category $A$,
objectwise localization induces equivalences of
$\infty$-categories
\[
 \begin{aligned}
 N(\mathcal M_T^A)[W_{\mathrm{obj}}^{-1}]
 &\simeq\operatorname{Fun}(N(A),\mathcal D(T;\bQ)),\\
 N(\mathcal M_{T,\mathrm{per}}^A)[W_{\mathrm{obj}}^{-1}]
 &\simeq\operatorname{Fun}(N(A),
                          \mathcal D_{\mathrm{per}}(T;\bQ)).
 \end{aligned}
\]
In each line $W_{\mathrm{obj}}$ is taken in the category on that
line. For $A$ the one-object category with only its identity,
these identify $N(\mathcal M_T)[W_{\mathrm{qis}}^{-1}]$ with
$\mathcal D(T;\bQ)$ and
$N(\mathcal M_{T,\mathrm{per}})[W_{\mathrm{qis}}^{-1}]$ with
$\mathcal D_{\mathrm{per}}(T;\bQ)$.
\end{proposition}
\begin{proof}
The category $\operatorname{Sh}(T;\bQ)$ is Grothendieck.
By \cite[Propositions~1.3.5.3 and~1.3.5.15]{HA},
$\mathcal M_T$ has the combinatorial injective model structure
whose cofibrations are degreewise monomorphisms and whose weak
equivalences are quasi-isomorphisms; its underlying
$\infty$-category is $\mathcal D(T;\bQ)$.

For periodic complexes use the unfolded sheaf differential
graded $\Lambda$-module model of
the paragraph on sheaf complexes in
Subsection~\ref{app:complexes-duality}, where
$\Lambda=\bQ[u,u^{-1}]$, $|u|=2$, and $du=0$.
Apply \cite[Proposition~5.1(1)]{LowenVanDenBergh} to
$\operatorname{Sh}(T;\bQ)$ and the differential graded category
with one object and endomorphism algebra $\Lambda$.
It gives the injective model structure with monomorphisms as
cofibrations and quasi-isomorphisms as weak equivalences.
The construction in \cite[Proposition~5.2 and its proof]{LowenVanDenBergh}
is cofibrantly generated. Its underlying category is Grothendieck,
hence locally presentable, so this model structure is
combinatorial. Unfolding identifies its weak equivalences with
the periodic quasi-isomorphisms used above.

In either case, give $\mathcal M^A$ its injective diagram model
structure, whose cofibrations and weak equivalences are tested
objectwise. Every diagram is cofibrant: its map from the zero
diagram is an objectwise monomorphism. Therefore the localization
of cofibrant diagrams in \cite[Proposition~1.3.4.25]{HA}
is exactly $N(\mathcal M^A)[W_{\mathrm{obj}}^{-1}]$.
That proposition gives both asserted equivalences.
\end{proof}

\begin{proposition}[Simultaneous representatives]
\label{prop:simultaneous-diagram-representatives}
Let $(\mathcal M,\mathcal D)$ be either
$(\mathcal M_T,\mathcal D(T;\bQ))$ or
$(\mathcal M_{T,\mathrm{per}},\mathcal D_{\mathrm{per}}(T;\bQ))$,
and let $\ell:N(\mathcal M)\to\mathcal D$ be the localization
functor. Let $A$ be a small ordinary category.
\begin{enumerate}[(i),wide=0pt]
\item Every coherent diagram $F:N(A)\to\mathcal D$ admits
a strict representative $K:A\to\mathcal M$ together with an
equivalence
\[
 \eta:\ell\circ N(K)\xrightarrow{\sim}F
 \quad\text{in }\operatorname{Fun}(N(A),\mathcal D).
\]
\item A natural transformation $\alpha:F\to G$ can be
represented jointly with its source and target by strict
diagrams $K,K':A\to\mathcal M$ and an ordinary natural
transformation $u:K\to K'$. The equivalences
$\eta_F:\ell\circ N(K)\simeq F$ and
$\eta_G:\ell\circ N(K')\simeq G$ are part of an equivalence
of the corresponding $A\times[1]$-diagrams, where
$[1]=\{0\to1\}$. In particular, this equivalence includes the homotopy
\[
 \eta_G\circ\ell(u)\simeq\alpha\circ\eta_F
\]
and all its compatibilities with the arrows of $A$; here
$\ell(u)$ denotes the natural transformation induced by $u$
after localization.
\item More generally, for any small ordinary category $J$,
a coherent $J$-indexed diagram of coherent $A$-indexed diagrams
admits such a simultaneous representative indexed by $A\times J$.
\end{enumerate}
These assertions allow replacement of the complexes; they do not
assert existence of the representatives with arbitrarily
prescribed cochain models at the vertices.
\end{proposition}
\begin{proof}
Part (i) is essential surjectivity of
Proposition~\ref{prop:diagram-model-comparison}.
For (iii), the nerve preserves products, and the functor-category
equivalence
\[
 \operatorname{Fun}\bigl(N(J),
             \operatorname{Fun}(N(A),\mathcal D)\bigr)
 \simeq\operatorname{Fun}(N(A\times J),\mathcal D)
\]
identifies the given family with one coherent diagram.
Apply (i) to $A\times J$. Taking $J=[1]$ gives (ii),
including the equivalence of the entire arrow diagrams.
\end{proof}

This replacement of a coherent diagram by a strict one is
\emph{strictification}. In the chosen strict diagram, the internal
composition homotopies and higher composition components can all
be taken to be zero. This does not assert vanishing of the homotopy
components of the specified equivalence $\eta$ with the original
diagram.

\subsection{Homotopy limits, colimits, and cofibers}
\label{app:homotopy-limits-colimits}

\begin{definition}[Homotopy limits]
\label{def:general-holim} Let \(I\) be a small category and let \(K:N(I)\to\mathcal D(\mathbb Q)\) be a coherent diagram indexed by \(I\).
 A \emph{coherent cone} from a complex $L$
to $K$ is a morphism $\operatorname{const}_I L\to K$ of
coherent diagrams, where $\operatorname{const}_I L$ assigns
$L$ to every object and the identity to every arrow of $I$.
Thus a cone includes maps $L\to K_i$ and their compatibility
homotopies at every level. The \emph{homotopy limit}, denoted
$H=\operatorname*{holim}_{i\in I}K_i$, is a complex with a
coherent cone to $K$ such that composition induces natural
equivalences
\[
 \operatorname{Map}_{\mathcal D(\bQ)}(L,H)
 \simeq
 \operatorname{Map}_{\operatorname{Fun}(I,\mathcal D(\bQ))}
       (\operatorname{const}_I L,K).
\]
Here $\operatorname{Fun}(I,\mathcal D(\bQ))$ is the
$\infty$-category of coherent $I$-diagrams, so the right-hand
side is the space of coherent cones from $L$. Equivalently,
$H$ is the limit of $K$ in $\mathcal D(\bQ)$.
The same universal property in the $\infty$-category of spaces defines homotopy limits of spaces, including mapping spaces.
\end{definition}

The  derived $\infty$-category $\mathcal D(\bQ)$
admits all small limits; see
\cite[Proposition~1.3.5.21]{HA}. With its cone, a homotopy
limit is unique up to a contractible space of choices. The
construction is functorial in coherent morphisms of diagrams
and takes objectwise equivalences to equivalences. These
properties follow from its universal property.

\begin{definition}[Homotopy colimits]\label{def:general-hocolim}
For a small category $I$ and a coherent diagram $K$ in $\mathcal D(\bQ)$,
its \emph{homotopy colimit}
$\hocolim_{i\in I}K_i$ is its colimit in $\mathcal D(\bQ)$.
It is characterized, together with its structure maps, by natural
equivalences
\[
 \operatorname{Map}_{\mathcal D(\bQ)}
       \!\left(\hocolim_{i\in I}K_i,L\right)
 \simeq
 \operatorname*{holim}_{i\in I^{\mathrm{op}}}
          \operatorname{Map}_{\mathcal D(\bQ)}(K_i,L).
\]
Here the right-hand side is the space of families of maps
$K_i\to L$ together with their compatibility homotopies at every
level. The homotopy limit of mapping spaces is understood as in
Definition~\ref{def:general-holim}. Thus a map from the homotopy
colimit to $L$ is equivalent to such a coherent family of maps
to $L$, including its homotopies at every level.
\end{definition}

For a cofinal functor $F:I\to J$ as defined in
Subsection~\ref{app:coherent-diagrams} and a coherent diagram
$K:N(J)\to\mathcal D(\bQ)$, the canonical map
\[
 \hocolim_{i\in I}K_{F(i)}
       \xrightarrow{\sim}\hocolim_{j\in J}K_j
\]
is an equivalence by \cite[Proposition~4.1.1.8]{HTT}.

\emph{Cofiber sequences of complexes.}
For a morphism $f:K\to L$ in $\mathcal D(\F)$, its \emph{cofiber}, denoted
$\operatorname{cofib}(f)$, is the homotopy pushout of
$L\leftarrow K\to0$. Thus
\[
 \operatorname{cofib}(f):=L\amalg_K^h0
                    \simeq\operatorname{Cone}(f).
\]
Here, on the right, choose a cochain-map representative of $f$;
its mapping cone is
\[
 \operatorname{Cone}(f)^q=L^q\oplus K^{q+1},\qquad
 d(l,k)=(d_Ll+f(k),-d_Kk).
\]
The resulting derived object is independent, up to equivalence,
of the representative. If $f$ is an inclusion of subcomplexes,
the natural map $\operatorname{Cone}(f)\to L/K$ is a
quasi-isomorphism, so $\operatorname{cofib}(f)\simeq L/K$.
For a morphism $f:K\to L$, its \emph{homotopy fiber} is
$\operatorname{fib}(f):=\operatorname{cofib}(f)[-1]$; it fits into
the fiber sequence $\operatorname{fib}(f)\to K\to L$.
A \emph{cofiber sequence} $K\xrightarrow{f}L\xrightarrow{g}M$
includes a specified null-homotopy $h$ of $gf$, satisfying
$gf=d_Mh+hd_K$, such that the induced cochain map
\[
 \operatorname{Cone}(f)\longrightarrow M,\qquad
 (l,k)\longmapsto g(l)+h(k),
\]
is a quasi-isomorphism. Here $h:K^q\to M^{q-1}$ has degree
$-1$. Equivalently, $M$ is the homotopy pushout of
$L\leftarrow K\to0$, with the maps and null-homotopy just
specified. A short exact sequence of complexes
$0\to K\to L\to M\to0$ gives such a cofiber sequence,
with the zero null-homotopy. A \emph{cofiber sequence of
diagrams} is a coherently compatible family of these sequences:
the maps and null-homotopies are compatible with the arrows
of the indexing category and their higher compositions.

For a diagram $F:J\to I$ in an ordinary category $I$, a
\emph{cocone} consists of an object $i\in I$ and arrows
$\lambda_j:F(j)\to i$ such that
\[
 \lambda_{j'}\circ F(a)=\lambda_j
       \qquad\text{for every arrow }a:j\to j'\text{ in }J.
\]
Thus all objects of the diagram map compatibly to one object
$i$. Requiring this for every finite diagram means allowing
every finite category $J$ and every functor $F:J\to I$.
\begin{proposition}\label{prop:hocolim-properties}
Every small coherent diagram in $\mathcal D(\bQ)$ has a homotopy
colimit. With its structure maps it is unique up to a contractible
space of choices, and the construction is functorial in coherent
morphisms of diagrams. It has the following properties.
\begin{enumerate}[(i)]
\item A coherent morphism which is a quasi-isomorphism at every
object induces an equivalence of homotopy colimits.
\item Homotopy colimits preserve cofiber sequences of diagrams
and commute with coproducts, iterated homotopy colimits, and tensoring
with a fixed complex over $\bQ$.
\item If $I$ is filtered, then
\[
 H^q\!\left(\hocolim_{i\in I}K_i\right)
       \cong\varinjlim_{i\in I}H^q(K_i).
\]
For a diagram of actual complexes and cochain maps indexed by $I$,
its degreewise ordinary colimit represents the homotopy colimit.
In particular, restriction of a sequence to any cofinal subsequence
does not change its homotopy colimit.
\end{enumerate}
Here filtered means that $I$ is nonempty and every finite diagram
in $I$ admits a cocone.
\end{proposition}

\begin{proof}[Guide to proof]
The existence assertion is the small-colimit property of the
 derived category, proved in
\cite[Proposition~1.3.5.21]{HA}. The universal property in
Definition~\ref{def:general-hocolim} gives uniqueness, functoriality,
and (i). A cofiber is a pushout over zero, so commutation of colimits
gives the cofiber assertion and the other colimit assertions in (ii).
Tensoring with a fixed complex is left adjoint to the derived mapping
complex functor and hence also preserves colimits.

For (iii), filtered colimits of vector spaces are exact. Consequently
degreewise filtered colimits preserve quasi-isomorphisms and compute
the derived colimit, which agrees with the above homotopy colimit
by \cite[Propositions~1.3.4.24--1.3.4.25]{HA}.
Exactness also gives the displayed cohomology formula, with no
boundedness assumption. Ordinary filtered colimits are unchanged
on cofinal subsequences, proving the last assertion.
\end{proof}

The same limit, colimit, and cofiber constructions apply in
$\mathcal D_2(\F)$. The proofs using split complexes and exact
filtered colimits work in each parity. Since
$\operatorname{Per}_2$ is a left adjoint, it preserves homotopy
colimits and cofibers. In particular, it preserves the telescope
models and their finite-stage filtrations constructed in
Subsection~\ref{app:filtered-diagrams}.


\subsection{\v Cech nerves, totalization, and realization}
\label{app:totalization-realization}

\begin{definition}[\v Cech nerve]
\label{def:cech-nerve}
Let $\Delta$ be the category of nonempty finite ordered sets
$[n]=\{0,\ldots,n\}$ and order-preserving maps, and write
$\Delta^{\mathrm{op}}$ for its opposite category.
For a continuous map $u:U\to S$, its \emph{augmented \v Cech
nerve} is the simplicial space $U_\bullet:\Delta^{\mathrm{op}}
\to\mathrm{Top}$, with its map to the constant simplicial space
$S$, given by
\[
 U_n=\underbrace{U\times_S\cdots\times_S U}_{n+1\text{ factors}}.
\]
The face map $d_i:U_n\to U_{n-1}$ omits the $i$th factor,
and the degeneracy $s_i:U_n\to U_{n+1}$ repeats that factor
by a diagonal. The augmentation is induced by $u$.
For an atlas $U=\coprod_{i\in I}U_i\to S$, one has
\[
 U_n=\coprod_{(i_0,\ldots,i_n)\in I^{n+1}}
        U_{i_0}\times_S\cdots\times_S U_{i_n}.
\]
Repeated indices are included. A finite atlas gives finitely
many terms in each degree, but the nerve has all degrees
$n\ge0$. The same definition applies to morphisms of schemes,
using scheme-theoretic fiber products.
\end{definition}

\begin{definition}[Homotopy totalization and realization]
\label{def:totalization-realization}
Let $\mathcal C$ be $\mathcal D(\F)$ or
$\mathcal D(T;\F)$ for a fixed topological space $T$, or one
of their two-periodic counterparts from
Subsection~\ref{app:complexes-duality}. These categories admit
all small limits and colimits. Let $\Delta$ be the category of
Definition~\ref{def:cech-nerve}.

A coherent \emph{cosimplicial object} of $\mathcal C$ is a
functor $C^\bullet:\Delta\to\mathcal C$. Its
\emph{homotopy totalization} is the limit in $\mathcal C$:
\[
 \Tot(C^\bullet)
   :=\operatorname*{holim}_{[q]\in\Delta}C^q
   \in\mathcal C.
\]
A coherent \emph{simplicial object} of $\mathcal C$ is a
functor $E_\bullet:\Delta^{\mathrm{op}}\to\mathcal C$.
Its \emph{homotopy realization} is the colimit in
$\mathcal C$:
\[
 |E_\bullet|
   :=\hocolim_{[q]\in\Delta^{\mathrm{op}}}E_q
   \in\mathcal C.
\]
Thus these constructions are functors
\[
 \Tot:\operatorname{Fun}(\Delta,\mathcal C)\to\mathcal C,
 \qquad
 |-|:\operatorname{Fun}(\Delta^{\mathrm{op}},\mathcal C)
          \to\mathcal C.
\]
Their inputs include the composition homotopies and higher
compatibilities of the diagrams. In a sheaf category, all terms
and their limit or colimit are complexes on the same space $T$,
and the structure maps are morphisms in that category.
\end{definition}

\begin{lemma}[Simultaneous cochain representatives after
totalization and realization]
\label{lem:totalization-realization-representatives}
Let $\mathcal C$ be one of the derived $\infty$-categories in
Definition~\ref{def:totalization-realization}, and let
$\mathcal M$ be the corresponding ordinary category of cochain
complexes, in the same grading and, for sheaves, on the same
fixed space. Write
$\ell:N(\mathcal M)\to\mathcal C$ for the localization functor.

Let $J$ be a small ordinary category, and let
\[
 \mathscr C:N(J)\longrightarrow
       \operatorname{Fun}(\Delta,\mathcal C),
 \qquad
 \mathscr E:N(J)\longrightarrow
       \operatorname{Fun}(\Delta^{\mathrm{op}},\mathcal C)
\]
be coherent families of cosimplicial and simplicial objects,
respectively. There exist strict diagrams of cochain complexes
$\mathsf T,\mathsf R:J\to\mathcal M$ and equivalences
\[
 \ell\circ N(\mathsf T)\simeq\Tot\circ\mathscr C,
 \qquad
 \ell\circ N(\mathsf R)\simeq|-|\circ\mathscr E
\]
in $\operatorname{Fun}(N(J),\mathcal C)$.
These are equivalences of entire coherent diagrams, including
their comparison morphisms and higher compatibilities.

In particular, taking $J=[1]=\{0\to1\}$, coherent natural
transformations
\[
 \alpha:C^\bullet\longrightarrow D^\bullet,
 \qquad
 \beta:E_\bullet\longrightarrow F_\bullet
\]
induce morphisms
\[
 \Tot(\alpha):\Tot(C^\bullet)\longrightarrow\Tot(D^\bullet),
 \qquad
 |\beta|:|E_\bullet|\longrightarrow|F_\bullet|
\]
that admit cochain representatives after choosing their source
and target complexes jointly. If $\alpha$ or $\beta$ is an
equivalence in every simplicial degree, the corresponding
cochain representative is a quasi-isomorphism.

All these choices allow replacement of the complexes; no
assertion is made for arbitrarily prescribed cochain models.
\end{lemma}

\begin{proof}
By Definition~\ref{def:totalization-realization},
$\Tot$ and $|-|$ are functors of $\infty$-categories.
Consequently $\Tot\circ\mathscr C$ and $|-|\circ\mathscr E$
are coherent $J$-indexed diagrams in $\mathcal C$.
Apply Proposition~\ref{prop:simultaneous-diagram-representatives}(i)
to these two diagrams to obtain $\mathsf T$, $\mathsf R$,
and the asserted equivalences.

An objectwise equivalence of cosimplicial or simplicial objects
is an equivalence in the corresponding functor
$\infty$-category. Functors preserve equivalences, so its
totalization or realization is an equivalence in $\mathcal C$.
A cochain map represents such an equivalence precisely when it
is a quasi-isomorphism.
\end{proof}

For either strict diagram $K=\mathsf T$ or $K=\mathsf R$,
every arrow $a:j\to j'$ of $J$ is represented by a cochain map
satisfying
\[
 d_{K(j')}\circ K(a)=K(a)\circ d_{K(j)}.
\]
Whenever these representatives are computed by \v Cech
total-complex models, the differentials in this identity are
the full total differentials, including the alternating maps
between overlap terms. The lemma requires the original
comparisons to be coherent morphisms of diagrams, with their
specified homotopies and higher compatibilities; independently
chosen maps on the individual terms do not supply this input.

\subsection{Filtered objects and telescope models}
\label{app:filtered-diagrams}

The definitions and results of this subsection also hold in
$\mathcal D_2(\F)$, with the same formulas and proofs read
modulo $2$.

For a sequence, $0\to1\to2\to\cdots$ denotes the
\emph{indexing category}: its objects are the nonnegative
integers, and it has one arrow $i\to j$ whenever $i\le j$.
The integers here label complexes; they are not themselves
complexes or cochain groups. Choose a strict representative
of a coherent sequential diagram: this consists of actual
cochain complexes $K_i$ and specified cochain maps
$u_i:K_i\to K_{i+1}$. The map assigned to $i\to j$, for $i<j$,
is then $u_{j-1}\circ\cdots\circ u_i$, and the identity arrow
is assigned the identity cochain map. Thus, once the successive
cochain maps are chosen, the maps for all longer intervals are
their composites and need no separate choices.

\begin{definition}[Filtered objects]
\label{def:filtered-derived-object}
A nonnegative increasing \emph{filtered object} over $\F$ is an
object $C_\bullet$ of
\[
 \mathcal D_{\mathrm{fil}}(\F)
   :=\operatorname{Fun}(\Z_{\ge0},\mathcal D(\F));
\]
thus it is a coherent diagram $C_0\to C_1\to\cdots$.
An \emph{exhaustive filtration} of $C$ additionally specifies an
equivalence $\hocolim_PC_P\simeq C$. We use $C_P$ for these
derived stages.
A morphism of filtered objects is a natural transformation of
the coherent diagrams: it includes the homotopies with all
transition maps and their higher compatibilities, as in
Subsection~\ref{app:coherent-diagrams}. It is an equivalence in
$\mathcal D_{\mathrm{fil}}(\F)$ precisely when each component is
an equivalence in $\mathcal D(\F)$.

An \emph{exhaustive cochain model} of $C_\bullet$ consists of a
cochain complex $T$, an increasing filtration by subcomplexes
\[
 0=F_{-1}T\subseteq F_0T\subseteq F_1T\subseteq\cdots,
 \qquad T=\bigcup_{P\ge0}F_PT,
\]
and a specified equivalence of entire diagrams
\[
 (F_PT)_{P\ge0}\simeq(C_P)_{P\ge0}
 \quad\text{in }\mathcal D_{\mathrm{fil}}(\F).
\]
The terms $F_PT$ in this equivalence are viewed in
$\mathcal D(\F)$. The equivalence includes compatibility with
all transition morphisms and their higher homotopies; a
collection of unrelated stagewise equivalences does not specify
it. Passing to homotopy colimits gives $T\simeq C$ when the
exhaustive filtration of $C$ has been specified. We use $F_PT$
for the actual subcomplex of the chosen model, not as an
alternative name for the derived stage $C_P$.

Put $C_{-1}=0$ and define
\[
 \gr_PC_\bullet:=\operatorname{cofib}(C_{P-1}\to C_P),
 \qquad
 \gr_PT:=F_PT/F_{P-1}T.
\]
The specified diagram equivalence induces
$\gr_PT\simeq\gr_PC_\bullet$ in $\mathcal D(\F)$.
\end{definition}

For these nonnegative filtrations, a morphism is an equivalence
if and only if it induces equivalences on all associated-graded
complexes, by induction on the cofiber sequences.

A homotopy of filtered morphisms is a homotopy of entire diagrams.
Thus equality in $h\mathcal D_{\mathrm{fil}}(\F)$ means that the
two morphisms are joined by such a homotopy. Separate homotopies
at individual stages do not supply this compatibility.

The following telescope construction gives an exhaustive
filtration by subcomplexes representing any filtered object,
after strictification of its sequential diagram by
Proposition~\ref{prop:simultaneous-diagram-representatives}(i).

\begin{proposition}\label{prop:sequential-telescope}
Let $K_0\xrightarrow{u_0}K_1\xrightarrow{u_1}\cdots$ be a
sequence of cochain complexes with specified cochain maps
$u_p:K_p\to K_{p+1}$: each $u_p$ has degree zero and commutes
with the differentials. The telescope below uses these actual
maps. Let $s$ send the $p$-th summand of
$\bigoplus_{p\ge0}K_p$ to the $(p+1)$-st by $u_p$. Then the
\emph{mapping telescope}
\[
 \operatorname{Tel}(K_\bullet)
  :=\operatorname{Cone}\!\left(1-s:\bigoplus_{p\ge0}K_p
                                 \longrightarrow\bigoplus_{p\ge0}K_p\right)
\]
represents $\hocolim_pK_p$. It has an exhaustive nonnegative
filtration by the subcomplexes
\[
 T_P:=\operatorname{Cone}\!\left(1-s:\bigoplus_{0\le p<P}K_p
                             \longrightarrow\bigoplus_{0\le p\le P}K_p\right).
\]
Writing $T:=\operatorname{Tel}(K_\bullet)$ and $F_PT:=T_P$
therefore gives a filtration by actual subcomplexes. There are
natural quasi-isomorphisms $T_P\to K_P$ commuting with the
transition maps; together they give an equivalence of filtered
diagrams $(F_PT)_P\simeq(K_P)_P$. No injectivity of the maps
$u_P$ is required.
\end{proposition}

\begin{proof}[Guide to proof]
On the summand $K_p$ in
$\bigoplus_{p\le P}K_p$, take the composite transition map to
$K_P$, including the identity for $p=P$. These maps give an
augmentation and an exact sequence of complexes
\[
 0\longrightarrow\bigoplus_{p<P}K_p
 \xrightarrow{1-s}\bigoplus_{p\le P}K_p
 \longrightarrow K_P\longrightarrow0.
\]
Its exactness follows successively from the components $p=0,1,\ldots,P$.
The resulting map $T_P\to K_P$, zero on the shifted source
summands, is therefore a quasi-isomorphism. These maps commute
with $T_P\hookrightarrow T_{P+1}$ and $u_P$.
Moreover $\bigcup_PT_P=\operatorname{Tel}(K_\bullet)$, and
filtered exactness gives a quasi-isomorphism from this union to
$\varinjlim_PK_P$. Proposition~\ref{prop:hocolim-properties}(iii)
identifies the latter with the homotopy colimit. This is the
classical telescope construction of
\cite[Definition~1.6.4]{Neeman}; compare
\cite[Tags~090Z, 093W, 0949, and~0CRK]{Stacks}.
\end{proof}


\begin{lemma}[Passage from stage operations to the filtered complex]
\label{lem:filtered-operations-colimit}
Let $C^{(1)}_\bullet,\ldots,C^{(r)}_\bullet,D_\bullet$ be
filtered objects, with $r\ge1$, and put
$C^{(i)}=\hocolim_PC^{(i)}_P$ and $D=\hocolim_PD_P$.
Suppose that
\[
 \theta_{\mathbf P}:\bigotimes_{i=1}^r C^{(i)}_{P_i}
        \longrightarrow D_{|\mathbf P|},
 \qquad |\mathbf P|=\sum_iP_i,
\]
is a morphism of coherent diagrams indexed by
$\mathbf P=(P_1,\ldots,P_r)\in\Z_{\ge0}^r$.
Define the \emph{sum filtration} of their tensor product to be
the filtered derived object $E_\bullet=(E_N)_{N\ge0}$, where
\[
 E_N:=\hocolim_{|\mathbf P|\le N}
               \bigotimes_i C^{(i)}_{P_i},
\]
the indexing tuples are ordered coordinatewise, and the maps
$E_N\to E_{N+1}$ are induced by inclusions of the indexing
posets. Then $\theta$ induces a morphism $E_\bullet\to D_\bullet$ of
filtered objects and a morphism of their underlying derived
objects
\begin{equation}\label{eq:filtered-operation-colimit}
 \bigotimes_i C^{(i)}
 \simeq\hocolim_{\mathbf P}\bigotimes_iC^{(i)}_{P_i}
 \xrightarrow{\hocolim\theta_{\mathbf P}}
       \hocolim_{\mathbf P}D_{|\mathbf P|}
 \simeq D.
\end{equation}
These constructions preserve coherent comparisons of the stage
operations.

Fix filtered cochain models of $C^{(i)}_\bullet$ and
$D_\bullet$: namely, cochain complexes $T^{(i)}$ and $U$
with increasing filtrations by subcomplexes satisfying
\[
 T^{(i)}=\bigcup_{P\ge0}F_PT^{(i)},
 \qquad U=\bigcup_{P\ge0}F_PU,
\]
together with specified equivalences
\[
 (F_PT^{(i)})_P\simeq(C^{(i)}_P)_P,
 \qquad (F_PU)_P\simeq(D_P)_P
\]
in $\mathcal D_{\mathrm{fil}}(\F)$; in particular one may fix
the finite-initial-telescope models. These are equivalences of
entire diagrams in the sense of
Definition~\ref{def:filtered-derived-object}.
The induced filtered morphism has a cochain representative
\[
 \theta_T:\bigotimes_iT^{(i)}\longrightarrow U,
 \qquad
 \theta_T\Bigl(\bigotimes_iF_{P_i}T^{(i)}\Bigr)
       \subseteq F_{|\mathbf P|}U,
\]
on these same complexes with their fixed differentials.
Its filtered chain-homotopy class depends only on the induced
filtered morphism. The assertions also hold with a fixed complex
$V$ tensored into every source and into the source of
\eqref{eq:filtered-operation-colimit}.
\end{lemma}

\begin{proof}
Compose each $\theta_{\mathbf P}$ with
$D_{|\mathbf P|}\to D_N$ for $|\mathbf P|\le N$.
The assumed coherences give a cocone, hence a map from the
displayed homotopy colimit $E_N$ to $D_N$, coherently in $N$.
Tensor products preserve homotopy colimits separately, by
Proposition~\ref{prop:hocolim-properties}(ii), giving the first
equivalence in \eqref{eq:filtered-operation-colimit}.
The addition functor $\Z_{\ge0}^r\to\Z_{\ge0}$ is cofinal:
for every $N$, the poset of tuples with $|\mathbf P|\ge N$
is nonempty and filtered under coordinatewise maxima.
This gives the last equivalence. The same argument shows that
the sum filtration is exhaustive.

For the specified cochain models, the derived sum filtration
$E_\bullet$ is represented by the actual subcomplexes
\[
 S_N:=F_N\Bigl(\bigotimes_iT^{(i)}\Bigr)
     :=\sum_{|\mathbf P|\le N}
              \bigotimes_iF_{P_i}T^{(i)}
       \ \subseteq\ \bigotimes_iT^{(i)}.
\]
Indeed, the diagram over $|\mathbf P|\le N$ can be built by
adjoining its objects in increasing total degree. At each step,
the union of the proper predecessor terms maps injectively into
the new tensor-product term. The needed intersection identities
follow by choosing degreewise vector-space splittings of the
input filtrations. Each such attachment therefore computes a
homotopy pushout, so its ordinary colimit $S_N$ computes the
required homotopy colimit. These comparisons are natural in
$N$ and give $(S_N)_N\simeq(E_N)_N$ in
$\mathcal D_{\mathrm{fil}}(\F)$.

It remains to represent the filtered morphism on $S_\bullet$
and $F_\bullet U$, without changing these models. A filtration
$A_\bullet$ by subcomplexes which is acyclic at every stage
admits a filtration-preserving contraction. To see this, each
$\gr_PA$ is acyclic and hence a direct sum of two-term identity
complexes of vector spaces. It is therefore projective in the
category of complexes. The sequences
$0\to A_{P-1}\to A_P\to\gr_PA\to0$ split as complexes;
induction identifies the filtration with the partial sums of
$\bigoplus_P\gr_PA$. Contract each summand.
Consequently, a stagewise quasi-isomorphism between filtrations
by subcomplexes has a filtration-preserving homotopy inverse,
since its filtered cone has such a contraction.

Proposition~\ref{prop:diagram-model-comparison}, with
$T=\mathrm{pt}$ and $A=\Z_{\ge0}$, identifies the
$\infty$-category of coherent sequential diagrams with the
$\infty$-categorical localization of the ordinary category of
strict diagrams at stagewise quasi-isomorphisms. Passing to
homotopy categories gives the ordinary localization used in the
following zigzag argument.
Replace intermediate strict diagrams in a representative zigzag
by their telescope filtrations and invert their comparison maps
using the preceding argument. This represents the given
filtered morphism by a cochain map $S_\bullet\to F_\bullet U$.
Filtered mapping cylinders show that homotopic maps become equal
in this localization; conversely, the filtered homotopy category
already inverts all its stagewise quasi-isomorphisms. Hence it is
the same localization, proving uniqueness up to filtered chain
homotopy. Taking unions gives $\theta_T$.
Tensoring with $V$ preserves the colimits, inclusions, and
intersection identities used above, proving the last assertion.

The comparison also preserves higher homotopies on these fixed
models. For exhaustively filtered cochain complexes $T,U$,
define the filtered mapping complex by
\[
 \begin{split}
 \underline{\Hom}_{\mathrm{fil}}^r(T,U)
  :=\{f\in\underline{\Hom}^r(T,U):
       f(F_PT)\subseteq F_PU\text{ for every }P\},\\
 d(f):=d_Uf-(-1)^r f d_T.
 \end{split}
\]
Here $\underline{\Hom}^r(T,U)$ is the vector space of
degree-$r$ graded linear maps defined in
Subsection~\ref{app:complexes-duality}; the containment means
$f((F_PT)^j)\subseteq(F_PU)^{j+r}$ for every $j$.
Let $\mathcal K_{\mathrm{fil}}$ be the dg nerve of the category
of these filtered cochain complexes with these mapping
complexes. Sending a cochain filtration to its derived diagram
defines a functor to $\mathcal D_{\mathrm{fil}}(\F)$ and hence
comparison maps
\[
 \operatorname{Map}_{\mathcal K_{\mathrm{fil}}}(T,U)
 \longrightarrow
 \operatorname{Map}_{\mathcal D_{\mathrm{fil}}(\F)}
                   (F_\bullet T,F_\bullet U).
\]
For $n\ge0$, the map on $\pi_n$ at zero identifies with
\[
 H^{-n}\underline{\Hom}_{\mathrm{fil}}(T,U)
 \longrightarrow
 \Hom_{h\mathcal D_{\mathrm{fil}}(\F)}
                   (F_\bullet T,F_\bullet(U[-n])).
\]
It is an isomorphism by the preceding fixed-model argument
applied to $U[-n]$; for $n=0$ it is exactly the existence and
uniqueness up to filtered chain homotopy already proved.
Translation in the additive mapping spaces gives the same
conclusion at every basepoint. Thus the mapping-space comparison
is an equivalence, not merely a bijection on homotopy classes.
It applies also to the sum-filtered tensor sources and the fixed
parameter complex $V$ used above.

Consequently the maps, comparison homotopies, and higher
compatibilities of a coherent diagram can be represented jointly
on the fixed models. Previously chosen endpoint maps may be
kept fixed: first use their specified comparison paths to the
specified derived maps, and then transport the comparison
homotopies relative to those endpoints. Equivalences of mapping
spaces induce equivalences on the homotopy fibers prescribing
these endpoints and, successively, the boundaries of higher
homotopies. This preserves the required coherences without
asserting that independently chosen cochain representatives
satisfy strict composition identities.
\end{proof}

Thus a product with $r$ inputs requires a coherent morphism on
all input bounds, not merely separately chosen maps
$\bigotimes_iC_{P_i}\to C_{|\mathbf P|}$.
The resulting cochain map on a fixed model represents those maps
under the specified comparisons; it need not restrict literally
to independently chosen representatives of them.

\phantomsection

\end{document}